\documentclass[twoside]{scrartcl}
\usepackage{scrlayer-scrpage}

\usepackage[utf8]{inputenc}
\usepackage[T1]{fontenc}
\usepackage[english]{babel}
\usepackage[unicode=true,plainpages=false]{hyperref}

\usepackage[left=2.5cm, right=3cm, top=2cm, bottom=3.5cm]{geometry}

\usepackage[shortlabels]{enumitem}

\usepackage{amsmath}
\usepackage{amssymb}
\usepackage{amsthm}
\usepackage{mathrsfs}
\usepackage{yhmath}
\usepackage{wasysym}
\usepackage{mathtools}

\usepackage[capitalise]{cleveref}

\usepackage{caption}
\usepackage{subcaption}

\usepackage{graphicx}
\usepackage{pgf,tikz}
\usepackage{tikz-cd}
\usepackage{rotating}

\usepackage{lmodern} \usepackage{tabularx} \usepackage[draft=false,babel,kerning=true,spacing=true]{microtype}

\allowdisplaybreaks

\let\cprod\times

\makeatletter
\newcommand*{\relrelbarsep}{.386ex}
\newcommand*{\relrelbar}{  \mathrel{    \mathpalette\@relrelbar\relrelbarsep
  }}
\newcommand*{\@relrelbar}[2]{  \raise#2\hbox to 0pt{$\m@th#1\relbar$\hss}  \lower#2\hbox{$\m@th#1\relbar$}}
\providecommand*{\rightrightarrowsfill@}{  \arrowfill@\relrelbar\relrelbar\rightrightarrows
}
\providecommand*{\leftleftarrowsfill@}{  \arrowfill@\leftleftarrows\relrelbar\relrelbar
}
\providecommand*{\xrightrightarrows}[2][]{  \ext@arrow 0359\rightrightarrowsfill@{#1}{#2}}
\providecommand*{\xleftleftarrows}[2][]{  \ext@arrow 3095\leftleftarrowsfill@{#1}{#2}}
\makeatother

\renewcommand{\times}{\cdot}
\DeclarePairedDelimiter{\abs}{\lvert}{\rvert}

 \newcommand{\noloc}{\nobreak\mskip6mu plus1mu\mathpunct{}\nonscript\mkern-\thinmuskip{:}\mskip2mu\relax}

\newcommand{\isom}{\cong}

\newcommand{\from}{\mathbin{\leftarrow}}
\newcommand{\xto}[1]{\mathbin{\xrightarrow{#1}}} \newcommand{\xfrom}[1]{\mathbin{\xleftarrow{#1}}} \newcommand{\xinjto}[1]{\mathbin{\xhookrightarrow{#1}}} \newcommand{\isoto}{\xto\sim}

\newcommand{\xtofrom}[1]{\mathbin{\xleftrightarrow{#1}}}
\newcommand{\isotofrom}{\xtofrom\sim}

\DeclareMathOperator{\Hom}{Hom}

\DeclareMathOperator{\Fun}{Fun}
\DeclareMathOperator{\Obj}{Obj}
\DeclareMathOperator*{\colim}{colim}
\DeclareMathOperator{\Tow}{Tow}  \DeclareMathOperator{\Mul}{Mul}

\newcommand{\setst}{\mathrel{|}}
\newcommand{\isect}{\mathbin{\cap}}
\newcommand{\bigisect}{\bigcap}
\newcommand{\union}{\mathbin{\cup}}
\newcommand{\bigunion}{\bigcup}
\newcommand{\dunion}{\mathbin{\sqcup}}
\newcommand{\bigdunion}{\bigsqcup}

\renewcommand{\emptyset}{\varnothing}

\renewcommand{\subset}{\subseteq}
\renewcommand{\supset}{\supseteq}

\newcommand{\comp}{\mathbin{\circ}}
\DeclareMathOperator{\id}{id}
\DeclareMathOperator{\ev}{ev}
\newcommand{\pr}[1]{\mathrm{pr}_{#1}}
\newcommand{\restrict}[2]{{#1}|_{#2}}

\newcommand{\injto}{\mathrel{\hookrightarrow}}
\newcommand{\surjto}{\mathrel{\twoheadrightarrow}}

\DeclareMathOperator{\Map}{Map}

\newcommand{\dsum}{\oplus}
\newcommand{\bsum}{\mathbin{\boxplus}}
\newcommand{\bigdsum}{\bigoplus}

\newcommand{\cd}{\mathrm{cd}}

\newcommand{\Z}{{\mathbb{Z}}}
\newcommand{\N}{\mathbb{N}}
\newcommand{\Fld}{\mathbb{F}}
\newcommand{\Q}{\mathbb{Q}}
\newcommand{\R}{{\mathbb{R}}}

\newcommand{\divides}{\mathbin{|}}

\newcommand{\GL}{\mathrm{GL}}

\DeclareMathOperator{\End}{End}

\DeclareMathOperator{\img}{im}

\DeclareMathOperator{\coker}{coker}

\newcommand{\units}{^\cprod}
\newcommand{\tensor}{\otimes}

\newcommand{\bigtensor}{\bigotimes}

\DeclareMathOperator{\Gal}{Gal}
\DeclareMathOperator{\Sym}{Sym}

\newcommand{\sep}{{\mathrm{sep}}}

\DeclareMathOperator{\Spec}{Spec}

\newcommand{\Spf}{\mathrm{Spf}\,}

\newcommand{\setP}{\mathbb{P}}

\DeclareMathOperator{\Cond}{Cond} \newcommand{\disc}{{\mathrm{disc}}}

\newcommand{\triv}{{\mathrm{triv}}}
\newcommand{\sm}{{\mathrm{sm}}}

\DeclareMathOperator{\catset}{Set}

\DeclareMathOperator{\catsimpmark}{Set^+_\Delta}

\newcommand{\Drepsld}[2]{\D_\solid(#1, #2)}

\newcommand{\Drepsldsm}[2]{\D^{\mathrm{sm}}_\solid(#1, #2)}
\newcommand{\Drepsldasm}[2]{\D^{\mathrm{sm},a}_\solid(#1, #2)}

\newcommand{\Drepsldasmp}[2]{\D^{\mathrm{sm},+a}_\solid(#1, #2)}

\DeclareMathOperator{\Shv}{Shv}
\newcommand{\catsh}[2]{\Shv_{#2}(#1)}
\newcommand{\cathsh}[2]{\Shv_{#2}(#1)\hat{\,}} \newcommand{\catshs}[1]{\operatorname{Shv}(#1)}

\DeclareMathOperator{\IHom}{\underline{\Hom}}\DeclareMathOperator{\smIHom}{\IHom^{\mathrm{sm}}}

 \newcommand{\an}{{\mathrm{an}}}  \newcommand{\et}{{\mathrm{et}}}

\newcommand{\proet}{{\mathrm{proet}}} \newcommand{\qproet}{{\mathrm{qproet}}}  \newcommand{\vsite}{{\mathrm{v}}}    
\newcommand{\cplt}{\hat{}\,}
\newcommand{\ri}{\mathcal O} \newcommand{\mm}{\mathfrak m} 
\DeclareMathOperator{\Spa}{Spa}

\newcommand{\opp}{{\mathrm{op}}}

\newcommand{\solid}{{\scalebox{0.5}{$\square$}}}
\DeclareMathOperator{\Mod}{Mod}

\newcommand{\Perf}{\mathrm{Perf}}
\newcommand{\Perfd}{\mathrm{Perfd}}
\newcommand{\AffPerfd}{\Perfd^{\mathrm{aff}}}
\newcommand{\AffPerfCoeff}{\Perf^{\mathrm{aff}}_\Lambda}

\DeclareMathOperator{\D}{\mathcal D}
\DeclareMathOperator{\Dqcohri}{\D^a_\solid}
\newcommand{\Dqcohriarg}[3]{\operatorname{\D^{#1a#2}_{\solid#3}}}
\newcommand{\Dqcohriq}{\Dqcohriarg?{}{}}
\newcommand{\Dqcohrip}{\Dqcohriarg+{}{}}
\newcommand{\Dqcohrim}{\Dqcohriarg{--}{}{}}
\newcommand{\Dqcohrib}{\Dqcohriarg b{}{}}

\newcommand{\Dqcohriwap}{\Dqcohriarg{}{w+}{}}
\newcommand{\Dqcohriwam}{\Dqcohriarg{}{w--}{}}
\newcommand{\Dqcohriwab}{\Dqcohriarg{}{wb}{}}
\newcommand{\Dqcohrile}[1]{\Dqcohriarg{}{}{\le#1}}
\newcommand{\Dqcohrige}[1]{\Dqcohriarg{}{}{\ge#1}}
\newcommand{\DqcohriX}[1]{\Dqcohri(\ri^+_{#1}/\pi)}
\newcommand{\DqcohriXZ}[1]{\Dqcohri(\ri^+_{#1}/\pi, \Z)}

\DeclareMathOperator{\Tot}{Tot}

\newcommand{\waperf}{{\mathrm{waperf}}}
\newcommand{\perf}{{\mathrm{perf}}}

\newcommand{\infcatinf}{\mathcal Cat_\infty}
\newcommand{\catcat}{\mathcal Cat}
\newcommand{\catS}{\mathcal S} \newcommand{\catPrL}{\mathcal Pr^L}
\newcommand{\catPrR}{\mathcal Pr^R}

\newcommand{\catsldmod}[1]{\mathrm{Mod}_\solid(#1)}
\newcommand{\catsldmoda}[1]{\mathrm{Mod}_\solid^a(#1)}

\DeclareMathOperator{\AnRing}{AnRing}
\DeclareMathOperator{\AnAssRing}{AnAssRing}

\DeclareMathOperator{\AnSpec}{AnSpec}
\DeclareMathOperator{\AnSpace}{AnSpace}
\DeclareMathOperator{\AlmSetup}{AlmSetup}

\DeclareMathOperator{\Ring}{Ring}
\DeclareMathOperator{\AssRing}{AssRing}
\DeclareMathOperator{\LMod}{LMod}
\DeclareMathOperator{\Ani}{Ani}

\newcommand{\oc}{{\mathrm{oc}}}
\DeclareMathOperator{\fib}{fib}
\DeclareMathOperator{\cofib}{cofib}

\newcommand{\catfiltcofin}[2]{\mathrm{Filt}^{#1}_{#2}}
\newcommand{\catPrLk}[1]{\mathrm{Pr}^L_{#1}} \newcommand{\catAcc}[1]{\mathrm{Acc}_{#1}}
\DeclareMathOperator{\Ind}{Ind} \DeclareMathOperator{\Pro}{Pro}

\DeclareMathOperator{\cosk}{cosk}

\DeclareMathOperator{\Corr}{Corr}
\DeclareMathOperator{\vStackspi}{vStack^\sharp_\pi}
\DeclareMathOperator{\vStackspip}{vStack^\sharp_{\pi\divides p}}
\DeclareMathOperator{\vStacksCoeff}{vStack_\Lambda}
\DeclareMathOperator{\vStacksCoeffVar}{vStack^\dagger_\Lambda}

\DeclareMathOperator{\catFinAst}{Fin_\ast}

\DeclareMathOperator{\opComm}{Comm^\tensor}
\DeclareMathOperator{\opAssoc}{Assoc^\tensor}
\DeclareMathOperator{\opLM}{LM^\tensor}
\DeclareMathOperator{\opTriv}{Triv^\tensor}
\newcommand{\act}{{\mathrm{act}}} \newcommand{\cart}{{\mathrm{cart}}}
\newcommand{\StkSch}{{\mathrm{StkSch}}}

\DeclareMathOperator{\trc}{tr.c} \newcommand{\dimtrg}{\mathrm{dim.trg}}
\newcommand{\std}{{\mathrm{std}}} \newcommand{\nuc}{{\mathrm{nuc}}} \newcommand{\lnuc}{{\mathrm{lnuc}}}
\DeclareMathOperator{\eq}{eq}
\DeclareMathOperator{\ExtrDisc}{ExtrDisc} \DeclareMathOperator{\Alg}{Alg} \DeclareMathOperator{\CAlg}{CAlg} 
\newcommand{\equalizer}[4]{\eq\big( \begin{tikzcd}[ampersand replacement=\&] #1 \arrow[r,shift left,"#3"] \arrow[r,shift right,swap,"#4"] \& #2 \end{tikzcd} \big)}

\newcommand{\sigobj}[1]{\text{$\sigma$-$#1$}}

\theoremstyle{plain}
\newtheorem{theorem}{Theorem}[subsection]
\newtheorem{theorem*}{Theorem}
\newtheorem{proposition}[theorem]{Proposition}
\newtheorem{proposition*}[theorem*]{Proposition}
\newtheorem{corollary}[theorem]{Corollary}
\newtheorem{lemma}[theorem]{Lemma}

\theoremstyle{definition}
\newtheorem{definition}[theorem]{Definition}
\newtheorem{definition*}[theorem*]{Definition}
\newtheorem{example}[theorem]{Example}

\newtheorem{remark}[theorem]{Remark}
\newtheorem{remarks}[theorem]{Remarks}
\newtheorem{warning}[theorem]{Warning}

\newtheorem{hypothesis*}[theorem*]{Hypothesis}

\newtheorem{setup}[theorem]{Setup}

\numberwithin{equation}{theorem}
\numberwithin{figure}{subsection}
\numberwithin{table}{subsection}

\newlist{thmenum}{enumerate}{1}
\setlist[thmenum]{label=(\roman*), ref=\thetheorem.(\roman*)}
\crefalias{thmenumi}{theorem}

\newlist{propenum}{enumerate}{1}
\setlist[propenum]{label=(\roman*), ref=\theproposition.(\roman*)}
\crefalias{propenumi}{proposition}

\newlist{corenum}{enumerate}{1}
\setlist[corenum]{label=(\roman*), ref=\thecorollary.(\roman*)}
\crefalias{corenumi}{corollary}

\newlist{lemenum}{enumerate}{1}
\setlist[lemenum]{label=(\roman*), ref=\thelemma.(\roman*)}
\crefalias{lemenumi}{lemma}

\newlist{exampleenum}{enumerate}{1}
\setlist[exampleenum]{label=(\alph*), ref=\theexamples.(\alph*)}
\crefalias{exampleenumi}{example}

\newlist{remarksenum}{enumerate}{1}
\setlist[remarksenum]{label=(\roman*), ref=\theremarks.(\roman*)}
\crefalias{remarksenumi}{remark}

\newlist{defenum}{enumerate}{1}
\setlist[defenum]{label=(\alph*), ref=\thedefinition.(\alph*)}
\crefalias{defenumi}{definition}

\newcommand{\thetitle}{A $p$-Adic 6-Functor Formalism\\in Rigid-Analytic Geometry}
\newcommand{\theauthor}{Lucas Mann}

\title{\thetitle}
\author{\theauthor}
\date{\today}

\begin{document}

\maketitle

\begin{abstract}
We develop a full 6-functor formalism for $p$-torsion étale sheaves in rigid-analytic geometry. More concretely, we use the recently developed condensed mathematics by Clausen--Scholze to associate to every small v-stack (e.g. rigid-analytic variety) $X$ with pseudouniformizer $\pi$ an $\infty$-category $\DqcohriX X$ of ``derived quasicoherent complete topological $\ri^+_X/\pi$-modules'' on $X$. We then construct the six functors $\tensor$, $\IHom$, $f^*$, $f_*$, $f_!$ and $f^!$ in this setting and show that they satisfy all the expected compatibilities, similar to the $\ell$-adic case. By introducing $\varphi$-module structures and proving a version of the $p$-torsion Riemann-Hilbert correspondence we relate $\ri^+_X/\pi$-sheaves to $\Fld_p$-sheaves. As a special case of this formalism we prove Poincaré duality for $\Fld_p$-cohomology on rigid-analytic varieties. In the process of constructing $\DqcohriX X$ we also develop a general descent formalism for condensed modules over condensed rings.
\end{abstract}

\tableofcontents

\section{Introduction}

\subsection{Poincaré Duality in Rigid-Analytic Geometry}

In complex analytic geometry we have the following well-known version of Poincaré duality: Given a compact complex manifold $X$ of dimension $d$ and any field $k$, then for all $i \in \Z$ there is a natural perfect pairing
\begin{align*}
	H^i(X, k) \tensor_k H^{2d-i}(X, k) \to k
\end{align*}
induced by the cup product and a trace map $H^{2d}(X, k) \to k$. Here $H^i(X, k)$ denotes the singular cohomology of $X$ with coefficients in $k$. It is natural to ask whether a similar duality result holds in other contexts. For example, it is a classical result that the above duality holds also in the case that $X$ is a proper smooth variety over some field $K$ and that $k = \Fld_\ell$ for some prime $\ell$ not equal to the characteristic of $K$, where we replace singular cohomology by étale cohomology.

In this thesis we study Poincaré duality in the setting of rigid-analytic geometry, i.e. the $p$-adic analog of complex analytic geometry. Then Poincaré duality takes the following form:

\begin{theorem} \label{rslt:intro-poincare-duality}
Let $K$ be an algebraically closed non-archimedean field of mixed characteristic $(0, p)$ and let $X$ be a proper smooth rigid-analytic variety of pure dimension $d$ over $K$. Then for all primes $\ell$ and all $i \in \Z$ there is a natural perfect pairing
\begin{align*}
	H^i_\et(X, \Fld_\ell) \tensor_{\Fld_\ell} H^{2d-i}_\et(X, \Fld_\ell) \to \Fld_\ell(-d)
\end{align*}
\end{theorem}

Here $\Fld_\ell(-d)$ denotes the $-d$-th Tate twist of $\Fld_\ell$, which parametrizes the choice of an $\ell$-th root of unity in $K$. Unlike in the case of schemes, where only the characteristic of $K$ was relevant, in the rigid-analytic world also the residue characteristic of $K$ (in this case $p$) plays a role in the proof of \cref{rslt:intro-poincare-duality}. Namely, in the case $\ell \ne p$ the result is classical, independently proved by Berkovich \cite{berkovich-etale-cohomology} and Huber \cite{huber-etale-cohomology}. The case $\ell = p$ remained unsolved for a long time and was only proved very recently by Zavyalov \cite{zavyalov-poincare-duality} by loosely following a proof strategy by Gabber. This proof relies on an extensive study of formal models.

The aim of this thesis is to give a new proof of \cref{rslt:intro-poincare-duality} in the case $\ell = p$ (see \cref{rslt:Fp-Poincare-duality-over-field}) by developing a general and powerful 6-functor formalism for $p$-adic sheaves on rigid-analytic varieties, much like its $\ell$-adic analog in \cite{etale-cohomology-of-diamonds}. This allows us to check \cref{rslt:intro-poincare-duality} locally on $X$, which reduces the computation to the case that $X$ is a torus, where it is easy to do by hand (and already essentially appeared in \cite{faltings-p-adic-hodge}). In that regard, our proof of the Poincaré duality in the case $\ell = p$ is very similar to the proof of classical Poincaré duality and also similar to the proof of the $\ell \ne p$ case.

\subsection{The 6-Functor Formalism}

The basic idea behind our 6-functor formalism is as follows: To every rigid-analytic variety $X$ over a fixed non-archimedean field $K$ we want to associate a derived category $D(\mathcal C_X)$ of certain ``complexes of étale $p$-torsion sheaves'' on $X$. The category $D(\mathcal C_X)$ comes equipped with a derived tensor product $\tensor$ and a derived internal hom, denoted $\IHom$, which are the first two functors of the 6-functor formalism. To every map $f\colon Y \to X$ of rigid-analytic varieties we then need to construct an adjoint pair of functors
\begin{align*}
	f^*\colon D(\mathcal C_X) \rightleftarrows D(\mathcal C_Y) \noloc f_*,
\end{align*}
where $f^*$ is the derived pullback functor and $f_*$ is the derived pushforward functor. This provides the next two functors of the 6-functor formalism. These four functors can usually be defined in any geometric setting, e.g. for sheaves on sites. The remaining two functors are a bit more special: For every map $f\colon Y \to X$ as above we need to construct another pair of adjoint functors
\begin{align*}
	f_!\colon D(\mathcal C_Y) \rightleftarrows D(\mathcal C_X) \noloc f^!.
\end{align*}
The functor $f_!$ is known from algebraic topology as the (derived) functor of ``sections with compact support''. It is characterized by the requirement that $f_! = f_*$ if $f$ is proper and $f^! = f^*$ if $f$ is an open immersion.

In order for the six functors $\tensor$, $\IHom$, $f^*$, $f_*$, $f_!$ and $f^!$ to form a 6-functor formalism, one requires them to satisfy several compatibilities. Among others, this includes functoriality of the last four functors in $f$ (i.e. $(f \comp g)^* = g^* \comp f^*$ etc.), proper base-change (i.e. base-change for $f_!$) and the projection formula for $f_!$. Having such a 6-functor formalism available, the proof of \cref{rslt:intro-poincare-duality} boils down to showing that if $f$ is smooth then $f^!$ is a twist of $f^*$ by some invertible sheaf -- then Poincaré duality follows formally by using that $f^!$ is right adjoint to $f_*$ if $f$ is proper. On the other hand, the computation of $f^!$ can be done locally on the source, which allows us to reduce to the case of a relative torus.

Having explained the general idea behind a 6-functor formalism in rigid-analytic geometry, let us now describe how we can implement it, i.e. how to define the categories $\mathcal C_X$ and the six functors concretely. In the case $\ell \ne p$ one can simply take $\mathcal C_X$ to be the category of étale $\Fld_\ell$-sheaves on $X$; this is done in \cite{etale-cohomology-of-diamonds}. In the case $\ell = p$ the situation becomes more subtle. Namely, a $p$-adic 6-functor formalism cannot work directly with étale $\Fld_p$-sheaves, at least not if one hopes to get similarly strong results as in the $\ell$-adic 6-functor formalism. For example, if $K = C$ is algebraically closed and $X$ is the closed ball over $C$, then $R\Gamma(X_\et, \Fld_p)$ ``sees'' the base-field $C$ (this is computed by the usual Artin-Schreier sequence, comparing it to the cohomology of $\mathbb G_a$). On the other hand $(\Spa C)_\et$ is independent of $C$, hence base-change for $\Fld_p$-cohomology along $\Spa C' \to \Spa C$, for any other algebraically closed field $C'$, does not hold.

To solve the problem with base-change, we instead propose to work with étale ``quasicoherent'' $\ri^+_X/\pi$-sheaves on the space $X$, for any pseudouniformizer $\pi$ in $K$ (i.e. any element $\pi \in K$ with $\abs\pi < 1$). Here the sheaf $\ri^+_X$ on $X$ is the sheaf of functions $f \in \ri_X$ such that $\abs{f(x)} \le 1$ for all $x \in X$. The category of ``quasicoherent $\ri^+_X/\pi$-sheaves'' not only knows the étale site of $X$, but also sees the structure sheaf of $X$. In particular, the pullback along $\Spa C' \to \Spa C$ amounts to the functor $- \tensor_{\ri_C/\pi} \ri_{C'}/\pi$, which makes a very general base-change possible. On the other hand, by the Primitive Comparison Theorem \cite[Theorem 5.1]{rigid-p-adic-hodge}, étale $\ri^+/\pi$-cohomology on a proper smooth rigid-analytic variety captures étale $\Fld_p$-cohomology (even better, we will construct a ``$p$-torsion Riemann-Hilbert correspondence'' which provides a strong link between $\Fld_p$-sheaves and $\ri^+_X/\pi$-sheaves, see \cref{rslt:intro-global-RH} below).

We thus need to define a category $\mathcal C_X$ of ``étale quasicoherent $\ri^+_X/\pi$-sheaves on $X$'' associated to every space $X$ and the six functors $\tensor$, $\IHom$, $f^*$, $f_*$, $f_!$ and $f^!$ in this setting. The category $\mathcal C_X$ should satisfy the following properties:
\begin{enumerate}[(i)]
	\item The objects of $\mathcal C_X$ should be quasicoherent, i.e. every object $\mathcal M \in \mathcal C_X$ is locally on $X$, say on some affinoid $U = \Spa(A, A^+)$, of the form $\mathcal M = \widetilde N$ for some $A^+/\pi$-module $N$. Here we use the language of adic spaces introduced by Huber \cite{huber-etale-cohomology} (as opposed to Berkovich spaces or Tate's rigid spaces) as it seems best fit to work with the sheaf $\ri^+_X$ and to define perfectoid spaces \cite{scholze-perfectoid-spaces}.

	\item We should work with almost coefficients throughout. This is necessary because for quasicoherent sheaves to form a good theory, we need the structure sheaf $\ri^+_X/\pi$ to be acyclic on affinoid spaces. With almost coefficients, this holds at least if $X$ is affinoid \emph{perfectoid}, cf. \cite[\S8]{etale-cohomology-of-diamonds}. We will provide a quick introduction to perfectoid geometry in \cref{sec:intro.diam} and to almost mathematics in \cref{sec:intro.almcond} below. In the following explanations it is enough to understand the superscript $(-)^a$ as ``we care about the object only up to things that are annihilated by $\mm$'', where $\mm$ is the ideal of topologically nilpotent elements; if we work over an algebraically closed field $C$ in characteristic $0$ then $\mm$ is generated by $p^{1/p^n}$ for $n \ge 0$.

	\item The objects in $\mathcal C_X$ should come equipped with some ``topological'' structure, as this is the only way that the lower-shriek functor $j_!\colon \mathcal C_U \to \mathcal C_X$ can exist for open immersions $j\colon U \injto X$ (much like in the case of quasicoherent sheaves on schemes, see \cite{condensed-mathematics}). Thus, the objects in $\mathcal C_X$ should be ``solid sheaves'' in a reasonable sense. We will provide a quick introduction to the solid framework in \cref{sec:intro.almcond} below, but for now the reader is invited to just assume solid modules over a ring $A$ to be some sort of ``complete topological'' $A$-modules; the category of such modules is denoted $\catsldmod A$ and we will often speak of ``$A_\solid$-modules'' to refer to objects in this category.
\end{enumerate}
It seems hard to directly define a category of sheaves $\mathcal C_X$ that satisfies all of the above properties, so we will instead define $\mathcal C_X$ by gluing: We know that $\mathcal C_X$ should be equal to the category of $(A^+/\pi)^a_\solid$-modules on an affinoid space $X = \Spa(A, A^+)$, so we can try to glue the categories $\catsldmoda{A^+/\pi}$. This strategy creates some new challenges, however:
\begin{enumerate}[(a)]
	\item The sheaf $\ri^{+a}_X/\pi$ is not acyclic on a general affinoid space $X = \Spa(A, A^+)$ -- in fact, this is almost never the case if $X$ is a rigid-analytic variety over some field. However, $\ri^{+a}_X/\pi$ is acyclic if $X$ is affinoid \emph{perfectoid}.

	\item For a map $\Spa(B, B^+) \to \Spa(A, A^+)$ of affinoid (perfectoid) spaces, the associated map $(A^+/\pi)^a_\solid \to (B^+/\pi)^a_\solid$ of solid rings is usually not flat (even in simple cases like open immersions).
\end{enumerate}
To overcome problem (a) we should glue $\catsldmoda{A^+/\pi}$ only on affinoid perfectoid spaces $X = \Spa(A, A^+)$. In order to get a theory on rigid-analytic spaces as well, this gluing must happen in the \emph{pro-étale} topology, i.e. the topology generated by cofiltered limits of étale maps (every rigid-analytic space admits a pro-étale cover by affinoid perfectoid spaces, see \cite[\S4]{rigid-p-adic-hodge} and \cite[\S15]{etale-cohomology-of-diamonds}, or \cref{sec:intro.diam} below). Problem (b) tells us that there is no way that we can actually glue $\catsldmoda{A^+/\pi}$ in any reasonable way. Instead, we need to pass to the derived $\infty$-category $\Dqcohri(A^+/\pi) := \D(\catsldmoda{A^+/\pi})$ and glue directly on the derived level. We provide a quick introduction to the language of $\infty$-categories and how to glue them in \cref{sec:intro.infty} below, but for now it is enough to view them as ordinary categories (even though the $\infty$-categorical structure is crucial for the gluing to work).

With all of the above considerations at hand, we now have a good strategy for defining the desired $\infty$-category $\DqcohriX X$ playing the role of $D(\mathcal C_X)$: The assignment $X \mapsto \DqcohriX X$ should be a pro-étale sheaf of $\infty$-categories which on every affinoid perfectoid space $X = \Spa(A, A^+)$ evaluates to $\Dqcohri(A^+/\pi)$. One of the main results of this thesis is that such a sheaf actually exists, at least if we only require $\DqcohriX X = \Dqcohri(A^+/\pi)$ on ``small'' affinoid perfectoid $X$. It turns out that this sheaf automatically satisfies v-descent, so we even get a good theory of $\DqcohriX X$ on every small v-stack $X$ (the v-topology and v-stacks are introduced in \cite{etale-cohomology-of-diamonds} and are quickly recalled in \cref{sec:intro.diam} below, but they are not relevant for the following; the reader who is unfamiliar with this formalism is invited to replace ``small v-stack'' by ``diamond'' or just ``analytic adic space over $\Z_p$'' as a good approximation). In order to formulate our theorem, we make one further abstraction: Instead of working with the base sheaf $\ri^+_X/\pi$, we work with a more general class of so-called ``integral torsion coefficients'' $\Lambda$. The precise definition of integral torsion coefficients is a bit subtle (see \cref{def:intergal-torsion-coefficients}), so in the following the reader is invited to always assume $\Lambda = \ri^+_X/\pi$ for some pseudouniformizer $\pi$. We denote by $\vStacksCoeff$ the category of pairs $(X, \Lambda)$ where $X$ is a small v-stack and $\Lambda$ is a system of integral torsion coefficients on $X$.

\begin{theorem} \label{rslt:intro-def-of-qcoh-Lambda-modules}
There is a unique hypercomplete sheaf of $\infty$-categories
\begin{align*}
	\vStacksCoeff^\opp \to \infcatinf, \qquad (X, \Lambda) \mapsto \Dqcohri(X, \Lambda)
\end{align*}
on the v-site of $\vStacksCoeff$ such that for every affinoid perfectoid space $X$ which is of weakly perfectly finite type (see \cref{def:perfectly-finite-type-in-AffPerf}) over some totally disconnected space we have
\begin{align*}
	\Dqcohri(X, \Lambda) = \Dqcohri(\Lambda(X)).
\end{align*}
\end{theorem}

For the definition of totally disconnected perfectoid spaces see \cite[\S7]{etale-cohomology-of-diamonds} or \cref{sec:intro.diam} below; the crucial property of these spaces is that they form a basis of the pro-étale topology. The proof of \cref{rslt:intro-def-of-qcoh-Lambda-modules} is essentially the combination of \cref{rslt:compute-Dqcohri-for-fin-type-over-tot-disc,rslt:v-hyperdescent-for-Dqcohri}, see \cref{rslt:def-of-qcoh-Lambda-modules}. These results rely on a general descent formalism for solid modules which we explain in \cref{sec:intro.descent,sec:intro.descperfd} below.

The $\infty$-category $\Dqcohri(X, \Lambda)$ is stable, hence the underlying ordinary category is triangulated, but in general it has no natural $t$-structure (this is expected by problem (b) above). It contains as a full subcategory the $\infty$-category $\D^a(X, \Lambda)_\omega \subset \Dqcohri(X, \Lambda)$ of discrete objects, which does have a $t$-structure and even admits a much more explicit definition in terms of actual sheaves on $X_\et$ (see \cref{sec:ri-pi.bd-and-disc}). Even though $\Dqcohri(X, \Lambda)$ does not admit a $t$-structure in general, one can still define full subcategories $\Dqcohrip(X, \Lambda)$, $\Dqcohrim(X, \Lambda)$ and $\Dqcohrib(X, \Lambda)$ of (left/right) bounded sheaves -- they can be defined by gluing the $\infty$-categories $\Dqcohriq(\Lambda(X))$ on affinoid perfectoid spaces $X$ as in \cref{rslt:intro-def-of-qcoh-Lambda-modules}, see \cref{rslt:bounded-Dqcohri-is-sheaf}.

Having defined the base $\infty$-category $\Dqcohri(X, \Lambda)$ of ``quasicoherent $\Lambda$-modules on $X$'' we want to work with, it is now time to define the promised six functors. We can immediately define four functors:

\begin{definition}
Let $f\colon Y \to X$ be a map in $\vStacksCoeff$.
\begin{defenum}
	\item The functor $f^*\colon \Dqcohri(X, \Lambda) \to \Dqcohri(Y, \Lambda)$ is the ``restriction'' map of the sheaf $\Dqcohri(-, \Lambda)$.
	\item The functor $f_*\colon \Dqcohri(Y, \Lambda) \to \Dqcohri(X, \Lambda)$ is a right adjoint of $f^*$.
	\item The $\infty$-category $\Dqcohri(X, \Lambda)$ inherits a closed symmetric monoidal structure from the closed symmetric monoidal structure on $\Dqcohri(\Lambda(U))$ for totally disconnected spaces $U$. This means that we have functors $\tensor$ and $\IHom$ on $\Dqcohri(X, \Lambda)$.
\end{defenum}
\end{definition}

With the above definitions, one easily obtains $p$-adic analogs of many $\ell$-adic results. For example, the functors $f^*$ and $f_*$ satisfy base-change on $\Dqcohrip$ along any qcqs map of small v-stacks (see \cref{rslt:base-change-for-bounded-Dqcohri}). This solves one of the main problems of $\Fld_p$-cohomology at the cost of a more involved formalism.

In order to complete the 6-functor formalism, it remains to define the shriek functors $f_!$ and $f^!$ along ``nice'' maps $f$ of small v-stacks. The correct definition for ``nice map'' is captured by so-called \emph{bdcs maps}, which are maps $f\colon Y \to X$ of small v-stacks that are locally compactifiable, representable in locally spatial diamonds and $p$-bounded (see \cref{def:bdcs-maps}). The notion of $p$-boundedness is a bit subtle to define, but roughly states that $f$ has universally finite $p$-cohomological dimension (see \cref{def:p-bounded}). While the definition of bdcs maps is a bit subtle, it is usually easy to check this condition in practice. For example, every map which is of finite type in some sense (e.g. any map of rigid-analytic varieties over some base field) is bdcs, and bdcs maps are stable under composition and base-change.

With the definition of bdcs maps at hand, we can finally define the shriek functors. In general, whenever we can factor $f\colon Y \to X$ as a composition of an étale map $j\colon Y \to Z$ and a proper map $g\colon Z \to X$ then we set $f_! = g_* \comp j_!$, where $j_!$ is a left adjoint of $j^*$ (this exists by a reduction to the case of schemes, see \cref{rslt:existence-of-lower-shriek-for-etale-maps}). Since bdcs maps are locally compactifiable and hence locally admit a factorization as above, we can define $f_!$ in general by extending from the compactifiable case via colimits (see \cref{def:shriek-functors}).

\begin{definition}
Let $f\colon Y \to X$ be a bdcs map in $\vStacksCoeff$.
\begin{defenum}
	\item We define $f_!\colon \Dqcohri(Y, \Lambda) \to \Dqcohri(X, \Lambda)$ as discussed above.

	\item We define $f^!\colon \Dqcohri(X, \Lambda) \to \Dqcohri(Y, \Lambda)$ as a right adjoint of $f_!$.
\end{defenum}
\end{definition}

With the shriek functors available, we have now defined all six functors of the 6-functor formalism. We can finally state the main result of this thesis (see \cref{rslt:main-6-functor-formalism-properties}):

\begin{theorem} \label{rslt:intro-six-functor-formalism}
The assignment $X \mapsto \Dqcohri(X, \Lambda)$ together with the six functors $\tensor$, $\IHom$, $f^*$, $f_*$, $f_!$ and $f^!$ define a 6-functor formalism, i.e. the following properties are satisfied:
\begin{thmenum}
	\item (Functoriality) If $f$ and $g$ are composable maps in $\vStacksCoeff$ then $(g \comp f)^* = f^* \comp g^*$ and $(g \comp f)_* = g_* \comp f_*$. If moreover both $f$ and $g$ are bdcs then $(g \comp f)_! = g_! \comp f_!$ and $(g \comp f)^! = f^! \comp g^!$.

	\item (Special Cases) Let $f$ be a map in $\vStacksCoeff$. If $f$ is étale then $f^! = f^*$ and if $f$ is proper then $f_! = f_*$.

	\item(Projection Formula) Let $f\colon Y \to X$ be a bdcs map in $\vStacksCoeff$, let $\mathcal M \in \Dqcohri(X, \Lambda)$ and $\mathcal N \in \Dqcohri(Y, \Lambda)$. Then there is a natural isomorphism
	\begin{align*}
		f_!(\mathcal N \tensor f^* \mathcal M) = f_! \mathcal N \tensor \mathcal M.
	\end{align*}

	\item (Proper Base-Change). Let
	\begin{center}\begin{tikzcd}
		Y' \arrow[r,"g'"] \arrow[d,"f'"] & Y \arrow[d,"f"]\\
		X' \arrow[r,"g"] & X
	\end{tikzcd}\end{center}
	be a cartesian square in $\vStacksCoeff$ and assume that $f$ is bdcs. Then there is a natural isomorphism
	\begin{align*}
		g^* \comp f_! = f'_! \comp g'^*
	\end{align*}
	of functors $\Dqcohri(Y, \Lambda) \to \Dqcohri(X', \Lambda)$.
\end{thmenum}
\end{theorem}

In our 6-functor formalism there is still one desirable component missing: Poincaré duality for smooth maps. In the world of diamonds and small v-stacks there is no general notion of smoothness, so instead of \emph{proving} Poincarè duality for smooth maps, we will \emph{define} smooth maps via Poincaré duality. This is similar to the $\ell$-adic case, where $\ell$-cohomologically smooth maps are defined. In this spirit, we say that a bdcs map $f\colon Y \to X$ of small v-stacks is \emph{$p$-cohomologically smooth} if its pullback to any totally disconnected space $X' = \Spa(A, A^+)$ satisfies Poincaré duality on $\DqcohriX{X'}$ with respect to some pseudouniformizer $\pi \in A^+$ (see \cref{def:p-cohomologically-smooth-map}). We then get the expected results (see \cref{rslt:p-cohom-smoothness-identities}):

\begin{theorem}
\begin{thmenum}
	\item Let $f\colon Y \to X$ be a bdcs and $p$-cohomologically smooth map in $\vStacksCoeff$. Then the natural morphism
	\begin{align*}
		f^! \Lambda^a \tensor f^* \isoto f^!
	\end{align*}
	of functors $\Dqcohri(X, \Lambda) \to \Dqcohri(Y, \Lambda)$ is an isomorphism and $f^! \Lambda^a$ is invertible.

	\item Let
	\begin{center}\begin{tikzcd}
		Y' \arrow[r,"g'"] \arrow[d,"f'"] & Y \arrow[d,"f"]\\
		X' \arrow[r,"g"] & X
	\end{tikzcd}\end{center}
	be a cartesian square in $\vStacksCoeff$ and assume that $f$ is bdcs and $p$-cohomologically smooth. Then the natural morphism
	\begin{align*}
		g'^* f^! \isoto f'^! g^*
	\end{align*}
	is an isomorphism of functors $\Dqcohri(X, \Lambda) \to \Dqcohri(Y', \Lambda)$.
\end{thmenum}
\end{theorem}

Having defined $p$-cohomologically smooth maps of small v-stacks, it is natural to ask which maps appearing in practice are $p$-cohomologically smooth. One would expect those maps to be $p$-cohomologically smooth which come from smooth maps of adic spaces via diamondification (see \cite[\S15]{etale-cohomology-of-diamonds} or \cref{sec:intro.diam} below), at least in characteristic $0$. This is indeed the case:

\begin{theorem} \label{rslt:intro-adic-smooth-maps-over-Qp-are-p-cohom-smooth}
Let $f\colon Y \to X$ be a smooth map of analytic adic spaces over $\Q_p$. Then the induced map $f^\diamond\colon Y^\diamond \to X^\diamond$ of diamonds is $p$-cohomologically smooth.
\end{theorem}

The proof is found in \cref{rslt:adic-smooth-maps-over-Qp-are-p-cohom-smooth}. It illustrates the benefits of having a full 6-functor formalism available. Namely, $p$-cohomological smoothness is v-local on $X$ and étale local on $Y$, so we can reduce to the case that $X = \Spa(A, A^+)$ is totally disconnected and $Y$ is a relative torus of dimension $1$. In this setting, the required computation essentially already appears in Faltings' work \cite{faltings-p-adic-hodge} and goes as follows. Note first that there is a canonical smooth formal model $\mathfrak Y$ of $Y$ (the torus over $\Spf A^+$). We show in \cref{sec:ri-pi.poincare} that $f_!\colon \DqcohriX Y \to \Dqcohri(A^+/\pi)$ factors as the composition of a ``pushforward'' $t_*\colon \DqcohriX Y \to \Dqcohri(\ri_{\mathfrak Y}/\pi)$ and the functor $\mathfrak f_!\colon \Dqcohri(\ri_{\mathfrak Y}/\pi) \to \Dqcohri(A^+/\pi)$. The latter functor $\mathfrak f_!$ is defined in \cite{condensed-mathematics}, noting that $(\abs{\mathfrak Y}, \ri_{\mathfrak Y}/\pi)$ is a scheme over $A^+/\pi$. Now the computation of $f^!$ reduces to the computation of $f_!$, which in turn reduces to computing $t_*$ and $\mathfrak f_!$:
\begin{itemize}
	\item For $t_*$ we note that since $\mathfrak Y$ is affine, $t_*$ essentially computes cohomology on $Y$. By considering the $\Z_p$-torsor $Y_\infty \to Y$ given by the perfectoid torus $Y_\infty$, we note that $\Dqcohri(Y, \Lambda)$ is the $\infty$-category of derived smooth $\Z_p$-representations (see \cref{rslt:compute-Dqcohri-on-stack-quotient-of-p-bounded-space}). Thus $t_*$ computes smooth $\Z_p$-cohomology, which is given by a Koszul complex.

	\item For $\mathfrak f_!$ we apply the Poincaré duality for coherent cohomology of schemes from \cite{condensed-mathematics} (see also \cref{rslt:scheme-6-functor}.
\end{itemize}
Combining both of the above computations eventually leads to a proof of \cref{rslt:intro-adic-smooth-maps-over-Qp-are-p-cohom-smooth}. In order to deduce $p$-adic Poincaré duality, i.e. \cref{rslt:intro-poincare-duality} in the case $\ell = p$, from \cref{rslt:intro-adic-smooth-maps-over-Qp-are-p-cohom-smooth} we need a way of passing back and forth between $\Fld_p$-sheaves and $\ri^{+a}_X/\pi$-sheaves. The Primitive Comparison Theorem \cite[Theorem 5.1]{rigid-p-adic-hodge} is a first step in that direction, but it is not enough for our purpose: Given a rigid-analytic variety $f\colon X \to \Spa K$ over some non-archimedean field $K$ and a pseudouniformizer $\pi \in K$ with $\pi \divides p$, we need to know that $\mathcal L := f^! \ri^a_K/\pi \in \DqcohriX X$ is of the form $\mathcal L = \mathcal L_0 \tensor_{\Fld_p} \ri^{+a}_X/\pi$ for some invertible complex of étale $\Fld_p$-sheaves $\mathcal L_0$ on $X$. By the computation in \cref{rslt:intro-adic-smooth-maps-over-Qp-are-p-cohom-smooth} this holds locally on $X$, but there is a priori no reason for it to hold globally.

In order to solve the issues raised in the previous paragraph, we introduce the notion of $\varphi$-modules: Given a ring $A$ in characteristic $p$, we denote by $\varphi\colon A \to A$, $x \mapsto x^p$ the Frobenius of $A$. Then a $\varphi$-module over $A$ is an $A$-module $M$ together with a $\varphi$-linear map $\varphi_M\colon M \to M$ such that the induced map $\varphi^* M \to M$ is an isomorphism. This definition naturally extends to the solid and derived setting and then also to solid quasicoherent $\ri^{+a}_X/\pi$-sheaves (provided that $\pi \divides p$), see \cref{def:phi-module-in-DqcohriX}. We denote the resulting $\infty$-category by $\DqcohriX X^\varphi$. Coming back to our problem of understanding $\mathcal L = f^!(\ri^a_K/\pi)$, we observe that $\mathcal L$ comes naturally equipped with a $\varphi$-module structure, induced by the Frobenius $\varphi\colon \ri^a_K/\pi \to \ri^a_K/\pi$. We can now use the following version of a $p$-torsion Riemann-Hilbert correspondence (see \cref{rslt:global-Riemann-Hilbert}):

\begin{theorem} \label{rslt:intro-global-RH}
Let $X$ be a small v-stack over $\Z_p$ with pseudouniformizer $\pi$ such that $\pi \divides p$. Let $\D_\et(X, \Fld_p)^\oc$ denote the $\infty$-category of overconvergent étale derived $\Fld_p$-sheaves on $X$ (see \cref{def:overconvergent-Fp-sheaves}). Then the functor
\begin{align*}
	- \tensor \ri^{+a}_X/\pi\colon \D_\et(X, \Fld_p)^\oc \injto \DqcohriX X^\varphi
\end{align*}
is fully faithful and induces an equivalence of perfect objects on both sides.
\end{theorem}

Applying \cref{rslt:intro-global-RH} to a smooth map $f\colon Y \to X$ of analytic adic spaces over $\Q_p$ we see that $f^!(\ri^{+a}_X/p) = \mathcal L_0 \tensor_{\Fld_p} \ri^{+a}_Y/p$ for some invertible complex $\mathcal L_0$ of étale $\Fld_p$-sheaves on $Y$. If $f$ is of pure dimension $d$ then the proof of \cref{rslt:intro-adic-smooth-maps-over-Qp-are-p-cohom-smooth} shows that $\mathcal L_0$ is concentrated in cohomological degree $-2d$ and is thus of the form $\mathcal L_0 = \omega_f[2d]$ for some $\Fld_p$-local system $\omega_f$ on $Y$. With this knowledge at hand, one can now perform a standard ``deformation to the normal cone'' argument to deduce the following explicit description of $f^!(\ri^{+a}_X/p)$ (see \cref{rslt:dualizing-complex-for-smooth-maps-over-Qp}):

\begin{theorem} \label{rslt:intro-dualizing-complex-for-smooth-maps-over-Qp}
Let $f\colon Y \to X$ be a smooth map of analytic adic spaces over $\Q_p$. If $f$ has pure dimension $d$ then there is a natural isomorphism
\begin{align*}
	f^{\diamond!}(\ri^{+a}_X/p) = \ri^{+a}_Y/p(d)[2d].
\end{align*}
\end{theorem}

It is an immediate corollary of \cref{rslt:intro-dualizing-complex-for-smooth-maps-over-Qp} that Poincaré duality for étale $\Fld_p$-cohomology holds on smooth rigid-analytic varieties (see \cref{rslt:Fp-Poincare-duality-over-field}):

\begin{corollary} \label{rslt:intro-Fp-Poincare-duality-with-omega-X}
Let $K$ be a non-archimedean field over $\Q_p$ and let $X$ be a proper smooth rigid-analytic space of pure dimension $d$ over $K$. Then for all $i \in \Z$ there is a natural perfect pairing
\begin{align*}
	H^i_\et(X, \Fld_p) \tensor_{\Fld_p} H^{2d - i}_\et(X, \Fld_p) \to \Fld_p(-d).
\end{align*}
\end{corollary}

See \cref{rslt:relative-F-p-Poincare-duality} for a relative version of \cref{rslt:intro-Fp-Poincare-duality-with-omega-X} (which can be deduced in the same way).

Poincaré duality for $\Fld_p$-cohomology is the first application of our theory, and was the main motivation for developing it. However, we expect many more applications to be possible, especially in the $p$-adic Langlands program -- after all, for any perfectoid field $K$ and any (locally) profinite group $G$, $\DqcohriX{\Spa K/G}$ is the derived $\infty$-category of smooth $G$-representations on solid almost $\ri_K/\pi$-modules.

\subsection{Overview of the Thesis}

This thesis is structured as follows:
\begin{itemize}
	\item In \cref{sec:andesc} we develop a general theory of analytic spaces (similar to \cite{scholze-analytic-spaces}) in the almost setting and study descent of quasicoherent sheaves. The most important examples for our purposes are schemes, in which case we obtain several descent results for solid quasicoherent sheaves. For example we show that solid quasicoherent sheaves descend along ind-fppf covers of affine schemes.

	\item \cref{sec:ri-pi} is the heart of the thesis: Here we define the $\infty$-category $\Dqcohri(X, \Lambda)$ and construct the six functors for it. The main technical result needed is \cref{rslt:intro-def-of-qcoh-Lambda-modules}, which relies heavily on the descent theory developed in \cref{sec:andesc}.

	\item In \cref{sec:infcat} we provide some technical results on $\infty$-categories which we have not found in the literature.
\end{itemize}
The following subsections of the introduction provide more explanations on the machinery that goes into the 6-functor formalism and sketches the proof of the central descent result \cref{rslt:intro-def-of-qcoh-Lambda-modules}. More concretely, \cref{sec:intro.diam,sec:intro.infty,sec:intro.almcond} explain the main tools used in this thesis: Firstly diamonds, which provide a powerful framework to study rigid-analytic varieties, secondly $\infty$-categories, which allows us to glue derived categories, and thirdly condensed mathematics, whose theory of solid modules makes the existence of lower shriek functors possible. The material covered in these sections is mainly part of the literature, with one exception: In \cref{sec:intro.almcond} we also discuss how to combine condensed mathematics with almost mathematics (which is done in detail in the first part of \cref{sec:andesc}). After discussing the preliminaries, \cref{sec:intro.descent} explains our descent formalism in the context of condensed almost mathematics, which is developed in the second part of \cref{sec:andesc}. Using this descent formalism, \cref{sec:intro.descperfd} sketches the proof of \cref{rslt:intro-def-of-qcoh-Lambda-modules}.

We encourage the reader to have a look at \cref{sec:intro.notation}, where we discuss the somewhat non-standard notation and conventions used in the thesis.

\subsection{Adic Spaces and Diamonds} \label{sec:intro.diam}

Even though our main application of the $p$-adic 6-functor formalism, namely Poincaré duality, is a statement about rigid-analytic varieties (and hence classical $p$-adic geometry), it is very convenient for the development of our theory to allow more general spaces than just rigid-analytic spaces. Among others, we need to work with the pro-étale topology, whose objects are not necessarily rigid-analytic varieties anymore. Scholze's framework of diamonds and small v-stacks, as developed in \cite{etale-cohomology-of-diamonds}, provides such a theory of very general ``spaces'', so we will use it throughout the thesis. In the following we give a quick introduction to the theory of diamonds, addressed at the reader who comes from classical rigid-analytic geometry.

The basic building blocks of diamonds are the perfectoid spaces, first defined by Scholze in \cite{scholze-perfectoid-spaces}. Roughly, a perfectoid space $X$ is an analytic adic space over $\Z_p$ (the reader who is unfamiliar with adic spaces may use Berkovich spaces as an approximation) whose coordinate sheaf $\ri_X$ admits very good approximations of $p$-th roots; more precisely, $X$ is perfectoid if it can be covered by affinoid spaces $\Spa(A, A^+)$ such that $A$ is a complete uniform Tate ring and the Frobenius $\varphi\colon A^\circ/p \to A^\circ/p$ is surjective. As a basic example of a perfectoid space, take any algebraically closed non-archimedean field $C$ and let
\begin{align*}
	X = \Spa(C\langle T^{1/p^\infty}\rangle, \ri_C\langle T^{1/p^\infty}\rangle),
\end{align*}
the perfectoid unit disc over $C$. Perfectoid spaces are usually not noetherian and are therefore not covered by classical results in rigid-analytic geometry, but the property of allowing $p$-th roots still allows a very good theory. To every perfectoid space $X$ we can associate its tilt $X^\flat$, which is a perfectoid space in characteristic $p$ and inherits the same étale site (see \cite[\S7]{scholze-perfectoid-spaces} or \cite[\S3.6]{relative-p-adic-hodge-1}). This allows to transfer many problems in mixed characteristic to the equal characteristic setting.

The notion of diamonds is motivated by the following observation: Given a rigid-analytic variety $X$ (in fact the following works for any analytic adic space $X$ over $\Z_p$) one can consider the pro-étale site $X_\proet$ of $X$ whose objects are roughly the cofiltered limits of objects in $X_\et$, similar to the pro-étale site for schemes discussed in \cite{proetale-topology}. One now notes that $X_\proet$ admits a basis consisting of perfectoid spaces, as is illustrated by the following example:

\begin{example}
Let $X = \Spa(C\langle T^{\pm1} \rangle, \ri_C \langle T^{\pm1} \rangle)$ be a torus over some algebraically closed field $C$ in characteristic $0$. Then one can obtain the perfectoid torus
\begin{align*}
	\tilde X = \Spa(C\langle T^{\pm1/p^\infty}\rangle, \ri_C\langle T^{\pm1/p^\infty}\rangle)
\end{align*}
by adjoining all $p^n$-th roots of the coordinate $T$. Since adjoining the $p^n$-th roots of $T$ yields an étale map, in the limit we obtain a pro-étale map $\tilde X \to X$.
\end{example}

In general we can choose a pro-étale cover $\tilde X \to X$ by some perfectoid space $\tilde X$. But then $X$ itself can be seen as the quotient of $\tilde X$ along some pro-étale equivalence relation $R \subset \tilde X \cprod \tilde X$ (this is similar to how a nice topological space can be interpreted as a quotient of its universal cover by the action of the fundamental group). We can thus replace $X$ by the pair $(\tilde X, R)$ without losing much information. This perspective has many advantages due to the excellent properties of perfectoid spaces. In particular, through this perspective we can leave the realm of noetherian spaces and for example form cofiltered limits of our spaces.

To formalize the idea of considering pairs $(\tilde X, R)$ as above (which we only want to consider up to some form of equivalence, of course), the theory of sheaves comes in handy. In fact, we can replace $(\tilde X, R)$ by the sheaf of sets $\tilde X/R$ on the pro-étale site $\Perfd_\proet$ of all perfectoid spaces (where we confuse $\tilde X$ with the sheaf it represents). Taking things one step further, we can make use of the tilting correspondence and instead consider the sheaf of sets $\tilde X^\flat/R^\flat$ on the site $\Perf$ of all perfectoid spaces \emph{of characteristic $p$}. By the tilting correspondence we do not lose much information by this step (e.g. all étale information is preserved), but gain a lot more power because perfectoid spaces in characteristic $p$ are particularly simple. We arrive at the following definition:

\begin{definition}[{cf. \cite[Definition 11.1]{etale-cohomology-of-diamonds}}]
A \emph{diamond} is a sheaf of sets $X$ on the pro-étale site on the category $\Perf$ of perfectoid spaces of characteristic $p$ such that $X$ can be written as the quotient $\tilde X/R$ of a perfectoid space $\tilde X \in \Perf$ (confused with the sheaf it represents) and a representable equivalence relation $R \subset \tilde X \cprod \tilde X$ such that the two projections $R \to \tilde X$ are pro-étale.
\end{definition}

We refer the reader to \cite[\S11]{etale-cohomology-of-diamonds} for some basic properties of diamonds, in particular their excellent stability properties inside the category of all sheaves on $\Perf$. Note that the definition of diamonds can roughly be seen as an analog of the definition of algebraic spaces in algebraic geometry.

Starting with a rigid-analytic variety $X$, we have now explained how to associated to $X$ a diamond $X^\diamond = \tilde X^\flat/R^\flat$ (see \cite[Definition 15.5]{etale-cohomology-of-diamonds} for the precise construction). It is now natural to ask what information of $X$ is captured by $X^\diamond$. We have the following:
\begin{itemize}
	\item To every diamond $X$ one can associate a topological space $\abs X$ (see \cite[Definition 11.14]{etale-cohomology-of-diamonds}) by taking the quotient $\abs X := \abs{\tilde X}/\abs R$ for any representation $X = \tilde X/R$ (of course one needs to verify that this is independent of the chosen representation for $X$). If $X$ is the diamond associated to some rigid-analytic variety, then $\abs X$ is equal to the underlying topological space of that variety (see \cite[Lemma 15.6]{etale-cohomology-of-diamonds}).

	\item We define a map $f\colon Y \to X$ of diamonds to be \emph{étale} if it is locally separated (this is a technical condition) and the pullback $\tilde f\colon \tilde Y \to \tilde X$ of $f$ along any map $\tilde X \to X$ from a perfectoid space $\tilde X \in \Perf$ is representable by an étale map of perfectoid spaces. If $X$ is the diamond associated to a rigid-analytic variety then the thus defined étale site $X_\et$ coincides with the étale site of that variety (see \cite[Lemma 15.6]{etale-cohomology-of-diamonds}).

	\item In general the structure sheaf $\ri_{X'}$ of a rigid-analytic variety $X'$ is not captured by the diamond $X = X'^\diamond$. The main obstacle for this is that in the process of constructing $X'^\diamond$ we passed to the tilts and hence lost any information about the characteristic we were in. We can recover this information by remembering the structure map $X' \to \Spa \Z_p$, which on the diamond side defines a map $X'^\diamond \to (\Spa \Z_p)^\diamond$. Here $(\Spa\Z_p)^\diamond$ is the sheaf that parametrizes untilts (see \cite[Lemma 15.1]{etale-cohomology-of-diamonds}); note that it is not a diamond because $\Spa\Z_p$ is not an \emph{analytic} adic space. In this thesis we will refer to a diamond $X$ together with a map $X \to (\Spa\Z_p)^\diamond$ as an \emph{untilted diamond} and often denote it by $X^\sharp$.

	Given an untilted diamond $X^\sharp$ we can define an associated structure sheaf $\ri_{X^\sharp}$ on $\abs X$ as follows: For every pro-étale map $\tilde X \to X$ from some perfectoid space $\tilde X \in \Perf$ we get an untilted diamond $\tilde X^\sharp$ (i.e. a map $\tilde X \to (\Spa\Z_p)^\diamond$), which in this case corresponds precisely to an untilt of $\tilde X$, i.e. a perfectoid space $\tilde X^\sharp$ whose tilt is $\tilde X$. We then define $\ri_{X^\sharp}(\tilde X) := \ri_{\tilde X^\sharp}(\tilde X^\sharp)$. One checks that this defines a pro-étale sheaf on the site of perfectoid spaces $\tilde X \in X_\proet$. As these form a basis of $X_\proet$ we can uniquely extend this sheaf to the whole of $X_\proet$ and in particular obtain the desired sheaf $\ri_{X^\sharp}$ on $\abs X$. One can similarly define the integral structure sheaf $\ri^+_{X^\sharp} \subset \ri_{X^\sharp}$.

	If $X$ is a locally noetherian analytic adic space over $\Q_p$ with associated untilted diamond $X^\diamond$ then $\ri^+_X/p = \ri^+_{X^\diamond}/p$ on $X_\et = X^\diamond_\et$. If $X$ is seminormal and of finite type over some non-archimedean field over $\Q_p$ then $\ri_X = \ri_{X^\diamond}$ and $\ri^+_X = \ri^+_{X^\diamond}$. Both of these results are summarized in \cite[Proposition 2.8]{mann-werner-simpson}, while the latter result is originally due to Kedlaya-Liu \cite[Theorem 8.2.3]{relative-p-adic-hodge-2}.
\end{itemize}
Using the topological space $\abs X$ associated to a diamond $X$ one can single out a special subclass of diamonds: A diamond $X$ is called \emph{spatial} if $X$ is qcqs and $\abs X$ admits a basis given by $\abs U$ for quasicompact open immersions $U \injto X$; $X$ is called \emph{locally spatial} if it admits an open cover by spatial diamonds (see \cite[Definition 11.17]{etale-cohomology-of-diamonds}). The diamond associated to an analytic adic space is always locally spatial, hence the locally spatial diamonds are the most important ones for our purpose.

By studying rigid-analytic varieties through their associated diamond, we have many new tools at our disposal. Most prominently we will make heavy use of pro-étale descent: Whenever we want to show that a property holds for all diamonds, it is enough to verify this property for a basis of $\Perf_\proet$ and then show that the property satisfies pro-étale descent. This technique is extremely powerful, partly because of the following basis of $\Perf_\proet$.

\begin{definition}[{cf. \cite[Definition 7.1]{etale-cohomology-of-diamonds}}]
A perfectoid space $X$ is called \emph{totally disconnected} if $X$ is qcqs and for every open cover $X = \bigunion_i U_i$ the map $\bigdunion_i U_i \to X$ has a section.
\end{definition}

On the one hand one can show that every perfectoid space admits a cover by totally disconnected perfectoid spaces (see \cite[Lemma 7.18]{etale-cohomology-of-diamonds}), hence the totally disconnected space do indeed form a basis of $\Perf_\proet$. On the other hand, totally disconnected spaces are of very simple shape: If $X$ is totally disconnected then it is affinoid and every connected component of $X$ is of the form $\Spa(K, K^+)$ for some perfectoid field $K$ and some open an bounded valuation subring $K^+ \subset K$. In other words, $X$ is assembled from a profinite set $\pi_0(X)$ and a collection of valuation rings $K^+$.

In a sense, the above explanations show that understanding $p$-adic geometry reduces to an understanding of valuation rings and profinite sets. One may compare this to the setting of schemes, where the pro-étale local theory (given by the w-local schemes or even w-contractible schemes, see \cite[\S2]{proetale-topology}) deals with profinite sets and henselian local rings, which are in general much more complicated than valuation rings. In that sense, $p$-adic geometry is easier than algebraic geometry, as it should be!

We close this section with a few final remarks about a generalization of diamonds. Namely, instead of working with the pro-étale topology on $\Perf$, one could use the much finer \emph{v-site} $\Perf_\vsite$. In the v-site, coverings are essentially given by jointly surjective families of maps (plus some quasicompactness requirements similar to the fpqc topology for schemes), see \cite[Definition 8.1]{etale-cohomology-of-diamonds}. One can show that every diamond is automatically a v-sheaf, which leads to the study of v-sheaves in general. In fact, we can make the following definition:

\begin{definition}[{cf. \cite[Definitions 12.1, 12.4]{etale-cohomology-of-diamonds}}]
A \emph{small v-sheaf} is a sheaf of sets on the v-site $\Perf_\vsite$ such that there is a surjective map of v-sheaves $\tilde X \surjto X$ from some perfectoid space $X$. A \emph{small v-stack} $Y$ is a stack on $\Perf_\vsite$ such that there is a surjective map of stacks $\tilde X \to Y$ for which $R = \tilde X \cprod_Y \tilde X$ is a small v-sheaf.
\end{definition}

Note that the smallness assumption on the v-sheaves and v-stacks is really just a set-theoretical subtlety coming from the fact that the site $\Perf_\vsite$ is not small. For example, every qcqs v-stack is automatically small.

The category of diamonds embedds fully faithfully into the category of small v-sheaves, which in turn embedds fully faithfully into the 2-category of small v-stacks. Many of the definitions and constructions about diamonds extend to small v-stacks. For example, for every small v-stack $X$ we get an associated topological space $\abs X$, a good notion of étale maps $Y \to X$ and if $X$ comes equipped with an untilt $X^\sharp$ (i.e. a map $X \to (\Spa\Z_p)^\diamond$ of small v-stacks) then we get the structure sheaves $\ri_{X^\sharp}$ and $\ri^+_{X^\sharp}$ on $X$. We denote by $X_\vsite$ the v-site of $X$, which consists of all small v-stacks over $X$ with covers being families of jointly surjective maps. By the above, $X_\vsite$ admits a basis consisting of totally disconnected perfectoid spaces. Hence whenever we want to show that a property holds for $X$, it is enough to verify this property on totally disconnected spaces and then show that the property satisfies descent along v-covers. This ``v-descent'' technique is used abundantly throughout the thesis.

\subsection{\texorpdfstring{$\infty$}{Infinity}-Categories} \label{sec:intro.infty}

Throughout this thesis we will make heavy use of the language of $\infty$-categories, as introduced by Lurie in \cite{lurie-higher-topos-theory,lurie-higher-algebra,lurie-spectral-algebraic-geometry} (see also Kerodon \cite{kerodon}, which is an analog of the Stacks Project for $\infty$-categories). While this language is very powerful and is arguably a much better model for derived categories than triangulated categories are (in particular it allows us to glue them, which is crucial for our constructions), it also takes some time to get used to. Since $\infty$-categories are not standard in arithmetic geometry yet, we decided to give a quick introduction to the topic (with a focus on what is needed in this thesis), so that the reader who is unfamiliar with $\infty$-categories can at least get a basic idea of the core concepts we are working with throughout the thesis.

Recall that a category $\mathcal C$ consists of a collection of objects and morphisms between them together with a rule of composing morphisms. These morphisms allow us to relate the objects of $\mathcal C$ and we usually only care about these objects up to isomorphism. The basic idea of $\infty$-categories is to apply a similar philosophy to the morphisms in $\mathcal C$. In other words, we want to introduce morphisms between the morphisms in $\mathcal C$, called the $2$-morphisms, and usually only care about morphisms in $\mathcal C$ up to $2$-isomorphism. Repeating the pattern, we then introduce $3$-morphisms (i.e. morphisms between $2$-morphisms) and so on, eventually leading to an $\infty$-category.

One of the biggest challenges when working with $\infty$-categories is keeping track of this infinite amount of data. This already starts with a precise definition of $\infty$-categories, which is not at all obvious how to do. One way to define $\infty$-categories is via simplicial sets. For the following definition, we denote by $\Delta$ the simplex category, i.e. the category of totally ordered finite sets, where morphisms are monotonous maps. Up to isomorphism the objects of $\Delta$ are given by the sets $[n] := \{ 1, \dots, n \}$ for all integers $n \ge 0$ (with the canonical ordering). Recall that a \emph{simplicial set} is a functor $X\colon \Delta^\opp \to \catset$ from $\Delta^\opp$ to the category of sets. For every $n \ge 0$ the set $X_n := X([n])$ is called the set of \emph{$n$-simplices} of $X$. We single out two particular types of simplicial sets: For every $n \ge 0$ we let $\Delta^n$ denote the simplicial set which corresponds to an $n$-dimensional tetraeder, i.e. $\Delta^n$ has a single $n$-simplex, $n + 1$ $(n-1)$-simplices (the faces) and so on, down to $n + 1$ $0$-simplices (the corners). For every $n \ge 0$ and every $0 \le i \le n$ we denote $\Lambda^n_i \subset \Delta^n$ the subsimplicial set where we removed the $n$-simplex and the $i$-th $(n-1)$-simplex. With this notation at hand, we can now define $\infty$-categories as follows:

\begin{definition}
An \emph{$\infty$-category} is a simplicial set $\mathcal C$ such that for all integers $0 < i < n$, every map $f_0\colon \Lambda^n_i \to \mathcal C$ of simplicial sets admits an extension $f\colon \Delta^n \to \mathcal C$.
\end{definition}

Given an $\infty$-category $\mathcal C$, the $0$-simplices $\mathcal C_0$ are the \emph{objects} of $\mathcal C$ and the $1$-simplices $\mathcal C_1$ are the \emph{morphisms}. Unlike in ordinary category theory, there is in general no composition rule in the $\infty$-category $\mathcal C$. Instead, given morphisms $f\colon X \to Y$ and $g\colon Y \to Z$ in $\mathcal C$ the definition of $\infty$-categories guarantees the existence of a $2$-simplex $\sigma \in \mathcal C_2$ corresponding to a diagram of the form
\begin{center}\begin{tikzcd}
	& Y \arrow[dr,"g"]\\
	X \arrow[ur,"f"] \arrow[rr,"h"] && Z
\end{tikzcd}\end{center}
We say that $\sigma$ exhibits the morphism $h\colon X \to Z$ as a composition of $f$ and $g$ and often abuse notation by writing $h \isom g \comp f$. In the case $Y = Z$ and $g = \id$ we get a notion of isomorphism $f \isom h$. By passing to isomorphism classes of morphisms (and forgetting all higher morphisms) we obtain an ordinary category underlying $\mathcal C$. Conversely every ordinary category can be seen as an $\infty$-category by adding exactly one higher morphism for every composition diagram.

One can generalize many notions from ordinary category theory to $\infty$-categories (see \cite{kerodon}). In particular one can define isomorphisms and initial and final objects in an $\infty$-category. Given an $\infty$-category $\mathcal C$ and any simplicial set $K$ one can show that $\Hom(K, \mathcal C)$ (which can naturally be seen as a simplicial set) is again an $\infty$-category. We usually denote it $\Fun(K, \mathcal C)$ and call it the \emph{$\infty$-category of functors} from $K$ to $\mathcal C$.

A \emph{diagram} in an $\infty$-category is a functor $I \to \mathcal C$ from some simplicial set $I$ (often $I$ is a category, but it does not need to be). Given a diagram $f\colon I \to \mathcal C$ one can then define the limit $\varprojlim_I f \in \mathcal C$ and the colimit $\varinjlim_I f \in \mathcal C$, if they exist. Roughly the limit (resp. colimit) is defined to be the universal extension of the diagram $f$ to a diagram $I^\triangleleft \to \mathcal C$ (resp. $I^\triangleright \to \mathcal C$), where $I^\triangleleft$ the simplicial set $I$ together with a new point $*$ that has a unique morphism to all objects in $I$ (and dually for $I^\triangleright$). When discussing derived $\infty$-categories below we will provide an example of a limit in the $\infty$-categorical context and explain the difference to limits in ordinary categories.

In the following we elaborate on some concepts in the $\infty$-categorical world that are used repeatedly throughout the thesis.

\paragraph{Derived $\infty$-Categories.} In this thesis most of the $\infty$-categories we deal with are derived $\infty$-categories or can be constructed from derived $\infty$-categories. Let us therefore explain in some detail how to work with them. Recall that classically, the derived category $D(\mathcal A)$ of a ``nice'' abelian category $\mathcal A$ is constructed in the following steps:
\begin{align*}
	\mathcal A \quad \leadsto \quad C(\mathcal A) \quad \leadsto \quad K(\mathcal A) \quad \leadsto \quad D(\mathcal A).
\end{align*}
Here $C(\mathcal A)$ denotes the category of chain complexes in $\mathcal A$, $K(\mathcal A)$ denotes the homotopy category and $D(\mathcal A)$ the derived category. In the step from $C(\mathcal A)$ to $K(\mathcal A)$ we mod out homotopies (i.e. we identify morphisms which differ by a chain homotopy) and in the step from $K(\mathcal A)$ to $D(\mathcal A)$ we invert quasi-isomorphisms. While the second step is not an issue, it turns out that in the first step we lose too much information. Instead of modding out homotopies we should add the homotopies as higher morphisms, leading to an $\infty$-category $\mathcal K(\mathcal A)$ whose underlying ordinary category is $K(\mathcal A)$. Inverting quasi-isomorphisms in $\mathcal K(\mathcal A)$ leads to the \emph{derived $\infty$-category} $\D(\mathcal A)$ (see \cite[\S1.3]{lurie-higher-algebra} for details).

It turns out that the $\infty$-category $\D(\mathcal A)$ contains all the information which was classically captured by the triangulated structure on $D(\mathcal A)$. For example, take any $M \in \D(\mathcal A)$ and let $X := \lim(0 \to M \from 0)$, so that we have a pullback diagram
\begin{center}\begin{tikzcd}
	X \arrow[r] \arrow[d] & 0 \arrow[d] \\
	0 \arrow[r] & M
\end{tikzcd}\end{center}
in $\D(\mathcal A)$. Taking this limit in the ordinary category $D(\mathcal A)$ would result in $X = 0$, but in the $\infty$-category $\D(\mathcal A)$ the situation is different. By definition of limits the above pullback diagram needs to be universal among all such squares in $\D(\mathcal A)$. Of course all the morphisms in such a square are determined: they need to be $0$. However, in the $\infty$-category $\D(\mathcal A)$ the square also contains a chain homotopy of the zero map $X \to M$ to itself (this is what the makes the diagram commute). One checks that the universal chain homotopy of the zero map to $M$ is given by $X = M[-1]$ with the obvious chain homotopy $M[-1] \to M$ (given by the identity maps on the complex). Thus the $\infty$-category $\D(\mathcal A)$ captures the shift operator on $D(\mathcal A)$.

The exact triangles in $D(\mathcal A)$ are also captured by $\D(\mathcal A)$: A triangle $X \to Y \to Z$ in $D(\mathcal A)$ is exact if and only if there is a pullback square of the form
\begin{center}\begin{tikzcd}
	X \arrow[r] \arrow[d] & Y \arrow[d]\\
	0 \arrow[r] & Z
\end{tikzcd}\end{center}
in $\D(\mathcal A)$. In particular this makes $X = \lim(Y \to Z \from 0)$ functorial in the morphism $Y \to Z$, a property that the classical perspective on derived categories lacks. We write
\begin{align*}
	X = \fib(Y \to Z) := \lim(Y \to Z \from 0)
\end{align*}
and call it the \emph{fiber} of the morphism $Y \to Z$. One remarkable property of $\D(\mathcal A)$ is that every pullback square in $\D(\mathcal A)$ is automatically a pushout square and vice-versa. In particular we can also recover $Z$ from the morphism $X \to Y$ via
\begin{align*}
	Z = \cofib(X \to Y) := \colim(0 \from X \to Y).
\end{align*}
We call $Z$ the \emph{cofiber} of $X \to Y$. With these definitions at hand, the whole triangulated structure of $D(\mathcal A)$ is captured by pullback and pushout squares in $\D(\mathcal A)$. From this point of view one sees that the triangulated structure of $D(\mathcal A)$ is a way of axiomatizing the behavior of finite limits and colimits in $\D(\mathcal A)$. In general limits and colimits behave poorly in $D(\mathcal A)$ (which is why the triangulated structure is needed), but $\D(\mathcal A)$ does not share this defect. In fact, we will make plenty use of large limits and colimits in $\D(\mathcal A)$ throughout the thesis.

Just like derived categories can be generalized to triangulated categories, there is a similar generalization of derived $\infty$-categories in the $\infty$-categorical setting:

\begin{definition}
A \emph{stable $\infty$-category} is an $\infty$-category $\mathcal C$ which admits a zero object $0 \in \mathcal C$ (i.e. $0$ is both initial and final), in which all finite limits and colimits exist and in which every square is a pullback square if and only if it is a pushout square.
\end{definition}

We refer the reader to \cite[\S1]{lurie-higher-algebra} for more information on stable $\infty$-categories (see \cite[Proposition 1.1.3.4]{lurie-higher-algebra} for our non-standard definition). Given a stable $\infty$-category $\mathcal C$, the underlying ordinary category of $\mathcal C$ comes naturally equipped with a triangulated structure in the same way as discussed for derived $\infty$-categories above. The class of stable $\infty$-categories is stable under many operations, including arbitrary limits of $\infty$-categories. It follows that most of the $\infty$-categories we deal with in this thesis are stable.

A functor $F\colon \mathcal C \to \mathcal D$ of stable $\infty$-categories is called \emph{exact} if it preserves all finite limits and colimits. It follows easily from the definition of stable $\infty$-categories that $F$ is exact as soon as it preserves either finite limits or finite colimits, which means that essentially all the functors we deal with are exact.

A \emph{$t$-structure} on a stable $\infty$-category $\mathcal C$ is a $t$-structure on the underlying ordinary triangulated category, i.e. it is a pair $(\mathcal C_{\ge0}, \mathcal C_{\le0})$ of full subcategories of $\mathcal C$ satisfying certain axioms (see \cite[\S1.2.1]{lurie-higher-algebra}). Given a $t$-structure on $\mathcal C$, the heart $\mathcal C^\heartsuit := \mathcal C_{\ge0} \isect \mathcal C_{\le0}$ is an ordinary abelian category. Moreover, for every $n \in \Z$ we have the truncation functors $\tau_{\ge n}\colon \mathcal C \to \mathcal C_{\ge n}$ and $\tau_{\le n}\colon \mathcal C \to \mathcal C_{\le n}$ and the homology functors $\pi_n\colon \mathcal C \to \mathcal C^\heartsuit$. Note that we use homological notation throughout (i.e. $\pi_n(M) = H^{-n}(M)$, where $H^\bullet$ denotes the cohomology functor associated with the $t$-structure).

\paragraph{Anima and Hom-Spaces.} In ordinary category theory the category $\catset$ of sets plays a special role because for any two objects $X$ and $Y$ in some category $\mathcal C$ the collection of morphisms $X \to Y$ forms a set $\Hom(X, Y)$ (up to set-theoretic issues which we do not care about here). In the $\infty$-categorical setting the category $\catset$ is replaced by the $\infty$-category $\Ani$ of \emph{anima} (anima are called ``spaces'' and $\Ani$ is denoted by $\catS$ in \cite{lurie-higher-topos-theory}). The $\infty$-category $\Ani$ can be constructed by starting with the category of CW complexes (e.g. modeled by certain simplicial sets) and adding homotopies as higher morphisms. In particular every $X \in \Ani$ has an underlying set $\pi_0(X)$ and for every $x \in \pi_0(X)$ and all $n \ge 1$ we have a homotopy group $\pi_n(X, x)$. These homotopy groups are conservative, i.e. a map of anima $f\colon X \to Y$ is an isomorphism as soon as $\pi_0(X) \to \pi_0(Y)$ and all $\pi_n(X, x) \to \pi_n(Y, f(x))$ are bijections. There is a natural fully faithful embedding $\catset \injto \Ani$ which identifies a set $S$ with the anima $S'$ with $\pi_0(S') = S$ and $\pi_n(S', x) = 0$ for $n > 0$ and all $x \in S$. We say that an anima $X \in \Ani$ is \emph{static} if $\pi_n(X, x) = 0$ for $n > 0$ and $x \in \pi_0(X)$, i.e. if $X$ lies in the image of $\catset \injto \Ani$.

Given an $\infty$-category $\mathcal C$ and objects $X, Y \in \mathcal C$, we get an associated anima $\Hom(X, Y) \in \Ani$ of morphisms $X \to Y$. Here the underlying set $\pi_0 \Hom(X, Y)$ can be identified with the isomorphism classes of morphisms $X \to Y$ in $\mathcal C$. The higher homotopy groups of $\Hom(X, Y)$ are in general harder to grasp. If $\mathcal C$ is a stable $\infty$-category then $\pi_n \Hom(X, Y)$ (which is independent of the base point) can be identified as the set of isomorphism classes of morphisms $X[n] \to Y$.

\paragraph{The $\infty$-Category of $\infty$-Categories.} One can construct an $\infty$-category $\infcatinf$ whose objects are the (small) $\infty$-categories and whose morphisms are the functors between $\infty$-categories (see \cite[Section 01YV]{kerodon}). In this thesis we often encounter the problem of constructing certain functors $K \to \infcatinf$ from some simplicial set $K$. This is in general hard to do because we need to take care of an infinite amount of higher coherences. Here the theory of Cartesian and coCartesian fibration comes in handy. Namely, functors $F\colon K \to \infcatinf$ correspond bijectively to coCartesian fibrations $\phi\colon X \to K$ of simplicial sets, see \cite[Section 01J2]{kerodon}. We do not want to give the precise definition of coCartesian fibrations here, but we do point out that for a map of simplicial sets to be a coCartesian fibration is a \emph{condition} rather than an additional datum. In the forward direction of the above correspondence, we associate to the functor $F$ the simplicial set $X$ as follows: An object of $X$ is a pair $(k, x)$ where $k$ is an object of $K$ and $x$ is an object of $F(k)$. A morphism $f\colon (k_x, x) \to (k_y, y)$ in $X$ is a pair $(\alpha, f_\alpha)$ where $\alpha\colon k_x \to k_y$ is a morphism in $K$ and $f_\alpha\colon F(\alpha)(x) \to y$ is a morphism in $F(k_y)$. Higher morphisms can be defined similarly. Passing back and forth between functors $K \to \infcatinf$ and coCartesian fibrations $X \to K$ is a very useful technique for constructing functors to $\infcatinf$.

The $\infty$-category $\infcatinf$ admits all small limits and colimits (see \cite[\S7.3]{kerodon}). Limits admit a particularly nice description using coCartesian fibrations. Namely, suppose we are given a diagram $F\colon I \to \infcatinf$ of $\infty$-categories and we want to compute $\varprojlim_I F \in \infcatinf$. Let $\phi\colon X \to I$ be the coCartesian fibration which classifies $F$. Then $\varprojlim_I F$ can be identified with the full subcategory of sections $s \in \Fun_I(I, X)$ such that for every morphism $\alpha\colon i \to j$ in $I$ the morphism $s(\alpha)$ in $X$ is $\phi$-coCartesian.

\paragraph{Gluing $\infty$-Categories.} The main object of study in this thesis, the $\infty$-category $\DqcohriX X$ of quasicoherent $\ri^{+a}/\pi$-sheaves on a small v-stack $X$, will be obtained by gluing the $\infty$-categories $\Dqcohri(A^+/\pi)$ of derived solid almost $A^+/\pi$-modules on affinoid perfectoid spaces $X = \Spa(A, A^+)$ along v-covers. Let us briefly describe what we mean by ``gluing'' in this context.

The general setup is as follows: Suppose we are given a site $\mathcal C$ (e.g. the big v-site of all small v-stacks) and we want to associate to every object $X \in \mathcal C$ an $\infty$-category $F(X)$ by ``gluing'' the predefined $\infty$-categories $F(B)$ on some basis $B \in \mathcal B \subset \mathcal C$. This idea can be formalized by seeking a certain functor $F\colon \mathcal C^\opp \to \infcatinf$ (where for every map $f\colon X \to Y$ in $\mathcal C$ the induced functor $F(f)$ is a pullback functor $f^*\colon \mathcal F(Y) \to \mathcal F(X)$) whose restriction to the basis $\mathcal B$ is predetermined. In order for such a functor $F$ to be unique we require it to satisfy descent along all covers in $\mathcal C$, i.e. we want $F$ to be a \emph{sheaf of $\infty$-categories} on $\mathcal C$.

Motivated by the previous paragraph, let us explain what a sheaf is in the $\infty$-categorical setting (see \cref{sec:infcat.sheaves} for more details). Suppose we are given a site $\mathcal C$ and some nice $\infty$-category $\mathcal D$ of values (e.g. $\mathcal D = \infcatinf$ or $\mathcal D = \Ani$). Then a \emph{$\mathcal D$-valued sheaf} on $\mathcal C$ is roughly a functor $F\colon \mathcal C^\opp \to \mathcal D$ which turns coproducts in $\mathcal C$ to products in $\mathcal D$ and satisfies the following descent condition: For every cover $Y \to X$ in $\mathcal C$ with associated Čech nerve $Y_\bullet \to X$ (i.e. $Y_n = Y \cprod_X \cdots \cprod_X Y$) the natural morphism
\begin{align*}
	F(X) \isoto \varprojlim_{[n]\in\Delta} F(Y_n)
\end{align*}
is an equivalence. To understand the latter condition, we claim that if $\mathcal D = \infcatinf$ then the right-hand limit computes the $\infty$-category of descent data of $F$ along the cover $Y \to X$. Namely, one checks that an object in the right-hand $\infty$-category can be described as follows: It consists of an object $M \in F(Y_0)$, an isomorphism of the two pullbacks of $M$ to $Y_1$, a homotopy in $F(Y_2)$ which exhibits the usual cocycle condition of the three pullbacks of the isomorphism on $Y_1$ to $Y_2$, a higher homotopy in $F(Y_3)$ which exhibits a cocycle condition for the homotopy on $F(Y_2)$ and so on. If all $F(Y_n)$ are ordinary categories then the cocycle homotopy in $F(Y_2)$ is uniquely determined and therefore automatically satisfies all higher cocycle conditions; thus in this case we arrive at the classical definition of descent data.

We say that a $\mathcal D$-valued sheaf $F$ on the site $\mathcal C$ is \emph{hypercomplete} if it satisfies the descent condition $F(X) \isoto \varprojlim_{[n]\in\Delta} F(Y_n)$ for all hypercovers $Y_\bullet \to X$. In the $\infty$-categorical setting this condition is usually stronger than just satisfying Čech descent.

The usual operations of sheaves, like pullback, pushforward and sheafification, extend to the $\infty$-categorical setting.

\paragraph{Totalizations in Derived $\infty$-Categories.} Let $\mathcal C$ be an $\infty$-category. A \emph{cosimplicial object in $\mathcal C$} is a functor $X^\bullet\colon \Delta \to \mathcal C$. The limit of this diagram is called the \emph{totalization} of $X^\bullet$ and denoted $\Tot(X^\bullet) \in \mathcal C$ (assuming it exists). By the discussion of sheaves of $\infty$-categories above it is unsurprising that we will often encounter totalizations of cosimplicial objects in $\infty$-categories throughout the thesis. In the following we make a few remarks on how to compute totalizations in the setting of derived $\infty$-categories (or more generally stable $\infty$-categories with $t$-structure). See \cite[\S1.2.4]{lurie-higher-algebra} for details.

Let $\mathcal C$ be a stable $\infty$-category and let $X^\bullet$ be a cosimplicial object in $\mathcal C$. For every $n \ge 0$ we denote by $\Tot_n(X^\bullet)$ the limit of the $n$-truncated simplicial object $\Delta_{\le n} \injto \Delta \xto{X^\bullet} \mathcal C$. By commuting limits with limits we see that $\Tot(X^\bullet) = \varprojlim_n \Tot_n(X^\bullet)$. This perspective is useful because any exact functor $F\colon \mathcal C \to \mathcal D$ of stable $\infty$-categories commutes with $\Tot_n$. Thus in order to show that $F$ commutes with $\Tot$ it only remains to check that $F$ commutes with the $\varprojlim_n$. One strategy to show such a statement is as follows: Suppose $\mathcal C$ comes equipped with a $t$-structure and that $X^n \in \mathcal C_{\le0}$ for all $n$. Then $\Tot_n$ and $\Tot$ agree on $\pi_i$ for all $i > -n$ (in fact $\Tot$ can be computed by a spectral sequence with $1$-page $\pi_p(X^q)$). In particular, if $\mathcal D$ also has a $t$-structure such that there is some $d \ge 0$ with $F(\mathcal C_{\le0}) \subset \mathcal D_{\le d}$ then $F(\Tot(X^\bullet)) = \Tot(F(X^\bullet))$.

\paragraph{Kan Extensions.} Kan extensions are an extremely useful tool for constructing functors in ordinary category theory and in $\infty$-category theory. In fact a surprising amount of constructions in ($\infty$-)category theory can be reduced to Kan extensions. This is even more helpful in the $\infty$-categorical context, where one can almost never construct something explicitly. The basic idea is as follows: Suppose we are given a functor $\alpha\colon \mathcal C_0 \to \mathcal C$ of $\infty$-categories. For any other $\infty$-category $\mathcal D$ we get an induced functor $\alpha^*\colon \Fun(\mathcal C, \mathcal D) \to \Fun(\mathcal C_0, \mathcal D)$. A \emph{left (right) Kan extension} along $\alpha$ is a left (right) adjoint of $\alpha^*$, i.e. a way of ``extending'' any functor $\mathcal C_0 \to \mathcal D$ to a functor $\mathcal C \to \mathcal D$ along $\alpha$. If $\mathcal D$ has enough colimits (resp. limits) then left (resp. right) Kan extensions exist and can be computed as follows: Given a functor $F\colon \mathcal C_0 \to \mathcal D$, its left Kan extension $\overline{\mathcal F}\colon \mathcal C \to \mathcal D$ is computed by
\begin{align*}
	\overline{\mathcal F}(X) = \varinjlim_{Y \in (\mathcal C_0)_{X/}} F(Y),
\end{align*}
where $(\mathcal C_0)_{X/}$ is the $\infty$-category whose objects are pairs $(Y, f)$ with $Y \in \mathcal C_0$ and $f\colon \alpha(Y) \to X$ a morphism. For example, if $\mathcal C = \mathcal C_0^\triangleright$ then the left Kan extension along $\mathcal C_0 \injto \mathcal C$ computes the colimit over a diagram $\mathcal C_0 \to \mathcal D$. We refer the reader to \cite[Section 02Y1]{kerodon} for more information on Kan extensions.

\paragraph{Presentable $\infty$-Categories.} The category of abelian groups is of particularly nice shape: It admits all small limits and colimits and moreover every abelian group can be written as a colimit of copies of $\Z$. Most of the $\infty$-categories we encounter in this thesis allow a similar description, which is formalized by the notion of \emph{presentable $\infty$-categories} and studied by Lurie in \cite[\S5.5]{lurie-higher-topos-theory}. Roughly, an $\infty$-category $\mathcal C$ is presentable if for some regular cardinal $\kappa$ there is a full subcategory $\mathcal C_0 \subset \mathcal C$ of $\kappa$-compact objects such that every object in $\mathcal C$ can be written as an iterated small colimit of objects in $\mathcal C_0$. Here an object $X \in \mathcal C$ is called \emph{$\kappa$-compact} if the functor $\Hom(X, -)\colon \mathcal C \to \Ani$ preserves all $\kappa$-filtered colimits. If we want to refer to the cardinal $\kappa$ used in the definition of presentable $\infty$-category we can say that $\mathcal C$ is \emph{$\kappa$-compactly generated}.

Presentable $\infty$-categories are particularly nice to deal with. Most prominently, the adjoint functor theorem tells us that a functor $F\colon \mathcal C \to \mathcal D$ of presentable $\infty$-categories admits a right adjoint iff it preserves all small colimits and it admits a left adjoint iff it preserves all small limits and is accessible (i.e. preserves $\kappa$-filtered colimits for some $\kappa$). This is yet another very useful way to construct functors in the $\infty$-categorical setting.

\paragraph{Higher Algebra.} Many of the $\infty$-categories we deal with in this thesis come equipped with a (symmetric) monoidal structure, i.e. a product operation $\tensor\colon \mathcal C \cprod \mathcal C \to \mathcal C$ which is coherently associative (and commutative). Encoding all the required coherences of such a product operation in the $\infty$-categorical setting is a non-trivial task (even in ordinary category theory this is a bit subtle). One way to do this is through the theory of $\infty$-operads which is studied in \cite[\S2]{lurie-higher-algebra}. Even though we use symmetric monoidal $\infty$-categories throughout the whole thesis, a good understanding of its precise definition is only required in \cref{sec:andesc} and not in \cref{sec:ri-pi}. We therefore refrain from attempting to give an introduction to $\infty$-operads and the theory of modules and algebras in that setting and instead refer the interested reader to \cite{lurie-higher-algebra}.

It is noteworthy that in this thesis we will make use of animated rings as opposed to $\mathbb E_\infty$-rings over $\Z$ (as a generalization of commutative rings to the $\infty$-categorical setting). The former can be defined by adding homotopies as higher morphisms to the category of simplicial commutative rings, while the latter are commutative algebra objects in the symmetric monoidal $\infty$-category $\D(\Z)$. The main difference between these two concepts is that the ring structure on an animated ring is defined to be commutative in a strict sense, while on an $\mathbb E_\infty$-ring it is only commutative up to coherent homotopies. There is a conservative limit- and colimit-preserving forgetful functor from animated rings to $\mathbb E_\infty$-rings over $\Z$ which induces isomorphisms on the associated $\infty$-categories of modules, so the difference between the two concepts often does not matter too much. The main reason why we chose to work with animated rings is that they seem to better extend classical algebraic geometry. For example, the free generators of the $\infty$-category of animated rings are the polynomial rings $\Z[T_1, \dots, T_n]$, while the free generators of the $\infty$-category of $\mathbb E_\infty$-rings over $\Z$ are more complicated (and not static).

\subsection{Condensed Almost Mathematics} \label{sec:intro.almcond}

In this thesis we will make heavy use of the recently developed theory of condensed mathematics by Clausen--Scholze \cite{condensed-mathematics,scholze-analytic-spaces}. In order to apply this theory to perfectoid geometry we need to combine it with almost mathematics, as introduced by Faltings and later studied by Gabber-Ramero \cite{almost-ring-theory}. The first part of \cref{sec:andesc} is devoted to providing a thorough treatment of the resulting ``Condensed Almost Mathematics''. In the following we provide an overview of this theory, with the aim of giving the reader enough intuition at hand so that they may skip \cref{sec:andesc} entirely if they are only interested in the $p$-adic 6-functor formalism.

For a start let us ignore the almost setting and instead consider the ``honest'' condensed setting introduced in \cite{condensed-mathematics,scholze-analytic-spaces}. At its core, condensed mathematics provides a new definition of topological spaces which behaves much better in the context of algebraic structures. This is achieved by redefining what a ``topology'' on a set should be and in particular allow non-trivial topological structures on the set with one element. This resolves one of the most annoying properties of topological spaces, namely that the category of topological abelian groups is not abelian. If it was then one would have a good notion of kernels and cokernels, but the map $\R_\disc \to \R$ (where $\R_\disc$ denotes the real numbers equipped with the discrete topology) is bijective on underlying sets and hence has vanishing kernel and cokernel, despite it not being an isomorphism. This issue is resolved by condensed mathematics, as the same map will have a non-trivial cokernel in the category of condensed abelian groups (i.e. a non-trivial ``topological'' abelian group with underlying set $\{ 0 \}$).

Let us come to the definition of condensed sets. The basic observation is that all ``nice'' topological spaces $X$ are determined by the continuous maps $Y \to X$ they allow from all compact Hausdorff spaces $Y$. For example, let $Y = \N_\infty := \N \union \{ \infty \}$ be the one-point compactification of $\N$. Then a continuous map $\N_\infty \to X$ is the same as a converging sequence in $X$, hence the set of continuous maps $\N_\infty \to X$ determines the convergence of sequences in $X$. This motivates to define a condensed set $X$ to be a sheaf of sets on the site of compact Hausdorff spaces (where covers are finite jointly surjective families of continuous maps), where we view its value $X(Y)$ on any compact Hausdorff space $Y$ as the set of maps ``$Y \to X$'' (by confusing $Y$ with the sheaf it represents, this statement becomes literally true by Yoneda). We now note that the site of compact Hausdorff spaces admits a basis by the so-called \emph{extremally disconnected sets}. These are the profinite sets $Y$ such that every surjection $Z \to Y$ from a compact Hausdorff space $Z$ admits a section. In particular, a sheaf on the site of compact Hausdorff spaces is the same as a sheaf on the site of extremally disconnected sets. Due to the described splitting behavior, the latter category of sheaves admits a very simple description:

\begin{definition} \label{def:intro-condensed-objects}
Let $\mathcal C$ be an $\infty$-category. For every strong limit cardinal $\kappa$ we define $\Cond(\mathcal C)_\kappa$ to be the $\infty$-category of contravariant functors from the category of $\kappa$-small extremally disconnected sets to $\mathcal C$ which map finite disjoint unions to products. We define
\begin{align*}
	\Cond(\mathcal C) := \varinjlim_\kappa \Cond(\mathcal C)_\kappa,
\end{align*}
where the transition functors are the pullback functors. We call $\Cond(\mathcal C)$ the $\infty$-category of \emph{condensed objects in $\mathcal C$} and $\Cond(\mathcal C)_\kappa \subset \Cond(\mathcal C)$ the full subcategory of \emph{$\kappa$-condensed objects in $\mathcal C$}.
\end{definition}

The need of using strong limit cardinals in \cref{def:intro-condensed-objects} is a set-theoretical annoyance coming from the fact that the category of extremally disconnected sets is not small. Essentially, $\Cond(\mathcal C)$ is the $\infty$-category of $\mathcal C$-valued sheaves on the site of extremally disconnected sets.

Setting $\mathcal C = \catset$ we obtain the category $\Cond(\catset)$ of condensed sets, which provides the promised replacement of topological spaces. In this thesis we are mainly interested in condensed algebraic structures:

\begin{definition}
Let $A$ be a classical ring. We denote by $\D(A)$ the derived $\infty$-category of condensed $A$-modules. For every strong limit cardinal $\kappa$ we denote by $\D(A)_\kappa \subset \D(A)$ the full subcategory of $\kappa$-condensed objects. In particular we obtain $\D(A)_\omega$, the $\infty$-category of \emph{discrete $A$-modules}.
\end{definition}

The objects of $\D(A)$ are complexes of sheaves of $A$-modules on the site of extremally disconnected sets; intuitively one should see these as ``derived topological $A$-modules''. The $\infty$-category $\D(A)_\omega$ is the classical derived $\infty$-category of $A$-modules. One observes that the $\infty$-category $\D(A)$ is compactly generated (up to set-theoretic issues) by the free generators $A[S]$ for extremally disconnected sets $S$. In other words, every object in $\D(A)$ can be written as a colimit of copies of $A[S]$ for varying $S$. Here $A[S]$ is the sheafification of the presheaf $T \mapsto A[\Hom(T, S)]$.

We have now defined a good notion of ``(derived) topological $A$-modules'', but the story does not end there. In practice one often needs to pass to a full subcategory of $\D(A)$ consisting of ``complete'' $A$-modules. One reason why one needs this comes from the application to algebraic geometry: One of the main applications in \cite{condensed-mathematics} is a full 6-functor formalism for (condensed, derived) quasicoherent sheaves on schemes. For such a 6-functor formalism to exist, one needs to construct, for every open immersion $j\colon \Spec A[1/f] \injto \Spec A$ a functor $j_!\colon \D(A[1/f]) \to \D(A)$ which is left adjoint to the pullback $j^* = - \tensor_A A[1/f]$. However, such a functor $j_!$ cannot exist in general because usually $j^*$ does not preserve limits. To solve this issue we can pass to the full subcategory of solid modules:

\begin{definition}
Let $A$ be a classical ring. For every profinite set $S = \varprojlim_i S_i$ we let
\begin{align*}
	A_\solid[S] := \varinjlim_{A' \subset A} \varprojlim_i A'[S_i],
\end{align*}
where $A'$ ranges over all subrings of $A$ which are of finite type over $\Z$. Note that there is a canonical map $A[S] \to A_\solid[S]$. An $A$-module $M \in \D(A)$ is called \emph{solid} if for all profinite sets $S$ the map $\IHom(A_\solid[S], M) \isoto \IHom(A[S], M)$ is an isomorphism in $\D(A)$. Here $\IHom$ denotes the (derived) internal hom in $\D(A)$. We denote by
\begin{align*}
	\D_\solid(A) \subset \D(A)
\end{align*}
the full subcategory spanned by the solid $A$-modules.
\end{definition}

One of the main results of \cite{condensed-mathematics} is that the stable $\infty$-category $\D_\solid(A)$ enjoys surprisingly nice properties. Among others it admits all small limits and colimits and the inclusion $\D_\solid(A) \injto \D(A)$ is stable under them. Also, an $A$-module $M \in \D(A)$ is solid if and only if all $\pi_n M$ are solid; in fact, $\D_\solid(A)$ is the derived $\infty$-category of its heart. The inclusion $\D_\solid(A) \injto \D(A)$ admits a left adjoint
\begin{align*}
	- \tensor_A A_\solid\colon \D(A) \to \D_\solid(A)
\end{align*}
sending $A[S] \mapsto A_\solid[S]$. Given any map $A \to B$ of classical rings we get a base-change functor
\begin{align*}
	- \tensor_{A_\solid} B_\solid\colon \D_\solid(A) \to \D_\solid(B)
\end{align*}
which is left adjoint to the forgetful functor $\D_\solid(B) \to \D_\solid(A)$. In fact we can compute $M \tensor_{A_\solid} B_\solid = (M \tensor_A B) \tensor_B B_\solid$. If $B = A[1/f]$ then the functor $j^* = - \tensor_{A_\solid} B_\solid$ preserves all small limits and does indeed admit a left adjoint $j_!$. It is now mostly formal to glue the $\infty$-categories $\D_\solid(A)$ in order to get the $\infty$-category $\D_\solid(X)$ of \emph{solid quasicoherent sheaves} on a scheme $X$, and to construct the 6-functor formalism for $\D_\solid(X)$. See \cite[\S11]{condensed-mathematics} or \cref{rslt:scheme-6-functor} for a precise formulation.

It is sometimes convenient to have a slightly more general version of solid modules at hand:

\begin{definition}
Let $A \to B$ be a map of classical rings. We let
\begin{align*}
	\D_\solid(B, A) \subset \D(B)
\end{align*}
denote the full subcategory spanned by those $M \in \D(B)$ such that $M$ is solid as an $A$-module via the forgetful functor $\D(B) \to \D(A)$.
\end{definition}

One checks that $\D_\solid(B, A)$ depends only on the integral closure of the image of $A$ in $B$. In fact, given any classical discrete Huber pair $(A, A^+)$ (i.e. $A$ is a classical ring and $A^+ \subset A$ is an integrally closed subring) then one should view $\D_\solid(A, A^+)$ as the $\infty$-category of solid quasicoherent sheaves on the discrete adic space $\Spa(A, A^+)$. The category of discrete adic spaces contains the category of schemes as a full subcategory, and it is often handy to have the more general theory of discrete adic spaces available. We can glue $\D_\solid(A, A^+)$ to get the $\infty$-category of solid quasicoherent sheaves $\D_\solid(X)$ on any discrete adic space $X$, and the 6-functor formalism extends to that case (see \cref{rslt:scheme-6-functor}). This finishes the overview of \cite{condensed-mathematics}.

Let us now introduce almost mathematics, i.e. combine \cite{condensed-mathematics} with \cite{almost-ring-theory}. The general paradigm of almost mathematics is that we care about objects and elements only up to things that are annihilated by some ideal $\mm$. Following \cite{almost-ring-theory} this idea can be formalized as follows:

\begin{definition}
An \emph{almost setup} is a pair $(V,\mm)$ where $V$ is a classical ring and $\mm \subset V$ is an ideal such that $\mm^2 = \mm$ and such that $\widetilde\mm := \mm \tensor_V \mm$ is flat over $V$ (here we use the non-derived tensor product).
\end{definition}

\begin{definition}
Let $(V,\mm)$ be an almost setup. A $V$-module $M \in \D(V)$ is called \emph{almost zero} if multiplication by $\varepsilon$ is zero on $M$ for all $\varepsilon \in \mm$. We say that a map $M \to N$ in $\D(V)$ is an \emph{almost isomorphism} if $\fib(M \to N)$ is almost zero. Let $W$ be the collection of almost isomorphisms in $\D(V)$ and define
\begin{align*}
	\D(V,\mm) := \D(V)[W^{-1}],
\end{align*}
the $\infty$-category of \emph{(condensed derived) $(V,\mm)$-modules}. If $\mm$ is clear from context then we also write $\D(V^a)$ for $\D(V,\mm)$ and refer to it as the $\infty$-category of \emph{(condensed derived) almost $V$-modules}.
\end{definition}

We study the $\infty$-category $\D(V^a)$ in \cref{sec:andesc.almmath} and generalize many of the results in \cite{almost-ring-theory} to the condensed derived setting. For example we get the following result (see \cref{rslt:properties-of-derived-almost-V-modules}):

\begin{proposition}
Fix an almost setup $(V,\mm)$
\begin{propenum}
	\item There is a natural ``almostification'' functor
	\begin{align*}
		(-)^a\colon \D(V) \to \D(V^a), \qquad M \mapsto M^a
	\end{align*}
	which is $t$-exact and preserves all small limits and colimits and the symmetric monoidal operation $\tensor$.

	\item The functor $(-)^a$ admits a right adjoint
	\begin{align*}
		(-)_*\colon \D(V^a) \to \D(V), \qquad M \mapsto M_*,
	\end{align*}
	which is left $t$-exact and fully faithful. Moreover, for any $M \in \D(V)$ there is a natural isomorphism $(M^a)_* = \IHom_V(\widetilde\mm, M)$.

	\item The functor $(-)^a$ admits a left adjoint
	\begin{align*}
		(-)_!\colon \D(V^a) \to \D(V), \qquad M \mapsto M_!,
	\end{align*}
	which is $t$-exact and fully faithful. It can be computed as $M_! = \widetilde\mm \tensor_V M_*$.
\end{propenum}
\end{proposition}

From now on we fix an almost setup $(V,\mm)$. The functors $(-)^a$, $(-)_*$ and $(-)_!$ allow us to move back and forth between the almost world and the non-almost world, thereby generalizing many results from \cite{condensed-mathematics,scholze-analytic-spaces} to the almost setting. For example, it follows immediately that every $M \in \D(V^a)$ can be written as a colimit of copies of the free generators $V^a[S] := V[S]^a$. However, it is in general no longer true that the generators $V^a[S]$ are compact (because $(-)_*$ does not necessarily commute with colimits). In particular $\D(V^a)$ is not compactly generated (but still presentable up to set-theoretic issues).

Given a classical $V$-algebra $A$ we get an associated $V^a$-algebra $A^a$ (i.e. a commutative algebra object in the heart of $\D(V^a)$) and thus also the $\infty$-category $\D(A^a)$ of \emph{almost $A$-modules}. Similar remarks as for $\D(V^a)$ apply, for example we have the almostification $(-)^a\colon \D(A) \to \D(A^a)$ with right adjoint $(-)_*$ and left adjoint $(-)_!$. One can also generalize solid modules:

\begin{definition}
Let $(V,\mm)$ be an almost setup and let $A \to B$ be a morphism of classical rings, where $B$ is a $V$-algebra. Then we define
\begin{align*}
	\D_\solid(B^a, A) \subset \D(B^a)
\end{align*}
to be the essential image of $\D_\solid(B, A)$ under the almostification functor $\D(B) \to \D(B^a)$.
\end{definition}

We show that $\D_\solid(B^a, A)$ satisfies similarly nice properties as $\D_\solid(B, A)$ and sometimes also denote it by $\D_\solid^a(B, A)$. Be aware that while this $\infty$-category only depends on $B$ via its almostification $B^a$, it certainly does depend on the non-almost version of $A$ (in fact, since $A$ is not required to live over $V$, there is in general no almost version $A^a$ of $A$). In the case $A = B$ we also abbreviate $\Dqcohri(B) = \D^a_\solid(B, B)$. This finally provides a precise definition of the $\infty$-category $\Dqcohri(A^+/\pi)$ associated to an affinoid perfectoid space $\Spa(A, A^+)$ with pseudouniformizer $\pi$ (where the almost setup is $(A^+, A^{\circ\circ})$ with $A^{\circ\circ} \subset A^+$ being the ideal of topologically nilpotent elements).

With the above explanations the reader should have enough knowledge in order to understand the solid framework used in \cref{sec:ri-pi}. However, in \cref{sec:andesc} we push the theory a bit further by also generalizing the analytic geometry developed in \cite{scholze-analytic-spaces}. Let us briefly touch upon what this entails. The basic idea comes from the observation that the $\infty$-category $\D_\solid(A)$ defined above is built from a simple recipe: For every extremally disconnected set $S$ we specify a solid version of the generator $A[S]$, called $A_\solid[S]$, and then we let $\D_\solid(A) \subset \D(A)$ be the full subcategory of those $M \in \D(A)$ such that $\IHom(A_\solid[S], M) = \IHom(A[S], M)$. This can be generalized:

\begin{definition}
Let $(V,\mm)$ be an almost setup.
\begin{defenum}
	\item An \emph{analytic ring over $V^a$} is a pair $\mathcal A = (\underline{\mathcal A}, \mathcal M_{\mathcal A})$, where $\underline{\mathcal A}$ is a condensed animated ring over $V^a$ (see \cref{def:almost-rings}) and $\mathcal M_{\mathcal A}$ is a functor from the category of extremally disconnected sets to $\D_{\ge0}(\underline{\mathcal A})$ together with a natural transformation $\underline{\mathcal A}[S] \to \mathcal M_{\mathcal A}(S)$, satisfying certain properties as in \cref{def:analytic-ring}.

	We will usually denote $\mathcal A[S] := \mathcal M_{\mathcal A}(S)$, removing $\mathcal M_{\mathcal A}$ from the notation.

	\item Let $\mathcal A$ be an analytic ring over $V^a$. We let
	\begin{align*}
		\D(\mathcal A) \subset \D(\underline{\mathcal A})
	\end{align*}
	denote the full subcategory spanned by those $M \in \D(\underline{\mathcal A})$ such that for all extremally disconnected sets $S$ the natural map $\IHom(\mathcal A[S], M) \isoto \IHom(\underline{\mathcal A}[S], M)$ is an isomorphism.
\end{defenum}
\end{definition}

For every map of classical rings $A \to B$ such that $B$ is a $V$-algebra the above definitions provide an analytic ring $(B^a, A)_\solid$ over $V^a$ with $(B^a, A)_\solid[S] = B^a \tensor_{A^a} A_\solid[S]^a$. The general theory of analytic rings is a vast generalization of this concept:
\begin{itemize}
	\item We allow the base ring $\underline{\mathcal A}$ to be a derived ring, i.e. it might have nontrivial $\pi_n(A)$ for $n > 0$. This is already useful in the solid context: Given a diagram of classical rings $B \to A \from C$ the (derived) tensor product $B \tensor_A C$ is in general a non-static ring. The above definition of analytic rings allows us to define the associated analytic ring $(B \tensor_A C)_\solid$.

	\item We allow the base ring $\underline{\mathcal A}$ to have a non-discrete topology, so we could take rings like $\Z_p$ or $\R$ as our basis. In this thesis this will only be relevant marginally (in the setting of $p$-adically complete rings), but in \cite{scholze-analytic-spaces} a large effort is spent on defining an analytic ring structure on $\R$ which ultimately allows to fit the classical theory of Banach spaces and manifolds in the condensed world.
\end{itemize}
In \cref{sec:andesc} we study the $\infty$-category $\AnRing_{V^a}$ of analytic rings over $V^a$ and show that it behaves in a very similar way as the $\infty$-category $\AnRing := \AnRing_\Z$ introduced in \cite{scholze-analytic-spaces}. One can then ``glue'' analytic rings in order to define analytic spaces over $V^a$. Restricting to analytic rings of the form $(A^a, A^+)_\solid$ (where $A^a$ is an \emph{animated} discrete ring over $V^a$) we obtain a theory of (derived) discrete adic spaces over $V^a$. In \cref{rslt:scheme-6-functor,rslt:lower-shriek-commutes-with-almost-localization} we provide a 6-functor formalism for solid quasicoherent sheaves on these spaces, thereby extending the results in \cite{condensed-mathematics}.

Even though in our 6-functor formalism for $\ri^{+a}_X/\pi$-sheaves we only need to work with analytic rings of the form $(A^+/\pi)^a_\solid$, it is useful to have the more general theory of analytic rings at our disposal. We believe that the theory of analytic spaces over $V^a$ (and the descent formalism explained in \cref{sec:intro.descent} below) will also have applications outside our theory of $\ri^{+a}_X/\pi$-modules.

\subsection{Descent Theory for Solid Modules} \label{sec:intro.descent}

The second half of \cref{sec:andesc} is devoted to developing a general descent formalism for modules over analytic almost rings (in particular for solid modules over classical rings) which is based on ideas by Mathew \cite[\S3]{akhil-galois-group-of-stable-homotopy} and Bhatt--Scholze \cite[\S11]{bhatt-scholze-witt}. This formalism provides the necessary tools for proving the descent of $(A^+/\pi)^a_\solid$-modules on affinoid perfectoid spaces (i.e. \cref{rslt:intro-def-of-qcoh-Lambda-modules}) which lies at the heart of our 6-functor formalism. We now give an overview of our descent formalism.

In the following we formulate everything in terms of general analytic rings, which were quickly introduced at the end of \cref{sec:intro.almcond}. However, our main applications lie in the world of solid modules, so the reader is free to assume that all the appearing analytic rings are of the form $(B^a, A)_\solid$ (but be aware that all tensor products are derived, so in general we need to allow $B^a$ to be a derived ring).

The basic question we want to answer is as follows: Under what assumptions on a map $f\colon \mathcal A \to \mathcal B$ of analytic rings over some almost setup $(V,\mm)$ can we guarantee that modules descent along $f$? In order to answer this question we first need a proper definition of ``descent''. Following the explanations on gluing $\infty$-categories in \cref{sec:intro.infty} we arrive at the following definition:

\begin{definition}
Let $(V,\mm)$ be an almost setup and let $f\colon \mathcal A \to \mathcal B$ be a map of analytic rings over $V^a$. We say that \emph{modules descent along $f$} if the natural functor
\begin{align*}
	\D(\mathcal A) \isoto \varprojlim_{[n]\in\Delta} \D(\mathcal B^{\tensor n+1})
\end{align*}
is an equivalence. Here $\mathcal B^{\tensor n}$ denotes the $n$-fold tensor product of $\mathcal B$ over $\mathcal A$.
\end{definition}

The general strategy for showing that modules descend along a given map $f\colon \mathcal A \to \mathcal B$ is as follows: Note that there is a natural pair of adjoint functors
\begin{align*}
	F\colon \D(\mathcal A) \rightleftarrows \varprojlim_{[n]\in\Delta} \D(\mathcal B^{\tensor n+1}) \noloc G,
\end{align*}
where $F$ is given by $F(M) = (M \tensor_{\mathcal A} \mathcal B^{\tensor \bullet+1})$ and $G$ is given by $G(N^\bullet) = \Tot(N^\bullet) := \lim_{[n]\in\Delta} N^n$. We want to show that $F$ and $G$ are equivalences, which by adjunctions reduces to the following two claims:
\begin{enumerate}[(a)]
	\item $F$ is conservative, i.e. if $M \in \D(\mathcal A)$ satisfies $F(M) = 0$ then $M = 0$. This boils down to the claim that $- \tensor_{\mathcal A} \mathcal B$ is conservative, which is usually not too hard to check (this is essentially saying that the associated map of analytic spaces is surjective).

	\item $G$ is fully faithful, i.e. the counit $FG \isoto \id$ is an isomorphism. This boils down to showing that for all $N^\bullet \in \varprojlim_{[n]\in\Delta} \D(\mathcal B^{\tensor n+1})$ the natural map $\Tot(N^\bullet) \tensor_{\mathcal A} \mathcal B \to N^0$ is an isomorphism. It follows easily from the fact that the map $\mathcal B \to \mathcal B \tensor_{\mathcal A} \mathcal B$ is split that we have $N^0 = \Tot(N^\bullet \tensor_{\mathcal A} \mathcal B)$. Hence the fully faithfulness of $G$ reduces to showing that the map
	\begin{align*}
		\Tot(N^\bullet) \tensor_{\mathcal A} \mathcal B \isoto \Tot(N^\bullet \tensor_{\mathcal A} \mathcal B)
	\end{align*}
	is an isomorphism, i.e. that $- \tensor_{\mathcal A} \mathcal B$ commutes with the totalization of $N^\bullet$.
\end{enumerate}
We thus need to find a good condition on $f$ under which it is guaranteed that $- \tensor_{\mathcal A} \mathcal B$ commutes with the totalization of $N^\bullet$ as in (b). We know that $- \tensor_{\mathcal A} \mathcal B$ is an exact functor and hence commutes with all finite limits. The crucial observation by Mathew \cite[\S3]{akhil-galois-group-of-stable-homotopy} is that one can often guarantee that the limit $\Tot(N^\bullet)$, while being an infinite limit by definition, is ``essentially finite''. Namely, the cosimplicial object $N^\bullet \tensor_{\mathcal A} \mathcal B$ is split (because $\mathcal B \to \mathcal B \tensor_{\mathcal A} \mathcal B$ splits), so that its totalization is essentially finite. It is therefore enough to ensure that $N^\bullet$ can be built from $N^\bullet \tensor_{\mathcal A} \mathcal B$ in finitely many simple steps. This leads to the following definition:

\begin{definition}[{cf. \cref{def:descendable-morphism-of-analytic-spaces}}] \label{def:intro.descendable-morphism}
Let $(V,\mm)$ be an almost setup. A map $f\colon \mathcal A \to \mathcal B$ of analytic rings over $V^a$ is ``\emph{descendable}'' if $f$ is steady (see \cref{def:steady-morphism-of-analytic-rings}; this is always satisfied in the solid discrete setting) and the functor $\id\colon \D(\mathcal A) \to \D(\mathcal A)$ can be obtained by a finite sequence of compositions, finite (co)limits and retracts from the functor $- \tensor_{\mathcal A} \mathcal B\colon \D(\mathcal A) \to \D(\mathcal A)$.
\end{definition}

\begin{remark}
In \cref{def:intro.descendable-morphism} we put ``descendable'' in quotes because the provided definition is not exactly the right one, even though it captures the main idea. The problem is that we should not work with the $\infty$-category $\Fun(\D(\mathcal A), \D(\mathcal A))$ of functors $\D(\mathcal A) \to \D(\mathcal A)$ but instead consider $\D(\mathcal A)$-enriched functors $\D(\mathcal A) \to \D(\mathcal A)$. Only then it is true that descendable morphisms are stable under base-change (see below). The precise definitions of the enriched functor $\infty$-categories are a bit subtle; we have devoted the whole \cref{sec:andesc.endofun} to dealing with these subtleties.
\end{remark}

Note that the descendability condition is a condition that only depends on finite limits and retracts of the functors $- \tensor_{\mathcal A} \mathcal B^{\tensor n}$ and thus does not need us to deal with any infinite limits. It does the expected job:

\begin{proposition} \label{rslt:intro.properties-of-descendability}
\begin{propenum}
	\item Modules descend along every descendable morphism.

	\item Descendable morphisms are stable under composition and base-change.
\end{propenum}
\end{proposition}

Part (i) of \cref{rslt:intro.properties-of-descendability} was essentially proved above (see \cref{rslt:descendable-implies-limit-of-categories} for details). Part (ii) is mostly straightforward but the stability under base-change heavily relies on the correct definition of endofunctors (see \cref{rslt:descendability-stable-under-composition,rslt:descendable-stable-under-base-change} for the proofs).

It is time for an example. Probably the first thing one might look for is faithfully flat descent. It is currently unclear whether faithfully flat maps of rings are descendable (even in Mathew's version), but it is true in the finitely presented case (see \cref{rslt:fppf-cover-of-rings-is-descendable}):

\begin{proposition} \label{rslt:intro.fppf-descent-for-rings}
Let $A \to B$ be a faithfully flat and finitely presented map of classical rings. Then the map $A_\solid \to B_\solid$ is descendable.
\end{proposition}

To prove \cref{rslt:intro.fppf-descent-for-rings} we first reduce to the case that $A$ and $B$ are of finite type over $\Z$, in which case the solid generators $A_\solid[S]$ and $B_\solid[S]$ are isomorphic to products $\prod_I A$ and $\prod_I B$. Using this fact one can reduce the descendability condition to the similar condition of Mathew and hence deduce the claim from Mathew's result on fppf covers \cite[Proposition 3.31]{akhil-galois-group-of-stable-homotopy}.

To formulate the next example, we note that the above descendability condition generalizes to maps $f\colon Y \to X$ of analytic spaces over $(V,\mm)$ (just replace $\D(\mathcal A)$ by $\D(X)$ and $\D(\mathcal B)$ by $\D(Y)$ and be careful with the correct notion of endofunctors), so in particular we get a notion of descendable morphisms of schemes and discrete adic spaces. We have (see \cref{rslt:h-cover-is-descendable}):

\begin{theorem} \label{rslt:intro.h-cover-is-descendable}
Let $f\colon Y \to X$ be an h-cover of qcqs noetherian classical schemes. Then $f$ is descendable, so that solid quasicoherent sheaves satisfy descent along $f$.
\end{theorem}

We refer the reader to \cite[\S2]{bhatt-scholze-witt} for an introduction to the h-topology of schemes, but remark that every fppf cover and every proper surjective map is an h-cover. The proof of \cref{rslt:intro.h-cover-is-descendable} goes along the lines of the proof of \cite[Proposition 11.25]{bhatt-scholze-witt} using \cref{rslt:intro.fppf-descent-for-rings}.

One checks easily that the descendability condition is stable under almostification. Thus, using \cref{rslt:intro.h-cover-is-descendable} one should be able to deduce that the map $(A^+/\pi)^a_\solid \to (B^+/\pi)^a_\solid$ is descendable for every \emph{étale} map $\Spa(A, A^+) \to \Spa(B, B^+)$ of affinoid perfectoid spaces by observing that such an étale map is given by blow-ups and fppf covers on the formal model. This is not enough for our purposes however, because we need pro-étale descent. We thus need to understand how the descendability condition behaves under filtered colimits. In general a filtered colimit of descendable maps is not descendable anymore, but the situation becomes much nicer if we put a bound on the ``complexity'' of the descendability. Namely, to every descendable map $f\colon \mathcal A \to \mathcal B$ of analytic rings one can define the \emph{index of descendability} of $f$, which is a positive integer measuring how complicated it is to get from $- \tensor_{\mathcal A} \mathcal B$ to $\id$ (see \cref{def:index-of-descendable-map}).

\begin{example} \label{rslt:intro-examples-of-index-of-descendability}
If $f\colon \mathcal A \to \mathcal B$ splits then it is descendable of index $1$. If $A \to B$ is a faithfully flat and finitely presented map of classical rings then $A_\solid \to B_\solid$ is descendable of index $\le 2$.
\end{example}

\begin{theorem}[{cf. \cref{rslt:weakly-descendable-maps-properties}}] \label{rslt:intro-modules-descend-along-weakly-descendable-map}
Let $(V,\mm)$ be an almost setup such that $\mm$ is countably generated. Suppose $f\colon \mathcal A \to \mathcal B$ is a map of analytic rings over $V^a$ which can be written as a filtered colimit of descendable maps of bounded index. Then modules descend along $f$.
\end{theorem}

The proof of \cref{rslt:intro-modules-descend-along-weakly-descendable-map} is done in two steps. One first shows that the claim is true for \emph{countable} filtered colimits. Namely in this case $f$ is even descendable, as an explicit computation shows (here one uses that $R^n \varprojlim$ vanishes for countable cofiltered limits of abelian groups and $n \ge 2$). This allows one to reduce the claim to the case that $f$ is an $\omega_1$-filtered colimit of maps of analytic rings along each of which modules descend. In the second step one looks at the $\omega_1$-compact generators inside the module $\infty$-categories and uses that $\omega_1$-filtered colimits of $\infty$-categories commute with countable limits like $\varprojlim_{[n]\in\Delta}$.

Motivated by \cref{rslt:intro-modules-descend-along-weakly-descendable-map} we make the following definition:

\begin{definition}
Let $(V,\mm)$ be an almost setup such that $\mm$ is countably generated. A map $f\colon \mathcal A \to \mathcal B$ of analytic rings over $V^a$ is called \emph{weakly descendable} if it is an iterated $\omega_1$-filtered colimit of descendable maps.
\end{definition}

By \cref{rslt:intro-modules-descend-along-weakly-descendable-map} every filtered colimit of descendable maps of bounded index is weakly descendable. In particular we get the following corollary of \cref{rslt:intro-examples-of-index-of-descendability}:

\begin{proposition} \label{rslt:into-ind-fppf-cover-is-weakly-descendable}
Let $f\colon \Spec B \to \Spec A$ be an ind-fppf cover of affine schemes. Then $A_\solid \to B_\solid$ is weakly descendable, so in particular solid modules descend along this map.
\end{proposition}

\subsection{Descent on Affinoid Perfectoid Spaces} \label{sec:intro.descperfd}

With the preparations of the previous subsections at hand, we can finally sketch the proof of the central descent result \cref{rslt:intro-def-of-qcoh-Lambda-modules}. Here we will work with the case $\Lambda = \ri^+_X/\pi$ throughout, as the case of general $\Lambda$ follows from the stability of descendable maps under base-change (see \cref{rslt:intro.properties-of-descendability}). To shorten notation, let us denote by $\AffPerfd_\pi$ the category of pairs $(X, \pi)$ where $X = \Spa(A, A^+)$ is an affinoid perfectoid space and $\pi \in A^+$ is a pseudouniformizer. We will often omit $\pi$ from the notation.

Before we start applying descent results, we need to get one detail out of the way: Suppose we are given a map $Y = \Spa(B, B^+) \to X = \Spa(A, A^+)$ of affinoid perfectoid spaces and we can show that modules descend along the map $(A^+/\pi)^a_\solid \to (B^+/\pi)^a_\solid$ of analytic rings. Then this descent involves modules over the tensor products $(B^+/\pi)^a_\solid \tensor_{(A^+/\pi)^a_\solid} \dots \tensor_{(A^+/\pi)^a_\solid} (B^+/\pi)^a_\solid$, which a priori might have nothing to do with sheaves on $\Spa(B^n, B^{n+}) := Y \cprod_X \dots \cprod_X Y$, i.e. with modules over $(B^{n+}/\pi)^a_\solid$. One can show that (see \cref{rslt:base-change-for-affinoid-perfectoid})
\begin{align*}
	(B^{n+}/\pi)^a_\solid = (B^+/\pi)^a_\solid \tensor_{(A^+/\pi)^a_\solid} \dots \tensor_{(A^+/\pi)^a_\solid} (B^+/\pi)^a_\solid,
\end{align*}
so there is no problem here. This is one of the places where we really need to work with almost coefficients.

The first descent result is the following (see \cref{rslt:explicit-description-of-DqcohriX}):

\begin{proposition} \label{rslt:intro-descent-on-tot-disc-spaces}
There is a unique pro-étale sheaf of $\infty$-categories $X \mapsto \DqcohriX X$ on $\AffPerfd_\pi$ such that for every totally disconnected space $X = \Spa(A, A^+) \in \AffPerfd_\pi$ we have $\DqcohriX X = \Dqcohri(A^+/\pi)$.
\end{proposition}

Up to some details, the proof of \cref{rslt:intro-descent-on-tot-disc-spaces} essentially reduces to showing that if $f\colon Y = \Spa(B, B^+) \to X = \Spa(A, A^+)$ is a pro-étale map of totally disconnected perfectoid spaces with pseudouniformizer $\pi$ then modules descend along $(A^+/\pi)^a_\solid \to (B^+/\pi)^a_\solid$. In fact we claim that this holds for any affinoid perfectoid space $Y$ and any v-cover $f$. Namely, by some general observations on descendable morphisms (see \cref{sec:andesc.fsdesc}) we can essentially reduce to the case that $X$ is connected, i.e. of the form $\Spa(K, K^+)$ for some perfectoid field $K$ and some open and bounded valuation subring $K^+ \subset K$. But the map $K^+ \to B^+$ is torsion-free, so since $K^+$ is a valuation ring this map is automatically ind-fppf. We conclude by \cref{rslt:into-ind-fppf-cover-is-weakly-descendable}.

Next up we improve on \cref{rslt:intro-descent-on-tot-disc-spaces} in the following way (see \cref{rslt:compute-Dqcohri-for-fin-type-over-tot-disc}):

\begin{theorem} \label{rslt:intro-descent-on-fin-type-over-tot-disc}
Suppose $X = \Spa(A, A^+) \in \AffPerfd_\pi$ is of weakly perfectly finite type (see \cref{def:perfectly-finite-type-in-AffPerf}) over some totally disconnected space. Then $\DqcohriX X = \Dqcohri(A^+/\pi)$, where $\DqcohriX X$ is defined as in \cref{rslt:intro-descent-on-tot-disc-spaces}
\end{theorem}

The proof of \cref{rslt:intro-descent-on-fin-type-over-tot-disc} is a bit subtle; it is the key step in proving \cref{rslt:intro-def-of-qcoh-Lambda-modules}. For simplicity let us assume that $X$ is of weakly perfectly finite type over $\Spa(K, K^+)$ for some perfectoid field $K$ (the general case can essentially be reduced to this case by passing to connected components). In this case we argue in the following two steps:
\begin{enumerate}[1.]
	\item We observe that $A^+$ has essentially finite global dimension (in the usual sense of commutative algebra). It follows from the definition of descendable maps that for every cover $Y = \Spa(B, B^+) \surjto X$ by some totally disconnected space the map $(A^+/\pi)^a_\solid \to (B^+/\pi, A^+/\pi)^a_\solid$ is descendable (see \cref{rslt:descendability-for-fin-global-dimension} for the general criterion). It is easy to see that $(B^+/\pi, A^+/\pi)^a = (B'^+/\pi)^a_\solid$, where $\Spa(B', B'^+) := \overline Y^{/X}$ is the relative compactification (see \cite[\S18]{etale-cohomology-of-diamonds}). It now follows formally that \cref{rslt:intro-descent-on-fin-type-over-tot-disc} reduces to showing that $\DqcohriX{\overline Y^{/X}} = \Dqcohri(B'^+/\pi)$. This claim in turn can essentially be reduced to the case that $Y$ is connected, i.e. of the form $Y = \Spa(K', K'^+)$

	\item By the previous step we are now in the following position: Let $Y' = \Spa(K', K'^+) \to Y = \Spa(K, K^+)$ be a v-cover, where $K$ and $K'$ are perfectoid fields with open and bounded valuation subrings $K^+ \subset K$ and $K'^+ \subset K'$, and $K'$ has finite topological transcendence degree over $K^+$; we need to show that $\DqcohriX{\overline Y'^{/Y}} = \Dqcohri(K'^+/\pi, K^+/\pi)$. To prove this, we let $Z = \Spa(C, C^+) \to \overline Y'^{/Y}$ be the w-localization (see \cite[Proposition 7.12]{etale-cohomology-of-diamonds}). This is a pro-étale cover by the totally disconnected space $Z$, which reduces the claim formally to showing that modules descend along $(K'^+/\pi, K^+/\pi)^a_\solid \to (C^+/\pi)^a_\solid$. In fact we claim that the latter map is weakly descendable; this follows from a general descendability result on open covers of Zariski-Riemann spaces (see \cref{rslt:descendability-of-open-cover-of-ZR-space}) which in turn is based on the fppf descent result \cref{rslt:intro.fppf-descent-for-rings}.
\end{enumerate}
The final step in the proof of \cref{rslt:intro-def-of-qcoh-Lambda-modules} is v-descent (see \cref{rslt:v-hyperdescent-for-Dqcohri}):

\begin{theorem} \label{rslt:intro-v-hyperdescent}
The assignment $X \mapsto \DqcohriX X$ from \cref{rslt:intro-descent-on-tot-disc-spaces} satisfies descent along all v-hypercovers in $\AffPerfd_\pi$.
\end{theorem}

The key observation in the proof of \cref{rslt:intro-v-hyperdescent} is that if $\Spa(B, B^+) \to X = \Spa(A, A^+)$ is a map of affinoid perfectoid spaces such that $X$ is totally disconnected, then the map $(A^+/\pi)^a_\solid \to (B^+/\pi)^a_\solid$ of analytic rings has Tor dimension $\le 1$. This result can by \cref{rslt:intro-descent-on-fin-type-over-tot-disc} be reduced to the following abstract result (see \cref{rslt:flat-map-of-val-rings-has-sld-Tor-dim-1}):

\begin{theorem} \label{rslt:intro-val-rings-Tor-dim}
Let $V \to V'$ be a flat map of valuation rings. Then $V_\solid \to V'_\solid$ has Tor dimension $\le 1$.
\end{theorem}

The proof of \cref{rslt:intro-val-rings-Tor-dim} is a fun application of the 6-functor formalism for quasicoherent sheaves on schemes. We encourage the reader to have a look at \cref{rslt:map-of-fields-is-analytically-flat}, where the general idea is discussed in the case that $V$ and $V'$ are fields (in which case $V_\solid \to V'_\solid$ is even flat).

\subsection{Notation and Conventions} \label{sec:intro.notation}

In this thesis we fully embrace the language of $\infty$-categories and condensed mathematics, as this is the most natural framework in which to formulate our results. In particular, by a \emph{ring} we mean a condensed animated ring and similarly a module over a ring $A$ is an object $M \in \D(A)$ in the derived $\infty$-category of condensed $A$-modules. Moreover, by a scheme we mean a derived scheme (see \cref{sec:andesc.scheme}). In many cases we put more assumptions on the rings and modules we deal with, according to the following terminology:
\begin{itemize}
	\item A ring/module is called \emph{discrete} if it is discrete in the sense of condensed mathematics.

	\item A ring/module is called \emph{static} if it is concentrated in degree $0$, i.e. has no ``derived'' data. In a similar fashion, a module is called \emph{connective} if it is concentrated in homological degree $\ge 0$.

	\item A ring/module is called \emph{classical} if it is discrete and static. Similarly, a scheme is called classical if it can be covered by affine schemes of the form $\Spec A$ for classical (i.e. non-derived) rings $A$.
\end{itemize}
All the functors we consider will always be their derived version (unless explicitly stated otherwise), so for example $\tensor$ denotes the derived tensor product of modules (even if we plug in classical modules on both sides). Moreover, we usually use homological language instead of cohomological language, i.e. we refer to the homotopy groups $\pi_n(M) = H^{-n}(M)$ of a module $M \in \D(A)$.

Regarding the $\infty$-categorical language, we deviate slightly from Lurie's books and denote by $\Ani$ (the $\infty$-category of \emph{anima}) what Lurie denotes $\catS$ (the $\infty$-category of spaces). Moreover, we denote by $\infcatinf^\tensor$ the $\infty$-category of symmetric monoidal $\infty$-categories and by $\infcatinf^\ocircle$ the $\infty$-category of monoidal $\infty$-categories.

\subsection{Acknowledgements}

I heartily thank Peter Scholze, my thesis advisor, for suggesting the topic and being a great mentor throughout. It has been a wonderful time being a Ph.D. student under his care and this thesis would not have been possible without his deep insights and continuous support. I would also like to thank Dustin Clausen for explaining some of the subtleties of condensed mathematics and $\infty$-categories to me, as well as pointing me to the $\infty$-category of correspondences and its usefulness in constructing a 6-functor formalism. I am grateful to K\k{e}stutis \v{C}esnavi\v{c}ius for his insights on valuation rings and resolutions of singularities, to Jacob Lurie for explaining some subtleties about enriched $\infty$-categories to me and to David Hansen for many helpful discussions about possible applications of this thesis as well as finding numerous typos. In addition I am very thankful to Eugen Hellmann and the organizers of the RAMpAGe seminar for letting me present my work. Lastly, I would like to thank the university of Bonn and specifically the members of the arithmetic geometry group for providing such an enjoyable environment during my Ph.D. studies, both mathematically and socially.

\clearpage
\section{Analytic Spaces and Descent} \label{sec:andesc}

The goal of this section is twofold: The first goal is to generalize the theory of condensed mathematics and analytic geometry developed by Clausen--Scholze in \cite{condensed-mathematics,scholze-analytic-spaces} to the setting of almost mathematics. The second goal is to develop a descent formalism in this setup and prove various descent results (on schemes and discrete adic spaces) which will later be applied to the setting of perfectoid spaces.

The section is structured as follows. We start in \cref{sec:andesc.cond} by introducing the $\infty$-category $\Cond(\mathcal C)$ of condensed objects in a given $\infty$-category $\mathcal C$ and discuss some set-theoretic subtleties arising from this definition. In \cref{sec:andesc.almmath} we introduce the $\infty$-category $\D(V,\mm)$ of condensed derived almost $V$-modules associated to an almost setup $(V,\mm)$ and generalize many of the constructions from \cite{almost-ring-theory} to this setting; in particular we discuss animated almost rings and provide tools to circumvent the fact that $\D(V,\mm)$ is not compactly generated (unlike $\D(V)$). In \cref{sec:andesc.anring} we introduce analytic rings in the almost setting, thereby generalizing most of the definitions and results in the second half of \cite{scholze-analytic-spaces}. In \cref{sec:andesc.anspace} we construct analytic spaces over an almost setup, generalizing the construction in \cite[\S13]{scholze-analytic-spaces}. This completes the general theory of condensed almost mathematics. Afterwards, in \cref{sec:andesc.endofun} we provide some technical constructions in order to properly speak about ``enriched endofunctors'' $\D(\mathcal A) \to \D(\mathcal A)$ associated to an analytic ring $\mathcal A$, which we use in \cref{sec:andesc.descmorph} to adapt Mathew's notion of \emph{descendability} \cite{akhil-galois-group-of-stable-homotopy} to the setting of analytic almost geometry. In \cref{sec:andesc.filtcolim} we investigate the behavior of descendability under filtered colimits and introduce the notion of \emph{weakly descendable} maps. In \cref{sec:andesc.fsdesc} we show that descent can be checked after passing to certain covers of the base (e.g. connected components of a scheme) by carefully analyzing convergence speeds of totalizations. Combining this with our results on filtered colimits we are lead to the notion of \emph{weakly fs-descendable} morphisms of analytic rings, which provide the most general form of descent discusses in this thesis and are well-suited to be applied to the setting of perfectoid spaces. This finishes the discussion of the general descent theory, so we continue with examples: In \cref{sec:andesc.scheme} we introduce solid analytic rings and glue them to schemes and discrete adic spaces in the setting of analytic almost geometry; among others, we construct a full 6-functor formalism for solid quasicoherent sheaves on discrete adic spaces based on \cite[\S11]{condensed-mathematics}. In \cref{sec:andesc.schemedesc} we apply our general descent formalism to the setting of schemes and discrete adic spaces in order to prove various descent results, in particular ind-fppf descent for solid modules. In \cref{sec:andesc.valrings} we show that for a flat map of valuation rings $V \to V'$ the associated map of analytic rings $V_\solid \to V'_\solid$ has Tor dimension $\le 1$, which is a crucial step in the proof of v-descent for $\ri^{+a}_X/\pi$-modules. Finally, \cref{sec:andesc.adic-sigobj} introduces two somewhat unrelated concepts: adically complete rings and modules, and modules equipped with a semilinear endomorphism; both of these concepts are used in the context of our $p$-torsion Riemann-Hilbert correspondence for $\ri^{+a}_X/\pi$-modules.

\subsection{Condensed Objects} \label{sec:andesc.cond}

The basis of analytic geometry is condensed mathematics, as introduced in \cite{condensed-mathematics}. It acts as a replacement of topology and removes many of the defects of topological spaces in the context of algebraic structures. In the following we introduce the main definitions and prove some basic properties; in particular we deal with set-theoretic issues which arise from working with sheaves on a big site. Our approach for addressing set-theoretic issues is that from \cite[Appendix to Lecture II]{condensed-mathematics} -- while it poses some challenges, it allows us to avoid the choice of any cutoff cardinals and does not need to assume the existence of strongly inaccessible cardinals. A large part of this subsection is devoted to showing that internal $\Hom$-objects exist in the category of condensed sets.

Let us start with the main definition: condensed objects in an $\infty$-category $\mathcal C$. Letting $\mathcal C$ be the category of sets produces the category of condesed sets and hence a replacement of topological spaces. Since we are interested in an algebraic theory, we will usually let $\mathcal C$ be an $\infty$-category of algebraic structures, most prominently modules over some ring or the derived category thereof.

\begin{definition}
Let $\mathcal C$ be an $\infty$-category which has all small colimits.
\begin{defenum}
	\item \label{def:Cond-C-kappa} For every strong limit cardinal $\kappa$ we let $\Cond(\mathcal C)_\kappa$ be the $\infty$-category of contravariant functors from $\kappa$-small extremally disconnected sets to $\mathcal C$ which take finite coproducts to products. We call $\Cond(\mathcal C)$ the $\infty$-category of \emph{$\kappa$-condensed objects in $\mathcal C$}.

	\item Via the fully faithful left-adjoints to the forgetful functors (given by Kan extension) we define
	\begin{align*}
		\Cond(\mathcal C) := \varinjlim_\kappa \Cond(\mathcal C)_\kappa,
	\end{align*}
	where the colimit on the right ranges over all strong limit cardinals $\kappa$. We call $\Cond(\mathcal C)$ the $\infty$-category of \emph{condensed objects in $\mathcal C$}.

	\item For every strong limit cardinal $\kappa$ we denote
	\begin{align*}
		(-)_\kappa\colon \Cond(\mathcal C) \to \Cond(\mathcal C)_\kappa, \qquad X \mapsto X_\kappa
	\end{align*}
	the restriction functor to $\kappa$-small extremally disconnected sets. We say that an object $X \in \Cond(\mathcal C)$ is \emph{discrete} if $X \in \Cond(\mathcal C)_\omega$.
\end{defenum}
\end{definition}

\begin{remark} \label{rmk:condensed-objects-are-sheaves}
In \cref{def:Cond-C-kappa}, if $\mathcal C$ has all small limits, then $\Cond(\mathcal C)_\kappa$ is the $\infty$-category of hypercomplete $\mathcal C$-valued sheaves on the site of $\kappa$-small profinite sets (where covers are finite jointly surjective families of maps). This follows from the fact that $\kappa$-small extremally disconnected sets form a basis of the site of $\kappa$-small profinite sets and that they are weakly contractible (see \cref{rslt:sheaves-on-basis-equiv-sheaves-on-whole-site} and \cite[\S2]{condensed-mathematics}).
\end{remark}

The definition of $\Cond(\mathcal C)_\kappa$ works well for strong limit cardinals $\kappa$ because of the existence of enough $\kappa$-small extremally disconnected sets in that case (as just explained in \cref{rmk:condensed-objects-are-sheaves}). However, there does not seem to be a good criterion for detecting whether a condensed object $X \in \Cond(\mathcal C)$ lies in $\Cond(\mathcal C)_\kappa \subset \Cond(\mathcal C)$. In particular this makes it a priori hard to verify that a given contravariant functor from extremally disconnected sets to $\mathcal C$ lies in $\Cond(\mathcal C)$, i.e. is ``small''. Such criteria do exist if we work with regular cardinals $\kappa$ instead of strong limit cardinals (see \cref{rslt:essential-image-of-pullback-of-regular-condensed-objects} below), which we will show next. By passing to larger strong limit cardinals, we can then deduce criteria for containment in $\Cond(\mathcal C)$.

Let us now study $\Cond(\mathcal C)_\kappa$ for regular cardinals $\kappa$, with the above mentioned application in mind. The $\infty$-categories $\Cond(\mathcal C)_\kappa$ for regular $\kappa$ will not be used anymore after we have proved the desired application to $\Cond(\mathcal C)$, so the reader may skip the details of the following definitions and results and just take \cref{rslt:characterize-condensed-objects-via-limit-property} for granted. We start with the definition:

\begin{definition}
Let $\kappa$ be a regular cardinal.
\begin{defenum}
	\item A profinite set $S$ is called \emph{$\kappa$-cocompact} if for all $\kappa$-cofiltered systems $(T_i)_i$ of profinite sets we have
	\begin{align*}
		\Hom(\varprojlim_i T_i, S) = \varinjlim_i \Hom(T_i, S).
	\end{align*}

	\item Let $\mathcal C$ be an $\infty$-category which has all small limits and colimits. Then we denote by $\Cond(\mathcal C)_\kappa$ the $\infty$-category of hypercomplete $\mathcal C$-valued sheaves on the site of $\kappa$-cocompact profinite sets (where covers are given by finite families of jointly surjective maps).
\end{defenum}
\end{definition}

\begin{remark}
Since the category of profinite sets is the $\Pro$-category of the category of finite sets, which in turn contains all finite limits, it follows that this category is ``cocompactly generated'', i.e. satisfies all the dual properties of compactly generated categories. For example, for every regular cardinal $\kappa$, the subcategory of $\kappa$-cocompact profinite sets is stable under $\kappa$-small limits and every profinite set is a $\kappa$-cofiltered limit of $\kappa$-cocompact profinite sets.
\end{remark}

In order to study $\Cond(\mathcal C)_\kappa$ for regular cardinals $\kappa$ we need a good understanding of $\kappa$-cocompact profinite sets. The following results show that every hypercover of profinite sets can be written as a $\kappa$-cofiltered limit of hypercovers of $\kappa$-cocompact sets (see \cref{rslt:hypercover-of-profinite-sets-is-lim-of-cocompact-hypercovers} below).

\begin{lemma}
Let $\kappa$ be a regular cardinal and $S$ a profinite set.
\begin{lemenum}
	\item \label{rslt:cocompact-profin-set-is-small-cofiltered-limit} $S$ is $\kappa$-cocompact if and only if it is a $\kappa$-small cofiltered limit of finite sets.

	\item \label{rslt:closed-subset-of-kappa-cocompact-profin-is-kappa-cocompact} Suppose that $S$ is $\kappa$-cocompact and let $Z \subset S$ a closed subset. Then $Z$ is $\kappa$-cocompact.
\end{lemenum}
\end{lemma}
\begin{proof}
In general, if $\mathcal C$ is any compactly generated ($\infty$-)category and $\kappa$ a regular cardinal, then the $\kappa$-compact object in $\mathcal C$ are precisely the retracts of $\kappa$-small filtered colimits of compact objects in $\mathcal C$. Namely, let $\mathcal C' \subset \mathcal C$ be the full subcategory spanned by these objects. Since $\kappa \gg \omega$, the proof of \cite[Proposition 5.4.2.11]{lurie-higher-topos-theory} implies that $\mathcal C'$ generates $\mathcal C$ under $\kappa$-filtered colimits, hence the proof of \cite[Proposition 5.4.2.2]{lurie-higher-topos-theory} implies $\mathcal C = \Ind_\kappa(\mathcal C')$, so that \cite[Lemma 5.4.2.4]{lurie-higher-topos-theory} implies $\mathcal C' = \mathcal C^\kappa$. Applying the dual statement to the category of profinite sets (which is cocompactly generated by the finite sets) we deduce that a profinite set $S$ is $\kappa$-cocompact if and only if it is a retract of a $\kappa$-small cofiltered limit of finite sets.

Now (i) and (ii) reduce to the following claim: Suppose $S = \varprojlim_i S_i$ is a $\kappa$-small cofiltered limit of finite sets and $Z \subset S$ is a closed subsect; then $Z$ is a $\kappa$-small cofiltered colimit of finite sets. This can be shown as follows: For each $i$ let $Z_i \subset S_i$ denote the image of $Z$ under the projection $S \to S_i$; it is enough to see that $Z = \varprojlim_i Z_i$. ``$\subset$'' is clear. For ``$\supset$'' note that for every $x \in S \setminus Z$ we can find a decomposition $S = S' \isect S''$ into clopen subsets such that $Z \subset S'$ and $x \in S''$. This decomposition must come via preimage from some $S_{i_0}$. In particular, $Z$ and $x$ get mapped to disjoint subsets of $S_{i_0}$ and hence $x \not\in \varprojlim_i Z_i$.
\end{proof}

\begin{lemma} \label{rslt:surjection-to-kappa-cofiltered-lim-is-kappa-filtered-lim-of-surjections}
Let $\kappa$ be a regular cardinal, let $(T_i)_{i\in I}$ be a $\kappa$-cofiltered system of $\kappa$-cocompact profinite sets with limit $T := \varprojlim_i T_i$ and let $S \surjto T$ be any surjective map from a profinite set $S$. Then there is a $\kappa$-cofiltered category $I'$ with a functor $\alpha\colon I' \to I$ and a diagram $(S_{i'} \surjto T_{\alpha(i')})_{i'\in I'}$ of surjections of $\kappa$-cocompact profinite sets such that:
\begin{lemenum}
	\item $\alpha$ is $\kappa$-cofinal, i.e. for every $\kappa$-small diagram $f\colon K \to I'$ and any $i \in I_{/\alpha f}$ there is some $i' \in I'_{/f}$ and a morphism $\alpha(i') \to i$ in $I_{/\alpha f}$.

	\item If $I$ admits $\kappa$-small cofiltered limits, then $I'$ can be chosen to have the same property and $\alpha$ can be chosen to preserve these limits.

	\item We have $(S \to T) = \varprojlim_{i\in I'} (S_{i'} \to T_{\alpha(i')})$.
\end{lemenum}
\end{lemma}
\begin{proof}
By the dual version of \cite[Proposition 5.3.5.15]{lurie-higher-topos-theory} there is a $\kappa$-filtered index set $J$ and a diagram $(S_j \to T_j)_{j\in J}$ of $\kappa$-cocompact profinite sets such that $(S \to T) = \varprojlim_{j\in J} (S_j \to T_j)$. In fact we can take $J$ to be the category of commutative diagrams
\begin{center}\begin{tikzcd}
	S \arrow[r] \arrow[d] & T \arrow[d]\\
	S_j \arrow[r] & T_j
\end{tikzcd}\end{center}
where $S_j$ and $T_j$ are $\kappa$-cocompact. Then $J$ admits a surjective projection to the category $J_0$ of all maps $T \to T_j$ for $T_j$ $\kappa$-cocompact. One checks easily that the map $I \to J_0, i \mapsto (T \to T_i)$ is cofinal, so we can form $J' := J \cprod_{J_0} I$, which is still $\kappa$-cofiltered and admits a surjective map $J' \to I$ and a cofinal map $J' \to J$. Let $I' \subset J'$ be the full subcategory of those $j$ where the map $S_j \to T_j$ is surjective. We claim that $I'$ is cofinal in $J'$. Given any diagram as above in $J'$, let $T_j'$ be the image of the map $S_j \to T_j$. By \cref{rslt:closed-subset-of-kappa-cocompact-profin-is-kappa-cocompact} $T_j'$ is still $\kappa$-cocompact and thus defines an element of $J$. Pick any $i \in I$ which admits a map $T_i \to T_j'$ under $T$, and let $S_i := S_j \cprod_{T_j'} T_i$. This defines an element of $I'$ with a map to the original diagram in $J'$. As $I'$ is cofinal in $J'$, it is clear that (iii) is satisfied and one easily checks that (i) also holds. Part (ii) follows easily from the construction, as $\kappa$-small limits of $\kappa$-cocompact profinite sets are $\kappa$-cocompact and cofiltered limits of surjections of compact Hausdorff spaces are surjections.
\end{proof}

\begin{lemma} \label{rslt:hypercover-of-profinite-sets-is-lim-of-cocompact-hypercovers}
Let $S_\bullet \to S$ be a hypercover of profinite sets and let $\kappa$ be an uncountable regular cardinal. Then there is a $\kappa$-cofiltered family $(T_{i,\bullet} \to T_i)_{i\in I}$ of hypercovers of $\kappa$-cocompact profinite sets such that
\begin{align*}
	(S_\bullet \to S) = \varprojlim_i (T_{i,\bullet} \to T_i).
\end{align*}
\end{lemma}
\begin{proof}
Using \cref{rslt:surjection-to-kappa-cofiltered-lim-is-kappa-filtered-lim-of-surjections} we can inductively construct a chain $\dots \to I_n \to \dots \to I_0$ of $\kappa$-cofiltered categories which admit $\kappa$-small cofiltered limits and whose transition maps $\alpha_n\colon I_n \to I_{n-1}$ preserve these limits, and for every $n$ a diagram $(T_{i,\bullet\le n} \to T_i)_{i\in I_n}$ of $n$-truncated hypercovers of $\kappa$-cocompact profinite sets such that $\varprojlim_{i\in I_n} (T_{i,\bullet\le n} \to T_i) = (S_{\bullet\le n} \to S)$ and such that $\alpha_n$ identifies the $(n-1)$-truncated parts. Then $I := \varprojlim_n I_n$ is still $\kappa$-cofiltered (this can be seen directly in this case; see \cref{rslt:catfiltcofin-stable-under-tau-small-limits} for the general statement), which produces the desired diagram.
\end{proof}

We now arrive at the promised criterion for checking containment in $\Cond(\mathcal C)_\kappa$ for regular cardinals $\kappa$. The following result is a generalization of the observation that a condensed set $X$ is discrete if and only if for every profinite set $S = \varprojlim_i S_i$ we have $X(S) = \varinjlim_i X(S_i)$.

\begin{proposition}
Let $\mathcal C$ be an $\infty$-category which admits all small limits and colimits. Let $\kappa' \ge \kappa$ be uncountable regular cardinals and assume that $\kappa$-filtered colimits commute with totalizations in $\mathcal C$. Then:
\begin{propenum}
	\item The pullback functor $\Cond(\mathcal C)_\kappa \to \Cond(\mathcal C)_{\kappa'}$ is fully faithful.

	\item \label{rslt:essential-image-of-pullback-of-regular-condensed-objects} An object $M \in \Cond(\mathcal C)_{\kappa'}$ lies in the essential image of the pullback functor if and only if for every $\kappa$-cofiltered limit $S = \varprojlim_i S_i$ of profinite sets such that $S$ and all $S_i$ are $\kappa'$-cocompact we have
	\begin{align*}
		M(S) = \varinjlim_i M(S_i).
	\end{align*}
\end{propenum}
\end{proposition}
\begin{proof}
Let $N \in \Cond(\mathcal C)_\kappa$ with pullback $M \in \Cond(\mathcal C)_{\kappa'}$. Let $M'$ be the left Kan extension of $N$ to $\kappa'$-cocompact profinite sets; then $M$ is the sheafification of $M'$. We claim that $M'$ is already a hypercomplete sheaf, i.e. $M = M'$. By construction of Kan extensions we have $M'(S) = \varinjlim_{S \to T} N(T)$ for every profinite set $S$, where the colimit is taken over the category of maps $S \to T$ with $T$ being $\kappa$-cocompact. Given any hypercover $S_\bullet \to S$ of $\kappa'$-cocompact profinite sets, by \cref{rslt:hypercover-of-profinite-sets-is-lim-of-cocompact-hypercovers} we can write $(S_\bullet \to S) = \varprojlim_i (T_{i,\bullet} \to T_i)$ for a $\kappa$-cofiltered diagram $(T_{i,\bullet} \to T_i)_i$ of hypercovers of $\kappa$-cocompact profinite sets. For each $n \ge 0$ the diagram $(T_{i,n})_i$ is cofinal among the diagram of maps $S_i \to T$ with $T$ $\kappa$-cocompact, hence
\begin{align*}
	M'(S) = \varinjlim_i N(T_i) = \varinjlim_i \varprojlim_{n\in\Delta} N(T_{i,n}) = \varprojlim_{n\in\Delta} \varinjlim_i N(T_{i,n}) = \varprojlim_{n\in\Delta} M'(S_n),
\end{align*}
proving that $M'$ is indeed a hypercomplete sheaf. Thus $M = M'$, which implies (i).

To prove (ii), let $S = \varprojlim_i S_i$ be given as in the claim. First assume that $M$ lies in the essential image of the pullback, i.e. $M = M'$ as above. Since colimits commute with colimits, all we have to show is that the category of maps $S \to T$ (with $T$ being $\kappa$-cocompact) is the filtered colimit of the categories of maps $S_i \to T$. This follows immediately from the definition of $\kappa$-cocompact profinite sets. Conversely, if $M \in \Cond(\mathcal C)_{\kappa'}$ satisfies the described limit property then in particular $M(S) = \varinjlim_{S \to T} M(T)$ (with $T$ ranging over $\kappa$-cocompact profinite sets), as this limit of $\kappa$-cofiltered; thus $M \in \Cond(\mathcal C)_\kappa$, as desired.
\end{proof}

We will now apply the above study of $\Cond(\mathcal C)_\kappa$ for regular $\kappa$ to get the promised containment criterion for when a contravariant functor from extremally disconnected sets to $\mathcal C$ lies in $\Cond(\mathcal C)$. The following result is the only application of the categories $\Cond(\mathcal C)_\kappa$ with regular $\kappa$; they will not be needed anymore afterwards.

\begin{proposition} \label{rslt:characterize-condensed-objects-via-limit-property}
Let $\mathcal C$ be an $\infty$-category which has all small limits and colimits and assume that there is a regular cardinal $\kappa$ such that $\kappa$-filtered colimits commute with totalizations in $\mathcal C$. Let $M$ be a contravariant functor from extremally disconnected sets to $\mathcal C$ which maps finite coproducts to products. Then $M \in \Cond(\mathcal C)$ if and only if there is a regular cardinal $\kappa$ such that for every $\kappa$-cofiltered limit $S = \varprojlim_i S_i$ of extremally disconnected sets we have
\begin{align*}
	M(S) = \varinjlim_i M(S_i).
\end{align*}
In this case, the above identity is also true for every $\kappa$-cofiltered limit $S = \varprojlim_i S_i$ of profinite sets.
\end{proposition}
\begin{proof}
The ``if'' part is easy: Assuming the existence of $\kappa$ such that $M$ has the given limit property, pick any strong limit cardinal $\kappa'$ whose cofinality is at least $\kappa$. Then for every extremally disconnected set $S$, the category of maps $S \to T$ with $T$ a $\kappa'$-small extremally disconnected set is $\kappa$-cofiltered, hence $M(S) = \varinjlim_{S \to T} M(T)$, showing that $M \in \Cond(\mathcal C)_{\kappa'}$.

We now prove the ``only if'' part together with the final remark. Assume $M \in \Cond(\mathcal C)_{\kappa'}$ for some strong limit cardinal $\kappa'$. Pick any regular cardinal $\kappa \ge \kappa'$ such that $\kappa$-filtered colimits commute with totalizations in $\mathcal C$. Let $S = \varprojlim_i S_i$ be any $\kappa$-cofiltered limit of profinite sets. Choose a regular cardinal $\kappa'' \ge \kappa$ such that $S$ and all $S_i$ are $\kappa''$-cocompact and admit hypercovers by $\kappa''$-cocompact extremally disconnected sets. Pick any strong limit cardinal $\kappa''' \ge \kappa''$ such that all $\kappa''$-cocompact profinite sets are $\kappa'''$-small. We have the following pullback functors:
\begin{align*}
	\Cond(\mathcal C)_{\kappa'} \xto{\alpha^*} \Cond(\mathcal C)_\kappa \xto{\beta^*} \Cond(\mathcal C)_{\kappa''} \xto{\gamma^*} \Cond(\mathcal C)_{\kappa'''}.
\end{align*}
Note that for any $N \in \Cond(\mathcal C)_{\kappa''}$ and any $\kappa''$-cocompact extremally disconnected set $T$ we have $(\gamma^*N)(T) = N(T)$. Namely, $\gamma^* N$ is the sheafification of the left Kan extension of $N$ to all $\kappa'''$-small profinite sets; but the left Kan extension does not change the value of $N$ on any $\kappa''$-cocompact profinite set and the sheafification does not change the value of $N$ on any $\kappa'''$-small extremally disconnected sets as these form a basis of the site of $\kappa'''$-small profinite sets (by \cref{rslt:sheaves-on-basis-equiv-sheaves-on-whole-site} sheafification on a site can be computed on any basis of the site followed by right Kan extension). By hypercompleteness and the choice of $\kappa''$ we deduce $(\gamma^*N)(S) = N(S)$ and $(\gamma^*N)(S_i) = N(S_i)$ for all $i$. Now the claim follows by applying \cref{rslt:essential-image-of-pullback-of-regular-condensed-objects} to $N = \beta^*(\alpha^* M)$.
\end{proof}

\begin{corollary} \label{rslt:condensed-objects-stable-under-pushforward}
Let $\mathcal C$ be an $\infty$-category which has all small limits and colimits and assume that there is a regular cardinal $\kappa$ such that $\kappa$-filtered colimits commute with totalizations in $\mathcal C$. Then for every $M \in \Cond(\mathcal C)$ and every profinite set $S$ (with projection $j_S\colon S \to *$) the functor
\begin{align*}
	j_{S*} M\colon T \mapsto M(S \cprod T)
\end{align*}
defines an object of $\Cond(\mathcal C)$.
\end{corollary}
\begin{proof}
Clearly $j_{S*} M$ is a hypercomplete sheaf on all profinite sets, i.e. satisfies descent along all hypercovers. By \cref{rslt:characterize-condensed-objects-via-limit-property} there is some regular cardinal $\kappa$ such that for all $\kappa$-cofiltered limits $T = \varprojlim_i T_i$ of profinite sets we have $M(T) = \varinjlim_i M(T_i)$. Since $T \cprod S = \varprojlim_i (T_i \cprod S)$, $j_{S*} M$ has the same property. But then \cref{rslt:characterize-condensed-objects-via-limit-property} implies that $j_{S*}M \in \Cond(\mathcal C)$, as desired.
\end{proof}

We now arrive at the main result of the present subsection. It states that if $\mathcal C$ is a ``nice'' $\infty$-category then $\Cond(\mathcal C)$ inherits the nice structure of $\mathcal C$.

\begin{proposition} \label{rslt:condensed-objects-in-presentable-monoidal-cat}
Let $\tau$ be a regular cardinal and let $\mathcal C$ be a $\tau$-compactly generated $\infty$-category. Then:
\begin{propenum}
	\item $\Cond(\mathcal C)$ has all small limits and colimits. For every strong limit cardinal $\kappa$, the full subcategory $\Cond(\mathcal C)_\kappa \subset \Cond(\mathcal C)$ is presentable and stable under all colimits. It is $\tau$-compactly generated by $\tau$-compact objects $N[S]$ for $\tau$-compact $N \in \mathcal C$ and $\kappa$-small extremally disconnected sets $S$ such that
	\begin{align*}
		\Hom(N[S], M) = \Hom(N, M(S))
	\end{align*}
	for all $M \in \Cond(\mathcal C)$. Moreover, for any given cardinal $\lambda$, there is a cofinal class of strong limit cardinals $\kappa$ such that $\Cond(\mathcal C)_\kappa$ is stable under $\lambda$-small limits in $\Cond(\mathcal C)$.

	\item A closed symmetric monoidal structure on $\mathcal C$ naturally induces a closed symmetric monoidal structure on $\Cond(\mathcal C)$. For every strong limit cardinal $\kappa$, the full subcategory $\Cond(\mathcal C)_\kappa \subset \Cond(\mathcal C)$ is stable under the monoidal structure.
\end{propenum}
\end{proposition}
\begin{proof}
Fix a regular cardinal $\tau$ such that $\mathcal C$ is $\tau$-compactly generated. For every strong limit cardinal $\kappa$ let $\Cond(\mathcal C)'_\kappa$ be the $\infty$-category of all contravariant functors from $\kappa$-small extremally disconnected sets to $\mathcal C$. By \cite[Proposition 5.5.3.6]{lurie-higher-topos-theory} $\Cond(\mathcal C)'_\kappa$ is presentable, hence for every $\kappa$-small extremally disconnected set $S$ the functor $\Cond(\mathcal C)'_\kappa \to \mathcal C$, $M \mapsto M(S)$ admits a left adjoint $N \mapsto N[S]'$. Then $\Cond(\mathcal C)_\kappa \subset \Cond(\mathcal C)'_\kappa$ is the localization along the set $W$ of morphisms $N[S_1]' \dunion N[S_2]' \to N[S_1 \dunion S_2]'$ for $\kappa$-small extremally disconnected sets $S_1, S_2$ and $\tau$-compact $N \in \mathcal C$. In particular by \cite[Proposition 5.5.4.15]{lurie-higher-topos-theory} $\Cond(\mathcal C)_\kappa$ is presentable. Let us denote $N[S]$ the sheafification of $N[S]'$, i.e. the image of $N[S]'$ under the localization $\Cond(\mathcal C)'_\kappa \to \Cond(\mathcal C)_\kappa$. By \cref{rslt:compact-generators-implies-compactly-generated} it follows that $\Cond(\mathcal C)_\kappa$ is $\tau$-compactly generated and that every object in $\Cond(\mathcal C)_\kappa$ is an iterated colimit of the objects $N[S]$ for varying $\kappa$-extremally disconnected sets $S$ and $\tau$-compact $N \in \mathcal C$ (and the same is true for $\Cond(\mathcal C)'_\kappa$ and the objects $N[S]'$).

If $\kappa' \ge \kappa$ is another strong limit cardinal then the pullback functor $\Cond(\mathcal C)_\kappa \to \Cond(\mathcal C)_{\kappa'}$ is given by left Kan extension. It is clear that this pullback preserves small colimits, which implies that $\Cond(\mathcal C)$ admits all small colimits and that $\Cond(\mathcal C)_\kappa$ is stable under small colimits.

To prove the claims about limits, let the cardinality $\lambda$ be given. Pick a strong limit cardinal $\kappa$ large enough so that the cofinality $\lambda'$ of $\kappa$ satisfies $\lambda' \gg \tau$ and $\lambda' \ge \lambda$ (there is a cofinal class of such $\kappa$). Then $\lambda'$-filtered colimits commute with $\lambda'$-small limits in $\mathcal C$ by \cref{rslt:filtered-colim-preserve-small-lim-in-compactly-generated-cat}, so by the same argument as in \cite[Proposition 2.9]{condensed-mathematics}, for every strong limit cardinal $\kappa' \ge \kappa$ the pullback $\Cond(\mathcal C)_\kappa \to \Cond(\mathcal C)_{\kappa'}$ commutes with $\lambda'$-small limits. This finishes the proof of (i).

Before continuing with the proof of (ii), let us mention the following explicit calculation of $N[S]'$: For every $N \in \mathcal C$, and all $\kappa$-small extremally disconnected sets $S$ and $T$ we have a natural isomorphism
\begin{align*}
	N[S]'(T) = \bigdunion_{\pi_0\Hom(T, S)} N
\end{align*}
in $\mathcal C$. This follows from the fact that the functor $N[S]'$ is computed as a left Kan extension of $N \in \mathcal C = \Fun(\{ S \}, \mathcal C)$ along the inclusion of $S$ into the category of ($\kappa$-small) extremally disconnected sets.

We now prove (ii), so assume that $\mathcal C$ is equipped with a closed symmetric monoidal structure $\tensor$. For any strong limit cardinal $\kappa$ we get an induced symmetric monoidal structure on $\Cond(\mathcal C)'_\kappa$. To show that this induces a symmetric monoidal structure on the localization $\Cond(\mathcal C)_\kappa$, we employ \cite[Proposition 4.1.7.4]{lurie-higher-algebra}: Recall the definition of the set of morphisms $W$ above; we need to show that for every $f \in W$ and every $M \in \Cond(\mathcal C)'_\kappa$ the map $f \tensor M$ becomes an isomorphism after localization. By transfinite induction we can assume that $M = N[S]'$ for some ($\tau$-compact) $N \in \mathcal C$ and some $\kappa$-small profinite set $S$. By the above explicit computation of $N[S]'$ we deduce
\begin{align*}
	N_1[S_1]' \tensor N_2[S_2]' = (N_1 \tensor N_2)[S_1 \cprod S_2]'
\end{align*}
for all $N_1, N_2 \in \mathcal C$ and all $\kappa$-small extremally disconnected sets $S_1$ and $S_2$. From this we deduce immediately that $f \tensor N[S]' \in W$. This finishes the construction of the symmetric monoidal structure on $\Cond(\mathcal C)_\kappa$.

A similar computation as in the previous paragraph shows that the pullback functors $\Cond(\mathcal C)_\kappa \to \Cond(\mathcal C)_{\kappa'}$ are symmetric monoidal, hence we also obtain a symmetric monoidal structure on $\Cond(\mathcal C)$. It remains to show that the symmetric monoidal structure on $\Cond(\mathcal C)$ is closed. On each $\Cond(\mathcal C)_\kappa$ this is true by the adjoint functor theorem, so we need to see that for any given $M, M' \in \Cond(\mathcal C)$ there is some strong limit cardinal $\kappa$ such that $M, M' \in \Cond(\mathcal C)_\kappa$ and for all strong limit cardinals $\kappa' \ge \kappa$ we have
\begin{align*}
	\IHom_{\Cond(\mathcal C)_\kappa}(M, M') = \IHom_{\Cond(\mathcal C)_{\kappa'}}(M, M').
\end{align*}
To see this, write $M$ as an iterated small colimit of objects $N[S]$ for varying $N \in \mathcal C$ and $\kappa$-small extremally disconnected sets $S$. We can pull this colimit out of the $\IHom$'s (after which it turns into an iterated small limit), so by the claims about limits in (i) we can reduce to the case $M = N[S]$. In this case choose a regular cardinal $\kappa_0 \ge \tau$ such that $M'(T)$ transforms $\kappa_0$-cofiltered limits in $T$ into $\kappa_0$-filtered colimits (see \cref{rslt:characterize-condensed-objects-via-limit-property}). We claim that the same is then true for $\IHom(N[S], M')(T)$ (computed in $\Cond(\mathcal C)_\kappa$ for large $\kappa$). Namely, for any $\tau$-compact $N' \in \mathcal C$ and any $\kappa_0$-cofiltered limit $T = \varprojlim_i T_i$,
\begin{align*}
	\Hom(N', \IHom(N[S], M')(\varprojlim_i T_i)) &= \Hom(N'[\varprojlim_i T_i], \IHom(N[S], M'))\\
	&= \Hom(N[S], \IHom(N'[\varprojlim_i T_i], M'))\\
	&= \Hom(N[S], \varinjlim_i \IHom(N'[T_i], M'))\\
	&= \varinjlim_i \Hom(N', \IHom(N[S], M')(T_i))\\
	&= \Hom(N', \varinjlim_i \IHom(N[S], M')(T_i)).
\end{align*}
This proves (ii).
\end{proof}

\begin{remark}
Most of \cref{rslt:condensed-objects-stable-under-pushforward} is rather formal, but the existence of $\IHom$ in $\Cond(\mathcal C)$ is surprisingly subtle: Its proof relies on the characterization of condensed objects in \cref{rslt:characterize-condensed-objects-via-limit-property}, which itself required a good understanding of $\kappa$-cocompact profinite sets. There may be more direct approaches to proving the existence of $\IHom$, but we found the presented one the most conceptional. Note that even in the case that $\mathcal C$ is the category of sets or the category of abelian groups, the existence of $\IHom$ in $\Cond(\mathcal C)$ has not been justified in the literature before.
\end{remark}

We close this section with the following result, showing that forming condensed objects and forming derived categories commute with each other.

\begin{proposition} \label{rslt:Cond-commutes-with-D}
Let $\mathcal A$ be a Grothendieck abelian category and $\kappa$ a strong limit cardinal. Then $\Cond(\mathcal A)_\kappa$ is a Grothendieck abelian category and there is a natural equivalence of $\infty$-categories
\begin{align*}
	\D(\Cond(\mathcal A)_\kappa) = \Cond(\D(\mathcal A))_\kappa.
\end{align*}
\end{proposition}
\begin{proof}
Note first that $\Cond(\mathcal A)_\kappa$ is an accessible localization of the category $\mathcal C$ of $\mathcal A$-valued presheaves on $\kappa$-small profinite sets and in particular a presentable abelian category. As colimits in $\Cond(\mathcal A)_\kappa$ are computed pointwise on extremally disconnected sets, filtered colimits in $\Cond(\mathcal A)_\kappa$ are exact, so that this category is indeed Grothendieck abelian.

There is a natural functor $F\colon \D(\Cond(\mathcal A)_\kappa) \to \Cond(\D(\mathcal A))_\kappa$ with $F(M)(S) = \Gamma(S, M)$ for all $M \in \D(\Cond(\mathcal A)_\kappa)$ and extremally disconnected sets $S$. As this functor commutes with cohomology, it is easily seen to be conservative. Also $F$ commutes with all colimits and hence admits a right adjoint $G$. For every extremally disconnected set $S$ and every $M \in \mathcal A$ denote $M[S] \in \Cond(\mathcal A)_\kappa$ the sheafification of the presheaf $T \mapsto M[\Hom(T, S)]$. We can view $M[S]$ as an object of $\D(\Cond(\mathcal A)_\kappa)$ and as an object of $\Cond(\D(\mathcal A))_\kappa$. In both cases it satisfies $\Hom(M[S], N) = \Hom(M, N(S))$ for all other objects $N$ (for $\D(\Cond(\mathcal A)_\kappa)$ use K-injective resolutions to reduce to the abelian level; for this it was shown in the proof of \cref{rslt:condensed-objects-in-presentable-monoidal-cat}). Thus, for every $N \in \Cond(\D(\mathcal A))_\kappa$ we have
\begin{align*}
	\Hom(M, F(G(N))(S)) = \Hom(M[S], G(N)) = \Hom(M, N(S)).
\end{align*}
By Yoneda (using that $\D(\mathcal A)$ is generated under colimits by $\mathcal A$) this implies $F(G(N)) = N$ for all $N$, i.e. $G$ is fully faithful. Together with the fact that $F$ is conservative we deduce that $F$ and $G$ are equivalences.
\end{proof}

\subsection{Almost Mathematics} \label{sec:andesc.almmath}

We now introduce a general framework to do almost mathematics by adapting and generalizing many of the basic definitions in \cite{almost-ring-theory}. To every almost setup $(V,\mm)$ (to be defined shortly) consisting of a classical ring $V$ and an ideal $\mm$ we associate the $\infty$-category $\D(V,\mm)$ of condensed derived almost $V$-modules. Intuitively, the objects of $\D(V,\mm)$ can be seen as complexes of topological $V$-modules, where we force every map $f\colon M \to M'$ to be an isomorphism whenever for all $n \in \Z$ both kernel and cokernel of $\pi_n(f)\colon \pi_n(M) \to \pi_n(M')$ are annihilated by $\mm$. We will then also consider almost rings over $(V,\mm)$ which should be seen as ``derived topological $V$-modules'' with a multiplication map that is associative and commutative up to $\mm$-torsion.

Let us start by providing a precise definition of almost setups $(V,\mm)$. The following definition is taken from \cite[\S2.1.1]{almost-ring-theory}.

\begin{definition} \label{def:almost-setups}
\begin{defenum}
	\item An \emph{almost setup} is a pair $(V, \mm)$, where $V$ is a classical ring and $\mm \subset V$ is an ideal such that $\mm^2 = \mm$ and such that
	\begin{align*}
		\widetilde\mm := \pi_0(\mm \tensor_V \mm)
	\end{align*}
	is flat over $V$.

	\item A \emph{morphism of almost setups} $\varphi\colon (V, \mm) \to (V', \mm')$ is a morphism $\varphi\colon V \to V'$ of rings such that $\mm' \subset \varphi(\mm) V'$. We say that $\varphi$ is \emph{strict} if $\mm' = \varphi(\mm) V'$ and we say that $\varphi$ is \emph{localizing} if $V = V'$.

	\item We denote the category of almost setups by $\AlmSetup$.
\end{defenum}
\end{definition}

\begin{remarks}
\begin{remarksenum}
	\item Given an almost setup $(V, \mm)$ and any morphism $\varphi\colon V \to V'$ of discrete rings, there is a unique choice of ideal $\mm' \subset V'$ such that $(V', \mm')$ is an almost setup and $(V, \mm) \to (V', \mm')$ is a strict morphism. Indeed, we must take $\mm' = \varphi(\mm) V'$. Then it is clear that $\mm'^2 = \mm'$, and by \cite[Remark 2.1.4.(ii)]{almost-ring-theory} $\widetilde \mm'$ is flat over $V'$. It follows that every morphism of almost setups factors as a composition of a strict morphism and a localizing morphism.

	\item If $(V, \mm) \to (V', \mm')$ is a strict morphism of almost setups then $\widetilde\mm \tensor_V V' = \widetilde \mm'$ by \cite[Remark 2.1.4.(ii)]{almost-ring-theory}.

	\item In our applications it will often be the case that $\mm$ is flat over $V$, in which case $\mm = \widetilde\mm$. However, the flatness of $\mm$ is not preserved under base-change (in the sense of (i)), which is why we do not want to assume it in general, similar to \cite{almost-ring-theory}.
\end{remarksenum}
\end{remarks}

In our application to perfectoid spaces, the almost setups $(V,\mm)$ we work with will always satisfy a stronger property: The ideal $\mm \subset V$ is generated by countably many elements in $V$. While this additional property is irrelevant for most of our results about almost mathematics, it plays an important role when dealing with filtered colimits of descendable maps later on. We therefore introduce a special terminology for it:

\begin{definition}
Let $\kappa$ be a regular cardinal. We say that an almost setup $(V,\mm)$ is \emph{$\kappa$-compact} if $\widetilde\mm$ is $\kappa$-compact in the derived category of classical $V$-modules.
\end{definition}

\begin{lemma} \label{rslt:m-countably-generated-implies-omega-1-compact}
Let $(V,\mm)$ be an almost setup such that $\mm$ is countably generated as a $V$-module. Then $(V,\mm)$ is $\omega_1$-compact.
\end{lemma}
\begin{proof}
By \cite[Theorem 2.1.12.(ii)(a)]{almost-ring-theory} $\widetilde\mm$ is countably presented. Thus by the proof of \cite[Lemma 2.1.16]{almost-ring-theory} there is a short exact sequence $0 \to L \to L \to \widetilde\mm \to 0$ of $V$-modules where $L$ is a countable direct sum of copies of $V$. Then $L$ is $\omega_1$-compact and hence the same follows for $\widetilde\mm$.
\end{proof}

We now come to the definition of (condensed) almost modules over an almost setup. The following definition is very similar to \cite[\S2.2.2]{almost-ring-theory} except that we work with condensed $V$-modules throughout, thereby allowing ``topological'' structures on our modules.

\begin{definition}
Let $(V, \mm)$ be an almost setup. Denote $\Mod_V^{\mathrm{class}}$ the abelian category of classical $V$-modules and denote $\Cond(V) := \Cond(\Mod_V^{\mathrm{class}})$ the category of static $V$-modules.
\begin{defenum}
	\item A condensed $V$-module $M \in \Cond(V)$ is called \emph{almost zero} if for all $\varepsilon \in \mm$ the multiplication-by-$\varepsilon$ map $M \xto{\varepsilon} M$ is zero.

	\item A morphism $f\colon M \to N$ in $\Cond(V)$ is called an \emph{almost isomorphism} if $\ker f$ and $\coker f$ are almost zero.

	\item Let $\kappa$ be a strong limit cardinal. The full subcategory $\Sigma_\kappa \subset \Cond(V)_\kappa$ of almost zero $\kappa$-condensed $V$-modules is a Serre subcategory, so we can define the quotient
	\begin{align*}
		\Cond(V, \mm)_\kappa := \Cond(V)_\kappa/\Sigma_\kappa,
	\end{align*}
	the category of \emph{static almost $(V, \mm)$-modules}. Taking the colimit over all $\kappa$ we obtain
	\begin{align*}
		\Cond(V,\mm) = \varinjlim_\kappa \Cond(V,\mm)_\kappa = \Cond(V)/\Sigma,
	\end{align*}
	where $\Sigma \subset \Cond(V)$ is the full subcategory of almost zero $V$-modules. The localization functor $\Cond(V) \to \Cond(V, \mm)$ is denoted by $(V, \mm)^*$. If $(V, \mm)$ is clear from context, we also denote the localization functor by $M \mapsto M^a$.

	\item For every profinite set $S$ we let $V[S] \in \Cond(V)$ denote the free generator on $S$, i.e. the sheafification of $T \mapsto V[\Hom(T, S)]$. We denote
	\begin{align*}
		V^a[S] := V[S]^a \in \Cond(V,\mm).
	\end{align*}
\end{defenum}
\end{definition}

Note that in the case $\mm = V$ we get $\Cond(V, V) = \Cond(V)$, hence almost mathematics generalizes standard algebra. We will now deduce some basic properties about almost modules, generalizing many of the results of \cite[\S2.2]{almost-ring-theory} to the condensed world.

\begin{lemma} \label{rslt:compute-almost-Hom-over-V}
Let $(V, \mm)$ be an almost setup. Then for any $M, N \in \Cond(V)$ we have
\begin{align*}
	\Hom(M^a, N^a) = \Hom(\widetilde\mm \tensor_V M, N).
\end{align*}
\end{lemma}
\begin{proof}
We argue as in \cite[\S2.2.2]{almost-ring-theory}. Let $\mathcal C$ be the category of almost isomorphisms $P \to M$ in $\Cond(V)$. We claim that $\widetilde\mm \tensor_V M \to M$ is an initial object of $\mathcal C$. Indeed, given any $\phi\colon P \to M$ in $\mathcal C$, consider the diagram
\begin{center}\begin{tikzcd}
	\widetilde\mm \tensor_V P \arrow[r,"\sim"] \arrow[d] & \widetilde\mm \tensor_V M \arrow[d]\\
	P \arrow[r,"\phi"] & M
\end{tikzcd}\end{center}
One sees easily that the top horizontal map is an isomorphism (this reduces to the claim that if $Q \in \Cond(V)$ is almost zero, then $\widetilde\mm \tensor_V Q = 0$). This induces a map $\psi\colon \widetilde \mm \tensor_V M \to P$ over $M$. It remains to see that $\psi$ is unique. But if $\psi'\colon \widetilde\mm \tensor_V M \to P$ is another morphism over $M$ then $\img(\psi - \psi') \subset \ker(\phi)$ is almost zero and hence zero, as $\mm(\widetilde\mm \tensor_V M) = \widetilde\mm \tensor_V M$.

By construction of the quotient category, $\Hom(M^a, N^a)$ is the colimit over $\mathcal C$ of $\Hom(P, N)$, which implies the claim.
\end{proof}

\begin{lemma} \label{rslt:properties-of-almost-V-modules}
Let $(V, \mm)$ be an almost setup.
\begin{lemenum}
	\item \label{rslt:V-almost-localization-functor-properties} The localization functor $(V,\mm)^*\colon M \mapsto M^a$ on $\Cond(V)$ maps injective objects to injective objects and preserves all limits and colimits.

	\item \label{rslt:V-right-adjoint-of-almost-localization-properties} The localization functor $(V,\mm)^*$ admits a fully faithful right adjoint $(V,\mm)_*\colon N \mapsto N_*$ which satisfies
	\begin{align*}
		(M^a)_* = \IHom_V(\widetilde\mm, M)
	\end{align*}
	for all $M \in \Cond(V)$.

	\item \label{rslt:V-left-adjoint-of-almost-localization-properties} The localization functor $(V,\mm)^*$ admits an exact fully faithful left adjoint defined as
	\begin{align*}
		(V,\mm)_!\colon \Cond(V,\mm) \to \Cond(V), \qquad M \mapsto M_! := \widetilde\mm \tensor_V M_*.
	\end{align*}

	\item There is a natural closed symmetric monoidal structure on $\Cond(V, \mm)$ such that $M \mapsto M^a$ is symmetric monoidal.

	\item \label{rslt:properties-of-kappa-condensed-almost-V-modules} For every strong limit cardinal $\kappa$, $\Cond(V, \mm)_\kappa$ is a Grothendieck abelian category and we have $\Cond(V,\mm)_\kappa = \Cond(\mathcal C)_\kappa$, where $\mathcal C \subset \Cond(V,\mm)$ is the full subcategory of discrete almost $V$-modules. There is a cofinal class of $\kappa$'s such that both $(V,\mm)_*$ and $(V,\mm)_!$ restrict to functors $\Cond(V,\mm)_\kappa \to \Cond(V)_\kappa$.
\end{lemenum}
\end{lemma}
\begin{proof}
For every strong limit cardinal $\kappa$, the subcategory $\Sigma_\kappa \subset \Cond(V)_\kappa$ is stable under colimits and $\Cond(V)_\kappa$ is a Grothendieck abelian category, hence by \cite[Proposition 4.6.3]{popescu-abelian-categories} $\Sigma_\kappa$ is a localizing subcategory, so by \cite[Proposition 4.6.2]{popescu-abelian-categories} the quotient $\Cond(V, \mm)_\kappa = \Cond(V)_\kappa/\Sigma_\kappa$ is a Grothendieck abelian category and the localization functor $(V, \mm)^*$ admits a fully faithful right adjoint $(V, \mm)_{\kappa*}\colon \Cond(V,\mm)_\kappa \to \Cond(V)_\kappa$. Note that a static $V$-module $M$ is almost zero if and only if $M(S)$ is almost zero for all profinite sets $S$; this implies $\Cond(V,\mm)_\kappa = \Cond(\mathcal C)_\kappa$ as in (v).

For every $M \in \Cond(V)_\kappa$ and every profinite set $S$ we use \cref{rslt:compute-almost-Hom-over-V} to compute
\begin{align*}
	(M^a)_{\kappa*}(S) &= \Hom(V[S], (M^a)_{\kappa*}) = \Hom(V[S]^a, M^a) = \Hom(\widetilde\mm \tensor_V V[S], M)\\
	&= \IHom_{V,\kappa}(\widetilde\mm, M)(S).
\end{align*}
Let $\lambda$ be a regular cardinal such that $\widetilde\mm$ is a $\lambda$-small colimit of copies of $V$. By \cref{rslt:condensed-objects-in-presentable-monoidal-cat} there is a cofinal class of strong limit cardinals $\kappa$ such that $\Cond(V)_\kappa \subset \Cond(V)$ is stable under $\lambda$-small limits. Pulling the colimit in $\widetilde\mm$ out of the $\IHom$, for this cofinal class of $\kappa$ we see that $\IHom_{V,\kappa}(\widetilde\mm, M)$ is independent of $\kappa$ and hence also equals $\IHom_V(\widetilde\mm, M)$. This proves the existence of $(V,\mm)_*$ and finishes the proof of (v) and (ii).

To prove (iii), note that it follows easily from \cref{rslt:compute-almost-Hom-over-V} that the given functor $(V,\mm)_!$ is left adjoint to $(V,\mm)^*$. In particular it is right exact. One easily verifies that the unit of the adjunction $\id \to (V,\mm)^* (V,\mm)_!$ is an isomorphism, so that $(V,\mm)_!$ is fully faithful. It follows from the adjunction that $(V,\mm)_!$ is right exact. But it is also left exact, because $M \to M_*$ is left exact (by adjunction) and $\widetilde\mm \tensor_V -$ is exact by flatness of $\widetilde\mm$.

Part (i) is a formal consequence of (ii) and (iii). It remains to prove (iv). We are forced to define $M^a \tensor N^a := (M \tensor N)^a$ for any $M, N \in \Cond(V)$ and one checks easily that this works.
\end{proof}

We have now defined and studied the category of \emph{static} almost modules, i.e. those without any derived structure. Following the convention in this thesis, a general almost module will be an object in the derived $\infty$-category of $\Cond(V,\mm)$. We will now introduce this $\infty$-category and show its basic properties.

\begin{definition}
Let $(V, \mm)$ be an almost setup. For every strong limit cardinal $\kappa$ we let $\D(V,\mm)_\kappa := \D(\Cond(V,\mm)_\kappa)$, the derived $\infty$-category of $\Cond(V,\mm)_\kappa$. We define
\begin{align*}
	\D(V, \mm) := \varinjlim_\kappa \D(V,\mm)_\kappa,
\end{align*}
the derived $\infty$-category of $\Cond(V, \mm)$. The objects in $\D(V,\mm)$ are called \emph{$(V,\mm)$-modules}. If the ideal $\mm$ is clear from context, we will also call them \emph{almost $V$-modules} or \emph{$V^a$-modules}.
\end{definition}

\begin{remark}
By \cref{rslt:Cond-commutes-with-D} we have $\D(V,\mm) = \Cond(\D(V,\mm)_\omega)$.
\end{remark}

\begin{proposition} \label{rslt:properties-of-derived-almost-V-modules}
Let $(V,\mm)$ be an almost setup.
\begin{propenum}
	\item $\D(V, \mm)$ is a stable closed symmetric monoidal $\infty$-category which has all small limits and colimits. It is naturally equipped with a left-complete $t$-structure which is compatible with products and filtered colimits.

	\item \label{rslt:almost-localization-on-V-modules} The localization functor $\Cond(V) \to \Cond(V,\mm)$ induces a symmetric monoidal $t$-exact functor
	\begin{align*}
		(V,\mm)^*\colon \D(V) \to \D(V,\mm), \qquad M \mapsto M^a
	\end{align*}
	which preserves all small limits and colimits. Moreover, for any $M, N \in \D(V)$ there is a natural isomorphism
	\begin{align*}
		\IHom_{V^a}(M^a, N^a) = \IHom_V(M, N)^a.
	\end{align*}

	\item \label{rslt:lower-star-on-almost-V-modules} The functor $(V,\mm)^*$ admits a right adjoint
	\begin{align*}
		(V,\mm)_*\colon \D(V,\mm) \to \D(V), \qquad M \mapsto M_*,
	\end{align*}
	which is left $t$-exact and fully faithful. Moreover, for any $M \in \D(V)$ there is a natural isomorphism
	\begin{align*}
		(M^a)_* = \IHom_V(\widetilde\mm, M).
	\end{align*}

	\item The functor $(V,\mm)^*$ admits a left adjoint
	\begin{align*}
		(V,\mm)_!\colon \D(V,\mm) \to \D(V), \qquad M \mapsto M_!,
	\end{align*}
	which is $t$-exact and fully faithful. It can be computed as
	\begin{align*}
		M_! = \widetilde\mm \tensor_V M_*.
	\end{align*}

	\item Let $\tau$ be a regular cardinal such that $(V,\mm)$ is $\tau$-compact. For every strong limit cardinal $\kappa$, $\D(V,\mm)_\kappa$ is a stably presentable $\infty$-category which is stable under all colimits and the symmetric monoidal structure in $\D(V,\mm)$. It is $\tau$-compactly generated by the $\tau$-compact objects $V^a[S]$ for $\kappa$-small extremally disconnected sets $S$. For every cardinal $\lambda$ there is a cofinal class of strong limit cardinals $\kappa$ such that $\D(V,\mm)_\kappa$ is stable under $\lambda$-small limits in $\D(V,\mm)$. Moreover, there is also a cofinal class of $\kappa$'s such that both $(V,\mm)_*$ and $(V,\mm)_!$ restrict to functors $\D(V,\mm)_\kappa \to \D(V)_\kappa$.
\end{propenum}
\end{proposition}
\begin{proof}
By \cref{rslt:properties-of-kappa-condensed-almost-V-modules,rslt:condensed-objects-in-presentable-monoidal-cat,rslt:Cond-commutes-with-D} $\D(V,\mm) = \Cond(\D(\mathcal C))$ is a stable closed symmetric monoidal $\infty$-category which has all limits and colimits and everything in (v) except the final statement are clear (for the claim about $\tau$-compact generators, note that $V^a$ is a $\tau$-compact generator in the derived category of almost $V$-modules by \cref{rslt:compute-almost-Hom-over-V}). To show that the $t$-structure on $\D(V,\mm)$ is left-complete, it is by \cite[Proposition 1.2.1.19]{lurie-higher-algebra} enough to show that the $t$-structure is compatible with (countable) products, i.e. that products in $\Cond(V,\mm)$ are exact. This boils down to the statement that a product of almost surjective maps in $\Cond(V)$ is almost surjective, which is clear. This finishes the proof of (i).

To prove (ii), note that since the almost localization functor $\Cond(V) \to \Cond(V,\mm)$ is exact, it can easily be upgraded to a $t$-exact functor $\D(V) \to \D(V,\mm)$. This functor is in particular exact, hence to show that it preserves limits and colimits, it is enough to show that it preserves products and coproducts. As the $t$-structures on $\D(V)$ and $\D(V,\mm)$ are compatible with products and coproducts, this can be checked on homotopy groups and hence reduces to the abelian level, where it is \cref{rslt:V-almost-localization-functor-properties}. To prove the identity on inner hom's, note first that by adjunctions we get a natural morphism $\IHom_V(M, N)^a \to \IHom_{V^a}(M^a, N^a)$ in $\D(V^a)$. We can pull out colimits in $M$ on both sides and hence reduce to the case that $M = V[S]$ for some extremally disconnected set $S$. Choosing a K-injective complex representing $N$ (and using \cref{rslt:V-almost-localization-functor-properties}) we reduce the claim to the abelian level. Now the claim follows easily using Yoneda and \cref{rslt:compute-almost-Hom-over-V}.

To prove (iii), note that for every strong limit cardinal $\kappa$, the adjoint functor theorem provides a right adjoint $(V,\mm)_{\kappa*}\colon \D(V,\mm)_\kappa \to \D(V)_\kappa$ of $(V,\mm)^*$. It can also be described as the right derived functor of the functor $M \mapsto M_{\kappa*}$ on the abelian level. By choosing K-injective resolutions it follows that for the cofinal class of $\kappa$'s in \cref{rslt:properties-of-kappa-condensed-almost-V-modules}, $M_{\kappa*}$ is independent of $\kappa$, so that we get the desired functor $M \to M_*$ on $\D(V,\mm) \to \D(V)$. As $M \mapsto M^a$ preserves injectives (by \cref{rslt:V-almost-localization-functor-properties}) the functor $M \mapsto (M^a)_{\kappa*}$ is the right derived functor of the corresponding functor on abelian level, hence the identity $(M^a)_{\kappa*} = \IHom_{V,\kappa}(\widetilde\mm, M)$ follows from the similar identity in \cref{rslt:V-right-adjoint-of-almost-localization-properties}. To check that $M \mapsto M_{\kappa*}$ is fully faithful we need to see that the counit $(M_*)^a \isoto M$ is an isomorphism for all $M \in \D(V,\mm)$. But given $M$, we can find a K-injective complex $I^\bullet$ representing $M$ and then $(M_*)^a = (I^\bullet_*)^a = I^\bullet$, because $(-)_*$ is fully faithful on the abelian level.

To prove (iv), note that for every strong limit cardinal $\kappa$ the adjoint functor theorem provides a left adjoint $(V,\mm)_{\kappa!}\colon \D(V,\mm)_\kappa \to \D(V)_\kappa$ of $(V,\mm)^*$. It can also be described as the left derived functor of the functor $M \mapsto M_{\kappa!}$ on the abelian level. The latter functor is exact (see \cref{rslt:V-left-adjoint-of-almost-localization-properties}), hence $(V,\mm)_{\kappa!}$ is $t$-exact. Thus it can be computed directly on homology groups, so for the cofinal class of $\kappa$'s in \cref{rslt:properties-of-kappa-condensed-almost-V-modules} $(V,\mm)_{\kappa!}$ is indenpendent of $\kappa$ and hence defines the desired functor $(V,\mm)_!$. As both $(V,\mm)_!$ and $(V,\mm)^*$ are $t$-exact, the full faithfulness of $(V,\mm)_!$ reduces easily to the abelian level, where it is shown in \cref{rslt:V-left-adjoint-of-almost-localization-properties}. The identity $M_! = \widetilde\mm \tensor_V M_*$ can be checked by representing $M$ by some K-injective object as in the proof of (iii) (note that there is a natural map $M_! \to \widetilde\mm \tensor_V M_*$ by adjunctions).
\end{proof}

With $\D(V,\mm)$ at hand we now have a good basis for almost mathematics. Next up we introduce $(V,\mm)$-algebras to get a derived condensed version of almost ring theory. In the case of associative (i.e. not necessarily commutative) rings the definitions are straightforward. For the following definition, recall the associative $\infty$-operad $\opAssoc$ (see \cite[Definition 4.1.1.3]{lurie-higher-algebra}).

\begin{definition}
Let $(V,\mm)$ be an almost setup.
\begin{defenum}
	\item An \emph{associative $(V,\mm)$-algebra} is an $\opAssoc$-algebra object in the symmetric monoidal $\infty$-category $\D_{\ge0}(V,\mm)$. We denote
	\begin{align*}
		\AssRing_{(V,\mm)} := \Alg_{\opAssoc}(\D_{\ge0}(V,\mm))
	\end{align*}
	the $\infty$-category of associative $(V,\mm)$-algebras. In the case $\mm = V$ we simply write $\AssRing_V$ for $\AssRing_{(V,V)}$. If moreover $V = \Z$ then we abbreviate $\AssRing := \AssRing_\Z$ and call its object \emph{associative rings}.

	\item For every $(V,\mm)$-algebra $A$ we let
	\begin{align*}
		\D(A) := \LMod_A(\D(V,\mm))
	\end{align*}
	denote the stable $\infty$-category of $A$-modules in $\D(V,\mm)$. For every strong limit cardinal $\kappa$ such that $A \in \D(V,\mm)_\kappa$ we let $\D(A)_\kappa \subset \D(A)$ denote the full subcategory of those $A$-modules which lie in $\D(V,\mm)_\kappa$.

	\item \label{def:mm-A-and-A-S-for-ass-ring-A} For every $V$-algebra $A$ and every profinite set $S$ we denote the following objects in $\D(A)$:
	\begin{align*}
		A[S] := V[S] \tensor A, \qquad \mm_A := \mm \tensor A, \qquad \widetilde\mm_A := \widetilde\mm \tensor A.
	\end{align*}
\end{defenum}
\end{definition}

The almost localization functor and its right adjoint induce analogous functors on associative rings, as follows:

\begin{lemma} \label{rslt:every-ass-almost-algebra-is-A-a-of-actual-algebra}
Let $(V,\mm)$ be an almost setup. The localization functor $(V,\mm)^*\colon \D_{\ge0}(V) \to \D_{\ge0}(V,\mm)$ induces an essentially surjective functor
\begin{align*}
	(V,\mm)^*\colon \AssRing_V \to \AssRing_{(V,\mm)}, \qquad A \mapsto A^a.
\end{align*}
It admits a fully faithful right adjoint
\begin{align*}
	(V,\mm)_{**}\colon \AssRing_{(V,\mm)} \injto \AssRing_V, \qquad A \mapsto A_{**}
\end{align*}
which satisfies $A_{**} = \tau_{\ge0} A_*$ on the underlying modules.
\end{lemma}
\begin{proof}
We will implicitly make use of \cref{rslt:properties-of-derived-almost-V-modules}. Since $(V,\mm)^*$ is symmetric monoidal, it induces a functor $\AssRing_V \to \AssRing_{(V,\mm)}$. Let $\AssRing_{(V,\mm),\kappa}$ denote the full subcategory of those associative $(V,\mm)$-algebras which lie in $\D(V,\mm)_\kappa$. For a cofinal class of $\kappa$ the functor $\tau_{\ge0}(V,\mm)_*$ restricts to a functor $\D_{\ge0}(V,\mm)_\kappa \to \D_{\ge0}(V)_\kappa$ and thus provides a right adjoint of $(V,\mm)^*$ on $\D_{\ge0}(-)_\kappa$. Then by \cref{rslt:automatic-lax-monoidal-structure-of-right-adjoint} $\tau_{\ge0} (V,\mm)_*$ induces a functor
\begin{align*}
	(V,\mm)_*\colon \AssRing_{(V,\mm),\kappa} \to \AssRing_{V,\kappa}, \qquad A \mapsto A_{**}
\end{align*}
which is right adjoint to $(V,\mm)^*\colon \AssRing_{V,\kappa} \to \AssRing_{(V,\mm),\kappa}$. Taking the colimit over all $\kappa$ we arrive at the desired right adjoint of $(V,\mm)^*$ on $\AssRing_V$. For every $(V,\mm)$-algebra $A$ the natural map $(A_{**})^a \isoto A$ is an isomorphism: Since the forgetful functor $\AssRing_{(V,\mm)} \to \D_{\ge0}(V,\mm)$ is conservative this can be checked on the underlying connective $(V,\mm)$-modules, where it follows from the full faithfulness of $\tau_{\ge0} (V,\mm)_*$.
\end{proof}

\begin{remark}
The functor $(V,\mm)^*\colon \AssRing_V \to \AssRing_{(V,\mm)}$ should also admit a left adjoint $A \mapsto A_{!!}$ as in \cite[Proposition 2.2.27]{almost-ring-theory}, but we do not need it.
\end{remark}

\begin{remark}
Given an almost setup $(V,\mm)$ and an associative $V$-algebra $A$ with $A \in \D(V)_\kappa$ for some strong limit cardinal $\kappa$, the $\infty$-category $\D(A)$ is generated by the compact projective objects $A[S]$ for $\kappa$-small extremally disconnected sets $S$. This implies that every $M \in \D_{\ge0}(A)$ is a sifted colimit of the $A[S]$'s, in an essentially unique way -- this will be used implicitly below. Note that the same is wrong for $\D(A^a)$: The objects $A^a[S]$ are neither compact nor projective in general. Thus, while it is still true that every $M \in \D_{\ge0}(A)$ is a sifted colimit of $A^a[S]$'s (this follows from \cref{rslt:general-properties-of-almost-localization-over-ass-alg} below), this sifted colimit is not essentially unique anymore.
\end{remark}

We now study the basic properties of associative almost rings and their associated modules. We observe that the properties of $\D(V,\mm)$ listed in \cref{rslt:properties-of-derived-almost-V-modules} generalize to similar properties of $\D(A^a)$ for every associative $(V,\mm)$-algebra $A^a$:

\begin{proposition} \label{rslt:general-properties-of-almost-localization-over-ass-alg}
Let $\tau$ be a regular cardinal, $(V, \mm)$ a $\tau$-compact almost setup and $A \in \AssRing_V$.
\begin{propenum}
	\item $\D(A^a)$ is a stable $\infty$-category which has all small limits and colimits. It comes equipped with a left complete $t$-structure which is compatible with products and filtered colimits.

	\item \label{rslt:construction-of-almost-localization-of-modules-over-ass-V-m-alg} There is a natural localization functor
	\begin{align*}
		(A,\mm_A)^*\colon \D(A) \to \D(A^a), \qquad M \mapsto M^a,
	\end{align*}
	which reduces to $(V,\mm)^*$ on underlying $V$-modules. The functor $(A,\mm_A)^*$ is $t$-exact and preserves all small limits and colimits.

	\item \label{rslt:construction-of-right-adjoint-for-modules-over-ass-V-m-alg} The functor $(A,\mm_A)^*$ admits a left $t$-exact fully faithful right adjoint
	\begin{align*}
		(A,\mm_A)_*\colon \D(A^a) \to \D(A), \qquad M \mapsto M_*
	\end{align*}
	which reduces to $(V,\mm)_*$ on underlying $(V,\mm)$-modules.

	\item The functor $(A,\mm_A)^*$ admits a $t$-exact fully faithful left adjoint
	\begin{align*}
		(A,\mm_A)_!\colon \D(A^a) \to \D(A), \qquad M \mapsto M_!,
	\end{align*}
	which reduces to $(V,\mm)_!$ on underlying $(V,\mm)$-modules.

	\item For every strong limit cardinal $\kappa$ such that $A \in \D(V)_\kappa$, $\D(A^a)_\kappa$ is a stably presentable $\infty$-category which is stable under all colimits. It is $\tau$-compactly generated by the $\tau$-compact objects $A^a[S]$ for $\kappa$-small extremally disconnected sets $S$. For every cardinal $\lambda$ there is a cofinal class of strong limit cardinals $\kappa$ such that $\D(A^a)_\kappa$ is stable under $\lambda$-small limits in $\D(A^a)$. Moreover, there is also a cofinal class of $\kappa$'s such that both $(A,\mm_A)_*$ and $(A,\mm_A)_!$ restrict to functors $\D(A^a)_\kappa \to \D(A)_\kappa$.
\end{propenum}
\end{proposition}
\begin{proof}
We will implicitly make use of \cref{rslt:properties-of-derived-almost-V-modules}. Part (i) and everything except the last statement of part (v) follow immediately from the corresponding statements about $\D(V,\mm)$ (using that the forgetful functor $\D(A^a) \to \D(V,\mm)$ is conservative and preserves all small limits and colimits).

To prove (ii), note first that since $(V,\mm)^*$ is symmetric monoidal, it induces a functor $(A,\mm_A)^*\colon \D(A) \to \D(A^a)$ (e.g. take the restriction to $\LMod_A(\D(V))$ of the functor $\LMod(\D(V)) = \Alg_{\opLM}(\D(V)^\tensor) \to \Alg_{\opLM}(\D(V,\mm)^\tensor) = \LMod(\D(V,\mm))$). The forgetful functors $\D(A) \to \D(V)$ and $\D(A^a) \to \D(V,\mm)$ preserve all small limits and colimits and are conservative. Thus, as $(V,\mm)^*$ preserves limits and colimits, so does $(A,\mm_A)^*$.

To prove (iii), note that the right adjoint $(A,\mm_A)_*$ exists by the adjoint functor theorem (a priori only on $\kappa$-condensed objects, but as shown below $(A,\mm_A)_*$ is just $(V,\mm)_*$ on underlying $V^a$-modules and hence independent of $\kappa$ for large $\kappa$). Moreover, by definition of $(A,\mm_A)^*$ there is a natural equivalence of functors
\begin{align*}
	(A,\mm_A)^* \comp (- \tensor_V A) \isom (- \tensor_{V^a} A^a) \comp (V,\mm)^*.
\end{align*}
Passing to right adjoints shows that $(A,\mm_A)_*$ is just $(V,\mm)_*$ on underlying $V^a$-modules. But the forgetful functor $\D(A^a) \to \D(V^a)$ is conservative, hence $(A,\mm_A)^*(A,\mm_A)_* M \isoto M$ is an isomorphism for all $M \in \D(A^a)$ (because the same is true for $(V,\mm)_*$), i.e. $(A,\mm_A)_*$ is fully faithful.

To prove (iv), fix a strong limit cardinal $\kappa$ such that $A \in \D(V)_\kappa$. Then by the adjoint functor theorem there is a left adjoint $(A,\mm_A)_{\kappa!}\colon \D(A^a)_\kappa \to \D(A)_\kappa$. Let us denote $r\colon \D(V)_\kappa \to \D(A)_\kappa$ and $r^a\colon \D(V^a)_\kappa \to \D(A^a)_\kappa$ the right adjoints to the forgetful functors. By adjunctions there is a natural morphism
\begin{align*}
	(A,\mm_A)^* \comp r \to r^a \comp (V,\mm)^*
\end{align*}
of functors $\D(V)_\kappa \to \D(A^a)_\kappa$. We claim that this is an isomorphism. This can be checked after applying the forgetful functor $\D(A^a)_\kappa \to \D(V^a)_\kappa$. But note that the composition of $r^a$ with the forgetful functor is just the functor $\IHom_{V^a,\kappa}(A^a, -)$ (and similarly for $r$), so the claimed isomorphism reduces to showing that the natural morphism
\begin{align*}
	\IHom_{V,\kappa}(A, -)^a \isoto \IHom_{V^a,\kappa}(A^a, (-)^a)
\end{align*}
is an isomorphism. But this was shown in \cref{rslt:almost-localization-on-V-modules}. Passing to left adjoints in the above isomorphism of functors shows that $(A,\mm_A)_{\kappa!}$ is computed as $(V,\mm)_{\kappa!}$ on underlying $V^a$-modules. In particular, on a cofinal class of $\kappa$'s, $(A,\mm_A)_{\kappa!}$ is independent of $\kappa$ and hence defines the desired functor $(A,\mm_A)_!$. As in the proof of (iii) it follows that $(A,\mm_A)_!$ is fully faithful.

The last statement of part (v) follows from the proofs of (iii) and (iv).
\end{proof}

We now turn our focus on (commutative) $(V,\mm)$-algebras. They are more subtle to define than associative $(V,\mm)$-algebras, for several reasons. The first reason is that there are two possible ways to define commutative rings in the $\infty$-categorical setting: We could either take commutative algebra objects in $\D_{\ge0}(\Z)$, or we could animate (cf. \cite[\S11]{scholze-analytic-spaces}) the category of (static) commutative rings in $\D(\Z)^\heartsuit$. Both approaches have their advantages and disadvantages. The first approach is somewhat easier to setup, mainly because commutative algebra objects in a symmetric monoidal $\infty$-category are usually easier to construct (and appear more ``naturally'') than animated rings. However, animated rings behave more like classical rings when it comes to algebraic geometry, and therefore seem to provide a better theory of schemes; most prominently, the polynomial rings $\Z[T]$ have the expected universal property in the $\infty$-category of animated rings, while the same is false in the $\infty$-category of commutative algebra objects (see also \cite[\S25]{lurie-spectral-algebraic-geometry} for more remarks on the difference between both approaches). As we are interested in algebraic geometry, we chose to work with animated rings throughout.

Let us start with the definition of (commutative) rings:

\begin{definition} \label{def:rings}
\begin{defenum}
	\item A \emph{(commutative) ring} is an animated condensed ring. That is, if we denote $\Ring^\heartsuit_\omega$ the ordinary category of classical commutative rings, then we define the $\infty$-category $\Ring$ as
	\begin{align*}
		\Ring := \Ani(\Cond(\Ring^\heartsuit_\omega)) = \Cond(\Ani(\Ring^\heartsuit_\omega))
	\end{align*}
	(see \cite[Lemma 11.8]{scholze-analytic-spaces} for the right-hand identity). The $\infty$-category $\Ring$ is freely generated under small sifted colimits by retracts of the rings $\Z[\N[S]]$ for extremally disconnected sets $S$. In particular, if for every strong limit cardinal $\kappa$ we denote
	\begin{align*}
		\Ring_\kappa := \Cond(\Ani(\Ring^\heartsuit_\omega))_\kappa
	\end{align*}
	then $\Ring = \varinjlim_\kappa \Ring_\kappa$ and each $\Ring_\kappa$ is compactly generated.

	\item There is a natural forgetful functor
	\begin{align*}
		(-)^\circ\colon \Ring \to \CAlg(\D_{\ge0}(\Z)), \qquad A \mapsto A^\circ
	\end{align*}
	which is conservative and preserves all small limits and colimits (by the same argument as in \cite[Proposition 25.1.2.2]{lurie-spectral-algebraic-geometry}). In particular, to every ring $A \in \Ring$ we have an associated symmetric monoidal stable $\infty$-category
	\begin{align*}
		\D(A) := \Mod_{A^\circ}(\D(\Z))
	\end{align*}
	of $A$-modules, which is functorial in $A$ (see \cite[Theorem 4.5.3.1]{lurie-higher-algebra}). The $\infty$-category $\D(A)$ depends only on the underlying associative algebra of $A$. In particular we obtain $\mm_A, \widetilde\mm_A, A[S] \in \D(A)$ from \cref{def:mm-A-and-A-S-for-ass-ring-A}.

	\item Let $A \in \Ring$ be a ring. An $A$-algebra is an object of $\Ring_A := \Ring_{A/}$. Note that if $A = V$ is a classical ring then $\Ring_V$ can also be defined as
	\begin{align*}
		\Ring_V = \Cond(\Ani(\Ring_{V,\omega}^\heartsuit))
	\end{align*}
	(see \cite[Proposition 25.1.4.2]{lurie-spectral-algebraic-geometry}) hence all the definitions and remarks from (i) and (ii) also apply to $\Ring_V$.
\end{defenum}
\end{definition}

Next up we introduce (commutative) almost rings. Since we are working with animated rings, the notion of almost rings is a bit subtle. We start with the definition:

\begin{definition} \label{def:almost-rings}
Let $(V,\mm)$ be an almost setup.
\begin{defenum}
	\item We define
	\begin{align*}
		\Ring_{(V,\mm)} \subset \Ring_V
	\end{align*}
	to be the full subcategory spanned by those $V$-algebras $A$ whose underlying $V$-module lies in the essential image of $\tau_{\ge0}(V,\mm)_*\colon \D_{\ge0}(V,\mm) \to \D_{\ge0}(V)$. The objects of $\Ring_{(V,\mm)}$ are called the \emph{$(V,\mm)$-algebras}. For every strong limit cardinal $\kappa$ we denote
	\begin{align*}
		\Ring_{(V,\mm),\kappa} \subset \Ring_{(V,\mm)}.
	\end{align*}
	the full subcategory spanned by those $(V,\mm)$-algebras $A$ such that $A^a \in \D(V,\mm)_\kappa$, where $A^a$ denotes the almost localization of the underlying $V$-module of $A$.

	\item There is a natural forgetful functor
	\begin{align*}
		(-)^\circ\colon \Ring_{(V,\mm)} \to \CAlg(\D_{\ge0}(V,\mm)), \qquad A \mapsto A^\circ
	\end{align*}
	given as the composition
	\begin{align*}
		\Ring_{(V,\mm)} \injto \Ring_V \xto{(-)^\circ} \CAlg(\D_{\ge0}(V)) \xto{(-)^a} \CAlg(\D_{\ge0}(V,\mm)).
	\end{align*}
	In particular, for every $(V,\mm)$-algebra $A$ we get an associated symmetric monoidal $\infty$-category
	\begin{align*}
		\D(A) := \Mod_{A^\circ}(\D(V,\mm))
	\end{align*}
	which is functorial in $A$.

	\item For every $(V,\mm)$-algebra $A$ we denote $\Ring_A := (\Ring_{(V,\mm)})_{A/}$ and call its objects the \emph{$A$-algebras}.
\end{defenum}
\end{definition}

While defining $\Ring_{(V,\mm)}$ is straightforward, it is a priori unclear that this $\infty$-category behaves as expected. For example, we would expect it to admit all small colimits and there should also exist a localization functor $\Ring_V \to \Ring_{(V,\mm)}$. Our next goal is to prove these statements. They rely on a general construction of animated rings introduced in the proof of \cite[Proposition 12.26]{scholze-analytic-spaces}. We extracted the argument and transformed it into the following result. To understand the statement, recall that for every ring $A$ and every prime power $p^m$ there is a natural \emph{Frobenius homomorphism} $\varphi_{p^m}\colon A \to A/p$ of rings, which intuitively acts as $x \mapsto x^{p^m}$. Namely, this homomorphism can be constructed as a natural transformation on the ordinary category of compact projective generators of $\Ring$, where the explicit formula provides the desired map.

\begin{lemma} \label{rslt:construct-left-adjoint-on-animated-rings}
Let $A$ be a ring and let $\mathcal C \subset \D_{\ge0}(A)$ be a full subcategory satisfying the following properties:
\begin{enumerate}[(a)]
	\item $\mathcal C$ is stable under cofibers in $\D_{\ge0}(A)$ and the inclusion $\mathcal C \injto \D_{\ge0}(A)$ admits a left adjoint $L\colon \D_{\ge0}(A) \to \mathcal C$.

	\item For every $M, Q \in \D_{\ge0}(A)$ such that $L(Q) = 0$ we have $L(Q \tensor M) = 0$.

	\item For every prime power $p^m$ and every $Q \in \D_{\ge0}(A)$ such that $L(Q) = 0$ we have $L(\varphi_{p^m}^*Q) = 0$, where the map $\varphi_{p^m}\colon A \to A/p$ is the Frobenius $x \mapsto x^{p^m}$ and we view $\varphi_{p^m}^*Q \in \D_{\ge0}(A/p)$ as an object in $\D_{\ge0}(A)$ via the forgetful functor.
\end{enumerate}
Let $\Ring_{\mathcal C} \subset \Ring_A$ be the full subcategory spanned by those $A$-algebras $B$ such that the underlying object in $\D_{\ge0}(A)$ lies in $\mathcal C$. Then the inclusion $\Ring_{\mathcal C} \injto \Ring_A$ admits a left adjoint $L\colon \Ring_A \to \Ring_{\mathcal C}$ which acts as the functor $L\colon \D_{\ge0}(A) \to \mathcal C$ on the underlying $A$-modules.
\end{lemma}
\begin{proof}
We argue as in the proof of \cite[Proposition 12.26]{scholze-analytic-spaces}. Namely, note that there is an adjunction
\begin{align*}
	F\colon \D_{\ge0}(A) \rightleftarrows \Ring_A \noloc G,
\end{align*}
where $G$ is the forgetful functor (to construct $F$, restrict to $\kappa$-condensed objects on both sides, where $\kappa$ is large enough such that $A \in \Ring_\kappa$, and apply the adjoint functor theorem; then take the colimit over $\kappa$). The functor $G$ preserves all sifted colimits because so does the forgetful functor $\CAlg(\D_{\ge0}(A)) \to \D_{\ge0}(A)$. By the Barr-Beck theorem (see \cite[Theorem 4.7.3.5]{lurie-higher-algebra}) it follows that the adjunction of $F$ and $G$ is monadic, i.e. there is a monad $T$ on $\D_{\ge0}(A)$ which acts as $T(M) = GF(M)$ such that $\Ring_A \isom \LMod_T(\D_{\ge0}(A))$.

We claim that $T$ restricts to a monad $T_{\mathcal C}$ on $\mathcal C$ with $T_{\mathcal C}(M) = LT(M)$. This boils down to showing that if $M \to N$ is a map in $\D_{\ge0}(A)$ such that $L(M) \isoto L(N)$ is an isomorphism, then $LT(M) \to LT(N)$ is an isomorphism. By \cite[Construction 25.2.2.6]{lurie-spectral-algebraic-geometry} (which easily generalizes to the condensed setting) we have $T(M) = \bigdsum_{n\ge0} \Sym^n_A(M)$, so it is enough to show that the map $L(\Sym^n_A(M)) \to L(\Sym^n_A(N))$ is an isomorphism. This can be seen as in the proof of \cite[Lemma 12.27]{scholze-analytic-spaces}, which reduces the claim to showing that if $Q \in \D_{\ge0}(A)$ satisfies $L(Q) = 0$ then for every prime power $p^m$ we have $L(\varphi_{p^m}^* Q) = 0$. This is true by assumption.

It follows formally that $T$-modules in $\D_{\ge0}(A)$ whose underlying $A$-module lies in $\mathcal C \subset \D_{\ge0}(A)$ are the same as $T_{\mathcal C}$-modules, i.e. we have $\Ring_{\mathcal C} \isom \LMod_{T_{\mathcal C}}(\mathcal C)$. The functor $L\colon \mathcal D_{\ge0}(A) \to \mathcal C$ induces a functor $L\colon \LMod_T(\D_{\ge0}(A)) \to \LMod_{T_{\mathcal C}}(\mathcal C)$ which defines the desired left adjoint of $\Ring_{\mathcal C} \injto \Ring_A$.
\end{proof}

With \cref{rslt:construct-left-adjoint-on-animated-rings} at hand, we can now prove the desired properties of the $\infty$-category $\Ring_{(V,\mm)}$ of $(V,\mm)$-algebras:

\begin{proposition} \label{rslt:almost-localization-for-rings}
Let $(V,\mm)$ be an almost setup. Then the almost localization functor $(V,\mm)^*\colon \D_{\ge0}(V) \to \D_{\ge0}(V,\mm)$ induces an essentially surjective functor
\begin{align*}
	(V,\mm)^*\colon \Ring_V \to \Ring_{(V,\mm)}, \qquad A \mapsto A^a,
\end{align*}
It is right adjoint to the inclusion
\begin{align*}
	(V,\mm)_{**}\colon \Ring_{(V,\mm)} \injto \Ring_V, \qquad A \mapsto A_{**},
\end{align*}
which acts as $A_{**} = \tau_{\ge0} A_*$ on the underlying modules.
\end{proposition}
\begin{proof}
To show the existence of the functor $L = (V,\mm)^*\colon \Ring_V \to \Ring_{(V,\mm)}$ we will apply \cref{rslt:construct-left-adjoint-on-animated-rings} to $A = V$ and $\mathcal C = \D_{\ge0}(V,\mm) \subset \D_{\ge0}(V)$. We thus have to check that $\mathcal C = \D_{\ge0}(V,\mm)$ satisfies conditions (a), (b) and (c). It is clear that (a) is satisfied, and (b) follows immediately from the fact that $(V,\mm)^*$ is symmetric monoidal (see \cref{rslt:almost-localization-on-V-modules}).

It remains to prove (c), so fix a prime power $p^m$ with associated Frobenius $\varphi\colon V \to V/p$ (where the quotient $V/p$ is understood in the derived sense) and any $Q \in \D_{\ge0}(V)$ such that $Q^a = 0$. Denote $Q' := \varphi^* Q \in \D_{\ge0}(V/p)$. Then we know that $\varphi_* Q' \in \D_{\ge0}(V)$ is almost zero (because almost localization is symmetric monoidal and $\varphi_* Q'$ is just $Q$ tensored with the $V$-module $V/p$), i.e. for every $n \ge 0$ we have $\mm \cdot \pi_n(\varphi_* Q') = 0$. Here $\varphi_*$ is the forgetful functor, which is $t$-exact and hence $0 = \mm \cdot \pi_n(\varphi_* Q') = \mm \cdot \varphi_* \pi_n (Q')$. Now $\pi_n Q'$ is a static $\pi_0(V/p)$-module, so one sees easily that $0 = \mm \cdot \varphi_* \pi_n (Q') = \mm' \cdot \pi_n(Q')$, where $\mm' \subset \pi_0(V/p)$ is the ideal generated by the image of $\pi_0 \varphi\colon V \to \pi_0(V/p)$. In order to show that $Q'^a = 0$, i.e. $\mm \cdot \pi_n(Q') = 0$, we are thus reduced to showing that $\mm'$ is the same as the ideal in $\pi_0(V/p)$ generated by the image of $\mm$ under the projection $V \to \pi_0(V/p)$. This amounts to saying that $\mm$ is generated by the $p^m$-th powers of the elements in $\mm$, which is shown in \cite[Proposition 2.1.7.(ii)]{almost-ring-theory}.
\end{proof}

\begin{corollary} \label{rslt:basic-properties-of-almost-rings}
Let $(V,\mm)$ be an almost setup.
\begin{corenum}
	\item The $\infty$-category $\Ring_{(V,\mm)}$ admits all small limits and colimits and the localization $(V,\mm)^*\colon \Ring_V \to \Ring_{(V,\mm)}$ preserves them.

	\item The forgetful functor $(-)^\circ\colon \Ring_{(V,\mm)} \to \CAlg(\D_{\ge0}(V,\mm))$ preserves all small limits and colimits and is conservative.

	\item For every strong limit cardinal $\kappa$ the $\infty$-category $\Ring_{(V,\mm),\kappa}$ is presentable.
\end{corenum}
\end{corollary}
\begin{proof}
It is clear that $\Ring_{(V,\mm)}$ admits all small limits. The existence of the localization $(V,\mm)^*$ (by \cref{rslt:almost-localization-for-rings}) implies that it also admits all small colimits (because $\Ring_V$ has all small colimits). To show that $(V,\mm)^*$ preserves small limits, we have to show that the endofunctor $A \mapsto (A^a)_{**}$ of $\Ring_V$ preserves small limits. But this functor acts as $M \mapsto \tau_{\ge0} (M^a)_*$ on the underlying $V$-modules, which preserves small limits. This proves (i).

To prove (ii), note first that it is clear that the forgetful functor $(-)^\circ$ preserves small limits, as it is defined as a composition of functors which preserve small limits. To check that $(-)^\circ$ is conservative, let $f\colon A \to B$ be a map in $\Ring_{(V,\mm)}$ such that $A^\circ \isoto B^\circ$ is an isomorphism. In order to verify that $f$ is an isomorphism, it is enough to show that $f_{**}\colon A_{**} \to B_{**}$ is an isomorphism in $\Ring_V$. This can be checked after applying the forgetful functor $\Ring_V \to \D_{\ge0}(V)$. The underlying $V$-modules of $A_{**}$ and $B_{**}$ lie in the essential image of $\tau_{\ge0}(V,\mm)_*\colon \D_{\ge0}(V,\mm) \to \D_{\ge0}(V)$, hence the isomorphism can be checked after applying $(-)^a$. But the resulting map $(A_{**})^a \to (B_{**})^a$ is the same as the map obtained by applying the forgetful functor $F\colon \CAlg(\D_{\ge0}(V,\mm)) \to \D_{\ge0}(V,\mm)$ to $A^\circ \isoto B^\circ$. By a similar reasoning one checks that $(-)^\circ$ preserves all sifted colimits (using that $F$ preserves sifted colimits). Using that almost localization is symmetric monoidal we deduce that $(-)^\circ$ preserves pushouts as well, hence all small colimits.

It remains to prove (iii) so let $\kappa$ be given. We know that $\Ring_{V,\kappa}$ is presentable (even compactly generated), so it is enough to show that $\Ring_{(V,\mm),\kappa}$ is an accessible localization of $\Ring_{V,\kappa}$. Note that $(-)_\kappa \comp (V,\mm)_{**}$ provides a fully faithful embedding $\Ring_{(V,\mm),\kappa} \injto \Ring_{V,\kappa}$ (here $(-)_\kappa\colon \Ring_V \to \Ring_{V,\kappa}$ is the obvious restriction functor) which is left adjoint to the almost localization $(V,\mm)^*$. Thus it only remains to check that $\Ring_{(V,\mm),\kappa} \injto \Ring_{V,\kappa}$ is stable under $\tau$-filtered colimits for some regular cardinal $\tau$. We can in fact take any $\tau$ such that $(V,\mm)$ is $\tau$-compact because then $\D_{\ge0}(V,\mm) \subset \D_{\ge0}(V)$ is stable under $\tau$-filtered colimits, which follows from the explicit computation of $(-)_*$ in \cref{rslt:lower-star-on-almost-V-modules} and the $\tau$-compactness of $\widetilde\mm$.
\end{proof}

In the following we will usually use \cref{rslt:almost-localization-for-rings,rslt:basic-properties-of-almost-rings} without explicitly mentioning these results. Regarding the associated symmetric monoidal $\infty$-category of modules, we obtain the following properties extending \cref{rslt:general-properties-of-almost-localization-over-ass-alg}:

\begin{proposition} \label{rslt:general-properties-of-almost-localization-over-alg}
Let $(V, \mm)$ be an almost setup and $A \in \Ring_V$. The $\infty$-category $\D(A^a)$ depends only on the underlying associative algebra of $A^a$, hence all the results from \cref{rslt:general-properties-of-almost-localization-over-ass-alg} apply. We moreover get the following results on the symmetric monoidal structure:
\begin{propenum}
	\item $\D(A^a)$ is a closed symmetric monoidal $\infty$-category.

	\item \label{rslt:properties-of-almost-localization-over-alg} The almost localization functor $(A,\mm_A)^*\colon \D(A) \to \D(A^a)$ can be upgraded to a symmetric monoidal functor. Moreover for all $M, N \in \D(A)$ there is a natural isomorphism
	\begin{align*}
		\IHom_{A^a}(M^a, N^a) = \IHom_A(M, N)^a.
	\end{align*}

	\item \label{rslt:computation-of-almost-lower-star-over-alg} For all $M \in \D(A)$ there is a natural isomorphism
	\begin{align*}
		(M^a)_* = \IHom_A(\widetilde\mm_A, M) \in \D(A).
	\end{align*}

	\item \label{rslt:properties-of-almost-lower-shriek-over-alg} For all $M \in \D(A^a)$ we have
	\begin{align*}
		M_! = \widetilde\mm_A \tensor_A M_* \in \D(A).
	\end{align*}
	In particular, for all $N \in \D(A)$ we have $(N^a)_! = \widetilde\mm_A \tensor_A N$.

	\item For every strong limit cardinal $\kappa$ such that $A \in \D(V)_\kappa$, $\D(A^a)_\kappa$ is stable under the symmetric monoidal structure on $\D(A^a)$.
\end{propenum}
\end{proposition}
\begin{proof}
We make implicit use of \cref{rslt:properties-of-derived-almost-V-modules}. Part (v) follows from the corresponding statement for $\D(V,\mm)$. In (i) the only non-trivial statement is that $\D(A^a)$ admits internal homs. Since it admits all small limits, we reduce to the existence of $\IHom_{A^a}(A^a[S], M)$ for extremally disconnected sets $S$ and all $M \in \D(A^a)$. This object exists by the adjoint functor theorem when restricted to $\D(A^a)_\kappa$ and does not depend on large enough $\kappa$ because its underlying $V^a$-module is just $\IHom_{V^a}(V^a[S], M)$.

We now prove (ii). Let us first show that $(A,\mm_A)^*$ can be upgraded to a symmetric monoidal functor. Let $\mathcal C'^\tensor \to \Delta^1 \cprod \opComm$ denote the coCartesian family of symmetric monoidal $\infty$-categories classifying the functor $(V,\mm)^*\colon \D(V)^\tensor \to \D(V,\mm)^\tensor$. By \cite[Remark 3.3.3.16]{lurie-higher-algebra} there is a fibration of generalized $\infty$-operads $\Mod(\mathcal C^\tensor)^\tensor \to \CAlg(\mathcal C^\tensor) \cprod \opComm$. Here $\Alg(\mathcal C^\tensor)$ classifies the functor $\CAlg(\D_{\ge0}(V)) \to \CAlg(\D_{\ge0}(V,\mm))$; in particular it contains a coCartesian morphism $f\colon A^\circ \to A^{\circ a}$ over $\Delta^1$. Let $\mathcal M^\tensor := \Mod(\mathcal C^\tensor)^\tensor_f$. Then we have a fibration of generalized $\infty$-operads $p\colon \mathcal M^\tensor \to \Delta^1 \cprod \opComm$ whose fiber over $0 \in \Delta^1$ is $\D(A)^\tensor$ and whose fiber over $1 \in \Delta^1$ is $\D(A^a)^\tensor$. We claim that $p$ is a coCartesian fibration. We first check that it is locally coCartesian, so suppose that $e\colon (0, x) \to (1, y)$ is an edge in $\Delta^1 \cprod \opComm$ (where we can w.l.o.g. assume that $e$ lies over the edge $0 \to 1$ in $\Delta^1$). We can reduce to the case that $y = \langle 1 \rangle$, $x = \langle n \rangle$ and $x \to y$ is the unique active morphism (by the same arguments as in \cite[\S4.5.3]{lurie-higher-algebra}). In this case $\mathcal M^\tensor_e$ classifies the functor $\D(A)^n \to \D(A^a)$ given by $(M_1, \dots, M_n) \mapsto (A,\mm_A)^* (M_1 \tensor_A \dots \tensor_A M_n)$. This shows that $p$ is a locally coCartesian fibration. To show that it is a Cartesian fibration, we use \cite[Remark 2.4.2.9]{lurie-higher-topos-theory}, which reduces the claim to showing that for all $M_1, M_2 \in \D(A)$ the natural morphism $(A,\mm_A)^*(M_1 \tensor_A M_2) \isoto (A,\mm_A)^*M_1 \tensor_{A^a} (A,\mm_A)^*M_2$ (induced by $p$) is an isomorphism. By computing $\tensor_A$ in terms of a bar resolution we reduce to the claim that $(V,\mm)^*$ is symmetric monoidal. But this we know to be true.

We now prove the identity on inner hom's in (ii). First note that by adjunctions we get a natural morphism $\IHom_A(M, N)^a \to \IHom_{A^a}(M^a, N^a)$ in $\D(A^a)$. We can pull out colimits in $M$ on both sides and hence reduce to the case that $M = A[S]$ for some extremally disconnected set $S$. On underlying $V^a$-modules we get
\begin{align*}
	\IHom_A(A[S], N)^a = \IHom_V(V[S], N)^a, \qquad \IHom_{A^a}(A[S]^a, N^a) = \IHom_{V^a}(V[S]^a, N^a).
\end{align*}
Using that the forgetful functor $\D(A^a) \to \D(V^a)$ is conservative, we are now reduced to the case $A = V$, where we have already shown it. This finishes the proof of (ii).

We now prove (iii). The fact that $M \mapsto M^a$ commutes with inner hom's shows that $\IHom_A(\widetilde\mm_A, M)^a = M^a$ for any $M \in \D(A)$. This produces a natural map $\IHom_A(\widetilde\mm_A, M) \to (M^a)_*$ by adjunction. To prove that this map is an isomorphism, we can first apply the forgetful functor $\D(A) \to \D(V)$; but then it follows from the similar claim over $V$.

It remains to prove (iv). The claimed identity $M_! = \widetilde\mm_A \tensor_A M_*$ follows from the similar identity over $V$ by applying the forgetful functor $\D(A) \to \D(V)$ to the natural morphism $M_! \to \widetilde\mm_A \tensor_A M_*$. Now if $M = N^a$ for some $N \in \D(A)$ then we have a natural almost isomorphism $N \to (N^a)_*$. Thus it becomes an isomorphism after applying $\widetilde\mm_A \tensor_A -$, i.e. we get $\widetilde\mm_A \tensor_A N = \widetilde\mm_A \tensor_A (N^a)_* = (N^a)_!$.
\end{proof}

One would expect the functor $(A,\mm_A)^*$ defined in \cref{rslt:general-properties-of-almost-localization-over-ass-alg} to be functorial in $A$ and in $(V,\mm)$. This is correct, as the following statements show. In the following we denote by $\infcatinf^\tensor$ the $\infty$-category of symmetric monoidal $\infty$-categories and symmetric monoidal functors.

\begin{lemma} \label{rslt:functoriality-of-almost-modules}
There is a natural functor
\begin{align*}
	\AlmSetup \to \infcatinf^\tensor, \qquad (V,\mm) \mapsto \D(V,\mm).
\end{align*}
It satisfies the following properties:
\begin{lemenum}
	\item Let $\varphi\colon (V, V) \to (V, \mm)$ be a localizing morphism of almost setups. Then the induced symmetric monoidal functor $\D(V, V) \to \D(V, \mm)$ coincides with $(V,\mm)^*$.

	\item \label{rslt:functoriality-of-almost-modules-along-strict} Let $\varphi\colon (V,\mm) \to (V',\mm')$ be a strict morphism of almost setups. Denote $V^a := (V,\mm)^*V$ and $V'^a := (V,\mm)^* V'$. Then $\D(V',\mm') = \D(V'^a)$ and the induced symmetric monoidal functor $\D(V,\mm) \to \D(V',\mm')$ coincides with $- \tensor_{V^a} V'^a$.
\end{lemenum}
\end{lemma}
\begin{proof}
Viewing every classical ring $V$ as an algebra in $\D(\Z)$, we deduce from \cite[Theorem 4.5.3.1]{lurie-higher-algebra} that there is a natural functor $\AlmSetup \to \infcatinf^\tensor$ mapping $(V, \mm) \mapsto \D(V)$. Let $\mathcal C'^\tensor \to \AlmSetup \cprod \opComm$ be the corresponding coCartesian family of symmetric monoidal $\infty$-categories. Then let $\mathcal C^\tensor \subset \mathcal C'^\tensor$ be the full subcategory where we only allow objects in $\D(V,\mm) \subset \D(V)$ over each $(V,\mm) \in \AlmSetup$; here the inclusion $\D(V,\mm) \injto \D(V)$ is given by $(V,\mm)_*$ (see \cref{rslt:lower-star-on-almost-V-modules}). We claim that $p\colon \mathcal C^\tensor \to \AlmSetup \cprod \opComm$ is still a coCartesian family of symmetric monoidal $\infty$-categories and hence classifies the desired functor $\AlmSetup \to \infcatinf^\tensor$. To see this, we first show that $p$ is a locally coCartesian fibration, so let $e\colon ((V,\mm), x) \to ((V',\mm'), y)$ be a morphism in $\AlmSetup \cprod \opComm$. We need to show that $\mathcal C_e$ is a coCartesian fibration. We can w.l.o.g. assume that $x = \langle n \rangle$ and $y = \langle 1 \rangle$ and that the map $x \to y$ is active. In this case the fiber of $\mathcal C$ over $((V,\mm), x)$ is $\D(V,\mm)^n$ and the fiber of $\mathcal C$ over $((V',\mm'), y)$ is $\D(V',\mm')$. We know that $\mathcal C'_e$ is a coCartesian fibration, classifying the functor $f'\colon \D(V)^n \to \D(V')$, $(M_1, \dots, M_n) \mapsto (M_1 \tensor_V \dots \tensor_V M_n) \tensor_V V'$. It follows easily that $\mathcal C_e$ is a coCartesian fibration, classifying the functor $f_{\varphi,n}\colon \D(V,\mm)^n \to \D(V',\mm')$ obtained by restricting $f'$ to $\D(V,\mm)^n \subset \D(V)^n$ and then composing it with $(V',\mm')^*\colon \D(V') \to \D(V',\mm')$. This shows that $p$ is indeed a locally coCartesian fibration. To show that it is a coCartesian fibration we use \cite[Remark 2.4.2.9]{lurie-higher-topos-theory} which reduces the claim to the following observations (similar to the proofs in \cite[\S4.5.3]{lurie-higher-algebra}):
\begin{itemize}
	\item For every morphism $\varphi\colon (V,\mm) \to (V',\mm')$ of almost setups and all $M, N \in \D(V,\mm)$, the natural morphism $f_{\varphi,1}(M \tensor N) \to f_{\varphi,1}(M) \tensor f_{\varphi,1}(N)$ (induced by $p$) is an isomorphism. To see this, we can factor $\varphi = \varphi_l \comp \varphi_s$, where $\varphi_s\colon (V,\mm) \to (V',\varphi(\mm) V')$ is strict and $\varphi_l\colon (V', \varphi(\mm) V') \to (V', \mm')$ is the ``localization from $\varphi(\mm) V'$ to $\mm'$''. Then $f_{\varphi,1} = f_{\varphi_s,1} \comp f_{\varphi_l,1}$, so it is enough to prove the above claim for $f_{\varphi_s,1}$ and $f_{\varphi_l,1}$ separately. We have $f_{\varphi_s,1} = - \tensor_{V^a} V'^a$ and this functor is clearly symmetric monoidal. Moreover, from $\mm' \subset \varphi(\mm)V'$ it follows that the localilzation functor $(V',\mm')$ factors over the localization functor $(V',\varphi(\mm)V')$; this easily implies that $f_{\varphi_l,1}$ is symmetric monoidal.

	\item For every composition $\varphi = \varphi' \comp \varphi''$ of morphisms in $\AlmSetup$, the induced natural transformation $f_{\varphi,1} \to f_{\varphi',1} \comp f_{\varphi'',1}$ (induced by $p$) is an isomorphism. This follows in a similar fashion as the previous step by factoring $\varphi$, $\varphi'$ and $\varphi''$ as a composition of a strict morphism and a localizing morphism.
\end{itemize}
This finishes the proof that $p$ is a coCartesian fibration, so that we obtain the claimed functor. The claims (i) and (ii) follow immediately from the construction.
\end{proof}

\begin{lemma} \label{rslt:functoriality-of-almost-rings}
There are natural functors
\begin{align*}
	\AlmSetup \to \infcatinf, \qquad &(V,\mm) \mapsto \AssRing_{(V,\mm)},\\
	\AlmSetup \to \infcatinf, \qquad &(V,\mm) \mapsto \Ring_{(V,\mm)}.
\end{align*}
They satisfy the following properties:
\begin{lemenum}
	\item Let $\varphi\colon (V, V) \to (V, \mm)$ be a localizing morphism of almost setups. Then the induced functors $\AssRing_{(V,V)} \to \AssRing_{(V,\mm)}$ and $\Ring_{(V,V)} \to \Ring_{(V,\mm)}$ coincide with the functor $(V,\mm)^*$.

	\item Let $\varphi\colon (V,\mm) \to (V',\mm')$ be a strict morphism of almost setups. Denote $V^a := (V,\mm)^* V$ and $V'^a := (V,\mm)^* V'$. Then $\AssRing_{(V',\mm')} = (\AssRing_{(V,\mm)})_{V'^a/}$ and the induced functor $\AssRing_{(V,\mm)} \to \AssRing_{(V',\mm')}$ coincides with $- \tensor_{V^a} V'^a$. The same holds for $\Ring$ in place of $\AssRing$.
\end{lemenum}
\end{lemma}
\begin{proof}
We only prove the statement for $\Ring$; the argument for $\AssRing$ is analogous. We start with the functor $\AlmSetup \to \infcatinf$, $(V,\mm) \mapsto \Ring_V$ where transition functors are given by base-change. Let $\mathcal C' \to \AlmSetup$ be the biCartesian fibration classified by this functor. Let $\mathcal C \subset \mathcal C'$ be the full subcategory where we only allow objects in $\Ring_{(V,\mm)} \subset \Ring_V$ in the fiber over $(V,\mm)$. Then $p\colon \mathcal C \to \AlmSetup$ is still a Cartesian fibration: This follows from the fact that for every morphism $\varphi\colon (V,\mm) \to (V',\mm')$ of almost setups and every $A \in \Ring_{(V',\mm')} \subset \Ring_{V'}$, the ring $\varphi_*A \in \Ring_V$ lies in $\Ring_{(V,\mm)} \subset \Ring_V$, which in turn follows from \cref{rslt:functoriality-of-almost-modules} by passing to right adjoints of the transition functors. We claim that $p$ is also a coCartesian fibration and thus classifies the desired functor. To see this, by \cite[Corollary 5.2.2.5]{lurie-higher-topos-theory} we need to see that for every morphism of almost setups $\varphi\colon (V,\mm) \to (V',\mm')$ the forgetful functor $\varphi_*\colon \Ring_{(V',\mm')} \to \Ring_{(V,\mm)}$ has a left adjoint. If $\varphi$ is localizing, i.e. $V = V'$, then this left adjoint is given by the almost localization $(V,\mm')^*$ from \cref{rslt:almost-localization-for-rings}. Now suppose that $\varphi$ is strict. Then it follows from \cref{rslt:functoriality-of-almost-modules-along-strict} that a $V'$-algebra $A \in \Ring_{V'}$ lies in $\Ring_{(V',\mm')}$ if and only if $\varphi_* A \in \Ring_V$ lies in $\Ring_{(V,\mm)}$. Thus we have $\Ring_{(V',\mm')} = (\Ring_{(V,\mm)})_{V'^a/}$ and the left adjoint of $\varphi_*$ is given by the functor $(V',\mm')^* (- \tensor_V V') = (- \tensor_{V^a} V'^a)$.
\end{proof}

\begin{proposition} \label{rslt:functoriality-of-almost-modules-over-rings}
Let $\AssRing_{(-)} \to \AlmSetup$ and $\Ring_{(-)} \to \AlmSetup$ denote the coCartesian fibrations classifying the functors from \cref{rslt:functoriality-of-almost-rings}. Then there are natural functors
\begin{align*}
	\AssRing_{(-)} \to \infcatinf, \qquad &A \mapsto \D(A), \\
	\Ring_{(-)} \to \infcatinf^\tensor, \qquad &A \mapsto \D(A).
\end{align*}
If $\varphi\colon (V, V) \to (V,\mm)$ is a localizing morphism of almost setups and $A$ is a (associative) $V$-algebra, then the induced functor $\D(A) \to \D(A^a)$ coincides with $(A,\mm_A)^*$.
\end{proposition}
\begin{proof}
We only construct the second functor; the construction of the first functor is analogous. Let $\mathcal C'^\tensor \to \Ring_{(-)} \cprod \opComm$ be the coCartesian family of $\infty$-operads classifying the functor $A \mapsto \D(A_{**})$, which we can see as the composite functor
\begin{align*}
	\Ring_{(-)} \xto{(-)_{**}} \Ring \xto{(-)^\circ} \CAlg(\D_{\ge0}(\Z)) \to \infcatinf^\tensor
\end{align*}
where the last functor is the one from \cite[Theorem 4.5.3.1]{lurie-higher-algebra}. Let $\mathcal C^\tensor \subset \mathcal C'^\tensor$ be the full subcategory where we only allow objects in $\D(A) \subset \D(A_{**})$ in the fiber over each $A \in \Ring_{(-)}$ (where the inclusion $\D(A) \injto \D(A_{**})$ is the functor $(A_{**},\mm_{A_{**}})_*$. We claim that the induced map $p\colon \mathcal C^\tensor \to \Ring_{(-)} \cprod \opComm$ is a coCartesian family of symmetric monoidal $\infty$-categories and hence defines the desired functor. To see this, we argue as in \cite[\S4.5.3]{lurie-higher-algebra}: It follows from the fact that $(A_{**}, \mm_{A_{**}})_*$ is computed as $(V,\mm)_*$ on underlying $(V,\mm)$-modules (see \cref{rslt:construction-of-right-adjoint-for-modules-over-ass-V-m-alg}) and from \cref{rslt:functoriality-of-almost-modules} that the fiber of $p$ over $\langle 1 \rangle \in \opComm$ is a Cartesian fibration. Since all the left adjoints of all forgetful functors exist (they can be composed of functors of the form $- \tensor_A B$ and $(A,\mm_A)^*$), this fiber is also a coCartesian fibration. Since all the left adjoints are moreover symmetric monoidal in a natural way, $p$ is indeed a coCartesian family of symmetric monoidal $\infty$-categories.
\end{proof}

\begin{definition}
For every morphism $\varphi\colon (V,\mm) \to (V',\mm')$ of almost setups we denote
\begin{align*}
	&\varphi^*\colon \D(V,\mm) \to \D(V',\mm'),\\
	&\varphi^*\colon \AssRing_{(V,\mm)} \to \AssRing_{(V',\mm')},\\
	&\varphi^*\colon \Ring_{(V,\mm)} \to \Ring_{(V',\mm')},
\end{align*}
the functors constructed in \cref{rslt:functoriality-of-almost-modules,rslt:functoriality-of-almost-rings}. Given any (associative) $(V,\mm)$-algebra $A$ we similarly denote
\begin{align*}
	\varphi^*\colon \D(A) \to \D(\varphi^* A)
\end{align*}
the functor constructed in \cref{rslt:functoriality-of-almost-modules-over-rings}. It admits a right adjoint, which we denote by
\begin{align*}
	\varphi_*\colon \D(\varphi^* A) \to \D(A).
\end{align*}
If $\varphi$ is strict then $\varphi_*$ is a forgetful functor. If $\varphi$ is localizing then $\varphi_*$ is of the form $(A,\mm_A)_*$ and $\varphi^*$ is of the form $(A,\mm_A)^*$.
\end{definition}

Given an almost setup $(V,\mm)$ and a $(V,\mm)$-algebra $A^a$, one frequently finds themselves in the position that an object $M \in \D(A^a)$ satisfies a given property only up to some $\varepsilon \in \mm$, i.e. it can be ``$\varepsilon$-approximated'' by some $M_\varepsilon \in \D(A)$ having the desired property. The following definition formalizes this phenomenon.

\begin{definition} \label{def:approx-epsilon}
Let $(V,\mm)$ be an almost setup and $A^a$ an associative $(V,\mm)$-algebra. Let $M, N \in \D(A^a)$ and $\varepsilon \in \mm$.
\begin{defenum}
	\item We write $M \approx_\varepsilon N$ if there are maps $f_\varepsilon\colon M \to N$ and $g_\varepsilon\colon N \to M$ such that $f_\varepsilon g_\varepsilon$ and $g_\varepsilon f_\varepsilon$ are isomorphic to the multiplication-by-$\varepsilon$ map. We write $M \approx N$ if $M \approx_\varepsilon N$ holds for all $\varepsilon \in \mm$.

	\item We say that $M$ is an $\varepsilon$-retract of $N$ if there are maps $f_\varepsilon\colon M \to N$ and $g_\varepsilon\colon N \to M$ such that $g_\varepsilon f_\varepsilon \isom \varepsilon \id$.
\end{defenum}
\end{definition}

Throughout the thesis we often encounter $\infty$-categories $\mathcal C$ such that for every two objects $X, Y \in \mathcal C$ the anima $\Hom(X, Y)$ can naturally be upgraded to an almost $V$-module (for some almost setup $(V,\mm)$). This phenomenon is best modeled by saying that $\mathcal C$ is \emph{enriched} over $\D_{\ge0}(V^a)_\omega$ (see \cref{sec:infcat.enriched} for a quick introduction to enriched $\infty$-categories; in the following we only need a vague understanding of that concept). Concretely, we introduce the following notation:

\begin{definition} \label{def:V-a-enriched-categories}
Let $(V, \mm)$ be an almost setup and $\mathcal C$ an $\infty$-category. We say that $\mathcal C$ is \emph{enriched over $V^a$} if it is enriched over $\D_{\ge0}(V^a)_\omega$. In this case we denote the enriched homomorphism functor by
\begin{align*}
	\Hom(-, -)^a\colon \mathcal C^\opp \cprod \mathcal C \to \D_{\ge0}(V^a)_\omega.
\end{align*}
Note that for every $\varepsilon \in \mm$ and every morphism $f\colon X \to Y$ in $\mathcal C$ there is a canonical morphism $\varepsilon f\colon X \to Y$ in $\mathcal C$. This allows us in particular to define $\varepsilon$-retracts and the relations $\approx_\varepsilon$ and $\approx$ for objects in $\mathcal C$ as in \cref{def:approx-epsilon}.
\end{definition}

One problem we often encounter in the world of almost mathematics is that $\infty$-categories which feel like being compactly generated are not actually compactly generated. For example this already holds for the $\infty$-category $\D(V^a)_\omega$ of discrete almost $V$-modules (for a given almost setup $(V,\mm)$): The object $V^a \in \D(V^a)_\omega$ is in general not compact because $\Hom(V^a, M^a) = \Hom(\widetilde\mm, M)$ for all $M \in \D(V)_\omega$, and $\widetilde\mm$ is usually not compact in $\D(V)_\omega$. Still, $V^a$ is \emph{almost compact}, meaning that $\Hom(V^a, -)$ commutes with colimits up to almost isomorphism. This intuition is made more precise in the following definition.

\begin{definition} \label{def:almost-compact-object}
Let $(V,\mm)$ be an almost setup and $\mathcal C$ a $V^a$-enriched $\infty$-category. An object $P \in \mathcal C$ is called \emph{almost compact} if the functor $\Hom(P, -)^a\colon \mathcal C \to \D_{\ge0}(V^a)_\omega$ preserves all filtered colimits. We denote $\mathcal C^{a\omega} \subset \mathcal C$ the full subcategory spanned by the almost compact objects.
\end{definition}

We get the following analog of \cref{rslt:compact-generators-implies-compactly-generated}, showing that almost compact ``generators'' behave as nicely as one may hope.

\begin{proposition} \label{rslt:almost-compact-generators-implies-generated}
Let $(V, \mm)$ be an almost setup, $\mathcal C$ a $V^a$-enriched $\infty$-category which has all small colimits, and $(P_i)_{i\in I}$ a small family of almost compact objects of $\mathcal C$ such that the family of functors $\Hom(P_i, -)^a\colon \mathcal C \to \D_{\ge0}(V^a)_\omega$ is conservative. Let $\mathcal C_0 \subset \mathcal C$ be the smallest full subcategory which contains all $P_i$ and is stable under finite colimits and retracts. Then:
\begin{propenum}
	\item \label{rslt:almost-compact-generators-fully-faithful-into-Ind} The natural functor $\Ind(\mathcal C_0) \to \mathcal C$ admits a fully faithful right adjoint. In particular $\mathcal C$ is presentable and every object in $\mathcal C$ is a filtered colimit of objects in $\mathcal C_0$.

	\item \label{rslt:characterization-of-almost-compact-objects-in-general-C} An object $P \in \mathcal C$ is almost compact if and only if for every $\varepsilon \in \mm$ there is some $P_\varepsilon \in \mathcal C_0$ such that $P$ is an $\varepsilon$-retract of $P_\varepsilon$.
\end{propenum}
\end{proposition}
\begin{proof}
We start with the proof of (i). By \cite[Proposition 5.3.5.10]{lurie-higher-topos-theory} there is a natural functor $F\colon \Ind(\mathcal C_0) \to \mathcal C$ which preserves all filtered colimits. By \cite[Lemma 5.4.2.4]{lurie-higher-topos-theory}, the compact objects of $\Ind(\mathcal C_0)$ are precisely $\mathcal C_0$, hence by \cite[Proposition 5.5.1.9]{lurie-higher-topos-theory} $F$ preserves all small colimits. Then by \cite[Remark 5.5.2.10]{lurie-higher-topos-theory} $F$ admits a right-adjoint $G\colon \mathcal C \to \Ind(\mathcal C_0)$. We will show that $G$ is fully faithful, i.e. that for every $M \in \mathcal C$ the natural map $F(G(M)) \isoto M$ is an isomorphism.

For every $P \in \mathcal C_0$ the functor $\Hom(P, -)^a\colon \mathcal C_0 \to \D_{\ge0}(V^a)_\omega$ extends uniquely to a functor $\Hom(P, -)^a\colon \Ind(\mathcal C_0) \to \D_{\ge0}(V^a)_\omega$ which preserves filtered colimits. It agrees with the functor $\Hom(P, F(-))^a\colon \Ind(\mathcal C_0) \to \D_{\ge0}(V^a)_\omega$, because both functors preserve filtered colimits and agree on $\mathcal C_0$. For every $M \in \mathcal C$ there is a natural map
\begin{align}
	\Hom(P, G(M))^a = \Hom(P, F(G(M)))^a \to \Hom(P, M)^a, \label{eq:almost-compact-generators-eq1}
\end{align}
which we claim to be an isomorphism. This does not follow formally from the adjunction of $F$ and $G$ because $\Hom(P, -)^a$ does not necessarily define an enrichment of $\Ind(\mathcal C_0)$. We can circumvent this issue as follows. Consider the functor $\Hom(P, -)_V\colon \mathcal C_0 \to \D_{\ge0}(V)_\omega$ which is computed as $\Hom(P, -)_V := \tau_{\ge0}(\Hom(P, -)^a)_*$. As before, this functor extends uniquely to a functor $\Hom(P, -)_V\colon \Ind(\mathcal C_0) \to \D_{\ge0}(V)_\omega$ which preserves all filtered colimits. By uniqueness we have a natural isomorphism $\Hom(P, -)^a = (\Hom(P, -)_V)^a$ of functors $\Ind(\mathcal C_0) \to \D_{\ge0}(V)_\omega$. Moreover there is a natural morphism $\Hom(P, -)_V \to \Hom(P, F(-))_V$ of functors $\Ind(\mathcal C_0) \to \D_{\ge0}(V)_\omega$ coming from the fact that both functors agree on $\mathcal C_0$ and $\Hom(P, -)_V$ is the left Kan extension of its restriction to $\mathcal C_0$. We thus get a morphism
\begin{align}
	\Hom(P, G(M))_V \to \Hom(P, F(G(M)))_V \to \Hom(P, M)_V \label{eq:almost-compact-generators-eq2}
\end{align}
of objects in $\D_{\ge0}(V)_\omega$, for every $M \in \mathcal C$. To show that this is an isomorphism, it is enough to do so after applying $\Hom(V, -)$. But on both sides we have
\begin{align*}
	\Hom(V, \Hom(P, -)_V) = \Hom(P, -)
\end{align*}
(For $\Hom(P, -)_V\colon \mathcal C \to \D_{\ge0}(V)_\omega$ this is clear and for $\Hom(P, -)_V\colon \Ind(\mathcal C_0) \to \D_{\ge0}(V)_\omega$ this follows because $\Hom(V, -)$ preserves filtered colimits.) Thus it follows from the adjunction of $F$ and $G$ that the composition of morphisms in \cref{eq:almost-compact-generators-eq2} is an isomorphism. By applying $(-)^a$ we deduce that the morphism in \cref{eq:almost-compact-generators-eq1} is an isomorphism. But by assumption the family functors $\Hom(P, -)^a$, for varying $P \in \mathcal C_0$, is conservative, which implies that $F(G(M)) \isoto M$ is indeed an isomorphism. This finishes the proof of (i).

We now prove (ii), so let $P \in \mathcal C$ be given. First assume that $P$ is almost compact. By (i) we can write $P$ as a filtered colimit $P = \varinjlim_i P_i$ of objects $P_i \in \mathcal C_0$. Then by assumption on $P$ the map
\begin{align*}
	\varphi\colon \varinjlim_i \pi_0 \Hom(P, P_i)_V \isoto \pi_0 \Hom(P, P)_V
\end{align*}
is an almost isomorphism of classical $V$-modules, whose underlying set is the set of homomorphisms in $\mathcal C$. In particular, for any given $\varepsilon \in \mm$, $\varepsilon \id_P$ lies in the image of $\varphi$, i.e. there is an index $i$ and a map $f_i\colon P \to P_i$ such that the composition of $f_i$ with the natural map $g_i\colon P_i \to P$ is multipliciation by $\varepsilon$. This proves the ``only if'' direction of (ii).

It remains to prove the ``if'' direction of (ii), so assume that $P$ satisfies the condition involving the $P_\varepsilon$ and let $(M_i)_{i\in I}$ be a filtered system of objects in $\mathcal C$. Fix any $\varepsilon \in \mm$ and choose $P_\varepsilon \in \mathcal C_0$ such that $P$ is an $\varepsilon$-retract of $P_\varepsilon$, i.e. there are maps $f_\varepsilon\colon P \to P_\varepsilon$, $g_\varepsilon\colon P_\varepsilon \to P$ such that $g_\varepsilon f_\varepsilon \isom \varepsilon \id$. Then $P_\varepsilon$ is almost compact, so we obtain the following diagram of objects in $\D_{\ge0}(V^a)_\omega$:
\begin{center}\begin{tikzcd}
	\varinjlim_i \Hom(P, M_i)^a \arrow[r] \arrow[d,shift left] & \Hom(P, \varinjlim_i M_i)^a \arrow[d,shift left]\\
	\varinjlim_i \Hom(P_\varepsilon, M_i)^a \arrow[r,"\sim"] \arrow[u,shift left] & \Hom(P_\varepsilon, \varinjlim_i M_i)^a \arrow[u,shift left]
\end{tikzcd}\end{center}
This diagram commutes if we restrict the vertical morphisms to either both upwards or both downwards. The composition of either downwards arrow with its corresponding upwards arrow is isomorphic to multiplication by $\varepsilon$. A quick diagram chase shows that both kernel and cokernel on homotopy groups of the upper horizontal map are killed by $\varepsilon$. Since $\varepsilon$ was arbitrary, the upper horizontal map is an (almost) isomorphism, as desired.
\end{proof}

\subsection{Analytic Rings} \label{sec:andesc.anring}

Throughout this subsection, we fix an almost setup $(V,\mm)$. We now introduce the notion of analytic rings, as defined in \cite[\S12]{scholze-analytic-spaces}. The main difference to the reference is that we work more generally with \emph{almost} $V$-modules (instead of honest abelian groups).

Let us start with the basic definitions. We first introduce associative analytic rings.

\begin{definition}
\begin{defenum}
	\item An \emph{uncompleted\footnote{In accordance with the authors of \cite{scholze-analytic-spaces} we renamed ``normalized analytic rings'' into ``complete analytic rings''.} pre-analytic associative ring} over $(V,\mm)$ is a pair $\mathcal A = (\underline{\mathcal A}, \mathcal M_{\mathcal A})$, where
	\begin{itemize}
		\item $\underline{\mathcal A} \in \AssRing_{(V,\mm)}$ is an associative $(V,\mm)$-algebra,

		\item $\mathcal M_{\mathcal A}$ is a functor
		\begin{align*}
			\mathcal M_{\mathcal A}\colon \ExtrDisc \to \D_{\ge0}(\underline{\mathcal A}), \qquad S \mapsto \mathcal M_{\mathcal A}[S]
		\end{align*}
		from the category of extremally disconnected sets to connective $\underline{\mathcal A}$-modules which maps finite coproducts to finite direct sums, together with a natural transformation $\underline{\mathcal A}[S] \to \mathcal M_{\mathcal A}[S]$ of functors $\ExtrDisc \to \D_{\ge0}(\underline{\mathcal A})$.
	\end{itemize}
	We will usually denote $\mathcal A[S] := \mathcal M_{\mathcal A}[S]$, making the functor $\mathcal M_{\mathcal A}$ obsolete in the notation.

	\item An \emph{uncompleted analytic associative ring} over $(V, \mm)$ is an uncompleted pre-analytic associative ring $\mathcal A$ over $(V, \mm)$ such that for every object $M \in \D(\underline{\mathcal A})$ which is a sifted colimit of the objects $\mathcal A[S]$, the natural map
	\begin{align*}
		\IHom_{\underline{\mathcal A}}(\mathcal A[S], M) \isoto \IHom_{\underline{\mathcal A}}(\underline{\mathcal A}[S], M)
	\end{align*}
	is an isomorphism in $\D(V^a)$. Here $\IHom_{\underline{\mathcal A}}$ denotes the $\D(V^a)$-enriched hom in $\D(\underline{\mathcal A})$ (note that $\D(\underline{\mathcal A})$ is tensored over $\D(V^a)$ and thus easily seen to be enriched; see \cref{sec:infcat.enriched} for more on enriched $\infty$-categories).

	\item \label{def:D-of-ass-analytic-ring} Let $\mathcal A$ be an uncompleted analytic associative ring over $(V, \mm)$. We denote
	\begin{align*}
		\D(\mathcal A) \subset \D(\underline{\mathcal A})
	\end{align*}
	the full subcategory spanned by all $M \in \D(\underline{\mathcal A})$ such that for every extremally disconnected set $S$ the map $\IHom_{\underline{\mathcal A}}(\mathcal A[S], M) \isoto \IHom_{\underline{\mathcal A}}(\underline{\mathcal A}[S], M)$ is an isomorphism in $\D(V^a)$.

	\item A \emph{morphism of uncompleted analytic associative rings} $\mathcal A \to \mathcal B$ is a morphism $\underline{\mathcal A} \to \underline{\mathcal B}$ of the underlying $(V,\mm)$-algebras such that the forgetful functor $\D(\mathcal B) \to \D(\underline{\mathcal A})$ takes image in $\D(\mathcal A)$. The $\infty$-category of uncompleted analytic associative rings over $(V, \mm)$ is denoted
	\begin{align*}
		\AnAssRing^\dagger_{(V,\mm)}.
	\end{align*}
\end{defenum}
\end{definition}

\begin{proposition} \label{rslt:properties-and-functoriality-of-D-of-analytic-ass-ring}
\begin{propenum}
	\item \label{rslt:properties-of-D-of-analytic-ass-ring} Let $\mathcal A$ be an uncompleted analytic associative ring over $(V,\mm)$. The $\infty$-category $\D(\mathcal A)$ is stable under all limits and colimits in $\D(\underline{\mathcal A})$ and generated under colimits by the objects $\mathcal A[S]$. The $t$-structure on $\D(\underline{\mathcal A})$ restricts to a $t$-structure on $\D(\mathcal A)$. The inclusion $\D(\mathcal A) \injto \D(\underline{\mathcal A})$ admits a left adjoint
	\begin{align*}
		- \tensor_{\underline{\mathcal A}} \mathcal A\colon \D(\underline{\mathcal A}) \to \D(\mathcal A)
	\end{align*}
	sending $\underline{\mathcal A}[S]$ to $\mathcal A[S]$.

	\item \label{rslt:functorialiy-of-D-A-for-ass-ring} For every morphism $\mathcal A \to \mathcal B$ of uncompleted analytic associative rings over $(V,\mm)$, the forgetful functor $\D(\mathcal B) \to \D(\mathcal A)$ admits a left adjoint
	\begin{align*}
		- \tensor_{\mathcal A} \mathcal B\colon \D(\mathcal A) \to \D(\mathcal B).
	\end{align*}
	This defines a functor from $\AnAssRing^\dagger_{(V,\mm)}$ to $\infty$-categories, sending $\mathcal A$ to $\D(\mathcal A)$ and $\mathcal A \to \mathcal B$ to $- \tensor_{\mathcal A} \mathcal B$.
\end{propenum}
\end{proposition}
\begin{proof}
We first prove (i), so let $\mathcal A$ be given. We will implicitly make use of \cref{rslt:general-properties-of-almost-localization-over-ass-alg}. By \cref{rslt:every-ass-almost-algebra-is-A-a-of-actual-algebra} we can pick an associative $V$-algebra $A'$ such that $\underline{\mathcal A} = A'^a$. Let $\mathcal C' \subset \D_{\ge0}(A')$ be the full subcategory spanned by the objects $A'[S]$ for extremally disconnected sets $S$. There is a natural functor $F_0\colon \mathcal C' \to \D_{\ge0}(\underline{\mathcal A})$ sending $A'[S]$ to $\mathcal A[S]$, which can be constructed as follows: Let $p'\colon X' \to \Delta^1$ denote the coCartesian fibration classifying the functor $\D(A') \to \D(\underline{\mathcal A})$ and let $X \subset X'$ be the full subcategory where we only allow objects isomorphic to some $A'[S]$ in the fiber over $0 \in \Delta^1$ and objects which are isomorphic to some $\mathcal A[S]$ in the fiber over $1 \in \Delta^1$. We claim that the map $p\colon X \to \Delta^1$ is still a coCartesian fibration. It is clearly a categorical fibration, so we only need to show the following: Given an extremally disconnected set $S$, there is a $p$-coCartesian edge $A'[S] \to M$ over the non-degenerate edge in $\Delta^1$. By definition of $X'$ there is a $p'$-coCartesian edge $A'[S] \to \underline{\mathcal A}[S]$ in $X'$ (where $\underline{\mathcal A}[S]$ lies in the fiber over $1 \in \Delta^1$) and by definition of uncompleted analytic associative rings there is also a natural morphism $\underline{\mathcal A}[S] \to \mathcal A[S]$. Composing these morphisms we get a morphism $f\colon A'[S] \to \mathcal A[S]$ in $X$ which we claim to be $p$-coCartesian. By \cite[Proposition 2.4.4.3]{lurie-higher-topos-theory} this boils down to the claim that for every extremally disconnected set $S'$ the natural map $\Hom_X(\mathcal A[S], \mathcal A[S']) \to \Hom_X(A'[S], \mathcal A[S'])$ induced by $f$ is an isomorphism of anima. Since the functor $\D(A') \to \D(\underline{\mathcal A})$ maps $A'[S] \mapsto \underline{\mathcal A}[S]$ we have $\Hom_X(A'[S], \mathcal A[S']) = \Hom(\underline{\mathcal A}[S], \mathcal A[S])$. Thus the claimed isomorphism of anima follows directly from the definition of uncompleted analytic associative rings. Now the desired functor $F_0\colon \mathcal C' \to \D(\mathcal A)$ is induced by the coCartesian fibration $p$. In a similar vein we can construct a natural transformation $(A',\mm_{A'})^* \to F_0$ by viewing it as a functor $\mathcal C' \to \Fun(\Delta^1, \D(\underline{\mathcal A}))$; we leave the details to the reader.

Since $\D_{\ge0}(A')$ is freely generated by $\mathcal C'$ under sifted colimits, the functor $F_0$ extends uniquely to a functor $F_0\colon \D_{\ge0}(A') \to \D_{\ge0}(\underline{\mathcal A})$ which preserves all small colimits (cf. \cite[Proposition 5.5.8.15]{lurie-higher-topos-theory}). Let us denote
\begin{align*}
	F := F_0 \comp (A',\mm_{A'})_!\colon \D_{\ge0}(\underline{\mathcal A}) \to \D_{\ge0}(\underline{\mathcal A}).
\end{align*}
This functor still preserves all colimits and inherits the natural transformation $\id \to F$. Let $\D_{\ge0}(\mathcal A) \subset \D_{\ge0}(\underline{\mathcal A})$ be the full subcategory spanned by those objects $M \in \D_{\ge0}(\underline{\mathcal A})$ such that for all extremally disconnected sets $S$ the natural map
\begin{align*}
	\tau_{\ge0} \IHom_{\underline{\mathcal A}}(\mathcal A[S], M) \isoto \tau_{\ge0} \IHom_{\underline{\mathcal A}}(\underline{\mathcal A[S]}, M)
\end{align*}
is an isomorphism in $\D_{\ge0}(V^a)$ (it will turn out that $\D_{\ge0}(\mathcal A) = \D(\mathcal A) \isect \D_{\ge0}(\underline{\mathcal A})$, but this is a priori not clear). The essential image of $F$ is generated under sifted colimits by the objects $\mathcal A[S]$, which by the definition of uncompleted analytic associative rings is contained in $\D_{\ge0}(\mathcal A)$. Given any $M \in \D_{\ge0}(\underline{\mathcal A})$ and $N \in \D_{\ge0}(\mathcal A)$ the natural map $\Hom(F(M), N) \isoto \Hom(M, N)$ is an isomorphism (we can pull out sifted colimits in $M$ to assume that $M = \underline{\mathcal A}[S]$ for some extremally disconnected $S$). By \cite[Proposition 5.2.2.8]{lurie-higher-topos-theory} $F$ is left adjoint to the inclusion $\D_{\ge0}(\mathcal A) \injto \D_{\ge0}(\underline{\mathcal A})$. It follows that $\D_{\ge0}(\mathcal A)$ is stable under all colimits. It is clearly stable under limits as well.

Next we claim that $\D_{\ge0}(\mathcal A) \subset \D(\mathcal A)$. Namely, if $M \in \D_{\ge0}(\mathcal A)$ then $M[n] \in \D_{\ge0}(\mathcal A)$ for all $n \ge 0$, as $\D_{\ge0}(\mathcal A)$ is stable under limits and colimits. But $\tau_{\ge0} \IHom_{\underline{\mathcal A}}(\mathcal A[S], M[n]) = \tau_{\ge-n} \IHom_{\underline{\mathcal A}}(\mathcal A[S], M)$ (and similarly for $\underline{\mathcal A[S]}$ in place of $\mathcal A[S]$), which implies $M \in \D(\mathcal A)$.

We can now show that $\D(\mathcal A) \subset \D(\underline{\mathcal A})$ is stable under all colimits: It is clearly stable under finite colimits, so it is enough to show that it is stable under filtered colimits. But note that for any $M \in \D(\mathcal A)$ we have $\tau_{\ge0}M \in \D_{\ge0}(\mathcal A)$ and $\tau_{\ge0}\colon \D(\underline{\mathcal A}) \to \D_{\ge0}(\underline{\mathcal A})$ commutes with filtered colimits. Applying this to all $M[n]$ and using that $\D_{\ge0}(\mathcal A) \subset \D_{\ge0}(\underline{\mathcal A})$ is stable under all colimits, we deduce that $\D(\mathcal A)$ is stable under filtered colimits as desired.

Clearly $\D(\mathcal A)$ is stable under limits. To construct the functor $- \tensor_{\underline{\mathcal A}} \mathcal A$, note that when restricted to $\D_{\ge0}(\underline{\mathcal A})$ this is just the functor $F$. Now $\D(\underline{\mathcal A})$ is just the $\infty$-category of spectrum objects of $\D_{\ge0}(\underline{\mathcal A})$ which allows us to extend $F$ to a functor $F\colon \D(\underline{\mathcal A}) \to \D(\underline{\mathcal A})$ (see \cite[Corollary 1.4.4.5]{lurie-higher-algebra}; the reference only works with presentable $\infty$-categories, but one can reduce to that case by choosing a cofinal class of strong limit cardinals $\kappa' \ge \kappa$ such that $F$ maps $\D_{\ge0}(\underline{\mathcal A})_\kappa$ to $\D(\underline{\mathcal A})_{\kappa'}$, i.e. $\kappa'$ is such that $\mathcal A[S] \in \D(\underline{\mathcal A})_{\kappa'}$ for all $\kappa$-small extremally disconnected sets $S$). Then $F$ preserves all colimits, so writing any $M \in \D(\underline{\mathcal A})$ as $M = \varinjlim_n \tau_{\ge-n} M$ we deduce that $F$ factors over $\D(\mathcal A)$. The natural transformation $\id \to F$ of functors extends to $\D(\underline{\mathcal A})$ and by the same argument as in the connective case we see that $F$ is indeed left adjoint to the inclusion $\D(\mathcal A) \injto \D(\underline{\mathcal A})$ and hence defines $- \tensor_{\underline{\mathcal A}} \mathcal A$.

Part (ii) is formal: The functor $- \tensor_{\mathcal A} \mathcal B$ can be defined as $M \mapsto (M \tensor_{\underline{\mathcal A}} \underline{\mathcal B}) \tensor_{\underline{\mathcal B}} \mathcal B$, which is indeed left adjoint to the forgetful functor. The functoriality of this construction follows by looking at the associated (co)Cartesian fibrations and using \cite[Corollary 5.2.2.5]{lurie-higher-topos-theory}); see the proof of \cref{rslt:functoriality-of-D-A} below where do this more explicitly.
\end{proof}

It is often convenient to have a presentable version of $\D(\mathcal A)$. Similar to the $\infty$-categories $\D(A)_\kappa$ for $(V,\mm)$-algebras $A$, we can make the following definition.

\begin{definition} \label{rslt:kappa-modules-over-analytic-ring}
Let $\mathcal A$ be an uncompleted analytic associative ring over $(V,\mm)$. Let $\kappa$ be a strong limit cardinal such that $\underline{\mathcal A} \in \D(V,\mm)_\kappa$. Then we define
\begin{align*}
	\D(\mathcal A)_\kappa \subset \D(\mathcal A)
\end{align*}
to be the smallest full subcategory which is closed under small colimits and contains all objects $\mathcal A[S]$ for $\kappa$-small extremally disconnected sets $S$.
\end{definition}

\begin{proposition}
Let $\tau$ be a regular cardinal such that $(V,\mm)$ is $\tau$-compact.
\begin{propenum}
	\item Let $\mathcal A$ be an uncompleted analytic associative ring over $(V,\mm)$ and let $\kappa$ be a strong limit cardinal such that $\underline{\mathcal A} \in \D(V,\mm)_\kappa$. Then $\D(\mathcal A)_\kappa$ is a stably presentable $\infty$-category which is stable under all colimits in $\D(\mathcal A)$. It is $\tau$-compactly generated by the $\tau$-compact objects $\mathcal A[S]$ for $\kappa$-small extremally disconnected sets $S$.
	\item Let $\mathcal A \to \mathcal B$ be a morphism of uncompleted analytic associative rings over $(V,\mm)$. Then for every strong limit cardinal $\kappa$ such that $\underline{\mathcal A}, \underline{\mathcal B} \in \D(V,\mm)_\kappa$ the functor $- \tensor_{\mathcal A} \mathcal B$ restricts to a symmetric monoidal functor $\D(\mathcal A)_\kappa \to \D(\mathcal B)_\kappa$.
\end{propenum}
\end{proposition}
\begin{proof}
We first prove (i), so let $\mathcal A$ and $\kappa$ be given. Then all $V^a[S] \in \D(V,\mm)$ are $\tau$-compact and thus all $\mathcal A[S] \in \D(\mathcal A)$ are $\tau$-compact. Let $\mathcal C \subset \D(\mathcal A)$ be the full subcategory spanned by all iterated $\tau$-small colimits of the objects $\mathcal A[S]$ for $\kappa$-small extremally disconnected sets $S$. Then $\mathcal C$ is small (at most $\tau$ steps are necessary in its construction) and consists of $\tau$-compact objects in $\D(\mathcal A)$. By \cite[Proposition 5.3.5.11]{lurie-higher-topos-theory} there is a fully faithful functor $F\colon \Ind_\tau(\mathcal C) \injto \D(\mathcal A)$. Since $\restrict F{\mathcal C}$ preserves $\tau$-small colimits, $F$ preserves all small colimits (analogous to the proof of \cite[Proposition 5.5.8.15.(3)]{lurie-higher-topos-theory}). Consequently the image of $F$ is stable under all small colimits. On the other hand this image is contained in $\D(\mathcal A)_\kappa$ and hence must be equal to $\D(\mathcal A)_\kappa$. It follows that $F$ induces an equivalence $\Ind_\tau(\mathcal C) \isoto \D(\mathcal A)_\kappa$, hence $\D(\mathcal A)_\kappa$ is presentable. It is clear that $\D(\mathcal A)_\kappa$ is a stable $\infty$-category and stable under all colimits in $\D(\mathcal A)$.

Part (ii) follows directly from the fact that $\mathcal A[S] \tensor_{\mathcal A} \mathcal B = \mathcal B[S]$ for every extremally disconnected set $S$.
\end{proof}

Given an uncompleted analytic associative ring $\mathcal A$ over $(V,\mm)$ and an $\underline{\mathcal A}$-algebra $B$, there is a natural induced uncompleted analytic ring structure $\mathcal B$ on $B$ such that $\D(\mathcal B) \subset \D(B)$ consists precisely of those $B$-modules whose underlying $\underline{\mathcal A}$-module lies in $\D(\mathcal A)$. More precisely, we have the following result:

\begin{lemma} \label{rslt:existence-of-induced-unnormalized-analytic-ass-ring-structure}
Let $\mathcal A$ be an uncompleted analytic associative ring over $(V,\mm)$. The forgetful functor $\AnAssRing^\dagger_{\mathcal A/} \to \AssRing_{\underline{\mathcal A}}$ admits a fully faithful left adjoint
\begin{align*}
	L\colon \AssRing_{\underline{\mathcal A}} \to \AnAssRing^\dagger_{\mathcal A/}.
\end{align*}
\end{lemma}
\begin{proof}
This is stated in \cite[Proposition 12.8]{scholze-analytic-spaces}, but as pointed out by Longke Tang and Yixiao Li the argument is wrong (it does work if $\underline{\mathcal A}$ comes from an $\mathbb E_\infty$-ring, though). In the following we sketch their fixed version of the proof.

Given an associative $\underline{\mathcal A}$-algebra $B$ we define $L(B)$ to be the following uncompleted pre-analytic associative ring with $\underline{L(B)} := B$: Given any $M \in \D(B)$ we denote by $F(M)$ the pushout of $B \tensor_{\underline{\mathcal A}} (\mathcal A \tensor_{\underline{\mathcal A}} M) \from B \tensor_{\underline{\mathcal A}} M \to M$. Then $F$ defines an endofunctor of $\D(B)$ and comes with a natural transformation $\id \to F$. For every $n \ge 0$ we have a natural morphism $F^n = \id \comp F^n \to F \comp F^n = F^{n+1}$ and we define $G := \varinjlim_n F^n$. We now claim that we can take $L(B)[S] := G(B[S])$. Namely, one first checks that $G(M) = \varinjlim_n (\mathcal A \tensor_{\underline{\mathcal A}} F^n(M))$ on the underlying $\underline{\mathcal A}$-module, so that $G(M) \in \D(\mathcal A)$ for all $M \in \D(B)$. In particular the functor $G$ defines a functor $G\colon \D(B) \to \D$, where $\D \subset \D(B)$ is the full subcategory spanned by those $B$-modules whose image in $\D(\underline{\mathcal A})$ lies in $\D(\mathcal A)$. One then checks that $G$ is in fact a left adjoint to the inclusion $\D \injto \D(B)$, from which it follows that the above defined $L(B)$ does indeed define an uncompleted associative analytic ring and we have $\D(L(B)) = \D$. This finishes the construction of $L(B)$.

It is clear that the assignment $B \mapsto L(B)$ defines a functor $L\colon \AssRing_{\underline{\mathcal A}} \to \AnAssRing^\dagger_{\mathcal A/}$ and it is also clear that $L$ is fully faithful. Clearly, $L$ followed by the forgetful functor $\mathrm{For}$ is the identity, which gives a natural transformation $\id \to \mathrm{For} \comp L$ of endofunctors of $\AssRing_{\underline{\mathcal A}}$. To show that $L$ is left adjoint to $\mathrm{For}$, by \cite[Proposition 5.2.2.8]{lurie-higher-topos-theory} it is enough to show that for all $B \in \AssRing_{\underline{\mathcal A}}$ and all uncompleted analytic associative rings $\mathcal B$ over $\mathcal A$, the natural map
\begin{align*}
	\Hom(L(B), \mathcal B) \isoto \Hom(B, \mathrm{For}(\mathcal B))
\end{align*}
is an isomorphism. But this is just saying that every map $B \to \underline{\mathcal B}$ of $\underline{\mathcal A}$-algebras is automatically a map of analytic rings, i.e. each $\mathcal B[S]$ lands in $\D(L(B))$ via the forgetful functor $\D(\mathcal B) \to \D(B)$. This is clear by the above description of $\D(L(B))$.
\end{proof}

We will now introduce (commutative) analytic rings. Having already studied analytic associative rings, the definitions are mostly straightforward:

\begin{definition}
\begin{defenum}
	\item \label{def:unnormalized-analytic-ring} An \emph{uncompleted analytic ring} over $(V,\mm)$ is a pair $\mathcal A = (\underline{\mathcal A}, \mathcal M_{\mathcal A})$, where $\underline{\mathcal A}$ is a $(V,\mm)$-algebra and $\mathcal M_{\mathcal A}$ induces an uncompleted analytic associative ring structure on the underlying associative $(V,\mm)$-algebra of $\underline{\mathcal A}$, such that for every prime $p$ the Frobenius $\varphi_p\colon \underline{\mathcal A} \to \underline{\mathcal A}/p$ induces a map of uncompleted analytic associative rings $\mathcal A \to \mathcal A/p$ (where we equip $\mathcal A/p$ with the induced uncompleted analytic ring structure from \cref{rslt:existence-of-induced-unnormalized-analytic-ass-ring-structure} along the projection $\underline{\mathcal A} \to \underline{\mathcal A}/p$).

	\item Let $\mathcal A$ be an uncompleted analytic ring over $(V,\mm)$. The $\infty$-category $\D(\mathcal A)$ of \emph{$\mathcal A$-modules} is the $\infty$-category defined in \cref{def:D-of-ass-analytic-ring} associated to the underlying uncompleted analytic associative ring.

	\item A \emph{morphism of uncompleted analytic rings} $\mathcal A \to \mathcal B$ is a morphism of $(V,\mm)$-algebras $\underline{\mathcal A} \to \underline{\mathcal B}$ which induces a morphism of uncompleted analytic associative rings. The $\infty$-category of uncompleted analytic rings over $(V,\mm)$ is denoted
	\begin{align*}
		\AnRing^\dagger_{(V,\mm)}.
	\end{align*}
\end{defenum}
\end{definition}

\begin{remark} \label{rmk:Frobenius-condition-in-analytic-ring-def}
The condition on the Frobenius in \cref{def:unnormalized-analytic-ring} seems to come out of nowhere, but is needed to define the completion functor below (see \cref{rslt:existence-of-initial-ring-over-unnormalized-analytic-ring}). This observation was made by Scholze, cf. the remarks following \cite[Definition 12.10]{scholze-analytic-spaces}. The condition is automatic in a number of cases, e.g. when all $\mathcal A[S]$ are $m$-truncated for a fixed $m \ge 0$ (see \cite[Proposition 12.24]{scholze-analytic-spaces}).
\end{remark}

Given an uncompleted analytic ring $\mathcal A$ over $(V,\mm)$, the results of \cref{rslt:properties-and-functoriality-of-D-of-analytic-ass-ring} can be extended to also include a symmetric monoidal structure on $\D(\mathcal A)$:

\begin{proposition}
\begin{propenum}
	\item \label{rslt:properties-of-D-of-analytic-ring} Let $\mathcal A$ be an uncompleted analytic ring over $(V,\mm)$. Then all the properties of \cref{rslt:properties-of-D-of-analytic-ass-ring} apply to $\D(\mathcal A)$. Moreover, there is a natural closed symmetric monoidal structure on $\D(\mathcal A)$ making $- \tensor_{\underline{\mathcal A}} \mathcal A$ symmetric monoidal. For every $M \in \D(\underline{\mathcal A})$ and $N \in \D(\mathcal A)$ we have
	\begin{align*}
		\IHom_{\underline{\mathcal A}}(M, N) \in \D(\mathcal A).
	\end{align*}
	In particular, internal hom's in $\D(\mathcal A)$ agree with those in $\D(\underline{\mathcal A})$.

	\item \label{rslt:functoriality-of-D-A} For every morphism $\mathcal A \to \mathcal B$ of uncompleted analytic rings over $(V,\mm)$ the forgetful functor $\D(\mathcal B) \to \D(\mathcal A)$ admits a symmetric monoidal left adjoint
	\begin{align*}
		- \tensor_{\mathcal A} \mathcal B\colon \D(\mathcal A) \to \D(\mathcal B).
	\end{align*}
	This defines a functor from $\AnRing^\dagger_{(V,\mm)}$ to symmetric monoidal $\infty$-categories, sending $\mathcal A$ to $\D(\mathcal A)$ and $\mathcal A \to \mathcal B$ to $- \tensor_{\mathcal A} \mathcal B$.
\end{propenum}
\end{proposition}
\begin{proof}
To construct the symmetric monoidal structure on $\D(\mathcal A)$ we employ \cite[Proposition 4.1.7.4]{lurie-higher-algebra}: Let $W$ be the set of morphisms $M \to N$ in $\D(\underline{\mathcal A})$ such that $M \tensor_{\underline{\mathcal A}} \mathcal A \isoto N \tensor_{\underline{\mathcal A}} \mathcal A$ is an isomorphism. We have to show that given $(M \to N) \in W$ and any $P \in \D(\underline{\mathcal A})$, we have $(M \tensor P \to N \tensor P) \in W$. Note that $W$ is generated by the morphisms $M \to M \tensor_{\underline{\mathcal A}} \mathcal A$, so we can assume that $(M \to N) = (M \to M \tensor_{\underline{\mathcal A}} \mathcal A)$. Writing $M$ as a colimit of $\underline{\mathcal A}[S]$'s we can reduce to the case $M = \underline{\mathcal A}[S]$ for some extremally disconnected set $S$. It is now enough to show that for every $Q \in \D(\mathcal A)$ the natural map
\begin{align*}
	\Hom(\mathcal A[S] \tensor P, Q) \isoto \Hom(\underline{\mathcal A}[S] \tensor P, Q)
\end{align*}
is an isomorphism (then the claim follows from Yoneda). But by the adjunction of $\tensor$ and $\IHom$, this reduces to the claim that $\IHom_{\underline{\mathcal A}}(P, Q)$ is in $\D(\mathcal A)$. This follows formally from the adjunction of $\tensor$ and $\IHom$ (here it is important to use the \emph{internal} hom's in the definition of analytic rings). This finishes the proof of (i).

Part (ii) is formal: Use \cite[Proposition 4.1.7.4]{lurie-higher-algebra} to get the symmetric monoidal structure on the functors $- \tensor_{\mathcal A} \mathcal B$. To get the desired functoriality, first note that by \cite[Theorem 4.5.3.1]{lurie-higher-algebra} we have a natural functor from $\AnRing^\dagger_{(V,\mm)}$ to the $\infty$-category of symmetric monoidal $\infty$-categories sending $\mathcal A \mapsto \D(\underline{\mathcal A})$. This functor is classified by a coCartesian family of symmetric monoidal $\infty$-categories $\mathcal C^\tensor \to \opComm \cprod \AnRing^\dagger_{(V,\mm)}$ (see \cite[Definition 4.8.3.1]{lurie-higher-algebra}). Now let $\mathcal D^\tensor \subset \mathcal C^\tensor$ be the full subcategory spanned by those objects which belong to $\D(\mathcal A)$ in each fiber (instead of merely belonging to $\D(\underline{\mathcal A})$). The existence of the symmetric monoidal functors $- \tensor_{\mathcal A} \mathcal B$ implies that $\mathcal D^\tensor \to \opComm \cprod \AnRing^\dagger_{(V,\mm)}$ is still a coCartesian fibration (e.g. use the criterion in \cite[Proposition 2.4.4.3]{lurie-higher-topos-theory}) and is thus classified by the desired functor $\AnRing^\dagger_{(V,\mm)} \to \infcatinf^\tensor$.
\end{proof}

\begin{proposition}
Let $\tau$ be a regular cardinal such that $(V,\mm)$ is $\tau$-compact.
\begin{propenum}
	\item Let $\mathcal A$ be an uncompleted analytic ring over $(V,\mm)$ and let $\kappa$ be a strong limit cardinal such that $\underline{\mathcal A} \in \D(V,\mm)_\kappa$. Then $\D(\mathcal A)_\kappa$ is stable under the symmetric monoidal structure in $\D(\mathcal A)$ and thus becomes itself symmetric monoidal.

	\item Let $\mathcal A \to \mathcal B$ be a morphism of uncompleted analytic rings over $(V,\mm)$. Then for every strong limit cardinal $\kappa$ such that $\underline{\mathcal A}, \underline{\mathcal B} \in \D(V,\mm)_\kappa$ the functor $- \tensor_{\mathcal A} \mathcal B$ restricts to a symmetric monoidal functor $\D(\mathcal A)_\kappa \to \D(\mathcal B)_\kappa$.
\end{propenum}
\end{proposition}
\begin{proof}
Part (ii) follows immediately from (i). To prove (i), let $\mathcal A$ and $\kappa$ be given. To show that $\D(\mathcal A)_\kappa$ is stable under the symmetric monoidal structure, we use the fact that this symmetric monoidal structure preserves all small colimits in both arguments to reduce the claim to showing that for all $\kappa$-small exremally disconnected sets $S$ and $T$ we have $\mathcal A[S] \tensor_{\mathcal A} \mathcal A[T] \in \D(\mathcal A)_\kappa$. We have
\begin{align*}
	\mathcal A[S] \tensor_{\mathcal A} \mathcal A[T] = (\underline{\mathcal A}[S] \tensor_{\underline{\mathcal A}} \underline{\mathcal A}[T]) \tensor_{\underline{\mathcal A}} \mathcal A.
\end{align*}
But clearly $\underline{\mathcal A}[S] \tensor_{\underline{\mathcal A}} \underline{\mathcal A}[T] \in \D(\underline{\mathcal A})_\kappa$ (because this category is stable under the symmetric monoidal structure), hence this object is a (sifted) colimit of objects $\mathcal A[S']$ for $\kappa$-small extremally disconnected sets $S'$. As $- \tensor_{\underline{\mathcal A}} \mathcal A$ preserves all colimits and maps $\underline{\mathcal A}[S']$ to $\mathcal A[S']$ we deduce $\mathcal A[S] \tensor_{\mathcal A} \mathcal A[T] \in \D(\mathcal A)_\kappa$, as desired.
\end{proof}

In the definition of uncompleted analytic rings it is possible that $\mathcal A[*]$ is not equal to $\mathcal A$. This behavior is undesirable in practice, so from now on we restrict to the case that $\mathcal A[*] = \mathcal A$. This leads to the following definition, which is the main definition of \cref{sec:andesc}.

\begin{definition} \label{def:analytic-ring}
An \emph{analytic (associative) ring} over $(V,\mm)$ is an uncompleted analytic (associative) ring $\mathcal A$ over $(V,\mm)$ such that the map $\underline{\mathcal A} \to \mathcal A[*]$ is an isomorphism. We denote
\begin{align*}
	\AnRing_{(V,\mm)}, \qquad \AnAssRing_{(V,\mm)},
\end{align*}
the full subcategories of $\AnRing^\dagger_{(V,\mm)}$ and $\AnAssRing^\dagger_{(V,\mm)}$ spanned by the analytic (associative) rings. We also abbreviate
\begin{align*}
	\AnRing := \AnRing_{(\Z,\Z)},
\end{align*}
the $\infty$-category of analytic rings.
\end{definition}

Our next step is to construct a completion functor $\mathcal A \mapsto \mathcal A^=$ which takes an uncompleted analytic (associative) ring $\mathcal A$ and produces an analytic (associative) ring $\mathcal A^=$. Basically, $\mathcal A^=$ will be defined as $\underline{\mathcal A^=} = \mathcal A[*]$ and $\mathcal A^=[S] = \mathcal A[S]$. To make this construction work, we need to check two things: Firstly, we need to equip $\mathcal A[*]$ with the structure of a (associative) ring and secondly we need to verify that the induced analytic ring structure on $\mathcal A[*]$ is indeed an analytic ring structure. This will be the content of the next results.

\begin{lemma} \label{rslt:existence-of-initial-ring-over-unnormalized-analytic-ring}
Let $\mathcal A$ be an uncompleted analytic (associative) ring over $(V,\mm)$. Consider the $\infty$-category of (associative) $\underline{\mathcal A}$-algebras $B$ such that $B$ lies in $\D(\mathcal A)$. This $\infty$-category has an initial object, whose underlying object in $\D(\mathcal A)$ is $\mathcal A[*]$.
\end{lemma}
\begin{proof}
We only treat the case of commutative rings. In the associative case, if $\underline{\mathcal A}$ comes from an $\mathbb E_\infty$-algebra (but $B$ may not) then one can argue similar to the commutative case, cf. the first part of the proof of \cite[Proposition 12.26]{scholze-analytic-spaces}). If $\underline{\mathcal A}$ is not commutative however, the argument in loc. cit. is wrong, as pointed out by Longke Tang and Yixiao Li. They suggested a fix, but we decided not to include it here as we do not really need this case.

From now on we work in the commutative case. Let $\Ring_{\mathcal A}$ be the described $\infty$-category of $\underline{\mathcal A}$-algebras. By definition it embedds into $\Ring_{\underline{\mathcal A}}$. It is enough to construct a left adjoint
\begin{align*}
	L\colon \Ring_{\underline{\mathcal A}} \to \Ring_{\mathcal A}
\end{align*}
of this embedding which acts as $- \tensor_{\underline{\mathcal A}} \mathcal A$ on the underlying $\underline{\mathcal A}$-modules; then $L(\underline{\mathcal A})$ is the desired initial object. If $(V,\mm) = (V,V)$ then we can directly apply \cref{rslt:construct-left-adjoint-on-animated-rings} to $A = \underline{\mathcal A}$ and $\mathcal C = \D_{\ge0}(\mathcal A)$: The only non-trivial condition to check is (c), but this follows easily from our assumption that the Frobenius is a morphism of analytic rings.

For general $(V,\mm)$ we can argue in the same way by generalizing \cref{rslt:construct-left-adjoint-on-animated-rings} to almost rings $A \in \Ring_{(V,\mm)}$. Namely, the relevant monad $T$ on $\D_{\ge0}(A)$ is given by $T(M) = \bigdsum_{n\ge0} \Sym^n_{A_{**}}(\tau_{\ge0} M_*)^a$ for $M \in \D_{\ge0}(A)$. Note that by the proof of \cref{rslt:almost-localization-for-rings}, if $M \to N$ is an almost isomorphism in $\D_{\ge0}(A_{**})$ then $\Sym^n_{A_{**}}(M) \isoto \Sym^n_{A_{**}}(N)$ is also an almost isomorphism. We therefore get a good notion of $\Sym^n_A(M) \in \D_{\ge0}(A)$ for every $M \in \D_{\ge0}(A)$; e.g. $\Sym^n_A(M) := \Sym^n_{A_{**}}(\tau_{\ge0} M_*)^a$. We need to show that if $M \to N$ is a map in $\D_{\ge0}(A)$ such that $L(M) \isoto L(N)$ is an isomorphism (notation as in \cref{rslt:construct-left-adjoint-on-animated-rings}) then $L(\Sym^n_A(M)) \to L(\Sym^n_A(N))$ is an isomorphism. Since almost localization preserves all the relevant constructions, we can argue as in the proof of \cite[Lemma 12.27]{scholze-analytic-spaces}.
\end{proof}

\begin{proposition} \label{rslt:existence-of-normalization-of-analytic-rings}
The inclusions $\AnAssRing_{(V,\mm)} \injto \AnAssRing^\dagger_{(V,\mm)}$ and $\AnRing_{(V,\mm)} \injto \AnRing^\dagger_{(V,\mm)}$ admit left adjoints
\begin{align*}
	(-)^=\colon \AnAssRing^\dagger_{(V,\mm)} \to \AnAssRing_{(V,\mm)}, \qquad &\mathcal A \mapsto \mathcal A^=,\\
	(-)^=\colon \AnRing^\dagger_{(V,\mm)} \to \AnRing_{(V,\mm)}, \qquad &\mathcal A \mapsto \mathcal A^=.
\end{align*}
We naturally have $\D(\mathcal A^=) = \D(\mathcal A)$ via the induced functor.
\end{proposition}
\begin{proof}
We argue only for $\AnRing$; the proof for $\AnAssRing$ is analogous. Given an uncompleted analytic ring $\mathcal A \in \AnRing^\dagger_{(V,\mm)}$, we claim that the $\infty$-category $\AnRing_{\mathcal A/}$ has an initial object. Namely, let $\Ring_{\mathcal A} \subset \Ring_{\underline{\mathcal A}}$ be the full subcategory of those $\underline{\mathcal A}$-algebras $B$ such that $B \in \D(\mathcal A)$. Note that the forgetful functor $\mathrm{For}\colon \AnRing_{\mathcal A/} \to \Ring_{\underline{\mathcal A}}$ factors over $\Ring_{\mathcal A}$ and that the induced uncompleted analytic ring structure from \cref{rslt:existence-of-induced-unnormalized-analytic-ass-ring-structure} provides a left adjoint $L\colon \Ring_{\mathcal A} \to \AnRing_{\mathcal A/}$ of $\mathrm{For}$. But by \cref{rslt:existence-of-initial-ring-over-unnormalized-analytic-ring} $\Ring_{\mathcal A}$ has an initial object $B$, hence $L(B)$ is an initial object of $\AnRing_{\mathcal A/}$.

Now define $(-)^=$ as the functor which associates to each $\mathcal A \in \AnRing^\dagger_{(V,\mm)}$ an initial object of $\AnRing_{\mathcal A/}$; this functor can be constructed as follows. Start with the coCartesian fibration $X' \to \Delta^1$ classifying the identity functor $\AnRing^\dagger_{(V,\mm)} \to \AnRing^\dagger_{(V,\mm)}$, then pass to the full subcategory $X' \subset X$ where we allow only (complete) analytic rings in the fiber over $1 \in \Delta^1$. The induced fibration $X \to \Delta^1$ is still coCartesian: For every $\mathcal A \in X$ lying in the fiber over $0 \in \Delta^1$ pick $\mathcal A \to \mathcal A^=$ the initial edge lying over $0 \to 1$ in $\Delta^1$. To see that this is a coCartesian edge, by \cite[Proposition 2.4.4.3]{lurie-higher-topos-theory} is is enough to verify that for every analytic ring $\mathcal B$ the induced map
\begin{align*}
	\Hom(\mathcal A^=, \mathcal B) \isoto \Hom(\mathcal A, \mathcal B)
\end{align*}
is an isomorphism of anima. But as $\mathcal A^=$ is initial in $\AnRing_{\mathcal A/}$ we have $\AnRing_{\mathcal A/} \isom \AnRing_{\mathcal A^=/}$, which immediately implies the claim. Of course $X \to \Delta^1$ is also a Cartesian fibration, classifying the inclusion $\AnRing_{(V,\mm)} \to \AnRing^\dagger_{(V,\mm)}$. This shows that the functor $\mathcal A \to \mathcal A^=$ is indeed left adjoint to this inclusion.
\end{proof}

\begin{definition}
\begin{defenum}
	\item For every uncompleted analytic (associative) ring $\mathcal A$ over $(V,\mm)$, the analytic (associative) ring $\mathcal A^=$ defined in \cref{rslt:existence-of-normalization-of-analytic-rings} is called the \emph{completion of $\mathcal A$}.

	\item \label{def:induced-analytic-ring-structure} Let $\mathcal A$ be an analytic (associative) ring over $(V,\mm)$. We define the functor
	\begin{align*}
		(-)_{\mathcal A/}\colon \AssRing_{\underline{\mathcal A}} \to \AnAssRing_{\mathcal A/}, \qquad B \mapsto B_{\mathcal A/},
	\end{align*}
	resp.
	\begin{align*}
		(-)_{\mathcal A/}\colon \Ring_{\underline{\mathcal A}} \to \AnRing_{\mathcal A/}, \qquad B \mapsto B_{\mathcal A/}
	\end{align*}
	to be the composition of the functor $L$ from \cref{rslt:existence-of-induced-unnormalized-analytic-ass-ring-structure} and the completion functor $(-)^=$. We call $B_{\mathcal A/}$ the \emph{induced analytic (associative) ring} on $B$ over $\mathcal A$.
\end{defenum}
\end{definition}

We have established a good basis for the theory of analytic rings over $(V,\mm)$, so we now start studying some of its properties. First up, we show that analytic ring structures on a $(V,\mm)$-algebra $A$ essentially only depend on the static ring $\pi_0(A)$:

\begin{proposition} \label{rslt:topological-invariance-of-analytic-ring-structures}
Let $A$ be a (associative) $(V,\mm)$-algebra. Then there is a bijective correspondence between analytic ring structures $\mathcal A$ on $A$ and analytic ring structures $\mathcal A_0$ on $\pi_0 A$. In the forward direction, this takes $\mathcal A$ to $\mathcal A_0 = (\pi_0 A)_{\mathcal A/}$.
\end{proposition}
\begin{proof}
By \cite[Proposition 12.21]{scholze-analytic-spaces} this holds in the case $(V,\mm) = (\Z, \Z)$ and for associative rings, but the same argument applies to general $(V,\mm)$ and commutative rings (the condition on Frobenius is clearly automatic in both directions). Also note that loc. cit. deals with uncompleted analytic ring structures, but for any (complete) analytic ring structure $\mathcal A$ on $A$, we have $\pi_0 A \in \D(\mathcal A)$ so that $(\pi_0 A)_{\mathcal A/}$ does not require any completions.
\end{proof}

Next up we study colimits of analytic rings over $(V,\mm)$. It turns out that all small colimits exist and admit a fairly explicit description:

\begin{proposition} \label{rslt:colimits-of-analytic-rings}
The $\infty$-category $\AnRing_{(V,\mm)}$ has all small colimits. More precisely:
\begin{propenum}
	\item The initial object of $\AnRing_{(V,\mm)}$ is the ring $V^a$.

	\item \label{rslt:sifted-colimits-of-analytic-rings} Given a sifted colimit $\mathcal A = \varinjlim_i \mathcal A_i$ in $\AnRing_{(V,\mm)}$, we have $\underline{\mathcal A} = \varinjlim_i \underline{\mathcal A_i}$ and $\mathcal A[S] = \varinjlim_i \mathcal A_i[S]$ for all extremally disconnected sets $S$.

	\item \label{rslt:pushouts-of-analytic-rings} Given a diagram $\mathcal B \from \mathcal A \to \mathcal C$ in $\AnRing_{(V,\mm)}$, the pushout $\mathcal E = \mathcal B \tensor_{\mathcal A} \mathcal C$ is the completion of the uncompleted analytic ring structure on $\underline{\mathcal B} \tensor_{\underline{\mathcal A}} \underline{\mathcal C}$ such that $\D(\mathcal E) \subset \D(\underline{\mathcal B} \tensor_{\underline{\mathcal A}} \underline{\mathcal C})$ consists precisely of those objects $M \in \D(\underline{\mathcal B} \tensor_{\underline{\mathcal A}} \underline{\mathcal C})$ with $M \in \D(\mathcal B)$ and $M \in \D(\mathcal C)$ via the forgetful functors. For every extremally disconnected set $S$, $\mathcal E[S]$ is the colimit over the repeated application of $(- \tensor_{\underline{\mathcal B}} \mathcal B) \tensor_{\underline{\mathcal C}} \mathcal C$ to $(\underline{\mathcal B} \tensor_{\underline{\mathcal A}} \underline{\mathcal C})[S]$.
\end{propenum}
\end{proposition}
\begin{proof}
This follows from the same proof as in \cite[Proposition 12.12]{scholze-analytic-spaces}.
\end{proof}

While sifted colimits of analytic rings over $(V,\mm)$ are very explicit, pushouts become much more subtle than what one might be used to from classical ring theory. This is because of the potentially countable colimit of repeated applications of $(- \tensor_{\underline{\mathcal B}} \mathcal B) \tensor_{\underline{\mathcal C}} \mathcal C$ involved in the computation of the generators $\mathcal E[S] = (\mathcal B \tensor_{\mathcal A} \mathcal C)[S]$. This means in particular that in general, for $M \in \D(\mathcal C)$ it might happen that $M \tensor_{\mathcal A} \mathcal B \ne M \tensor_{\mathcal C} (\mathcal B \tensor_{\mathcal A} \mathcal C)$, i.e. base-change does not hold in general. This phenomenon is not new: Similar subtleties arise for pushouts of Huber pairs along non-adic maps. While in general these subtleties around pushouts of analytic rings are an expected feature for a theory of ``complete topological rings'', it is useful to single out a special case of maps $\mathcal A \to \mathcal B$ where pushouts behave as nicely as possible. The following definition can be seen as an analog of adic maps.

\begin{definition} \label{def:steady-morphism-of-analytic-rings}
A map $\mathcal A \to \mathcal B$ of analytic rings over $(V,\mm)$ is called \emph{steady} if for all maps $\mathcal A \to \mathcal C$ of analytic rings over $(V,\mm)$, the functor $- \tensor_{\mathcal A} \mathcal B$ preserves the full subcategory $\D_{\ge0}(\mathcal C)$ of $\D_{\ge0}(\underline{\mathcal C})$.
\end{definition}

\begin{remark}
Our \cref{def:steady-morphism-of-analytic-rings} is taken straight from \cite[Definition 12.13]{scholze-analytic-spaces} (with a minor fix, in accordance with Scholze). However it is still an open question whether this definition for steadiness is really the ``right one''. The issue is that the steadiness of a map $f\colon \mathcal A \to \mathcal B$ depends on the whole $\infty$-category $\AnRing_{(V,\mm)}$ rather than being a more intrinsic property of the map $f$ itself. Also, there are some open questions regarding the interplay of steady maps and nuclear modules, cf. \cref{rmk:is-D-nuc-a-sheaf}.
\end{remark}

There are several equivalent possible definitions for steady maps, which are summarized by the following result.

\begin{proposition}
Let $f\colon \mathcal A \to \mathcal B$ be a map of analytic rings over $(V,\mm)$. Then the following are equivalent:
\begin{propenum}
	\item $f$ is steady.

	\item \label{rslt:steady-map-implies-easy-calculation-of-generators} For every map $\mathcal A \to \mathcal C$ of analytic rings and every extremally disconnected set $S$, the natural map
	\begin{align*}
		\mathcal C[S] \tensor_{\mathcal A} \mathcal B \isoto (\mathcal C \tensor_{\mathcal A} \mathcal B)[S]
	\end{align*}
	is an isomorphism of objects in $\D(\underline{\mathcal C} \tensor_{\underline{\mathcal A}} \underline{\mathcal B})$.

	\item \label{rslt:steady-map-of-analytic-rings-equiv-base-change} For all pushout diagrams
	\begin{center}\begin{tikzcd}
		\mathcal B' & \mathcal A' \arrow[l,"f'"]\\
		\mathcal B \arrow[u,"g'"] & \mathcal A \arrow[l,"f"] \arrow[u,"g"]
	\end{tikzcd}\end{center}
	in $\AnRing_{(V,\mm)}$, the natural map
	\begin{align*}
		f^* g_* \isoto g'_* f'^*
	\end{align*}
	of functors $\D(\mathcal A') \to \D(\mathcal B)$ is an isomorphism. Here $g_*\colon \D(\mathcal A') \to \D(\mathcal A)$ is the forgetful functor and $f^*\colon \D(\mathcal A) \to \D(\mathcal B)$ is the functor $- \tensor_{\mathcal A} \mathcal B$ (and similarly for $g'_*$ and $f'^*$).
\end{propenum}
\end{proposition}
\begin{proof}
We first show that (i) implies (ii), so assume that $f$ is steady and pick $\mathcal C$ and $S$ as in (ii). By definition of steadiness the restriction of $M := \mathcal C[S] \tensor_{\mathcal A} \mathcal B$ to $\D(\underline{\mathcal C})$ lies in $\D(\mathcal C)$. Of course also the restriction of $M$ to $\D(\underline{\mathcal B})$ lies in $\D(\mathcal B)$, hence by \cref{rslt:pushouts-of-analytic-rings} we have $M \in \D(\mathcal C \tensor_{\mathcal A} \mathcal B)$. One checks immediately that $M$ satisfies the universal property of $(\mathcal C \tensor_{\mathcal A} \mathcal B)[S]$, which implies the claimed isomorphism.

We now show that (ii) implies (iii), so let $\mathcal A'$ and $\mathcal B' = \mathcal A' \tensor_{\mathcal A} \mathcal B$ be given. As all functors in the base-change diagram preserve colimits, we are reduced to showing the isomorphism $f^* g_* M \isoto g'_* f'^* M$ for $M = \mathcal A'[S]$, where $S$ is any extremally disconnected set. But then this is just (ii) with $\mathcal C = \mathcal A'$.

It is clear that (ii) is just a special case of (iii), hence (ii) and (iii) are equivalent. It is similarly easy to see that (iii) implies (i).
\end{proof}

In practice, to show that a map $f\colon \mathcal A \to \mathcal B$ of analytic rings over $(V,\mm)$ is steady, one usually tries to break $f$ into simple pieces for which steadiness is easier to verify. This is possible by the following stability properties enjoyed by the class of steady morphisms:

\begin{proposition}
\begin{propenum}
	\item \label{rslt:stability-of-steady-maps} The class of steady maps in $\AnRing_{(V,\mm)}$ is stable under composition, base-change and all colimits.

	\item \label{rslt:steadyness-satisfies-2-out-of-3} Let $f\colon \mathcal A \to \mathcal B$ and $g\colon \mathcal B \to \mathcal C$ be maps of analytic rings over $(V,\mm)$. If $f$ and $g \comp f$ are steady, then so is $g$.
\end{propenum}
\end{proposition}
\begin{proof}
It follows easily from the base-change criterion in \cref{rslt:steady-map-of-analytic-rings-equiv-base-change} that steady maps are stable under composition and base-change. To show stability under colimits, it is enough to consider sifted colimits and finite coproducts. Suppose $(\mathcal A \to \mathcal B) = \varinjlim_i (\mathcal A_i \to \mathcal B_i)$ is a sifted colimit of steady maps in $\AnRing_{(V,\mm)}$. Replacing all maps $\mathcal A_i \to \mathcal B_i$ with their base-change to $\mathcal A$, we can assume $\mathcal A_i = \mathcal A$ for all $i$. Then $- \tensor_{\mathcal A} \mathcal B = \varinjlim_i - \tensor_{\mathcal A} \mathcal B_i$ by \cref{rslt:sifted-colimits-of-analytic-rings}, so the steadiness of $\mathcal A \to \mathcal B$ follows easily from the base-change criterion \cref{rslt:steady-map-of-analytic-rings-equiv-base-change}. It remains to handle the case of coproducts, so assume that $\mathcal A_1 \to \mathcal B_1$ and $\mathcal A_2 \to \mathcal B_2$ are steady. Then
\begin{align*}
	\mathcal B_1 \tensor \mathcal B_2 = (\mathcal B_1 \tensor \mathcal A_2) \tensor_{(\mathcal A_1 \tensor \mathcal A_2)} (\mathcal A_1 \tensor \mathcal B_2).
\end{align*}
Now $(\mathcal B_1 \tensor \mathcal A_2) \to (\mathcal A_1 \tensor \mathcal A_2)$ and $(\mathcal A_1 \tensor \mathcal B_2) \to (\mathcal A_1 \tensor \mathcal A_2)$ are steady by stability under base-change, hence $(\mathcal B_1 \tensor \mathcal B_2) \to (\mathcal A_1 \tensor \mathcal A_2)$ is the composition of two steady maps and hence itself steady. This finishes the proof of (i).

We now prove (ii), so let $f$ and $g$ be given with $f$ and $g \comp f$ steady. Pick any morphism $\mathcal B \to \mathcal B'$ of analytic rings and any $M \in \D(\mathcal B')$. We need to show that $M \tensor_{\mathcal B} \mathcal C$, considered as an object of $\D(\underline{\mathcal B'})$, lies in $\D(\mathcal B')$. To show this, recall that the bar construction allows us to write $\underline{\mathcal B} = \varinjlim_{n\in\Delta} \underline{\mathcal B}^{\tensor_{\mathcal A} n+2}$ in $\D(\mathcal A)$. Applying $- \tensor_{\underline{\mathcal B}} \mathcal B$ and $- \tensor_{\mathcal B} M$ we obtain $M = \varinjlim_{n\in\Delta} M \tensor_{\mathcal A} \mathcal B^{\tensor_{\mathcal A} n+1}$. Thus we see that
\begin{align*}
	M \tensor_{\mathcal B} \mathcal C = \varinjlim_{n\in\Delta} M \tensor_{\mathcal A} \mathcal B^{\tensor_{\mathcal A}n} \tensor_{\mathcal A} \mathcal C.
\end{align*}
By steadiness of $f$, $- \tensor_{\mathcal A} \mathcal B$ preserves $\D(\mathcal B')$ and by steadiness of $g \comp f$, $- \tensor_{\mathcal A} \mathcal C$ preserves $\D(\mathcal B')$. As $\D(\mathcal B')$ is stable under colimits, we deduce $M \tensor_{\mathcal B} \mathcal C \in \D(\mathcal B')$, as desired.
\end{proof}

We now introduce nuclear $\mathcal A$-modules, which form a special subcategory of $\D(\mathcal A)$. Nuclearity of modules can be seen as an analogous property to steadiness of algebras; for example it implies a similar base-change property (see \cref{rslt:nucelar-modules-are-steady} below). Moreover, nuclear modules are useful to build steady maps of analytic rings via induced structures (see \cref{rslt:nuclear-induced-ring-is-steady}).

\begin{definition} \label{def:nuclear-module}
Let $\mathcal A$ be an analytic ring over $(V,\mm)$. An object $N \in \D(\mathcal A)$ is called \emph{nuclear} if for all extremally disconnected sets $S$ and all $M \in \D(\mathcal A)$ the natural map
\begin{align*}
	(\IHom_{\mathcal A}(\mathcal A[S], M) \tensor_{\mathcal A} N)(*) \isoto (M \tensor_{\mathcal A} N)(S)
\end{align*}
is an isomorphism in $\D(V,\mm)_\omega$. We denote
\begin{align*}
	\D(\mathcal A)^\nuc \subset \D(\mathcal A)
\end{align*}
the full subcategory spanned by the nuclear objects.
\end{definition}

\begin{remark}
In the case $\mm = V$ the condition on $N$ in \cref{def:nuclear-module} can be weakened significantly without changing the notion of nuclearity, see \cite[Definition 13.10]{scholze-analytic-spaces} and \cite[Proposition 13.14]{scholze-analytic-spaces}. We do not know if the same works for general $(V,\mm)$.
\end{remark}

\begin{proposition}
Let $f\colon \mathcal A \to \mathcal B$ be a map of analytic rings over $(V,\mm)$.
\begin{propenum}
	\item \label{rslt:nuclear-modules-stable-under-colim} The subcategory $\D(\mathcal A)^\nuc \subset \D(\mathcal A)$ is stable under all small colimits.

	\item \label{rslt:nucelar-modules-are-steady} Let $N \in \D(\mathcal A)^\nuc$. Then for all $M \in \D(\mathcal B)$ the natural morphism
	\begin{align*}
		M \tensor_{\mathcal A} N \isoto M \tensor_{\mathcal B} f^* N
	\end{align*}
	is an isomorphism in $\D(\mathcal A)$.

	\item \label{rslt:nuclear-modules-stable-under-pullback} The pullback functor $f^*\colon \D(\mathcal A) \to \D(\mathcal B)$ preserves nuclearity.
\end{propenum}
\end{proposition}
\begin{proof}
Part (i) is immediate from the definition. We now prove (ii), so let $M$ and $N$ be given as in the claim. We have to show that the object $L := M \tensor_{\mathcal A} N \in \D(\underline{\mathcal B}_{\mathcal A/})$ lies in $\D(\mathcal B)$. Fix some large strong limit cardinal $\kappa$ and denote $\mathcal C_0 \subset \D(\underline{\mathcal B}_{\mathcal A/})_\kappa$ and $\mathcal C'_0 \subset \D(\mathcal B)_\kappa$ the full subcategories generated under retracts, shifts and finite colimits by the objects $\underline{\mathcal B}_{\mathcal A/}[S]$ and $\mathcal B[S]$ for $\kappa$-small extremally disconnected profinite sets $S$, respectively. There is a natural commuting diagram
\begin{center}\begin{tikzcd}
	\Ind(\mathcal C'_0) \arrow[r] \arrow[d] & \D(\mathcal B)_\kappa \arrow[d]\\
	\Ind(\mathcal C_0) \arrow[r] & \D(\underline{\mathcal B}_{\mathcal A/})_\kappa
\end{tikzcd}\end{center}
Here the horizontal maps are the unique colimit-preserving extensions of the inclusions $\mathcal C_0 \subset \D(\underline{\mathcal B}_{\mathcal A/})_\kappa$ and $\mathcal C'_0 \subset \D(\mathcal B)$. The right vertical map is the forgetful functor $\D(\mathcal B)_\kappa \to \D(\underline{\mathcal B}_{\mathcal A/})$ followed by the restriction to $\kappa$-condensed modules. The left vertical map is the one induced by the functor $- \tensor_{\underline{\mathcal B}_{\mathcal A/}} \mathcal B\colon \mathcal C_0 \to \mathcal C'_0$ via viewing $\Ind(-)$ as a full subcategory of the presheaf $\infty$-category $\mathcal P(-)$. To see that the diagram commutes, first note that all maps preserve filtered colimits, so it is enough to check that the map commutes on objects $P \in \mathcal C'_0 \subset \Ind(\mathcal C'_0)$, where it is easy to see (e.g. use the fully faithful right adjoint of the bottom horizontal map which exists by \cref{rslt:almost-compact-generators-fully-faithful-into-Ind}). After choosing $\kappa$ large enough, we have $L \in \D(\underline{\mathcal B}_{\mathcal A/})_\kappa$. Via picking a fully faithful right adjoint of the lower horizontal map (see \cref{rslt:almost-compact-generators-fully-faithful-into-Ind}) we can equivalently view $L$ as an object in $\Ind(\mathcal C_0)$; namely, it is the one with $L(P) = \Hom(P, L)$ for $P \in \mathcal C_0^\opp$ (viewing $\Ind(\mathcal C_0) \subset P(\mathcal C_0)$). On the other hand, for every extremally disconnected set $S$ the nuclearity of $N$ implies
\begin{align*}
	&L(\underline{\mathcal B}_{\mathcal A/}[S]) = \Hom(\underline{\mathcal B}_{\mathcal A/}[S], M \tensor_{\mathcal A} N) = (M \tensor_{\mathcal A} N)(S) =\\&\qquad= (\IHom_{\mathcal A}(\mathcal A[S], M) \tensor_{\mathcal A} N)(*) = (\IHom_{\mathcal B}(\underline{\mathcal B}_{\mathcal A/}[S] \tensor_{\underline{\mathcal B}_{\mathcal A/}} \mathcal B, M) \tensor_{\mathcal A} N)(*).
\end{align*}
Since all the involved functors preserve finite colimits and retracts, the same holds for any $P \in \mathcal C_0$ in place of $\underline{\mathcal B}_{\mathcal A/}[S]$. Therefore, $L$ is the image of the object
\begin{align*}
	L'\colon \mathcal (C_0')^\opp \to \Ani, \qquad P' \mapsto (\IHom_{\mathcal B}(P, M) \tensor_{\mathcal A} N)(*)
\end{align*}
under the map $\Ind(\mathcal C_0') \to \Ind(\mathcal C_0)$. By the above commuting diagram this implies that $L$ comes from some object in $\D(\mathcal B)_\kappa$. As this is true for all $\kappa$ we can deduce that $L \in \D(\mathcal B)$. This finishes the proof of (ii).

We now prove (iii), so pick any nuclear $N \in \D(\mathcal A)^\nuc$. Using (ii) we compute for every $M \in \D(\mathcal B)$ and every extremally disconnected set $S$
\begin{align*}
	(\IHom_{\mathcal B}(\mathcal B[S], M) \tensor_{\mathcal B} f^* N)(*) &= (\IHom_{\mathcal B}(\mathcal B[S], M) \tensor_{\mathcal A} N)(*)\\
	&= (\IHom_{\mathcal A}(\mathcal A[S], M) \tensor_{\mathcal A} N)(*)\\
	&= (M \tensor_{\mathcal A} N)(S)\\
	&= (M \tensor_{\mathcal B} f^* N)(S),
\end{align*}
proving that $f^*N$ is indeed nuclear in $\D(\mathcal B)$.
\end{proof}

\begin{corollary} \label{rslt:nuclear-induced-ring-is-steady}
Let $\mathcal A$ be an analytic ring over $(V,\mm)$ and $B$ an $\underline{\mathcal A}$-algebra. If $B$ lies in $\D(\mathcal A)$ and is nuclear as an $\mathcal A$-module, then the map of analytic rings $\mathcal A \to B_{\mathcal A/}$ is steady.
\end{corollary}
\begin{proof}
This follows directly from \cref{rslt:nucelar-modules-are-steady} and the definition of steadiness.
\end{proof}

Having studied the $\infty$-category of analytic rings over a fixed almost setup $(V,\mm)$, we now discuss functoriality of this construction with respect to a change of almost setups. In particular, for every morphism $\varphi\colon (V,\mm) \to (V',\mm')$ of almost setups we get a functor $\AnRing_{(V,\mm)} \to \AnRing_{(V',\mm')}$ mapping an analytic ring $\mathcal A$ over $(V,\mm)$ to an analytic ring $\varphi^* \mathcal A$ over $(V',\mm')$. This also induces a symmetric monoidal functor $\varphi^*\colon \D(\mathcal A) \to \D(\varphi^* \mathcal A)$ and all of these functors are completely functorial in $\varphi$ and $\mathcal A$. More precisely, we have the following results.

\begin{lemma} \label{rslt:functoriality-of-analytic-rings-over-AlmSetup}
There are natural functors
\begin{align*}
	\AlmSetup \to \infcatinf, \qquad &(V,\mm) \mapsto \AnAssRing_{(V,\mm)},\\
	\AlmSetup \to \infcatinf, \qquad &(V,\mm) \mapsto \AnRing_{(V,\mm)}.
\end{align*}
Moreover, for every morphism $\varphi\colon (V,\mm) \to (V',\mm')$ of almost setups we have the following properties:
\begin{lemenum}
	\item For every analytic (associative) ring $\mathcal A$ over $(V,\mm)$, the induced analytic (associative) ring $\varphi^* \mathcal A$ over $(V',\mm')$ satisfies $\underline{\varphi^* \mathcal A} = \varphi^* \underline{\mathcal A}$ and $(\varphi^* \mathcal A)[S] = \varphi^* (\mathcal A[S])$ for all extremally disconnected sets $S$.

	\item \label{rslt:functoriality-of-analytic-rings-along-strict-morphism} Suppose $\varphi$ is strict. Then we have $\AnRing_{(V',\mm')} = (\AnRing_{(V,\mm)})_{/V'^a}$ and the induced functor $\AnRing_{(V,\mm)} \to \AnRing_{(V',\mm')}$ coincides with $- \tensor_{V^a} V'^a$. The same holds for associative rings.
\end{lemenum}
\end{lemma}
\begin{proof}
We only prove the claim for $\AnRing$; the case of $\AnAssRing$ is similar. Let $\Ring_{(-)} \to \AlmSetup$ be the coCartesian fibration classifying the functor from \cref{rslt:functoriality-of-almost-rings}. We can explicitly construct a new simplicial set $\mathcal C$ with a map $\mathcal C \to \Ring_{(-)}$ as follows. The objects of $\mathcal C$ are the analytic rings over every almost setup $(V,\mm)$. Given two objects $\mathcal A, \mathcal A' \in \mathcal C$ then a morphism $\mathcal A \to \mathcal A'$ in $\mathcal C$ is the same as a morphism $f\colon \underline{\mathcal A} \to \underline{\mathcal A'}$ in $\Ring_{(-)}$ such that the restriction of the induced functor $f_*\colon \D(\underline{\mathcal A'}) \to \D(\underline{\mathcal A})$ to $\D(\mathcal A')$ factors over $\D(\mathcal A)$. Higher homotopies in $\mathcal C$ are the same as higher homotopies in $\Ring_{(-)}$ via $\mathcal A \mapsto \underline{\mathcal A}$. Then it is clear that $\mathcal C \to \Ring_{(-)}$ is a categorical fibration, hence the same follows for the map $p\colon \mathcal C \to \AlmSetup$. Moreover, it is clear that the fiber of $p$ over any almost setup $(V,\mm)$ is precisely $\AnRing_{(V,\mm)}$. To construct the desired functor $\AlmSetup \to \infcatinf$ it is now enough to show that $p$ is a coCartesian fibration (then $p$ classifies the desired functor).

We first show that $p$ is a locally coCartesian fibration, so let $\varphi\colon (V,\mm) \to (V',\mm')$ be any morphism in $\AlmSetup$ and denote $p_\varphi\colon \mathcal C_\varphi \to \Delta_1$ the pullback of $p$. Pick any analytic ring $\mathcal A \in \mathcal C_{(V,\mm)} = \AnRing_{(V,\mm)}$. We need to construct a $p_\varphi$-coCartesian morphism $f\colon \mathcal A \to \mathcal A'$ in $\mathcal C$ which lies over $\varphi$. Since $q\colon \Ring_{(-)} \to \AlmSetup$ is a coCartesian fibration, there is a $q$-coCartesian morphism $f'\colon \underline{\mathcal A} \to \varphi^* \underline{\mathcal A}$ in $\Ring_{(-)}$ lying over $\varphi$. Define the pre-analytic ring $\varphi^* \mathcal A$ as in claim (i), i.e. with $\underline{\varphi^* \mathcal A} = \varphi^* \underline{\mathcal A}$ and $(\varphi^* \mathcal A)[S] = \varphi^* (\mathcal A[S])$ for all extremally disconnected sets $S$. We claim that $\varphi^* \mathcal A$ is an analytic ring. Factoring $\varphi$ as a composition of a strict morphism and a localizing morphism we can treat these two cases separately. If $\varphi$ is strict then we have $\varphi^* \mathcal A = \mathcal A \tensor_{V^a} V'^a$ (by the explicit computation of the right-hand side in \cref{rslt:pushouts-of-analytic-rings}) so the claim is clear.
We can thus assume that $\varphi$ is localizing, i.e. of the form $\varphi\colon (V,\mm) \to (V,\mm')$. Fix an object $M \in \D(\varphi^* \underline{\mathcal A})$ which is a sifted colimit of the objects $\varphi^*\mathcal A[S]$. Denote $A := (V,\mm)_{**} \underline{\mathcal A} \in \Ring_V$. Then we have $\underline{\mathcal A} = (A,\mm_A)^* A$ and there is a natural equivalence $\varphi^* \comp (A,\mm_A)^* = (A,\mm'_A)^*$ of functors $\D(A) \to \D(\varphi^*\underline{\mathcal A})$ (where we use the localization functors from \cref{rslt:general-properties-of-almost-localization-over-ass-alg}). Hence if we define $M' := (A,\mm_A)_! (A,\mm'_A)^* M \in \D(\underline{\mathcal A})$ then we have $\varphi^* M' = M$. To prove that $\varphi^*\mathcal A$ is analytic, we need to show that the natural morphism
\begin{align}
	\IHom_{\varphi^*\underline{\mathcal A}}(\varphi^*\mathcal A[S], M) \isoto \IHom_{\varphi^*\underline{\mathcal A}}(\varphi^*\underline{\mathcal A}[S], M) \label{eq:almost-localization-of-analytic-ring-is-analytic}
\end{align}
is an isomorphism in $\D(\varphi^*\underline{\mathcal A})$. Note that it follows from \cref{rslt:properties-of-almost-localization-over-alg} that for all $N_1, N_2 \in \D(\underline{\mathcal A})$ we have a natural isomorphism
\begin{align*}
	\varphi^* \IHom_{\underline{\mathcal A}}(N_1, N_2) = \IHom_{\varphi^*\underline{\mathcal A}}(\varphi^* N_1, \varphi^* N_2)
\end{align*}
in $\D(\varphi^* \underline{\mathcal A})$ (namely, we can apply the $\IHom$-identity in \cref{rslt:properties-of-almost-localization-over-alg} to $(A,\mm_A)_* N_i$ with respect to both $(A,\mm_A)^*$ and $(A,\mm'_A)^*$). Thus the claimed isomorphism \cref{eq:almost-localization-of-analytic-ring-is-analytic} reduces to showing that the natural map
\begin{align*}
	\IHom_{\underline{\mathcal A}}(\mathcal A[S], M') \isoto \IHom_{\underline{\mathcal A}}(\underline{\mathcal A}[S], M')
\end{align*}
is an isomorphism; in other words we need to see that $M' \in \D(\mathcal A)$. By assumption $M'$ is a sifted colimit of copies of $\varphi^* \mathcal A[S]$ and since both $(A,\mm'_A)_!$ and $(A,\mm_A)^*$ preserve all small colimits it follows that $M'$ is a sifted colimit of copies of $\mathcal E[S] := (A,\mm_A)^* (A,\mm'_A)_! \varphi^*\mathcal A[S]$.  Therefore it is enough to show that each $\mathcal E[S]$ lies in $\D(\mathcal A)$. Fix the extremally disconnected set $S$ and pick any $E \in \D(A)$ such that $(A,\mm_A)^*E = \mathcal A[S]$. Using the various computation rules in \cref{rslt:general-properties-of-almost-localization-over-alg} we compute
\begin{align*}
	&\mathcal E[S] = (A,\mm_A)^* (A,\mm'_A)_! \varphi^*\mathcal A[S] = (A,\mm_A)^* (A,\mm'_A)_! \varphi^* (A,\mm_A)^* E =\\&\qquad= (A,\mm_A)^* (A,\mm'_A)_! (A,\mm'_A)^* E = (A,\mm_A)^* (\widetilde\mm'_A \tensor_A \IHom_A(\widetilde\mm'_A, E)) =\\&\qquad= (A,\mm_A)^* \widetilde\mm'_A \tensor_{\underline{\mathcal A}} (A,\mm_A)^* \IHom_A(\widetilde\mm'_A, E).
\end{align*}
Now $\widetilde\mm'_A$ is a colimit of copies of $A$ (since it is discrete), and as $\D(\mathcal A)$ is stable under colimits it is enough to show that the second factor in the above tensor product lies in $\D(\mathcal A)$. But by \cref{rslt:properties-of-almost-localization-over-alg} we have
\begin{align*}
	(A,\mm_A)^* \IHom_A(\widetilde\mm'_A, E) = \IHom_{\underline{\mathcal A}}((A,\mm_A)^* \widetilde\mm'_A, \mathcal A[S]).
\end{align*}
Since $\mathcal A[S]$ lies in $\D(\mathcal A)$, the same is true for the above $\IHom$ (see \cref{rslt:properties-of-D-of-analytic-ring}). This finishes the proof that $\varphi^*\mathcal A$ is an analytic ring (the required condition regarding the Frobenius on $\varphi^* \mathcal A$ follows easily from the corresponding property for $\mathcal A$).

Getting back to the construction of the desired $p_\varphi$-coCartesian edge $f\colon \mathcal A \to \mathcal A'$ in $\mathcal C$ over $\varphi$, we can now choose $f$ to be the morphism $f\colon \mathcal A \to \varphi^*\mathcal A$ induced by $f'\colon \underline{\mathcal A} \to \varphi^*\underline{\mathcal A}$. To see that $f$ is indeed a morphism in $\mathcal C$, i.e. that $\varphi_*\colon \D(\varphi^*\mathcal A) \to \D(\underline{\mathcal A})$ factors over $\D(\mathcal A)$, we note that it follows formally from adjunctions that for every $M \in \D(\varphi^* \mathcal A)$ we have
\begin{align*}
	\IHom_{\mathcal A}(\mathcal A[S], \varphi^* M) = \varphi_* \IHom_{\varphi^*\mathcal A}(\varphi^*\mathcal A[S], M) = \varphi_* \IHom(\varphi^* \underline{\mathcal A}[S], M) = \IHom_{\mathcal A}(\underline{\mathcal A}[S], M),
\end{align*}
as desired. We observe that even something stronger is true: For $M \in \D(\varphi^*\underline{\mathcal A})$ we have $\varphi_* M \in \D(\mathcal A)$ if \emph{and only if} $M \in \D(\varphi^* \mathcal A)$ (reverse the above computation and use that $\varphi_*$ is conservative).

We have now constructed the morphism $f\colon \mathcal A \to \varphi^* \mathcal A$ in $\mathcal C$. To show that $f$ is $p_\varphi$-coCartesian we employ \cite[Proposition 2.4.4.3]{lurie-higher-topos-theory} which reduces the claim to showing that for any analytic ring $\mathcal B'$ over $(V',\mm')$ the natural map
\begin{align*}
	\Hom_{\mathcal C}(\varphi^*\mathcal A, \mathcal B') \isoto \Hom_{\mathcal C}(\mathcal A, \mathcal B')
\end{align*}
induced by $f$ is an isomorphism of anima. We know that the similar statement is true for the underlying almost rings because $f'$ is coCartesian. Therefore the claim reduces to the following: Let $g\colon \varphi^*\underline{\mathcal A} \to \underline{\mathcal B'}$ be a morphism of $(V',\mm')$-algebras such that the composition $\underline{\mathcal A} \xto{f'} \varphi^* \underline{\mathcal A} \xto{g} \underline{\mathcal B'}$ lies in $\mathcal C$; then $g$ lies in $\mathcal C$. But this follows immediately from the above observation that for given $M \in \D(\varphi^* \underline{\mathcal A})$ we have $\varphi_* \in \D(\mathcal A)$ if and only if $M \in \D(\varphi^*\mathcal A)$. This finishes the proof that $p$ is a locally coCartesian fibration.

Finally, to see that $p$ is a coCartesian fibration we use \cite[Remark 2.4.2.9]{lurie-higher-topos-theory} which reduces the claim to showing that for every composition $\varphi = \varphi' \comp \varphi''\colon (V,\mm) \to (V',\mm')$ in $\AlmSetup$, the natural morphism $\varphi^* \isoto \varphi'^* \comp \varphi''^*$ (induced by $p$) is an isomorphism of functors $\AnRing_{(V,\mm)} \to \AnRing_{(V',\mm')}$. But this is clear because all these functors are essentially computed on the underlying algebras and modules, where the commutativity holds by \cref{rslt:functoriality-of-almost-rings,rslt:functoriality-of-almost-modules-over-rings}. We have thus constructed the claimed functor $\AlmSetup \to \infcatinf$ mapping $(V,\mm) \mapsto \AnRing_{(V,\mm)}$. The claims (i) and (ii) follow from the construction.
\end{proof}

\begin{definition}
For every morphism $\varphi\colon (V,\mm) \to (V',\mm')$ of almost setups we denote
\begin{align*}
	&\varphi^*\colon \AnAssRing_{(V,\mm)} \to \AnAssRing_{(V',\mm')},\\
	&\varphi^*\colon \AnRing_{(V,\mm)} \to \AnRing_{(V',\mm')}
\end{align*}
the functors constructed in \cref{rslt:functoriality-of-analytic-rings-over-AlmSetup}. In the case that $\varphi$ is of the form $(V, V) \to (V,\mm)$ we also denote $\varphi^* = (V,\mm)^*$. If $\mm$ is clear from context then we abbreviate $(V,\mm)^*$ by $\mathcal A \mapsto \mathcal A^a$.
\end{definition}

\begin{proposition} \label{rslt:functoriality-of-modules-over-analytic-rings-over-AlmSetup}
Let $\AnAssRing_{(-)} \to \AlmSetup$ and $\AnRing_{(-)} \to \AlmSetup$ denote the coCartesian fibrations classifying the functors from \cref{rslt:functoriality-of-analytic-rings-over-AlmSetup}. Then there are natural functors
\begin{align*}
	\AnAssRing_{(-)} \to \infcatinf, \qquad &\mathcal A \mapsto \D(\mathcal A),\\
	\AnRing_{(-)} \to \infcatinf^\tensor, \qquad &\mathcal A \mapsto \D(\mathcal A).
\end{align*}
Moreover, for every morphism $\varphi\colon (V,\mm) \to (V',\mm')$ of almost setups and every analytic (associative) ring $\mathcal A$ over $(V,\mm)$, the induced functor $\varphi^*\colon \D(\mathcal A) \to \D(\varphi^* \mathcal A)$ satisfies the following properties:
\begin{propenum}
	\item $\varphi^*$ is the restriction of the functor $\varphi^*\colon \D(\underline{\mathcal A}) \to \D(\varphi^*\underline{\mathcal A})$ to $\D(\mathcal A)$.

	\item Suppose $\mathcal A$ is a (commutative) analytic ring. Then $\varphi^*$ preserves nuclear modules.
\end{propenum}
\end{proposition}
\begin{proof}
We only treat the case of (commutative) rings; the case of associative rings is similar. Composing the forgetful functor $\AnRing_{(-)} \to \Ring(-)$, $\mathcal A \mapsto \underline{\mathcal A}$ with the functor $\Ring(-) \to \infcatinf^\tensor$ from \cref{rslt:functoriality-of-almost-modules-over-rings} we obtain a functor $\AnRing_{(-)} \to \infcatinf^\tensor$ mapping $\mathcal A \mapsto \D(\underline{\mathcal A})$. Let $\mathcal C'^\tensor \to \AnRing_{(-)} \cprod \opComm$ denote the corresponding coCartesian family of symmetric monoidal $\infty$-categories and let $\mathcal C^\tensor \subset \mathcal C'^\tensor$ denote the full subcategory where we only allow objects in $\D(\mathcal A) \subset \D(\underline{\mathcal A})$ in the fiber over any $\mathcal A \in \AnRing_{(-)}$. It follows from the existence of symmetric monoidal left adjoints $\D(\underline{\mathcal A}) \to \D(\mathcal A)$ to the inclusions $\D(\mathcal A) \injto \D(\underline{\mathcal A})$ that $p\colon \mathcal C^\tensor \to \AnRing_{(-)} \cprod \opComm$ is a coCartesian family of symmetric monoidal $\infty$-operads and thus classifies the desired functor.

It remains to prove the claims (i) and (ii), so let $\varphi\colon (V,\mm) \to (V',\mm')$ and $\mathcal A$ be given. To prove (i) we have to show that given any $M \in \D(\mathcal A)$ the object $\varphi^* M \in \D(\varphi^* \underline{\mathcal A})$ lies in $\D(\varphi^*\mathcal A)$. But this is clear: We can write $M$ as a colimit of copies of $\mathcal A[S]$ for varying extremally disconnected sets $S$; then $\varphi^* M$ is a colimit of copies of $\varphi^*\mathcal A[S]$ and hence lies in $\D(\varphi^* \mathcal A)$. To prove (ii) we can argue separately in the cases that $\varphi$ is strict and that it is localizing. In the former case the claim follows from \cref{rslt:nuclear-modules-stable-under-pullback}. In the latter case $\varphi^*$ is an almost localization which commutes with $\tensor$ and $\IHom$ by \cref{rslt:properties-of-almost-localization-over-alg}, so preservation of nuclearity is formal.
\end{proof}

\begin{definition} \label{def:almost-pullback-on-analytic-modules}
For every morphism $\varphi\colon (V,\mm) \to (V',\mm')$ and every analytic ring $\mathcal A$ over $(V,\mm)$ we denote $\varphi^*\colon \D(\mathcal A) \to \D(\varphi^*\mathcal A)$ the symmetric monoidal functor constructed in \cref{rslt:functoriality-of-modules-over-analytic-rings-over-AlmSetup}. It admits a right adjoint $\varphi_*\colon \D(\varphi^*\mathcal A) \to \D(\mathcal A)$ which is the restriction of $\varphi_*\colon \D(\varphi^*\underline{\mathcal A}) \to \D(\underline{\mathcal A})$ to $\D(\varphi^*\mathcal A)$.
\end{definition}

For the convenience of the reader we formulate the following straightforward generalization of \cref{rslt:general-properties-of-almost-localization-over-ass-alg}.

\begin{proposition} \label{rslt:general-properties-of-almost-localization-over-analytic-rings}
Let $\mathcal A$ be an analytic ring over $V$.
\begin{propenum}
	\item Denote $\mm_{\mathcal A} := \mm_{\underline{\mathcal A}}$ and $\widetilde\mm_{\mathcal A} := \widetilde\mm_{\underline{\mathcal A}}$. Then $\mm_{\mathcal A}, \widetilde\mm_{\mathcal A} \in \D(\mathcal A)$.

	\item \label{rslt:construction-of-almost-localization-of-modules-over-analytic-ring} There is a natural localization functor
	\begin{align*}
		(\mathcal A,\mm_{\mathcal A})^*\colon \D(\mathcal A) \to \D(\mathcal A^a), \qquad M \mapsto M^a,
	\end{align*}
	which reduces to $(V,\mm)^*$ on the underlying $V$-modules. The functor $(\mathcal A,\mm_{\mathcal A})^*$ is $t$-exact, symmetric monoidal and preserves all limits and colimits. Moreover, for all $M, N \in \D(\mathcal A)$ there is a natural isomorphism
	\begin{align*}
		\IHom_{\mathcal A^a}(M^a, N^a) = \IHom_{\mathcal A}(M, N)^a.
	\end{align*}

	\item The functor $(\mathcal A,\mm_{\mathcal A})^*$ admits a left $t$-exact fully faithful right adjoint
	\begin{align*}
		(\mathcal A,\mm_{\mathcal A})_*\colon \D(\mathcal A^a) \to \D(\mathcal A), \qquad M \mapsto M_*.
	\end{align*}
	Moreover, for all $M \in \D(\mathcal A)$ there is a natural isomorphism
	\begin{align*}
		(M^a)_* = \IHom_{\mathcal A}(\widetilde\mm_{\mathcal A}, M).
	\end{align*}

	\item \label{rslt:left-adjoint-of-almost-localization-over-analytic-ring} The functor $(\mathcal A,\mm_{\mathcal A})^*$ admits a $t$-exact fully faithful left adjoint
	\begin{align*}
		(\mathcal A,\mm_{\mathcal A})_!\colon \D(\mathcal A^a) \to \D(\mathcal A), \qquad M \mapsto M_!,
	\end{align*}
	It can be computed as
	\begin{align*}
		M_! = \widetilde\mm_{\mathcal A} \tensor_{\mathcal A} M_*.
	\end{align*}

	\item For every strong limit cardinal $\kappa$ such that $A \in \D(V)_\kappa$ the functor $(\mathcal A, \mm_{\mathcal A})^*$ restricts to a functor $\D(\mathcal A)_\kappa \to \D(\mathcal A^a)_\kappa$.
\end{propenum}
\end{proposition}
\begin{proof}
Part (i) follows from the fact that $\mm$ is discrete, so that $\mm_{\underline{\mathcal A}} = \mm \tensor_V \underline{\mathcal A}$ is a colimit of copies of $\mathcal A$. For part (ii) we let $(\mathcal A, \mm_{\mathcal A})^* := \varphi^*$, where $\varphi\colon (V, V) \to (V,\mm)$ is the obvious morphism of almost setups and $\varphi^*$ is the functor defined in \cref{def:almost-pullback-on-analytic-modules}. By \cref{rslt:functoriality-of-modules-over-analytic-rings-over-AlmSetup} $(\mathcal A, \mm_{\mathcal A})^*$ is just the restriction of the functor $(\underline{\mathcal A}, \mm_{\underline{\mathcal A}})^*$, so the claimed properties follow from \cref{rslt:properties-of-almost-localization-over-alg} and the fact that $\IHom_{\mathcal A} = \IHom_{\underline{\mathcal A}}$. Part (iii) follows similarly with $(\mathcal A, \mm_{\mathcal A})_* = \varphi_*$. For part (iv) we note that the functor $(\underline{\mathcal A}, \mm_{\underline{\mathcal A}})_!\colon \D(\underline{\mathcal A}^a) \to \D(\underline{\mathcal A})$ maps $\D(\mathcal A^a)$ to $\D(\mathcal A)$: this follows from the formula $(\underline{\mathcal A}, \mm_{\underline{\mathcal A}})_! M = \widetilde\mm_{\underline{\mathcal A}} \tensor_{\underline{\mathcal A}} M_*$, noting that since $\widetilde\mm_{\underline{\mathcal A}}$ is a colimit of copies of $\mathcal A$, the tensor product $\widetilde\mm_{\underline{\mathcal A}} \tensor_{\underline{\mathcal A}} -$ preserves $\D(\mathcal A)$. It follows that we can (and have to) define $(\mathcal A,\mm_{\mathcal A})_!$ by restricting $(\underline{\mathcal A}, \mm_{\underline{\mathcal A}})_!$ to $\D(\mathcal A^a)$. Part (v) is obvious from the definitions.
\end{proof}

Given a morphism $\varphi\colon (V,\mm) \to (V',\mm')$ of almost setups and an analytic ring $\mathcal A$ over $(V,\mm)$, the induced map $\mathcal A \to \varphi^*\mathcal A$ in $\AnRing_{(-)}$ behaves much like a steady morphism. Namely, it satisfies arbitrary base-change:

\begin{proposition} \label{rslt:almost-localization-is-steady}
Let $\varphi\colon (V,\mm) \to (V',\mm')$ be a morphism of almost setups and let $f\colon \mathcal A \to \mathcal B$ be a map of analytic rings over $(V,\mm)$. Then the diagram
\begin{center}\begin{tikzcd}
	\D(\varphi^* \mathcal B) \arrow[d,swap,"(\varphi^*f)_*"] & \D(\mathcal B) \arrow[l,"\varphi^*"] \arrow[d,"f_*"]\\
	\D(\varphi^*\mathcal A) & \D(\mathcal A) \arrow[l,"\varphi^*"]
\end{tikzcd}\end{center}
is naturally coherent, i.e. the natural base-change morphism
\begin{align*}
	\varphi^* f_* \isoto (\varphi^*f)_* \varphi^*
\end{align*}
is an isomorphism of functors $\D(\mathcal B) \to \D(\varphi^*\mathcal A)$.
\end{proposition}
\begin{proof}
We can show this separately in the cases that $\varphi$ is strict and that it is localizing. If $\varphi$ is strict then it is just the map of analytic rings $V^a \to V'^a$ (with the canonical analytic structure) which is easily seen to be steady (e.g. by \cref{rslt:nuclear-induced-ring-is-steady}). If $\varphi$ is localizing the claim follows from the fact that $\varphi^*$ is computed by $(V,\mm)^*$ on the underlying $V$-modules by \cref{rslt:construction-of-almost-localization-of-modules-over-analytic-ring}.
\end{proof}

A natural question one might ask is whether the almost localization functor $\AnRing_V \to \AnRing_{(V,\mm)}$, $\mathcal A \mapsto \mathcal A^a$ admits a (fully faithful) right adjoint $\mathcal A \mapsto \mathcal A_{**}$. This seems to fail due to set-theoretic issues, but we can get close to it:

\begin{proposition} \label{rslt:almost-localization-of-AnRing-is-surjective}
The almost localization functor
\begin{align*}
	(V,\mm)^*\colon \AnRing_V \to \AnRing_{(V,\mm)}, \qquad \mathcal A \mapsto \mathcal A^a
\end{align*}
is essentially surjective and preserves all small colimits. More precisely, for every analytic ring $\mathcal A$ over $(V,\mm)$ and every $V$-algebra $A'$ with $A'^a = \underline{\mathcal A}$ there is an analytic ring structure $\mathcal A'$ on $A'$ such that $M \in \D(A')$ satisfies $M \in \D(\mathcal A')$ if and only if $M^a \in \D(\mathcal A)$.
\end{proposition}
\begin{proof}
It follows easily from the explicit computation of colimits in \cref{rslt:colimits-of-analytic-rings} that $(V,\mm)^*$ preserves all small colimits.

Now let $\mathcal A \in \AnRing_{(V,\mm)}$ and $A' \in \Ring_V$ with $A'^a = \underline{\mathcal A}$ be given. Let us denote $A := \underline{\mathcal A}$. We will construct the desired analytic ring structure on $A'$. To do this, fix any strong limit cardinal $\kappa$ such that $A'$ is $\kappa$-condensed (i.e. lies in $\D(V)_\kappa$). We define $\D_{\ge0}(\mathcal A)'_\kappa \subset \D_{\ge0}(A)_\kappa$ to be the full subcategory of those objects $M \in \D_{\ge0}(A)_\kappa$ such that for all extremally disconnected sets $S$ the natural map
\begin{align*}
	\tau_{\ge0} \IHom_A(\mathcal A[S]_\kappa, M)_\kappa \isoto \tau_{\ge0} \IHom_A(A[S], M)_\kappa
\end{align*}
is an isomorphism (i.e. $\D_{\ge0}(\mathcal A)'_\kappa$ is a $\kappa$-condensed version of $\D_{\ge0}(\mathcal A)$). By the same argument as in \cref{rslt:properties-of-D-of-analytic-ass-ring} $\D_{\ge0}(\mathcal A)'_\kappa$ is stable under small limits and colimits in $\D_{\ge0}(A)_\kappa$ and the inclusion $\D_{\ge0}(\mathcal A)'_\kappa \injto \D_{\ge0}(A)_\kappa$ has a left adjoint $- \tensor_A \mathcal A_\kappa$. Now define $\D_{\ge0}(\mathcal A')'_\kappa \subset \D_{\ge0}(A')_\kappa$ to be the full subcategory spanned by those $M \in \D_{\ge0}(A')_\kappa$ such that $M^a \in \D_{\ge0}(\mathcal A)'_\kappa$.

Clearly $\D_{\ge0}(\mathcal A')'_\kappa$ is stable under limits and colimits in $\D_{\ge0}(A')_\kappa$. We claim that it is presentable. To see this, by \cref{rslt:compact-generators-implies-compactly-generated} it is enough to find a small family $(X_i)_{i\in I}$ of objects $X_i \in \D_{\ge0}(\mathcal A')'_\kappa$ with the following property: Given any $M \in \D_{\ge0}(\mathcal A')'_\kappa$ such that $\Hom(X_i, M) = *$ for all $i$ then $M = 0$. We choose the family $(X_i)_i$ to consist of the objects
\begin{align*}
	(\mathcal A[S]_\kappa)_![n], \qquad \cofib(\widetilde\mm_{A'} \to A') \tensor_{A'} A'[S][n]
\end{align*}
for $n \ge 0$ and $\kappa$-small extremally disconnected sets $S$. Now let any $M \in \D_{\ge0}(\mathcal A')'_\kappa$ be given and suppose that $\Hom(X_i, M) = *$ for all $X_i$. Taking $X_i = (\mathcal A[S]_\kappa)_![n]$ we see
\begin{align*}
	&* = \Hom((\mathcal A[S]_\kappa)_![n], M) = \Hom((\mathcal A[S]_\kappa)_!, \tau_{\ge0}(M[-n])) =\\&\qquad= \Hom(\mathcal A[S]_\kappa, \tau_{\ge0}(M[-n])^a) = \Hom(A[S], \tau_{\ge0}(M[-n])^a)
\end{align*}
Ranging over all $S$ we deduce $\tau_{\ge0}(M[-n])^a = 0$. Then ranging over all $n$ we deduce $M^a = 0$. Thus by \cref{rslt:computation-of-almost-lower-star-over-alg} we have $\IHom_{A'}(\widetilde\mm_{A'}, M) = 0$. It follows that
\begin{align*}
	* = \Hom(\cofib(\widetilde\mm_{A'} \to A') \tensor_{A'} A'[S][n], M) = \Hom(A'[S][n], M)
\end{align*}
for all $\kappa$-small $S$ and $n \ge 0$, hence $M = 0$ as desired. We have shown that $\D_{\ge0}(\mathcal A')'_\kappa$ is presentable. By the adjoint functor theorem the inclusion $\D_{\ge0}(\mathcal A')'_\kappa \injto \D_{\ge0}(A')_\kappa$ admits a left adjoint $L_\kappa$.

Now fix any $M \in \D_{\ge0}(A')$. We claim that for large enough $\kappa$ the value $L_\kappa(M) \in \D_{\ge0}(A')$ does not depend on $\kappa$. Since $L$ preserves all small colimits this reduces to the case $M = A'[S]$ for some fixed extremally disconnected set $S$. Choose a strong limit cardinal $\kappa_0$ such that $\mathcal A[S]$ lies in $\D(A)_{\kappa_0}$ and pick any strong limit cardinal $\kappa \ge \kappa_0$. Then $L_{\kappa_0}(A'[S])^a = \mathcal A[S] = \mathcal A[S]_\kappa$ and thus $L_{\kappa_0}(A'[S])$ lies in $\D_{\ge0}(\mathcal A')'_\kappa$. In particular we get a natural map $L_\kappa(A'[S]) \to L_{\kappa_0}(A'[S])$ which we claim to be an isomorphism. To see this, pick any $N \in \D_{\ge0}(\mathcal A_0)'_\kappa$. Let $N' := (N_{\kappa_0}^a)_* \in \D_{\ge0}(A')_{\kappa_0}$. Then
\begin{align*}
	&\Hom(L_{\kappa_0}(A'[S]), N') = \Hom(L_{\kappa_0}(A'[S])^a, N_{\kappa_0}^a) = \Hom(\mathcal A[S], N_{\kappa_0}^a) = \Hom(A[S], N_{\kappa_0}^a) =\\&\qquad= \Hom(A'[S], N').
\end{align*}
Now let $N'' := \cofib(N_{\kappa_0} \to N')$. Then $N''^a = 0$ and therefore $N'' \in \D_{\ge0}(\mathcal A')'_{\kappa_0}$. It follows that
\begin{align*}
	\Hom(L_{\kappa_0}(A'[S]), N'') = \Hom(A'[S], N'').
\end{align*}
The last two computations combine to
\begin{align*}
	&\Hom(L_{\kappa_0}(A'[S]), N) = \Hom(L_{\kappa_0}(A'[S]), N_{\kappa_0}) = \Hom(A'[S], N_{\kappa_0}) =\\&\qquad= \Hom(A'[S], N) = \Hom(L_\kappa(A'[S]), N)
\end{align*}
and thus $L_{\kappa_0}(A'[S]) = L_\kappa(A'[S])$, as desired. Now let $\D_{\ge0}(\mathcal A') \subset \D_{\ge0}(A')$ denote the full subcategory spanned by those objects $M \in \D_{\ge0}(A')$ with $M^a \in \D_{\ge0}(\mathcal A)$. Note that every $M \in \D_{\ge0}(\mathcal A')$ lies in $\D_{\ge0}(\mathcal A')'_\kappa$ for large enough $\kappa$. Since $L_\kappa(M)$ is independent of $\kappa$ for large $\kappa$ (as shown above) we deduce that the functors $L_\kappa$ induce a functor
\begin{align*}
	L\colon \D_{\ge0}(A') \to \D_{\ge0}(\mathcal A')
\end{align*}
which is left adjoint to the inclusion $\D_{\ge0}(\mathcal A') \injto \D_{\ge0}(A')$. Note moreover that for every extremally disconnected set $S$ the $\infty$-category $\D_{\ge0}(\mathcal A')$ is stable under the functor $\IHom_{A'}(A'[S], -)$ (this follows immediately from \cref{rslt:computation-of-almost-lower-star-over-alg}). It follows formally that the functor $S \mapsto \mathcal A'[S] := L(A'[S])$ defines an analytic ring structure on $A'$ (cf. \cite[Proposition 12.20]{scholze-analytic-spaces} and also note that the Frobenius condition for $\mathcal A_0$ follows immediately from the Frobenius condition for $\mathcal A$). \end{proof}

\begin{remark}
From \cref{rslt:almost-localization-of-AnRing-is-surjective} it should follow that the almost localization functor $\mathcal A \mapsto \mathcal A^a$ of analytic rings admits a left adjoint $\mathcal A \mapsto \mathcal A_{!!}$. However, it does not seem to admit a right adjoint $\mathcal A \mapsto \mathcal A_{**}$ due to set-theoretic issues: One needs to define $\mathcal A_{**}$ as the analytic ring structure on $(\underline{\mathcal A})_{**}$ which is obtained as the colimit $\mathcal A_{**} = \varinjlim_i \mathcal A_i$, where $\mathcal A_i$ ranges over all analytic ring structures on $(\underline{\mathcal A})_{**}$ which satisfy $\mathcal A_i^a = \mathcal A$. This colimit is filtered but usually not small, so that its existence is unclear.
\end{remark}

From the previous result we can deduce that steadiness is preserved under a change of almost setup, as the next result shows. This is crucial in order to define almost localizations of analytic spaces later on.

\begin{proposition} \label{rslt:almost-localization-preserves-steadiness}
Let $\varphi\colon (V,\mm) \to (V',\mm')$ be a morphism of almost setups. Then the functor $\varphi^*\colon \AnRing_{(V,\mm)} \to \AnRing_{(V',\mm')}$ preserves all small colimits and steady morphisms.
\end{proposition}
\begin{proof}
If $\varphi$ is strict then $\varphi^* = - \tensor_{V^a} V'^a$, so the claim follows from \cref{rslt:stability-of-steady-maps}. We can thus assume that $\varphi$ is localizing, i.e. $V' = V$ and $\mm' \subset \mm$. Preservation of colimits is shown in \cref{rslt:almost-localization-of-AnRing-is-surjective}, so we only need to show preservation of steadiness. Let $f\colon \mathcal A \to \mathcal B$ be a steady morphism of analytic rings over $(V,\mm)$ and let $\varphi^*\mathcal A \to \mathcal C$ be a morphism of analytic rings over $(V,\mm')$. Let $\mathcal B'$ and $\mathcal C'$ denote the analytic ring structures on $\underline{\mathcal B}$ and $(\underline{\mathcal C})_*$ constructed in \cref{rslt:almost-localization-of-AnRing-is-surjective}. Then there are natural maps $\mathcal B' \to \mathcal B$ and $\mathcal B' \to \mathcal C'$, so we can define $\mathcal C'' := \mathcal B \tensor_{\mathcal B'} \mathcal C'$. Then $\varphi^* \mathcal C'' = \mathcal C$ and the natural map $\mathcal B \to \mathcal C''$ gets transformed to the map $\varphi^*\mathcal B \to \mathcal C$ by applying $\varphi^*$. Then the base-change criterion (see \cref{rslt:steady-base-change-holds}) for $\varphi^*\mathcal A \to \varphi^* \mathcal B$ along $\varphi^* \mathcal B \to \mathcal C$ follows from the base-change criterion for $\mathcal A \to \mathcal B$ along $\mathcal B \to \mathcal C''$ by applying $\varphi^*$ everywhere and using that $\varphi^*$ commutes with pullbacks (by \cref{rslt:functoriality-of-modules-over-analytic-rings-over-AlmSetup}) and pushforwards (by \cref{rslt:almost-localization-is-steady}).
\end{proof}

For later use, let us record the following definition of almost bounded objects in $\D(\mathcal A)$ for an analytic ring $\mathcal A$ over $(V,\mm)$. Note that $\D(\mathcal A)$ is $V^a$-enriched, so that by \cref{def:V-a-enriched-categories} we have a notion of $\approx_\varepsilon$ and of $\varepsilon$-retracts in $\D(\mathcal A)$.

\begin{definition} \label{def:almost-bounded-objects-in-D-A}
Let $\mathcal A$ be an analytic ring over $(V, \mm)$. An object $M \in \D(\mathcal A)$ is called \emph{weakly almost (left/right) bounded} if for every $\varepsilon \in \mm$ there is some (left/right) bounded object $M_\varepsilon \in \D(\mathcal A)$ such that $M$ is an $\varepsilon$-retract of $M_\varepsilon$. We denote by
\begin{align*}
	\D^{w+}(\mathcal A), \D^{w-}(\mathcal A), \D^{wb}(\mathcal A)
\end{align*}
the respective full subcategories.
\end{definition}

\subsection{Analytic Geometry} \label{sec:andesc.anspace}

Fix an almost setup $(V,\mm)$. In the previous subsection we have defined and studied analytic rings over $(V,\mm)$. In the present subsection we will glue these analytic rings to get more general analytic spaces over $(V,\mm)$, much like classical schemes are obtained by gluing classical rings. The following constructions are analogous to \cite[\S13]{scholze-analytic-spaces} but we work in the more general setting of almost mathematics.

The gluing procedure is rather formal once we have defined the coordinate rings and the open immersions. The coordinate rings for analytic geometry are the analytic rings, and the open immersions will be steady localizations. We now introduce the latter.

\begin{definition}
A map $f\colon \mathcal A \to \mathcal B$ of analytic rings over $(V,\mm)$ is called a \emph{localization} if the forgetful functor $f_*\colon \D(\mathcal B) \to \D(\mathcal A)$ is fully faithful.
\end{definition}

\begin{proposition} \label{rslt:properties-of-localizations}
\begin{propenum}
	\item Localizations of analytic rings over $(V,\mm)$ are stable under composition and base-change. Steady localizations are additionally stable under all colimits.

	\item Let $f\colon \mathcal A \to \mathcal B$ and $g\colon \mathcal B \to \mathcal C$ be maps of analytic rings over $(V,\mm)$. If $f$ and $g \comp f$ are (steady) localizations, then so is $g$.

	\item \label{rslt:steady-localization-equiv-diagonal-is-isom} Let $f\colon \mathcal A \to \mathcal B$ be a map of analytic rings over $(V,\mm)$. If $f$ is a localization then the induced map $\mathcal B \isoto \mathcal B \tensor_{\mathcal A} \mathcal B$ is an isomorphism. Conversely, if $f$ is steady and $\mathcal B \to \mathcal B \tensor_{\mathcal A} \mathcal B$ is an isomorphism, then $f$ is a localization.
\end{propenum}
\end{proposition}
\begin{proof}
Part (ii) is obvious (using \cref{rslt:steadyness-satisfies-2-out-of-3} for the claim about steadiness). It is also clear that localizations are stable under composition. Next we show stability under base-change, so let $f\colon \mathcal A \to \mathcal B$ be any localization, $\mathcal A \to \mathcal A'$ any map and $\mathcal B' = \mathcal B \tensor_{\mathcal A} \mathcal A'$. Note that $\mathcal B$ is the completion of the uncompleted analytic ring $\tilde{\mathcal B}$ with $\underline{\tilde{\mathcal B}} = \underline{\mathcal A}$ and $\tilde{\mathcal B}[S] = \mathcal B[S]$. Since completion preserves all colimits (and in particular pushouts), $\mathcal B'$ is the completion of the uncompleted pushout $\tilde{\mathcal B'} := \tilde{\mathcal B} \tensor_{\mathcal A} \mathcal A'$. We have $\underline{\tilde{\mathcal B'}} = \underline{\mathcal A} \tensor_{\underline{\mathcal A}} \underline{\mathcal A'} = \underline{\mathcal A'}$, so all in all we deduce that $\mathcal B'$ is the completion of an uncompleted analytic ring with underlying $(V,\mm)$-algebra $\underline{\mathcal A'}$. It follows immediately that $\mathcal A' \to \mathcal B'$ is a localization. This proves the first statement in (i).

We now prove (iii), so let $f\colon \mathcal A \to \mathcal B$ be a map of analytic rings over $(V,\mm)$. Suppose first that $f$ is a localization. By the previous paragraph, $\mathcal B' := \mathcal B \tensor_{\mathcal A} \mathcal B$ is the completion of an uncompleted analytic ring $\tilde{\mathcal B'}$ with $\underline{\tilde{\mathcal B'}} = \underline{\mathcal B}$. By the explicit formula in \cref{rslt:pushouts-of-analytic-rings} we deduce $\tilde{\mathcal B'}[S] = \mathcal B[S]$ for all extremally disconnected sets $S$ (use that $\mathcal B[S] \tensor_{\mathcal A} \mathcal B = \mathcal B[S]$ by full faithfulness of the forgetful functor $\D(\mathcal B) \to \D(\mathcal A)$; this object lies in $\D(\mathcal B)$ under both forgetful functors $\D(\underline{\mathcal B} \tensor_{\underline{\mathcal A}} \underline{\mathcal B}) \to \D(\underline{\mathcal B})$, hence further applications of $- \tensor_{\underline{\mathcal B}} \mathcal B$ have no effect). Thus $\mathcal B' = \mathcal B$. Conversely, assume that $f$ is steady and $\mathcal B \isoto \mathcal B \tensor_{\mathcal A} \mathcal B$ is an isomorphism. Then for all $M \in \D(\mathcal B)$ we have, using steadiness (see \cref{rslt:steady-map-of-analytic-rings-equiv-base-change}),
\begin{align*}
	M \tensor_{\mathcal A} \mathcal B = M \tensor_{\mathcal B} (\mathcal B \tensor_{\mathcal A} \mathcal B) = M \tensor_{\mathcal B} \mathcal B = M,
\end{align*}
hence $\D(\mathcal B) \to \D(\mathcal A)$ is fully faithful, as desired. This proves (iii).

It remains to prove the second part of (i), i.e. that steady localizations are stable under all colimits. Let $(\mathcal A \to \mathcal B) = \varinjlim_i (\mathcal A_i \to \mathcal B_i)$ be a colimit of steady localizations. By (iii), for each $i$ the map $\mathcal B_i \isoto \mathcal B_i \tensor_{\mathcal A_i} \mathcal B_i$ is an isomorphism. Passing to the colimit we deduce that $\mathcal B \isoto \mathcal B \tensor_{\mathcal A} \mathcal B$ is an isomorphism, hence $\mathcal A \to \mathcal B$ is a steady localization by (iii) (here we also use that steady maps are stable under colimits, see \cref{rslt:stability-properties-of-steadiness}).
\end{proof}

We now come to the desired gluing procedure. Instead of directly defining analytic spaces over $(V,\mm)$, let us abstract the process a bit further by providing a general recipe of how to glue spaces from rings. This will avoid some repetitiveness when we later introduce schemes and discrete adic spaces in the setting of analytic geometry.

\begin{definition} \label{def:geometry-blueprint}
A \emph{geometry blueprint} is a pair $G = (R_G, L_G)$, where $R_G$ is a full subcategory of $\AnRing_{(V,\mm)}$ and $L_G$ is a class of morphisms in $R_G$ satisfying the following properties:
\begin{enumerate}[(i)]
	\item $L_G$ contains all equivalences and is stable under composition. In particular, if two morphisms $f$ and $f'$ in $R_G$ are homotopic, then $f \in L_G$ iff $f' \in L_G$.

	\item $R_G$ is stable under pushouts in $\AnRing_{(V,\mm)}$ and $L_G$ is stable under base-change.
	\item $R_G$ is stable under finite products in $\AnRing_{(V,\mm)}$. In particular $R_G$ contains the final object $0$.

	\item Every morphism in $L_G$ is a steady localization.
\end{enumerate}
If $G$ is a geometry blueprint over $(V,\mm)$, we call the objects of $R_G$ the \emph{$G$-analytic rings} and the objects of $L_G$ the \emph{$G$-localizations}.
\end{definition}

\begin{remark}
Let $G = (R_G, L_G)$ be a geometry blueprint over $(V,\mm)$. Then $L_G$ satisfies the following additional property: Let $f\colon \mathcal A \to \mathcal B$ and $g\colon \mathcal B \to \mathcal C$ be morphisms in $R_G$ with composition $h$ and assume that $h$ and $f$ lie in $L_G$. Then $g$ is in $L_G$. Namely, by \cref{rslt:steady-localization-equiv-diagonal-is-isom} we can write $g$ as the map $\mathcal B \isoto \mathcal C \tensor_{\mathcal A} \mathcal B \to \mathcal B$, where the second map is the base-change of $h$.
\end{remark}

\begin{definition}
Let $G = (R_G, L_G)$ be a geometry blueprint over $(V,\mm)$.
\begin{defenum}
	\item The objects of $R_G^\opp$ are denoted by $\AnSpec_G \mathcal A$ for $\mathcal A \in R_G$. The \emph{$G$-analytic site} on $R_G^\opp$ is the Grothendieck site whose coverings are generated by finite families of $G$-localizations $(\AnSpec_G \mathcal A_i \to \AnSpec_G \mathcal A)_i$ such that the pullback functor $\D(\mathcal A) \to \prod_i \D(\mathcal A_i)$ is conservative.

	\item We denote $\Shv(R_G^\opp)$ the $\infty$-category of sheaves (of anima) on the site $R_G^\opp$, i.e. the $\infty$-category of functors from $R_G$ to anima which satisfy descent along analytic Čech covers (see \cref{rslt:can-apply-cech-descent-to-G-analytic-site}).
\end{defenum}
\end{definition}

\begin{remarks}
\begin{remarksenum}
	\item Given a geometry blueprint $G$ over $(V,\mm)$, the $\infty$-category $R_G^\opp$ is usually not small, so the category of sheaves on it is not a topos (it is not presentable). In particular, one needs to be careful when dealing with these sheaves, for example sheafification may not exist in general (it may happen that the sheafification of a presheaf requires ``large'' anima as values). Nevertheless, many of the good properties of topoi are retained: For example there is a good notion of injective (monomorphisms) and surjective (effective epimorphisms) maps of sheaves -- to see this, one can embedd the category $\Shv(R_G^\opp)$ into the category of sheaves of large anima, which is a topos (a priori this argument requires the existence of Grothendieck universes, but this should be avoidable).

	\item \label{rslt:can-apply-cech-descent-to-G-analytic-site} For a geometry blueprint $G$ over $(V,\mm)$ the analytic site on $R_G^\opp$ is not directly an explicit covering site (see \cref{def:explicit-covering-site}) but very close to it. In particular it is not hard to see that the proof of \cref{rslt:sheaves-on-explicit-covering-site} still applies, so that sheaves on $R_G^\opp$ can be identified with functors satisfying Čech descent.
\end{remarksenum}
\end{remarks}

\begin{lemma}
Let $G = (R_G, L_G)$ be a geometry blueprint over $(V,\mm)$.
\begin{lemenum}
	\item For every analytic ring $\mathcal A \in R_G$ the functor
	\begin{align*}
		\AnSpec_G \mathcal A\colon R_G \to \Ani, \qquad \mathcal B \mapsto \Hom(\mathcal A, \mathcal B)
	\end{align*}
	is a sheaf for the $G$-analytic site. In particular, we can naturally view $R_G^\opp$ as a full subcategory of $\Shv(R_G^\opp)$.

	\item \label{rslt:qcoh-sheaves-is-sheaf-of-infty-categories-on-affine-analytic-spaces} The functor
	\begin{align*}
		\D\colon R_G \to \infcatinf^\tensor, \qquad \mathcal A \mapsto \D(\mathcal A)
	\end{align*}
	is a sheaf for the $G$-analytic site.
\end{lemenum}
\end{lemma}
\begin{proof}
We first prove (ii). Let $(\AnSpec_G \mathcal B_i \to \AnSpec_G \mathcal B)_{i\in I}$ be a $G$-analytic covering and for every non-empty subset $J \subset I$ denote $\mathcal B_J := \bigtensor_{j\in J} \mathcal B_j$ (via $\AnSpec_G$, this corresponds to the fiber product of the $\AnSpec_G \mathcal B_j$ over $\AnSpec_G \mathcal B$). We need to see that the natural functor
\begin{align*}
	\D(\mathcal B) \isoto \varprojlim_{J \subset I} \D(\mathcal B_J)
\end{align*}
is an isomorphism. This follows formally from the fact that base-change holds along steady localizations (see \cref{rslt:steady-map-of-analytic-rings-equiv-base-change}), see the proof of \cite[Proposition 10.5]{condensed-mathematics}.

We now prove (i), so let $\mathcal A$ be given as in the claim and pick $\mathcal B$ and $\mathcal B_i$ as before. We need to see that the natural map
\begin{align*}
	\Hom(\mathcal A, \mathcal B) \isoto \varprojlim_{J \subset I} \Hom(\mathcal A, \mathcal B_J)
\end{align*}
is an isomorphism. This boils down to showing that $\mathcal B \isoto \varprojlim_J \mathcal B_J$ is an isomorphism of analytic rings, which in turn reduces to showing that $\mathcal B[S] \isoto \varprojlim_J \mathcal B_J[S]$ is an isomorphism in $\D(\mathcal B)$ for all extremally disconnected sets $S$ (setting $S = *$ then shows that $\underline{\mathcal B} \isoto \varprojlim_J \underline{\mathcal B_J}$ as $(V,\mm)$-algebras, because the forgetful functor from algebras to modules is conservative). But this follows directly from (ii). The second part of the claim is an instance of the Yoneda lemma.
\end{proof}

We now introduce the gluing procedure. For a given geometry blueprint $G$ this will result in an $\infty$-category of $G$-analytic spaces, which are glued from the affine spaces in $R_G^\opp$ along $G$-localizations.

\begin{definition}
Let $G = (R_G, L_G)$ be a geometry blueprint over $(V,\mm)$. An \emph{affine $G$-analytic space} is a sheaf $X \in \Shv(R_G^\opp)$ which is isomorphic to $\AnSpec_G \mathcal A$ for some analytic ring $\mathcal A \in R_G$. \end{definition}

\begin{definition}
Let $G = (R_G, L_G)$ be a geometry blueprint over $(V,\mm)$.
\begin{defenum}
	\item Let $X = \AnSpec_G \mathcal A$ be an affine $G$-analytic space.
		A \emph{$G$-analytic subspace} of $X$ is a subsheaf $U \subset X$ which admits a cover $\bigdunion_i U_i \surjto U$ by a small family of affine subsheaves $U_i = \AnSpec_G \mathcal B_i \subset U$ such that the induced maps $\mathcal A \to \mathcal B_i$ are $G$-localizations. Note that by \cref{rslt:steady-localization-equiv-diagonal-is-isom} every $G$-localization gives rise to a subsheaf and hence to an affine $G$-analytic subspace.

	\item Let $X$ be a sheaf on $R_G^\opp$. A \emph{$G$-analytic subspace} of $X$ is a subsheaf $U \subset X$ such that for every map $Y \to X$ from an affine $G$-analytic space $Y$, the pullback $U \cprod_X Y \subset Y$ is a $G$-analytic subspace.
\end{defenum}
\end{definition}

\begin{definition}
Let $G = (R_G, L_G)$ be a geometry blueprint over $(V,\mm)$. A \emph{$G$-analytic space} is a sheaf $X$ on $R_G^\opp$ which admits a cover by a small family of affine $G$-analytic subspaces. We denote by
\begin{align*}
	\AnSpace_G \subset \Shv(R_G^\opp)
\end{align*}
the full subcategory spanned by the $G$-analytic spaces.
\end{definition}

As a special case of the above gluing procedure we obtain the $\infty$-category of all analytic spaces over $(V,\mm)$:

\begin{definition}
An \emph{analytic space over $(V,\mm)$} is a $G_0$-analytic space for the geometry blueprint $G_0 = (R_{G_0}, L_{G_0})$ where $R_{G_0} = \AnRing_{(V,\mm)}$ and $L_{G_0}$ is the class of all steady localizations. We write
\begin{align*}
	\AnSpace_{(V,\mm)}
\end{align*}
for the $\infty$-category of analytic spaces over $(V,\mm)$. We also abbreviate $\AnSpace := \AnSpace_{(\Z,\Z)}$, the $\infty$-category of analytic spaces. A \emph{steady subspace} of an analytic space over $(V,\mm)$ is a $G_0$-analytic subspace. We also write $\AnSpec \mathcal A = \AnSpec_{G_0} \mathcal A$ for every analytic ring $\mathcal A$ over $(V,\mm)$.
\end{definition}

Checking that a subsheaf $U \subset X$ of a sheaf $X$ on $R_G^\opp$ is a $G$-analytic subspace can sometimes be difficult due to the fact that it is a requirement on \emph{all} objects in $R_G^\opp$. In particular this can make it hard to verify that $X$ is a $G$-analytic space in general. To remedy this, we provide the following criterion for $X$ to be a $G$-analytic space. Roughly, it says that everything that can be glued from $G$-analytic spaces is again a $G$-analytic space.

\begin{lemma} \label{rslt:obtain-analytic-space-from-gluing}
Let $G = (R_G, L_G)$ be a geometry blueprint over $(V,\mm)$ and let $X \in \Shv(R_G^\opp)$. Suppose there is a cover $X = \bigunion_{i\in I} U_i$ of $X$ by subsheaves $U_i \subset X$ such that each $U_i$ is a $G$-analytic space and for all $i, j \in I$ the map $U_i \isect U_j \injto U_i$ is a $G$-analytic subspace. Then $X$ is a $G$-analytic space and the maps $U_i \injto X$ are $G$-analytic subspaces.
\end{lemma}
\begin{proof}
By refining the cover $(U_i)_i$ we can assume that all $U_i$ are affine. Then $X$ is a $G$-analytic space as soon as we can show that all $U_i \subset X$ are $G$-analytic subspaces. To verify that the latter condition is satisfied, pick any affine $G$-analytic space $Y = \AnSpec_G \mathcal B$ with a map $Y \to X$. We need to see that the subsheaves $V_i := U_i \cprod_X Y \subset Y$ are $G$-analytic subspaces. By \cref{rslt:cover-of-sheaves-characterized-by-sections} we can find a cover $(Y_j \to Y)_{j\in J}$ such that each $Y_j = \AnSpec_G \mathcal B_j$ for some $G$-localization $\mathcal B \to \mathcal B_j$ and such that each map $Y_j \to X$ factors over some $U_i$, so in particular $Y_j \subset V_i$. Since the $Y_j$'s form a cover of $Y$ by affine $G$-analytic subspaces, it is now enough to verify that for all $i \in I$ and $j \in J$ the subsheaf $V_i \isect Y_j \subset Y_j$ is a $G$-analytic subspace. Fix $i$ and $j$ and pick an index $i' \in I$ such that $Y_j \subset V_{i'}$. Then it is enough to show that $V_i \isect V_{i'} \subset V_{i'}$ is a $G$-analytic subspace. But this is the pullback of $U_i \isect U_{i'} \subset U_{i'}$ along $Y \to X$ and the latter is a $G$-analytic subspace by assumption.
\end{proof}

We next introduce some basic terminology regarding $G$-analytic spaces, like quasicompact and quasiseparated spaces and the associated symmetric monoidal $\infty$-category $\D(X)$ of quasicoherent sheaves.

\begin{definition}
Let $G = (R_G, L_G)$ be a geometry blueprint over $(V,\mm)$.
\begin{defenum}
	\item The \emph{$G$-analytic site} on $\AnSpace_G$ is the site whose coverings are small families $(U_i \injto X)_i$ of jointly surjective maps such that each $U_i \subset X$ is a $G$-analytic subspace.

	\item For every $G$-analytic space $X$ we denote $X_\an$ the site whose objects are the $G$-analytic subspaces of $X$ and whose covers are as in (a).

	\item A $G$-analytic space $X$ is called \emph{quasi-compact} (resp. \emph{quasi-separated}) if it is so as an object of the site $\AnSpace_G$. Concretely, $X$ is quasi-compact if every cover $(U_i \injto X)_i$ of $X$ by $G$-analytic subspaces admits a finite subcover, and $X$ is quasi-separated if for any two quasi-compact $G$-analytic spaces $Y_1, Y_2$ with maps $Y_i \to X$, the fiber product $Y_1 \cprod_X Y_2$ is quasi-compact.

	A map $f\colon Y \to X$ of $G$-analytic spaces is \emph{quasi-compact} (resp. \emph{quasi-separated}) if for every quasi-compact (resp. quasi-separated) $G$-analytic space $Z$ with a map $Z \to X$, the pullback $Z \cprod_X Y$ is quasi-compact (resp. quasi-separated).

	We abbreviate the property of being ``quasi-compact and quasi-separated'' by \emph{qcqs}.
\end{defenum}
\end{definition}

\begin{definition} \label{def:quasicoherent-sheaves-on-G-analytic-space}
Let $G$ be a geometry blueprint over $(V,\mm)$.
\begin{defenum}
	\item We define
	\begin{align*}
		\D\colon \AnSpace_G^\opp \to \infcatinf^\tensor, \qquad X \mapsto \D(X)
	\end{align*}
	to be the sheaf of symmetric monoidal $\infty$-categories on the analytic site of $\AnSpace_G$ whose restriction to affine $G$-analytic spaces is the sheaf in \cref{rslt:qcoh-sheaves-is-sheaf-of-infty-categories-on-affine-analytic-spaces} (using \cref{rslt:sheaves-on-basis-equiv-sheaves-on-whole-site}).\footnote{A priori there might be set-theoretic issues with defining the sheaf $\D(-)$ since $\AnSpace_G$ is usually a large $\infty$-category. They can easily be circumvented by using that every $G$-analytic space $X$ admits a \emph{small} covering by affine spaces and noting that $\D(X)$ only depends on the values of $\D(-)$ on any analytic cover of $X$.} The objects of $\D(X)$ are called the \emph{quasicoherent sheaves} on the $G$-analytic space $X$.

	\item For every map $f\colon Y \to X$ of $G$-analytic spaces we let
	\begin{align*}
		f^*\colon \D(X) \to \D(Y)
	\end{align*}
	denote the ``restriction'' functor of the sheaf $\D$ and call it the \emph{pullback} along $f$.

	\item For every qcqs map $f\colon Y \to X$ of $G$-analytic spaces we let
	\begin{align*}
		f_*\colon \D(Y) \to \D(X)
	\end{align*}
	denote a right-adjoint of $f^*$ (see \cref{rslt:pushforward-of-qcoh-on-analytic-spaces-exists} below) and call it the \emph{pushforward} along $f$.
\end{defenum}
\end{definition}

\begin{remark}
It is reasonable to only define the pushforward $f_*\colon \D(Y) \to \D(X)$ along \emph{qcqs} maps $f\colon Y \to X$ because quasicoherent sheaves only behave well under qcqs pushforward.
\end{remark}

\begin{lemma} \label{rslt:pushforward-of-qcoh-on-analytic-spaces-exists}
Let $G$ be a geometry blueprint over $(V,\mm)$ and $f\colon Y \to X$ a qcqs map of $G$-analytic spaces. Then $f^*\colon \D(X) \to \D(Y)$ admits a right adjoint $f_*\colon \D(Y) \to \D(X)$. Moreover, $f_*$ commutes with restriction, i.e. for every $G$-analytic subspace $U \subset X$ and every $\mathcal M \in \D(Y)$ we have $\restrict{(f_* \mathcal M)}U = (\restrict{f_*}{f^{-1}U}) \restrict{\mathcal M}{f^{-1}U}$.
\end{lemma}
\begin{proof}
We will use implicitly that if the right adjoint $f_*\colon \D(Y) \to \D(X)$ exists then it can be made functorial in $f$ (see \cite[Corollary 5.2.2.5]{lurie-higher-topos-theory}). First consider the case that $X$ is affine and $Y$ is a $G$-analytic subspace of an affine space. Cover $Y$ by a (finite) collection of affine $G$-analytic subspaces $W_i \subset Y$, $i \in I$. For each non-empty subset $J \subset I$ let $W_J := \bigisect_{j\in J} W_j$, which is still affine because $Y$ is a subspace of an affine space. Letting $f_J\colon W_J \to X$ denote the projection, each pullback functor $f_J^*$ admits a right adjoint $f_{J*}$ (namely the forgetful functor). Then $f_* := \varprojlim_J f_{J*}$ is a right adjoint of $f^*$.

Now assume that $X$ is still affine but $Y$ is arbitrary. We can again pick $W_i \subset Y$, $W_J$ and $f_J\colon W_J \to X$ as above. Then each $W_J$ is a $G$-analytic subspace of some $W_i$, hence by the previous case each $f_J^*$ admits a right adjoint $f_{J*}$ and we can again choose $f_* = \varprojlim_J f_{J*}$ as the right adjoint of $f^*$. Note that the claim about restrictions is easily checked in the above two cases, because everything commutes with finite limits, so we can reduce to the affine case, which follows from \cref{rslt:steady-map-of-analytic-rings-equiv-base-change}.

Now assume that $X$ is a $G$-analytic subspace of an affine space (and $Y$ is arbitrary). Pick a cover $X = \bigunion_{i\in I} U_i$ by $G$-analytic subspaces $U_i \subset X$. For each nonempty finite subset $J \subset I$ denote $U_J := \bigisect_{j\in J} U_j$ and let $W_J := f^{-1} U_J \subset Y$ with projection $f_J\colon W_J \to U_J$. Since $X$ is a subspace of an affine space, all $U_J$ are affine, hence all $f_{J*}$ exist by the above. We have $\D(X) = \varprojlim_J \D(U_J)$ and $\D(Y) = \varprojlim_J \D(W_J)$. Using the fact that all $f_{J*}$ commute with restrictions, we can thus define $f_* (\mathcal M_J)_J := (f_{J*} \mathcal M_J)_J$ (this construction preserves coCartesian edges). Then $f_*$ is a right adjoint of $f^*$ and the claim about restrictions is immediate from the construction.

Finally, if $X$ and $Y$ are arbitrary, we can again choose $U_i$, $U_J$, $W_J$ as before. Now each $U_J$ is a $G$-analytic subspace of an affine subspace (namely some $U_i$), hence by the previous case we can construct $f_*$ as before. It is again clear that the claim about restrictions is true.
\end{proof}

In the $\infty$-category $\AnRing_{(V,\mm)}$ of analytic rings over $(V,\mm)$ we singled out a special class of morphisms, namely the steady ones. They are defined by the property that they enable very general base-change. We can generalize the notion of steadiness to analytic spaces, as follows.

\begin{definition} \label{def:steady-morphism-of-analytic-spaces}
Let $G$ be a geometry blueprint over $(V,\mm)$. A morphism $f\colon Y \to X$ of $G$-analytic spaces is called \emph{steady} if for all affine $G$-analytic subspaces $U \subset X$ and all affine $G$-analytic subspaces $V \subset f^{-1}(U)$ the morphism $V \to U$ is induced by a steady morphism of analytic rings.
\end{definition}

\begin{lemma} \label{rslt:steadiness-can-be-checked-on-cover}
Let $f^\sharp\colon \mathcal A \to \mathcal B$ be a morphism of analytic rings over $(V,\mm)$ and assume that there exists a finite family of steady morphisms $g_i^\sharp\colon \mathcal B \to \mathcal C_i$ such that each $g_i^\sharp \comp f^\sharp$ is steady and the base-change morphisms $\D(\mathcal B) \to \D(\mathcal C_i)$ are jointly conservative. Then $f^\sharp$ is steady.
\end{lemma}
\begin{proof}
Follows from the base-change criterion for steadiness (see \cref{rslt:steady-map-of-analytic-rings-equiv-base-change}).
\end{proof}

\begin{lemma} \label{rslt:equivalent-characerizations-of-steady-morphisms}
Let $G$ be a geometry blueprint over $(V,\mm)$ and let $f\colon Y \to X$ be a morphism of $G$-analytic spaces. Then the following are equivalent:
\begin{lemenum}
	\item $f$ is steady.
	\item There exists a cover $\bigdunion_{i\in I} U_i \surjto X$ by $G$-analytic subspaces $U_i \subset X$ such that each $f^{-1}(U_i) \to U_i$ is steady.
	\item There exists a cover $\bigdunion_{i\in I} U_i \surjto X$ by affine $G$-analytic subspaces $U_i \subset X$ and for each $i$ a cover $\bigdunion_{j\in J_i} V_j \surjto f^{-1}(U_i)$ by affine $G$-analytic subspaces $V_j \subset f^{-1}(U_i)$ such that all the morphisms $V_j \to U_i$ are induced by steady morphisms of analytic rings.
\end{lemenum}
In particular, a morphism of affine $G$-analytic spaces is steady if and only if the corresponding morphism of analytic rings is steady.
\end{lemma}
\begin{proof}
It is clear that (i) implies (iii). To show that (iii) implies (ii) we can assume that $X = U_i$ for some $i$, in particular that $X$ is affine. Given any $G$-analytic subspace $V \subset Y$, it follows from the preservation of steadiness under composition and base-change that each morphism $V \isect V_j \to X$ is induced by a steady morphism of analytic rings. Letting $X = \AnSpec \mathcal A$, $V = \AnSpec \mathcal B$ and $V_j = \AnSpec \mathcal C_j$, it follows from \cref{rslt:steadiness-can-be-checked-on-cover} that $\mathcal A \to \mathcal B$ is steady, as desired. The proof that (ii) implies (i) is similar.
\end{proof}

\begin{proposition} \label{rslt:stability-properties-of-steadiness}
Let $G$ be a geometry blueprint over $(V,\mm)$.
\begin{propenum}
	\item Steady morphisms of $G$-analytic spaces are stable under composition and arbitrary base-change.

	\item Let $f\colon Y \to X$ and $g\colon Z \to Y$ be maps of $G$-analytic spaces. If $f$ and $f \comp g$ are steady, then so is $g$.
\end{propenum}
\end{proposition}
\begin{proof}
Follows easy from the corresponding properties for steady morphisms of analytic rings (see \cref{rslt:stability-of-steady-maps}) using \cref{rslt:equivalent-characerizations-of-steady-morphisms}.
\end{proof}

In the following we often make implicit use of \cref{rslt:equivalent-characerizations-of-steady-morphisms} and \cref{rslt:stability-properties-of-steadiness}. The next result shows that the steady maps defined above satisfy the expected base-change:

\begin{proposition} \label{rslt:steady-base-change-holds}
Let $G$ be a geometry blueprint over $(V,\mm)$ and let
\begin{center}\begin{tikzcd}
	Y' \arrow[r,"g'"] \arrow[d,"f'"] & Y \arrow[d,"f"]\\
	X' \arrow[r,"g"] & X
\end{tikzcd}\end{center}
be a cartesian square of $G$-analytic spaces such that $g$ is steady and $f$ is qcqs. Then the base change morphism
\begin{align*}
	g^* f_* \isoto f'_* g'^*
\end{align*}
is an isomorphism of functors from $\D(Y)$ to $\D(X')$.
\end{proposition}
\begin{proof}
Since pushforward commutes with restriction (see \cref{rslt:pushforward-of-qcoh-on-analytic-spaces-exists}) we can reduce to the case that $X$ and $X'$ are affine, so in particular $Y$ and $Y'$ are qcqs. Let $Y = \bigunion_{i\in I} W_i$ be a finite covering by affine $G$-analytic subspaces. For each non-empty subset $J \subset I$ let $W_J := \bigisect_{j \in J} W_j$. Then $\D(Y) = \varprojlim_J \D(W_J)$, so for every $\mathcal M = (\mathcal M_J)_J \in \D(Y)$ we have $f_* \mathcal M = \varprojlim_J f_{J*} \mathcal M_J$, where $f_{J*}\colon W_J \to X$ is the map induced by $f$. We get a similar representation for $\D(Y')$ and $f'_*$ in terms of the affine steady covering $\bigdunion_i f'^{-1}(W_i) \surjto Y'$. Since $g^*$ is exact and hence preserves finite limits, the claim reduces to the cases $Y = W_J$.

We are now reduced to the case that $Y$ is a $G$-analytic subspace of an affine space $\tilde Y$. In particular, for any two affine $G$-analytic subspaces $U, V \subset Y$ we get that $U \isect V = U \cprod_Y V = U \cprod_{\tilde Y} V$ is still affine. Thus, if we again choose an affine steady covering $Y = \bigunion_i W_i$ then all $W_J$ are affine. This reduces the claim to the case that $Y$ is affine. Then $Y'$ is also affine and we are done by \cref{rslt:steady-map-of-analytic-rings-equiv-base-change}.
\end{proof}

We also get a generalization of the notion of nuclear modules to the setting of analytic spaces. The following definition is sensible by \cref{rslt:nuclear-modules-stable-under-pullback}.

\begin{definition}
Let $G$ be a geometry blueprint over $(V,\mm)$, let $X$ be a $G$-analytic space and let $\mathcal M \in \D(X)$ be a quasicoherent sheaf on $X$. We say that $\mathcal M$ is \emph{locally nuclear} if there is an analytic cover $X = \bigunion_i U_i$ of $X$ by affine $G$-analytic subspaces $U_i = \AnSpec_G \mathcal B_i$ such that each $\restrict{\mathcal M}{U_i} \in \D(\mathcal B_i)$ is nuclear. We denote
\begin{align*}
	\D(X)^\lnuc \subset \D(X)
\end{align*}
the full subcategory of locally nuclear sheaves on $X$.
\end{definition}

\begin{remark} \label{rmk:is-D-nuc-a-sheaf}
It would be desirable for the functor $\D(-)^\lnuc$ to be a sheaf of $\infty$-categories on the analytic site of ($G$-)analytic spaces (in which case we could name its objects ``nuclear sheaves'' instead of ``\emph{locally} nuclear sheaves''). This is currently unknown however; it might depend on choosing the ``right'' definition of steady localizations.
\end{remark}

Given a geometry blueprint $(V,\mm)$, we would in general not expect the $\infty$-category $\AnSpace_G$ of $G$-analytic spaces to embedd fully faithfully into the $\infty$-category $\AnSpace_{(V,\mm)}$ of all analytic spaces over $(V,\mm)$. However, there is still a comparison functor, as follows.

\begin{lemma} \label{rslt:properties-of-analytification}
Let $G$ be a geometry blueprint over $(V,\mm)$. There is a natural functor
\begin{align*}
	(-)^\an\colon \AnSpec_G \to \AnSpec_{(V,\mm)}
\end{align*}
which satisfies the following properties:
\begin{lemenum}
	\item For every $G$-analytic ring $\mathcal A$ we have $(\AnSpec_G \mathcal A)^\an = \AnSpec \mathcal A$.

	\item The functor $(-)^\an$ preserves disjoint unions and fiber products.

	\item Given a $G$-analytic space $X$, if $X$ is quasicompact or quasiseparated then so is $X^\an$.

	\item Given a map $f\colon Y \to X$ of $G$-analytic spaces, if $f$ is quasicompact, quasiseparated, steady or a cover, then so is $f^\an\colon Y^\an \to X^\an$. Moreover, if $U \subset X$ is a $G$-analytic subspace then $U^\an \subset X^\an$ is a steady subspace of $X^\an$.

	\item For every $G$-analytic space $X$ there is a natural equivalence $\D(X) = \D(X^\an)$ of symmetric monoidal $\infty$-categories, compatible with the functors $f^*$ and $f_*$. Under this equivalence we have $\D(X)^\lnuc \subset \D(X^\an)^\lnuc$.
\end{lemenum}
\end{lemma}
\begin{proof}
We ignore all set-theoretic issues arising from the fact that we work with big sites; they can be resolved by temporarily working with sheaves of \emph{big} anima. The embedding $R_G \injto \AnRing_{(V,\mm)}$ defines a morphism of sites $\nu\colon\AnRing_{(V,\mm)}^\opp \to R_G^\opp$ and hence an associated exact pullback functor $\nu^*\colon \Shv(R_G^\opp) \to \Shv(\AnRing_{(V,\mm)}^\opp)$ (see \cite[Proposition 6.2.3.20]{lurie-higher-topos-theory}). We define the desired functor $(-)^\an$ as the restriction of $\nu^*$ to the full subcategory $\AnSpace_G \subset \Shv(R_G^\opp)$. To make this work, we observe the following properties of $\nu^*$:
\begin{itemize}
	\item For every $G$-analytic ring $\mathcal A$ we have $\nu^* \AnSpec_G \mathcal A = \AnSpec \mathcal A$. Namely, $\nu^*$ is computed as the composition of a left Kan extension along $i\colon R_G^\opp \injto \AnRing_{(V,\mm)}^\opp$ followed by sheafification. The left Kan extension is left adjoint to the functor $\mathcal P(\AnRing_{(V,\mm)}) \to \mathcal P(R_G)$ given by composition with $i$, hence by \cite[Proposition 5.2.6.3]{lurie-higher-topos-theory} it maps $\AnSpec_G \mathcal A$ to $\AnSpec \mathcal A$. But the latter is already a sheaf, so no sheafification is needed.

	\item Let $\mathcal A$ be a $G$-analytic ring and $U \subset X := \AnSpec_G \mathcal A$ a $G$-analytic subspace. Then $\nu^* U \subset \nu^* X = \AnSpec \mathcal A$ is a steady subspace. Namely, by definition of $G$-analytic subspaces (and the fact that $\nu^*$ preserves covers) we reduce to the case that $U$ is an affine $G$-analytic subspace of $X$, i.e. given by a $G$-localization $f\colon \mathcal A \to \mathcal B$. Then $f$ is in particular a steady localization and hence defines a steady subspace of $\nu^* X$.

	\item Let $X$ be a $G$-analytic space. Then $\nu^* X$ is an analytic space. Namely, pick any covering $X = \bigunion_i U_i$ by affine $G$-analytic subspaces. Then $\nu^* X = \bigunion_i \nu^* U_i$. By the previous observations each $\nu^* U_i$ is an affine analytic space and each $\nu^* U_i \isect \nu^* U_j \subset \nu^* U_i$ is a steady subspace. Thus \cref{rslt:obtain-analytic-space-from-gluing} implies that $\nu^*X$ is indeed an analytic space.
\end{itemize}
We have now constructed the functor $(-)^\an\colon \AnSpace_G \to \AnSpace_{(V,\mm)}$ as the restriction of $\nu^*$. Moreover, the functor $(-)^\an$ satisfies (i) and (ii) because so does $\nu^*$. It remains to prove claims (iii), (iv) and (v).

For (iii), let $X$ be a given $G$-analytic space. If $X$ is quasicompact then we can cover it by a finite number of affine $G$-analytic subspaces. Then $X^\an$ is covered by a finite number of affine steady subspaces and thus quasicompact (because affine analytic spaces are quasicompact). Similarly, $X$ is quasiseparated if and only if it admits a cover $X = \bigunion_i U_i$ by affine $G$-analytic subspaces $U_i \subset X$ such that each intersection $U_i \isect U_j$ can be covered by finitely many affine steady subspaces. If this is the case then $X^\an$ clearly has the same property (with ``$G$-analytic subspace'' replaced by ``steady subspace'').

To prove (iv), let $f\colon Y \to X$ be given. Similar to the proof of (iii) one can show that if $f$ is quasicompact or quasiseparated, then so is $f^\an$ (cf. \cite[Lemmas 01K4, 01KO]{stacks-project}). Moreover, if $f$ is steady then it comes from steady maps of analytic rings on covers of $X$ and $Y$ by $G$-analytic subspaces; the same is then true for $f^\an$ so that $f^\an$ is also steady (see \cref{rslt:equivalent-characerizations-of-steady-morphisms}). Finally, if $U \subset X$ is a $G$-analytic subspace then $U^\an \subset X^\an$ is a steady subspace, because both can be checked after passing to a cover of $X$ (resp. $X^\an$) by affine $G$-analytic (resp. steady) subspaces, where it is true by the above observation for $\nu^*$.

It remains to prove (v). Let $\mathcal C \to \Delta^1$ be the Cartesian fibration classifying the functor $(-)^\an$. Then $\mathcal C$ can naturally be equipped with an analytic site, of which the full subcategory $\mathcal C_a \subset \mathcal C$ spanned by the affine ($G$-)analytic spaces is a basis. Clearly $\D(-)$ defines a sheaf $\mathcal C_a \to \infcatinf^\tensor$ and hence extends to a sheaf $\D(-)\colon \mathcal C^\opp \to \infcatinf^\tensor$ by \cref{rslt:sheaves-on-basis-equiv-sheaves-on-whole-site}. For every $G$-analytic space $X$ this induces a symmetric monoidal functor $\D(X) \to \D(X^\an)$ compatible with pullbacks. That this functor is an equivalence can now be checked on any $G$-analytic cover of $X$, so that w.l.o.g. $X$ is affine. But then the claim is clear. Moreover, since the isomorphism $\D(X) = \D(X^\an)$ is compatible with pullbacks, it is automatically compatible with pushforwards as well. The claim about locally nuclear sheaves is obvious.
\end{proof}

\begin{definition} \label{def:analytification}
Let $G$ be a geometry blueprint over $(V,\mm)$ and $X$ a $G$-analytic space. The analytic space $X^\an$ over $(V,\mm)$ constructed in \cref{rslt:properties-of-analytification} is called the \emph{analytification} of $X$.
\end{definition}

Like in the previous subsections we will finish this subsections by discussing functoriality of the above definitions with respect to a change of almost setup. We get the following results.

\begin{lemma} \label{rslt:functoriality-of-analytic-spaces}
There is a natural functor
\begin{align*}
	\AlmSetup \to \infcatinf, \qquad (V,\mm) \mapsto \AnSpace_{(V,\mm)}.
\end{align*}
For every morphism $\varphi\colon (V,\mm) \to (V',\mm')$ of almost setups, the induced functor
\begin{align*}
	\varphi^*\colon \AnSpace_{(V,\mm)} \to \AnSpace_{(V',\mm')}
\end{align*}
satisfies the following properties:
\begin{lemenum}
	\item For every analytic ring $\mathcal A$ over $(V,\mm)$ we have $\varphi^* \AnSpec \mathcal A = \AnSpec \varphi^* \mathcal A$.

	\item $\varphi^*$ preserves disjoint unions and fiber products.

	\item Given an analytic space $X$ over $(V,\mm)$, if $X$ is quasicompact or quasiseparated then so is $\varphi^* X$.

	\item Given a map $f\colon Y \to X$ of $G$-analytic spaces, if $f$ is quasicompact, quasiseparated, steady or a cover, then so is $\varphi^* f\colon \varphi^* Y \to \varphi^* X$. Moreover, if $U \subset X$ is a steady subspace of some analytic space $X$ over $(V,\mm)$, then $\varphi^* U \subset \varphi^* X$ is a steady subspace of $\varphi^* X$.
\end{lemenum}
\end{lemma}
\begin{proof}
Starting with the functor $(V,\mm) \mapsto \AnRing_{(V,\mm)}$ from \cref{rslt:functoriality-of-analytic-rings-over-AlmSetup} we can construct a functor
\begin{align*}
	\AlmSetup^\opp \to \infcatinf, \qquad (V,\mm) \mapsto \Fun((\AnRing_{(V,\mm)})^\opp, \Ani)
\end{align*}
Let $p''\colon \mathcal C'' \to \AlmSetup$ be the Cartesian fibration classifying this functor. Let $\mathcal C' \subset \mathcal C''$ be the full subcategory where in the fiber over every $(V,\mm)$ we only allow sheaves for the analytic site. Since by \cref{rslt:almost-localization-preserves-steadiness} the functors $\varphi^*\colon \AnRing_{(V,\mm)} \to \AnRing_{(V',\mm')}$ preserve steady localizations, the map $p'\colon \mathcal C' \to \AlmSetup$ is still a Cartesian fibration. It is also a coCartesian fibration (where we temporarily pass to sheaves of \emph{large} anima): By \cite[Corollary 5.2.2.5]{lurie-higher-topos-theory} it is enough to show that for every morphism $\varphi\colon (V,\mm) \to (V',\mm')$ of almost setups, the pushforward
\begin{align*}
	\varphi_*\colon \Shv((\AnRing_{(V',\mm')})^\opp) \to \Shv((\AnRing_{(V,\mm)})^\opp)
\end{align*}
admits a left adjoint $\varphi^*$. This follows from \cite[Proposition 6.2.3.20]{lurie-higher-topos-theory} provided we can show that the functor $\varphi^*\colon \AnRing_{(V,\mm)} \to \AnRing_{(V',\mm')}$ preserves analytic covers. This boils down to the following statement: Given a map $f\colon \mathcal A \to \mathcal B$ of analytic rings over $(V,\mm)$ such that $f^*\colon \D(\mathcal A) \to \D(\mathcal B)$ is conservative, then also $(\varphi^* f)^*\colon \D(\varphi^* \mathcal A) \to \D(\varphi^* \mathcal B)$ is conservative. This can be checked separately in the case that $\varphi$ is strict and that $\varphi$ is localizing. The former case is easy because then $\varphi^*$ is just a base-change. In the latter case the functor $\varphi^*\colon \D(-) \to \D(\varphi^* -)$ admits a fully faithful left adjoint $\varphi_!$ by \cref{rslt:left-adjoint-of-almost-localization-over-analytic-ring}. Then if $M \in \D(\varphi^* \mathcal A)$ satisfies $(\varphi^* f)^* M = 0$, we deduce from \cref{rslt:almost-localization-is-steady} (by passing to left adjoints) that $0 = \varphi_! (\varphi^* f)^* M = f^* \varphi_! M$. Thus $\varphi_! M = 0$ because $f^*$ is conservative, and consequently $M = 0$ as desired. This finishes the proof that $p'$ is a coCartesian fibration.

Now let $\mathcal C \subset \mathcal C'$ be the full subcategory where in the fiber over $(V,\mm)$ we only allow the analytic spaces over $(V,\mm)$. We claim that the map $p\colon \mathcal C \to \AlmSetup$ is still a coCartesian fibration (and hence defines the desired functor). To see this it is enough to show that the sheaf pullback functor $\varphi^*$ constructed above preserves analytic spaces. This follows in the same way as in the proof of \cref{rslt:properties-of-analytification}. We similarly deduce properties (i) to (iv).
\end{proof}

\begin{proposition} \label{rslt:functoriality-of-sheaves-on-analytic-spaces-over-AlmSetup}
Let $\AnSpace_{(-)} \to \AlmSetup$ be the coCartesian fibration classifying the functor from \cref{rslt:functoriality-of-analytic-spaces}. Then there is a natural functor
\begin{align*}
	\AnSpace_{(-)} \to \infcatinf^\tensor, \qquad X \mapsto \D(X).
\end{align*}
Moreover, for every morphism $\varphi\colon (V,\mm) \to (V',\mm')$ of almost setups, the induced functor $\D(X) \to \D(\varphi^* X)$ preserves locally nuclear sheaves.
\end{proposition}
\begin{proof}
The desired functor can be defined as a right Kan extension of the functor $\AnRing_{(-)} \to \infcatinf^\tensor$, $\mathcal A \mapsto \D(\mathcal A)$ defined in \cref{rslt:functoriality-of-modules-over-analytic-rings-over-AlmSetup}. It preserves locally nuclear sheaves by the same result.
\end{proof}

\begin{definition} \label{def:functoriality-of-analytic-spaces-and-sheaves}
Let $\varphi\colon (V,\mm) \to (V',\mm')$ be a morphism of almost setups. We denote
\begin{align*}
	\varphi^*\colon \AnSpace_{(V,\mm)} \to \AnSpace_{(V',\mm')}, \qquad X \mapsto \varphi^* X
\end{align*}
the functor from \cref{rslt:functoriality-of-analytic-spaces}. If $\varphi$ is of the form $(V, V) \to (V,\mm)$ then we also denote $\varphi^*$ as $(V,\mm)^*$ and abbreviate it as $X \mapsto X^a$.

For every analytic space $X$ over $(V,\mm)$ we denote
\begin{align*}
	\varphi^*\colon \D(X) \to \D(\varphi^*X)
\end{align*}
the functor induced from \cref{rslt:functoriality-of-sheaves-on-analytic-spaces-over-AlmSetup}. By construction it commutes with pullbacks along maps of analytic spaces. If $\varphi$ is of the form $(V, V) \to (V,\mm)$ then we also denote $\varphi^*$ as $(V,\mm)^*$ and abbreviate it as $\mathcal M \mapsto \mathcal M^a$.
\end{definition}

\begin{remark} \label{rmk:almost-localization-of-G-analytic-spaces}
If $\varphi\colon (V,\mm) \to (V',\mm')$ is a morphism of almost setups and $G = (R_G, L_G)$ is a geometry blueprint over $(V,\mm)$ then one can define the geometry blueprint $\varphi^* G := (\varphi^* R_G, \varphi^* L_G)$ over $(V',\mm')$. Like in the proof of \cref{rslt:functoriality-of-analytic-spaces} one also obtains a functor
\begin{align*}
	\varphi^*\colon \AnSpace_G \to \AnSpace_{\varphi^* G}.
\end{align*}
Of course this observation can be made functorial in $\varphi$ and $G$ (with an appropriate definition of morphisms of geometry blueprints).
\end{remark}

\subsection{Steady Endofunctors} \label{sec:andesc.endofun}

As before, we fix an almost setup $(V,\mm)$. In the previous subsections we have constructed the $\infty$-category of analytic spaces over $(V,\mm)$. We now want to build a powerful theory of descent for quasicoherent sheaves on analytic spaces, following Mathew's theory of descendable algebras in \cite[\S3]{akhil-galois-group-of-stable-homotopy}. In the affine setting this means that we want to define when a map $f\colon \mathcal A \to \mathcal B$ of analytic rings over $(V,\mm)$ exhibits $\mathcal B$ as a ``descendable'' algebra over $\mathcal A$. However, Mathew's notion of descendability requires the $\mathcal A$-algebra $\mathcal B$ to sit inside a surrounding stable $\infty$-category of ``$\mathcal A$-modules''. This is easily satisfied in the classical setting of affine schemes (with just the usual derived $\infty$-category of discrete modules), but in the setting of analytic rings it is much less clear what this $\infty$-category of $\mathcal A$-modules should be: The analytic ring $\mathcal B$ carries a notion of completion (parametrized by the objects $\mathcal B[S]$ for extremally disconnected sets $S$) which is not captured by the underlying object $\underline{\mathcal B} \in \D(\mathcal A)$. The present section solves this problem by introducing a new notion of the $\infty$-category of ``$\mathcal A$-modules'' $\mathcal E(\mathcal A)$ which does indeed allows us to view any (steady) analytic $\mathcal A$-algebra $\mathcal B$ as an object inside of $\mathcal E(\mathcal A)$.

The main idea to construct this new $\infty$-category of ``$\mathcal A$-modules'' is as follows: We can characterize any $\mathcal A$-algebra $\mathcal B$ by its action $- \tensor_{\mathcal A} \mathcal B$ on $\D(\mathcal A)$, i.e. as a functor $\D(\mathcal A) \to \D(\mathcal A)$. Thus as a first approximation to $\mathcal E(\mathcal A)$ we can take the $\infty$-category $\Fun(\D(\mathcal A), \D(\mathcal A))$ of endofunctors of $\D(\mathcal A)$. This $\infty$-category comes equipped with the composition monoidal structure, which can be seen as a replacement for the tensor product on $\D(\mathcal A)$. This already allows us to define a reasonable notion of descendable $\mathcal A$-algebras by viewing them as objects in $\Fun(\D(\mathcal A), \D(\mathcal A))$. However, this approach has one major downside: Given a map $f\colon \mathcal A \to \mathcal A'$ to some other analytic ring over $(V,\mm)$, there is no good ``pullback'' functor
\begin{align*}
	f^\natural\colon \Fun(\D(\mathcal A), \D(\mathcal A)) \to \Fun(\D(\mathcal A'), \D(\mathcal A')).
\end{align*}
Such a pullback functor should be monoidal and also satisfy $f^\natural(- \tensor_{\mathcal A} \mathcal B) = - \tensor_{\mathcal A'} \mathcal B'$, where $\mathcal B' := \mathcal B \tensor_{\mathcal A} \mathcal A'$. The presence of such a pullback functor would allow us to conclude that if $\mathcal B$ is descendable as an $\mathcal A$-algebra then $\mathcal B'$ is descendable as an $\mathcal A'$-algebra -- a property that is certainly required to hold for any useful definition of descendability! The closest approximation to $f^\natural$ one can get is the functor $F \mapsto f^* \comp F \comp f_*$, but this functor clearly does not satisfy the above requirements; in fact it does not even satisfy $f^\natural \id = \id$ unless $f$ is a localization.

To fix the pullback issue explained in the previous paragraph, we need to work with a more involved formalism. Namely, we observe that for every $\mathcal A$-algebra $\mathcal B$, the functor $F := - \tensor_{\mathcal A} \mathcal B\colon \D(\mathcal A) \to \D(\mathcal A)$ comes equipped with an additional ``$\mathcal A$-linear'' structure. Roughly this means that for all $\mathcal A$-modules $M, N \in \D(\mathcal A)$ there is a canonical map $M \tensor F(N) \to F(M \tensor N)$. More concretely, $F$ can naturally be seen as a \emph{$\D(\mathcal A)$-enriched} functor $\D(\mathcal A) \to \D(\mathcal A)$ (see \cref{sec:infcat.enriched} for an introduction to enriched $\infty$-categories). This suggests that we should replace $\Fun(\D(\mathcal A), \D(\mathcal A))$ by the $\infty$-category $\Fun_{\D(\mathcal A)}(\D(\mathcal A), \D(\mathcal A))$ of enriched endofunctors of $\D(\mathcal A)$. For set-theoretical reasons we will instead work with a slightly smaller full subcategory, which we will denote $\End(\D(\mathcal A))$ (see \cref{def:analytic-endofunctors} below). In the first part of this subsection we study the basic properties of $\End(\D(\mathcal A))$ and in particular define the desired pullback functor
\begin{align*}
	f^\natural\colon \End(\D(\mathcal A)) \to \End(\D(\mathcal A')).
\end{align*}
Unfortunately we are still not done: While this new pullback functor does indeed satisfy $f^\natural(- \tensor_{\mathcal A} \mathcal B) = - \tensor_{\mathcal A'} \mathcal B'$, it is still not monoidal with respect to the composition monoidal structure. Moreover, the thus defined functor $\AnRing_{(V,\mm)} \to \infcatinf$, $\mathcal A \to \End(\D(\mathcal A))$ is not a sheaf for the analytic topology, so if we define descendability via $\End(\D(\mathcal A))$ then it is not clear that this notion can be checked on analytic covers -- which is certainly desirable.

The above mentioned problems of $\End(\D(\mathcal A))$ can be resolved by passing to a full subcategory $\End(\D(\mathcal A))^\std \subset \End(\D(\mathcal A))$ of so-called \emph{steady} endofunctors (see \cref{def:steady-object-of-End} below). They are defined to satisfy a base-change property which is very similar to the definition of steady morphisms of analytic rings. In particular, for every steady map $\mathcal A \to \mathcal B$ the endofunctor $- \tensor_{\mathcal A} \mathcal B$ is steady and hence an object of $\End(\D(\mathcal A))^\std$. The pullback functors restrict to monoidal pullback functors
\begin{align*}
	f^\natural\colon \End(\D(\mathcal A))^\std \to \End(\D(\mathcal A'))^\std
\end{align*}
and the thus defined functor $\End(\D(-))^\std$ defines an analytic sheaf. In particular it glues and therefore provides a good notion of a stable monoidal $\infty$-category $\End(\D(X))^\std$ for every analytic space $X$ over $(V,\mm)$.

The $\infty$-category $\End(\D(X))^\std$ works very well and enables a good notion of descendability. It still has two downsides however, which would make some of its applications a bit awkward: Firstly, evaluating a steady endofunctor $F \in \End(\D(\mathcal A))^\std$ on some module $M \in \D(\mathcal A)$ does not commute with pullbacks along steady localizations and hence does not glue to define an evaluation map $\End(\D(X))^\std \cprod \D(X) \to \D(X)$. Secondly, the monoidal structure on $\End(\D(X))^\std$ is not functorial in the almost setup $(V,\mm)$. Both of these problems can be solved by restricting to the full subcategory $\mathcal E(X) \subset \End(\D(X))^\std$ which is generated by endofunctors of the form $- \tensor_{\mathcal A} \mathcal B$. We thus obtain the desired stable monoidal $\infty$-category $\mathcal E(X)$ which provides a suitable notion of ``generalized nuclear $\mathcal A$-modules''.

After the above introductory remarks, we now start with the rigorous construction of $\mathcal E(X)$. Let us begin with the definition of $\End(\D(\mathcal A))$. In the following, note that whenever we have a map $f\colon \mathcal A \to \mathcal B$ of analytic rings over $(V,\mm)$, we naturally get an enrichment of $\D(\mathcal B)$ over $\D(\mathcal A)$ via $f^*$, as in \cref{ex:enriched-category-via-monoidal-functor}. Under this enrichment, the $\D(\mathcal A)$-enriched $\Hom$'s in $\D(\mathcal B)$ are just the internal $\Hom$-objects in $\D(\mathcal B)$ viewed as $\mathcal A$-modules via the forgetful functor $f_*$.

\begin{definition}
\begin{defenum}
	\item Let $f\colon \mathcal A \to \mathcal B$ be a map of analytic rings over $(V,\mm)$. We denote
	\begin{align*}
		\Fun_f(\D(\mathcal B), \D(\mathcal A)) \subset \Fun_{\D(\mathcal A)}(\D(\mathcal B), \D(\mathcal A))
	\end{align*}
	the full subcategory of those $\D(\mathcal A)$-enriched functors $\D(\mathcal B) \to \D(\mathcal A)$ whose underlying functor $F\colon \D(\mathcal B) \to \D(\mathcal A)$ is exact and satisfies the following property: For all large enough strong limit cardinals $\kappa'$ there is a strong limit cardinal $\kappa$ such that the restriction of $F$ to $\D(\mathcal B)_{\kappa'}$ factors through $\D(\mathcal A)_\kappa$. Here $\D(\mathcal B)$ is enriched over $\D(\mathcal A)$ via $f$, as explained above.

	Let $\kappa$ and $\kappa'$ be strong limit cardinals such that the restriction of $f_*\colon \D(\mathcal B) \to \D(\mathcal A)$ to $\D(\mathcal B)_{\kappa'}$ factors over $\D(\mathcal A)_\kappa$, we similarly define
	\begin{align*}
		\Fun_f(\D(\mathcal B)_{\kappa'}, \D(\mathcal A)_\kappa) \subset \Fun_{\D(\mathcal A)}(\D(\mathcal B)_{\kappa'}, \D(\mathcal A)_{\kappa})
	\end{align*}
	to be the full subcategory spanned by the exact functors.

	\item \label{def:analytic-endofunctors} Let $\mathcal A$ be an analytic ring over $(V,\mm)$. We denote
	\begin{align*}
		\End(\D(\mathcal A)) := \Fun_{\id}(\D(\mathcal A), \D(\mathcal A)).
	\end{align*}
\end{defenum}
\end{definition}

At the end we will only be interested in the $\infty$-category $\End(\D(\mathcal A))$, but for setting up the theory it is helpful to have the more general definition of $\Fun_f(\D(\mathcal B), \D(\mathcal A))$ at hand.

\begin{remark} \label{rslt:analytic-enriched-functors-are-lim-colim-of-presentable-case}
Let $f\colon \mathcal A \to \mathcal B$ be a map of analytic rings over $(V,\mm)$ and let $\kappa'$ be a strong limit cardinal. Then by definition, for every $F \in \Fun_f(\D(\mathcal B), \D(\mathcal A))$ and every large enough strong limit cardinal $\kappa'$, the restriction of $F$ to $\D(\mathcal B)_{\kappa'}$ lies in $\Fun_f(\D(\mathcal B)_{\kappa'}, \D(\mathcal A)_\kappa)$ for some strong limit cardinal $\kappa$. It follows that
\begin{align*}
	\Fun_f(\D(\mathcal B), \D(\mathcal A)) = \varprojlim_{\kappa'} \varinjlim_{\kappa \in I_{\kappa'}} \Fun_f(\D(\mathcal B)_{\kappa'}, \D(\mathcal A)_\kappa),
\end{align*}
where for each strong limit cardinal $\kappa'$ we choose $I_{\kappa'}$ to be a cofinal class of strong limit cardinals $\kappa$ such that the forgetful functor $\D(\mathcal B)_{\kappa'} \to \D(\mathcal A)$ factors over $\D(\mathcal A)_\kappa$. Moreover, for all $\kappa$ and $\kappa'$ as above we have
\begin{align*}
	\Fun_f(\D(\mathcal B)_{\kappa'}, \D(\mathcal A)_\kappa) \subset \Fun_{\D(\mathcal A)_\kappa}(\D(\mathcal B)_{\kappa'}, \D(\mathcal A)_\kappa),
\end{align*}
i.e. we can work with $\D(\mathcal A)_\kappa$-enriched functors instead of $\D(\mathcal A)$-enriched functors. These observations are very useful because they often allow us to reduce our questions to the case of presentable $\infty$-categories.
\end{remark}

\begin{remark}
Let $f\colon \mathcal A \to \mathcal B$ be a map of analytic rings over $(V,\mm)$. Then $\Fun_f(\D(\mathcal B), \D(\mathcal A))$ is stable and admits all small limits and colimits. This follows easily from the fact that the same is true for the underlying functor category $\Fun(\D(\mathcal B), \D(\mathcal A))$ using \cref{rslt:enriched-forgetful-functor-properties}.
\end{remark}

\begin{example}
Let $f\colon \mathcal A \to \mathcal B$ be a map of analytic rings over $(V, \mm)$. Then by \cref{rslt:weak-enriched-yoneda}, for every $M \in \D(\mathcal A)$ and $N \in \D(\mathcal B)$ there is the object
\begin{align*}
	\IHom_{\mathcal B}(N, -) \tensor_{\mathcal A} M \in \Fun_f(\D(\mathcal B), \D(\mathcal A)).
\end{align*}
\end{example}

The main technical reason for working with enriched functors instead of ordinary functors is the following variant of \cref{rslt:invariance-of-enriched-presheaves}:

\begin{lemma} \label{rslt:invariance-of-analytic-enriched-presheaves}
Let $\mathcal A \xto{f} \mathcal A' \xto{g} \mathcal B$ be maps of analytic rings over $(V,\mm)$ with composition $h = g \comp f$. Assume that $\mathcal A'$ is equipped with the induced analytic ring structure from $\mathcal A$. Then the natural forgetful functor
\begin{align*}
	\Fun_g(\D(\mathcal B), \D(\mathcal A')) \isoto \Fun_h(\D(\mathcal B), \D(\mathcal A))
\end{align*}
is an equivalence.
\end{lemma}
\begin{proof}
The functor is the one from \cref{rslt:forgetful-functor-of-enriched-presheaves} (using \cref{rslt:analytic-enriched-functors-are-lim-colim-of-presentable-case}). To show that it is an equivalence, we can by \cref{rslt:analytic-enriched-functors-are-lim-colim-of-presentable-case} reduce to showing that for all suitable strong limit cardinals $\kappa'$ and $\kappa$ the functor
\begin{align*}
	\Fun_{\D(\mathcal A')_\kappa}(\D(\mathcal B)_{\kappa'}, \D(\mathcal A)_\kappa) \isoto \Fun_{\D(\mathcal A)_\kappa}(\D(\mathcal B)_{\kappa'}, \D(\mathcal A)_\kappa)
\end{align*}
is an equivalence. This follows immediately from \cref{rslt:invariance-of-enriched-presheaves}: Clearly the forgetful functor $f_*\colon \D(\mathcal A')_\kappa \to \D(\mathcal A)_\kappa$ is conservative and preserves small colimits (note that since $\mathcal A'$ has the induced analytic ring structure from $\mathcal A$, the forgetful functor $f_*$ preserves $\kappa$-condensed objects); it also satisfies the projection formula because $\mathcal A'$ has the induced analytic ring structure from $\mathcal A$.
\end{proof}

With \cref{rslt:invariance-of-analytic-enriched-presheaves} at hand we can now prove the following functoriality results of the $\infty$-categories $\Fun_f(\D(\mathcal B), \D(\mathcal A))$:

\begin{lemma} \label{rslt:lower-upper-natural-functors-for-Fun-of-analytic-rings}
Let
\begin{center}\begin{tikzcd}
	\mathcal B' & \mathcal B \arrow[l,"g'"]\\
	\mathcal A' \arrow[u,"f'"] & \mathcal A \arrow[l,"g"] \arrow[u,"f"]
\end{tikzcd}\end{center}
be a commuting diagram of analytic rings over $(V,\mm)$ (not necessarily cartesian). Then there is a natural pair of adjoint functors
\begin{align*}
	(f',f)^\natural\colon \Fun_g(\D(\mathcal A'), \D(\mathcal A)) & \rightleftarrows \Fun_{g'}(\D(\mathcal B'), \D(\mathcal B))\noloc (f',f)_\natural
\end{align*}
satisfying the following properties:
\begin{lemenum}
	\item \label{rslt:properties-of-lower-natural} $(f',f)_\natural$ is conservative and preserves small limits and colimits. It acts as $F \mapsto f_* \comp F \comp f'^*$ on the underlying functors.

	\item \label{rslt:upper-natural-on-generators} For all $N \in \D(\mathcal A')$ and all $M \in \D(\mathcal A)$ we have
	\begin{align*}
		(f',f)^\natural \big(\IHom_{\mathcal A'}(N, -) \tensor_{\mathcal A} M \big) = \IHom_{\mathcal B'}(f'^*N, -) \tensor_{\mathcal A'} f^*M.
	\end{align*}

	\item \label{rslt:computation-of-upper-natural-if-f-localization} If $f$ is a localization then $(f',f)^\natural$ acts as $F \mapsto f^* \comp F \comp f'_*$ on the underlying functors.

	\item \label{rslt:lower-natural-is-fully-faithful-if-f'-localization} If $f'$ is a localization then $(f',f)_\natural$ is fully faithful.
\end{lemenum}
\end{lemma}
\begin{proof}
We first construct the functor $(f', f)_\natural$. By \cref{rslt:forgetful-functor-of-enriched-presheaves} (together with \cref{rslt:analytic-enriched-functors-are-lim-colim-of-presentable-case}) there is a natural functor $\theta_1\colon \Fun_{g'}(\D(\mathcal B'), \D(\mathcal B)) \to \Fun_{g'f}(\D(\mathcal B'), \D(\mathcal A))$ which acts as $F \mapsto f_* \comp F$ on the underlying functors. Moreover, we can view the functor $f'^*\colon \D(\mathcal A') \to \D(\mathcal B')$ as a morphism of $\D(\mathcal A)$-modules and in particular as a $\D(\mathcal A)$-enriched functor. Thus, composition with $f'^*$ defines a functor $\theta_2\colon \Fun_{f'g}(\D(\mathcal B'), \D(\mathcal A)) \to \Fun_g(\D(\mathcal A'), \D(\mathcal A))$. This functor acts as $F \mapsto F \comp f'^*$ on the underlying functor categories. Noting that $f'g = g'f$, we can define $(f',f)_\natural = \theta_2 \comp \theta_1$. Then $(f',f)_\natural$ acts as $F \mapsto f_* \comp F \comp f'^*$ on the underlying functor categories and is thus conservative and preserves all small limits and colimits. This proves (i).

We now prove the existence of the left adjoint $(f', f)^\natural$. First assume that $\mathcal B = \mathcal A$ and $f = \id$. We can view the map $f'^*\colon \D(\mathcal A') \to \D(\mathcal B')$ as a map of $\D(\mathcal A)$-modules and in particular as a $\D(\mathcal A)$ enriched functor. We can then construct $(f', f)^\natural \Fun_g(\D(\mathcal A'), \D(\mathcal A)) \to \Fun_{g'}(\D(\mathcal B'), \D(\mathcal A))$ as an operadic left Kan extension along $f'^*$ (see \cite[Corollary 3.1.3.4]{lurie-higher-algebra} and use similar reductions as in the proof of \cref{rslt:weak-enriched-yoneda}); it exists by the existence of the right adjoint $f_*$ of $f'^*$. One also sees that $(f', f)^\natural$ acts as $F \mapsto F \comp f'_*$ in this case.

For the proof of the existence of $(f', f)^\natural$ we can from now on assume that $\mathcal A' = \mathcal B'$ and $f' = \id$. Now factor $f$ as $\mathcal B \xto{f_1} \underline{\mathcal B}_{\mathcal A/} \xto{f_2} \mathcal A$, which allows us to reduce the construction of $(f', f)^\natural$ to the cases $f = f_1$ and $f = f_2$. Note that $(\id, f_2)_\natural$ is an equivalence by \cref{rslt:invariance-of-analytic-enriched-presheaves}, which only leaves the case $f = f_1$, i.e. we can assume that $f$ is a localization. Let $\mathcal C'^{\tensor\tensor} \subset \mathcal C^{\tensor\tensor} \to \opLM$ denote the fibrations of $\infty$-operads exhibiting the self-enrichments of $\D(\mathcal B)$ and $\D(\mathcal A)$, respectively. Let furthermore $\mathcal D^\tensor \to \opLM$ denote the fibration of $\infty$-operads exhibiting the $\D(\mathcal B)$-enrichment of $\D(\mathcal B')$. Then by definition $\Fun_{g'}(\D(\mathcal B'), \D(\mathcal B))$ consists of maps of $\infty$-operads $\mathcal D^\tensor \to \mathcal C'^{\tensor\tensor}$ over $\opLM$ which are the identity over $\mathfrak a$, while by the proof of \cref{rslt:forgetful-functor-of-enriched-presheaves} we can identify $\Fun_g(\D(\mathcal B'), \D(\mathcal A))$ with maps of $\infty$-operads $\mathcal D^\tensor \to \mathcal C^{\tensor\tensor}$ which restrict to the inclusion $\D(\mathcal B) \injto \D(\mathcal A)$ over $\mathfrak a$. With this description, $(\id,f)_\natural$ amounts to the composition with $f_*\colon \mathcal C'^{\tensor\tensor} \to \mathcal C^{\tensor\tensor}$. We can thus define a left adjoint $(\id, f)^\natural$ of $(\id, f)_\natural$ by mapping $F\colon \mathcal D^\tensor \to \mathcal C^{\tensor\tensor}$ to $f^* \comp F$, where $f^*\colon \mathcal C^{\tensor\tensor} \to \mathcal C'^{\tensor\tensor}$ is the functor induced from the symmetric monoidal base-change functor $\D(\mathcal A) \to \D(\mathcal B)$. This definition works because if a map $F\colon \mathcal D^\tensor \to \mathcal C^{\tensor\tensor}$ restricts to $f_*$ over $\mathfrak a$ then $(\id, f)^\natural F$ restricts to the identity over $\mathfrak a$ (since $f^* f_* = \id$) and thus lands in the correct category. This finishes the construction of the functor $(f', f)^\natural$.

Claim (iii) follows immediately from the above construction of $(f', f)^\natural$.

We now prove (iv), so assume that $f'$ is a localization. Note that we can factor the functor $(f',f)_\natural$ as
\begin{align*}
	\Fun_{g'}(\D(\mathcal B'), \D(\mathcal B)) \xto{(\id,f)_\natural} \Fun_{g'f}(\D(\mathcal B'), \D(\mathcal A)) \xto{(f',\id)_\natural} \Fun_g(\D(\mathcal A'), \D(\mathcal A)).
\end{align*}
Using (i), (iii) and $f'^* f'_* \isoto \id$, one checks immediately that the counit $(f',\id)^\natural \comp (f',\id)_\natural \to \id$ is an isomorphism, i.e. that $(f',\id)_\natural$ is fully faithful. It thus remains to show that $(\id,f)_\natural$ is fully faithful, so we can assume that $\mathcal B' = \mathcal A'$ and $f' = \id$ from now on. Factoring $f$ as $\mathcal B \xto{f_1} \underline{\mathcal B}_{\mathcal A/} \xto{f_2} \mathcal A$ allows us to reduce to the cases $f = f_1$ and $f = f_2$. The case $f = f_2$ is immediately handled by \cref{rslt:invariance-of-analytic-enriched-presheaves}. We are thus left with the case $f = f_1$, i.e. $f$ is a localization. Applying (i), (iii) and $f^* f_* \isoto \id$ we deduce that $(\id, f)^\natural \comp (\id, f)_\natural \isoto \id$ is an isomorphism, proving that $(\id, f)_\natural$ is fully faithful as desired.

It remains to prove (ii), so let $N \in \D(\mathcal A')$ and $M \in \D(\mathcal A)$ be given. Then by (i) the left one of the following diagrams commutes:
\begin{center}
	\begin{tikzcd}
		\Fun_g(\D(\mathcal A'), \D(\mathcal A)) \arrow[r,"\ev_N"] & \D(\mathcal A)\\
		\Fun_{g'}(\D(\mathcal B'), \D(\mathcal B)) \arrow[r,"\ev_{f'^*N}"] \arrow[u,"{(f',f)_\natural}"] & \D(\mathcal B) \arrow[u,"f_*"]
	\end{tikzcd}
	\qquad\qquad
	\begin{tikzcd}
		\Fun^L_{g_*}(\D(\mathcal A'), \D(\mathcal A)) \arrow[d,swap,"{(f',f)^\natural}"] & \D(\mathcal A) \arrow[l,swap,"\ell_N"] \arrow[d,swap,"f^*"]\\
		\Fun^L_{g'_*}(\D(\mathcal B'), \D(\mathcal B))  & \D(\mathcal B) \arrow[l,swap,"\ell_{f'^*N}"]
	\end{tikzcd}
\end{center}
Passing to left adjoints produces the right diagram (with $\ell_N$ as in \cref{rslt:weak-enriched-yoneda}), which is therefore still homotopy coherent. Plugging in $M$ immediately produces the claimed identity in (ii).
\end{proof}

\begin{definition} \label{def:lower-upper-natural-for-End-of-analytic-ring}
Let $f\colon \mathcal A \to \mathcal B$ be a map of analytic rings. We define the adjoint pair of functors
\begin{align*}
	f^\natural\colon \End(\D(\mathcal A)) \rightleftarrows \End(\D(\mathcal B)) \noloc f_\natural
\end{align*}
as $f^\natural := (f,f)^\natural$ and $f_\natural := (f,f)_\natural$ using \cref{rslt:lower-upper-natural-functors-for-Fun-of-analytic-rings}.
\end{definition}

Our next goal is to glue the $\infty$-categories $\End(\D(\mathcal A))$ along analytic covers in order to obtain an $\infty$-category ``$\End(\D(X))$'' for every analytic space $X$ over $(V,\mm)$.
The gluing procedure is formal if we can show that $\End(\D(-))$ is conservative and satisfies base-change along analytic covers. Base-change is easy to show (using steady base-change for $\D(-)$ and the explicit description of the functors $f^\natural$ and $f_\natural$), but it seems to be wrong that $\End(\D(-))$ is conservative along analytic covers (i.e. if $f\colon \mathcal A \to \mathcal B_1 \cprod \dots \mathcal B_n$ is a map of analytic rings over $(V,\mm)$ such that each $\mathcal A \to \mathcal B_i$ is a steady localization and $f^*$ is conservative, it does not follow that $f^\natural$ is conservative). We explored several different avenues for solving this issue and came up with the following solution:

\begin{definition} \label{def:steady-object-of-End}
Let $\mathcal A$ be an analytic ring over $(V,\mm)$. An endofunctor $F \in \End(\D(\mathcal A))$ is called \emph{steady} if for all maps $g\colon \mathcal A \to \mathcal B$ of analytic rings the object $(g,\id)^\natural F \in \Fun_g(\D(\mathcal B), \D(\mathcal A))$ lies in the essential image of $(\id,g)_\natural\colon \End(\D(\mathcal B)) \injto \Fun_g(\D(\mathcal B), \D(\mathcal A))$. We denote by
\begin{align*}
	\End(\D(\mathcal A))^\std \subset \End(\D(\mathcal A))
\end{align*}
the full subcategory of steady endofunctors.
\end{definition}

More intuitively, an endofunctor $F \in \End(\D(\mathcal A))$ is steady if for every map $g\colon \mathcal A \to \mathcal B$ the functor $F \comp g_*\colon \D(\mathcal B) \to \D(\mathcal A)$ factors as $\D(\mathcal B) \to \D(\mathcal B) \xto{g_*} \D(\mathcal A)$, i.e. for every $M \in \D(\mathcal B)$, $F(M) \in \D(\mathcal A)$ comes ``naturally'' equipped with a $\mathcal B$-module structure. The reason for choosing \cref{def:steady-object-of-End} is the following list of properties:

\begin{proposition} \label{rslt:steady-objects-of-End-enum}
\begin{propenum}
	\item \label{rslt:properties-of-End-std} For every analytic ring $\mathcal A$ over $(V,\mm)$, $\End(\D(\mathcal A))^\std$ is a stable $\infty$-category and it is stable under colimits in $\End(\D(\mathcal A))$. Moreover, for all $M, N \in \D(\mathcal A)$ with $M$ locally nuclear we have $\IHom(N, -) \tensor M \in \End(\D(\mathcal A))^\std$.

	\item \label{rslt:lower-and-upper-natural-restrict-to-End-std} For every map $f\colon \mathcal A \to \mathcal B$ of analytic rings over $(V,\mm)$, $f^\natural$ preserves steadiness and hence induces a colimit-preserving functor
	\begin{align*}
		f^\natural\colon \End(\D(\mathcal A))^\std \to \End(\D(\mathcal B))^\std.
	\end{align*}
	Moreover, $f$ is steady if and only if $f_\natural \id \in \End(\D(\mathcal A))$ is steady. In this case, $f_\natural$ preserves steady objects and hence induces a conservative, limit- and colimit-preserving functor
	\begin{align*}
		f_\natural\colon \End(\D(\mathcal B))^\std \to \End(\D(\mathcal A))^\std
	\end{align*}
	which is right adjoint to $f^\natural$.

	\item \label{rslt:base-change-for-End-std} Let
	\begin{center}\begin{tikzcd}
		\mathcal B' & \mathcal B \arrow[l,"g'"]\\
		\mathcal A' \arrow[u,"f'"] & \mathcal A \arrow[l,"g"] \arrow[u,"f",swap]
	\end{tikzcd}\end{center}
	be a pushout square of analytic rings over $(V,\mm)$ and assume that $f$ is steady. Then the base-change morphism
	\begin{align*}
		g^\natural f_\natural \isoto f'_\natural g'^\natural
	\end{align*}
	is an equivalence of functors from $\End(\D(\mathcal B))^\std$ to $\End(\D(\mathcal A'))^\std$.

	\item \label{rslt:conservativity-for-End-std} Let $(f_i\colon \mathcal A \to \mathcal B_i)_{i\in I}$ be a finite family of steady localizations which form an analytic cover. Then if $F \in \End(\D(\mathcal A))^\std$ satisfies $f_i^\natural F \isom 0$ for all $i$ then $F \isom 0$.
\end{propenum}
\end{proposition}
\begin{proof}
We start with (i), so let $\mathcal A$ be given. Given any map $g\colon \mathcal A \to \mathcal B$ of analytic rings over $(V,\mm)$, it follows from \cref{rslt:properties-of-lower-natural,rslt:lower-natural-is-fully-faithful-if-f'-localization} that $\End(\D(\mathcal B)) \subset \Fun_g(\D(\mathcal B), \D(\mathcal A))$ is stable under small colimits. As $(g, \id)^\natural$ preserves small colimits, it follows that $\End(\D(\mathcal A))^\std \subset \End(\D(\mathcal A))$ is stable under small colimits. Since $\End(\D(\mathcal A))$ is a stable $\infty$-category, the same follows for $\End(\D(\mathcal A))^\std$.

To finish the proof of (i) we need to show that $e_{M,N} := \IHom(N, -) \tensor M \in \End(\D(\mathcal A))$ is steady for all $M, N \in \D(\mathcal A)$ with $M$ nuclear. Fix such $M$ and $N$. Note that we need to show that for any given $g\colon \mathcal A \to \mathcal B$, the natural morphism $(g,\id)^\natural e_{M,N} \to (\id,g)_\natural g^\natural e_{M,N}$ of objects in $\Fun_g(\D(\mathcal B), \D(\mathcal A))$ is an isomorphism. This can be checked on underlying functors (by \cref{rslt:enriched-forgetful-functor-is-conservative}), so it is enough to check this after plugging in any $L \in \D(\mathcal B)$. By the explicit formulas in \cref{rslt:lower-upper-natural-functors-for-Fun-of-analytic-rings} we compute
 \begin{align*}
 	((g,\id)^\natural e_{M,N})(L) &= \IHom_{\mathcal A}(N, L) \tensor_{\mathcal A} M\\
 	&= \IHom_{\mathcal B}(N \tensor_{\mathcal A} \mathcal B, L) \tensor_{\mathcal A} M\\
 	((\id,g)_\natural g^\natural e_{M,N})(L) &= \IHom_{\mathcal B}(N \tensor_{\mathcal A} \mathcal B, L) \tensor_{\mathcal B} (M \tensor_{\mathcal A} \mathcal B)
\end{align*}
If $M$ is nuclear then by \cref{rslt:nucelar-modules-are-steady} the two functors $- \tensor_{\mathcal A} M$ and $- \tensor_{\mathcal B} (M \tensor_{\mathcal A} \mathcal B)$ from $\D(\mathcal B)$ to $\D(\mathcal A)$ are isomorphic, proving the desired identity. One checks immediately that this isomorphism can be checked on an analytic cover so that it also holds if $M$ is only \emph{locally} nuclear. This finishes the proof of (i).

To prove (ii), let $f\colon \mathcal A \to \mathcal B$ be given. To see that $f^\natural$ preserves steadiness, choose any $F \in \End(\D(\mathcal A))^\std$ and any map $g\colon \mathcal B \to \mathcal C$ of analytic rings over $(V,\mm)$. We have to show that $(g,\id)^\natural (f^\natural F)$ lies in the essential image of $(\id, g)_\natural$. But
\begin{align*}
	(g,\id_{\mathcal B})^\natural (f^\natural F) = (gf,f)^\natural F = (\id_{\mathcal C},f)^\natural (gf,\id_{\mathcal A})^\natural F.
\end{align*}
By steadiness of $F$, $(gf,\id_{\mathcal A})^\natural F$ lies in the essential image of $(\id_{\mathcal C},gf)_\natural = (\id_{\mathcal C},f)_\natural (\id_{\mathcal C},g)_\natural$. Thus the claim follows from the fact that $(\id_{\mathcal C},f)_\natural$ is fully faithful (see \cref{rslt:lower-natural-is-fully-faithful-if-f'-localization}). This finishes the proof that $f^\natural$ preserves steadiness. It follows immediately from (i) that $f^\natural$ preserves all small colimits.

Now assume that $f$ is steady. We want to show that $f_\natural$ preserves steadiness. Given any map $g\colon \mathcal A \to \mathcal A'$ of analytic rings, let $g'\colon \mathcal B \to \mathcal B'$ denote the base-change along $f$ and let $f'\colon \mathcal A' \to \mathcal B'$ denote the remaining map in the Cartesian diagram (as in (iii)). We obtain the following diagram of functors:
\begin{center}\begin{tikzcd}
	\End(\D(\mathcal B')) \arrow[r,"{(\id,g')_\natural}"] \arrow[d,swap,"f'_\natural"] & \Fun_{g'_*}(\D(\mathcal B'), \D(\mathcal B)) \arrow[d,"{(f',f)_\natural}"] & \End(\D(\mathcal B)) \arrow[l,swap,"{(g',\id)^\natural}"] \arrow[d,"f_\natural"]\\
	\End(\D(\mathcal A')) \arrow[r,swap,"{(\id,g)_\natural}"] & \Fun_{g_*}(\D(\mathcal A'), \D(\mathcal A)) & \End(\D(\mathcal A)) \arrow[l,"{(g,\id)^\natural}"]
\end{tikzcd}\end{center}
This diagram is naturally homotopy coherent: For the left square this is evident; for the right square this amounts to saying that the natural adjunction morphism $(g,\id)^\natural f_\natural \to (f',f)_\natural (g',\id)^\natural$ is an isomorphism. Using the explicit descriptions of all the functors (see \cref{rslt:computation-of-upper-natural-if-f-localization}) this reduces immediately to the claim that the natural adjunction morphism $f^* g_* \to g'_* f'^*$ is an isomorphism, which follows from the steadiness of the map $f$ (see \cref{rslt:steady-map-of-analytic-rings-equiv-base-change}). It follows that if $F \in \End(\D(\mathcal B))$ is steady then $(g,\id)^\natural (f_\natural F) = (f',f)_\natural (g',\id)^\natural F$ lies in the image of $(\id,g)_\natural$. Since $g$ was arbitrary, this proves that $f_\natural F$ is steady, as desired. It is clear that the restriction of $f_\natural$ to steady objects is still right adjoint to $f^\natural$, hence the restricted $f_\natural$ preserves limits (this is a priori not clear because steady objects are probably not preserved by limits, hence limits in $\End(\D(\mathcal A))^\std$ differ from limits in $\End(\D(\mathcal A))$).

To finish the proof of (ii), it remains to show that if $f_\natural \id \in \End(\D(\mathcal A))$ is steady, then $f$ is steady. Assume that $f_\natural \id$ is steady and let $g\colon \mathcal A \to \mathcal A'$ be any map of analytic rings. Then $(g,\id)^\natural f_\natural \id = - \tensor_{\mathcal A} \mathcal B$ as a functor $\D(\mathcal A') \to \D(\mathcal A)$. This functor clearly factors over $\D(\underline{\mathcal A'})$ and by the steadiness of $f_\natural \id$ it factors (necessarily uniquely) over $\D(\mathcal A')$. In other words, for every $M \in \D(\mathcal A')$ we have $M \tensor_{\mathcal A} \mathcal B \in \D(\mathcal A')$, so $f$ is steady by definition.

Part (iii) follows immediately from the homotopy coherent diagram in the proof of (ii), using that the left horizontal functors are fully faithful.

It remains to prove (iv). Let $(f_i)_i$ and $F \in \End(\D(\mathcal A))^\std$ be given as in the claim. For every non-empty subset $J \subset I$ let $\mathcal B_J := \bigtensor_{i\in J} \mathcal B_i$ and let $f_J\colon \mathcal A \to \mathcal B_J$ be the associated map. Then $f_J^\natural F \isom 0$ for all $J$. By the sheafiness of $\D(-)$ (see \cref{rslt:qcoh-sheaves-is-sheaf-of-infty-categories-on-affine-analytic-spaces}) we have $M = \varprojlim_J (M \tensor_{\mathcal A} \mathcal B_J)$ for every $M \in \D(\mathcal A)$. On the other hand, since $F$ is steady, $f_J^\natural F \isom 0$ implies $(f_J,\id)^\natural F \isom 0$, i.e. $F(N) = 0$ for any $N \in \D(\mathcal B_J)$. In particular,
\begin{align*}
	F(M) = F(\varprojlim_J (M \tensor_{\mathcal A} \mathcal B_J)) = \varprojlim_J F(M \tensor_{\mathcal A} \mathcal B_J) = \varprojlim_J 0 = 0.
\end{align*}
Thus $F \isom 0$, as desired.
\end{proof}

As a last preparation for gluing $\End(\D(\mathcal A))^\std$ we need to make it functorial in $\mathcal A$. We also want to maintain the composition monoidal structure on the glued $\infty$-categories $\End(\D(X))$, so we need to see that $\End(\D(\mathcal A))^\std$ is even functorial as a monoidal $\infty$-category. Here it is important that we work with \emph{steady} endofunctors, as they guarantee the pullback functors $f^\natural$ to be monoidal.

\begin{lemma}
\begin{lemenum}
	\item Let $\mathcal A$ be an analytic ring over $(V,\mm)$. Then $\End(\D(\mathcal A))^\std$ is stable under the composition monoidal structure on $\End(\D(\mathcal A))$ and hence becomes itself a monoidal $\infty$-category. The monoidal structure is exact in both arguments.

	\item Let $f\colon \mathcal A \to \mathcal B$ be a map of analytic rings over $(V,\mm)$. Then for all steady endofunctors $F_1, F_2 \in \End(\D(\mathcal A))^\std$ there is a natural isomorphism $f^\natural F_1 \comp f^\natural F_2 = f^\natural (F_1 \comp F_2)$.

	\item \label{rslt:functoriality-of-analytic-steady-endofunctors} The assignment $\mathcal A \mapsto \End(\D(\mathcal A))^\std$, $f \mapsto f^\natural$ defines a functor
	\begin{align*}
		\End(\D(-))\colon \AnRing_{(V,\mm)} \to \infcatinf^{\ocircle},
	\end{align*}
	where $\infcatinf^\ocircle$ denotes the $\infty$-category of monoidal $\infty$-categories.
\end{lemenum}
\end{lemma}
\begin{proof}
By the same argument as in the proof of \cref{rslt:functoriality-of-enriched-endofunctors} with \cref{rslt:invariance-of-analytic-enriched-presheaves} in place of \cref{rslt:invariance-of-enriched-presheaves} we obtain a functor
\begin{align*}
	\AnRing_{(V,\mm)} \to \infcatinf^\ocircle, \qquad \mathcal A \mapsto \End(\D(\underline{\mathcal A})).
\end{align*}
This functor maps $f\colon \mathcal A \to \mathcal B$ to $\underline f^\natural$, so that in particular $\underline f^\natural$ is monoidal. Let us denote $p''\colon \mathcal E''^\ocircle \to \AnRing_{(V,\mm)} \cprod \opAssoc$ the coCartesian family of monoidal $\infty$-categories classifying the above functor (cf. \cite[Example 4.8.3.3]{lurie-higher-algebra}). Then let $\mathcal E'^\ocircle \subset \mathcal E''^\ocircle$ be the full subcategory where we only allow objects in $\End(\D(\mathcal A)) \subset \End(\D(\underline{\mathcal A}))$ in the fiber over each analytic ring $\mathcal A$ (cf. \cref{rslt:lower-natural-is-fully-faithful-if-f'-localization}). We claim that the map $p'\colon \mathcal E'^\ocircle \to \AnRing_{(V,\mm)} \cprod \opAssoc$ is a locally coCartesian fibration. Since we know that $p''$ is a coCartesian fibration, this follows from the following two observations:
\begin{itemize}
	\item For every analytic ring $\mathcal B$ over $(V,\mm)$ the inclusion $\End(\D(\mathcal B)) \subset \End(\D(\underline{\mathcal B}))$ admits a left adjoint. Namely, this left adjoint is given by $j^\natural$, where $j\colon \underline{\mathcal B} \to \mathcal B$ is the natural map.

	\item For every analytic ring $\mathcal A$ over $(V,\mm)$ the subcategory $\End(\D(\mathcal A)) \subset \End(\D(\underline{\mathcal A}))$ is stable under the composition monoidal structure on $\End(\D(\underline{\mathcal A}))$. In fact, the composition monoidal structure on $\End(\D(\mathcal A))$ comes via restriction from the composition monoidal structure on $\End(\D(\underline{\mathcal A}))$: this follows easily from the direct description of $i_\natural$ in \cref{rslt:properties-of-lower-natural}, where $i\colon \underline{\mathcal A} \to \mathcal A$ is the natural map.
\end{itemize}
Namely, it follows from these observations that for every edge $e\colon (\mathcal A, \langle n \rangle) \to (\mathcal B, \langle 1 \rangle)$ in $\AnRing_{(V,\mm)} \cprod \opAssoc$, where $\langle n \rangle \to \langle 1 \rangle$ is the active morphism corresponding to the natural ordering on $\langle n \rangle^\circ$, the map $\mathcal E'^\ocircle_e \to \Delta^1$ classifies the functor
\begin{align*}
	\End(\D(\mathcal A))^n \to \End(\D(\mathcal B)), \qquad (F_1, \dots, F_n) \mapsto j^\natural \underline f^\natural (F_1 \comp \dots \comp F_n) = f^\natural (F_1 \comp \dots \comp F_n).
\end{align*}
Now let $\mathcal E^\ocircle \subset \mathcal E'^\ocircle$ be the full subcategory where we only allow objects in $\End(\D(\mathcal A))^\std \subset \End(\D(\mathcal A))$ in the fiber over every analytic ring $\mathcal A$. We claim that the induced map $p\colon \mathcal E^\ocircle \to \AnRing_{(V,\mm)} \cprod \opAssoc$ is a locally coCartesian fibration. As the same is true for $p'$, this follows from the following observations:
\begin{itemize}
	\item For every map $f\colon \mathcal A \to \mathcal B$ of analytic rings over $(V,\mm)$ the functor $f^\natural\colon \End(\D(\mathcal A)) \to \End(\D(\mathcal B))$ restricts to a functor $\End(\D(\mathcal A))^\std \to \End(\D(\mathcal B))^\std$. This follows from \cref{rslt:lower-and-upper-natural-restrict-to-End-std}.

	\item For every analytic ring $\mathcal A$ over $(V,\mm)$, the subcategory $\End(\D(\mathcal A))^\std \subset \End(\D(\mathcal A))$ is stable under the composition monoidal structure. To see this, let $F_1, F_2 \in \End(\D(\mathcal A))^\std$ be given and pick any map $g\colon \mathcal A \to \mathcal B$ of analytic rings over $(V,\mm)$. We need to show that $F_1 \comp F_2$ is steady. The locally coCartesian fibration $p'$ induces a natural map $g^\natural (F_1 \comp F_2) \to g^\natural F_1 \comp g^\natural F_2$ (cf. \cite[Remark 2.4.2.9]{lurie-higher-topos-theory}). It is now enough to show that the natural map
	\begin{align}
		(g,\id)^\natural (F_1 \comp F_2) \to (\id,g)_\natural g^\natural (F_1 \comp F_2) \to (\id,g)_\natural (g^\natural F_1 \comp g^\natural F_2) \label{eq:upper-natural-is-monoidal-on-std}
	\end{align}
	is an isomorphism. This can be checked on the underlying functors and in particular after plugging in any $N \in \D(\mathcal B)$. By steadiness of $F_1$ and $F_2$ we have $g_* (g^\natural F_i)(N) = F_i(g_* N)$ for $i = 1, 2$ (this amounts to saying $(\id,g)_\natural g^\natural F_i = (g,\id)^\natural F_i$). Hence
	\begin{align*}
		&(\id,g)_\natural (g^\natural F_1 \comp g^\natural F_2)(N) = g_* (g^\natural F_1)((g^\natural F_2)(N)) = F_1(g_* (g^\natural F_2)(N)) = F_1(F_2(g_* N)) =\\&\qquad= ((g,\id)^\natural (F_1 \comp F_2))(N),
	\end{align*}
	where we used \cref{rslt:computation-of-upper-natural-if-f-localization} in the last step. This proves \cref{eq:upper-natural-is-monoidal-on-std}.
\end{itemize}
Namely, from these observations we deduce that for every edge $e\colon (\mathcal A, \langle n \rangle) \to (\mathcal B, \langle 1 \rangle)$ in $\AnRing_{(V,\mm)} \cprod \opAssoc$, where $\langle n \rangle \to \langle 1 \rangle$ is the active morphism corresponding to the natural ordering on $\langle n \rangle^\circ$, the map $\mathcal E^\ocircle_e \to \Delta^1$ classifies the functor
\begin{align*}
	(\End(\D(\mathcal A))^\std)^n \to \End(\D(\mathcal B))^\std, \qquad (F_1, \dots, F_n) \mapsto f^\natural (F_1 \comp \dots \comp F_n).
\end{align*}
To show that $p$ is a coCartesian fibration, we use \cite[Remark 2.4.2.9]{lurie-higher-topos-theory} which reduces this claim to showing that for every map $f\colon \mathcal A \to \mathcal B$ of analytic rings over $(V,\mm)$ and all $F_1, F_2 \in \End(\D(\mathcal A))^\std$ the natural map $f^\natural(F_1 \comp F_2) \to f^\natural F_1 \comp f^\natural F_2$ (induced by the locally coCartesian fibration $p$) is an isomorphism. But this follows immediately from \cref{eq:upper-natural-is-monoidal-on-std} for $g = f$ by applying the functor $(\id,g)^\natural$. This finishes the proof that $p$ is a coCartesian fibration and hence defines a coCartesian family of monoidal $\infty$-categories. By \cite[Example 4.8.3.3]{lurie-higher-algebra} it classifies a functor $\AnRing_{(V,\mm)} \to \infcatinf^\ocircle$, proving (iii).

Parts (i) and (ii) were shown in the process of proving (iii). For (i), note that the monoidal structure is exact in both arguments because our endofunctors are required to be exact functors by definition.
\end{proof}

With the functoriality of $\End(\D(-))^\std$ at hand we can finally perform the promised gluing procedure.

\begin{proposition} \label{rslt:End-is-analytic-sheaf-on-AnRing}
The presheaf $\mathcal A \mapsto \End(\D(\mathcal A))^\std$ on the $\infty$-category $\AnRing_{(V,\mm)}$ of affine analytic spaces over $(V,\mm)$ is a sheaf of monoidal $\infty$-categories for the analytic topology.
\end{proposition}
\begin{proof}
The presheaf exists by \cref{rslt:functoriality-of-analytic-steady-endofunctors}. To show that it is a sheaf we can employ the same argument as in \cite[Proposition 10.5]{condensed-mathematics}: Condition (i) of that reference follows from $\natural$-base-change (see \cref{rslt:base-change-for-End-std}) and \cref{rslt:lower-natural-is-fully-faithful-if-f'-localization}. Condition (ii) follows from \cref{rslt:conservativity-for-End-std}.
\end{proof}

\begin{definition}
\begin{defenum}
	\item We define the sheaf
	\begin{align*}
		(\AnSpace_{(V,\mm)})^\opp \to \infcatinf^\ocircle, \qquad X \mapsto \End(\D(X))^\std
	\end{align*}
	on the analytic site of analytic spaces to be the unique sheaf of monoidal $\infty$-categories which extends the sheaf from \cref{rslt:End-is-analytic-sheaf-on-AnRing} (see \cref{rslt:sheaves-on-basis-equiv-sheaves-on-whole-site}).\footnote{A priori there might be set-theoretic issues with this definition, which can be avoided by noting that every analytic space admits a \emph{small} cover by affine spaces.} We denote the monoidal structure on $\End(\D(X))^\std$ by $\comp$ and often abbreviate $FG = F \comp G$.

	\item For a map $f\colon Y \to X$ of analytic spaces over $(V,\mm)$, we denote
	\begin{align*}
		f^\natural\colon \End(\D(X))^\std \to \End(\D(Y))^\std
	\end{align*}
	the ``restriction'' functor of the sheaf $\End(\D(-))^\std$. It is a monoidal functor by construction. If $f$ is steady and qcqs then $f^\natural$ admits a right adjoint
	\begin{align*}
		f_\natural\colon \End(\D(Y))^\std \to \End(\D(X))^\std
	\end{align*}
	given by gluing the affine version of $f_\natural$ along analytic covers (analogous to \cref{rslt:pushforward-of-qcoh-on-analytic-spaces-exists}).
\end{defenum}
\end{definition}

Since $\End(\D(X))^\std$ is in general defined by some large limit, it can be hard to explicitly describe objects inside it. However, by \cref{rslt:upper-natural-on-generators} we can construct a large class of examples in $\End(\D(X))^\std$:

\begin{lemma} \label{rslt:construct-generators-in-End-std-of-space}
There is a natural transformation
\begin{align*}
	h\colon \D(-)^\lnuc \cprod \D(-)^\opp \to \End(\D(-))^\std
\end{align*}
of functors $(\AnSpace_{(V,\mm)})^\opp \to \infcatinf$ such that for every analytic ring $\mathcal A$ over $(V,\mm)$ and all $M \in \D(\mathcal A)^\lnuc$ and $N \in \D(\mathcal A)$ we have $h(M, N) = \IHom(N, -) \tensor M \in \End(\D(\mathcal A))^\std$.
\end{lemma}
\begin{proof}
By \cref{rslt:sheaves-on-basis-equiv-sheaves-on-whole-site} $\End(\D(-))^\std$ is the right Kan extension of its restriction to $\AnRing_{(V,\mm)}$, so it is enough to construct the functor $h$ on affine spaces. In other words, we need to see that the assignment $(M, N) \mapsto \IHom(N, -) \tensor M$ can be upgraded to a natural transformation $\D(-)^\lnuc \cprod \D(-)^\opp \to \End(\D(-))^\std$ of functors $\AnRing_{(V,\mm)} \to \infcatinf$. By abuse of notation we denote $\D(-)^\lnuc \to S$, $\D(-)^\opp \to S$ and $\End(\D(-))^\std \to S$ the coCartesian fibrations over $S = \AnRing_{(V,\mm)}$ classifying the respective functors. By the construction in the proofs of \cref{rslt:functoriality-of-enriched-endofunctors} and \cref{rslt:functoriality-of-analytic-steady-endofunctors}, the coCartesian fibration $\End(\D(-))^\std \to S$ factors over a monomorphism $\End(\D(-)) \injto \End_S(\D(-)^\tensor)$, where $\D(-)^\tensor \to S$ is the coCartesian fibration such that the fiber over each $A \in S$ is the $\infty$-operad $\D(\mathcal A)^\tensor$ exhibiting the self-enrichment of $\D(\mathcal A)$. Here $\End_S(-)$ is defined by the universal property in \cref{eq:def-of-End-S}. In particular, there is a natural map $\D(-)^\tensor \cprod_S \End_S(\D(-)^\tensor) \to \D(-)^\tensor$ of simplicial sets over $S$, which restricts to a map $\D(-) \cprod_S \End(\D(-)) \to \D(-)$ of simplicial sets over $S$. We thus obtain a map of simplicial sets
\begin{align*}
	\D(-) \cprod_S \D(-)^\opp \cprod_S \End(\D(-)) \to \D(-) \cprod_S \D(-)^\opp \to \Ani,
\end{align*}
where the second map comes from the Yoneda embedding of $\D(-)$. We can now argue as in the proof of \cref{rslt:weak-enriched-yoneda-is-functorial} to derive a map
\begin{align*}
	\alpha\colon \D(-) \cprod_S \D(-)^\opp \to \End(\D(-))
\end{align*}
of coCartesian fibrations over $S$, which over each $\mathcal A \in S$ restricts to the map $(M, N) \mapsto \IHom(N, -) \tensor M$. It is now enough to show that $\alpha$ maps coCartesian edges over $S$ to coCartesian edges over $S$, because then the desired natural transformation is obtained by applying the straightening construction to $\alpha$ (and restricting to locally nuclear objects and steady endofunctors using \cref{rslt:properties-of-End-std}). The claim on $\alpha$ amounts to saying that for every map $f\colon \mathcal A \to \mathcal B$ in $S$ and all $M, N \in \D(\mathcal A)$ the natural morphism
\begin{align*}
	f^\natural\big(\IHom_{\mathcal A}(N, -) \tensor_{\mathcal A} M\big) \isoto \IHom_{\mathcal B}(f^* N, -) \tensor_{\mathcal B} f^*M
\end{align*}
is an isomorphism in $\End(\D(\mathcal B))$. But this is precisely the statement of \cref{rslt:upper-natural-on-generators}.
\end{proof}

\begin{definition}
Let $X$ be an analytic space over $(V,\mm)$. For every $\mathcal M \in \D(X)^\lnuc$ and $\mathcal N \in \D(X)$ we denote the object $h(\mathcal M, \mathcal N)$ constructed in \cref{rslt:construct-generators-in-End-std-of-space} by
\begin{align*}
	\IHom(\mathcal N, -) \tensor \mathcal M \in \End(\D(X))^\std.
\end{align*}
In the case that $\mathcal N = \ri_X$ we will often denote $h(\mathcal M, \ri_X) =: \mathcal M \in \End(\D(X))^\std$.
\end{definition}

\begin{remark} \label{rmk:D-nuc-into-E-X-in-general-not-fully-faithful}
For every analytic space $X$ the functor
\begin{align*}
	\D(X)^\opp \injto \End(\D(X))^\std, \qquad \mathcal N \mapsto \IHom(\mathcal N, -)
\end{align*}
is fully faithful: Since both $\D(-)^\opp$ and $\End(\D(-))^\std$ are sheaves, this can be checked on affine spaces, where it reduces to \cref{rslt:weak-enriched-yoneda-is-functorial}. It is unclear if also the functor
\begin{align*}
	\D(X)^\lnuc \to \End(\D(X))^\std, \qquad \mathcal M \mapsto - \tensor \mathcal M
\end{align*}
is fully faithful: this boils down to the question whether $\D(-)^\lnuc$ is a sheaf (cf. \cref{rmk:is-D-nuc-a-sheaf}).
\end{remark}

It is straightforward to generalize the previous results from the affine case to general analytic spaces over $(V,\mm)$:

\begin{proposition} \label{rslt:results-on-End-of-analytic-space-enum}
\begin{propenum}
	\item \label{rslt:properties-of-End-of-analytic-space} For every analytic space $X$ over $(V,\mm)$, $\End(\D(X))^\std$ is a stable monoidal $\infty$-category which has all small colimits. The monoidal structure $\comp$ is exact in both arguments.

	\item \label{rslt:End-analytic-space-upper-natural-properties} For a map $f\colon Y \to X$ of analytic spaces over $(V,\mm)$ the functor $f^\natural\colon \End(\D(X))^\std \to \End(\D(Y))^\std$ is monoidal and preserves all small colimits.

	\item \label{rslt:End-analytic-space-natural-base-change} Let
	\begin{center}\begin{tikzcd}
		Y' \arrow[r,"g'"] \arrow[d,"f'"] & Y \arrow[d,"f"]\\
		X' \arrow[r,"g"] & X
	\end{tikzcd}\end{center}
	be a Cartesian square of analytic spaces over $(V,\mm)$ such that $f$ is steady and qcqs. Then the base-change morphism
	\begin{align*}
		g^\natural f_\natural \isoto f'_\natural g'^\natural
	\end{align*}
	is an isomorphism of functors from $\End(\D(Y))^\std$ to $\End(\D(X'))^\std$.

	\item \label{rslt:computation-of-f-lower-nat-f-upper-nat-for-affine-target} Let $f\colon Y \to X$ be a qcqs steady map of analytic spaces over $(V,\mm)$. Then for every $\mathcal F \in \End(\D(X))^\std$ there is a natural isomorphism $\mathcal F \comp (f_\natural \ri_Y) = f_\natural f^\natural \mathcal F$. Moreover, if $X = \AnSpec \mathcal A$ is affine then the underlying functor of $f_\natural \ri_Y$ is $f_* f^*$.
\end{propenum}
\end{proposition}
\begin{proof}
Parts (i) and (ii) follow directly from \cref{rslt:properties-of-End-std,rslt:functoriality-of-analytic-steady-endofunctors}. Part (iii) follows from \cref{rslt:base-change-for-End-std} by an analogous argument as in \cref{rslt:steady-base-change-holds}.

To prove (iv) we note that by adjunction there is a natural map $\mathcal F \comp (f_\natural \ri_Y) \to f_\natural f^\natural \mathcal F$. To show that it is an isomorphism we can now assume that $X = \AnSpec \mathcal A$ is affine. Now note that since everything commutes with finite limits and $f$ is qcqs, we can w.l.o.g. assume that $Y = \AnSpec \mathcal B$ is affine as well. By definition of steady objects we have $f_\natural f^\natural \mathcal F = (f,\id)_\natural (f,\id)^\natural \mathcal F$ and hence the claim follows immediately from the explicit computations in \cref{rslt:lower-upper-natural-functors-for-Fun-of-analytic-rings}.
\end{proof}

For every analytic space $X$ one would expect the $\infty$-category $\End(\D(X))^\std$ to act on $\D(X)$ via ``evaluation'', i.e. for every $\mathcal F \in \End(\D(X))^\std$ and $\mathcal M$ there should be a functorial way of computing $\mathcal F(\mathcal M) \in \D(X)$. We would in particular expect $\restrict{\mathcal F(\mathcal M)}U = (\restrict{\mathcal F}U)(\restrict{\mathcal M}U)$ for every steady subspace $U$. This seems to be false for steady endofunctors, but we can fix it by restricting to an even smaller subcategory:

\begin{definition}
\begin{defenum}
	\item Let $\mathcal A$ be an analytic ring over $(V,\mm)$. We define
	\begin{align*}
		\mathcal E(\mathcal A) \subset \End(\D(\mathcal A))^\std
	\end{align*}
	to be the smallest full subcategory which is stable under small colimits and $\comp$ and which contains all objects of the form $f_\natural M$ for steady maps $f\colon \mathcal A \to \mathcal B$ and nuclear $M \in \D(\mathcal B)^\nuc$.

	\item Let $X$ be an analytic space over $(V,\mm)$. We define
	\begin{align*}
		\mathcal E(X) \subset \End(\D(X))^\std
	\end{align*}
	to be the full subcategory of those $\mathcal F \in \End(\D(X))^\std$ such that for every map $f\colon Y \to X$ from an affine analytic space $Y = \AnSpec \mathcal B$ we have $f^\natural \mathcal F \in \mathcal E(\mathcal B)$.
\end{defenum}
\end{definition}

\begin{proposition} \label{rslt:properties-of-E-of-X}
\begin{propenum}
	\item For every analytic space $X$ over $(V,\mm)$, $\mathcal E(X)$ is a stable monoidal $\infty$-category which has all small colimits. The monoidal structure $\comp$ preserves all small colimits in both arguments. If $X = \AnSpec \mathcal A$ is affine then $\mathcal E(X) = \mathcal E(\mathcal A)$.

	\item \label{rslt:E-of-X-stability-under-natural-functors} Let $f\colon Y \to X$ be a morphism of analytic spaces over $(V,\mm)$. Then $f^\natural$ restricts to a functor $f^\natural\colon \mathcal E(X) \to \mathcal E(Y)$. If $f\colon U \injto X$ is qcqs and identifies $U$ as a steady subspace of $X$, then $f_\natural$ restricts to a functor $f_\natural\colon \mathcal E(U) \to \mathcal E(X)$.
	\item Let $f\colon Y \to X$ be a qcqs steady morphism of analytic spaces over $(V,\mm)$. Then for all $\mathcal G \in \End(\D(Y))^\std$ and all $\mathcal F \in \mathcal E(X)$ the natural morphisms
	\begin{align*}
		\mathcal F \comp f_\natural \mathcal G \isoto f_\natural(f^\natural \mathcal F \comp \mathcal G), \qquad f_\natural \mathcal G \comp \mathcal F \isoto f_\natural(\mathcal G \comp f^\natural \mathcal F)
	\end{align*}
	are isomorphisms in $\End(\D(X))^\std$. In particular there is a natural isomorphism $\mathcal F \comp f_\natural \ri_Y = f_\natural \ri_Y \comp \mathcal F$.

	\item \label{rslt:E-of-X-is-a-sheaf} The assignment $X \mapsto \mathcal E(X)$ is a sheaf of monoidal $\infty$-categories on the analytic site of analytic spaces over $(V,\mm)$.

	\item \label{rslt:evaluation-map-on-E-of-X} There is a natural coCartesian family of $\opLM$-monoidal $\infty$-categories
	\begin{align*}
		ED^\tensor_{(V,\mm)} \to \AnSpace^\std_{(V,\mm)} \cprod \opLM
	\end{align*}
	whose fiber over every analytic space $X$ over $(V,\mm)$ exhibits $\D(X)$ as tensored over $\mathcal E(X)$. If $X = \AnSpec \mathcal A$ is affine then the induced functor $\mathcal E(\mathcal A) \cprod \D(\mathcal A) \to \D(\mathcal A)$ is the evaluation functor $(F, M) \mapsto F(M)$. Here $\AnSpace^\std_{(V,\mm)}$ denotes the $\infty$-category of analytic spaces over $(V,\mm)$ where we only allow steady morphisms.
					\end{propenum}
\end{proposition}
\begin{proof}
In part (i) it is clear that $\mathcal E(X)$ is a stable monoidal $\infty$-category and has all small colimits. We now show that $\comp$ preserves all small colimits in both arguments. To see this, we can pass to an analytic cover to assume that $X = \AnSpec \mathcal A$ is affine. In this case it is enough to observe that for every $F \in \mathcal E(X)$ the underlying functor $\D(\mathcal A) \to \D(\mathcal A)$ preserves all small colimits. Since this property is stable under small colimits and compositions of functors, we reduce to the case $F = f_\natural N$ for some steady map $f\colon \mathcal A \to \mathcal B$ and some nuclear module $N \in \D(\mathcal B)^\nuc$. Then \cref{rslt:properties-of-lower-natural} implies that the underlying functor of $f_\natural N$ is $M \mapsto f_* (N \tensor_{\mathcal B} f^* M)$, which indeed preserves all small colimits.

To finish the proof of (i) it remains to show that if $X = \AnSpec \mathcal A$ is affine then $\mathcal E(X) = \mathcal E(\mathcal A)$. In other words, we need to show that for every $F \in \mathcal E(\mathcal A)$ and every map $f\colon \mathcal A \to \mathcal B$ of analytic rings over $(V,\mm)$, we have $f^\natural F \in \mathcal E(\mathcal B)$. Since $f^\natural$ is monoidal and preserves all small colimits, we can reduce to the case that $F = g_\natural N$ for some steady map $g\colon \mathcal A \to \mathcal C$ and some $N \in \D(\mathcal C)^\nuc$. Denote $f'$ and $g'$ the base-changes of $f$ and $g$ along $g$ and $f$, respectively. Then by \cref{rslt:base-change-for-End-std} we have $f^\natural F = g'_\natural(f'^* N)$, so we conclude by \cref{rslt:nuclear-modules-stable-under-pullback}.

For part (ii), the claim about $f^\natural$ is obvious from the definitions. Now assume that $f\colon U \injto X$ is qcqs and identifies $U \subset X$ as a steady subspace. Given any $\mathcal F \in \mathcal E(U)$ we need to show that $f_\natural \mathcal F \in \mathcal E(X)$. By base-change (see \cref{rslt:End-analytic-space-natural-base-change}) this reduces immediately to the case that $X = \AnSpec \mathcal A$ is affine. Since $\mathcal E(X)$ is stable under finite limits, we can further reduce to the case that $U = \AnSpec \mathcal B$ is affine and hence $f$ is given by a steady localization $\mathcal A \to \mathcal B$. Now the property $f_\natural \mathcal F \in \mathcal E(X)$ holds for the generators $\mathcal F = g_\natural N$ and is stable under small colimits in $\mathcal F$, since $f_\natural$ preserves small colimits. It remains to show that the property is stable under composition, i.e. if $F, G \in \mathcal E(\mathcal B)$ satisfy $f_\natural F, f_\natural G \in \mathcal E(\mathcal A)$, then $f_\natural(F \comp G) \in \mathcal E(\mathcal A)$. But $f_\natural(F \comp G) = f_\natural F \comp f_\natural G$ (there is a natural morphism from right to left so it can be checked on underlying functors, where it follows from $f^* f_* = \id$).

We now prove (iii), so let $f\colon Y \to X$ be steady and qcqs and let $\mathcal F \in \mathcal E(X)$ and $\mathcal G \in \End(\D(Y))^\std$ be given. By \cref{rslt:End-analytic-space-natural-base-change} we can reduce to the case that $X = \AnSpec \mathcal A$ is affine. Since everything commutes with finite limits we can further reduce to the case that $Y = \AnSpec \mathcal B$ is also affine. Moreover, since $f$ is qcqs, $f_\natural$ preserves all small colimits (this can be checked on a finite affine cover of $Y$; but in the affine case it follows from \cref{rslt:lower-and-upper-natural-restrict-to-End-std}). Using also that $f^\natural$ is monoidal it follows easily that the set of $\mathcal F \in \mathcal E(X)$ which satisfy the claimed isomorphisms (for all $\mathcal G$) is closed under small colimits and under compositions. We can therefore reduce to the case that $\mathcal F = g_\natural N$ for some steady map $g\colon \mathcal A \to \mathcal C$ and some $N \in \D(\mathcal C)^\nuc$. Now the claimed isomorphisms can be checked on underlying functors $\D(\mathcal A) \to \D(\mathcal A)$ and in particular after plugging in any $M \in \D(\mathcal A)$. Let $f'\colon \mathcal C \to \mathcal B \tensor_{\mathcal A} \mathcal C$ and $g'\colon \mathcal B \to \mathcal B \tensor_{\mathcal A} \mathcal C$ denote the base-changes of $f$ and $g$. The second isomorphism follows from the following chain of equations:
\begin{align*}
	&(f_\natural \mathcal G \comp \mathcal F)(M) = f_* \mathcal G f^* g_* (N \tensor g^* M) = f_* \mathcal G g'_* f'^* (N \tensor g^* M) = f_* \mathcal G g'_* (f'^* \mathcal N \tensor f'^* g^* M) =\\&\qquad= f_* \mathcal G g'_* (f'^* N \tensor g'^* f^* M) = f_\natural(\mathcal G \comp g'_\natural(f'^\natural N))(M) = f_\natural(\mathcal G \comp f^\natural \mathcal F)(M).
\end{align*}
The first isomorphism can be proved similarly using \cref{rslt:nucelar-modules-are-steady}.

Part (iv) follows formally from the sheafiness of $\End(\D(-))^\std$ using that $f_\natural$ preserves $\mathcal E(-)$ if $f$ is a qcqs steady subspace by (ii).

It remains to prove (v). By the sheafiness of $\mathcal E(-)$ and $\D(-)$ the desired coCartesian family of $\opLM$-monoidal $\infty$-categories is automatically a sheaf of $\opLM$-monoidal $\infty$-categories for the analytic site. Therefore by \cref{rslt:sheaves-on-basis-equiv-sheaves-on-whole-site} it can be constructed on the basis of affine analytic spaces, so we only need to construct the coCartesian family $ED^\tensor_{(V,\mm)} \to S \cprod \opLM$ with $S = \AnRing^\std_{(V,\mm)}$. As in the proof of \cref{rslt:functoriality-of-analytic-steady-endofunctors} we apply \cite[Corollary 4.7.1.40]{lurie-higher-algebra} to the $\infty$-category $(\infcatinf)_{/S}$, viewed as enriched over itself via the Cartesian monoidal structure. We then obtain a categorical fibration $ED'^\tensor \to S \cprod \opLM$ whose fiber over every $\mathcal A \in S$ exhibits $\D(\mathcal A)^\tensor$ as tensored over $\Fun(\D(\mathcal A)^\tensor, \D(\mathcal A)^\tensor)$ via the composition monoidal structure. Let $ED^\tensor_{(V,\mm)} \subset ED'^\tensor$ be the subcategory where we restrict to $\mathcal E(-)$ over $\mathfrak a \in \opLM$ and to $\D(-)$ over $\mathfrak m \in \opLM$. It remains to show that the induced categorical fibration $p\colon ED^\tensor_{(V,\mm)} \to S \cprod \opLM$ is a coCartesian fibration. By the standard arguments (as in the proofs above) this boils down to showing the following: Given a steady map $f\colon \mathcal A \to \mathcal B$ of analytic rings over $(V,\mm)$, some $F \in \mathcal E(\mathcal A)$ and any $M \in \D(\mathcal A)$, the natural map $f^* F(M) \isoto (f^\natural F)(f^* M)$ induced by $p$ is an isomorphism. This condition is stable under small colimits and compositions in $F$, so we can assume that $F$ is a generator of $\mathcal E(\mathcal A)$, i.e. of the form $F = g_\natural \mathcal N$ for some steady map $g\colon \mathcal A \to \mathcal C$ and some $N \in \D(\mathcal C)^\nuc$. Letting $f'$ and $g'$ be the base-changes of $f$ and $g$ along $g$ and $f$ respectively, we have $f^\natural F = g'_\natural (f'^* N)$. Thus
\begin{align*}
	&(f^\natural F)(f^* M) = g'_*(f'^* N \tensor g'^* f^* M) = g'_*(f'^* N \tensor f'^* g^* M) = g'_* f'^* (N \tensor g^* M) =\\&\qquad= f^* g_* (N \tensor g^* M) = f^* F(M),
\end{align*}
as desired.
\end{proof}

\begin{definition} \label{def:evaluation-map-on-E-of-X}
For every analytic space $X$ over $(V,\mm)$, we denote
\begin{align*}
	\ev\colon \mathcal E(X) \cprod \D(X) \to \D(X), \qquad (\mathcal F, \mathcal M) \mapsto \mathcal F(\mathcal M)
\end{align*}
the induced functor from \cref{rslt:evaluation-map-on-E-of-X} and call it the \emph{evaluation functor on $X$}. It is compatible with the monoidal structure on $\mathcal E(X)$, preserves all small colimits in both arguments and commutes with pullback along steady maps of analytic spaces, in particular with restriction to steady subspaces.
\end{definition}

\begin{remark}
Given a geometry blueprint $G$ over $(V,\mm)$ (see \cref{def:geometry-blueprint}), all of the above constructions and results apply verbatim to $G$-analytic spaces instead of analytic spaces over $(V,\mm)$ by applying the analytification functor from \cref{def:analytification}. In particular, for every $G$-analytic space $X$ we define $\mathcal E(X) := \mathcal E(X^\an)$.
\end{remark}

We have now developed the theory of endofunctors on analytic spaces. We close this section by discussing the functoriality of this construction with respect to the almost setup $(V,\mm)$.

\begin{proposition} \label{rslt:functoriality-of-steady-endofunctors-over-AlmSetup}
Let $\AnSpace_{(-)} \to \AlmSetup$ be the coCartesian fibration classifying the functor from \cref{rslt:functoriality-of-analytic-spaces}. There is a natural coCartesian family of $\opLM$-monoidal $\infty$-categories
\begin{align*}
	ED^\tensor \to \AnSpace^\std_{(-)} \cprod \opLM
\end{align*}
whose fiber over every analytic space $X$ over $(V,\mm)$ exhibits $\D(X)$ as tensored over $\mathcal E(X)$. Moreover, for every morphism of almost setups $\varphi\colon (V,\mm) \to (V',\mm')$ we have:
\begin{propenum}
	\item For every analytic ring $\mathcal A$ over $(V,\mm)$, the induced functor $\varphi^\natural\colon \mathcal E(\mathcal A) \to \mathcal E(\varphi^*\mathcal A)$ acts as $F \mapsto \varphi^* \comp F \comp \varphi_*$. It is monoidal and preserves all small colimits.

	\item For every steady map $f\colon \mathcal A \to \mathcal B$ of analytic rings over $(V,\mm)$ and every $M \in \D(\mathcal B)^\nuc$ the natural morphism $\varphi^\natural f_\natural M \isoto (\varphi^* f)_\natural (\varphi^* M)$ is an isomorphism.
\end{propenum}
\end{proposition}
\begin{proof}
As in the proof of \cref{rslt:evaluation-map-on-E-of-X} it is enough to construct the desired coCartesian fibration over $\AnRing_{(-)}^\std$, where it can be done by a very similar construction. We leave the $\infty$-categorical details to the reader and only check the relevant compatibilities, which are as follows:
\begin{enumerate}[(a)]
	\item For every morphism $\varphi\colon (V,\mm) \to (V',\mm')$ and every analytic ring $\mathcal A$ over $(V,\mm)$ with induced morphism $\varphi_{\mathcal A}\colon \mathcal A \to \varphi^*\mathcal A$ in $\AnRing_{(-)}$, there is a natural functor $\varphi^\natural\colon \mathcal E(\mathcal A) \to \mathcal E(\varphi^*\mathcal A)$ which acts as $F \mapsto \varphi^* \comp F \comp \varphi_*$ on the underlying functors.

	\item The functors $\varphi^\natural\colon \mathcal E(\mathcal A) \to \mathcal E(\varphi^*\mathcal A)$ constructed in (a) are monoidal and commute with the functors $f^\natural$ in the obvious sense. Note that all the required commutativity relations can be checked on the $1$-categorical level and then even on the underlying functors in $\mathcal E(\mathcal A)$ because they come from natural transformations induced by $p$.

	\item For every $F \in \mathcal E(\mathcal A)$ and every $M \in \D(\mathcal A)$ we have $(\varphi^\natural F)(\varphi^* M) = \varphi^* F(M)$.
\end{enumerate}
To prove (a) we can construct $\varphi^\natural$ (and its right adjoint $\varphi_\natural$) as in \cref{rslt:lower-upper-natural-functors-for-Fun-of-analytic-rings}; in fact $\varphi^*\colon \D(\mathcal A) \to \D(\varphi^*\mathcal A)$ behaves just like a localization. It is then easy to see that this functor $\varphi^\natural$ commutes with the functors $f^\natural$ in the obvious sense (e.g. by looking at the right adjoints, which are given very explicitly). To check that $\varphi^\natural$ is monoidal, it suffices to check that for every $F \in \mathcal E(\mathcal A)$ and every $M \in \D(\mathcal A)$ we have $(\varphi^\natural F)(\varphi^* M) = \varphi^* F(M)$ via the natural morphism (as in the proof of \cref{rslt:evaluation-map-on-E-of-X}). This condition is stable under colimits and compositions in $F$, so we can assume that $F = g_\natural N$ for some steady map $g\colon \mathcal A \to \mathcal B$ and some $N \in \D(\mathcal B)^\nuc$. Using that $\varphi^*$ behaves like a steady morphism (see \cref{rslt:almost-localization-is-steady}) we can use exactly the same computation as in the proof of \cref{rslt:evaluation-map-on-E-of-X} provided that (ii) holds (here we implicitly use \cref{rslt:almost-localization-preserves-steadiness} to see that the statement makes sense because $\varphi^* g$ is steady). In order to see (ii) we can assume that $\varphi$ is localizing and we compute for all $M \in \mathcal E(\varphi^* \mathcal A)$ (again using \cref{rslt:almost-localization-is-steady})
\begin{align*}
	&(\varphi^\natural f_\natural N)(M) = \varphi^* f_*(N \tensor f^* \varphi_* M) = (\varphi^*f)_* \varphi^* (N \tensor f^* \varphi_* M) =\\&\qquad= (\varphi^* f)_* (\varphi^* N \tensor \varphi^* f^* \varphi_* M) = (\varphi^* f)_* (\varphi^* N \tensor (\varphi^*f)^* M) =\\&\qquad= (\varphi^* f)_* \varphi^* (N \tensor \varphi_* (\varphi^*f)^* M) = (\varphi^*f)_\natural(\varphi^* N)(M),
\end{align*}
as desired. Claim (c) follows by the same computation as in the proof of \cref{rslt:evaluation-map-on-E-of-X} using \cref{rslt:almost-localization-is-steady}.
\end{proof}

\begin{definition}
For every morphism $\varphi\colon (V,\mm) \to (V',\mm')$ of almost setups and every analytic space $X$ over $(V,\mm)$ we denote
\begin{align*}
	\varphi^\natural\colon \mathcal E(X) \to \mathcal E(\varphi^* X)
\end{align*}
the functor constructed in \cref{rslt:functoriality-of-steady-endofunctors-over-AlmSetup}. It is monoidal, preserves all small colimits and commutes with the evaluation functor from \cref{def:evaluation-map-on-E-of-X}.
\end{definition}

\subsection{Descendable Morphisms} \label{sec:andesc.descmorph}

Fix an almost setup $(V,\mm)$. We now introduce descendable morphisms of analytic spaces over $(V,\mm)$, which form the basis of the descent theory we wish to develop in subsequent sections. The following definitions and results can all be seen as analogs of the notion of descendable algebras in \cite[\S3]{akhil-galois-group-of-stable-homotopy}. In fact, the construction of the steady endofunctor $\infty$-categories $\mathcal E(X)$ in \cref{sec:andesc.endofun} allows us to more or less directly carry over the definitions and results from loc. cit. Let us start with the core definitions:

\begin{definition} \label{def:descendable-morphism-of-analytic-spaces}
Let $f\colon Y \to X$ be a qcqs steady morphism of analytic spaces over $(V,\mm)$. We denote
\begin{align*}
	\langle f \rangle \subset \mathcal E(X)
\end{align*}
the smallest full subcategory which contains $f_\natural \ri_Y$ and is stable under finite (co)limits, retracts and composition. We say that $f$ is \emph{descendable} if $\langle f \rangle$ contains $\ri_X$.

A steady morphism $f^\sharp\colon \mathcal A \to \mathcal B$ of analytic rings over $(V,\mm)$ is called \emph{descendable} if the corresponding morphism $f\colon \AnSpec \mathcal B \to \AnSpec \mathcal A$ is descendable.
\end{definition}

\begin{definition}
Given a morphism $f\colon Y \to X$ of analytic spaces over $(V,\mm)$, we say that \emph{quasicoherent sheaves descend along $f$} if, denoting by $Y_\bullet \to X$ the Čech nerve of $f$, the natural functor
\begin{align*}
	\D(X) \isoto \varprojlim_{n\in\Delta} \D(Y_n)
\end{align*}
is an equivalence of $\infty$-categories. The same definition applies to a morphism of analytic rings.
\end{definition}

Of course, for the terminology to make sense, we need to show that quasicoherent sheaves descend along descendable morphisms:

\begin{proposition} \label{rslt:descendable-implies-limit-of-categories}
Let $f\colon Y \to X$ be a descendable morphism of analytic spaces over $(V,\mm)$. Then quasicoherent sheaves descend along $f$.
\end{proposition}
\begin{proof}
We can formally reduce to the case that $X = \AnSpec \mathcal A$ is affine, noting that by \cref{rslt:properties-of-E-of-X,rslt:End-analytic-space-natural-base-change} the descendability of $f$ remains unchanged if we pass to an analytic cover. Now we argue as in \cite[Proposition 3.22]{akhil-galois-group-of-stable-homotopy}, i.e. we apply Lurie's Beck-Chevalley condition \cite[Corollary 4.7.5.3]{lurie-higher-algebra} to the augmented cosimplicial $\infty$-category $\D(X)^\opp \to \D(Y_\bullet)^\opp$.

We first check condition (1) of loc. cit. Let $M^\bullet$ be an $f^*$-split cosimplicial object of $\D(\mathcal A)$. Then $(f_\natural \ri_Y)(M^\bullet) = f_* f^* (M^\bullet)$ is split, in particular the Tot tower of $(f_\natural \ri_Y)(M^\bullet)$ is constant. It follows easily that the Tot tower of $F(M^\bullet)$ is constant for every $F \in \langle f \rangle$. Choosing $F = \id_{\mathcal A} \in \langle f \rangle$ (by definition of descendable morphisms) we deduce that the Tot tower of $M^\bullet$ is constant. But then it is clear that $f^*$ commutes with the totalization of $M^\bullet$, as desired.

Condition (2) of \cite[Corollary 4.7.5.3]{lurie-higher-algebra} follows from \cref{rslt:steady-base-change-holds} using that $f$ (and hence any base-change of an $n$-fold composition of $f$ with itself) is steady and qcqs.

It remains to show that the functor $f^*$ is conservative. For that it is enough to show that $f_\natural \ri_Y = f_*f^*$ is conservative. If $M \in \D(\mathcal A)$ is such that $(f_\natural \ri_Y)(M) \isom 0$, then it follows immediately that $F(M) \isom 0$ for all $F \in \langle f \rangle$. Choosing $F = \id_{\mathcal A}$ yields $M \isom 0$, as desired.
\end{proof}

Next up we provide a few different characterizations of descendable morphisms which come in handy in many situations. Among others, this also allows us to define a ``complexity'' of descendable morphisms, called the \emph{index}.

\begin{definition}
Let $f\colon Y \to X$ be a qcqs steady morphism of analytic spaces over $(V,\mm)$. Then a morphism $\mathcal F_1 \to \mathcal F_2$ in $\mathcal E(X)$ is called \emph{$f$-zero} if the induced morphism $f_\natural f^\natural \mathcal F_1 \to f_\natural f^\natural \mathcal F_2$ is zero.
\end{definition}

\begin{proposition} \label{rslt:descendable-criteria}
Let $f\colon Y \to X$ be a qcqs steady morphism of analytic spaces over $(V,\mm)$. Then the following are equivalent:
\begin{propenum}
	\item $f$ is descendable.

	\item \label{rslt:descendable-nilpotent-criterion} There exists an integer $s \ge 0$ such that the composition of any $s$ consecutive $f$-zero morphisms in $\mathcal E(X)$ is zero.

	\item \label{rslt:descendable-Un-criterion} Letting
	\begin{align*}
		\mathcal K_f := \fib(\ri_X \to f_\natural \ri_Y),
	\end{align*}
	there is some $n \ge 0$ such that the natural morphism $\mathcal K_f^n \to \ri_X$ is zero.
\end{propenum}
\end{proposition}
\begin{proof}
The following proof is mostly analogous to the proof of \cite[Proposition 3.27]{akhil-galois-group-of-stable-homotopy}.

We first show that (i) implies (ii). For every $\mathcal G \in \mathcal E(X)$ let $\mathcal I_{\mathcal G}$ be the collection of morphisms $\mathcal F_1 \to \mathcal F_2$ in the homotopy category of $\mathcal E(X)$ such that $\mathcal F_1 \mathcal G \to \mathcal F_2 \mathcal G$ is zero. By \cref{rslt:computation-of-f-lower-nat-f-upper-nat-for-affine-target}, $\mathcal I_{f_\natural \ri_Y}$ consists precisely of the $f$-zero morphisms, so we want to show that $\mathcal I_{f_\natural \ri_Y}^s = 0$ for some $s \ge 0$. We note the following:
\begin{itemize}
	\item For $\mathcal G, \mathcal G' \in \mathcal E(X)$ we have $\mathcal I_{\mathcal G'} \subset \mathcal I_{\mathcal G' \mathcal G}$.

	\item If $\mathcal G'$ is a retract of $\mathcal G$ in $\mathcal E(X)$ then $\mathcal I_{\mathcal G} \subset \mathcal I_{\mathcal G'}$.

	\item If $\mathcal G' \to \mathcal G \to \mathcal G''$ is a cofiber sequence in $\mathcal E(X)$ then $\mathcal I_{\mathcal G'} \mathcal I_{\mathcal G''} \subset \mathcal I_{\mathcal G}$. To prove this let $\varphi''\colon \mathcal F_1 \to \mathcal F_2$ be a $\mathcal G''$-zero morphism and $\varphi\colon \mathcal F_2 \to \mathcal F_3$ a $\mathcal G'$-zero morphism in $\mathcal E(X)$. We need to show that the composition of $\varphi'$ and $\varphi''$ is $\mathcal G$-zero. We have a homotopy coherent diagram
	\begin{center}\begin{tikzcd}
		\mathcal F_1 \mathcal G' \arrow[r] \arrow[d] & \mathcal F_2 \mathcal G' \arrow[r] \arrow[d] & \mathcal F_3 \mathcal G' \arrow[d] \\
		\mathcal F_1 \mathcal G \arrow[r] \arrow[d] & \mathcal F_2 \mathcal G \arrow[r] \arrow[d] & \mathcal F_3 \mathcal G \arrow[d] \\
		\mathcal F_1 \mathcal G'' \arrow[r] & \mathcal F_2 \mathcal G'' \arrow[r] & \mathcal F_3 \mathcal G''
	\end{tikzcd}\end{center}
	in which the vertical arrows are cofiber sequences. From the facts that $\mathcal F_1 \mathcal G'' \to \mathcal F_2 \mathcal G''$ is zero and that the middle column in the above diagram is a cofiber sequence we deduce that the morphism $\mathcal F_1 \mathcal G \to \mathcal F_2 \mathcal G$ factors over $\mathcal F_1 \mathcal G \to \mathcal F_2 \mathcal G'$. But then the composition $\mathcal F_1 \mathcal G \to \mathcal F_2 \mathcal G \to \mathcal F_3 \mathcal G$ factors as $\mathcal F_1 \mathcal G \to \mathcal F_2 \mathcal G' \xto{0} \mathcal F_3 \mathcal G' \to \mathcal F_3 \mathcal G$ and hence is zero.
\end{itemize}
It follows that for any $\mathcal G \in \mathcal E(X)$ and any $\overline{\mathcal G}$ in the full $\infty$-subcategory of $\mathcal E(X)$ generated by $\mathcal G$ under finite (co)limits, retracts and compositions we have $\mathcal I_{\mathcal G}^s \subset \mathcal I_{\overline{\mathcal G}}$ for some $s \ge 0$. Letting $\mathcal G = f_\natural \ri_Y$ and $\overline{\mathcal G} = \ri_X$ (which is possible by the assumption that $f$ is descendable) we deduce $\mathcal I_{f_\natural \ri_Y}^s \subset \mathcal I_{\ri_X} = 0$ for some $s \ge 0$, as desired.

To prove that (ii) implies (iii) note that it is clearly enough to show that $\mathcal K_f \to \ri_X$ is $f$-zero (then the same is true for $\mathcal K_f^k \to \mathcal K_f^{k-1}$ for all $k > 0$ so that $\mathcal K_f^n \to \ri_X$ is the composition of $n$ consecutive $f$-zero maps). To see that $\mathcal K_f \to \ri_X$ is $f$-zero, it is enough to show that $f^\natural \mathcal K_f \to \ri_Y$ is zero in $\mathcal E(Y)$. Letting $g\colon Y \cprod_X Y \to Y$ denote the projection to the first factor, we see (using steady base-change, \cref{rslt:End-analytic-space-natural-base-change}) that applying $f^\natural$ to the map $\ri_X \to f_\natural \ri_Y$ yields the map $\ri_Y \to g_\natural \ri_{Y\cprod_X Y}$. It follows that $f^\natural \mathcal K_f = \mathcal K_g$. Thus we can replace $f$ by $g$ to reduce to the following statement: If $f$ is split then $K_f \to \ri_X$ is zero. But $f$ being split implies that $\ri_X \to f_\natural \ri_Y$ is split, hence the claim follows.

We now prove that (iii) implies (i). Clearly $\mathcal K_f^n = \fib(\mathcal K_f^{n-1} \to \mathcal K_f^{n-1} f_\natural \ri_Y)$ for all $n > 0$. By induction on $n$ it follows that $\mathcal K_f^n f_\natural \ri_Y \in \langle f \rangle$ for all $n \ge 0$. Thus for all $n > 0$ we have $\cofib(\mathcal K_f^n \to \mathcal K_f^{n-1}) \in \langle f \rangle$ because this cofiber equals $\mathcal K_f^{n-1} f_\natural \ri_Y$. By chasing cofiber sequences we deduce that $\cofib(\mathcal K_f^n \to \ri_X) \in \langle f \rangle$. But by assumption $\mathcal K_f^n \to \ri_X$ is zero, so that $\cofib(\mathcal K_f^n \to \ri_X) = \ri_X \dsum \mathcal K_f^n[1]$. Thus $\ri_X$ is a retract of $\cofib(\mathcal K_f^n \to \ri_X)$ and consequently lies in $\langle f \rangle$ as well.
\end{proof}

\begin{remarks}
\begin{remarksenum}
			\item The criterion \cref{rslt:descendable-Un-criterion} is analogous to \cite[Definition 11.18]{bhatt-scholze-witt}.

	\item \label{rslt:f-zero-only-depends-on-index} From the proof of \cref{rslt:descendable-criteria} it follows that the integer $n$ can be chosen as $n = s$. The proof also implies that for every $n \ge 0$ there is some $s(n) \ge 0$ (only depending on $n$) such that if $f\colon Y \to X$ satisfies \cref{rslt:descendable-Un-criterion} then $f$ satisfies \cref{rslt:descendable-nilpotent-criterion} with $s = s(n)$.

	More precisely, we can construct $s(n)$ as follows: By the last part of the proof of \cref{rslt:descendable-criteria} there is a representation of $\ri_X \in \mathcal E(X)$ as a retract of a finite limit of the functors $(f_\natural \ri_Y)^k$ for $k \ge 1$ and $0$. Moreover, the shape of this finite limit depends only on the index of descendability $n$. We can view this limit as a successive chain of cofiber sequences $\mathcal F_1 \to \mathcal F_2 \to \mathcal F_3$ in $\mathcal E(X)$, where both $\mathcal F_1$ and $\mathcal F_3$ are either equal to $(f_\natural \ri_Y)^k$ for some $k \ge 1$ or have been constructed by a previous cofiber sequence. By $s((f_\natural \ri_Y)^k) := 1$ and $s(\mathcal F_2) := s(\mathcal F_1) + s(\mathcal F_3)$ we can inductively define $s(\mathcal F)$ for every object appearing in the limit construction. Then the proof of \cref{rslt:descendable-nilpotent-criterion} shows that the composition of $s(\mathcal F)$ consecutive $f$-zero morphisms is $\mathcal F$-zero. At the end $\ri_X$ is the retract of some object $\mathcal F$ and we set $s(n) := s(\mathcal F)$.
\end{remarksenum}
\end{remarks}

\begin{definition} \label{def:index-of-descendable-map}
A morphism $f\colon Y \to X$ of analytic spaces over $(V,\mm)$ is \emph{descendable of index $n$} if $f$ is descendable and $n$ is the smallest integer such that $\mathcal K_f^n \to \ri_X$ is zero (see \cref{rslt:descendable-Un-criterion}). The same definition applies to a morphism of analytic rings.
\end{definition}

\begin{example}
Let $f\colon Y \to X$ be a qcqs steady morphism of analytic spaces over $(V,\mm)$. If $f$ is split then it is descendable of index $1$. Indeed, if $f$ is split then $\ri_X \to f_\natural \ri_Y$ is split and hence $\mathcal K_f \to \ri_X$ is zero.
\end{example}

Note that our definition of the index of descendable maps is analogous to \cite[Definition 11.18]{bhatt-scholze-witt}. One nice feature of descendable morphisms is that they enjoy many stability properties. We now discuss several of them. Later in \cref{sec:andesc.filtcolim} we investigate how descendable morphisms behave under filtered colimits of analytic rings.

\begin{lemma} \label{rslt:descendable-stable-under-base-change}
Let
\begin{center}\begin{tikzcd}
	Y' \arrow[r] \arrow[d,"f'"] & Y \arrow[d,"f"]\\
	X' \arrow[r] & X
\end{tikzcd}\end{center}
be a Cartesian square of analytic spaces over $(V,\mm)$. If $f$ is descendable (of index $n$) then $f'$ is descendable (of index $\le n$).
\end{lemma}
\begin{proof}
Assume that $f$ is descendable of index $n$, so in particular it is steady and qcqs. By base-change for the $\natural$-functors (see \cref{rslt:End-analytic-space-natural-base-change}) the map $\ri_X \to f_\natural \ri_Y$ is transformed under $g^\natural$ to $\ri_{X'} \to f'_\natural \ri_{Y'}$. It follows that the map $\mathcal K_f \to \ri_X$ is transformed under $g^\natural$ to $\mathcal K_{f'} \to \ri_{X'}$. Since $g^\natural$ is monoidal with respect to the composition monoidal structure (see \cref{rslt:End-analytic-space-upper-natural-properties}) it follows that the map $\mathcal K_f^n \to \ri_X$ is transformed under $g^\natural$ to $\mathcal K_{f'}^n \to \ri_{X'}$. The claim follows.
\end{proof}

\begin{lemma} \label{rslt:stability-of-descendable-maps-regarding-composition}
Let $g\colon Z \to Y$ and $f\colon Y \to X$ be two morphisms of analytic spaces over $(V,\mm)$.
\begin{lemenum}
	\item \label{rslt:descendability-stable-under-composition} If $f$ is descendable (of index $n$) and $g$ is descendable (of index $m$) then $f \comp g$ is descendable (of index $\le nm$).

	\item \label{rslt:descendability-of-composition-implies-desc-of-first-map} If $f \comp g$ is descendable (of index $n$) and $f$ is qcqs and steady, then $f$ is descendable (of index $\le n$).
\end{lemenum}
\end{lemma}
\begin{proof}
Let $h := f \comp g$. Then $h_\natural \ri_Z = f_\natural (g_\natural \ri_Z)$, so there is a cofiber sequence $f_\natural \mathcal K_g \to f_\natural \ri_X \to h_\natural \ri_Z$. Moreover, the counit $\ri_X \to h_\natural \ri_Z$ factors as $\ri_X \to f_\natural \ri_Y \to h_\natural \ri_Z$, where $\ri_X \to f_\natural \ri_Y$ is the counit and $f_\natural \ri_Y \to h_\natural \ri_Z$ is obtained via $f_\natural$ from the counit $\ri_Y \to g_\natural \ri_Z$. Completing this factorization to cofiber sequences, we obtain the following homotopy coherent diagram in $\mathcal E(X)$:
\begin{center}\begin{tikzcd}
	\mathcal K_f \arrow[r] \arrow[d] & \mathcal K_f \arrow[r] \arrow[d] & 0 \arrow[d]\\
	\mathcal K_h \arrow[r] \arrow[d,dashed] & \ri_X \arrow[r] \arrow[d] & h_\natural \ri_Z \arrow[d]\\
	f_\natural \mathcal K_g \arrow[r] & f_\natural \ri_Y \arrow[r] & h_\natural \ri_Z
\end{tikzcd}\end{center}
(All horizontal and vertical lines in this diagram are cofiber sequences and the dashed morphism is obtained through the diagram.) It follows that the morphism $\mathcal K_f \to \ri_X$ factors as $\mathcal K_f \to \mathcal K_h \to \ri_X$. Consequently, $\mathcal K_f^n \to \ri_X$ factors as $\mathcal K_f^n \to \mathcal K_h^n \to \ri_X$ for all $n \ge 0$, which implies (ii).

To prove (i), let the indices $n$ and $m$ be as in the claim. Consider the following diagram in $\mathcal E(X)$:
\begin{center}\begin{tikzcd}
	& \mathcal K_f \arrow[d]\\
	\mathcal K_h^m \arrow[r] \arrow[d] & \ri_X \arrow[d]\\
	(f_\natural \mathcal K_g)^m \arrow[r] & f_\natural \ri_Y
\end{tikzcd}\end{center}
Here the morphism $\mathcal K_h^m \to (f_\natural \mathcal K_g)^m$ is obtained from the morphism $\mathcal K_h \to f_\natural \mathcal K_g$ in the previous diagram and the morphism $(f_\natural \mathcal K_g)^m \to f_\natural \ri_Y$ is obtained as the composition of the obvious morphism $(f_\natural \mathcal K_g)^m \to (f_\natural \ri_Y)^m$ and the morphism $(f_\natural \ri_Y)^m \to f_\natural \ri_Y$ (given by adjunction).
Note that the natural map $\ri_X \to f_\natural \ri_Y$ factors as $\ri_X = \ri_X^m \to (f_\natural \ri_Y)^m \to f_\natural \ri_Y$ (e.g. use that $(f_\natural\ri_Y)^m = (f_\natural f^\natural)^m \ri_Y$ via \cref{rslt:computation-of-f-lower-nat-f-upper-nat-for-affine-target} and then use the properties of the monad $T = f_\natural f^\natural$). Using also the homotopy coherence of the lower left square of the first diagram, it follows that the second diagram is homotopy coherent.

One checks easily that the map $(f_\natural \mathcal K_g)^m \to f_\natural \ri_Y$ factors as $(f_\natural \mathcal K_g)^m \to f_\natural(\mathcal K_g^m) \to f_\natural \ri_Y$, where the second morphism is the one induced from the morphism $\mathcal K_g^m \to \ri_Y$ in $\mathcal E(Y)$ by applying $f_\natural$. By the descendability of $g$ this second morphism is zero, consequently $(f_\natural \mathcal K_g)^m \to f_\natural \ri_Y$ is zero. It follows that the composition $\mathcal K_h^m \to \ri_X \to f_\natural \ri_Y$ is zero and because the vertical column in the second diagram is a cofiber sequence we deduce that the map $\mathcal K_h^m \to \ri_X$ factors as $\mathcal K_h^m \to \mathcal K_f \to \ri_X$. This implies that the map $\mathcal K_h^{mn} \to \ri_X$ factors as $\mathcal K_h^{mn} \to \mathcal K_f^n \to \ri_X$ and is therefore zero by the descendability of $f$.
\end{proof}

\begin{lemma} \label{rslt:steady-covering-is-descendable}
Let $f\colon \bigdunion_{i=1}^n U_i \surjto X$ be a finite qcqs cover of an analytic space $X$ over $(V,\mm)$ by steady subspaces $U_i \subset X$. Then $f$ is descendable with index bounded by a constant $c(n)$ depending only on $n$.
\end{lemma}
\begin{proof}
For each nonempty subset $I \subset \{ 1, \dots, n \}$ let $U_I := \bigisect_{i\in I} U_i$ and let $f_I\colon U_I \to X$ denote the embedding. Then the natural functor $\mathcal E(X) \to \varprojlim_I \mathcal E(U_I)$, $\mathcal F \mapsto (f_I^\natural \mathcal F)_I$ is an equivalence (this is just saying that $\mathcal E(-)$ satisfies descent along steady covers, which is true by \cref{rslt:E-of-X-is-a-sheaf}). This equivalence has a right adjoint $\varprojlim_I \mathcal E(U_I) \to \mathcal E(X)$ given by $(\mathcal F_I)_I \mapsto \varprojlim_I f_{I\natural} \mathcal F_I$ (cf. \cref{rslt:E-of-X-stability-under-natural-functors}). In particular the natural morphism $\ri_X \to \varprojlim_I f_{I\natural} \ri_{U_I}$ is an isomorphism. Now each $f_{\{i\}\natural} \ri_{U_{\{i\}}}$ is a retract of $f_\natural \ri_{\bigdunion_i U_i}$ and any $f_{I\natural} \ri_{U_I}$ is obtained from that by composition, which proves that $f$ is descendable. Moreover, the structure of the limit $\ri_X = \varprojlim_I f_{I\natural} \ri_{U_I}$ does not depend on $X$ and the $U_i$'s, which implies that the index of descendability can be bounded by a constant only depending on $n$.
\end{proof}

\begin{corollary} \label{rslt:descendability-is-local}
Let $f\colon Y \to X$ be a map of analytic spaces over $(V,\mm)$.
\begin{corenum}
	\item \label{rslt:descendability-is-qcqs-local-on-target} Let $\bigdunion_{i=1}^n U_i \to X$ be a finite qcqs covering by steady subspaces $U_i \subset X$. If each restricted map $f^{-1}(U_i) \to U_i$ is descendable of index $\le m$, then $f$ is descendable of index $\le c(n) m$.

	\item Let $\bigdunion_{j\in J} V_j \to Y$ be a covering by steady subspaces $V_j \subset Y$. If the composed map $\bigdunion_{j\in J} V_j \to Y \to X$ is descendable (of index $n$) and $f$ is qcqs then $f$ is descendable (of index $\le n$).
\end{corenum}
\end{corollary}
\begin{proof}
Part (ii) is just a special case of \cref{rslt:descendability-of-composition-implies-desc-of-first-map} (using \cref{rslt:steadiness-can-be-checked-on-cover} to see that $f$ is steady). To prove (i), note that it follows from the assumptions that the map $\bigdunion_i f^{-1}(U_i) \to \bigdunion_i U_i$ is descendable of index $\le m$. Composing this map with the map $\bigdunion_i U_i \to X$ and using \cref{rslt:steady-covering-is-descendable} and \cref{rslt:descendability-stable-under-composition} we are reduced to (ii).
\end{proof}

\begin{remark}
Let $G$ be a geometry blueprint over $(V,\mm)$ (see \cref{def:geometry-blueprint}). Then all of the above definitions and results can be adapted to $G$-analytic spaces via the analytification functor from \cref{def:analytification}. For example we say that a qcqs steady map $f\colon Y \to X$ of $G$-analytic spaces is descendable (of index $n$) if this is true for $f^\an\colon Y^\an \to X^\an$.
\end{remark}

We close this subsection by discussing how descendability behaves under a change of almost setup $(V,\mm)$. In fact, this works as nicely as one can hope:

\begin{proposition}
Let $\varphi\colon (V,\mm) \to (V',\mm')$ be a morphism of almost setups and let $f\colon \mathcal A \to \mathcal B$ be a morphism of analytic rings over $(V,\mm)$. If quasicoherent sheaves descend along $f$, then they also descend along $\varphi^* f$.
\end{proposition}
\begin{proof}
Assume that quasicoherent sheaves descend along $f$ and denote $\mathcal A \to \mathcal B^\bullet$ the co-Čech nerve of $f$. We have natural pairs of adjoint functors
\begin{align*}
	L\colon \D(\mathcal A) \rightleftarrows \varprojlim_{n\in\Delta} \D(\mathcal B^n) \noloc R, \qquad L'\colon \D(\varphi^*\mathcal A) \rightleftarrows \varprojlim_{n\in\Delta} \D(\varphi^*\mathcal B^n) \noloc R'.
\end{align*}
Here $L(M) = (M \tensor_{\mathcal A} \mathcal B^\bullet)$ and $R(M^\bullet) = \varprojlim_{n\in\Delta} M = \Tot(M^\bullet)$ (and similarly for $L'$ and $R'$). By assumption both $L$ and $R$ are equivalences, and we need to show that the same is true for $L'$ and $R'$. This can be checked separately in the cases that $\varphi$ is strict and that $\varphi$ is localizing.

First assume that $\varphi$ is strict. Then $\varphi$ corresponds to the morphism $V^a \to V'^a$ of analytic rings over $(V,\mm)$ and $\varphi^* \mathcal A$ is the induced analytic ring structure on $\underline{\mathcal A} \tensor_V V'$ from $\mathcal A$ (and similarly for $\mathcal B^n$ in place of $\mathcal A)$. In particular, $\varphi^*\colon \D(\mathcal A) \to \D(\varphi^*\mathcal A)$ has a right adjoint $\varphi_*$ (similarly for $\mathcal B^n$), which induces a functor
\begin{align*}
	\varphi_*\colon \varprojlim_{n\in\Delta} \D(\varphi^*\mathcal B^n) \to \varprojlim_{n\in\Delta} \D(\mathcal B^n), \qquad (M^\bullet) \mapsto (\varphi_* M^\bullet).
\end{align*}
(For this to work one needs to check that $\varphi_*$ commutes with pullbacks in the diagram $\D(\varphi^* \mathcal B^\bullet)$ which is obvious by the above description of $\mathcal B^\bullet$ in terms of the induced analytic ring structure.) Note that $\varphi_* L' = L \varphi_*$ and $\varphi_* R' = R \varphi_*$. In particular $\varphi_* L'R' = LR \varphi_* = \varphi_*$ and $\varphi_* R'L' = \varphi_*$. Since $\varphi_*$ is conservative this implies $L'R' = \id$ and $R'L' = \id$, hence $L'$ and $R'$ are equivalences.

Now assume that $\varphi$ is localizing. Then by \cref{rslt:left-adjoint-of-almost-localization-over-analytic-ring} the functor $\varphi^*\colon \D(\mathcal A) \to \D(\varphi^* \mathcal A)$ admits a fully faithful left adjoint $\varphi_!$ (and similarly for $\mathcal B^n$ in place of $\mathcal A$). By passing to left adjoints in \cref{rslt:almost-localization-is-steady} we deduce that $\varphi_!$ commutes with pullbacks in the diagram $\D(\varphi^* \mathcal B^\bullet)$ and hence defines a functor
\begin{align*}
	\varphi_!\colon \varprojlim_{n\in\Delta} \D(\varphi^*\mathcal B^n) \to \varprojlim_{n\in\Delta} \D(\mathcal B^n), \qquad (M^\bullet) \mapsto (\varphi_! M^\bullet).
\end{align*}
Then $\varphi_! L' = L \varphi_!$ and $R' = \varphi^* R \varphi_!$. Therefore $R'L' = \varphi^* R \varphi_! L' = \varphi^* RL \varphi_! = \varphi^* \varphi_! = \id$. By functoriality of pullbacks (see \cref{rslt:functoriality-of-modules-over-analytic-rings-over-AlmSetup}) we also have $L' \varphi^* = \varphi^* L$ and hence $L' R' = L' \varphi^* R \varphi_! = \varphi^* LR \varphi_! = \varphi^* \varphi_! = \id$. This shows that $L'$ and $R'$ are equivalences, as desired.
\end{proof}

\begin{proposition} \label{rslt:descendability-stable-under-almost-localization}
Let $\varphi\colon (V,\mm) \to (V',\mm')$ be a morphism of almost setups and let $f\colon Y \to X$ be a map of analytic spaces over $(V,\mm)$. If $f$ is descendable (of index $n$) then $\varphi^* f$ is descendable (of index $\le n$).
\end{proposition}
\begin{proof}
Assume that $f$ is descendable of some index $n$. By \cref{rslt:functoriality-of-analytic-spaces} $\varphi^* f$ is qcqs and steady. By definition of descendability the morphism $\mathcal K_f^n \to \ri_X$ in $\mathcal E(X)$ is zero. By \cref{rslt:functoriality-of-steady-endofunctors-over-AlmSetup} the morphism $\mathcal K_{\varphi^* f}^n = \varphi^\natural \mathcal K_f^n \to \varphi^\natural \ri_X = \ri_{\varphi^* X}$ in $\mathcal E(\varphi^* X)$ is zero, hence $\varphi^* f$ is descendable of index $\le n$.
\end{proof}

\begin{example}
Let $f\colon \mathcal A \to \mathcal B$ be a map of analytic rings over $V$. If quasicoherent sheaves descend along $f$ then they also descend along $f^a\colon \mathcal A^a \to \mathcal B^a$. Moreover, if $f$ is descendable (of index $n$) then $f^a$ is descendable (of index $\le n$).
\end{example}

\subsection{Descendability of Filtered Colimits} \label{sec:andesc.filtcolim}

As before we fix an almost setup $(V,\mm)$. In our applications we are often interested in the descendability of filtered colimits of maps of analytic rings over $(V,\mm)$ each of which is of a nice simple form; most prominently, pro-étale maps are locally given as a filtered colimit of étale maps. Unfortunately descendability is in general not preserved under filtered colimits (cf. \cite[Example 11.21]{bhatt-scholze-witt}). However, quasicoherent sheaves still descend along a filtered colimit of descendable maps of \emph{bounded index}. Proving this statement is the main goal of the present subsection (see \cref{rslt:weakly-descendable-maps-properties}). The proof consists of two steps: The first step is to show that countable filtered colimits of descendable maps of bounded index are still descendable (similar to \cite[Lemma 11.22]{bhatt-scholze-witt}) and the second step is to show that the descent of quasicoherent sheaves is stable under $\omega_1$-filtered colimits of analytic rings.

We start with the first step, i.e. that descendability is preserved by countable filtered colimits. This will be the content of \cref{rslt:descendable-stable-under-countable-colimit} below, for which we first need to do a straightforward computation in the $\infty$-category of spectra. Recall that this $\infty$-category comes equipped with a natural $t$-structure whose heart is the category of abelian groups. Moreover, the $\infty$-category of spectra is symmetric monoidal via the smash product $\tensor$.

\begin{lemma}
Let $I$ be a countable filtered category.
\begin{lemenum}
	\item \label{rslt:countable-limit-of-spectra-SES} Let $(X_i)_{i\in I}$ be a diagram of spectra. Then for every $n \in \Z$ there is a natural short exact sequence of abelian groups
	\begin{align*}
		0 \to R^1 \varprojlim_i \pi_{n+1} X_i \to \pi_n \varprojlim_i X_i \to \varprojlim_i \pi_n X_i \to 0.
	\end{align*}
	Here the limit on the right denotes the non-derived limit of abelian groups.

	\item \label{rslt:countable-limit-of-spectra-pair-to-R2-lim} Let $(X_i)_{i\in I}$ and $(X'_i)_{i\in I}$ be diagrams of spectra. Then for all $n, k \in \Z$ the natural map
	\begin{align*}
		(R^1 \varprojlim_i \pi_{n+1} X_i) \tensor (R^1 \varprojlim_i \pi_{k+1} X'_i) \to (\pi_n \varprojlim_i X_i) \tensor (\pi_k \varprojlim_i X'_i) \to \pi_{n+k} \varprojlim_i (X_i \tensor X_i')
	\end{align*}
	is zero. Here the second map comes from the natural map $\varprojlim_i X_i \tensor \varprojlim_i X'_i \to \varprojlim_i X_i \tensor X'_i$.
\end{lemenum}
\end{lemma}
\begin{proof}
We first prove (i), so let $(X_i)_i$ and $n \in \Z$ be given. There is a natural cofiber sequence $(\tau_{\ge n+1} X_i)_i \to (X_i)_i \to (\tau_{\le n} X_i)_i$ of diagrams of spectra. By applying $\varprojlim$ and taking the associated homology long exact sequence we obtain an exact sequence of abelian groups
\begin{align*}
	&\pi_{n-1} \varprojlim_i \tau_{\le n} X_i \to \pi_n \varprojlim_i \tau_{\ge n+1} X_i \to \pi_n \varprojlim_i X_i \to \pi_n \varprojlim_i \tau_{\le n} X_i \to \pi_n \varprojlim_i \tau_{\ge n+1} X_i \to\\&\qquad\to \pi_{n-1} \varprojlim_i \tau_{\ge n+1} X_i
\end{align*}
Now $\varprojlim_i \tau_{\le n} X_i$ is concentrated in homological degrees $\le n$, so the first term of the above sequence is $0$. Moreover, since $I$ is countable, $\varprojlim_i \tau_{\ge n+1} X_i$ is concentrated in homological degrees $\ge n$ (this boils down to the fact that $R^2\varprojlim_i$ vanishes on systems of abelian groups), hence the last term of the above sequence is $0$. This proves (i).

For (ii), first note that the natural map $\varprojlim_i X_i \tensor \varprojlim_i X'_i \to \varprojlim_i X_i \tensor X'_i$ is obtained by noticing that $\varprojlim_i$ is right adjoint to the functor $X \mapsto (X)_i$ from spectra to $I$-indexed diagrams of spectra (indeed, $\varprojlim_i$ is just a right Kan extension). Since the functor $X \mapsto (X)_i$ is clearly symmetric monoidal, its right adjoint is lax symmetric monoidal. Now apply the functors $\varprojlim$ and $\tensor$ to the cofiber sequence
\begin{align*}
	(\tau_{\ge n+1} X_i, \tau_{\ge k+1} X'_i)_i \to (X_i, X'_i)_i \to (\tau_{\le n} X_i, \tau_{\le k} X'_i)_i
\end{align*}
of $I \cprod I$-indexed diagrams of spectra. Looking at the associated homology long exact sequences we see that the described map in (ii) factors over $\pi_{n+k} \varprojlim_i (\tau_{\ge n+1} X_i \tensor \tau_{\ge k+1} X'_i)$. But $\tau_{\ge n+1} X_i \tensor \tau_{\ge k+1} X'_i$ is concentrated in homological degrees $\ge n + k + 2$, hence the $\varprojlim$ of that lies in homological degrees $\ge n + k + 1$ (again because $I$ is countable and therefore $R^2\varprojlim_i = 0$). Thus the map in (ii) factors over $0$, as desired.
\end{proof}

\begin{proposition} \label{rslt:descendable-stable-under-countable-colimit}
Let $I$ be a countable filtered category and let $(f_i\colon \mathcal A_i \to \mathcal B_i)_{i \in I}$ be a diagram of maps of analytic rings over $(V,\mm)$. Assume that there is some $n \ge 0$ such that each $f_i$ is descendable of index $\le n$. Then the colimit
\begin{align*}
	f := \varinjlim_i f_i\colon \varinjlim_i \mathcal A_i \to \varinjlim_i \mathcal B_i
\end{align*}
is descendable of index $\le 2n$.
\end{proposition}
\begin{proof}
Let $\mathcal A := \varinjlim_i \mathcal A_i$ and let $\mathcal B'_i := \mathcal B_i \tensor_{\mathcal A_i} \mathcal A$ with morphism $f'_i \colon \mathcal A \to \mathcal B'_i$ for every $i \in I$. Then each $f'_i$ is descendable of index $\le n$ by \cref{rslt:descendable-stable-under-base-change} and $f = \varinjlim_i f'_i$. We can thus replace $\mathcal A_i$ by $\mathcal A$ and $\mathcal B_i$ by $\mathcal B'_i$ to assume that all $f_i$ have the same source $\mathcal A$.

The rest of the proof is analogous to \cite[Lemma 11.22]{bhatt-scholze-witt}: Note first that since $\mathcal E(\mathcal A)$ is stable, it is naturally enriched over the $\infty$-category of spectra (cf. \cite[Example 7.4.14]{gepner-haugseng-enriched}), so in the following we can view all $\Hom$'s as spectra. We have $f_\natural \mathcal B = \varinjlim_i f_{i\natural} \id_{\mathcal B_i}$ in $\mathcal E(\mathcal A)$ and consequently $\mathcal K_f = \varinjlim_i \mathcal K_{f_i}$. By assumption the morphisms $\mathcal K_{f_i}^n \to \id_{\mathcal A}$ are zero for all $i$. We will show that $\mathcal K_f^{2n} \to \id_{\mathcal A}$ is zero. By \cref{rslt:countable-limit-of-spectra-SES} we have a short exact sequence of abelian groups
\begin{align*}
	0 \to R^1 \varprojlim_i \pi_1\Hom(\mathcal K_{f_i}^n, \id_{\mathcal A}) \to \pi_0 \Hom(\mathcal K_f^n, \id_{\mathcal A}) \to \varprojlim_i \pi_0 \Hom(\mathcal K_{f_i}^n, \id_{\mathcal A}) \to 0.
\end{align*}
It follows that the map $\mathcal K_f^n \to \id_{\mathcal A}$ lies in the subgroup
\begin{align*}
	R^1 \varprojlim_i \pi_1\Hom(\mathcal K_{f_i}^n, \id_{\mathcal A}) \subset \pi_0 \Hom(\mathcal K_f^n, \id_{\mathcal A}).
\end{align*}
In order to conclude, it is enough to observe that map on $\pi_0$ induced by the obvious map
\begin{align*}
	\Hom(\mathcal K_f^n, \id_{\mathcal A}) \cprod \Hom(\mathcal K_f^n, \id_{\mathcal A}) \to \Hom(\mathcal K_f^{2n}, \id_{\mathcal A})
\end{align*}
kills the subgroups $R^1 \varprojlim_i \pi_1\Hom(\mathcal K_{f_i}^n, \id_{\mathcal A})$ (note that the morphism $\mathcal K_f^{2n} \to \id_{\mathcal A}$ is the image of the morphism $\mathcal K_f^n \to \id_{\mathcal A}$ under the above map). This follows from the fact that the above map pairs the $R^1\varprojlim$'s to an $R^2\varprojlim = 0$, see \cref{rslt:countable-limit-of-spectra-pair-to-R2-lim}.
\end{proof}

We now come to the second step in our study of descent along filtered colimits. While descendability is not retained under filtered colimits, we now show that descent of quasicoherent sheaves is retained under $\omega_1$-filtered colimits. This is mostly a formal consequence of the fact that $\omega_1$-filtered colimits commute with countable limits in $\infcatinf$ (see \cref{rslt:filtered-colim-preserve-finite-lim-in-infcatinf}) and that countable limits of $\omega_1$-compactly generated $\infty$-categories are still $\omega_1$-compactly generated (see \cref{rslt:tau-small-limits-in-PrL-tau}). We could not find these two results in the literature, so we included proofs in the appendix. With these results at hand, we can argues as follows (see \cref{rslt:descent-of-qcoh-is-stable-under-w1-filtered-colim} for the final result).

\begin{lemma} \label{rslt:description-of-filtered-colim-of-infcat}
Let $\mathcal C = \varinjlim_{i\in I} \mathcal C_i$ be a filtered colimit in $\infcatinf$. Then up to equivalence, $\mathcal C$ can be described as follows:
\begin{lemenum}
	\item The objects of $\mathcal C$ are given by $\Obj(\mathcal C) = \varinjlim_i \Obj(\mathcal C_i)$.
	\item Given two objects $X, Y \in \mathcal C$, choose $i \in I$ so that $X, Y$ lie in the image of $\mathcal C_i \to \mathcal C$. For every morphism $i \to j$ in $I$ let $X_j$ and $Y_j$ be the images of $X$ and $Y$ along $\mathcal C_i \to \mathcal C_j$. Then
	\begin{align*}
		\Hom_{\mathcal C}(X, Y) = \varinjlim_{j\in I_{i/}} \Hom_{\mathcal C_j}(X_j, Y_j).
	\end{align*}
\end{lemenum}
\end{lemma}
\begin{proof}
A close inspection of the construction in the proof of \cref{rslt:infcatinf-is-compactly-generated} reveals that the functor $\infcatinf \equiv S^{-1} \mathcal P(\Delta) \injto \mathcal P(\Delta)$ is given by sending $\mathcal C \in \infcatinf$ to the functor $[n] \mapsto \Hom_{\infcatinf}([n], \mathcal C)$. As we showed that the morphisms in $S$ are between compact objects, the proof of \cref{rslt:localization-compactly-generated-criterion} shows that $\infcatinf \subset \mathcal P(\Delta)$ is stable under filtered colimits. As colimits in $\mathcal P(\Delta)$ are computed pointwise, the claim follows.
\end{proof}

\begin{lemma} \label{rslt:filtered-colim-in-AnRing-implies-colim-of-compact}
Let $\kappa$ be a regular cardinal such that $(V,\mm)$ is $\kappa$-compact. Let $\mathcal A = \varinjlim_{i\in I} \mathcal A_i$ be a $\kappa$-filtered colimit of analytic rings over $(V,\mm)$. Then there is a natural equivalence
\begin{align*}
	\D(\mathcal A)^\kappa = \varinjlim_i \D(\mathcal A_i)^\kappa.
\end{align*}
\end{lemma}
\begin{proof}
Let us abbreviate $\mathcal C_i = \D(\mathcal A_i)$ and $\mathcal C = \D(\mathcal A)$. Let further $\mathcal C' := \varinjlim_i \mathcal C_i^\kappa$. We will implicitly make use of the description from \cref{rslt:description-of-filtered-colim-of-infcat} for $\mathcal C'$.

First note that since pullback functors preserve $\kappa$-compact objects, there is a natural functor $f\colon \mathcal C' \to \mathcal C^\kappa$. We first show that this functor is fully faithful, so let $X', Y' \in \mathcal C'$ be given. Choose $i \in I$ large enough so that $X'$ and $Y'$ come from objects $X_i, Y_i \in \mathcal C_i^\kappa$. For each $j \in I_{i/}$ let $X_j, Y_j \in \mathcal C_j^\kappa$ be the pullback of $X_i$ and $Y_i$ to $\mathcal C_j$. Finally, let $X$ and $Y$ be the images of $X'$ and $Y'$ in $\mathcal C^\kappa$. By the adjunction of pullback and pushforward (the latter being a forgetful functor) and the $\kappa$-compactness of $X_i$ we get
\begin{align*}
	&\Hom_{\mathcal C}(X, Y) = \Hom_{\mathcal C_i}(X_i, Y) = \Hom_{\mathcal C_i}(X_i, \varinjlim_{j \in I_{i/}} Y_j) = \varinjlim_{j \in I_{i/}} \Hom_{\mathcal C_i}(X_i, Y_j) =\\&\qquad= \varinjlim_{j \in I_{i/}} \Hom_{\mathcal C_j}(X_j, Y_j) = \Hom_{\mathcal C'}(X, Y).
\end{align*}
This proves full faithfulness of the functor $f\colon \mathcal C' \injto \mathcal C^\kappa$. To show essential surjectivity, first note that all $\mathcal A[S] \in \mathcal C^\kappa$ lie in the image of $f$, where $S$ runs through extremally disconnected sets (because they come via pullback from $\mathcal A_i[S] \in \mathcal C_i^\kappa$ for any $i \in I$). By the full faithfulness of $f$ we know that the essential image of $f$ is stable under retracts and $\kappa$-small colimits (for the latter use \cite[Proposition 5.5.7.11]{lurie-higher-topos-theory}). But by \cref{rslt:compact-generators-implies-compactly-generated} every $X \in \mathcal C^\kappa$ is a retract of an iterated $\kappa$-small colimit of $\mathcal A[S]$'s (more precisely, apply \cref{rslt:compact-generators-implies-compactly-generated} to the presentable $\infty$-categories $\D(\mathcal A)_{\kappa'}$ for increasing strong limit cardinals $\kappa'$). This finishes the proof.
\end{proof}

\begin{proposition} \label{rslt:descent-of-qcoh-is-stable-under-w1-filtered-colim}
Let $\kappa$ be an uncountable regular cardinal such that $(V,\mm)$ is $\kappa$-compact and let $(f_i\colon \mathcal A \to \mathcal B_i)_{i\in I}$ be a $\kappa$-filtered family of morphisms of analytic rings over $(V,\mm)$. If quasicoherent sheaves descend along each $f_i$, then they also descend along the colimit $f\colon \mathcal A \to \varinjlim_i \mathcal B_i$.
\end{proposition}
\begin{proof}
Let $\mathcal B = \varinjlim_i \mathcal B_i$. For each $i \in I$ and $n \ge 0$ let $\mathcal C_i^n := \D(\mathcal B_i^{\tensor n+1})$, where $\mathcal B_i^{\tensor k}$ is the $k$-fold tensor product of $\mathcal B_i$ over $\mathcal A$. Similarly define $\mathcal C^n := \D(\mathcal B^{\tensor n+1})$ and let also $\mathcal C := \D(\mathcal A)$. By assumption the functor
\begin{align*}
	\mathcal C \isoto \varprojlim_{n\in\Delta} \mathcal C_i^n
\end{align*}
is an equivalence for all $i \in I$. We have to show that $\mathcal C \to \varprojlim_{n\in\Delta} \mathcal C^n$ is an equivalence as well.

Let us denote $\mathcal C' := \varprojlim_{n\in\Delta} \mathcal C^n$. Working with $\D(-)_{\kappa'}$ instead of $\D(-)$ in all of the above definitions (where $\kappa'$ is a suitably large strong limit cardinal), all $\mathcal C_i^n$ and $\mathcal C^n$ are $\kappa$-compactly generated, so that by \cref{rslt:tau-small-limits-in-PrL-tau} all $\infty$-categories introduced in this proof are $\kappa$-compactly generated and the functor $\mathcal C \to \mathcal C'$ preserves $\kappa$-compact objects. Now working with $\D(-)$ instead of $\D(-)_{\kappa'}$ again, it is enough to show that the induced functor $\mathcal C^\kappa \to \mathcal C'^\kappa$ is an equivalence.

From \cref{rslt:tau-small-limits-in-PrL-tau} we get $\mathcal C'^\kappa = \varprojlim_n (\mathcal C^n)^\kappa$ and for fixed $i$, $\mathcal C^\kappa = \varprojlim_n (\mathcal C_i^n)^\kappa$. On the other hand, for fixed $n$ we have $(\mathcal C^n)^\kappa = \varinjlim_i (\mathcal C_i^n)^\kappa$ by \cref{rslt:filtered-colim-in-AnRing-implies-colim-of-compact}. By \cref{rslt:filtered-colim-preserve-finite-lim-in-infcatinf} countable limits commute with $\kappa$-filtered colimits in $\infcatinf$, so we get
\begin{align*}
	\mathcal C'^\kappa = \varprojlim_n (\mathcal C^n)^\kappa = \varprojlim_n \varinjlim_i (\mathcal C_i^n)^\kappa = \varinjlim_i \varprojlim_n (\mathcal C_i^n)^\kappa = \varinjlim_i \mathcal C^\kappa = \mathcal C^\kappa,
\end{align*}
as desired.
\end{proof}

We have prove the promised stability of descent under $\kappa$-filtered colimits. This suggests the following definition, which generalizes the notion of descendability while keeping many of its stability properties:

\begin{definition} \label{def:weakly-descendable-morphism-of-analytic-rings}
A morphism $\mathcal A \to \mathcal B$ of analytic rings over $(V,\mm)$ is called \emph{weakly descendable} if it is an iterated $\kappa$-filtered colimit of descendable maps, where $\kappa$ is any uncountable regular cardinal such that $(V,\mm)$ is $\kappa$-compact.
\end{definition}

\begin{remark}
One could introduce an index of descendability on weakly descendable morphisms by adding a bound on the index of descendability in the $\kappa$-filtered family of morphisms. However, the main motivation to introduce the index of descendability was to show descendability of countable colimits (cf. \cref{rslt:descendable-stable-under-countable-colimit}) and this motivation does not apply to weakly descendable morphisms.
\end{remark}

We now arrive at the main result of the present subsection, which summarizes the above work by showing that weakly descendable morphisms form a rather nice class of morphisms for studying descent:

\begin{theorem} \label{rslt:weakly-descendable-maps-properties}
\begin{thmenum}
	\item Quasi-coherent sheaves descend along weakly descendable morphisms.
	\item Every descendable morphism is weakly descendable and weakly descendable morphisms are stable under base-change and $\kappa$-filtered colimits, where $\kappa$ is any uncountable regular cardinal such that $(V,\mm)$ is $\kappa$-compact.
	\item \label{rslt:filtered-colim-of-bounded-desc-is-weakly-desc} Suppose that $(V,\mm)$ is $\omega_1$-compact (e.g. $\mm$ is countably generated). Then a filtered colimit of descendable morphisms of bounded index is weakly descendable.
\end{thmenum}
\end{theorem}
\begin{proof}
Part (i) follows from \cref{rslt:descent-of-qcoh-is-stable-under-w1-filtered-colim,rslt:descendable-implies-limit-of-categories}. The only non-trivial part of (ii) is the stability under base-change, which follows from \cref{rslt:descendable-stable-under-base-change}.

It remains to prove part (iii). Note that if $\mm$ is countably generated then $(V,\mm)$ is indeed $\omega_1$-compact by \cref{rslt:m-countably-generated-implies-omega-1-compact}. Let $(f_i\colon \mathcal A \to \mathcal B_i)_{i\in I}$ be a filtered family of morphisms of analytic rings over $(V,\mm)$ and assume that there is some integer $d \ge 0$ such that all $f_i$ are descendable of index $\le d$. We want to show that the colimit $f\colon \mathcal A \to \varinjlim_i \mathcal B_i$ is weakly descendable. By \cite[Proposition 5.3.1.16]{lurie-higher-topos-theory} we can assume that $I$ is a filtered partially ordered set. Let $J$ be the filtered partially ordered set of all functors $j\colon \N \to I$, where we define $j_1 \le j_2$ iff $j_1(k) \le j_2(k)$ for $k \gg 0$. Then $J$ is in fact $\omega_1$-filtered, i.e. for every countable subset $\{ j_1, j_2, \dots \} \subset J$ there is some $j \in J$ with $j_k \le j$ for all $k$; just choose $j(k)$ inductively to be greater than $j_1(k), \dots, j_k(k)$ and $j(k-1)$. For every $j \in J$ let $\mathcal B_j := \varinjlim_{k\in\N} \mathcal B_{j(k)}$ with associated morphism $f_j\colon \mathcal A \to \mathcal B_j$. We then obtain the $\omega_1$-filtered system $(f_j)_{j\in J}$ of morphisms in $\AnRing_{(V,\mm)}$. By \cref{rslt:descendable-stable-under-countable-colimit} each $f_j$ is descendable. This implies the claim.
\end{proof}

\subsection{Descendability over Flat Covers} \label{sec:andesc.fsdesc}

As before we fix an almost setup $(V,\mm)$. In the present subsection we study a weaker version of descendability which still implies descent of quasicoherent sheaves. Roughly, instead of asking for a map $f\colon \mathcal A \to \mathcal B$ of analytic rings to be descendable, we only require that $f$ be descendable (of bounded index) over a ``flat cover'' of $\AnSpec \mathcal A$. In our application this ``flat cover'' will be given by the connected components of an affine scheme $\Spec A$. Note that this cover is not an fpqc cover in general -- we make no assumption on quasicompactness.

The main ingredient for the weaker version of descendability is a thorough study of convergence speeds. In the proof of \cref{rslt:descendable-implies-limit-of-categories} the crucial input from the descendability condition on a map $f\colon Y \to \AnSpec \mathcal A$ of analytic spaces over $(V,\mm)$ was the fact that if a cosimplicial object $M^\bullet$ in $\D(\mathcal A)$ becomes split in $\D(Y)$ after applying $f^*$ then the associated $\Tot$-tower $(\Tot_n(M^\bullet))_n$ is pro-constant. In the following we will refine this statement by showing that the index of descendability of $f$ gives a bound on the speed of convergence of $(\Tot_n(M^\bullet))_n$.

The following definitions and results are similar to \cite[\S3]{akhil-thick-subcategory-theorem} but we keep close attention to the speed of convergence throughout.

\begin{definition}
Let $\mathcal C$ be a stable $\infty$-category.
\begin{defenum}
	\item Let $\Tow(\mathcal C) := \Fun((\Z_{\ge0})^\opp, \mathcal C)$ denote the category of \emph{towers} in $\mathcal C$.
	\item A tower $(X_n)_n \in \Tow(\mathcal C)$ is called \emph{nilpotent of index $\le h$} if for all $n \ge h$ the map $X_{n+h} \to X_n$ is zero.
	\item A tower $(X_n)_n \in \Tow(\mathcal C)$ is called \emph{strongly constant of index $h$} if $X := \varprojlim_n X_n$ exists and the cofiber of $X \to (X_n)_n$ is nilpotent of index $h$ (where we identify $X$ with the constant tower with value $X$).
	\item A cosimplicial object $M^\bullet$ in $\mathcal C$ is called \emph{strongly constant of index $h$} if the associated tower $(\Tot_n(M^\bullet))_n$ is strongly constant of index $h$.
\end{defenum}
\end{definition}

\begin{remarks}
\begin{remarksenum}
	\item A strongly constant tower is automatically pro-constant, but the converse is far from true (cf. \cite[Remark 3.6]{akhil-thick-subcategory-theorem}). Evidently a tower $(X_n)_n$ is nilpotent (of index $h$) if and only if it is strongly constant (of index $h$) and $\varprojlim_n X_n \isom 0$.
	\item \label{rmk:refined-nilpotent-tower-index} One could refine the definition of nilpotent towers by saying that $(X_n)_n$ is nilpotent of index $(h, n_0)$ if for all $n \ge n_0$ the map $X_{n+h} \to X_n$ is zero. However, one is usually only concerned with upper bounds on both $h$ and $n_0$ so that we can savely assume $n_0 = h$ to simplify the notation.
\end{remarksenum}
\end{remarks}

\begin{lemma} \label{rslt:f*M-strongly-constant-implies-M-strongly-constant-if-f-descendable}
Let $f\colon Y \to X$ be a morphism of analytic spaces over $(V,\mm)$ which is descendable of index $m$ and let $(\mathcal M_n)_n \in \Tow(\D(X))$. If $(f^*M_n)_n$ is strongly constant of index $h$ then $(\mathcal M_n)_n$ is strongly constant of index $s(m) h$, where $s(m)$ is the constant from \cref{rslt:f-zero-only-depends-on-index}.
\end{lemma}
\begin{proof}
Suppose that $(f^* \mathcal M_n)_n \in \Tow(\D(Y))$ is strongly constant of index $h$. Then the same is true for $(f_* f^* \mathcal M_n)_n = ((f_\natural \ri_Y)(\mathcal M_n))_n \in \Tow(\D(X))$ (where we use the evaluation functor $\mathcal E(X) \cprod \D(X) \to \D(X)$ from \cref{def:evaluation-map-on-E-of-X}). By definition of descendability, $\ri_X \in \mathcal E(X)$ can be obtained as a finite sequence of compositions, finite (co)limits and retracts from $f_\natural \ri_Y$. One checks easily that strongly constant towers are stable under these operations as well, proving that $(\mathcal M_n)_n$ is indeed strongly constant. To get the bound on the index, we can argue as follows: We can reduce to the case that $((f_\natural \ri_Y)(\mathcal M_n))_n$ is nilpotent of index $h$, i.e. for every $n \ge h$ the map $(f_\natural \ri_Y)(\mathcal M_{n+h}) \to (f_\natural \ri_Y)(\mathcal M_n)$ is zero. For every $\mathcal G \in \mathcal E(X)$ let $\mathcal I_{\mathcal G}$ be the collection of morphisms $\mathcal M \to \mathcal N$ in $\D(X)$ such that $\mathcal G(\mathcal M) \to \mathcal G(\mathcal N)$ is zero. Arguing as in the proof of \cref{rslt:descendable-nilpotent-criterion} we see that $\mathcal I_{f_\natural \ri_Y}^{s(m)} = 0$. Since each map $\mathcal M_{n+h} \to \mathcal M_n$ lies in $\mathcal I_{f_\natural \ri_Y}$ we deduce that $M_{n+s(m)h} \to M_n$ is zero and hence that $(\mathcal M_n)_n$ is nilpotent of index $s(m) h$.
\end{proof}

In the following we will provide an analogue of \cref{rslt:f*M-strongly-constant-implies-M-strongly-constant-if-f-descendable} in terms of the convergence of homotopy groups and will relate this to the associated spectral sequence.

\begin{definition} \label{def:spectral-sequence-of-tower}
Let $\mathcal C$ be a stable $\infty$-category with $t$-structure and let $(X_n)_n \in \Tow(\mathcal C)$. After passing to the dual category $\mathcal C^\opp$, \cite[Lemma 1.2.2.4 and Construction 1.2.2.6]{lurie-higher-algebra} associates to $(X_n)_n$ a homological spectral sequence $E^r_{p,q}$ in $(\mathcal C^\heartsuit)^\opp$, which is the same as a cohomological spectral sequence $\tilde E_r^{p,q}$ in $\mathcal C^\heartsuit$. As we want to keep homological notation in $\mathcal C$ it is convenient to flip the $q$-coordinate and define $E_r^{p,q} := \tilde E_r^{p,-q}$ to be the \emph{spectral sequence associated to $(X_n)_n$}. Be aware that the differentials $d_r\colon E_r^{*,*} \to E_r^{*,*}$ have bidegree $(r, r-1)$.
\end{definition}

The spectral sequence associated to a tower has the following alternative description in terms of exact couples:

\begin{lemma} \label{rslt:spectral-sequence-comes-from-exact-couple}
Let $\mathcal C$ be a stable $\infty$-category with $t$-structure and let $(X_n)_n \in \Tow(\mathcal C)$. For $p, q \in \Z$ let $D_1^{p,q} := \pi_{q-p}(X_p)$ and let $E_1^{p,q} := \pi_{q-p}(\fib(X_p \to X_{p-1}))$. There is an exact couple
\begin{center}\begin{tikzcd}
	D_1^{*,*} \arrow[rr,"g"] && D_1^{*,*} \arrow[ld,"h"]\\
	& E_1^{*,*} \arrow[ul,"f"]
\end{tikzcd}\end{center}
Here $f$ is induced by $\fib(X_p \to X_{p-1}) \to X_p$ and has bidegree $(0, 0)$, $g$ is induced by $X_p \to X_{p-1}$ and has bidegree $(-1, -1)$, and $h$ is induced by $X_p \to \fib(X_{p+1} \to X_p)[1]$ and has bidegree $(1, 0)$.

Moreover, the spectral sequence associated to the above exact couple is naturally isomorphic to the spectral sequence associated to $(X_n)_n$ as defined in \cref{def:spectral-sequence-of-tower}.
\end{lemma}
\begin{proof}
It follows immediately from the long exact sequence associated to $\fib(X_p \to X_{p-1}) \to X_p \to X_{p-1}$ that the given maps $f$, $g$ and $h$ define an exact couple in $\mathcal C^\heartsuit$. Let $C_r\colon D_r^{*,*} \to E_r^{*,*} \to D_r^{*,*}$ denote the $r$-th derived exact couple, so that $E_*^{*,*}$ is the spectral sequence associated to the exact couple. Let us temporarily denote by $\tilde E_*^{*,*}$ the spectral sequence associated to $(X_n)_n$ as in \cref{def:spectral-sequence-of-tower}. We want to show that $E_*^{*,*} \isom \tilde E_*^{*,*}$.

By definition (see \cite[Construction 1.2.2.6]{lurie-higher-algebra}) we have
\begin{align*}
	\tilde E_r^{p,q} = \img(\pi_{q-p} X(p+r-1, p-1) \to \pi_{q-p} X(p, p-r)),
\end{align*}
where $X(i, j) \in \mathcal C$ is functorial in $i, j \in (\Z \union \{ -\infty \})^\opp$ such that $X(i,-\infty) = X_i$ for all $i \ge 0$, $X(i, -\infty) = 0$ for $i < 0$ and $X(i, j) \to X(i, k) \to X(j, k)$ is a cofiber sequence for all $i \le j \le k$. In particular we have $X(i, j) \isom \fib(X_i \to X_j)$.

For all $p$, $q$ and $r$ we have a commutative diagram
\begin{center}\begin{tikzcd}
	\pi_{q-p+1} X_{p-1} \arrow[r] \arrow[d] & \pi_{q-p} X(p+r-1, p-1) \arrow[r] \arrow[d] & \pi_{q-p} X_{p+r-1} \arrow[d]\\
	\pi_{q-p+1} X_{p-r} \arrow[r] & \pi_{q-p} X(p, p-r) \arrow[r] & \pi_{q-p} X_p
\end{tikzcd}\end{center}
Note that $\img(\pi_{q-p+1} X_{p-1} \to \pi_{q-p+1} X_{p-r}) = D_r^{p-r,q-r+1}$ and $\img(\pi_{q-p} X_{p+r-1} \to \pi_{q-p} X_p) = D_r^{p, q}$, so that by taking images along the vertical maps in the above diagram we obtain a couple $\tilde C_r\colon D_r^{*,*} \xto{\tilde h_r} \tilde E_r^{*,*} \xto{\tilde f_r} D_r^{*,*}$. To finish the proof it is enough to show that $C_r \isom \tilde C_r$. One checks easily that $C_1 = \tilde C_1$, so it is enough to show the following, for every $r \ge 1$: Assuming that $\tilde C_r$ is an \emph{exact} couple, then $\tilde C_{r+1}$ is the derived couple of $\tilde C_r$. This boils down to the following claims:
\begin{itemize}
	\item For every $r \ge 1$ the differential $\tilde d_r\colon \tilde E_r^{*,*} \to \tilde E_r^{*,*}$ is given by $\tilde d_r = \tilde h_r \comp \tilde f_r$. This follows easily from the functoriality of $X(*,*)$, e.g. the map $\delta\colon \pi_{q-p} X(p, p-r) \to \pi_{q-p-1} X(p+r, p)$ factors as $\pi_{q-p} X(p, p-r) \to \pi_{q-p} X_p \to \pi_{q-p-1} X(p+r, p)$, where the first map is used to define $\tilde f_r$ and the second one is used to define $\tilde h_r$.

	\item For every $r \ge 1$, the map $\tilde f_{r+1}$ is equal to $\tilde f_r$ after identifying $\tilde E_{r+1}^{*,*} \isom \ker(\tilde d_r) / \img(\tilde d_r)$. To see this, we use the following commutative diagram:
	\begin{center}\begin{tikzcd}
		\pi_{q-p+1} X_{p-1} \arrow[r] \arrow[d,equal] & \pi_{q-p} X(p+r, p-1) \arrow[r] \arrow[d] & \pi_{q-p} X_{p+r} \arrow[d]\\
		\pi_{q-p+1} X_{p-1} \arrow[r] \arrow[d] & \pi_{q-p} X(p+r-1, p-1) \arrow[r] \arrow[d] & \pi_{q-p} X_{p+r-1} \arrow[d]\\
		\pi_{q-p+1} X_{p-r} \arrow[r] \arrow[d,"g"] & \pi_{q-p} X(p, p-r) \arrow[r] \arrow[d] & \pi_{q-p} X_p \arrow[d,equal]\\
		\pi_{q-p+1} X_{p-r-1} \arrow[r] & \pi_{q-p} X(p, p-r-1) \arrow[r] & \pi_{q-p} X_p
	\end{tikzcd}\end{center}
	The image of the middle vertical map is $\tilde E_r^{p,q}$ and the image of the vertical composition of maps is $E_{r+1}^{p,q}$. The diagram provides the desired identification $\tilde E_{r+1}^{*,*} \isom \ker(\tilde d_r) / \img(\tilde d_r)$ (cf. the proof of \cite[Proposition 1.2.2.7]{lurie-higher-algebra}). The maps $\tilde f_r$ and $\tilde f_{r+1}$ are obtained from the right horizontal maps by taking images along the middle vertical maps and the composition of vertical maps, respectively. The claim about $\tilde f_{r+1}$ and $\tilde f_r$ thus follows directly from the fact that the lower right vertical map in the diagram is the identity.

	\item For every $r \ge 1$, assuming that $\tilde C_r$ is exact, the map $\tilde h_{r+1}$ is equal to $\tilde h_r \comp g^{-1}$ after identifying $\tilde E_{r+1}^{*,*} \isom \ker(\tilde d_r) / \img(\tilde d_r)$. Here $g^{-1}$ denotes the process of choosing any preimage under $g$, using the fact that the choice does not matter. To prove this one can equivalently show that $\tilde h_r$ is equal to $\tilde h_{r+1} \comp g$. This follows from the above diagram in a similar way as for $\tilde f_{r+1}$ and $\tilde f_r$.
\end{itemize}
This finishes the proof by induction on $r$.
\end{proof}

\begin{corollary} \label{rslt:spectral-sequence-in-exact-sequence}
Let $\mathcal C$ be a stable $\infty$-category with $t$-structure and let $(X_n)_n \in \Tow(\mathcal C)$ with associated spectral sequence $E_*^{*,*}$ as defined in \cref{def:spectral-sequence-of-tower}. For $r \ge 1$ and $p, q \in \Z$ let $D_r^{p,q} := \img(\pi_{q-p} X_{p+r-1} \to X_p)$ (where we set $X_p = 0$ for $p < 0$). Then for all $r$, $p$ and $q$ there is a long exact sequence
\begin{align*}
	\dots \to E_r^{p,q+1} \to D_r^{p+1,q+1} \xto{g} D_r^{p,q} \to E_r^{p+r-1,q+r-1} \to \dots
\end{align*}
where the map $g$ is induced from the maps $X_p \to X_{p-1}$ of the tower.
\end{corollary}
\begin{proof}
This is just rephrasing that $E_r^{*,*}$ fits into an exact couple with $D_r^{*,*}$ as shown in \cref{rslt:spectral-sequence-comes-from-exact-couple}.
\end{proof}

Before we can proceed further, we need a short detour on systems of objects in abelian categories (cf. \cite[\S3.2]{akhil-thick-subcategory-theorem}).

\begin{definition}
Let $\mathcal A$ be an abelian category and let $(A_n)_n = (\dots \to A_n \to \dots \to A_1 \to A_0)$ be an inverse system of objects of $\mathcal A$.
\begin{defenum}
	\item For any $r \ge 0$, the \emph{$r$-th derived inverse system} of $(A_n)_n$ is the inverse system $(\img(A_{n+r} \to A_n))_n$.
	\item We say that $(A_n)_n$ is \emph{constant} if all transition maps $A_{n+1} \to A_n$ are constant. We say that $(A_n)_n$ is \emph{constant from index $n_0 \ge 0$} if the inverse system $(A_n)_{n\ge n_0}$ is constant. In this case the \emph{value} of $(A_n)_n$ is $A := \varprojlim_n A_n$, which is naturally isomorphic to $A_n$ for $n \ge n_0$. We say $(A_n)_n$ is \emph{eventually constant} if $(A_n)_n$ is constant from index $n_0$ for some $n_0 \ge 0$.
\end{defenum}
\end{definition}

\begin{lemma} \label{rslt:short-exact-sequence-and-r-th-derived-systems}
Let $\mathcal A$ be an abelian category and let $(A_n)_n \to (B_n)_n \to (C_n)_n$ be an exact sequence of inverse systems in $\mathcal A$. Fix integers $r, n_0 \ge 0$.
\begin{lemenum}
	\item Suppose that $0 \to (A_n)_n \to (B_n)_n$ is exact and that both the $r$-th derived systems of $(B_n)_n$ and $(C_n)_n$ are constant from index $n_0$. Then the $r$-th derived system of $(A_n)_n$ is constant from index $n_0$.
	\item Suppose that $(B_n)_n \to (C_n)_n \to 0$ is exact and that both the $r$-th derived systems of $(A_n)_n$ and $(B_n)_n$ are constant from index $n_0$. Then the $r$-th derived system of $(C_n)_n$ is constant from index $n_0$.
	\item Suppose that $0 \to (A_n)_n \to (B_n)_n \to (C_n)_n \to 0$ is exact and that both the $r$-th derived systems of $(A_n)_n$ and $(C_n)_n$ are constant from index $n_0$. Then the $2r$-th derived system of $(B_n)_n$ is constant from index $n_0$.
\end{lemenum}
\end{lemma}
\begin{proof}
Apart from the claims about the index $n_0$ this is precisely \cite[Lemma 3.10]{akhil-thick-subcategory-theorem}. However, the proof in the reference immediately provides the additional claims about $n_0$ as well.
\end{proof}

\begin{lemma}[{cf. \cite[Lemma 3.11]{akhil-thick-subcategory-theorem}}] \label{rslt:long-exact-sequence-and-r-th-derived-systems}
Let $\mathcal A$ be an abelian category and let
\begin{align*}
	(A_n)_n \to (B_n)_n \to (C_n)_n \to (D_n)_n \to (E_n)_n
\end{align*}
be an exact sequence of inverse systems in $\mathcal A$. Suppose that the $r$-th derived systems of $(A_n)_n$, $(B_n)_n$, $(D_n)_n$ and $(E_n)_n$ are constant from index $n_0 \ge 0$. Then the $2r$-th derived system of $(C_n)_n$ is constant from index $n_0$.
\end{lemma}
\begin{proof}
This follows by applying \cref{rslt:short-exact-sequence-and-r-th-derived-systems} three times.
\end{proof}

We are now in the position to show that strongly constant towers automatically have a fast converging spectral sequence. In fact, a weaker version of ``strongly constant'' is sufficient. Let us make the following definitions:

\begin{definition}
Let $\mathcal C$ be a stable $\infty$-category with $t$-structure.
\begin{defenum}
	\item A tower $(X_n)_n \in \Tow(\mathcal C)$ is called \emph{$t$-nilpotent of index $\le h$} if for all $n \ge h$ the map $\pi_* X_{n+h} \to \pi_* X_n$ is zero.
	\item A tower $(X_n)_n \in \Tow(\mathcal C)$ is called \emph{strongly $t$-constant of index $h$} if $X := \varprojlim_n X_n \in \mathcal C$ exists and $\cofib(X \to (X_n)_n) \in \Tow(\mathcal C)$ is $t$-nilpotent of index $h$.
	\item A cosimplicial object $M^\bullet$ in $\mathcal C$ is called \emph{strongly $t$-constant of index $h$} if the associated tower $(\Tot_n(M^\bullet))_n$ is strongly $t$-constant of index $h$.
\end{defenum}
\end{definition}

\begin{definition}
Let $E_*^{*,*}$ be a spectral sequence in an abelian category $\mathcal A$. We say that $E_*^{*,*}$ has a \emph{horizontal vanishing line at height $h \ge 0$} if for all $p$, $q$ and $r$ with $p \ge h$ and $r \ge h+1$ we have $E_r^{p,q} = 0$.
\end{definition}

\begin{remarks}
\begin{remarksenum}
	\item In general, being strongly $t$-constant of index $h$ is much weaker than being strongly constant of index $h$. However, as we will see below, under mild assumptions on $\mathcal C$ a strongly $t$-constant tower $(X_n)_n$ still satisfies $\pi_*(\varprojlim_n X_n) = \varprojlim_n \pi_* X_n$. Note also that if there is a constant $d$ such that the composition of $d$ phantom maps in $\mathcal C$ is zero then if a tower $(X_n)_n \in \Tow(\mathcal C)$ is strongly $t$-constant of index $h$ it is automatically strongly constant of index $\le hd$. For example this is true if $\mathcal C$ is the category of spectra and $d = 2$.

	\item Similar to \cref{rmk:refined-nilpotent-tower-index} one can define a refined version of the horizontal vanishing line by making the page at which it occurs independent of $h$. We decided not to do that for the same reasons as for nilpotent towers.
\end{remarksenum}
\end{remarks}

The next proposition shows that for a tower to be strongly $t$-constant is equivalent to its associated spectral sequence having a horizontal vanishing line.

\begin{proposition}[{cf. \cite[Proposition 3.12]{akhil-thick-subcategory-theorem}}] \label{rslt:strongly-t-constant-equiv-vanishing-line}
Let $\mathcal C$ be a stable $\infty$-category with $t$-structure such that countable products exist in $\mathcal C$ and are compatible with the $t$-structure. Let $(X_n)_n \in \Tow(\mathcal C)$ with associated spectral sequence $E_*^{*,*}$ as defined in \cref{def:spectral-sequence-of-tower}.
\begin{propenum}
	\item If $(X_n)_n$ is strongly $t$-constant of index $h$ then $E_*^{*,*}$ has a horizontal vanishing line at height $3h$.
	\item If $E_*^{*,*}$ has a horizontal vanishing line at height $h$ then $(X_n)_n$ is strongly $t$-constant of index $2h$.
\end{propenum}
Moreover, if $(X_n)_n$ is strongly $t$-constant of index $h$ then there is a natural isomorphism
\begin{align*}
	\pi_*(\varprojlim_n X_n) = \varprojlim_n \pi_* X_n = \img(\pi_* X_{3h} \to \pi_* X_{h}).
\end{align*}
\end{proposition}
\begin{proof}
Let $X := \varprojlim_n X_n$ (which exists because countable products in $\mathcal C$ exist) and let $(N_n)_n := \cofib(X \to (X_n)_n) \in \Tow(\mathcal C)$. Note that there is an exact sequence of graded inverse systems in $\mathcal C^\heartsuit$,
\begin{align*}
	(\pi_* N_n)_n \to (\pi_* X)_n \to (\pi_* X_n)_n \to (\pi_* N_n)_n \to (\pi_* X)_n \to (\pi_* X_n)_n,
\end{align*}
where the first and the second-to-last map lower the grading by $1$ each. Obviously $(\pi_* X)_n$ is constant, as is any of its derived systems.

We first prove (i), so let us assume that $(X_n)_n$ is strongly $t$-constant of index $h$, i.e. that $(N_n)_n$ is $t$-nilpotent of index $h$. By definition of $t$-nilpotent towers the $h$-th derived system of $(\pi_* N_n)_n$ is constant from index $h$. Thus by \cref{rslt:long-exact-sequence-and-r-th-derived-systems} the above exact sequence implies that the $2h$-th derived system of $(\pi_* X_n)_n$ is constant from index $h$. From \cref{rslt:spectral-sequence-in-exact-sequence} we deduce that $E_{2h+1}^{p,q} = 0$ for $p \ge 3h$. It follows that $E_r^{p,q} = 0$ for $r \ge 2h+1$ and $p \ge 3h$, finishing the proof of (i).

Let us now prove the last part of the claim, so we again assume that $(X_n)_n$ is strongly $t$-constant of some index $h$. By the proof of (i) we deduce that the $2h$-th derived system of $(\pi_* X_n)_n$ is constant from index $h$. This already implies the second identity. It remains to prove the first one. As both limits in the claim are invariant under the transformation $(X_n)_n \leadsto (X_{cn+d})_n$ for some constants $c \ge 1$ and $d \ge 0$, we can replace $(X_n)_n$ by $(X_{2hn+h})_n$ and hence assume that the first derived system of $(\pi_* X_n)_n$ is constant, say of value $A_*$. Denoting $f_n\colon X_n \to X_{n-1}$ we deduce that for all $n \ge 1$, the map $X_n \to \img f_n \isoto \img f_{n+1}$ is split surjective, so that $X_n \isom \img f_{n+1} \dsum \ker f_n$. Now note that the inverse limit $X = \varprojlim_n X_n$ fits into a cofiber sequence $X \to \prod_n X_n \to \prod_n X_n$ where the second map is given by $(\id - f_{n+1})_n$. As countable products are compatible with the $t$-structure in $\mathcal C$ we have $\pi_*(\prod_n X_n) = \prod_n \pi_*X_n$. We thus obtain a long exact sequence
\begin{align*}
	\dots \to \prod_n \pi_{k+1}X_n \to \prod_n \pi_{k+1}X_n \to \pi_k X \to \prod_n \pi_k X_n \to \prod_n \pi_k X_n \to \dots
\end{align*}
From the above splitting property of $(\pi_* X_n)_n$ it follows easily that $\prod_n \pi_{k+1} X_n \to \prod_n \pi_{k+1} X_n$ is surjective. Hence $\pi_k X \injto \prod_n \pi_k X_n$ is injective, so that $\pi_k X$ is isomorphic to the kernel of $\prod_n \pi_k X_n \to \prod_n \pi_k X_n$. Using that the first derived inverse system of $(\pi_k(X_n))_n$ is constant, this kernel is easily seen to be isomorphic to $A_*$.

It remains to prove (ii), so assume that $E_*^{*,*}$ has a horizontal vanishing line at height $h$. By \cref{rslt:spectral-sequence-in-exact-sequence} we deduce that for all $p \ge h$ the map
\begin{align*}
	g\colon \img(\pi_* X_{p+h+1} \to \pi_* X_{p+1}) \isoto \img(\pi_* X_{p+h} \to \pi_* X_p)
\end{align*}
is an isomorphism. In other words the $h$-th derived system of $(\pi_* X_n)_n$ is constant from index $h$. From the exact sequence of inverse systems at the beginning of the proof and from \cref{rslt:long-exact-sequence-and-r-th-derived-systems} we deduce that the $2h$-th derived system of $(\pi_* N_n)_n$ is constant from index $h$. Clearly $\varprojlim_n N_n = 0$, so from the proof of the last part of the claim we deduce $\varprojlim_n \pi_*(N_n) = \pi_*(\varprojlim_n N_n) = 0$. Hence the constant value of the $2h$-th derived system of $(\pi_* N_n)_n$ must be $0$, i.e. for all $n \ge h$, $\pi_* N_{n+2h} \to \pi_* N_n$ is the zero map. Thus $(N_n)_n$ is $t$-nilpotent and hence $(X_n)_n$ is strongly $t$-constant.
\end{proof}

We now come back to the promised analogue of \cref{rslt:f*M-strongly-constant-implies-M-strongly-constant-if-f-descendable} in terms of homotopy groups. In contrast to strongly constant towers, the property of being strongly $t$-constant is not preserved by exact functors in general. Therefore the following lemma requires additional hypotheses compared to \cref{rslt:f*M-strongly-constant-implies-M-strongly-constant-if-f-descendable}.

\begin{lemma} \label{rslt:f*M-strongly-t-constant-implies-M-strongly-t-constant-if-f-descendable}
Let $f\colon Y \to X$ be a morphism of affine analytic spaces over $(V,\mm)$ which is descendable of index $m$ and let $(M_n)_n \in \Tow(\D(X))$. Assume that there is a constant $h \ge 0$ such that $((f_*f^*)^k M_n)_n \in \Tow(\D(X))$ is strongly $t$-constant of index $h$ for all $k$. Then $(M_n)_n$ is strongly $t$-constant of index $s(m) h$, where $s(m)$ is the constant from \cref{rslt:f-zero-only-depends-on-index}.
\end{lemma}
\begin{proof}
We can argue as in the proof of \cref{rslt:descendable-nilpotent-criterion}, noting that the tower $(M_n)_n$ can be constructed from the towers $((f_\natural \ri_Y)^k(M_n))_n$ using finite limits and retracts. More precisely, we have the following:
\begin{itemize}
	\item Let $(R_n)_n, (S_n)_n \in \Tow(\D(X))$ such that $(R_n)_n$ is a retract of $(S_n)_n$. If $(S_n)_n$ is strongly $t$-constant of index $h$ then $(R_n)_n$ is strongly $t$-constant of index $\le h$. To see this we can easily reduce to the case that that $(S_n)_n$ is $t$-nilpotent of index $h$. From $(R_n)_n$ being a retract of $(S_n)_n$ we deduce that $R_{n+h} \to R_n$ factors over $S_{n+h} \to S_n$ for all $n \ge 0$, which implies the claim.

	\item Let $(R_n)_n \to (S_n)_n \to (T_n)_n$ be a cofiber sequence in $\Tow(\D(X))$. If $(R_n)_n$ and $(T_n)_n$ are strongly $t$-constant of index $h_1$ and $h_2$ respectively then $(S_n)_n$ is strongly $t$-constant of index $\le h_1 + h_2$. To see this we can again reduce to the case that $(R_n)_n$ and $(T_n)_n$ are $t$-nilpotent. Then the claim follows similarly as in the proof of \cite[Proposition 3.5]{akhil-thick-subcategory-theorem} (this is essentially the same argument as in the proof of \cref{rslt:descendable-nilpotent-criterion}).
\end{itemize}
By construction of $s(m)$ this finishes the proof.
\end{proof}

It is often the case that the assumptions of \cref{rslt:f*M-strongly-t-constant-implies-M-strongly-t-constant-if-f-descendable} can be deduced from a stronger property of $(f^*M_n)_n \in \Tow(\D(Y))$. One example of this is when $(f^*M_n)_n$ is strongly constant, because this property is preserved by exact functors (however, this gives nothing new compared to \cref{rslt:f*M-strongly-constant-implies-M-strongly-constant-if-f-descendable}). Another case is the following, which is also the only case we are interested in:

\begin{corollary} \label{rslt:f*M-split-implies-M-split-and-hor-vanishing-line}
Let $f\colon Y \to X$ be a morphism of affine analytic spaces over $(V,\mm)$ which is descendable of index $m$ and let $M^\bullet$ be a cosimplicial object in $\D(X)$. Assume that $f^*M^\bullet$ is split. Then:
\begin{corenum}
	\item $M^\bullet$ is strongly $t$-constant of index $2s(m)$.
	\item The spectral sequence associated to $M^\bullet$ has a horizontal vanishing line at height $6s(m)$.
\end{corenum}
Here $s(m)$ is the constant from \cref{rslt:f-zero-only-depends-on-index}.
\end{corollary}
\begin{proof}
Being a split cosimplicial object is preserved under any functor, so in particular $(f_\natural \ri_Y)^k(M^\bullet)$ is a split cosimplicial object of $\D(X)$ for all $k \ge 1$. Using \cref{rslt:strongly-t-constant-equiv-vanishing-line} and \cref{rslt:f*M-strongly-t-constant-implies-M-strongly-t-constant-if-f-descendable} it is therefore enough to show the following: If $N^\bullet$ is a split simplicial object in $\D(X)$ then the associated spectral sequence $E_*^{*,*}$ has a horizontal vanishing line at height $1$. By the dual version of \cite[Remark 1.2.4.4]{lurie-higher-algebra}, for every $q \in \Z$, $E_1^{*,q}$ is the normalized chain complex associated to the simplicial object $\pi_q N^\bullet$ of $\D(X)^\heartsuit$. By \cite[Proposition 1.2.3.17]{lurie-higher-algebra} this normalized chain complex is quasi-isomorphic to the unnormalized chain complex
\begin{align*}
	\pi_q N^0 \to \pi_q N^1 \to \pi_q N^2 \to \dots
\end{align*}
where the differentials are given by the alternating sums over the coface maps. Since $N^\bullet$ is split by assumption, this unnormalized chain complex is split exact and hence has vanishing cohomology in degree greater zero. On the other hand, this cohomology is precisely $E_2^{*,q}$ so that $E_2^{p,q} = 0$ for $p \ge 1$. Therefore $E_*^{*,*}$ has a horizontal vanishing line at height $1$, as desired.
\end{proof}

\begin{remark}
One can also attempt to show that every split simplicial object is strongly constant, which would simplify the proof of \cref{rslt:f*M-split-implies-M-split-and-hor-vanishing-line}.
\end{remark}

We now come to the promised version of descendability. As mentioned above, we want to work with ``flat covers'' of $\AnSpec \mathcal A$ for an analytic ring $\mathcal A$ over $(V,\mm)$. Here is the precise definition of that notion:

\begin{definition}
Let $f\colon \mathcal A \to \mathcal B$ be a map of analytic rings over $(V,\mm)$.
\begin{defenum}
	\item A map $\mathcal A \to \mathcal A'$ of analytic rings is called an \emph{fs-morphism} if it is flat (i.e. $- \tensor_{\mathcal A} \mathcal A'$ is $t$-exact) and steady. An \emph{fs-cover} of $\mathcal A$ is a family $(\mathcal A \to \mathcal A_i)_i$ of fs-morphisms such that the family of functors $- \tensor_{\mathcal A} \mathcal A_i$ is conservative.

	\item $f$ is called \emph{fs-descendable of index $\le d$} if there is an fs-cover $(\mathcal A \to \mathcal A_i)_i$ such that all the induced maps $\mathcal A_i \to \mathcal B \tensor_{\mathcal A} \mathcal A_i$ are descendable of index $\le d$. We say that $f$ is \emph{fs-descendable} if it is fs-descendable of some index.

	\item $f$ is called \emph{weakly fs-descendable (of index $\le d$)} if it is an iterated $\kappa$-filtered colimit of fs-descendable maps (of index $\le d$), where $\kappa$ is any uncountable regular cardinal such that $(V,\mm)$ is $\kappa$-compact.
\end{defenum}
\end{definition}

\begin{theorem} \label{rslt:fs-descendability-properties}
\begin{thmenum}
	\item \label{rslt:pi-0-weakly-descendable-implies-descent} Quasi-coherent sheaves descend along weakly fs-descendable morphisms.

	\item \label{rslt:stability-of-pi-0-weakly-descendable-maps}  (Weakly) fs-descendable morphisms are stable under base-change, fs-descendable morphisms are stable under composition, and weakly fs-descendable morphisms are stable under $\kappa$-filtered colimits, where $\kappa$ is any uncountable regular cardinal such that $(V,\mm)$ is $\kappa$-compact.

	\item \label{rslt:filtered-colim-of-bounded-pi-0-desc-is-weakly-pi-0-desc} Suppose that $(V,\mm)$ is $\omega_1$-compact (e.g. $\mm$ is countably generated). Let $(\mathcal A \to \mathcal B_i)_{i\in I}$ be a filtered diagram of maps of analytic rings over $(V,\mm)$ which are fs-descendable of bounded index $\le d$ with respect to a common fs-cover of $\mathcal A$. Then $\mathcal A \to \varinjlim_i \mathcal B_i$ is weakly fs-descendable of index $\le 2d$.
\end{thmenum}
\end{theorem}
\begin{proof}
Since fs-morphisms are stable under base-change, it follows easily from \cref{rslt:descendable-stable-under-base-change} that fs-descendable morphisms are stable under base-chage. It is also clear that $\kappa$-filtered colimits of weakly fs-descendable morphisms are weakly fs-descendable. This proves (ii). Part (iii) follows from \cref{rslt:descendable-stable-under-countable-colimit} (as in the proof of \cref{rslt:filtered-colim-of-bounded-desc-is-weakly-desc}).

It remains to prove (i). By \cref{rslt:descent-of-qcoh-is-stable-under-w1-filtered-colim} we are reduced to showing that quasicoherent sheaves descend along fs-descendable morphisms. Let $f\colon \mathcal A \to \mathcal B$ be such a morphism and let $d$ be an integer and $(\mathcal A \to \mathcal A_i)_{i\in I}$ an fs-cover such that all $f_i\colon \mathcal A_i \to \mathcal B_i := \mathcal B \tensor_{\mathcal A} \mathcal A_i$ are descendable of index $\le d$. It follows easily that $f$ is steady (note that the family of functors $- \tensor_{\mathcal B} \mathcal B_i$ is conservative, now apply any criterion for steadiness, e.g. \cref{rslt:steady-map-of-analytic-rings-equiv-base-change}). As in the proof of \cref{rslt:descendable-implies-limit-of-categories} we can thus apply Lurie's \cite[Corollary 4.7.5.3]{lurie-higher-algebra} to see that we only have to show the following two properties of the base-change functor $f^*\colon \D(\mathcal A) \to \D(\mathcal B)$:
\begin{enumerate}[(a)]
	\item $f^*$ is conservative.
	\item Given a cosimplicial object $M^\bullet \in \D(\mathcal A)$ such that $f^* M^\bullet \in \D(\mathcal B)$ is split, then $f^*$ preserves the totalization of $M^\bullet$, i.e. $\Tot(f^* M^\bullet) = f^* \Tot(M^\bullet)$.
\end{enumerate}
For every $i \in I$ denote $g_i\colon \mathcal A \to \mathcal A_i$ and $g'_i\colon \mathcal B \to \mathcal B_i$. Then we have a commuting diagram
\begin{center}\begin{tikzcd}
	\D(\mathcal B_i) & \arrow[l,swap,"g_i'^*"] \D(\mathcal B)\\
	\D(\mathcal A_i) \arrow[u,"f_i^*"] & \arrow[l,swap,"g_i^*"] \D(\mathcal A) \arrow[u,swap,"f^*"]
\end{tikzcd}\end{center}
We first show that $f^*$ satisfies (a). Suppose we are given a morphism $M \to N$ in $\D(\mathcal A)$ such that $f^* M \isoto f^* N$ is an isomorphism in $\D(\mathcal B)$. We need to show that $M \to N$ is an isomorphism. It is enough to show that for all $i \in I$, $g_i^* M \to g_i^* N$ is an isomorphism. But we have $f_i^* g_i^* M = g_i'^* f^* M$ (and similarly for $N$), so that $f_i^* g_i^* M \isoto f_i^* g_i^* N$ is an isomorphism. The claim thus follows from the fact that all $f_i^*$ are conservative (because $f_i$ is descendable).

It remains to prove that $f^*$ satisfies (b), so let $M^\bullet \in \D(\mathcal A)$ be such that $f^* M^\bullet \in \D(\mathcal B)$ is split. Then, for every $i \in I$, $f_i^* g_i^* M = g_i'^* f^* M^\bullet \in \D(\mathcal B)$ is split. Since $f_i$ is descendable of index $\le d$, \cref{rslt:f*M-split-implies-M-split-and-hor-vanishing-line} implies that $g_i^* M^\bullet$ is strongly $t$-constant of index $h := 2s(d)$ (where $s(d)$ is a constant defined in \cref{rslt:f-zero-only-depends-on-index} and only depends on $d$). In other words, the $h$-th derived system of $(\pi_* \Tot_n(g_i^* M^\bullet))_n$ is constant from index $h$.

By assumption $g_i^*$ is $t$-exact and in particular commutes with homology. It follows that $g_i^*$ maps the $h$-th derived system of $(\pi_* \Tot_n(M^\bullet))_n$ to the $h$-th derived system of $(\pi_* \Tot_n(g_i^* M^\bullet))_n$. In other words, the $h$-th derived system of $(\pi_* \Tot_n(M^\bullet))_n$ is mapped to a system which is constant from index $h$, for all $i$. Since the functors $g_i^*$ form a conservative family, it follows that the $h$-th derived system of $(\pi_* \Tot_n(M^\bullet))_n$ is constant from index $h$, i.e. that $M^\bullet$ is strongly $t$-constant of index $h$. By \cref{rslt:strongly-t-constant-equiv-vanishing-line} we have $\pi_* \Tot(M^\bullet) = \varprojlim_n \pi_* \Tot_n(M^\bullet)$, and since $g_i^*$ commutes with finite limits and with homology groups we deduce
\begin{align*}
	g_i^* \Tot(M^\bullet) = \Tot(g_i^* M^\bullet)
\end{align*}
for every $i \in I$.

It is now easy to deduce that $f^*$ preserves the totalization of $M^\bullet$: There is a natural morphism $f^* \Tot(M^\bullet) \to \Tot(f^* M^\bullet)$ and since the forgetful functor $f_*$ is conservative, it is enough to to show that the induced map $g_i^* f_* f^* \Tot(M^\bullet) \to g_i^* f_* \Tot(f^* M^\bullet)$ is an isomorphism for all $i \in I$. By steady base-change (see \cref{rslt:steady-map-of-analytic-rings-equiv-base-change}) we have
\begin{align*}
	g_i^* f_* f^* \Tot(M^\bullet) &= f_{i*} g_i'^* f^* \Tot(M^\bullet) = f_{i*} f_i^* g_i^* \Tot(M^\bullet) = f_{i*} f_i^* \Tot(g_i^* M^\bullet) =\\
	&= f_{i*} \Tot(f_i^* g_i^* M^\bullet),\\
	g_i^* f_* \Tot(f^* M^\bullet) &= f_{i*} g_i'^* \Tot(f^* M^\bullet) = f_{i*} \Tot(g_i'^* f^* M^\bullet) = f_{i*} \Tot(f_i^* g_i^* M^\bullet).
\end{align*}
(In the first line we used that $f_i^*$ preserves $f_i^*$-split totalizations and in the second line we used that $f^* M^\bullet$ is split so that its totalization commutes with any exact functor, in particular with $g_i'^*$.) This finishes the proof.
\end{proof}

\begin{corollary}
Let $f\colon \mathcal A \to \mathcal B$ be a map of analytic rings over $(V, \mm)$ and $d \ge 0$ an integer. Then:
\begin{corenum}
	\item If $f$ is fs-descendable of index $d$ then every $f^*$-split cosimplicial object $M^\bullet$ of $\D(\mathcal A)$ is strongly $t$-constant of index $\le 2s(d)$.

	\item \label{rslt:weakly-fs-descendable-implies-right-boundedness-descends} There is a constant $b(d) \in \Z$ only depending on $d$ with the following property: If $f$ is weakly fs-descendable of index $d$ and $M \in \D(\mathcal A)$ satisfies $f^* M \in \D_{\ge0}(\mathcal B)$, then $M \in \D_{\ge b(d)}(\mathcal A)$.
\end{corenum}
\end{corollary}
\begin{proof}
Part (i) follows from the proof \cref{rslt:pi-0-weakly-descendable-implies-descent}. We now prove (ii), so assume that $f$ is weakly fs-descendable of index $d$ and suppose $M \in \D(\mathcal A)$ satisfies $f^* M \in \D_{\ge0}(\mathcal B)$. By \cref{rslt:pi-0-weakly-descendable-implies-descent} we have $\D(\mathcal A) = \varprojlim_{n\in\Delta} \D(\mathcal B^{\tensor_{\mathcal A} n+1})$, so for the cosimplicial object $M^\bullet = M \tensor_{\mathcal A} \mathcal B^{\tensor_{\mathcal A} \bullet+1}$ we have $M = \Tot(M^\bullet)$. Writing $f$ as a an iterated $\kappa$-filtered colimit of fs-descendable maps $f_i\colon \mathcal B_i \to \mathcal A$, we can similarly write $M^\bullet$ as an iterated $\kappa$-filtered colimit of the associated cosimplicial objects $M_i^\bullet$. Then each $M_i^\bullet$ is $f_i^*$-split, hence by (i) each $M_i^\bullet$ is strongly $t$-constant of index $\le h := 2s(d)$. Since $\pi_*$ and $\varprojlim_n$ commute with $\kappa$-filtered colimits, it follows directly from the definitions that $M^\bullet$ is strongly $t$-constant of index $\le h$. In particular, by \cref{rslt:strongly-t-constant-equiv-vanishing-line} we have
\begin{align*}
	\pi_*(M) = \pi_*(\Tot(M^\bullet)) = \img(\pi_*\Tot_{3h}(M^\bullet) \to \pi_*\Tot_h(M^\bullet)))
\end{align*}
But each $M^n$ lies in $\D_{\ge0}$ and $\Tot_h$, $\Tot_{3h}$ are finite limits of fixed form (only depending on $h$), which implies the claim.
\end{proof}

\subsection{Schemes and Discrete Adic Spaces} \label{sec:andesc.scheme}

In the previous subsections we have developed a very general theory of analytic rings and descent. Now we finally come to some examples of that theory by discussing how schemes (and discrete adic spaces) fit into that setting. To every discrete ring $A$ we will associate an analytic ring $A_\solid$, called the \emph{solid analytic ring} associated to $A$. In fact we can more generally associate a solid analytic ring $(A, A^+)_\solid$ to every discrete Huber pair $(A, A^+)$. The basic construction of $(A, A^+)_\solid$ follows \cite[\S8]{condensed-mathematics}, but we expand upon it by proving several additional properties of the $\infty$-category of solid analytic rings. We also generalize the construction slightly by combining it with almost mathematics. With a good theory of solid analytic rings at hand, we can then glue them to obtain the $\infty$-categories of schemes and of discrete adic spaces. At the end of this subsection we discuss a 6-functor formalism for (solid) quasicoherent sheaves on schemes and discrete adic spaces, following \cite[\S11]{condensed-mathematics}.

Let us start with the basic definitions. In order to make the exposition easier to follow, we first work exclusively in the non-almost setting and only later introduce an almost version (see \cref{def:solid-almost-ring}) below. Most of the results about non-almost solid rings generalize directly to the almost world by viewing a solid almost ring as the almostification of an honest solid ring.

\begin{definition}
\begin{defenum}
	\item Let $A$ be a classical ring. We associate to $A$ the pre-analytic ring $A_\solid$ with
	\begin{align*}
		A_\solid[S] := \varinjlim_{A' \to A} \varprojlim_i A'[S_i]
	\end{align*}
	for every profinite set $S = \varprojlim_i S_i$ (with all $S_i$ finite), where $A'$ ranges through all finite-type classical $\Z$-algebras with a map to $A$.
\end{defenum}
\end{definition}

\begin{lemma} \label{rslt:A-solid-is-analytic-ring}
For every classical ring $A$, the pre-analytic ring $A_\solid$ defined above is an analytic ring. If $A$ is of finite type over $\Z$ then for every profinite set $S$ there is some set $I$ such that $A_\solid[S] = \prod_I A$.
\end{lemma}
\begin{proof}
The second claim follows from \cite[Theorem 5.4]{condensed-mathematics} (as in \cite[Corollary 5.5]{condensed-mathematics}). For the first claim, by definition we have $A_\solid = \varinjlim_{A' \to A} A'_\solid$, where $A'$ ranges over finite-type classical $\Z$-algebras with a map to $A$. Hence we can w.l.o.g. assume that $A$ is of finite type over $\Z$. But then the claim is shown in \cite[Theorem 8.13.(i)]{condensed-mathematics} (use also \cite[Proposition 12.24]{scholze-analytic-spaces}).
\end{proof}

\begin{definition}
\begin{defenum}
	\item A \emph{discrete Huber pair} is a pair $(A, A^+)$, where $A$ is a discrete ring and $A^+$ is an integrally closed subring of $\pi_0 A$. A morphism $(A, A^+) \to (B, B^+)$ of discrete Huber pairs is a morphism $A \to B$ of rings such that the image of $A^+$ under the induced map $\pi_0 A \to \pi_0 B$ lies in $B^+$.

	\item Let $(A, A^+)$ be a discrete Huber pair. We let
	\begin{align*}
		(A, A^+)_\solid \in \AnRing
	\end{align*}
	denote the following analytic ring structure on $A$: If $A$ is static then there is a map $A^+ \injto A$ and we let $(A, A^+)_\solid := A_{A^+_\solid/}$ be the induced analytic ring. For general $A$ we employ \cref{rslt:topological-invariance-of-analytic-ring-structures} to define $(A, A^+)_\solid$ as the analytic ring structure induced from $(\pi_0 A, A^+)_\solid$.

	\item \label{def:abberviated-version-of-solid-rings} Given discrete rings $A$ and $A^+$ and a map $A^+ \to \pi_0 A$ we also denote $(A, A^+)_\solid := (A, \tilde A^+)_\solid$, where $\tilde A^+ \subset \pi_0 A$ is the integral closure of the image of the map $\pi_0 A^+ \to \pi_0 A$. We further abbreviate, for every discrete ring $A$,
	\begin{align*}
		A_\solid := (A, \pi_0 A)_\solid \in \AnRing.
	\end{align*}

	\item We abbreviate $\D_\solid(A) := \D(A_\solid)$ and $\D_\solid(A, A^+) = \D((A, A^+)_\solid)$. A \emph{solid $(A, A^+)$-module} is an object of $\D_\solid(A, A^+)$.
\end{defenum}
\end{definition}

The first thing we want to do is check that the abbreviated version of $(A, A^+)_\solid$ from \cref{def:abberviated-version-of-solid-rings} is actually the induced analytic ring structure from $A^+_\solid$ (and not only from $\tilde A^+_\solid$). This allows us to efficiently compute the generators $(A, A^+)_\solid[S]$ in many cases.

\begin{lemma} \label{rslt:solid-of-map-A+-A-is-induced-ring-structure}
Let $A^+ \to A$ be a map of discrete rings and assume that $A$ is static. Then
\begin{align*}
	(A, A^+)_\solid = A_{A^+_\solid/}.
\end{align*}
In other words, for every profinite set $S$ we have
\begin{align*}
	(A, A^+)_\solid[S] = A \tensor_{A^+} A^+_\solid[S].
\end{align*}
\end{lemma}
\begin{proof}
First note that the map $A^+ \to A = \pi_0 A$ factors over $\pi_0 A^+$, so we can easily reduce to the case that $A^+$ is static. Let $\tilde A^+$ be the integral closure of the image of $A^+ \to A$. Then by definition we have $(A, A^+)_\solid = A_{\tilde A^+_\solid/}$, so it is enough to verify that $\tilde A^+_\solid = \tilde A^+_{A^+_\solid/}$. We are thus reduced to the following statement: Let $B \to B'$ be an integral map of classical rings; then $B'_\solid = B'_{B_\solid/}$. To prove this, first note that by writing $B'$ as a filtered colimit of finite-type classical $B$-algebras and pulling out the colimit on both sides of the claimed identity we can reduce to the case that $B \to B'$ is of finite type. By the same argument we can then reduce to the case that $B \to B'$ is of finite presentation. Then it comes via base-change from a finite map $B_0 \to B_0'$ of finite-type classical $\Z$-algebras. Write $B = \varinjlim_i B_i$, where $B_i$ ranges over finite-type classical $\Z$-algebras with maps $B_0 \to B_i \to B$. Let $B' = \varinjlim_i B'_i$ be the base-change along $B_0 \to B_0'$. Then $B_\solid = \varinjlim_i (B_i)_\solid$ and $B'_{B_\solid/} = \varinjlim_i (B'_i)_{(B_i)_\solid/}$, so we can reduce to the case that $B$ and $B'$ are of finite type over $\Z$. Now by \cref{rslt:A-solid-is-analytic-ring} we need to see that for every set $I$ we have $B' \tensor_B \prod_I B = \prod_I B'$. In fact we claim that this is true for any finite classical $B$-module $M$ in place of $B'$: If $M$ is free then this is clear, in general we can resolve $M$ by finite free modules to reduce to the free case.
\end{proof}

A general map $A^+ \to A$ of classical rings can be factored into a composition of a finite map and a polynomial algebra (in possibly infinitely many variables). If $A^+ \to A$ is finite (or more generally integral) then $A_\solid = (A, A^+)_\solid$ by definition, thus by \cref{rslt:solid-of-map-A+-A-is-induced-ring-structure} $A_\solid$ possesses the induced analytic ring structure from $A^+_\solid$. In this case many questions about $A_\solid$ can easily be reduced to $A^+_\solid$. Therefore, in order to understand the map $A^+_\solid \to A_\solid$ in general it remains to get a good understanding of polynomial algebras. We have the following result:

\begin{lemma} \label{rslt:finite-type-solid-pullback-preserves-limits}
Let $A$ be a classical finite-type $\Z$-algebra. Then the functor
\begin{align*}
	- \tensor_{A_\solid} A[T]_\solid\colon \D_\solid(A) \to \D_\solid(A[T])
\end{align*}
is $t$-exact and preserves all limits.
\end{lemma}
\begin{proof}
Let $f\colon \Spec A[T] \to \Spec A$ be the obvious map. In the following we make use of the shriek functors $f_!$ and $f^!$ defined in \cite[Theorem 8.13]{condensed-mathematics}. For all $M \in \D_\solid(A)$ we have $f^! M = f^*M \tensor_{A[T]_\solid} f^! A$ and $f^! A = A[T][1]$ by (the proof of) \cite[Theorem 8.13.(iv)]{condensed-mathematics}, hence $f^! M = f^* M[1]$. Therefore
\begin{align*}
	M \tensor_{A_\solid} A[T]_\solid = f^*M = f^! M[-1] = \IHom_{A[T]}(A[T], f^! M[-1]) = \IHom_A(f_! A[T][1], M)
\end{align*}
But $f_! A[T] = A((T^{-1}))/A[T][-1]$ (see e.g. \cite[Observation 8.11]{condensed-mathematics}), which implies
\begin{align*}
	M \tensor_{A_\solid} A[T]_\solid = \IHom_A(A((T^{-1}))/A[T], M).
\end{align*}
Clearly this functor preserves limits in $M$ and is left $t$-exact (hence $t$-exact).
\end{proof}

We will now prove some useful basic results about the analytic rings $(A, A^+)_\solid$. Namely, we show that the association $(A, A^+) \mapsto (A, A^+)_\solid$ exhibits the $\infty$-category of discrete Huber pairs as a full subcategory of the $\infty$-category of all analytic rings, compatible with all small colimits. Moreover, we show that every map of solid discrete analytic rings is steady, so that unconditional base-change for modules holds, just like one is used to from classical algebra.

\begin{proposition}
\begin{propenum}
	\item \label{rslt:solidification-is-fully-faithful-functor} The assignment $(A, A^+) \mapsto (A, A^+)_\solid$ defines a fully faithful functor from the $\infty$-category of discrete Huber pairs to $\AnRing$.

	\item \label{rslt:colimits-of-solid-discrete-rings} The functor in (i) preserves all non-empty colimits. In particular, for any diagram $(B, B^+) \from (A, A^+) \to (C, C^+)$ of discrete Huber pairs we have
	\begin{align*}
		(B, B^+)_\solid \tensor_{(A, A^+)_\solid} (C, C^+)_\solid = (B \tensor_A C, B^+ \tensor_{A^+} C^+)_\solid,
	\end{align*}
	where on the right-hand side we use the abbreviated notation from \cref{def:abberviated-version-of-solid-rings}.
\end{propenum}
\end{proposition}
\begin{proof}
To prove (i) we need to verify the following: Let $(A, A^+)$ and $(B, B^+)$ be two discrete Huber pairs and let $f\colon A \to B$ be a morphism of rings. Then $f$ is a morphism $(A, A^+)_\solid \to (B, B^+)_\solid$ of analytic rings if and only if it is a morphism of Huber pairs. To prove this, let us first show the ``if'' direction, so assume that $f$ is a morphism of Huber pairs. Given any $M \in \D_\solid(B, B^+)$ we need to see that the image of $M$ under the forgetful functor $\D(B) \to \D(A)$ lies in $\D_\solid(A, A^+)$. This can be checked on homology groups, so we can assume that $M$ is concentrated in degree $0$. Then we can view $M$ as an object in $\D(\pi_0 B, B^+)$ and it is enough to see that its image in $\D(\pi_0 A)$ lies in $\D_\solid(\pi_0 A, A^+)$; in other words we can assume that both $A$ and $B$ are concentrated in degree $0$. In this case $M$ lies in $\D_\solid(B, B^+)$ if and only if its image under the forgetful functor $\D(B) \to \D(B^+)$ lies in $\D_\solid(B^+)$ (and similarly for $A$). Hence we can assume $B = B^+$ and $A = A^+$ (note that by the assumption on $f$, $f$ restricts to a map $A^+ \to B^+$), i.e. we want to show that $f$ is a morphism of analytic rings $A_\solid \to B_\solid$. Writing $f$ as a filtered colimit of maps of finite-type classical $\Z$-algebras we can reduce to the case that $A$ and $B$ are of finite type over $\Z$. Then for every extremally disconnected set $S$, $B_\solid[S] = \varprojlim_i B[S_i]$ is a limit of discrete (hence solid) $A$-modules, therefore its image under the forgetful functor $\D(B) \to \D(A)$ lies in $\D_\solid(A)$, as desired.

Conversely, assume that $f$ is a morphism $(A, A^+)_\solid \to (B, B^+)_\solid$ of analytic rings. Then $f$ also induces a morphism $(\pi_0 A, A^+)_\solid \to (\pi_0 B, B^+)_\solid$ of analytic rings, so we can assume that $A$ and $B$ are concentrated in degree $0$. By the same argument as in the proof of \cite[Proposition 13.16]{scholze-analytic-spaces}, $f$ maps $A^+$ to $B^+$.

We now prove (ii). It is easy to see that the functor preserves non-empty filtered colimits. We now handle pushouts, so let $(B, B^+) \from (A, A^+) \to (C, C^+)$ be a diagram of discrete Huber pairs. Clearly its pushout in the $\infty$-category of Huber pairs is $(B \tensor_A C, (B^+ \tensor_{A^+} C^+)\tilde{})$, where $(B^+ \tensor_{A^+} C^+)\tilde{}$ denotes the integral closure of the image of $B^+ \tensor_{A^+} C^+ \to \pi_0(B \tensor_A C)$. Thus in order to verify that our functor preserves this pushout, we must show the claimed identity about the pushout of the associated analytic rings. By \cref{rslt:topological-invariance-of-analytic-ring-structures} this can be checked after passage to the induced structures on $0$-truncated analytic rings, which allows us to assume that $A$, $B$ and $C$ are concentrated in degree $0$. Using \cref{rslt:solid-of-map-A+-A-is-induced-ring-structure} and the explicit computation of pushouts (see \cref{rslt:pushouts-of-analytic-rings}) we can further reduce to the case $A = A^+$, $B = B^+$ and $C = C^+$. Now we want to show
\begin{align*}
	B_\solid \tensor_{A_\solid} C_\solid = (B \tensor_A C)_\solid.
\end{align*}
Both sides commute with filtered colimits in $(B \from A \to C)$, so we can assume that $A$, $B$ and $C$ are of finite type over $\Z$. Choose a surjection $B' := A[T_1, \dots, T_n] \to B$ of $A$-algebras. By \cref{rslt:solid-of-map-A+-A-is-induced-ring-structure} we have $B_\solid = (B, B')_\solid$ hence by the same reduction as above we can reduce to the case $B = B'$. By induction on $n$ we reduce to $B = A[T]$. By \cref{rslt:pushouts-of-analytic-rings} $B_\solid \tensor_{A_\solid} C_\solid$ is the analytic ring on $B \tensor_A C$ such that $(B_\solid \tensor_{A_\solid} C_\solid)[S]$ is the colimit of the repeated application of $(- \tensor_C C_\solid) \tensor_B B_\solid$ to $(B \tensor_A C)[S]$. We have
\begin{align*}
	((B \tensor_A C)[S] \tensor_C C_\solid) \tensor_B B_\solid = (B \tensor_A C_\solid[S]) \tensor_B B_\solid = B_\solid \tensor_{A_\solid} C_\solid[S].
\end{align*}
By using \cref{rslt:finite-type-solid-pullback-preserves-limits} (for $B = A[T]$) and the explicit definition of $C_\solid[S]$, we obtain $B_\solid \tensor_{A_\solid} C_\solid[S] = (B \tensor_A C)_\solid[S]$. But this object is both $B$-solid and $C$-solid, hence subsequent applications of $- \tensor_B B_\solid$ and $- \tensor_C C_\solid$ leave this object invariant. We deduce $(B_\solid \tensor_{A_\solid} C_\solid)[S] = (B \tensor_A C)_\solid[S]$. This finishes the case of pushouts.

To finish the proof of (ii) it remains to see that the functor in (i) preserves coproducts. Coproducts in the category of discrete Huber pairs are just tensor products over $(\Z, \Z)$ and hence translate to tensor products over $\Z_\solid$ on the side of analytic rings. Coproducts in $\AnRing$ are tensor products over $\Z$. But by the description of pushouts in \cref{rslt:pushouts-of-analytic-rings}, both tensor products are easily seen to agree.
\end{proof}

\begin{proposition}
\begin{propenum}
	\item \label{rslt:solid-discrete-equiv-nuclear} Let $(A, A^+)$ be a discrete Huber pair. Then $P \in \D_\solid(A, A^+)$ is nuclear if and only if it is discrete.

	\item \label{rslt:map-of-discrete-rings-is-steady} Let $(A, A^+)$ and $(B, B^+)$ be discrete Huber pairs. Then every morphism $f\colon (A, A^+)_\solid \to (B, B^+)_\solid$ is steady.
\end{propenum}
\end{proposition}
\begin{proof}
We start with the proof of (i), so let $(A, A^+)$ and $P \in \D_\solid(A, A^+)$ be given. Suppose first that $P$ is nuclear. Pick some extremally disconnected set $S = \varprojlim_i S_i$, where all $S_i$ are finite. Then
\begin{align*}
	(A, A^+)_\solid[S]^\vee &= \IHom_A((A, A^+)_\solid[S], A) = \IHom_A(A[S], A) = \IHom(S, A) =\\
	&= \varinjlim_i \IHom(S_i, A) = \varinjlim_i A^{S_i}.
\end{align*}
Consequently, by definition of nuclearity we get
\begin{align*}
	P(S) &= ((A, A^+)_\solid[S]^\vee \tensor_{(A, A^+)_\solid} P)(*) = (\varinjlim_i A^{S_i} \tensor_A P)(*) = \varinjlim_i P^{S_i}(*) =\\
	&= \varinjlim_i P(S_i).
\end{align*}
This immediately implies that $P$ is discrete. Conversely, if $P$ is discrete then it is a colimit of copies of $A$. On the other hand, $A$ is obviously nuclear over $(A, A^+)_\solid$ and nuclear objects are stable under colimits by \cref{rslt:nuclear-modules-stable-under-colim}, hence $P$ is nuclear.

We now prove (ii), so let $(A, A^+)$ and $(B, B^+)$ be given. In the following we will treat $A^+$ and $B^+$ in the more general way of \cref{def:abberviated-version-of-solid-rings}. By \cref{rslt:solidification-is-fully-faithful-functor} $f$ is induced by a morphism $(A, A^+) \to (B, B^+)$ of Huber pairs. We can factor $f$ as $(A, A^+)_\solid \xto{g} (B, A^+)_\solid \xto{h} (B, B^+)_\solid$, so it is enough to show that $g$ and $h$ are steady. For $g$, note that $(B, A^+)_\solid$ has the induced analytic ring structure from $(A, A^+)_\solid$ (by \cref{rslt:solid-of-map-A+-A-is-induced-ring-structure}) and hence $g$ is steady by \cref{rslt:nuclear-induced-ring-is-steady}, because $B$ is nuclear in $\D_\solid(A, A^+)$. To prove steadiness of $h$, by the ``2-out-of-3-property'' \cref{rslt:steadyness-satisfies-2-out-of-3} we can assume $A^+ = \Z$, i.e. we want to show that the map $(B, \Z)_\solid \to (B, B^+)_\solid$ is steady. Writing $B^+$ as a filtered colimit of finite-type $\Z$-algebras we reduce to the case that $B^+$ is of finite type over $\Z$. Since $(B, B^+)_\solid$ only depends on the image of $B^+$ in $\pi_0 B$ we can assume that $B^+ = \Z[T_1, \dots, T_n]$ for some $n$. Now $h$ is the base-change of the map $(B^+, \Z)_\solid \to (B^+, B^+)_\solid$ along $(B^+, \Z)_\solid \to (B, \Z)_\solid$ (by \cref{rslt:colimits-of-solid-discrete-rings}), so we reduce to showing that the map
\begin{align*}
	(\Z[T_1, \dots, T_n], \Z)_\solid \to (\Z[T_1, \dots, T_n], \Z[T_1, \dots, T_n])_\solid
\end{align*}
is steady. By induction (using preservation of steadiness under base-change and composition) this reduces to the case $n = 1$. In this case the steadiness follows easily from an explicit computation of $- \tensor_{(\Z[T], \Z)_\solid} \Z[T]_\solid$ in terms of a $\IHom$ along the lines of the proof of \cref{rslt:finite-type-solid-pullback-preserves-limits} (see \cite[Example 13.15.(2)]{scholze-analytic-spaces}).
\end{proof}

\begin{corollary} \label{rslt:solid-base-change-for-discrete-rings}
Let $A' \from A \to B$ be a diagram of discrete rings and denote $B' := A' \tensor_A B$. Then for every $M \in \D_\solid(B)$ the natural map
\begin{align*}
	M \tensor_{A_\solid} A'_\solid \isoto M \tensor_{B_\solid} B'_\solid
\end{align*}
is an isomorphism. In other words, base change along the following diagram holds:
\begin{center}\begin{tikzcd}
	\D_\solid(B') \arrow[d] && \arrow[ll,swap,"- \tensor_{B_\solid} B'_\solid"] \arrow[d] \D_\solid(B)\\
	\D_\solid(A') && \arrow[ll,swap,"- \tensor_{A_\solid} A'_\solid"] \D_\solid(A)
\end{tikzcd}\end{center}
\end{corollary}
\begin{proof}
By \cref{rslt:colimits-of-solid-discrete-rings} we have $B'_\solid = A'_\solid \tensor_{A_\solid} B_\solid$ and by \cref{rslt:map-of-discrete-rings-is-steady} the map $A_\solid \to A'_\solid$ is steady. Hence the claim follows from base-change along steady maps, see \cref{rslt:steady-map-of-analytic-rings-equiv-base-change}.
\end{proof}

We now introduce an almost version of the above definitions and results by reducing them to the non-almost case over $\Z$.

\begin{definition} \label{def:solid-almost-ring}
Let $(V,\mm)$ be an almost setup.
\begin{defenum}
	\item A \emph{discrete Huber pair over $(V,\mm)$} is a pair $(A, A^+)$ where $A$ is a discrete $(V,\mm)$-algebra and $A^+$ is an integrally closed subring of $\pi_0 A_{**\omega}$ (where $(-)_\omega$ denotes the discretization functor). A morphism $(A, A^+) \to (B, B^+)$ of discrete Huber pairs over $(V,\mm)$ is a morphism $A \to B$ of $(V,\mm)$-algebras such that the image of $A^+$ under the induced map $\pi_0 A_{**\omega} \to \pi_0 B_{**\omega}$ lies in $B^+$.

	\item Let $(A, A^+)$ be a discrete Huber pair over $(V,\mm)$. We let
	\begin{align*}
		(A, A^+)_\solid \in \AnRing_{(V,\mm)}
	\end{align*}
	denote the following analytic ring structure on $A$: The pair $(A_{**\omega}, A^+)$ is a discrete Huber pair over $V$, and we let $(A, A^+)_\solid := (A_{**\omega}, A^+)^a$.

	\item \label{def:abberviated-version-of-solid-almost-rings} Given a discrete $(V,\mm)$-algebra $A$, a discrete ring $A^+$ and a map $A^+ \to \pi_0 A_*$ we also denote $(A, A^+)_\solid := (A, \tilde A^+)_\solid$, where $\tilde A^+ \subset \pi_0 A_{**\omega}$ is the integral closure of the image of the map $\pi_0 A^+ \to \pi_0 A_{**\omega}$.

	\item We abbreviate $\D_\solid(A, A^+) = \D((A, A^+)_\solid)$. The objects of $\D_\solid(A, A^+)$ are called the \emph{solid $(A, A^+)$-modules}.
\end{defenum}
\end{definition}

\begin{proposition}
Let $(V,\mm)$ be an almost setup.
\begin{propenum}
	\item The assignment $(A, A^+) \mapsto (A, A^+)_\solid$ defines a fully faithful functor from the $\infty$-category of discrete Huber pairs over $(V,\mm)$ to $\AnRing_{(V,\mm)}$.

	\item \label{rslt:colimits-of-solid-discrete-almost-rings} The functor in (i) preserves all non-empty colimits. In particular, for any diagram $(B, B^+) \from (A, A^+) \to (C, C^+)$ of discrete Huber pairs over $(V,\mm)$ we have
	\begin{align*}
		(B, B^+)_\solid \tensor_{(A, A^+)_\solid} (C, C^+)_\solid = (B \tensor_A C, B^+ \tensor_{A^+} C^+)_\solid,
	\end{align*}
	where on the right-hand side we use the abbreviated notation from \cref{def:abberviated-version-of-solid-almost-rings}.

	\item \label{rslt:almost-solid-discrete-equiv-nuclear} Let $(A, A^+)$ be a discrete Huber pair over $(V,\mm)$. Then an object $P \in \D_\solid(A, A^+)$ is nuclear if and only if it is discrete.

	\item \label{rslt:map-of-discrete-almost-rings-is-steady} Let $(A, A^+)$ and $(B, B^+)$ be discrete Huber pairs over $(V,\mm)$. Then every morphism $f\colon (A, A^+)_\solid \to (B, B^+)_\solid$ is steady.
\end{propenum}
\end{proposition}
\begin{proof}
The functor $(A, A^+) \mapsto (A_{**\omega}, A^+)$ defines a fully faithful functor from discrete Huber pairs over $(V,\mm)$ to discrete Huber pairs over $V$. One easily deduces (i) and (ii) from \cref{rslt:solidification-is-fully-faithful-functor,rslt:colimits-of-solid-discrete-rings}. Part (iii) can be proved in exactly the same way as \cref{rslt:solid-discrete-equiv-nuclear}. Part (iv) can be deduced from \cref{rslt:map-of-discrete-rings-is-steady,rslt:almost-localization-preserves-steadiness}.
\end{proof}

Before we continue, we need to deal with some set-theoretic issues. In general, if $\mathcal A$ is an analytic ring then we defined the presentable $\infty$-categories $\D(\mathcal A)_\kappa$ for large enough strong limit cardinals $\kappa$ such that $\D(\mathcal A) = \bigunion_\kappa \D(\mathcal A)_\kappa$. By definition, $\D(\mathcal A)_\kappa$ is generated by the objects $\mathcal A[S]$ for $\kappa$-small extremally disconnected sets $S$. However, it may in general be false that $\mathcal A[S]$ is itself $\kappa$-condensed and therefore $\D(\mathcal A)_\kappa$ may not be a subcategory of $\D(\underline{\mathcal A})_\kappa$. In particular, for a map $f\colon \mathcal A \to \mathcal B$ of analytic rings, the ``forgetful functor'' $f_{*\kappa}\colon \D(\mathcal B)_\kappa \to \D(\mathcal A)_\kappa$ is in general the composition of the actual forgetful functor $f_*\colon \D(\mathcal B)_\kappa \to \D(\mathcal A)$ and the restriction $\D(\mathcal A) \to \D(\mathcal A)_\kappa$. Therefore $f_{*\kappa}$ may not be conservative in general, even though $f_*$ is. One can often circumvent this issue by juggling strong limit cardinals, but it turns out that in the setting of solid analytic rings the whole situation becomes much easier:

\begin{definition} \label{def:solid-cutoff-cardinal}
A \emph{solid cutoff cardinal} is a strong limit cardinal $\kappa$ such that either $\kappa = \omega$ or $\kappa$ satisfies the following property: For all cardinals $\lambda < \kappa$ there is a strong limit cardinal $\kappa_\lambda < \kappa$ such that the cofinality of $\kappa_\lambda$ is larger than $\lambda$.
\end{definition}

\begin{lemma} \label{rslt:properties-of-solid-cutoff-cardinals}
\begin{lemenum}
	\item There is a cofinal class of solid cutoff cardinals.

	\item Let $\kappa$ be a solid cutoff cardinal and $(V,\mm)$ an almost setup. Then for every discrete Huber pair $(A, A^+)$ over $(V,\mm)$ we have $\D_\solid(A, A^+)_\kappa \subset \D(A)_\kappa$.

	In particular, for every map $f\colon (A, A^+) \to (B, B^+)$ of discrete Huber pairs over $(V,\mm)$, the forgetful functor $f_*\colon \D_\solid(B, B^+) \to \D_\solid(A, A^+)$ restricts to a forgetful functor $f_*\colon \D_\solid(B, B^+)_\kappa \to \D_\solid(A, A^+)_\kappa$.
\end{lemenum}
\end{lemma}
\begin{proof}
Part (i) is \cite[Lemma 4.1]{etale-cohomology-of-diamonds}. We now prove (ii), so let $\kappa$ and $(A, A^+)$ be given. If $\kappa = \omega$ then the claim is obvious, so from now on we assume that $\kappa$ is uncountable. Fix a $\kappa$-small extremally disconnected set $S$. We need to show that $(A, A^+)_\solid[S]$ is $\kappa$-condensed, i.e. lies in $\D(V,\mm)_\kappa$. Using that the almost localization functor $(V,\mm)^*\colon \D(V) \to \D(V,\mm)$ maps $\D(V)_\kappa$ to $\D(V,\mm)_\kappa$ we can reduce to the case that $(A, A^+)$ is a discrete Huber pair over $V$; we can further assume $V = \Z$. Since $\D(\Z)_\kappa$ is stable under all small colimits we can reduce to the case that $A = A^+$ and that this ring is of finite type over $\Z$. Then by definition we have $(A, A^+)_\solid[S] = A_\solid[S] = \varprojlim_i A[S_i]$, where $S = \varprojlim_{i\in I} S_i$ with all $S_i$ finite. We can arrange that $I$ is $\kappa$-small (e.g. we can choose $I$ to be the set of disjoint clopen covers of $S$, so that $\abs I \le \omega^{\kappa'} \le 2^{\omega \cprod \kappa'} < \kappa$, where $\kappa' < \kappa$ is the cardinality of $S$). Let $\lambda < \kappa$ denote the cardinality of $I$. By definition of solid cutoff cardinals, we can find a strong limit cardinal $\kappa_\lambda < \kappa$ such that the cofinality of $\kappa_\lambda$ is larger than $\lambda$. By (the proof of) \cref{rslt:condensed-objects-in-presentable-monoidal-cat} it follows that $\D(\Z)_{\kappa_\lambda}$ is stable under $\lambda$-small limits in $\D(\Z)$, which immediately implies that $\varprojlim_i A[S_i] \in \D(\Z)_{\kappa_\lambda} \subset \D(\Z)_\kappa$, as desired.
\end{proof}

We now glue the solid analytic rings studied above in order to obtain schemes and discrete adic spaces. Using geometry blueprints (see \cref{def:geometry-blueprint}) the gluing process is completely formal:

\begin{definition}
Let $(V,\mm)$ be an almost setup.
\begin{defenum}
	\item We let $G_{\omega-adic}$ denote the geometry blueprint over $\Z$ whose $G_{\omega-adic}$-analytic rings are precisely the analytic rings of the form $(A, A^+)_\solid$ for discrete Huber pairs $(A, A^+)$ and whose $G_{\omega-adic}$-localizations are generated by the two maps
	\begin{align*}
		(\Z[T], \Z)_\solid \to (\Z[T,T^{-1}], \Z)_\solid, \qquad (\Z[T], \Z)_\solid \to \Z[T]_\solid.
	\end{align*}
	We let $G^{(V,\mm)}_{\omega-adic}$ denote the geometry blueprint which is generated by the image of $G_{\omega-adic}$ under the almost localization $\AnRing_V \to \AnRing_{(V,\mm)}$.

	\item A \emph{discrete adic space over $(V,\mm)$} is a $G_{\omega-adic}^{(V,\mm)}$-analytic space. For every discrete Huber pair $(A, A^+)$ over $(V,\mm)$ we denote
	\begin{align*}
		\Spa(A, A^+) := \AnSpec_{G_{\omega-adic}^{(V,\mm)}} (A, A^+)_\solid.
	\end{align*}
	An \emph{open immersion} $U \injto X$ of discrete adic spaces over $(V,\mm)$ is a map which exhibits $U$ as a $G_{\omega-adic}^{(V,\mm)}$-analytic subspace of $X$. A \emph{discrete adic space} is a discrete adic space over $\Z$.

	\item A discrete adic space is called \emph{classical} if it can be covered by open sets of the form $\Spa(A, A^+)$ with $A$ static.
\end{defenum}
\end{definition}

\begin{definition}
\begin{defenum}
	\item We let $G_{sch}$ denote the geometry blueprint over $\Z$ whose $G_{sch}$-analytic rings are precisely the analytic rings of the form $A_\solid$ for discrete rings $A$ and whose $G_{sch}$-localizations are generated by the map $\Z[T]_\solid \to \Z[T, T^{-1}]_\solid$.

	\item A \emph{scheme} is a $G_{sch}$-analytic space. For every discrete ring $A$ we denote
	\begin{align*}
		\Spec A := \AnSpec_{G_{sch}} A_\solid.
	\end{align*}
	An \emph{open immersion} $U \injto X$ of schemes is a map which exhibits $U$ as a $G_{sch}$-analytic subspace of $X$.

	\item A scheme is called \emph{classical} if it can be covered by open sets of the form $\Spec A$ with $A$ static.
\end{defenum}
\end{definition}

\begin{remarks}
\begin{remarksenum}
	\item By \cref{rslt:solidification-is-fully-faithful-functor} a scheme (in the above sense) can be identified with a functor $\Alg_{\ge0}(\Z)_\omega \to \Ani$ which satisfies descent along open covers of affine schemes. It follows that the above defined $\infty$-category of schemes is equivalent to the $\infty$-category of derived schemes as mentioned in \cite[\S25]{lurie-spectral-algebraic-geometry} via the Yoneda embedding (use also \cref{rslt:sheaves-on-basis-equiv-sheaves-on-whole-site} to see that Zariski sheaves on the site of all derived schemes are the same as Zariski sheaves on affine derived schemes).

	\item There is an analytification functor $X \mapsto X^\an$ (see \cref{def:analytification}) associating to every discrete adic space over some almost setup $(V,\mm)$ an analytic space over $(V,\mm)$. By an unpublished argument of Scholze the restriction of this functor to schemes is fully faithful.

	\item There is a natural functor from schemes to discrete adic spaces. This functor can be shown to be fully faithful by the same argument as for the analytificiation functor.

	\item For every almost setup $(V,\mm)$ there is an almost localization functor $X \mapsto X^a$ associating to every discrete adic space $X$ over $V$ a discrete adic space $X^a$ over $(V,\mm)$ (see \cref{def:functoriality-of-analytic-spaces-and-sheaves}). It maps $\Spa(A, A^+) \mapsto \Spa(A^a, A^+)$. This allows to deduce many of the following results (which are only stated over $\Z$) also in the almost setting.
\end{remarksenum}
\end{remarks}

\begin{remark} \label{rslt:morphism-of-discrete-adic-spaces-is-steady}
It follows from \cref{rslt:map-of-discrete-almost-rings-is-steady} that every morphism of discrete adic spaces over any almost setup $(V,\mm)$ is steady.
\end{remark}

The above defined $\infty$-category of discrete adic spaces is a generalization of the category of classical discrete adic spaces where we additionally allow non-static structure sheaves. More concretely we have the following comparison:

\begin{lemma} \label{rslt:classical-disc-adic-spaces-embedd-into-adic-spaces}
There is a natural fully faithful functor from the category of classical discrete adic spaces to the $\infty$-category of discrete adic spaces. This functor has a right adjoint mapping $\Spa(A, A^+) \mapsto \Spa(\pi_0 A, A^+)$. The right adjoint preserves fiber products and open immersions.
\end{lemma}
\begin{proof}
The existence of the desired fully faithful functor essentially boils down to the fact that if $\Spa(B, B^+) \injto \Spa(A, A^+)$ is an open immersion of discrete adic spaces and $\Spa(A, A^+)$ is classical, then $\Spa(B, B^+)$ is classical, which is easy to see. The right adjoint is easily constructed (on the level of the surrounding $\infty$-categories of sheaves it is just a pushforward, which in this case is just a restriction) and has the desired properties (note that preservation of pullbacks is true for any right adjoint).
\end{proof}

\begin{definition}
Let $X$ be a discrete adic space. We call the associated classical discrete adic space from \cref{rslt:classical-disc-adic-spaces-embedd-into-adic-spaces} the \emph{underlying classical discrete adic space $X_0$} of $X$.\end{definition}

Next we want to introduce the underlying topological space $\abs X$ of a discrete adic space $X$. This can be done directly using the $\infty$-category of discrete adic spaces as defined above, e.g. the points of $\abs X$ are precisely the equivalence classes of maps $\Spa(K, K^+) \to X$, where $K$ is a field and $K^+ \subset K$ is a valuation ring with fraction field $K$. However, we take a shortcut by referring to the known topological space underlying a classical discrete adic space. To make this work, we need the following result:

\begin{lemma} \label{rslt:topological-cover-equiv-disc-adic-cover}
A finite family $((A, A^+)_\solid \to (B_i, B_i^+)_\solid)_i$ of $G_{\omega-adic}$-localizations is a $G_{\omega-adic}$-analytic cover if and only if the family $\abs{\Spa(\pi_0 B_i, B_i^+)} \to \abs{\Spa(\pi_0 A, A^+)}$ of open immersions forms a topological cover.
\end{lemma}
\begin{proof}
We first prove the ``if'' direction, so assume that the given maps form a topological cover. We need to show that the family of functors $- \tensor_{(A, A^+)_\solid} (B_i, B_i^+)_\solid$ is conservative. Note that each of the maps $(A, A^+)_\solid \to (B_i, B_i^+)_\solid$ identifies $(B_i, B_i^+)_\solid$ with an $(A, A^+)_\solid$-algebra of the form
\begin{align*}
	(A[1/f], A^+[g_1/f, \dots, g_n/f])_\solid
\end{align*}
for some $f, g_1, \dots, g_n \in A$ (one checks easily that the collection of such maps is stable under base-change and compositions and hence provides all $G_{\omega-adic}$-localizations). In other words, each of the maps $(A, A^+)_\solid \to (B_i, B_i^+)_\solid$ is given by a rational open subset of $\abs{\Spa(\pi_0 A, A^+)}$ and by assumption these rational open subsets form a cover of $\abs{\Spa(\pi_0 A, A^+)}$. We can now argue as in the proof of \cite[Lemma 10.3.(ii)]{condensed-mathematics}; the argument roughly goes as follows: One first refines the given topological cover by a topolocial cover of the form $(A[1/s_i], A^+[s_1/s_i, \dots, s_N/s_i])_\solid$ for certain $s_1, \dots, s_N \in \pi_0 A$ generating the unit ideal. Then one can argue by induction on $N$ which reduces the claim to the case of a topological cover by $(A, A^+[f])_\solid$ and $(A[1/f], A^+[1/f])_\solid$ for some $f \in \pi_0 A$. This cover is the base-change of the standard topological cover of $(\Z[T], \Z)_\solid$ given by $\Z[T]_\solid$ and $(\Z[T, T^{-1}], \Z[T^{-1}])_\solid$. Then a quick computation shows that this last topological cover is indeed a $G_{\omega-adic}$-analytic cover, which proves the claim.

We now prove the ``only if'' direction. Suppose that the given maps form a $G_{\omega-adic}$-analytic cover but not a topological cover. Then one can find a valuation ring $K^+$ with fraction field $K$ and a map $(A, A^+)_\solid \to (K, K^+)_\solid$ which does not factor over any $(B_i, B_i^+)_\solid$. In particular the maps $(K, K^+)_\solid \to (A, A^+)_\solid \tensor_{(K, K^+)_\solid} (B_i, B_i^+)_\solid$ do not form a topological cover. We can thus reduce to the case that $(A, A^+) = (K, K^+)$. In this case every $(B_i, B_i^+)_\solid$ is of the form $(K, K^+[f_1, \dots, f_n])_\solid$ and two such analytic rings agree precisely if the $f_1, \dots, f_n$ generate the same valuation subring of $K$. Since the valuation rings between $K$ and $K^+$ form a chain under inclusion, we reduce to the case of a single map $\alpha\colon (K, K^+)_\solid \to (K, K'^+)_\solid$ which is a $G_{\omega-adic}$-analytic cover, where $K'^+ \supset K^+$ is another valuation ring with fraction field $K$. But then $\alpha$ is a localization and base-change along $\alpha$ is an equivalence, which shows that $\alpha$ is an isomorphism. This implies $K'^+ = K^+$ (using the full faithfulness of the functor $(A, A^+) \mapsto (A, A^+)_\solid$ from \cref{rslt:solidification-is-fully-faithful-functor} which we implicitly used throughout the whole proof). Contradiction!
\end{proof}

\begin{definition}
For every discrete adic space $X$ we denote $\abs X$ the topological space glued from the topological spaces $\abs{\Spa(\pi_0 A, A^+)}$ on an open cover of $X$ by affinoid subspaces $\Spa(A, A^+)$. This topological space is well-defined by \cref{rslt:topological-cover-equiv-disc-adic-cover} and depends only on the underlying classical discrete adic space $X_0$.
\end{definition}

\begin{remark}
Even if one is only interested in classical discrete adic spaces one still needs to show \cref{rslt:topological-cover-equiv-disc-adic-cover} (at least for classical discrete Huber rings) to verify that our definition of classical discrete adic spaces coincides with the definition by Huber.
\end{remark}

With the underlying topological space of a discrete adic space at hand, one could now attempt to characterize discrete adic spaces as certain locally ringed spaces. We do not pursue this idea any further (but cf. \cite[Proposition 14.2]{scholze-analytic-spaces}). However, let us prove the following handy tools motivated by the picture of locally ringed spaces:

\begin{lemma} \label{rslt:disc-adic-space-top-subsets-equiv-open-subspaces}
Let $X$ be a discrete adic space. Then the assignment $U \mapsto \abs U$ defines a bijection between open subspaces $U \subset X$ and open subsets of $\abs X$.
\end{lemma}
\begin{proof}
This can be checked locally on $X$ so we can assume that $X = \Spa(A, A^+)$ is affine. Now every open subset of $\abs X$ can be covered by rational open subsets and hence comes from some open subspace $U \subset X$. Moreover, if $U, U' \subset X$ are two open subspaces such that $\abs U = \abs U'$ then we must have $U = U'$ by \cref{rslt:topological-cover-equiv-disc-adic-cover}.
\end{proof}

\begin{corollary} \label{rslt:disc-adic-space-map-factors-over-U-if-over-top-U}
Let $f\colon Y \to X$ be a map of discrete adic spaces and let $U \subset X$ be an open subspace such that $\abs f\colon \abs Y \to \abs X$ factors over $\abs U$. Then $f$ factors over $U$.
\end{corollary}
\begin{proof}
To show that $f$ factors over $U$, we can equivalently show that the map $U \cprod_X Y \injto Y$ is an isomorphism. By assumption on $f$ we have $\abs{U \cprod_X Y} = \abs Y$, hence the claim follows from \cref{rslt:disc-adic-space-top-subsets-equiv-open-subspaces}.
\end{proof}

\begin{lemma} \label{rslt:locally-mono-plus-inj-on-top-implies-globally-mono}
Let $f\colon Y \to X$ be a map of discrete adic spaces which satisfies the following properties:
\begin{enumerate}[(a)]
	\item There is a cover of $Y$ by open subspaces $V \subset Y$ such that each map $V \to X$ is a monomorphism.

	\item The map $\abs Y \injto \abs X$ is injective.
\end{enumerate}
Then $f$ is a monomorphism.
\end{lemma}
\begin{proof}
Let $Z = \Spa(A, A^+)$ be an affine discrete adic space. In order to show that $f$ is a monomorphism we need to see that the induced map $\Hom(Z, X) \to \Hom(Z, Y)$ is a monomorphism of anima. To see this, we can fix a map $h\colon Z \to X$ and now need to show that the anima $\Hom_X(Z, Y)$ is either empty or contractible. Assume that it is non-empty, i.e. that there is a map $g\colon Z \to X$ such that $h = f \comp g$. Using (a) we pick a cover $Y = \bigunion_i V_i$ of open subspaces $V_i \subset Y$ such that each map $V_i \to X$ is a monomorphism. Let $W_i := g^{-1}(V_i) \subset Z$, let $W' := \bigdunion_i W_i$ and let $W'_\bullet \to Z$ be the associated Čech cover. Then $\Hom_X(Z, Y) = \varprojlim_{n\in\Delta} \Hom_X(W'_n, Y)$, so it is enough to show that each $\Hom_X(W'_n, Y)$ is contractible. Each $W'_n$ is a disjoint union of finite intersections of $W_i$'s, so we are reduced to showing the following: Let $V \subset Y$ be an open subset such that the map $V \to X$ is a monomorphism and let $W := g^{-1}(V) \subset Z$; then $\Hom_X(W, Y)$ is contractible. We know that $\Hom_X(W, V)$ is contractible, so we only need to show that $\Hom_X(W, V) \isoto \Hom_X(W, Y)$ is an isomorphism; since this map is a monomorphism this boils down to showing that every $X$-morphism $W \to Y$ factors over $V$. By (b) it factors over $\abs V$, hence the claim follows from \cref{rslt:disc-adic-space-map-factors-over-U-if-over-top-U}.
\end{proof}

One big advantage of discrete adic spaces over schemes is that discrete adic spaces allow canonical compactifications. This in particular provides a useful notion of ``affine proper maps'', so that often questions about proper maps can be checked locally on the source. The canonical compactification can be constructed as follows:

\begin{lemma} \label{rslt:compactification-of-disc-adic-space}
Let $X$ be a discrete adic space. Then there is a discrete adic space $\overline X$ which is uniquely characterized by the property
\begin{align*}
	\Hom(\Spa(A, A^+), \overline X) = \Hom(\Spa(A, A), X)
\end{align*}
for all discrete Huber pairs $(A, A^+)$. The functor $X \mapsto \overline X$ preserves finite limits, monomorphisms, disjoint unions and covers.
\end{lemma}
\begin{proof}
Viewing $X$ as a sheaf on the analytic site of $G_{\omega-adic}$-analytic rings, the given $\Hom$-identity immediately defines $\overline X$ as a presheaf on that site. It follows from \cref{rslt:topological-cover-equiv-disc-adic-cover} that $\overline X$ satisfies descent along analytic covers: One needs to check that if $(\Spa(B_i, B_i^+) \to \Spa(A, A^+))_i$ is a cover of classical discrete adic spaces by rational open subsets then $(\Spa(B_i, B_i) \to \Spa(A, A))_i$ is still a cover; this is an easy exercise. Hence $\overline X$ is a sheaf on the analytic site of $G_{\omega-adic}$-analytic rings. It is clear that the thus defined functor $X \mapsto \overline X$ (seen as a functor of sheaves on the site of $G_{\omega-adic}$-analytic rings; we do not yet know that $\overline X$ is a discrete adic space) preserves finite limits and monomorphisms. To check that it preserves covers we employ \cref{rslt:cover-of-sheaves-characterized-by-sections}: Suppose we are given a cover $Y \to X$ of discrete adic spaces and a section $s\colon \Spa(C, C^+) \to \overline X$ for some discrete Huber pair $(C, C^+)$. Then $s$ corresponds to a map $s'\colon \Spa(C, C) \to X$ and thus there are a $G_{\omega-adic}$-analytic cover $((C, C) \to (C_i, C_i))_i$ (by the same argument as above) and sections $t_i'\colon \Spa(C_i, C_i) \to Y$ compatible with $s$. Then each $C_i$ is of the form $C[1/f_i]$ for some $f_i \in \pi_0 C$, so we can form the $G_{\omega-adic}$-analytic cover $((C, C^+) \to (C_i, C_i^+))_i$, where $C_i^+ = C^+[1/f_i]$. Then the sections $t'_i$ correspond to sections $t_i\colon \Spa(C_i, C_i^+) \to \overline Y$, compatible with $s$. This shows that $\overline Y \to \overline X$ is a cover. One can similarly show that the functor $X \mapsto \overline X$ preserves disjoint unions.

It remains to show that $\overline X$ is a discrete adic space over $X$. If $X = \Spa(B, B^+)$ is affine then one checks immediately that $\overline X = \Spa(B, \Z)$, which is indeed a discrete adic space. If $(B, B^+) \to (B', B'^+)$ is a $G_{\omega-adic}$-localization then so is $(B, \Z) \to (B', \Z)$ (e.g. one checks that the class of $G_{\omega-adic}$-localizations which satisfies this property after any base-change is stable under base-change and composition and contains the generating $G_{\omega-adic}$-localizations). Altogether this implies that if $U \subset X$ is an open subspace of the affine discrete adic space $X$ then $\overline U \subset \overline X$ is an open subspace of $\overline X$. We conclude by \cref{rslt:obtain-analytic-space-from-gluing} that $\overline X$ is a discrete adic space for general $X$.
\end{proof}

\begin{definition}
\begin{defenum}
	\item Let $X$ be a discrete adic space. We call the discrete adic space $\overline X$ from \cref{rslt:compactification-of-disc-adic-space} the \emph{compactification of $X$}. There is a canonical map $X \to \overline X$.

	\item Let $f\colon Y \to X$ be a map of discrete adic spaces. We denote
	\begin{align*}
		\overline f^{/X}\colon \overline Y^{/X} := \overline Y \cprod_{\overline X} X \to X
	\end{align*}
	and call it the \emph{compactification of $f$}. We also call $\overline Y^{/X}$ the \emph{relative compactification of $Y$ over $X$}. There is a natural map $Y \to \overline Y^{/X}$ of discrete adic spaces.
\end{defenum}
\end{definition}

\begin{example}
If $Y = \Spa(B, B^+)$ and $X = \Spa(A, A^+)$ are affine then $\overline Y^{/X} = \Spa(B, A^+)$.
\end{example}

In general we cannot say much about the map $Y \to \overline Y^{/X}$, but properties of this map have a strong connection to the following classical properties of $f$:

\begin{definition} \label{rslt:classical-properties-of-maps-of-disc-adic-spaces}
Let $f\colon Y \to X$ be a map of discrete adic spaces.
\begin{defenum}
	\item We say that $f$ is \emph{locally of $+$-finite type} if for all $y \in \abs Y$ there are affine open neighbourhoods $V = \Spa(B, B^+) \subset Y$ of $y$ and $U = \Spa(A, A^+) \subset X$ of $f(x)$ such that $f(V) \subset U$ and there exists a collection of finitely many elements $g_1, \dots, g_n \in B^+$ with $(B, B^+)_\solid = (B, A^+[g_1, \dots, g_n])_\solid$. We say that $f$ is \emph{of $+$-finite type} if $f$ is locally of $+$-finite type and quasicompact.

	\item We say that $f$ is \emph{separated} if for all valuation rings $K^+$ with fraction field $K$ and every commutative solid diagram
	\begin{center}\begin{tikzcd}
		\Spa(K, K) \arrow[r] \arrow[d] & Y \arrow[d]\\
		\Spa(K, K^+) \arrow[r] \arrow[ur,dashed] & X
	\end{tikzcd}\end{center}
	there is a dashed morphism making the diagram commute.

	\item We say that $f$ is \emph{partially proper} if it is separated, locally of $+$-finite type and for every diagram as in (ii), the dashed morphism is unique. We say that $f$ is \emph{proper} if $f$ is partially proper and quasicompact.
\end{defenum}
\end{definition}

Note that all of the properties in \cref{rslt:classical-properties-of-maps-of-disc-adic-spaces} only depend on the underlying map $f_0\colon Y_0 \to X_0$ of classical discrete adic spaces. In particular if $X$ and $Y$ are schemes then the above notion of separatedness agrees with its counterpart in \cite{lurie-spectral-algebraic-geometry}, while our notion of properness is slightly more general, as we allow maps which need not be of finite type (but are so up to an integral extension); see \cite[Remark 5.0.0.1]{lurie-spectral-algebraic-geometry}. It is easy to see that every map of affine spaces is separated and every separated map is quasiseparated. The notion ``of $+$-finite type'' seems to be new, but is motivated by the next results.

\begin{proposition} \label{rslt:characterization-of-compactification}
Let $f\colon Y \to X$ be a map of discrete adic spaces.
\begin{propenum}
	\item $f$ is separated if and only if the map $Y \to \overline Y^{/X}$ is a monomorphism.

	\item $f$ is locally of $+$-finite type if and only if the map $Y \to \overline Y^{/X}$ is locally on $Y$ an open immersion.

	\item \label{rslt:partially-proper-iff-equiv-to-compactification} $f$ is partially proper if and only if the map $Y \to \overline Y^{/X}$ is an isomorphism.
\end{propenum}
\end{proposition}
\begin{proof}
We first prove (i). The map $j\colon Y \to \overline Y^{/X}$ is always locally a monomorphism. Namely, this can be checked in the case that $X = \Spa(A, A^+)$ and $Y = \Spa(B, B^+)$ are affine, in which case $j$ becomes the map $\Spa(B, B^+) \to \Spa(B, A^+)$ which is a monomorphism. It follows directly from the definitions that $f$ is separated if and only if the map $\abs Y \to \abs{\overline Y^{/X}}$ is bijective (use that the points of these topological spaces are precisely the maps from $\Spa(K, K^+)$'s). Hence (i) reduces to \cref{rslt:locally-mono-plus-inj-on-top-implies-globally-mono}.

Part (ii) can be checked locally on $X$ and $Y$, so we can assume that $X = \Spa(A, A^+)$ and $Y = \Spa(B, B^+)$ are affine. Then $j$ is the map $\Spa(B, B^+) \to \Spa(B, A^+)$. If $\Spa(B, B^+) = \Spa(B, A^+[g_1, \dots, g_n])$ for some $g_1, \dots, g_n \in B^+$ then clearly $j$ is an open immersion. Conversely, if $j$ is an open immersion then clearly $j$ is of $+$-finite type (cover it by rational open subsets), which immediately implies that $f$ is of $+$-finite type.

It remains to prove (iii). By (i) and (ii), $j$ is an isomorphism if and only if $f$ is separated, locally of $+$-finite type and the map $\abs j\colon \abs Y \to \abs{\overline Y^{/X}}$ is a bijection (see \cref{rslt:disc-adic-space-top-subsets-equiv-open-subspaces}). This is exactly the definition of partial properness for $f$.
\end{proof}

\begin{proposition} \label{rslt:compactification-is-proper}
Let $f\colon Y \to X$ be a map of discrete adic spaces and let $\overline f^{/X}\colon \overline Y^{/X} \to X$ be its compactification. Then $\overline f^{/X}$ is partially proper. If $f$ is quasicompact then $\overline f^{/X}$ is proper.
\end{proposition}
\begin{proof}
The partial properness of $\overline f^{/X}$ follows immediately from its construction. If $f$ is quasicompact then after base-change to any affine open subspace $\Spa(A, A^+) \subset X$, $Y$ can be covered by finitely many affine open subspaces $\Spa(B_i, B_i^+)$. Then $\overline Y^{/X}$ can be covered by the affine open subspaces $\Spa(B_i, A^+)$ and thus is quasicompact. This shows that $\overline f^{/X}$ is quasicompact and thus proper.
\end{proof}

In the following we will often make implicit use of \cref{rslt:characterization-of-compactification,rslt:compactification-is-proper} when dealing with relative compactifications.

To every discrete adic space $X$ we have an associated $\infty$-category $\D_\solid(X)$ of solid quasicoherent sheaves on $X$ (see \cref{def:quasicoherent-sheaves-on-G-analytic-space}). If $X = \Spa(A, A^+)$ is affine then $\D_\solid(X) = \D_\solid(A, A^+)$. For every solid cutoff cardinal $\kappa$ we can define a similar $\infty$-category $\D_\solid(X)_\kappa$ glued from $\D_\solid(A, A^+)_\kappa$'s. One of the main results of \cite{condensed-mathematics} (and a big motivation for the theory of solid analytic rings) states that there is a full $6$-functor formalism for $\D_\solid(X)$ on schemes $X$. The following result makes this more precise.
\begin{proposition} \label{rslt:scheme-6-functor}
Let $f\colon Y \to X$ be a separated $+$-finite-type map of discrete adic spaces. Then there is a natural pair of adjoint functors
\begin{align*}
	f_!\colon \D_\solid(Y) \rightleftarrows \D_\solid(X) \noloc f^!
\end{align*}
satisfying the following properties:
\begin{propenum}
	\item (Functoriality) If $g\colon Z \to Y$ is another separated $+$-finite-type map of discrete adic spaces then $(f \comp g)_! \isom f_! \comp g_!$ and $(f \comp g)^! \isom g^! \comp f^!$.

	\item \label{rslt:scheme-6-functor-etale-and-proper-case} (Special Cases) If $f$ is proper then $f_! \isom f_*$. If $f$ is an étale map of classical schemes then $f^! \isom f^*$.

	\item \label{rslt:scheme-6-functor-projection-formula} (Projection Formula) For all $\mathcal M \in \D_\solid(X)$ and all $\mathcal N \in \D_\solid(Y)$ there is a natural isomorphism $f_!(\mathcal N \tensor f^* \mathcal M) \isom (f_! \mathcal N) \tensor \mathcal M$.

	\item \label{rslt:scheme-6-functor-proper-base-change} (Proper Base-Change) Let $g\colon X' \to X$ be any map of discrete adic spaces with base-change $Y' := Y \cprod_X X'$ and projections $f'\colon Y' \to X'$ and $g'\colon Y' \to Y$. Then there is a natural equivalence $g^* f_! \isom f'_! g'^*$ of functors $\D_\solid(Y) \to \D_\solid(X')$.

	\item \label{rslt:scheme-6-functor-poincare-duality} (Poincaré Duality) If $f$ is a smooth map of classical schemes and has pure dimension $d$ then $f^! \isom f^* \tensor \omega_{Y/X}[d]$, where $\omega_{Y/X} := \bigwedge^d \Omega^1_{Y/X}$.

	\item \label{rslt:scheme-6-functor-cutoff-cardinals} For every solid cutoff cardinal (see \cref{def:solid-cutoff-cardinal}) $\kappa$ we have $f_!(\D_\solid(Y)_\kappa) \subset \D_\solid(X)_\kappa$.
\end{propenum}
\end{proposition}
\begin{proof}
We will first prove everything in the case that $X = \Spa(A, A^+)$ and $Y = \Spa(B, B^+)$ are affine and $\Spa(B, B^+) = \Spa(A, A^+[g_1, \dots, g_n])$ for some $g_1, \dots, g_n \in B^+$. We factor $f$ into the maps $j\colon Y \to \overline Y^{/X}$ and $\overline f^{/X}\colon \overline Y^{/X} \to X$, where $\overline Y^{/X} = \Spa(B, A^+)$ is the relative compactification of $f$. As in \cite[Definition 11.3]{condensed-mathematics} we will define $f_! := (\overline f^{/X})_* \comp j_!$, where $j_!\colon \D_\solid(B, B^+) \to \D_\solid(B, A^+)$ is a left adjoint of $j^*$. To make this construction work, we need to show the existence of $j_!$. Fix a solid cutoff cardinal $\kappa$. We claim that $j^*\colon \D_\solid(B, B^+)_\kappa \to \D_\solid(B, A^+)_\kappa$ preserves all small limits. Using the elements $g_1, \dots, g_n$ we see that the map $(B, A^+)_\solid \to (B, B^+)_\solid$ is a base-change of the map
\begin{align*}
	(\Z[T_1, \dots, T_n], \Z)_\solid \to (\Z[T_1, \dots, T_n], \Z[T_1, \dots, T_n])_\solid,
\end{align*}
so it suffices to prove that pullback along the latter map preserves all small limits. This readily reduces to showing that the pullback $j^*\colon \D_\solid(\Z[T], \Z)_\kappa \to \D_\solid(\Z[T])_\kappa$ preserves all small limits. This follows from \cite[Observation 8.11]{condensed-mathematics}, which explicitly constructs a left-adjoint. We have finished the proof that $j^*\colon \D_\solid(B, A^+)_\kappa \to \D_\solid(B, B^+)_\kappa$ preserves small limits. As both $\D_\solid(B, A^+)_\kappa$ and $\D_\solid(B, B^+)_\kappa$ are presentable, this formally implies the existence of a left-adjoint
\begin{align*}
	j_{!\kappa}\colon \D_\solid(B, B^+)_\kappa \to \D_\solid(B, A^+)_\kappa.
\end{align*}
It follows formally from steady base-change (see \cref{rslt:steady-base-change-holds}) by passing to right-adjoints that $j_{!\kappa}$ satisfies base-change with respect to any morphism $(C, C^+)_\solid \to (B, A^+)_\solid$ of discrete Huber pairs. Now pick another solid cutoff cardinal $\kappa' \ge \kappa$. Then for every $\kappa$-small profinite set $S$ we have $j_{!\kappa'} (B, B^+)_\solid[S] \in \D_\solid(B, A^+)_\kappa$. Namely, by base-change for $j_{!\kappa'}$ we can reduce to the case $B = B^+ = \Z[T]$ and $A^+ = \Z$, where we have by \cite[Theorem 8.1]{condensed-mathematics}
\begin{align*}
	j_{!\kappa'} \Z[T]_\solid[S] = \Z[T]_\solid[S] \tensor_{(\Z[T], \Z)_\solid} j_{!\kappa'} \Z[T],
\end{align*}
and $j_{!\kappa'} \Z[T]$ lies in $\D_\solid(\Z[T], \Z)_\kappa$ for every solid cutoff cardinal $\kappa$ by the explicit construction in \cite[Observation 8.11]{condensed-mathematics}. Since $j_{!\kappa'}$ preserves small colimits, we deduce $j_{!\kappa'}(\D_\solid(B, B^+)_\kappa) \subset \D_\solid(B, A^+)_\kappa$, i.e. the computation of $j_{!\kappa}$ does not depend on $\kappa$. This shows that the functors $j_{!\kappa}$ combine to a functor
\begin{align*}
	j_!\colon \D_\solid(B, B^+) \to \D_\solid(B, A^+)
\end{align*}
which is left-adjoint to $j^*$. This finishes the construction of $f_! := (\overline f^{/X})_* \comp j_!$ (in the affine case).

It is easy to see that $f_!$ admits a right adjoint $f^!$: This reduces to the case $f = \overline f^{/X}$, where $f_! = f_*$ is just a forgetful functor. This forgetful functor preserves all small colimits and $\D(-)_\kappa$ (by \cref{rslt:properties-of-solid-cutoff-cardinals}), so one can easily construct right adjoint $f^{!\kappa}$. Then note that $f_* f^{!\kappa}$ is computed by an $\IHom$ with fixed first argument and can thus be made independent of $\kappa$ (for large enough $\kappa$).

We now verify that the above defined functors $f_!$ and $f^!$ satisfy properties (i) -- (vi): Part (iv) follows from base-change for $j_!$ (as discussed above) and base-change for $(\overline f^{/X})_*$ (see \cref{rslt:solid-base-change-for-discrete-rings}). Part (iii) can be separately shown for $j_!$ and $\overline f^{/X}$; the latter case is easy. For the former case, i.e. the projection formula for $j_!$, note first that there is a natural map $j_!(\mathcal N \tensor j^* \mathcal M) \to j_! \mathcal N \tensor \mathcal M$. Both sides commute with colimits in $\mathcal M$ and $\mathcal N$, so we can reduce to the case $\mathcal M = (B, A^+)_\solid[S]$ and $\mathcal N = (B, B^+)_\solid[T]$ for some profinite sets $S$ and $T$. Now we can use base-change as above to reduce to the case $B = B^+ = \Z[T]$ and $A^+ = \Z$. In this case we can argue as in the proof of the projection formula in \cite[Theorem 8.2]{condensed-mathematics}: The two sides of the projection formula become isomorphic after applying $j^*$, hence their cofiber is a $B_\infty = \Z((T^{-1}))/B$-module. This means that tensoring the cofiber with $B_\infty$ leaves it invariant (see \cite[Observation 8.7]{condensed-mathematics}); on the other hand, tensoring any side of the projection formula with $B_\infty$ produces $0$ (in general, $j_! M \tensor_{(B, A)_\solid} B_\infty = 0$, which follows from the triangle $j_!M \to M \to M \tensor B_\infty$ and \cite[Observation 8.7]{condensed-mathematics}). This implies that the cofiber is $0$, hence the projection formula holds. Part (v) follows from \cite[Theorem 11.6]{condensed-mathematics}. The proper part of (ii) is obvious and the étale part of (ii) is a special case of (v).
Part (vi) can be shown separately for $j_!$ (where it is true by the above construction) and for $(\overline f^{/X})_*$ (where it follows from \cref{rslt:properties-of-solid-cutoff-cardinals}).

To finish the affine case, it remains to prove (i). Given $f$ and $g$ as in the claim, one easily reduces the claim to showing that the following diagram commutes (cf. the proof of \cite[Proposition 22.9]{etale-cohomology-of-diamonds}):
\begin{center}\begin{tikzcd}
	\D(\overline Z^{/Y}) \arrow[r,"j'_!"] \arrow[d,"h'_*"] & \D(\overline Z^{/X}) \arrow[d,"h_*"] \\
	\D(Y) \arrow[r,"j_!"] & \D(\overline Y^{/X})
\end{tikzcd}\end{center}
This follows formally from base-change and the projection formulas (cf. the proof of \cref{rslt:etale-lower-shriek-compatible-with-pushforward}). This finishes the proof of all claims in the affine case.

Now let $f$ be general. We construct $f_!$ as in the affine case: Factor $f$ as a composition of $j\colon Y \to \overline Y^{/X}$ and $\overline f^{/X}\colon \overline Y^{/X} \to X$ and set $f_! := (\overline f^{/X})_* \comp j_!$. For this to work, we need to check the existence of the left-adjoint $j_!\colon \D_\solid(Y) \to \D_\solid(\overline Y^{/X})$ of $j^*$. We know that $j_!$ exists in the affine case and satisfies arbitrary base-change (and in particular base-change along open immersions), so it can formally be glued to a functor $j_!$ for general $X$ and $Y$ (write the solid quasicoherent sheaf $\infty$-categories as limits over the solid module categories on an affine cover and use the explicit description of limits of $\infty$-covers in terms of coCartesian sections; then observe that the base-change for $j_!$ implies that the componentwise $j_!$ preserves coCartesian morphisms and thus induces a functor on the limits). The existence of $f^!$ follows in a similar way as in the affine case and the claims (i) -- (vi) can also easily be reduced to the affine case.
\end{proof}

\begin{remark}
One can try to improve upon \cref{rslt:scheme-6-functor} in several ways: Firstly one may ask if a similar statement is true in the almost setting (see \cref{rslt:lower-shriek-commutes-with-almost-localization} below), secondly one may make the $6$-functor formalism much more functorial (analogous to \cref{sec:ri-pi.6-functor}) and thirdly one may improve upon the Poincaré duality statement by allowing non-classical schemes and even discrete adic spaces.
\end{remark}

\begin{lemma} \label{rslt:lower-shriek-commutes-with-almost-localization}
Let $(V,\mm)$ be an almost setup and let $f\colon Y \to X$ be a separated $+$-finite type map of discrete adic spaces over $V$ with associated map $f^a\colon Y^a \to X^a$ of discrete adic spaces over $(V,\mm)$. Then there is a natural pair of adjoint functors
\begin{align*}
	f^a_!\colon \D_\solid(Y^a) \rightleftarrows \D_\solid(X^a) \noloc f^{a!}
\end{align*}
satisfying $f^a_! \comp (-)^a = (-)^a \comp f_!$ and $f^{a!} \comp (-)^a = (-)^a \comp f^!$.
\end{lemma}
\begin{proof}
As in the proof of \cref{rslt:scheme-6-functor} the construction of $f^a_!$ and $f^{a!}$ formally reduces to the case that $Y = \Spa(B, B^+)$ and $X = \Spa(A, A^+)$ are affine. Then $f^a_!$ can be constructed via the relative compactification of $f$, which reduces its construction to the case that $A = B$. By \cref{rslt:scheme-6-functor} the map $f^*\colon \D_\solid(B, A^+) \to \D_\solid(B, B^+)$ has a left adjoint $f_!$. Note that $f^{a*}$ is the composition $f^{a*} = (-)^a \comp f^* \comp (-)_*$ and all of the occurring functors on the right-hand side admit left adjoints; hence $f^{a*}$ admits a left adjoint $f^a_!$. We now prove that this functor $f^a_!$ satisfies $f^a_! \comp (-)^a = (-)^a \comp f_!$. Passing to right adjoints we need to see that for all $N \in \D^a_\solid(B, A^+)$ we have
\begin{align*}
	(N \tensor_{(B^a, A^+)_\solid} (B^a, B^+)_\solid)_* = N_* \tensor_{(B, A^+)_\solid} (B, B^+)_\solid.
\end{align*}
Writing $N = N'^a$ for some $N' \in \D_\solid(B, A^+)$ and using that $N \tensor_{(B^a, A^+)_\solid} (B^a, B^+)_\solid = (N' \tensor_{(B, A^+)_\solid} (B, B^+)_\solid)^a$ this amounts to
\begin{align*}
	\IHom_B(\mm \tensor_V B, N' \tensor_{(B, A^+)_\solid} (B, B^+)_\solid) = \IHom_B(\mm \tensor_V B, N') \tensor_{(B, A^+)_\solid} (B, B^+)_\solid.
\end{align*}
Note that on both sides we can pull colimits in $\mm$ out, turning them into limits: On the left-hand side this is clear and on the right-hand side this follows because $- \tensor_{(B, A^+)_\solid} (B, B^+)_\solid$ preserves all small limits by the existence of its left adjoint $f_!$. But $\mm$ is a colimit of copies of $V$ (in $\D_\solid(V)$), hence the above identity reduces to the case $\mm = V$, where it is obvious. This finishes the construction of the functor $f^a_!$.

The construction of $f^{a!}$ can also be performed using relative compactifications: In the case $A = B$ we have $f^{a!} = f^{a*}$ by the above construction of $f^a_!$, so this case is easy. It remains to handle the case $A^+ = B^+$ (more precisely, $B^+$ is the integral closure of the image of $A^+$ in $B$, but this distinction does not matter). Now $f^a_! = f^a_*$ and we can define the right adjoint $f^{a!}$ by the adjoint functor theorem (using the usual arguments to reduce to the presentable case, cf. the proof of \cref{rslt:scheme-6-functor}). It remains to verify the identity $f^{a!} \comp (-)^a = (-)^a \comp f^!$. By passing to left adjoints this amounts to showing $(-)_! \comp f^a_* = f_* \comp (-)_!$, which is easy to see by \cref{rslt:properties-of-almost-lower-shriek-over-alg}.
\end{proof}

We end this section with a result about the flat dimension of the compact projective generators $(A, A^+)_\solid[S] \in \D_\solid(A, A^+)$, for a classical discrete Huber pair $(A, A^+)$. Even though $(A, A^+)_\solid[S]$ is projective, it is in general not flat. Namely, while the functor $- \tensor -$ on $\D_\solid(A, A^+)$ is the left derived functor of its counterpart on the heart, this does not imply the same for the functor $M \tensor -$ for any fixed $M \in \D_\solid(A, A^+)^\heartsuit$. In particular the fact that $(A, A^+)_\solid[S]$ is projective does not imply its flatness. However, if $A^+$ can be generated (over $\Z$) by finitely many elements, then $(A, A^+)_\solid[S]$ does at least have finite flat dimension (see \cref{rslt:finite-type-A+-implies-finite-flat-dimension-of-generators} below). The following argument was suggested by Scholze.

\begin{lemma} \label{rslt:j-shriek-embedding-of-D-solid-Z-into-D-Z}
There is a fully faithful colimit-preserving functor
\begin{align*}
	j_!\colon \D_\solid(\Z) \injto \D(\Z), \qquad \prod_I \Z \mapsto \prod_I \R/\Z[-1].
\end{align*}
The $t$-structure on $\D(\Z)$ restricts to a $t$-structure on the essential image of $\D_\solid(\Z)^\omega$ under $j_!$. Moreover, $j_!$ is left $t$-exact and maps $\D_{\solid,\ge0}(\Z)$ to $\D_{\ge-1}(\Z)$.
\end{lemma}
\begin{proof}
In the following we abbreviate $M^\solid := M \tensor_\Z \Z_\solid$ for a $\Z$-module $M \in \D(\Z)$. We also ignore set-theoretic issues, which can all be resolved by taking filtered colimits over strong limit cardinals $\kappa$.

First note that \cite[Corollary 6.1.(iii)]{condensed-mathematics} implies that for all $M \in \D_\solid(\Z)$ we have $\IHom(\R, M) = 0$. This shows that $\IHom(\prod_I \R, M) = \IHom_\R(\prod_I \R, \IHom(\R, M)) = 0$ and hence that $(\prod_I \R)^\solid = 0$ for every small set $I$. We deduce that
\begin{align*}
	(\prod_I \R/\Z)^\solid = \prod_I \Z[1].
\end{align*}
Now let $\mathcal C^\omega \subset \D(\Z)$ be the full subcategory generated under finite (co)limits by objects of the form $\prod_I \R/\Z$. By the above identity, the solidification functor $(-)^\solid$ restricts to a functor $\mathcal C^\omega \to \mathcal \D_\solid(\Z)^\omega$. This functor is an equivalence: To see this, it is enough to verify that the functor is fully faithful, which boils down to the following computation based on \cite[Theorem 4.3]{condensed-mathematics}:
\begin{align*}
	&\IHom(\prod_J \R/\Z, \prod_I \R/\Z) = \prod_I \IHom(\prod_J \R/\Z, \R/\Z) = \prod_I \IHom(\prod_J \R/\Z, \Z)[1] =\\&\qquad= \prod_I \IHom(\prod_J \Z[1], \Z)[1] = \IHom(\prod_J \Z, \prod_I \Z).
\end{align*}
We thus get a functor $j_!\colon \D_\solid(\Z)^\omega \to \D(\Z)$ with essential image $\mathcal C^\omega$ and this functor preserves all finite (co)limits. Since $\D_\solid(\Z)$ is freely generated by $\D_\solid(\Z)^\omega$ under filtered colimits, we immediately obtain a colimit-preserving functor $j_!\colon \D_\solid(\Z) \to \D(\Z)$ extending the functor on $\D_\solid(\Z)^\omega$. It remains to show that $j_!$ has the desired properties.

We first show that $j_!$ is fully faithful. This boils down to showing that for every filtered diagram $(M_j)_{j\in J}$ in $\mathcal C^\omega$ and every small set $I$ we have
\begin{align*}
	\IHom(\prod_I \R/\Z, \varinjlim_j M_j) = \varinjlim_j \IHom(\prod_I \R/\Z, M_j).
\end{align*}
Note that $\prod_I \R/\Z$ is pseudocoherent, i.e. admits a resolution by the compact projective generators $\Z[S]$ for extremally disconnected sets $S$ (this is true for every compact abelian group $K$, which by \cite[Theorem 4.5]{condensed-mathematics} reduces to showing that $\Z[K]$ is pseudocoherent, which follows by simplicially resolving $K$ by extremally disconnected sets). In other words we can write $\prod_I \R/\Z$ as a geometric realization of some cosimplicial object in $\D(\Z)$ all of whose terms are of the form $\Z[S]$. Plugging this into the $\IHom$ turns it into a totalization, which commutes with uniformly left-bounded filtered colimits (as it is computed by a spectral sequence). Thus the above identity on $\IHom$'s is true if $(M_j)_j$ is uniformly left-bounded, even for an arbitrary filtered system $(M_j)_j$ in $\D(\Z)$ (i.e. without requiring it to lie in $\mathcal C^\omega$). In general we can now write $(M_j)_j = \varprojlim_k (\tau_{\le k} M_j)_j$. This reduces the above identity to showing that
\begin{align*}
	\varinjlim_j \varprojlim_k \IHom(\prod_I \R/\Z, \tau_{\le k} M_j) = \varprojlim_k \varinjlim_j \IHom(\prod_I \R/\Z, \tau_{\le k} M_j).
\end{align*}
Since $J$ is filtered, this holds as soon as we can show that $\IHom(\prod_I \R/\Z, \tau_{\ge k+1} M_j)$ lies in homological degrees $\ge k + d$ for some fixed constant $d$. In fact we claim that one can take $d = 0$. Namely, this boils down to the following claim: For any map $f\colon \prod_{I_1} \R/\Z \to \prod_{I_2} \R/\Z$, the complex $\IHom(\prod_I \R/\Z, \ker f)$ is concentrated in homological degrees $\ge -1$. To see that this is true, note that $K := \ker f$ is a compact abelian group, so by applying Pontrjagin duality (see \cite[Theorem 4.1.(iii)]{condensed-mathematics}) we get $K = \mathbb D(\mathbb D(K))$ and $\mathbb D(K)$ is discrete. Choosing a $2$-term free resolution of $\mathbb D(K)$ by free $\Z$-modules, we obtain a $2$-term resolution of $K$ in terms of products of $\R/\Z$. But $\IHom(\prod_I \R/\Z, \prod_{I'} \R/\Z)$ is concentrated in degree $0$ (e.g. by the above explicit computation), so the claim follows. This finishes the proof that $j_!$ is fully faithful.

It remains to prove the claim about the $t$-structures. We first show that the $t$-structure on $\D(\Z)$ restricts to a $t$-structure on the essential image $\mathcal C \subset \D(\Z)$ of $j_!$. As every object in $\mathcal C$ is a filtered colimit of objects in $\mathcal C^\omega$, it is enough to show that the $t$-structure on $\D(\Z)$ restricts to a $t$-structure on $\mathcal C^\omega$, i.e. for every $M \in \mathcal C^\omega$ we have $\tau_{\ge0} M, \tau_{\le0} M \in \mathcal C^\omega$. This in turn reduces to showing that for any map $f\colon \prod_I \R/\Z \to \prod_J \R/\Z$, both $\ker f$ and $\coker f$ lie in $\mathcal C^\omega$. To see that this is true for $K := \ker f$, note that $K$ is a compact abelian group, so by the Pontrjagin duality argument employed above we can find a $2$-term resolution of $K$ in terms of products of $\R/\Z$. This shows that $\ker f \in \mathcal C^\omega$. It follows immediately that also $\coker f \in \mathcal C^\omega$ because $\mathcal C^\omega$ is stable under finite limits.

We now show that $j_! \D_{\solid,\ge0}(\Z) \subset \D_{\ge-1}(\Z)$. Given any $M \in \D_{\solid,\ge0}(\Z)$ we can represent $M$ by a complex in homological degree $\ge 0$ all of whose terms are of the form $\bigdsum \prod \Z$. Then $j_!$ maps $M$ to a similar complex, where $\bigdsum \prod \Z$ is transformed into $\bigdsum \prod \R/\Z$ and shifted one to the right. Hence $j_! M \in \D_{\ge-1}(\Z)$, as desired.

We now prove that $j_!$ is left $t$-exact, so let $M \in \D_{\solid,\le0}(\Z)$ be given. It is enough to show that for every $N \in \D_{\ge1}(\mathcal C)$ we have $\Hom(N, j_! M) = 0$. Writing $N = j_! N'$ for some $N' \in \D_\solid(\Z)$ we see that $N' = (j_! N')^\solid$ and thus $N' \in \D_{\solid,\ge1}(\Z)$. But then
\begin{align*}
	\Hom(N, j_! M) = \Hom(j_! N', j_! M) = \Hom(N', M) = 0,
\end{align*}
as desired.
\end{proof}

\begin{lemma} \label{rslt:Z-solid-S-has-flat-dim-le-1}
Let $S$ be a profinite set. Then the object $\Z_\solid[S] \in \D_\solid(\Z)$ has flat dimension $\le 1$, i.e. the functor
\begin{align*}
	- \tensor_{\Z_\solid} \Z_\solid[S]\colon \D_\solid(\Z) \to \D_\solid(\Z)
\end{align*}
maps $\D_{\le0}$ to $\D_{\le1}$.
\end{lemma}
\begin{proof}
By \cref{rslt:j-shriek-embedding-of-D-solid-Z-into-D-Z} the $t$-structure of $\D(\Z)$ restricts to a $t$-structure on the image $\mathcal C^\omega$ of $\D_\solid(\Z)^\omega$ under the embedding $j_!\colon \D_\solid(\Z) \to \D(\Z)$ and this $t$-structure differs from the $t$-structure on $\D_\solid(\Z)$ only by an error of $1$.

Now let $M \in \D_{\solid,\le0}(\Z)$ be given and pick a small set $I$ such that $\Z_\solid[S] = \prod_I \Z$. Write $j_! M = \varinjlim_j P'_j$ as a filtered colimit of objects $P'_j \in \mathcal C^\omega$, where by \cref{rslt:j-shriek-embedding-of-D-solid-Z-into-D-Z} and the above remark about the $t$-structure on $\mathcal C^\omega$ we can assume $P'_j \in \D_{\le0}(\Z)$ for all $j$. We obtain a corresponding colimit $M = \varinjlim_j P_j$ in $\D_\solid(\Z)$, where $P_j = P_j'^\solid$. Then $M \tensor_{\Z_\solid} \Z_\solid[S] = \varinjlim_j (P_j \tensor_{\Z_\solid} \Z_\solid[S])$, so it is enough to show that each $P_j \tensor_{\Z_\solid} \Z_\solid[S]$ lies in $\D_{\solid,\le1}(\Z)$. Since the $t$-structures on $\D_\solid(\Z)^\omega$ and on $\mathcal C^\omega$ differ only by an error of $1$ (via $j_!$), it is enough to show that $j_!(P_j \tensor_{\Z_\solid} \Z_\solid[S])$ lies in $\D_{\le0}(\Z)$. Taking a resolution of $P_j'$ in terms of objects of the form $\prod_J \R/\Z$, $j_!(P_j \tensor_{\Z_\solid} \Z_\solid[S])$ is computed by the same resolution but where we replace $\prod_J \R/\Z$ by $\prod_{I\cprod J} \R/\Z = \prod_I \prod_J \R/\Z$ (see \cite[Proposition 6.3]{condensed-mathematics}), i.e. $j_!(P_j \tensor_{\Z_\solid} \Z_\solid[S]) = \prod_I P'_j$ and this clearly lies in $\D_{\le0}(\Z)$.
\end{proof}

\begin{proposition} \label{rslt:finite-type-A+-implies-finite-flat-dimension-of-generators}
Let $A^+ \to A$ be a map of classical rings such that $A^+$ is of finite type over $\Z$. Then for every profinite set $S$, the object $(A, A^+)_\solid[S] \in \D_\solid(A, A^+)$ has finite flat dimension, i.e. the functor
\begin{align*}
	- \tensor_{(A, A^+)_\solid} (A, A^+)_\solid[S]\colon \D_\solid(A, A^+) \to \D_\solid(A, A^+)
\end{align*}
maps $\D_{\le0}$ to $\D_{\le n}$ for some $n \ge 0$. In fact, if there is a surjective map $\Z[T_1, \dots, T_m] \surjto A^+$ then one can chose $n = m + 1$.
\end{proposition}
\begin{proof}
By \cref{rslt:solid-of-map-A+-A-is-induced-ring-structure} we have $(A, A^+)_\solid = A_{A^+_\solid/}$ which implies that $- \tensor_{(A, A^+)_\solid} (A, A^+)_\solid[S] = - \tensor_{A^+_\solid} A^+_\solid[S]$. Thus we can assume $A = A^+$ from now on. Choose a surjection $\Z[T_1, \dots, Z_m] \surjto A$. Then $A_\solid = (A, \Z[T_1, \dots, T_m])_\solid$, so by the same argument as before we can reduce to the case $A = \Z[T_1, \dots, T_m]$. Using \cref{rslt:Z-solid-S-has-flat-dim-le-1} we are thus reduced to the following claim: Let $B$ be a classical finite-type $\Z$-algebra such that $B_\solid[S] \in \D_\solid(B)$ has flat dimension $\le d$; then $B[T]_\solid[S] \in \D_\solid(B[T])$ has flat dimension $\le d + 1$. To prove this claim, we make use of the fully faithful functor $j_!\colon \D_\solid(B[T]) \injto \D_\solid(B[T], B)$ defined in the proof of \cref{rslt:scheme-6-functor}, which is the left adjoint of $j^* = - \tensor_{(B[T], B)_\solid} B_\solid$. Using the projection formula for $j_!$ (which is also part of the proof of \cref{rslt:scheme-6-functor}) we get, for every $M \in \D_\solid(B[T])$,
\begin{align*}
	& M \tensor_{B[T]_\solid} B[T]_\solid[S] = j^*j_! (M \tensor_{B[T]_\solid} B[T]_\solid[S]) = j^* (M \tensor_{(B[T], B)_\solid} j_! B[T]_\solid[S]).
	\end{align*}
We note:
\begin{itemize}
	\item The functor $j^*\colon \D_\solid(B[T], B) \to \D_\solid(B[T])$ has Tor dimension $\le 1$. This follows directly from the adjunction of $j^*$ and $j_!$ using that $j_! B[T]$ is concentrated in homological degrees $\ge -1$ (cf. \cite[Example 13.15.(2)]{scholze-analytic-spaces}).

	\item The object $j_! B[T]_\solid[S] \in \D_\solid(B[T], B)$ has flat dimension $\le d$, i.e. tensoring with it maps $\D_{\le0}$ to $\D_{\le d}$. Namely, using the projection formula for $j_!$ we compute
	\begin{align*}
		j_! B[T]_\solid[S] &= (B[T], B)_\solid[S] \tensor_{(B[T], B)_\solid} j_! B\\
		&= (B[T], B)_\solid[S] \tensor_{(B[T], B)_\solid} B_\infty/B[-1],
	\end{align*}
	where $B_\infty = B((T^{-1}))$. We therefore have a cofiber sequence
	\begin{align*}
		j_! B[T]_\solid[S] \to (B[T], B)_\solid[S] \to (B[T], B)_\solid[S] \tensor B_\infty
	\end{align*}
	in $\D_\solid(B[T], B)$, which reduces the claimed flat dimension of $j_! B[T]_\solid[S]$ to showing that $(B[T], B)_\solid[S]$ has flat dimension $\le d$ and $(B[T], B)_\solid[S] \tensor B_\infty$ has flat dimension $\le d + 1$. The former follows from our assumption on $B$ and the latter follows from the former using the two-term resolution of $B_\infty$ in \cite[Observation 8.6]{condensed-mathematics}.
\end{itemize}
Combining the above two observations with the formula for $M \tensor_{B[T]_\solid} B[T]_\solid[S]$ beforehand we deduce that $B[T]_\solid[S]$ has flat dimension $\le d + 1$, as desired.
\end{proof}

\subsection{Descent on Schemes} \label{sec:andesc.schemedesc}

In the previous subsection we defined the $\infty$-category of schemes and more generally of discrete adic spaces. We will now prove various descent results on schemes, using the machinery developed on general analytic spaces. For simplicity we state all results only in the non-almost world, but one easily deduces very similar statements for discrete adic spaces over any almost setup $(V,\mm)$ by almost localization (see \cref{rslt:descendability-stable-under-almost-localization}).

In \cref{sec:andesc.descmorph} we defined a general notion of descendable morphisms for analytic spaces. This in particular applies to discrete adic spaces (over any almost setup), so that we have a powerful theory of descent at hand. The only thing that is lacking so far is a good amount of examples of descendable morphisms. To provide such examples, our first goal is to carry over similar descent results from the literature to our setting. Namely, given a morphism $f\colon Y \to X$ of schemes, Mathew defines a descendability condition for $f$ in \cite[\S3]{akhil-galois-group-of-stable-homotopy}, which is refined in \cite[\S11.2]{bhatt-scholze-witt} including a notion of index of descendability analogous to our version. The main difference to our notion of descendability is that Mathew's version caters only to \emph{discrete} quasicoherent sheaves (and in particular is a weaker condition in general). Let us first recall how discrete quasicoherent sheaves fit in our picture.

To every discrete adic space $X$ we associated a symmetric monoidal $\infty$-category $\D_\solid(X)$ of solid quasicoherent sheaves on $X$ by gluing the $\infty$-categories $\D_\solid(A, A^+)$ on affine pieces $\Spa(A, A^+)$. For every solid cutoff cardinal $\kappa$ (see \cref{def:solid-cutoff-cardinal}) we get a full subcategory $\D_\solid(X)_\kappa \subset \D_\solid(X)$ stable under the symmetric monoidal structure and small colimits. In particular we obtain the $\infty$-category
\begin{align*}
	\D(X)_\omega := \D_\solid(X)_\omega
\end{align*}
of \emph{discrete quasicoherent sheaves on $X$}. If $X$ is a classical scheme then $\D(X)_\omega$ coincides with the derived $\infty$-category of classical quasicoherent sheaves on $X$ (because this is true on affine schemes and thus follows in the general case by gluing).

\begin{remark}
By \cref{rslt:solid-discrete-equiv-nuclear} we also have $\D(X)_\omega = \D(X)^\lnuc$, the $\infty$-category of (locally) nuclear quasicoherent sheaves on $X$. We will stick with the intuition of discrete rather than (locally) nuclear objects for now, but on more general spaces $X$ (where the structure sheaves have a non-trivial topology), $\D(X)^\lnuc$ is the correct replacement for $\D(X)_\omega$.
\end{remark}

Now Mathew's notion of descendability (together with the notion of index of descendability from \cite[Definition 11.18]{bhatt-scholze-witt}) can be phrased as follows:

\begin{definition}
A map $f\colon Y \to X$ of discrete adic spaces is \emph{discretely descendable of index $\le n$} if $f$ is qcqs and the natural morphism $(f_* \ri_Y)^{\tensor n} \to \ri_X$ in $\D(X)_\omega$ is zero.
\end{definition}

It is not hard to relate the notion of discrete descendability to that of actual descendability. In fact, we have the following result:

\begin{proposition} \label{rslt:scheme-classical-desc-iff-canonical-compact-desc}
Let $f \colon Y \to X$ be a qcqs morphism of discrete adic spaces and $n \ge 0$ an integer. Then the following are equivalent:
\begin{propenum}
	\item $f$ is discretely descendable of index $n$.
	\item The relative compactification $\overline f^{/X}\colon \overline Y^{/X} \to X$ is descendable of index $n$.
	\item \label{rslt:discrete-desc-equiv-abs-compactification-desc} The absolute compactification $\overline f\colon \overline Y \to \overline X$ is descendable of index $n$.
\end{propenum}
\end{proposition}
\begin{proof}
We first show that (i) and (ii) are equivalent. Let us abbreviate $Y' = \overline Y^{/X}$ and $f' = \overline f^{/X}$. Clearly $f'_* \ri_{Y'} \in \D_\solid(X)$ is locally nuclear (this is true for any discrete quasicoherent sheaf in place of $\ri_{Y'}$ since $f$ is qcqs), so we can view it as an object of $\mathcal E(X)$, where $\mathcal E(X)$ was defined in \cref{sec:andesc.endofun} (note that since $\D_\solid(-)^\lnuc = \D(-)_\omega$ is actually a sheaf, the functor $\D_\solid(X)^\lnuc \injto \mathcal E(X)$ is fully faithful, cf. \cref{rmk:D-nuc-into-E-X-in-general-not-fully-faithful}). We claim that the natural morphism $f'_* \ri_{Y'} \to f'_\natural \ri_{Y'}$ is an isomorphism in $\mathcal E(X)$. This can be checked locally, so we can assume that $Y = \Spa(B, B^+)$ and $X = \Spa(A, A^+)$ are affine. Then $Y' = \Spa(B, A^+)$ and we have (using \cref{rslt:solid-of-map-A+-A-is-induced-ring-structure})
\begin{align*}
	f'_\natural \ri_{Y'} = - \tensor_{(A, A^+)_\solid} (B, A^+)_\solid = - \tensor_{(A, A^+)_\solid} B = f_* \ri_Y,
\end{align*}
as desired. The equivalence of (i) and (ii) follows easily by comparing the two definitions of descendability of index $n$.

Since descendability is stable under base-change (see \cref{rslt:descendable-stable-under-base-change}) it is clear that (iii) implies (ii). It is now enough to show that (i) implies (iii), so assume that $f$ is discretely descendable of index $n$. Note that the pullback functor $\D(\overline X)_\omega \isoto \D(X)_\omega$ is an equivalence (and similarly for $Y$ in place of $X$): This can be checked locally, hence on affines, where it amounts to the obvious statement $\D_\solid(A, \Z)_\omega \isoto \D_\solid(A, A^+)_\omega$. We deduce that $\overline f$ is discretely descendable of index $n$. Then by the equivalence of (i) and (ii), $\overline f$ (which is the relative compactification of itself) is descendable of index $n$.
\end{proof}

\begin{remarks}
\begin{remarksenum}
	\item If both $Y$ and $X$ are affine schemes in \cref{rslt:scheme-classical-desc-iff-canonical-compact-desc} then the statement reduces to the following: A morphism $A \to B$ of discrete rings is discretely descendable (of index $n$) if and only if the morphism $A_\solid \to (B, A)_\solid$ of analytic rings is descendable (of index $n$).

	\item \label{rslt:descendable-implies-discretely-descendable} Given a qcqs morphism $f\colon Y \to X$ of discrete adic spaces, it follows immediately from \cref{rslt:scheme-classical-desc-iff-canonical-compact-desc} and \cref{rslt:descendability-of-composition-implies-desc-of-first-map} that if $f\colon Y \to X$ is descendable of index $n$ then $f$ is discretely descendable of index $\le n$ (use that $f$ factors as $Y \to \overline Y^{/X} \to X$).
\end{remarksenum}
\end{remarks}

The converse of \cref{rslt:descendable-implies-discretely-descendable} does not seem to hold in general: Usually $f$ being descendable is stronger than $f$ being discretely descendable. There are some cases in which the descendability of $f$ can still be deduced from discrete descendability, as the following results show.

\begin{lemma} \label{rslt:descendability-of-proper-morphism}
Let $f\colon Y \to X$ be a proper morphism of discrete adic spaces and $n \ge 0$ an integer. Then $f$ is discretely descendable of index $n$ if and only if $f$ is descendable of index $n$.
\end{lemma}
\begin{proof}
Combine \cref{rslt:partially-proper-iff-equiv-to-compactification} with \cref{rslt:scheme-classical-desc-iff-canonical-compact-desc}.
\end{proof}

\begin{lemma} \label{rslt:fppf-cover-of-rings-is-descendable}
Let $f\colon A \to B$ be a faithfully flat and finitely presented morphism of classical rings. Then the induced morphism $f_\solid\colon A_\solid \to B_\solid$ of analytic rings is descendable of index $\le 2$.
\end{lemma}
\begin{proof}
By \cite[Lemma 034Y]{stacks-project} the morphism $A \to B$ comes via base-change from a faithfully flat morphism of finite-type $\Z$-algebras. By \cref{rslt:descendable-stable-under-base-change} (and \cref{rslt:colimits-of-solid-discrete-rings}) we can therefore assume that $A$ and $B$ are of finite type over $\Z$. Moreover, $f_\solid$ is steady by \cref{rslt:map-of-discrete-rings-is-steady} and $f$ is discretely descendable of index $\le 2$ by the proof of \cite[Proposition 3.31]{akhil-galois-group-of-stable-homotopy}.

Let $\End^L(\D_\solid(A)) \subset \End(\D_\solid(A))$ be the full subcategory spanned by those endofunctors which preserve all small colimits (here $\End(\D(-))$ is as in \cref{def:analytic-endofunctors}). Then via operadic Kan extensions (see \cite[Corollary 3.1.3.4]{lurie-higher-algebra}) one shows that
\begin{align*}
	\End^L(\D_\solid(A)) = \Fun^L_{\D_\solid(A)}(\D_\solid(A)^\omega, \D_\solid(A))
\end{align*}
In particular, for every solid cutoff cardinal $\kappa$, $\End^L(\D_\solid(A)_\kappa)$ is presentable by \cref{rslt:enriched-presheaves-are-presentable}. The evaluation-at-$A$ functor $\ev_{A}\colon \End^L(\D_\solid(A)_\kappa) \to \D_\solid(A)_\kappa$ preserves all small colimits and hence admits a right adjoint $r_\kappa\colon \D_\solid(A)_\kappa \to \End(\D_\solid(A)_\kappa)$. For any compact object $P \in \D_\solid(A)^\omega_\kappa$, the composed functor $\ev_P \comp r_\kappa\colon D_\solid(A)_\kappa \to D_\solid(A)_\kappa$ is right adjoint to the composed functor $\ev_A \comp \ell_P$ (where $\ell_P$ is the left adjoint of $\ev_P$). By enriched Yoneda (see \cref{rslt:weak-enriched-yoneda}) the latter functor is $\ev_A \comp \ell_P = - \tensor P^\vee\colon D_\solid(A)_\kappa \to D_\solid(A)_\kappa$, where $P^\vee = \IHom_A(P, A)$ is the dual of $P$ in $D_\solid(A)_\kappa$. It follows that $\ev_P \comp r_\kappa = \IHom_A(P^\vee, -)_\kappa$, i.e. for all $M \in D_\solid(A)_\kappa$ we have
\begin{align*}
	r_\kappa(M)(P) = \IHom_A(P^\vee, M)_\kappa.
\end{align*}
In particular for fixed $M$ and $P$ the value $r_\kappa(M)(P)$ becomes independent of $\kappa$ for large enough $\kappa$. Hence we can take the limit over $\kappa$ to obtain a functor $r\colon \D_\solid(A) \to \End^L(\D_\solid(A))$ which is right adjoint to $\ev_A$. Now if $P = \prod_I A$ for some set $I$ and $M \in \D(A)_\omega \subset \D_\solid(A)$ is discrete then
\begin{align*}
	r(M)(P) = \IHom_A(P^\vee, M) = \IHom_A(\bigdsum_I A, M) = \prod_I M.
\end{align*}
Given any finite-type map $g\colon A \to C$ of classical rings we claim that the natural morphism $g_{\solid\natural} C \to r(C)$ in $\End^L(\D_\solid(A))$ is an isomorphism. This can be checked by evaluating both sides on compact projective generators, i.e. on $\prod_I A$ for all (small) sets $I$. We know $(g_{\solid\natural} C)(\prod_I A) = (\prod_I A) \tensor_{A_\solid} C_\solid = \prod_I C$ and by the above discussion we have $r(C)(\prod_I A) = \prod_I C$, as desired.

For all $n \ge 0$ let $B^{\tensor n}$ denote the $n$-fold tensor product of $B$ over $A$ and let $f_n\colon A \to B^{\tensor n}$ be the obvious map. By \cref{rslt:colimits-of-solid-discrete-rings}, $(B^{\tensor n})_\solid$ is the $n$-fold tensor product of $B_\solid$ over $A_\solid$. The above discussion implies that $r(B^{\tensor n}) = f_{n\solid\natural} B^{\tensor n}$ for all $n$. Writing $F = \fib(A \to B)$ we deduce $r(F) = \mathcal K_{f_\solid}$. Applying $r$ to the cofiber sequence $F \tensor_A B \to B \to B^{\tensor2}$ we obtain the cofiber sequence $r(F \tensor_A B) \to f_{\solid\natural} B \to f_{2\solid\natural} B^{\tensor 2}$, which implies that $r(F \tensor_A B) = f_{\solid\natural} f_\solid^\natural \mathcal K_{f_\solid}$. If we now apply $r$ to the cofiber sequence $F^{\tensor2} \to F \to F \tensor_A B$ we deduce $r(F^{\tensor2}) = \mathcal K_{f_\solid}^2$.

Now we can conclude: Since $f\colon A \to B$ is discretely descendable of index $\le 2$, the map $F^{\tensor2} \to A$ is zero. Applying $r$ to this map we deduce that the map $\mathcal K_{f_\solid}^2 \to A$ is zero, proving that $f_\solid$ is descendable of index $\le 2$.
\end{proof}

\begin{corollary} \label{rslt:fppf-descendability-for-schemes}
Let $f\colon Y \to X$ be an fppf cover of qcqs classical schemes. Then $f$ is descendable. Moreover, if $X$ can be covered by $n$ affines then the index of $f$ is bounded by $2 c(n)$, where $c(n)$ is the constant from \cref{rslt:steady-covering-is-descendable}.
\end{corollary}
\begin{proof}
Follows easily from \cref{rslt:descendability-is-local} and \cref{rslt:fppf-cover-of-rings-is-descendable}.
\end{proof}

The work so far allows us to prove one of the main results on descendability of schemes from \cite{bhatt-scholze-witt} also in the analytic setting (cf. \cite[Proposition 11.25]{bhatt-scholze-witt}):

\begin{theorem} \label{rslt:h-cover-is-descendable}
Let $f\colon Y \to X$ be an $h$-cover of qcqs noetherian classical schemes. Then $f$ is descendable.
\end{theorem}
\begin{proof}
Just as in the proof of \cite[Proposition 11.25]{bhatt-scholze-witt} we can use \cref{rslt:stability-of-descendable-maps-regarding-composition}, \cref{rslt:descendability-is-local} and \cref{rslt:fppf-descendability-for-schemes} in order to reduce the claim to the case that $f$ is proper. Then the claim follows from \cite[Proposition 11.25]{bhatt-scholze-witt} and \cref{rslt:descendability-of-proper-morphism}.
\end{proof}

In our application to perfectoid geometry, we often want to show descent of a filtered colimit of maps of classical rings. In order to apply \cref{rslt:filtered-colim-of-bounded-desc-is-weakly-desc} to this problem, it is crucial to find bounds on the index of descendability of a morphism $f\colon Y \to X$ which depend only on the base space $X$. For example, \cref{rslt:fppf-descendability-for-schemes} is such a result: The index of descendability of any fppf cover $Y \to X$ is universally bounded by a constant only depending on $X$. This applies in particular to the case that $Y = \bigdunion_{i=1}^n U_i$ for a qcqs open cover $(U_i)_i$ of $X$ (a priori the index of such a cover is bounded only by the constant $c(n)$ depending on $n$). A similar result can be shown for Riemann-Zariski spaces, as follows (see \cref{rslt:descendability-of-open-cover-of-ZR-space} below).

\begin{lemma} \label{rslt:affine-cover-of-proj-scheme-over-valuation-ring}
Let $R$ be a local ring and let $X$ be a classical projective $R$-scheme whose special fiber has dimension $d$. Then there exists a covering of $X$ by $d+1$ affine subsets.
\end{lemma}
\begin{proof}
Choose a closed immersion $X \subset \setP_R^n$. Let $\mm \subset R$ be the maximal ideal, $\kappa = R/\mm$ the fraction field and $X_\mm \subset \setP_\kappa^n$ the special fiber of $X$. It is well-known (cf. \cite[11.3.C]{vakil-rising-sea}) that there exist $d + 1$ homogeneous polynomials $\overline f_0, \dots, \overline f_d \in \kappa[T_0, \dots, T_n]$ with the following property: Letting $\overline H_i \subset \setP_\kappa^n$ denote the vanishing locus of $\overline f_i$, then the intersection of the $\overline H_i$'s is disjoint from $X_\mm$.

We can assume that each $\overline f_i$ has a coefficient $1$. Now for each $i = 0, \dots, d$ choose some lift $f_i \in R[T_0, \dots, T_n]$ of $\overline f_i$ with one coefficient $1$. Let $H_i \subset \setP_R^n$ denote the zero locus of $f_i$. Then each $\setP^n_R \setminus H_i$ is affine: To see this we can pass to a Veronese embedding to assume that $f_i$ has degree $1$; then since one of the coefficients of $f_i$ is $1$, after a transformation of coordinates $f_i$ is of the form $f_i = T_k$ for some $k$, where the claim is clear.

Now define $X_i := X \isect (\setP_R^n \setminus H_i)$. This is a closed subscheme of the affine scheme $\setP_R^n \setminus H_i$ and hence affine. We have therefore constructed open affine subsets $X_i \subset X$. To finish the proof, we need to show that $X = \bigunion_{i=0}^d X_i$. By construction of the $f_i$'s we know that this is true in the special fiber, i.e. $X_\mm = \bigunion_{i=0}^d X_{i,\mm}$. Since each $X_i \subset X$ is open and in particular stable under generalizations, it is enough to show that every $x \in X$ specializes to a point $x_\mm$ in $X_\mm$. This follows from the fact that $X \to \Spec R$ is proper and in particular closed.
\end{proof}

\begin{proposition} \label{rslt:descendability-of-open-cover-of-ZR-space}
Let $V$ be a valuation ring with fraction field $K$ and let $K' \supset K$ be a field extension of finite transcendence degree $d$. Let $X := \Spa(K', V)$ and choose any covering $U := \bigdunion_{i=1}^n U_i \surjto X$ by qcqs open subsets $U_i \subset X$. Then $U \to X$ is descendable of index $\le 2c(d+1)$.
\end{proposition}
\begin{proof}
Consider the category $\mathcal C$ of (classical) integral projective $V$-schemes $Y$ with a dominant point $x \in Y(K')$ (where morphisms are $V$-morphisms mapping the dominant points to each other). This category is filtered and by \cite[Theorem VI.17.41]{zariski-commutative-algebra} there is a canonical homeomorphism $\abs f\colon \abs X \isoto \varprojlim_{Y \in \mathcal C} \abs Y$. This map is given as the inverse limit of closed continuous maps (cf. \cite[Lemma VI.17.4]{zariski-commutative-algebra}) $\abs{f_Y}\colon \abs X \to \abs Y$ for $Y \in \mathcal C$. More precisely, $\abs{f_Y}$ maps a valuation $v \in X$ to the unique point $\mathfrak p \in Y$ such that the valuation ring $V_v \subset K'$ of $v$ dominates the local ring $\ri_{\mathfrak p} \subset K'$ of $Y$ at $\mathfrak p$ (where the inclusion $\ri_{\mathfrak p} \subset K'$ is induced by the dominant point $x_Y \in Y(K')$).

We claim that for every $Y \in \mathcal C$ the topological map $\abs{f_Y}\colon \abs X \to \abs Y$ is induced by a map of discrete adic spaces $f_Y\colon X \to Y$. Indeed, choose any affine covering $Y = \bigunion_{i=1}^m Y_i$, $Y_i = \Spec A_i$. Then each $A_i$ is integral and the dominant point $x_Y \in Y(K')$ induces an inclusion $A_i \subset K'$. We obtain canonical morphisms $\Spa(K', A_i) \to Y_i = \Spa(A_i, A_i)$ for all $i$. By the existence of the map $\abs{f_Y}$ we have $X = \bigunion_{i=1}^m \Spa(K', A_i)$ and one checks easily that the morphisms $\Spa(K', A_i) \to Y_i$ glue to give the desired morphism $f_Y\colon X \to Y$.

The fact that $\abs f$ is a homeomorphism (and all $\abs{f_Y}$ are closed) implies that there is some $Y \in \mathcal C$ and a covering $W := \bigdunion_{i=1}^n W_i \surjto Y$, where each $W_i \subset Y$ is a qcqs open subset, such that $\abs U = \abs X \cprod_{\abs Y} \abs W$. It follows that $U = X \cprod_Y W$. The existence of the dominant point $x_Y \in Y(K')$ together with the transcendence degree of $K'/K$ imply that the generic fiber of $Y$ (which is a projective scheme over $K$) has dimension $\le d$. Moreover, the generic fiber is non-empty and hence all fibers are non-empty (because $Y \to \Spec V$ is proper and in particular closed). By \cite[Lemma 0B2J]{stacks-project} the special fiber of $Y$ over $V$ has the same dimension as the generic fiber, so in particular it has dimension $\le d$. Thus by \cref{rslt:affine-cover-of-proj-scheme-over-valuation-ring} $Y$ can be covered by $d + 1$ affine open subsets. It follows from \cref{rslt:fppf-descendability-for-schemes} that the map $W \to Y$ is descendable of index $\le 2c(d+1)$. By \cref{rslt:descendable-stable-under-base-change} the same is true for $U \to X$.
\end{proof}

There is another case where the index of descendability of analytic rings can be bounded in terms of the base ring, namely if that base ring has finite global dimension. In fact, in this case descendability is satisfied very easily:

\begin{proposition} \label{rslt:descendability-for-fin-global-dimension}
For every integer $d \ge 0$ there is a constant $g(d)$ with the following property: Let $A \to B$ be a morphism of classical rings; assume that
\begin{enumerate}[(a)]
	\item $A$ has finite global dimension $\le d$,
	\item for every $k \ge 1$ the $k$-fold tensor product $B^{\tensor_A k}$ lies in $\D_{\le d}(A)$, and
	\item the totalization of the cobar construction
	\begin{center}\begin{tikzcd}
		B \arrow[r, shift left] \arrow[r, shift right] & B \tensor^L_A B \arrow[r] \arrow[r, shift left=2] \arrow[r, shift right=2] & \dots
	\end{tikzcd}\end{center}
	is isomorphic to $A$ via the natural morphism (in the category $\D(A)$).
\end{enumerate}
Then $A_\solid \to (B, A)_\solid$ is descendable of index $\le g(d)$.
\end{proposition}
\begin{proof}
By \cref{rslt:scheme-classical-desc-iff-canonical-compact-desc} we have to show that the morphism $A \to B$ is discretely descendable with index bounded by a constant depending only on $d$. Let $B^\bullet$ denote the cobar construction of $A \to B$, so that $B^\bullet$ lies in $\D_{\le d}(A)$ and $\Tot(B^\bullet) = A$ by assumption. By \cite[Proposition 1.2.4.5]{lurie-higher-topos-theory} we know that $\tau_{\ge-d} \Tot_{2d+1}(B^\bullet) = \tau_{\ge-d} A = A$. Let $X \in \D(A)$ fill the cofiber sequence
\begin{align*}
	A = \Tot(B^\bullet) \to \Tot_{2d+1}(B^\bullet) \to X.
\end{align*}
Then $X \in \D_{\le-d}(A)_\omega$, so since $A$ has global dimension $\le d$, there are no non-zero maps from $X$ to an object of $\D_{\ge1}(A)_\omega$. In particular the induced map $X \to A[1]$ must be zero, so that the above cofiber sequence is split. It follows that $A$ is a retract of $\Tot_{2d+1}(B^\bullet)$. As the structure of the finite limit $\Tot_{2d+1}(B^\bullet)$ only depends on $d$, there is a constant $g(d)$ such that the composition of $g(d)$ consecutive $B$-zero morphisms in $\D(A)$ is zero (this follows from the proof of \cref{rslt:descendable-nilpotent-criterion}, also see the proof of \cite[Proposition 3.27]{akhil-galois-group-of-stable-homotopy} which is the same in Mathew's setting). By the proof of \cref{rslt:descendable-Un-criterion}, $g(d)$ bounds the index of descendability of $A \to B$.
\end{proof}

In the setting of totally disconnected perfectoid spaces (which we will frequently encounter in this thesis) it is very useuful to reduce descent checks to connected components. This can be done using the notion of fs-descendability (see \cref{rslt:fs-descendability-properties}) and the following results.

\begin{lemma} \label{rslt:compare-catsldmod-with-sheaves-on-pi-0}
Let $X = \Spec A$ be a classical affine scheme. Every $x \in \pi_0(X)$ can be canonically viewed as a closed affine subscheme $x = \Spec A_x \subset X$. Then the collection of maps $(A_\solid \to A_{x,\solid})_x$ is an fs-cover.
\end{lemma}
\begin{proof}
Let $\rho\colon X \to \pi_0(X)$ be the projection of topological spaces. For every qcqs open subset $U \subset \pi_0(X)$, $\rho^{-1}(U) \subset X$ is an open and closed subset and hence an affine subscheme; let $A_U$ be its coordinate ring, i.e. $\rho^{-1}(U) = \Spec A_U$. Then $A_\solid \to A_{U\solid}$ is an fs-morphism: The map $A \to A_U$ is finite by the closedness of $\Spec A_U \subset \Spec A$, so that $A_{U\solid} = (A_U, A)_\solid$ and thus it is enough to show that $A \to A_U$ is a flat map of rings; but this follows from the openness of $\Spec A_U \subset \Spec A$. Given $x \in \pi_0(X)$ we have $A_x = \varinjlim_{U \ni x} A_U$, where $U$ ranges over all qcqs open neighbourhoods of $x$ in $\pi_0(X)$. In particular, the flatness of all $A_\solid \to A_{U\solid}$ implies the flatness of $A_\solid \to A_{x\solid}$.

It remains to show that the collection of functors $- \tensor_{A_\solid} A_{x,\solid}$ is conservative. Let $M \in \D_\solid(A)$ such that $M \tensor_{A_\solid} A_{x,\solid} = 0$ for all $x \in \pi_0(X)$. We want to see that $M = 0$; by flatness we can assume that $M$ is concentrated in degree $0$. We can view $M$ as a sheaf of condensed $\rho_* \ri_X$-modules on $\pi_0(X)$ whose value on a qcqs open $U \subset \pi_0(X)$ is $M(U) = M \tensor_{A_\solid} A_{U\solid} = M \tensor_A A_U$. Under this identification, $- \tensor_{A_\solid} A_{x,\solid}$ is the stalk functor at $x \in \pi_0(X)$. But if all stalks of $M$ are zero then $M$ is zero, as desired.
\end{proof}

\begin{corollary} \label{rslt:fs-descendable-can-be-checked-on-connected-components}
Let $f\colon A \to B$ be a map of classical rings. Assume that there is a constant $d$ such that for all connected components $x \in \pi_0(A)$ the fiber $A_x \to A_x \tensor_A B$ is descendable of index $\le d$. Then $f$ is fs-descendable of index $\le d$.
\end{corollary}
\begin{proof}
This follows immediately from \cref{rslt:compare-catsldmod-with-sheaves-on-pi-0}.
\end{proof}

\subsection{Valuation Rings} \label{sec:andesc.valrings}

The local situation of small v-stacks is ruled by valuation rings. In order to develop our theory of quasicoherent sheaves on small v-stacks we thus need a good understanding of valuation rings, particularly of their associated solid analytic ring structure. The main goal of this subsection is to show that given a flat map $V \to V'$ of valuation rings (equipped with the discrete topology) then the associated map of analytic rings $V_\solid \to V'_\solid$ has Tor dimension $\le 1$. Note that bounding the Tor dimension of $V_\solid \to V'_\solid$ is much stronger than that of $V \to V'$ (which is $0$ in this case), as the former Tor dimension measures the exactness of tensoring \emph{topological} $V$-modules with $V'$.

Let us first introduce the relevant notions.

\begin{definition}
Let $A \to B$ be a map of classical rings.
\begin{defenum}
	\item \label{def:flat-map-of-solid-rings} We say that $A_\solid \to B_\solid$ has \emph{Tor dimension $\le d$} if for all $M \in \D_{\solid,\le0}(A)$ the complex $M \tensor_{A_\solid} B_\solid$ lies in $\D_{\solid,\le d}(B)$. We say that $A_\solid \to B_\solid$ is \emph{flat} if it is of Tor dimension $0$.

	\item A static $A_\solid$-module $F \in \D_\solid(A)^\heartsuit$ is called \emph{$B_\solid$-flat} if $F \tensor_{A_\solid} B_\solid$ is static. A $B_\solid$-flat resolution of $M$ is an exact sequence $\dots \to F_n \to \dots F_1 \to F_0 \to M \to 0$ in $\D_\solid(A)^\heartsuit$ such that all $F_n$ are $B_\solid$-flat.
\end{defenum}
\end{definition}

\begin{remark}
In contrast to classical commutative algebra, we cannot characterize the flatness of the map $A_\solid \to B_\solid$ by the flatness of $B$ as an $A_\solid$-module. Thus in order to understand flatness of maps $A_\solid \to B_\solid$ it is not enough to study flatness inside the category $\D_\solid(A)^\heartsuit$.
\end{remark}

\begin{lemma}
Let $A \to B$ be a map of classical rings. Then the following are equivalent:
\begin{lemenum}
	\item $A_\solid \to B_\solid$ has Tor dimension $d$.
	\item Every $A_\solid$-module admits a $B_\solid$-flat resolution of length $d$.
\end{lemenum}
\end{lemma}
\begin{proof}
It is clear that (ii) implies (i). To prove that (i) implies (ii), assume that $A_\solid \to B_\solid$ has Tor dimension $d$ and let $M \in \D_\solid(A)^\heartsuit$ be given. Choose any projective resolution $P_\bullet \to M$ of $M$ and define the resolution $F_\bullet \to M$ by $F_n = P_n$ for $n < d$, $F_n = 0$ for $n > d$ and $F_d = \ker(P_d \to P_{d-1})$. We claim that $F_\bullet \to M$ is a $B_\solid$-flat resolution of $M$, i.e. that $F_d$ is $B_\solid$-flat. Let $F' \in \D_{\solid,[0,d-1]}(A)$ be the complex $F_{d-1} \to \dots \to F_0$. Then $F' \tensor_{A_\solid} B_\solid \in \D_{\solid,[0,d-1]}(B)$. Thus applying $- \tensor_{A_\solid} B_\solid$ to the triangle $F' \to M \to F_d[d]$ in $D_\solid(A)$ and looking at the associated long exact sequence of homology groups proves that $F_d \tensor_{A_\solid} B_\solid$ is concentrated in degree $0$, as desired.
\end{proof}

We have the following nice result about field extensions. While this is only a special case of the maps of valuation rings we seek to understand, it already explains quite well the general strategy we want to pursue.

\begin{proposition} \label{rslt:map-of-fields-is-analytically-flat}
Let $K \subset K'$ be a field extension. Then $K_\solid \to K'_\solid$ is flat.
\end{proposition}
\begin{proof}
By \cite[Lemma 07BV]{stacks-project} we can write $K'$ as a filtered colimit of global complete intersection algebras $A$ over $K$. Then $- \tensor_{K_\solid} K'_\solid$ is the filtered colimit of the functors $- \tensor_{K_\solid} A_\solid$ so that it is enough to show that each $K_\solid \to A_\solid$ is flat. From now on we fix such a global complete intersection algebra $A$ over $K$.

Denote $f\colon \Spec A \to \Spec K$ and fix any static $K_\solid$-module $M$. Recall the functors $f_!$ and $f^!$ from \cref{rslt:scheme-6-functor}. By \cite[Proposition 11.4]{condensed-mathematics} we have
\begin{align*}
	f^! M = (M \tensor_{K_\solid} A_\solid) \tensor_{A_\solid} f^! K.
\end{align*}
We claim that $f^! K = A[d]$ in $\D_\solid(A)$, where $L \in \D(A)^\heartsuit_\omega$ is some line bundle on $\Spec A$ and $d$ is the Krull dimension of $A$. Indeed, by assumption we can write $A = K[T_1, \dots, T_n]/(x_1, \dots, x_r)$ such that $x_1, \dots, x_r$ is a regular sequence and $d = n - r$. This reduces the computation of $f^!$ to the case of polynomial algebras and modding out by non-zero divisors. For the former case use \cref{rslt:scheme-6-functor-poincare-duality} and for the latter case the standard argument works (use the adjunction of $f^!$ and $f_!$ and the fact that $f_! = f_*$ if $f$ is finite).

Plugging $f^! K = A[d]$ into the above equation for $f^! M$ we deduce
\begin{align*}
	M \tensor_{K_\solid} A_\solid = f^! M[-d].
\end{align*}
Thus it is enough to show that $f^! M \in \D_{\solid,\le d}(A)$. By the adjunction of $f^!$ and $f_!$ we have
\begin{align*}
	f_* f^! M = f_* \IHom_A(A, f^! M) = \IHom_K(f_! A, M)
\end{align*}
in $\D_\solid(K)$. This proves that it is enough to show that $f_! A \in \D_{\solid,\ge-d}(K)$. By Noether normalization, $f$ factors as
\begin{align*}
	f\colon \Spec A \xto{i} \Spec K[T_1, \dots, T_d] \xto{g} \Spec K
\end{align*}
where $i$ is finite. Then $f_! A = g_! i_! A = g_! (i_* A)$ and so it is enough to show that $g_!$ maps $\D_{\ge 0}$ to $\D_{\ge -d}$. This follows from the explicit computation of the lower-shriek functor along $R \to R[T]$ for any ring $R$, see \cite[Theorem 8.2]{condensed-mathematics} (plus base-change, see \cref{rslt:scheme-6-functor-proper-base-change}).
\end{proof}

When analyzing the proof of \cref{rslt:map-of-fields-is-analytically-flat} we see that it relies on the following key observations:
\begin{enumerate}[(a)]
	\item Every field extension can be written as a filtered colimit of (global) complete intersection algebras.

	\item For a (global) complete intersection algebra $A$ over a field $K$, the map $K_\solid \to A_\solid$ is flat.
\end{enumerate}
In order to move from fields to valuation rings, we need to generalize both arguments. Let us start by generalizing (a) to valuation rings. In order to do that it is handy to weaken the requirement ``global complete intersection algebra'' slightly and only require it to be a finitely presented, flat, local complete intersection algebra, which is also called a syntomic algebra (cf. \cite[Lemma 069K]{stacks-project}).

\begin{definition}
A map $A \to B$ of classical rings is called \emph{ind-syntomic} if $B$ can be written as a filtered colimit of syntomic $A$-algebras.
\end{definition}

\begin{lemma} \label{rslt:stability-of-ind-syntomic-maps}
Ind-syntomic maps are stable under filtered colimits, compositions, finite products and base-change.
\end{lemma}
\begin{proof}
It is clear that ind-syntomic maps are stable under finite products and base-change (for the latter see \cite[Lemma 00SN]{stacks-project}). To prove that ind-syntomic maps are stable under filtered colimits and compositions, we use the following characterization from \cite[Lemma 07C3]{stacks-project}: A map $A \to A'$ of classical rings is ind-syntomic if and only if for every finitely presented classical $A$-algebra $B$, every $A$-algebra map $B \to A'$ factors over some syntomic $A$-algebra. Now suppose that $A' \to A''$ is another ind-syntomic map. To show that the composition $A \to A''$ is ind-syntomic, fix an $A$-algebra map $B \to A''$ (where $B$ is finitely presented over $A$). Then the induced $A'$-algebra map $A' \tensor_A B \to A''$ factors over some syntomic $A'$-algebra $A''_0$. It is thus enough to show that the map $A \to A''_0$ is ind-syntomic, i.e. we can assume that $A' \to A''$ is syntomic. Write $A' = \varinjlim_i A_i$ with all $A \to A_i$ syntomic. Then the map $A' \to A''$ comes via base-change from a map $A_i \to A''_i$, which by \cite[Lemma 0C33]{stacks-project} can be assumed to be syntomic. But then $A \to A''$ is the filtered colimit of the syntomic maps $A \to A''_i \tensor_{A_i} A_j$. This proves that ind-syntomic maps are stable under composition. Stability under filtered colimits can be proved similarly and is left to the reader.
\end{proof}

\begin{lemma} \label{rslt:val-extension-to-alg-closure-is-ind-syn}
Let $V \to V'$ be a flat map of valuation rings and assume that the fraction field $K'$ of $V'$ is an algebraic closure of the fraction field $K$ of $V$. Then $V \to V'$ is ind-syntomic.
\end{lemma}
\begin{proof}
The following argument is based on suggestions by K\k{e}stutis \v{C}esnavi\v{c}ius.

Note that the flatness of $V \to V'$ implies that $K \subset K'$ canonically. Let $V_0 := K \isect V'$. Then $V_0$ is valuation ring of $K$ containing $V$ and hence a localization of $V$ (cf. \cite[(1.1.14)]{huber-etale-cohomology}) and in particular ind-syntomic (as open immersions are syntomic). Since $V \to V'$ factors as $V \to V_0 \to V'$, by \cref{rslt:stability-of-ind-syntomic-maps} it is enough to show that $V_0 \to V'$ is syntomic, so that we can replace $V$ by $V_0$ and hence assume that the map $V \to  V'$ is local.

Now let $V^h$ be the henselization of $V$. Clearly $V'$ is its own henselization so that by \cite[Lemma 04GS]{stacks-project} the map $V \to V'$ factors as $V \to V^h \to V'$. By \cite[Lemma 04GV]{stacks-project} $V^h$ is a filtered colimit of étale $V$-algebras, so in particular $V \to V^h$ is ind-syntomic. Hence by \cref{rslt:stability-of-ind-syntomic-maps} it is enough to show that $V^h \to V'$ is ind-syntomic. In other words we can replace $V$ by $V^h$ which allows us to assume that $V$ is henselian.

Let $K'/L/K$ be an intermediate extension which is finite Galois over $K$ and let $W := L \isect V'$. Then $W$ is a valuation ring of $L$ and hence a localization of the integral closure $\tilde V^L$ of $V$ in $L$ (see \cite[\S7, Theorem 12]{zariski-commutative-algebra}). By \cite[Theorem 6.2]{kato-thatte-ram-grp-for-val-rings} the map $V \to \tilde V^L$ is ind-syntomic, hence so is the map $V \to W$ (again using that localizations are ind-syntomic and \cref{rslt:stability-of-ind-syntomic-maps}). Taking the union of all such $L$ we obtain the separable closure $K^\sep \subset K'$ of $K$ in $K'$. Let $V^\sep := K^\sep \isect V'$. By the above discussion we see that $V \to V^\sep$ is a filtered colimit of ind-syntomic maps, so by \cref{rslt:stability-of-ind-syntomic-maps} it is itself an ind-syntomic map. We can thus replace $V$ by $V^\sep$ and hence assume that $K'/K$ is purely inseparable.

Let $p$ be the characteristic of $K$. If $p = 0$ then $K = K'$ and we are done. Otherwise let $K'/L/K$ be an intermediate extension such that $L^p \subset K$ and let $W = L \isect V'$. By \cite[Theorem A.3]{vorst-conj-and-valuation-appendix} the map $V \to W$ is ind-syntomic. Note that since $K'/K$ is purely inseparable, any finite extension $L$ of $K$ can be written as successive extensions of the form just described; hence the associated map of valuation rings is a composition of ind-syntomic maps and therefore ind-syntomic. As $K'$ is the union of all such $L$, $V \to V'$ is a filtered colimit of ind-syntomic maps and hence ind-syntomic as well (again using \cref{rslt:stability-of-ind-syntomic-maps}).
\end{proof}

\begin{proposition} \label{rslt:flat-val-ext-is-ind-syn-if-K'-alg-closed}
Let $V \to V'$ be a flat map of valuation rings and assume that the fraction field $K'$ of $V'$ is algebraically closed. Then $V \to V'$ is ind-syntomic.
\end{proposition}
\begin{proof}
In the following proof, all rings and schemes are classical. Let $K$ be the fraction field of $V$, which is a subfield of $K'$ by the flatness of $V \to V'$. Let $\overline K \subset K'$ be the algebraic closure of $K$ in $K'$ (which is algebraically closed since $K'$ is). Let $\overline V := \overline K \isect V'$, which is a valuation ring of $\overline K$. Then $V \to V'$ factors as $V \to \overline V \to V'$ and by \cref{rslt:val-extension-to-alg-closure-is-ind-syn} the map $V \to \overline V$ is ind-syntomic. Hence by \cref{rslt:stability-of-ind-syntomic-maps} it is enough to prove that $\overline V \to V'$ is ind-syntomic, so we can replace $V$ by $\overline V$ and assume that $K$ is algebraically closed and that $V = K \isect V'$.

Let $R$ be the ``base ring'' of $V$ and $V'$, defined as follows: Let $p$ be the characteristic of $K$; if $p > 0$ then $R = \Fld_p$, otherwise $R = V \isect \Q$ which is a valuation ring of $\Q$ so that either $R = \Q$ or $R = \Z_{(p)}$ for some prime $p$. Let $B \subset V'$ be a $V$-subalgebra of finite type over $V$. Since $V$ is the filtered union of such $B$'s it is enough to show that the inclusion $B \subset V'$ factors as $B \to B' \to V'$ for some syntomic $V$-algebra $B'$. Since $B$ is flat over $V$ (because it is torsion-free) and of finite type, it is automatically finitely presented over $V$ by \cite[Lemma 053E]{stacks-project}. Therefore there exist a finite-type $R$-subalgebra $A_0 \subset V$ and a finite-type $A_0$-algebra $B_0$ such that $B = B_0 \tensor_{A_0} V$.

By \cite[Theorem 1.2.9]{temkin-tame-distillation-alterations} there exist alterations $X'_0 \surjto \Spec A_0$ and $Y'_0 \surjto \Spec B_0$ and a map $f'\colon Y'_0 \to X'_0$ compatible with $\Spec B_0 \to \Spec A_0$ such that $X'_0$ and $Y'_0$ are regular schemes and there are strict normal crossing divisors $Z' \subset X'_0$, $f'^{-1}(Z') \subset Y'_0$ such that the induced map $(Y'_0, f'^{-1}(Z')) \to (X'_0, Z')$ of log schemes is log smooth. We can assume that $X_0'$ and $Y_0'$ are connected and hence integral (because they are regular). Since $X'_0 \surjto \Spec A_0$ is generically finite and $K$ is algebraically closed, the map $\Spec K \to \Spec A_0$ lifts uniquely to a map $\Spec K \to X'_0$. By properness of $X'_0 \surjto \Spec A_0$ we can then uniquely lift the map $\Spec V \to \Spec A_0$ to a map $\Spec V \to X'_0$. We can similarly lift the map $\Spec V' \to \Spec B_0$ uniquely to a map $\Spec V' \to Y'_0$. Let $Y' := Y_0' \cprod_{X_0'} \Spec V$. There are natural maps $\Spec V' \to Y' \to \Spec B$ of $V$-schemes. We claim that $Y'$ is syntomic over $V$; this finishes the proof as we can then choose $B'$ to be the coordinate ring of any affine open neighbourhood of the image of $\Spec V' \to Y'$.

To prove that $Y' \to \Spec V$ is syntomic, let $x \in X'_0$ be the image of the closed point of $\Spec V$ and let $y \in Y'_0$ be any preimage of $x$ under $f'$. The maps $\Spec V \to X'_0$ and $\Spec K' \to Y'_0$ induce inclusions $\ri_{X'_0,x} \subset V$ and $\ri_{Y'_0,y} \subset K'$ which are compatible with the necessarily injective map $\ri_{X'_0,x} \injto \ri_{Y'_0,y}$ by the uniqueness of choices. By the explicit description of log smoothness in our case (cf. \cite[\S1.2.7]{temkin-tame-distillation-alterations}) there are regular parameters $t_1, \dots, t_s \in \ri_{X'_0,x}$ and $z_1, \dots, z_r \in \ri_{Y'_0,y}$ which define the divisors $Z'$ and $f'^{-1}(Z')$ resp., such that there exist open neighbourhoods $U_x$ of $x$ and $U_y$ of $y$ inducing a smooth map
\begin{align*}
	U_y \to U_x \cprod_{\Spec \Z[t_1, \dots, t_s]} \Spec \Z[z_1, \dots, z_r].
\end{align*}
Now the strict normal crossing divisor $Z'$ on $X'_0$ is locally at $x$ given by the principal ideal $t_1 \dots t_s \ri_{X'_0,x} \subset \ri_{X'_0,x}$ and the strict normal crossing divisor $f'^{-1}(Z')$ on $Y'_0$ is locally at $y$ given by the principal ideal $t_1 \cdots t_s \ri_{Y'_0,y} = z_1 \cdots z_r \ri_{Y'_0,y} \subset \ri_{Y'_0,y}$. Therefore we have $t_1 \cdots t_s = u z_1 \cdots z_r$ for some unit $u \in \ri\units_{Y'_0,y}$. On the other hand by our choice of the $t_i$'s and $z_j$'s (cf. \cite[\S1.2.7]{temkin-tame-distillation-alterations}) we can write $t_i = \prod_{j=1}^r z_j^{l_{ij}}$ for certain exponents $l_{ij} \ge 0$. Thus we get
\begin{align*}
	u z_1 \cdots z_r = t_1 \cdots t_s = z_1^{l_{11} + \dots + l_{s1}} \cdots z_r^{l_{1r} + \dots + l_{sr}}.
\end{align*}
But since $z_1, \dots, z_r$ are part of a regular system of parameters for the regular local ring $\ri_{Y'_0,y}$ the exponents on both sides must agree (cf. \cite[Lemma 00NO]{stacks-project}), so that each sum $\sum_i l_{ij}$ is $1$. It follows that for fixed $j$ precisely one of the $l_{ij}$'s is equal to $1$ while all the others are $0$. Then up to reordering we obtain indices $n_1, \dots, n_{s-1}, n_s = r$ such that
\begin{align*}
	t_1 = z_1 \cdots z_{n_1}, \quad t_2 = z_{n_1+1} \cdots z_{n_2}, \quad \dots, \quad t_s = z_{n_{s-1}+1} \cdots z_r.
\end{align*}
We now claim that $\Spec V \cprod_{\Spec \Z[t_1, \dots, t_s]} \Spec \Z[z_1, \dots, z_r]$ is syntomic over $V$. Since syntomic maps are stable under composition and base-change it is enough to show for each $i$ that the (affine) scheme
\begin{align*}
	\Spec V \cprod_{\Spec \Z[t_i]} \Spec \Z[z_{n_i+1}, \dots, z_{n_{i+1}}] \isom \Spec V[T_{n_i+1}, \dots, T_{n_{i+1}}]/(t_i - T_{n_i+1} \cdots T_{n_{i+1}})
\end{align*}
is syntomic over $V$. This is clear.

The rest is easy: Let $U = U_y \cprod_{X'_0} V \subset Y'$ be the preimage of $U_y$ under $Y' \to Y'_0$. Since smoothness is stable under base-change, $U$ is smooth over
\begin{align*}
	\Spec V \cprod_{X'_0} (U_x \cprod_{\Spec \Z[t_1, \dots, t_s]} \Spec \Z[z_1, \dots, z_r]) = \Spec V \cprod_{\Spec \Z[t_1, \dots, t_s]} \Spec \Z[z_1, \dots, z_r],
\end{align*}
which itself is syntomic over $V$ by the above. Altogether we obtain that $U$ is syntomic over $V$. Ranging through all possible $y$, we obtain a cover of $Y'$ by open subsets $U \subset Y'$ which are syntomic over $V$; hence $Y'$ is syntomic over $V$.
\end{proof}

\begin{remarks}
\begin{remarksenum}
	\item One can weaken the hypothesis of \cref{rslt:flat-val-ext-is-ind-syn-if-K'-alg-closed} by only requiring that $K'$ has no non-trivial extensions of $p$-power degree, where $p$ is the residue characteristic of $V'$. Indeed, this follows by the same proof if we note additionally that the alterations $X'_0 \surjto \Spec A_0$ and $Y'_0 \surjto \Spec B_0$ are actually $p$-alterations, i.e. induce a $p$-power degree extension of residue fields (this is still \cite[Theorem 1.2.9]{temkin-tame-distillation-alterations}). Note that in particular if $p = 0$ (i.e. $V$ is a $\Q$-algebra) then we do not need to put any assumptions on $K'$.

	\item We expect \cref{rslt:flat-val-ext-is-ind-syn-if-K'-alg-closed} to hold without any assumption on $K'$. This would follow if one can strengthen the used alteration result to an actual resolution of singularities.
	\end{remarksenum}
\end{remarks}

Having covered part (a) of the above discussion in sufficient generality, let us now turn to part (b). In other words, given a valuation ring $V$ and a syntomic $V$-algebra, we want to bound the Tor dimension of $V_\solid \to A_\solid$. In the case that $V = K$ is a field, the proof relied on two facts: Denoting $f\colon \Spec A \to \Spec K$ and $d = \dim A$, the first required fact is that $f^! K = A[d]$ and the second required fact is that $f_! A \in \D_{\solid,\ge-d}(B)$. While the computation for the first fact works in sufficient generality to cover general syntomic maps (even though equality with $A[d]$ only holds locally then), the computation of $f_! A$ relied on Noether normalization, which we do not have available for general $V$. However, we can still show that $f_! A \in \D_{\solid,\ge-d-1}(B)$ for general $V$, which is enough for our purpose. More precisely, we have the following results.

\begin{lemma} \label{rslt:upper-shriek-for-syntomic-map}
Let $f\colon \Spec B \to \Spec A$ be a syntomic map of classical schemes of relative dimension $d$ (cf. \cite[Lemma 02K1]{stacks-project}). Then $f^! A = L[d]$, where $L \in \D(B)^\heartsuit_\omega$ is a line bundle on $B$.
\end{lemma}
\begin{proof}
This can be shown locally on $\Spec B$, so by \cite[Lemmas 00SY, 07D2]{stacks-project} we can assume that $B = A[T_1, \dots, T_n]/(g_1, \dots, g_r)$ such that $g_1, \dots, g_r$ is a Koszul regular sequence. Let us abbreviate $A_n = A[T_1, \dots, T_n]$, then write $f = i \comp f'$, where we denote $f'\colon \Spec A_n \to \Spec A$ and $i\colon \Spec B \injto \Spec A_n$. By \cite[Theorem 11.6]{condensed-mathematics} we have $f'^! A = A_n[n]$. Moreover, by the adjunction of $i^!$ and $i_! = i_*$ we have
\begin{align*}
	i^! A_n = \IHom_B(B, i^! A_n) = \IHom_{A_n}(B, A_n).
\end{align*}
The Koszul regularity of $g_1, \dots, g_r$ implies that the Koszul complex of $g_1, \dots, g_r$ provides a projective resolution of $B$ in $\D(A)_\omega$ and hence in $D_\solid(A)$. Computing $\IHom_{A_n}(B, A_n)$ via this resolution we obtain the dual Koszul complex. By the self-duality of Koszul complexes we deduce $i^! A_n = B[-r]$. In total we get $f^! A = f'^! i^! A = B[n-r]$. Note that $d = n-r$ by definition of relative global complete intersections (see \cite[Definition 00SP]{stacks-project}).
\end{proof}

\begin{lemma}
Let $A \to B$ be a map of classical rings such that the induced map $f\colon \Spec B \to \Spec A$ is an open immersion. Then:
\begin{lemenum}
	\item \label{rslt:open-immersion-of-affine-schemes-has-Tor-dim-1} $A_\solid \to B_\solid$ has Tor dimension $\le 1$.
	\item \label{rslt:lower-shriek-for-affine-open-immersion} $f_! B \in \D_{\solid,[-1,0]}(A)$.
\end{lemenum}
\end{lemma}
\begin{proof}
We first prove (i). For any $g \in A$ let $D(g) \subset \Spec A$ be the associated distinguished open subset. We can then find $g_1, \dots, g_k \in A$ such that $\Spec B = \bigunion_{i=1}^k D(g_i)$ (as subsets of $\Spec A$). For every $I \subset \{ 1, \dots, k \}$ let $g_I = \prod_{i\in I} g_i$, so that $D(g_I) = \bigisect_{i\in I} D(g_i)$, and denote $f_I\colon D(g_I) \injto \Spec A$. Then for every $M \in \D_\solid(A)$, $f^*$ is computed by gluing (i.e. by a totalization of) the objects $f_I^* M$. Since totalizations preserve $\D_{\solid,\le n}(B)$ for any $n$ it is enough to show that each $f_I$ has Tor dimension $1$. In other words, we are reduced to the case $\Spec B = D(g)$ for some $g \in A$, i.e. $B = A_g = A[T]/(gT - 1)$. Note that $A_\solid \to A[T]_\solid$ is flat by \cref{rslt:finite-type-solid-pullback-preserves-limits}. Moreover, since $A[T] \to A_g$ is finite (hence integral) we have $A_{g,\solid} = (A[g], A[T])_\solid$, so that $- \tensor_{A[T]_\solid} A_{g,\solid} = - \tensor_{A[T]} A_g$ (where the right-hand side is the tensor product of condensed $A[T]$-modules). But $A_g$ has a projective resolution of length $1$ over $A[T]$, hence $A[T]_\solid \to A_{g,\solid}$ has Tor dimension $\le 1$. Altogether we obtain that $A_\solid \to A_{g,\solid}$ has Tor dimension $\le 1$, as desired.

Now we prove (ii). Note first that $f_! B \in \D_{\solid,[-n,0]}(A)$ for some $n \ge 0$ because $B$ is finitely presented over $A$. Let $M = \pi_{-n}(f_! B)$. Then
\begin{align*}
	0 \ne \pi_0 \Hom_{\D_\solid(A)}(f_! B, M[-n]) = \pi_n \Hom_{\D_\solid(A)}(f_! B, M) = \pi_n \Hom_{\D_\solid(B)}(B, f^! M).
\end{align*}
But $f^! M = f^* M$ is concentrated in homological degree $\le 1$ by (i), hence $n \le 1$.
\end{proof}

\begin{lemma} \label{rslt:syntomic-map-over-val-ring}
Let $V$ be a valuation ring and let $A$ be a syntomic $V$-algebra of relative dimension $d$ (cf. \cite[Lemma 02K1]{stacks-project}). Denote $f\colon \Spec A \to \Spec V$ the induced map of schemes. Then:
\begin{lemenum}
	\item We have $f_! A \in \D_{\solid,[-d-1,0]}(V)$. Let
	\begin{align*}
		H^{d+1}_c(A/V) := H^{d+1}(f_! A) \ (= \pi_{-d-1}(f_! A))
	\end{align*}
	be the top degree ``cohomology with compact support''.

	\item \label{rslt:syntomic-map-over-val-ring-flatness-criterion} A static $V_\solid$-module $M \in \D_\solid(V)^\heartsuit$ is $A_\solid$-flat if and only if $\pi_0\IHom_V(H^{d+1}_c(A/V), M) = 0$.
\end{lemenum}
In particular, $V_\solid \to A_\solid$ has Tor dimension $\le 1$.
\end{lemma}
\begin{proof}
We start with the proof of (i). We do not have Noether normalization over $V$, so we have to argue differently than in the proof of \cref{rslt:map-of-fields-is-analytically-flat}. Write $X = \Spec A$, choose any embedding $X \injto \setP_V^n$ for some $n$ and let $\overline X \subset \setP_V^n$ be the closure of $X$ (cf. \cite[Lemma 01RG]{stacks-project}). We claim that $\overline X$ has relative dimension $d$ over $\Spec V$. This seems less obvious than one might think, so we provide the following argument (here we use that $V$ is a valuation ring and not just any local ring): $\overline X$ is obtained from $X$ by adding all specializations. In particular one sees immediately that the generic fiber of $\overline X$ has dimension $d$. We now show that the special fiber of $\overline X$ also has dimension $d$ (this is really all we care about, but one can deduce the same for all fibers similarly). For any $x \in \setP^n_V$ let $d(x)$ denote the length of the longest chain of specializations of $x$ in the fiber of $(\setP^n_V)_{f(x)}$ over $f(x)$ (i.e. $d(x)$ is the dimension of the closure of $x$ in the fiber). Since $\overline X$ is obtained by adding specializations, one checks that it is enough to show that if $x \leadsto x'$ is a specialization in $\setP^n_V$ then $d(x') \le d(x)$. Letting $Z := \overline{\{ x \}} \subset \setP^n_V$ denote the closure of $x$, we have to show that the dimension of the special fiber of $Z$ is at most the dimension of the fiber $Z_{f(x)}$ over $f(x)$. Letting $\mathfrak p \subset V$ be the prime ideal given by $f(x)$ and replacing $V$ by $V/\mathfrak p$ and $X$ by $X \cprod_{\Spec V} \Spec V/\mathfrak p$ we can assume that $x$ lies in the generic fiber. But $Z$ is proper (in particular $Z \to \Spec V$ is closed and hence surjective) and irreducible, hence \cite[Lemma 0B2J]{stacks-project} implies the claim.

Applying \cref{rslt:affine-cover-of-proj-scheme-over-valuation-ring} we deduce that there is a covering $\overline X = \bigunion_{i=0}^d X_i$ by open affine subspaces $X_i = \Spec B_i \subset \overline X$. Note that $f\colon X \to \Spec V$ factors as $f = \overline f \comp j$, where $j\colon X \injto \overline X$ is the open immersion and $\overline f\colon \overline X \to \Spec V$ the projection. Thus $f_! A = f_! \ri_X = (\overline f_! \comp j_!) \ri_X = \overline f_* (j_! \ri_X)$. For every $I \subset \{ 0, \dots, d \}$ let $X_I := \bigisect_{i\in I} X_i$ and let $j_I\colon X \isect X_I \injto X_I$ be the open immersion. By separatedness of $\overline X \to V$ all $X_I =: \Spec B_I$ and $X \isect X_I =: \Spec A_I$ are affine. Since $j_! \ri_X$ is glued from the objects $M_I := j_{I!} A_I$, $\overline f_* (j_! \ri_X)$ is the totalization of the associated cosimplicial object $M^\bullet$, where $M^k = \prod_{\abs I \le k} M_I$. But $M^\bullet$ is constant for $\bullet \ge d+1$, so that the associated spectral sequence $E_r^{p,q}$ vanishes for $p \ge d+1$. By \cref{rslt:lower-shriek-for-affine-open-immersion} the $M_I$'s are concentrated in homological degrees $[-1, 0]$, so that $E_r^{p,q}$ also vanishes for $q \not\in [-1, 0]$. It follows that the totalization $f_! A = \overline f_* (j_! A)$ is concentrated in homoloical degree $[-d-1,0]$, as desired.

Now we prove (ii), so let $M \in \D_\solid(V)^\heartsuit$ be given. By \cref{rslt:upper-shriek-for-syntomic-map} we have $f^! V = L[d]$ for some line bundle $L$ on $A$. Let $L^{-1}$ denote the inverse (i.e. dual) of $L$. By \cite[Proposition 11.4]{condensed-mathematics} we have $f^! M = (M \tensor_{V_\solid} A_\solid) \tensor_{A_\solid} L[d]$, hence
\begin{align*}
	M \tensor_{V_\solid} A_\solid = f^! M \tensor_A L^{-1}[-d]
\end{align*}
Since tensoring with $L^{-1}$ is $t$-exact, we see that $M$ is $A_\solid$-flat if and only if $f^! M$ is concentrated in homological degree $\le d$. This can be checked after applying $f_*$; we have
\begin{align*}
	f_* f^! M = f_* \IHom_A(A, f^! M) = \IHom_V(f_! A, M).
\end{align*}
Since $f_! A$ is concentrated in homological degree $\ge -d-1$, $f^! M$ is concentrated in homological degree $\le d+1$ and
\begin{align*}
	\pi_{d+1} f^! M &= \pi_{d+1} \IHom_V(f_! A, M) = \pi_0 \IHom_V(f_! A[d+1], M) =\\
	&= \pi_0 \IHom_V(\pi_0(f_! A[d+1]), M).
\end{align*}
(To see the last equality, apply $\IHom_V(-, M)$ to the triangle $\tau_{\ge1} (f_! A[d+1]) \to f_! A[d+1] \to \pi_0(f_! A[d+1])[0]$, write down the associated long exact sequence and use that in general $\pi_0 \IHom_V(N_1, N_2) = 0$ if $N_1 \in \D_{\solid,\ge 1}(V)$ and $N_2 \in \D_{\solid,\le 0}(V)$.) But $\pi_0(f_! A[d+1]) = H^{d+1}_c(A/V)$, so that altogether we obtain that $\pi_0\IHom_V(H^{d+1}_c(A/V), M) = 0$ if and only if $f^! M$ is concentrated in homological degree $\le d$ if and only $M$ is $A_\solid$-flat.

The claim about the Tor dimension of $V_\solid \to A_\solid$ now follows by observing that the flatness criterion of (ii) is preserved under passage to submodules. Alternatively, note that the above proof of (ii) shows the claimed Tor dimension directly.
\end{proof}

\begin{corollary} \label{rslt:ind-syntomic-over-val-ring-implies-Tor-dim-1}
Let $V$ be a valuation ring and let $V \to A$ be an ind-syntomic map. Then $V_\solid \to A_\solid$ has Tor dimension $\le 1$.
\end{corollary}
\begin{proof}
This follows directly from \cref{rslt:syntomic-map-over-val-ring} because being of Tor dimension $\le 1$ is preserved under filtered colimits.
\end{proof}

We are finally in the position to prove the main result of this subsection: extensions of valuation rings have analytic Tor dimension $\le 1$.

\begin{theorem} \label{rslt:flat-map-of-val-rings-has-sld-Tor-dim-1}
Let $V \to V'$ be a flat map of valuation rings. Then $V_\solid \to V'_\solid$ has Tor dimension $\le 1$.
\end{theorem}
\begin{proof}
If the residue field $K'$ of $V'$ is algebraically closed then the claim follows from \cref{rslt:flat-val-ext-is-ind-syn-if-K'-alg-closed} and \cref{rslt:ind-syntomic-over-val-ring-implies-Tor-dim-1}. We expect that once a better theory of resolution of singularities in characteristic $p$ is available, one can remove the algebraic closedness assumption on $K'$ in order to also conclude in the general case. However, as this is not available yet, we will instead use the following argument.

Let $\overline K'$ be any algebraic closure of $K'$ and let $\overline V' \subset \overline K'$ be a valuation ring of $\overline K'$ which dominates $V'$. For every $n \ge 0$ let $\overline V'^n := \overline V'^{\tensor_{V'} n+1}$, so that $\Spec \overline V'^\bullet \to \Spec V'$ is the Čech nerve of $\Spec \overline V' \to \Spec V'$. Since $V' \to \overline V'$ is ind-fppf (either by \cref{rslt:flat-val-ext-is-ind-syn-if-K'-alg-closed} or by \cite[Lemma 053E]{stacks-project}), \cref{rslt:fppf-cover-of-rings-is-descendable} and \cref{rslt:filtered-colim-of-bounded-desc-is-weakly-desc} imply that $V'_\solid \to \overline V'_\solid$ is weakly descendable. In particular, for every $M \in \catsldmod V$ we have
\begin{align*}
	M \tensor_{V_\solid} V'_\solid = \varprojlim_{n\in\Delta} M \tensor_{V_\solid} \overline V'^n_\solid
\end{align*}
But the totalization on the right preserves $\D_{\solid,\le1}(V')$, so that it is enough to show that $V_\solid \to \overline V'^n_\solid$ has Tor dimension $\le 1$ for all $n$. By \cref{rslt:ind-syntomic-over-val-ring-implies-Tor-dim-1} it is enough to show that $V \to \overline V'^n$ is ind-syntomic for all $n$. But by \cref{rslt:flat-val-ext-is-ind-syn-if-K'-alg-closed} both $V' \to \overline V'$ and $V' \to \overline V'$ are ind-syntomic, so the claim follows from the fact that ind-syntomic maps are preserved by base-change and composition (see \cref{rslt:stability-of-ind-syntomic-maps}).
\end{proof}

\subsection{Adic Completions and \texorpdfstring{$\sigma$}{Sigma}-Modules} \label{sec:andesc.adic-sigobj}

Fix an almost setup $(V,\mm)$. This section introduces two somewhat unrelated concepts regarding analytic rings and their modules, both of which will be applied to the study of $\varphi$-modules in \cref{sec:ri-pi.phi-mod}. The first concept is that of $I$-adic completeness and $I$-adic completion of $\mathcal A$-modules, where $\mathcal A$ is an analytic ring over $(V,\mm)$ with an ideal $I$. Our main result in that regard is that the solid tensor product preserves $I$-adically complete modules in many cases (see \cref{rslt:solid-tensor-product-preserves-adic-complete} below), which in particular provides an easy way to compute the solid tensor product in theses cases. The second concept we want to introduce is that of ``$\sigma$-modules'': Given an analytic ring $\mathcal A$ and an endomorphism $\sigma\colon \mathcal A \to \mathcal A$ we introduce the $\infty$-category $\D(\mathcal A)^\sigma$ of pairs $(M, \sigma_M)$ where $M$ is an $\mathcal A$-module and $\sigma_M\colon M \to M$ is a $\sigma$-linear map such that the induced map $\sigma^* M \isoto M$ is an isomorphism. The most important special case of this construction is when $\mathcal A$ is an $\Fld_p$-algebra and $\sigma$ is the Frobenius. The resulting $\infty$-category is the $\infty$-category of $\varphi$-modules over $\mathcal A$ and will be studied extensively in \cref{sec:ri-pi.phi-mod}; in the present section we provide the necessary foundation, mainly to get the $\infty$-categorical subtleties out of the way.

Let us start with the first of the above mentioned concepts: $I$-adic completeness. We first introduce some notation:

\begin{definition}
Let $\mathcal A$ be an analytic ring over $(V,\mm)$.
\begin{defenum}
	\item An \emph{element of $\mathcal A$} is an element $x \in \pi_0 \underline{\mathcal A}_{**\omega}$, where $(-)_\omega$ denotes the discretization. Equivalently, an element is a map of rings $\Z[T] \to \underline{\mathcal A}_**$ or a map of analytic rings $V^a[T] \to \mathcal A$ over $(V,\mm)$.

	\item An \emph{ideal of $\mathcal A$} is an ideal $I \subset \pi_0 \underline{\mathcal A}_{**\omega}$.
\end{defenum}
\end{definition}

\begin{definition}
Let $\mathcal A$ be an analytic ring over $(V,\mm)$.
\begin{defenum}
	\item For every element $x$ of $\mathcal A$, there is a natural map $\mathcal A \xto{x} \mathcal A$ of $\mathcal A$-modules defined as the base-change of the map $V^a[T] \xto{T} V^a[T]$ along the map $V^a[T] \to \mathcal A$ given by $x$. In particular for every $\mathcal A$-module $M \in \D(\mathcal A)$ there is a natural map $M \xto{x} M$.

	\item For every collection $x_1, \dots, x_k$ of elements of $\mathcal A$ and every $\mathcal A$-module $M \in \D(\mathcal A)$ we define
	\begin{align*}
		M/(x_1, \dots, x_k) := M \tensor_{V^a[T_1, \dots, T_k]} V^a[T_1, \dots, T_k]/(T_1, \dots, T_k) \in \D(\mathcal A),
	\end{align*}
	where the map $V^a[T_1, \dots, T_k] \to \mathcal A$ is given by $x_1, \dots, x_k$. If $k = 1$ we also write $M/x := M/(x) = \cofib(M \xto{x} M)$. Be warned that $M/(x_1, \dots, x_k)$ depends on the elements $x_1, \dots, x_k$ and not only on the ideal they generate. For example $M/0 = M \dsum M[1]$.

	\item For every element $x$ of $\mathcal A$ we define $\mathcal A[1/x] := A \tensor_{V^a[T]} V^a[T, T^{-1}]$, where the map $V^a[T] \to \mathcal A$ is given by $x$.
\end{defenum}
\end{definition}

We can now come to the definition and basic proeprties of $I$-adically complete modules over an analytic ring. The following exposition vaguely follows \cite[\S3.4]{proetale-topology} (what we call ``$I$-adically complete'' is called ``derived $I$-complete'' in the reference).

\begin{definition}
Let $\mathcal A$ be an analytic ring over $(V,\mm)$ with a finitely generated ideal $I$ and let $M \in \D(\mathcal A)$. We say that $M$ is \emph{$I$-adically complete} if for all $x \in I$ we have
\begin{align*}
	\varprojlim(\dots \xto{x} M \xto{x} M \xto{x} M) = 0.
\end{align*}
We denote by
\begin{align*}
	\D(\mathcal A)_{\hat I} \subset \D(\mathcal A)
\end{align*}
the full subcategory spanned by the $I$-adically complete $\mathcal A$-modules. If $I = (x)$ is generated by a single element $x$ then we also denote $\D(\mathcal A)_{\hat x} := \D(\mathcal A)_{\hat I}$.
\end{definition}

\begin{lemma} \label{rstl:characterizations-of-I-adically-complete-modules}
Let $\mathcal A$ be an analytic ring over $(V,\mm)$ with a finitely generated ideal $I$ and let $M \in \D(\mathcal A)$. Then the following are equivalent:
\begin{lemenum}
	\item $M$ is $I$-adically complete.

	\item For every $x \in I$ the natural map $M \isoto \varprojlim_n M/x^n$ is an isomorphism.

	\item For every $x \in I$ and every $N \in \D(\mathcal A[1/x])$ we have $\IHom_{\mathcal A}(N, M) = 0$.

	\item \label{rslt:I-adically-complete-equiv-limit-over-mod-generators-to-the-n} For some (resp. every) collection of generators $x_1, \dots, x_k \in I$ the natural map
	\begin{align*}
		M \isoto \varprojlim_n M/(x_1^n, \dots, x_k^n)
	\end{align*}
	is an isomorphism.

	\item \label{rslt:I-adically-complete-only-depends-on-forget-to-polynomial-ring} For some (resp. every) collection of generators $x_1, \dots, x_k \in I$, $M$ is $(T_1, \dots, T_k)$-adically complete when viewed as a $V^a[T_1, \dots, T_k]$-module via the map $V^a[T_1, \dots, T_k] \to \mathcal A$ given by $x_1, \dots, x_k$.

	\item \label{rslt:M-is-I-adically-complete-iff-all-pi-n-M-are-so} All $\pi_i(M)$ are $I$-adically complete.
\end{lemenum}
\end{lemma}
\begin{proof}
To prove the equivalence of (i) and (ii) we observe that for every $x \in I$ there is the following homotopy coherent diagram in $\D(\mathcal A)$:
\begin{center}\begin{tikzcd}
	\cdots \arrow[r,"x"] & M \arrow[r,"x"] \arrow[d,"x^4"] & M \arrow[r,"x"] \arrow[d,"x^3"] & M \arrow[r,"x"] \arrow[d,"x^2"] & M \arrow[d,"x"]\\
	\cdots \arrow[r,"\id"] & M \arrow[r,"\id"] & M \arrow[r,"\id"] & M \arrow[r,"\id"] & M
\end{tikzcd}\end{center}
The limit of the cofibers along the vertical maps is $\varprojlim_n M/x^n$. Since limits commute with cofibers, this is also the cofiber of the limits along the top and bottom row, i.e. we get
\begin{align*}
	\varprojlim_n M/x^n = \cofib(\varprojlim(\dots \xto{x} M \xto{x} M \xto{x} M) \to M).
\end{align*}
The equivalence of (i) and (ii) is now evident.

It is clear that (iii) implies (i): Just take $N = \mathcal A[1/x] = \varinjlim(\mathcal A \xto{x} \mathcal A \xto{x} \dots)$. For the converse note that $\IHom_{\mathcal A}(N, M) = \IHom_{\mathcal A[1/x]}(N, \IHom_{\mathcal A}(\mathcal A[1/x], M))$.

We now prove that (iv) implies (iii), so assume that the given isomorphism holds for some collection of generators $x_1, \dots, x_k \in I$. By plugging this isomorphism into $\dots \xto{x} M \xto{x} M$ and commuting limits we easily deduce that $\IHom_{\mathcal A}(\mathcal A[1/x_i], M) = \varprojlim(\dots \xto{x_i} M \xto{x_i} M) = 0$ for all $i$. Hence (iii) holds for all $x \in I$ with $x = x_i$ for some $i$. In particular it holds for all $x \in I$ of the form $x = ax_i$ for some element $a$ of $\mathcal A$. Now given any $x \in I$ we write $x = a_1 x_1 + \dots + a_k x_k$ and repeatedly apply \cite[Lemma 3.4.6]{proetale-topology}.

We now prove that (ii) implies (iv). This follows from the observation that
\begin{align*}
	M/(x_1^n, \dots, x_k^n) = ((M/x_1^n)/x_2^n)/\dots/x_k^n
\end{align*}
and that $M \mapsto M/x$ preserves limits.

It is clear that (iv) and (v) are equivalent. The equivalence of (vi) with the other properties follows in the same way as in \cite[Proposition 3.4.4]{proetale-topology}.
\end{proof}

Having introduced $I$-adically complete $\mathcal A$-modules, we now turn our focus on the $I$-adic completion functor. It exists by the following result:

\begin{lemma} \label{rslt:I-adic-completion-functor}
Let $\mathcal A$ be an analytic ring over $(V,\mm)$ with a finitely generated ideal $I$. Then the inclusion $\D(\mathcal A)_{\hat I} \injto \D(\mathcal A)$ admits a left adjoint
\begin{align*}
	(-)_{\hat I}\colon \D(\mathcal A) \to \D(\mathcal A)_{\hat I}, \qquad M \mapsto M_{\hat I}.
\end{align*}
If $I$ is generated by $x_1, \dots, x_k$ then for all $M \in \D(\mathcal A)$ we have
\begin{align*}
	M_{\hat I} = \varprojlim_n M/(x_1^n, \dots, x_k^n).
\end{align*}
\end{lemma}
\begin{proof}
Fix generators $x_1, \dots, x_k \in I$. We need to show that the functor $M \mapsto \hat M := \varprojlim_n M/(x_1^n, \dots, x_k^n)$ is left adjoint to the inclusion $\D(\mathcal A)_{\hat I} \injto \D(\mathcal A)$. There is a natural transformation $M \to \hat M$, so that for all $N \in \D(\mathcal A)_{\hat I}$ we get a natural map
\begin{align*}
	\Hom_{\mathcal A}(\hat M, N) \to \Hom_{\mathcal A}(M, N)
\end{align*}
of anima. We need to show that this map is an isomorphism. By \cref{rslt:I-adically-complete-equiv-limit-over-mod-generators-to-the-n} we have $N = \varprojlim_n N/(x_1^n, \dots, x_k^n)$. Denoting $\mathcal A_n = \mathcal A/(x_1^n, \dots, x_k^n)$ we deduce
\begin{align*}
	\Hom_{\mathcal A}(M, N) = \varprojlim_n \Hom_{\mathcal A}(M, N \tensor_{\mathcal A} \mathcal A_n) = \varprojlim_n \Hom_{\mathcal A_n}(M \tensor_{\mathcal A} \mathcal A_n, N \tensor_{\mathcal A} \mathcal A_n)
\end{align*}
and similarly for $\hat M$ in place of $M$. But clearly $M \tensor_{\mathcal A} \mathcal A_n = \hat M \tensor_{\mathcal A} \mathcal A_n$ (because $- \tensor_{\mathcal A} \mathcal A_n$ amounts to the composition of the functors $(-)/x_i^n$), proving the claim.
\end{proof}

\begin{definition}
Let $\mathcal A$ be an analytic ring over $(V,\mm)$ with a finitely generated ideal $I$. The functor $\D(\mathcal A) \to \D(\mathcal A)_{\hat I}$, $M \mapsto M_{\hat I}$ from \cref{rslt:I-adic-completion-functor} is called the \emph{$I$-adic completion}. If $I$ is clear from context we will also denote $\hat M = M_{\hat I}$.

We get an induced symmetric monoidal structure $\hat\tensor_{\mathcal A}$ on $\D(\mathcal A)_{\hat I}$ which is computed as $M \hat\tensor_{\mathcal A} N = (M \tensor_{\mathcal A} N)_{\hat I}$.
\end{definition}

From the above results we can easily deduce the following nice representation of the symmetric monoidal $\infty$-category $\D(\mathcal A)_{\hat I}$:

\begin{lemma} \label{rslt:aidc-complete-modules-equiv-limit-over-A-mod-I-n}
Let $\mathcal A$ be an analytic ring over $(V,\mm)$ with ideal $I$ generated by some elements $x_1, \dots x_k$. Then the natural functor
\begin{align*}
	\D(\mathcal A)_{\hat I} \isoto \varprojlim_n \D(\mathcal A/(x_1^n, \dots, x_k^n))
\end{align*}
is an equivalence of symmetric monoidal $\infty$-categories.
\end{lemma}
\begin{proof}
The given functor has a right adjoint given by $(M_n)_n \mapsto \varprojlim_n M_n$. Thus the claim follows easily from \cref{rslt:I-adically-complete-equiv-limit-over-mod-generators-to-the-n}.
\end{proof}

We now come to the preservation of $I$-adically complete modules under the solid tensor product, as promised at the beginning of this subsection. Let us first introduce the following notation for ``generalized Huber pairs'':

\begin{definition} \label{def:generalized-solid-analytic-ring}
Let $A$ be a $(V,\mm)$-algebra and let $A^+ \to A_{**}$ be a morphism of rings such that $A^+$ is a classical ring and $A_{**} \in \D_\solid(A^+)$. We let $(A_{**}, A^+)_\solid := (A_{**})_{A^+/}$ and then define $(A, A^+)_\solid := (A_{**}, A^+)^a$, which is an analytic ring over $(V,\mm)$. The associated $\infty$-category of modules is denoted $\D_\solid(A, A^+) := \D((A, A^+)_\solid)$. If $A^+ = A_{**\omega}$ then we also write $A_\solid = (A, A^+)_\solid$.
\end{definition}

\begin{lemma} \label{rslt:solid-tensor-product-preserves-adic-complete-over-fin-type-Z}
Let $A$ be a finite-type classical $\Z$-algebra with some ideal $I \subset A$. Let $M, N \in \D^-_\solid(A)$ be $I$-adically complete. Then $M \tensor_{A_\solid} N$ is $I$-adically complete.
\end{lemma}
\begin{proof}
By choosing projective resolutions we can write $M$ and $N$ as geometric realizations of simplicial objects $M_\bullet$ and $N_\bullet$ in $\D_\solid(A)$ such that each $M_n$ and $N_n$ is a direct sum of copies of $A_\solid[S]$ for varying profinite sets $S$. The $I$-adic completion functor is exact and maps $\D_{\ge0}$ to $\D_{\ge-1}$ and thus commutes with uniformly right-bounded geometric realizations (as these can be computed as a colimit of $n$-truncated geometric realizations and the $n$-truncated and $(n+1)$-truncated geometric realizations differ only in homotopy groups $\ge n$, see \cite[Proposition 1.2.4.5.(4)]{lurie-higher-algebra}). Therefore $M = \hat M = \varinjlim_{n\in\Delta} \hat M_n$ and similarly for $N$. It follows that $M \tensor_{A_\solid} N = \varinjlim_{n\in\Delta} \hat M_n \tensor_{A_\solid} \hat N_n$ and using again the fact that completion commutes with uniformly right-bounded geometric realizations we deduce that it is enough to show that each $\hat M_n \tensor_{A_\solid} \hat N_n$ is $I$-adically complete. Thus from now on we can assume that $M$ and $N$ are completed direct sums of copies of $A_\solid[S]$ for varying profinite sets $S$.

We can also assume that $A$ is a polynomial ring over $\Z$. Namely, choose a surjection $A' \surjto A$ by a polynomial ring $A'$ and let $I' \subset A'$ be the preimage under $A$. Then $A_\solid = (A, A')_\solid$, hence
\begin{align*}
	M \tensor_{A_\solid} N = \varinjlim_{n\in\Delta} M \tensor_{A'_\solid} A^{\tensor n-1} \tensor_{A'_\solid} N,
\end{align*}
which follows from the bar resolution of $A$ over $A'$. But $A$ is a finite $A'$-module and can hence be written as a geometric resolution of a simplicial object consisting of finite free $A'$-modules. It follows that $M \tensor_{A'_\solid} A^{\tensor n}$ is $I'$-adically complete and therefore (assuming the claim holds for $A'$) $M \tensor_{A_\solid} N$ is $I'$-adically complete as an $A'$-module. By \cref{rslt:I-adically-complete-only-depends-on-forget-to-polynomial-ring} it is also $I$-adically complete as an $A$-module. Finally, we can assume that $I$ is generated by a single element $x$, because $I$-adic completeness is a condition which can be checked on the elements of $I$.

As a next step we need to understand the $I$-adic completion of $A_\solid[S]$. Let $\hat A$ denote the $I$-adic completion of $A$. We can write it as
\begin{align*}
	\hat A = A[[T]]/(T - x)
\end{align*}
(one sees easily that the right-hand $A$-module is $I$-adically complete and reduces to $A/x^n$ modulo $x^n$, hence it is indeed isomorphic to $\hat A$). Now $A_\solid[S] \isom \prod_J A$ for some set $J$ (by \cref{rslt:A-solid-is-analytic-ring}) from which it follows that $\prod_J A \tensor_{A_\solid} \prod_K A = \prod_{J\cprod K} A$ for all sets $J$ and $K$. Thus the above description of $\hat A$ (as a finite quotient of $\prod_\N A$) implies that
\begin{align*}
	(\prod_J A)_{\hat I} \tensor_{A_\solid} (\prod_K A)_{\hat I} = \prod_J \hat A \tensor_{A_\solid} \prod_K \hat A = \prod_{J \cprod K} (\hat A \tensor_{A_\solid} \hat A) = \prod_{J \cprod K} \hat A = (\prod_{J \cprod K} A)_{\hat I}.
\end{align*}
Let us get back to $M$ and $N$. By the above reductions we can assume that $M$ and $N$ are completed direct sums of copies of $\prod_J \hat A$ for varying sets $J$. Since $I$-adic completions commute with $\omega_1$-filtered colimits (like all countable limits do) we can even reduce to the case that $M = \widehat\bigdsum_{m\ge0} \prod_{J_m} \hat A$ and $N = \widehat\bigdsum_{m\ge0} \prod_{K_m} \hat A$ are \emph{countable} completed direct sums. Note that $M$ and $N$ are static: This follows from the fact that $\hat A$ is static by the above description of it. Therefore $M$ can be seen as topological $A$-modules whose elements are sequences $(x_m)_m$ with $x_m \in \prod_{J_m} \hat A$ such that $x_m$ converges $I$-adically to zero. It follows that we can write $M$ as a union of topological $A$-modules of the form $X = \prod_{m\ge0} x^{\alpha_m} \prod_{J_m} \hat A$ for a monotonous sequence $\alpha_m \in \Z_{\ge0}$ converging to $\infty$. Since $A$ does not have zero-divisors (recall that we reduced to the case that $A$ is a polynomial algebra over $\Z$), each $X$ is isomorphic to $\prod_{m\ge0} \prod_{J_m} \hat A$. Commuting the tensor product $M \tensor_{A_\solid} N$ with the colimit over the $X$'s (for both $M$ and $N$) and using the above formula $\prod_J \hat A \tensor_{A_\solid} \prod_K \hat A = \prod_{J \cprod K} \hat A$ we arrive at
\begin{align*}
	M \tensor_{A_\solid} N = \widehat\bigdsum_{m, m' \ge 0} \prod_{J_m \cprod K_{m'}} \hat A,
\end{align*}
which is indeed $I$-adically complete.
\end{proof}

\begin{proposition} \label{rslt:solid-tensor-product-preserves-adic-complete}
Let $A$ be a $(V,\mm)$-algebra and let $A^+ \to A_{**}$ be a morphism of rings such that $A^+$ is a classical ring and $A_{**} \in \D_\solid(A^+)$. Let $I^+$ be a finitely generated ideal of $A^+$, which induces a finitely generated ideal $I$ of $A$. Suppose that $A$ is $I$-adically complete and let $M, N \in \D^-_\solid(A, A^+)$ be $I$-adically complete. Assume that one of the following assumptions is satisfied:
\begin{propenum}
	\item $A^+$ is of finite type over $\Z$.
	\item For some generators $x_1, \dots, x_k \in I^+$, $M/(x_1, \dots, x_k)$ and $A/(x_1, \dots, x_k)$ are discrete.
\end{propenum}
Then $M \tensor_{(A, A^+)_\solid} N$ is $I$-adically complete.
\end{proposition}
\begin{proof}
First assume that (i) is satisfied. Note that $M_*, N_* \in \D_\solid(A_{**}, A^+)$ are still $I$-adically complete, so we can replace $A$ by $A_{**}$ and thus assume $(V,\mm) = (\Z, \Z)$. Then
\begin{align*}
	M \tensor_{(A, A^+)_\solid} N = \varinjlim_{n\in\Delta} M \tensor_{A^+_\solid} A^{\tensor n-1} \tensor_{A^+_\solid} N,
\end{align*}
which follows from the bar resolution of $A$ over $A^+$. As in the proof of \cref{rslt:solid-tensor-product-preserves-adic-complete-over-fin-type-Z} the geometric realization on the right preserves $I^+$-adically complete objects in $\D_\solid(A^+)$, so it is enough to show that $M \tensor_{A^+_\solid} A^{\tensor n} \tensor_{A^+_\solid} N$ is $I$-adically complete. But this follows from \cref{rslt:solid-tensor-product-preserves-adic-complete-over-fin-type-Z}.

Now assume that (ii) is satisfied. Write $A^+ = \varinjlim_i A^+_i$ as a filtered colimit of finite-type classical $\Z$-algebras. Then
\begin{align*}
	M \tensor_{(A, A^+)_\solid} N = \varinjlim_i M \tensor_{(A, A^+_i)_\solid} N.
\end{align*}
By part (i) each $X_i := M \tensor_{(A, A^+_i)_\solid} N$ is $I$-adically complete. It follows from the condition on $M$ and \cref{rslt:I-adically-complete-equiv-limit-over-mod-generators-to-the-n} that $X_i$ actually lies in $\D_\solid(A, A^+)$ (tensoring with a discrete module is just a colimit of fixed shape). In particular all $X_i$ are isomorphic to each other and hence isomorphic to $M \tensor_{(A, A^+)_\solid} N$, which is therefore also $I$-adically complete.
\end{proof}

\begin{example}
Let $A$ be a ring with an element $\pi$ such that $A/\pi$ is discrete and $A$ is $\pi$-adically complete. Given $\pi$-adically complete $M, N \in \D^-_\solid(A)$ such that $M/\pi$ is discrete, then $M \tensor_{A_\solid} N$ is $\pi$-adically complete.
\end{example}

We have finished our exposition on $I$-adically complete modules and now turn our focus on the second concept discussed in the present subsection: $\sigma$-modules, i.e. modules equipped with a semilinear map. The basic notions can be introduced in a very general setting:

\begin{definition}
Let $\mathcal C$ be an $\infty$-category.
\begin{defenum}
	\item Let $*_\N$ denote the category with one object and endomorphism monoid $\N := \Z_{\ge0}$. We denote the endomorphism of $*$ in $*_\N$ corresponding to $n \in \N$ by $\sigma^n$.

	\item A \emph{$\sigma$-object of $\mathcal C$} is a functor $*_\N \to \mathcal C$. We denote
	\begin{align*}
		\sigobj{\mathcal C} := \Fun(*_\N, \mathcal C)
	\end{align*}
	the $\infty$-category of $\sigma$-objects of $\mathcal C$.

	\item \label{def:sigma-invariants-on-category} Let $X\colon *_\N \to \mathcal C$ be a $\sigma$-object of $\mathcal C$. Then we denote
	\begin{align*}
		X^\sigma := \lim X \in \mathcal C
	\end{align*}
	and call it the \emph{$\sigma$-invariants of $X$} (if this limit exists in $\mathcal C$).
\end{defenum}
\end{definition}

\begin{remark} \label{rslt:sigma-objects-equivalent-to-object-plus-endo}
Let $\mathcal C$ be an $\infty$-category. Let $E_2$ denote the category with two objects $0$ and $1$ and two morphisms $i, s\colon 0 \to 1$. There is a functor $\alpha\colon E_2 \to *_\N$ sending $i$ to $\id$ and $s$ to $\sigma$. It induces a functor
\begin{align*}
	\alpha^*\colon \sigobj{\mathcal C} \to \Fun(E_2, \mathcal C).
\end{align*}
By considering the right Kan extension along $\alpha$ one checks easily that $\alpha^*$ is fully faithful and its essential image consists of those functors $F\colon E_2 \to \mathcal C$ such that $F(i)$ is an isomorphism (the relevant limit is a limit over a countable zig-zag diagram, which can be computed inductively over increasing finite subdiagrams). Such a functor can be identified with the pair $(X, \sigma)$, where $X := F(0)$ and $\sigma := F(i)^{-1} \comp F(s)$. In particular we obtain the more intuitive (yet less precise) description
\begin{align*}
	\sigobj{\mathcal C} = \{ (X, \sigma) \setst \text{$X \in \mathcal C$ is an object, $\sigma\colon X \to X$ is a morphism} \}.
\end{align*}
We also deduce that for every $\sigma$-object $(X, \sigma) \in \sigobj{\mathcal C}$, the $\sigma$-invariants can be computed as
\begin{align*}
	X^\sigma = \equalizer XX\id\sigma
\end{align*}
In particular the $\sigma$-invariants exist as soon as $\mathcal C$ admits equalizers.
\end{remark}

\begin{lemma} \label{rslt:alpha-beta-sigma-invariants}
Let $\mathcal C$ be an $\infty$-category and let $\alpha\colon X \to Y$ and $\beta\colon Y \to X$ be two morphisms in $\mathcal C$. Assume that the two $\sigma$-objects $(X, \beta\alpha)$ and $(Y, \alpha\beta)$ admit $\sigma$-invariants in $\mathcal C$. Then $X^{\beta\alpha} = Y^{\alpha\beta}$ via the following diagram of equivalences:
\begin{center}\begin{tikzcd}
	X^{\beta\alpha} \arrow[r,"\alpha",shift left] & Y^{\alpha\beta} \arrow[l,"\beta",shift left]
\end{tikzcd}
\end{center}
\end{lemma}
\begin{proof}
Using \cref{rslt:sigma-objects-equivalent-to-object-plus-endo} we can write
\begin{align*}
	X^{\beta\alpha} = \equalizer XX\id{\beta\alpha}, \qquad Y^{\alpha\beta} = \equalizer YY\id{\alpha\beta}.
\end{align*}
It is then clear that the claimed morphisms between $X^{\beta\alpha}$ and $Y^{\alpha\beta}$ exist. Their composition is $\beta\alpha$ respectively $\alpha\beta$ which is evidently isomorphic to the identity on the equalizer.
\end{proof}

With the general abstract notion of $\sigma$-objects at hand, we now apply this definition to the setting of analytic rings. In the following, let $\AnRing_{(-)} \to \AlmSetup$ denote the coCartesian fibration classifying the functor from \cref{rslt:functoriality-of-analytic-rings-over-AlmSetup}.

\begin{definition}
\begin{defenum}
	\item An \emph{analytic $\sigma$-ring over $(V,\mm)$} is a $\sigma$-object $(\mathcal A, \sigma)$ of $\AnRing_{(-)}$ such that $\mathcal A$ lies over $(V,\mm)$.

	We say that $(\mathcal A, \sigma)$ is \emph{perfect} if the induced map $\sigma_V\colon (V,\mm) \to (V,\mm)$ is strict and the map $\sigma\colon \mathcal A \to \sigma_{V*} \mathcal A$ is an isomorphism of analytic rings over $(V,\mm)$.

	\item Let $(\mathcal A, \sigma)$ be an analytic $\sigma$-ring over $(V,\mm)$. Via the functor $\mathcal A \mapsto \D(\mathcal A)$ from \cref{rslt:functoriality-of-modules-over-analytic-rings-over-AlmSetup} we get an induced $\sigma$-object in $\infcatinf^\tensor$. In particular we get the associated symmetric monoidal $\infty$-category
	\begin{align*}
		\D(\mathcal A)^\sigma = \equalizer{\D(\mathcal A)}{\D(\mathcal A)}\id{\sigma^*}
	\end{align*}
	(cf. \cref{rslt:sigma-objects-equivalent-to-object-plus-endo}). A \emph{$\sigma$-module over $\mathcal A$} is an object in $\D(\mathcal A)^\sigma$.

	\item Let $(\mathcal A, \sigma)$ be an analytic $\sigma$-ring over $(V,\mm)$. Let $\D_\sigma(\mathcal A) \to *_\N$ be the coCartesian fibration classifying the induced functor $*_\N \to \infcatinf$. We denote
	\begin{align*}
		\D(\mathcal A[T_\sigma]) := \Fun_{*_\N}(*_\N, \D_\sigma(\mathcal A)).
	\end{align*}
	A \emph{lax $\sigma$-module over $\mathcal A$} is an object of $\D(\mathcal A[T_\sigma])$.
\end{defenum}
\end{definition}

\begin{remark}
Let $(\mathcal A, \sigma)$ be an analytic $\sigma$-ring over $(V,\mm)$. Then a $\sigma$-module $M$ over $\mathcal A$ is an $\mathcal A$-module $M \in \D(\mathcal A)$ together with an isomorphism $\sigma_M\colon \sigma^* M \isoto M$. Via adjunctions we can equivalently view $\sigma_M$ as a map $M \to \sigma_* M$, which is the same as a $\sigma$-linear map $M \to M$.

Similarly a lax $\sigma$-module over $\mathcal M$ is an $\mathcal A$-module $M \in \D(\mathcal A)$ together with a map $\sigma^* M \to M$, equivalently a $\sigma$-linear map $M \to M$.
\end{remark}

\begin{remark} \label{rmk:our-sigma-modules-compared-to-Bhatt-Lurie}
In our definition of $\sigma$-modules we require the map $\sigma^* M \to M$ to be an isomorphism. There is a ``dual'' version of this definition, where one instead requires the adjoint map $M \to \sigma_* M$ to be an isomorphism (the latter version appears for example in \cite[\S3.2]{riemann-hilbert-mod-p}). If $(\mathcal A, \sigma)$ is perfect then both definitions are equivalent, but in general they differ.
\end{remark}

The following result provides some basic properties of the $\infty$-category of $\sigma$-modules over an analytic $\sigma$-ring $\mathcal A$:

\begin{lemma}
Let $(\mathcal A, \sigma)$ be an analytic $\sigma$-ring over $(V,\mm)$.
\begin{lemenum}
	\item $\D(\mathcal A)^\sigma$ and $\D(\mathcal A[T_\sigma])$ are stable symmetric monoidal $\infty$-categories which admit all small limits and colimits.

	\item \label{rslt:properties-of-sigma-module-forgetful-functors} There are natural symmetric monoidal functors
	\begin{align*}
		\D(\mathcal A)^\sigma \injto \D(\mathcal A[T_\sigma]) \to \D(\mathcal A)
	\end{align*}
	The first functor is fully faithful and preserves all small colimits. The second functor is conservative and preserves all small limits and colimits. Moreover, if $\sigma^*\colon \D(\mathcal A) \to \D(\mathcal A)$ preserves small limits, then so does the functor $\D(\mathcal A)^\sigma \to \D(\mathcal A[T_\sigma])$.

	\item \label{rslt:sigma-module-dualizable-iff-underlying-module-is-so} An object $P \in \D(\mathcal A)^\sigma$ is dualizable if and only if its underlying object in $\D(\mathcal A)$ is dualizable.
\end{lemenum}
\end{lemma}
\begin{proof}
It is clear that $\D(\mathcal A)^\sigma$ is symmetric monoidal. The symmetric monoidal structure on $\D(\mathcal A[T_\sigma])$ can be constructed using the observation that the coCartesian fibration $\D_\sigma(\mathcal A) \to *_\N$ can be upgraded to a coCartesian family of symmetric monoidal $\infty$-categories $\D_\sigma(\mathcal A) \to *_\N \cprod \opComm$. It is clear that $\D(\mathcal A[T_\sigma])$ admits all small limits and colimits, since so does $\D(\mathcal A)$. The same follows for $\D(\mathcal A)^\sigma$ by writing it as $\varinjlim_\kappa \D(\mathcal A)^\sigma_\kappa$ (cf. \cref{rslt:filtered-colim-preserve-finite-lim-in-infcatinf}) and noting that each $\D(\mathcal A)^\sigma_\kappa$ is presentable. This proves (i).

To prove (ii) we note that $\D(\mathcal A)^\sigma = \varprojlim_{*_\N} \D(-)$ can be understood as the full subcategory of $\Fun_{*_\N}(*_\N, \D_\sigma(\mathcal A))$ consisting of those sections $s\colon *_\N \to \D_\sigma(\mathcal A)$ such that $s(f)$ is coCartesian for all morphisms $f$ in $*_\N$ (see \cite[Proposition 3.3.3.1]{lurie-higher-topos-theory}). This provides the fully faithful functor $\D(\mathcal A)^\sigma \injto \D(\mathcal A[T_\sigma])$ which evidently preserves all small colimits (because so does $\sigma^*$) and can easily be made symmetric monoidal. The forgetful functor $\D(\mathcal A[T_\sigma]) \to \D(\mathcal A)$ is given by the restriction along $* \to *_\N$. It clearly satisfies all the claimed properties.

It remains to prove (iii). Since the forgetful functor $\D(\mathcal A)^\sigma \to \D(A)$ is symmetric monoidal, the ``only if'' part of (iii) is obvious. To prove the converse, let $\mathcal C^\sigma \subset \D(\mathcal A)^\sigma$ denote the full subcategory of those object $P \in \D(\mathcal A)^\sigma$ whose underlying object in $\D(\mathcal A)$ is dualizable and let $\mathcal C \subset \D(\mathcal A)$ be the full subcategory of dualizable objects. Then $\mathcal C$ and $\mathcal C^\sigma$ are symmetric monoidal. We claim that $\IHom$'s exist in $\mathcal C^\sigma$ and that their underlying object in $\mathcal C$ is computed by the $\IHom$ in $\mathcal C$. Namely, this follows formally from the fact that $\sigma^*\colon \mathcal C \to \mathcal C$ commutes with $\IHom$. But now for every $P \in \mathcal C^\sigma$ the natural map $P \tensor \IHom(P, \mathcal A) \isoto \IHom(P, P)$ is an isomorphism, as this can be checked after applying the forgetful functor. We conclude by \cref{rslt:equivalent-defs-for-dualizable}.
\end{proof}

\begin{remark} \label{rslt:weak-sigma-modules-equiv-modules-over-A-T-sigma}
Let $(\mathcal A, \sigma)$ be an analytic $\sigma$-ring over $(V,\mm)$. Using left Kan extensions one checks that the forgetful functor $\D(\mathcal A[T_\sigma]) \to \D(\mathcal A)$ admits a left adjoint which maps $M \in \D(\mathcal A)$ to $\bigdsum_{n\ge0} \sigma^{n*} M$. This provides a monadic adjunction $\D(\mathcal A) \rightleftarrows \D(\mathcal A[T_\sigma])$, so by the Barr-Beck theorem (see \cite[Theorem 4.7.3.5]{lurie-higher-algebra}) we can identify $\D(\mathcal A[T_\sigma])$ with the $\infty$-category of $T_\sigma$-modules in $\D(\mathcal A)$, where $T_\sigma$ is the monad mapping $M \mapsto \bigdsum_{n\ge0} \sigma^{n*} M = (\bigdsum_{n\ge0} \sigma^{n*} \mathcal A) \tensor_{\mathcal A} M$. We obtain an induced associative $\mathcal A$-algebra structure on $A[T_\sigma] := \bigdsum_{n\ge0} \sigma^{n*} \mathcal A$ which allows us to identify $\D(\mathcal A[T_\sigma])$ with the $\infty$-category of $\mathcal A[T_\sigma]$-modules. This justifies our notation.
\end{remark}

Next up we study functoriality of the $\infty$-category of $\sigma$-modules over analytic $\sigma$-rings. We get the following result:

\begin{lemma}
There are natural functors
\begin{align*}
	&\sigobj{\AnRing_{(-)}} \to \infcatinf^\tensor, \qquad (\mathcal A, \sigma) \mapsto \D(\mathcal A)^\sigma,\\
	&\sigobj{\AnRing_{(-)}} \to \infcatinf^\tensor, \qquad (\mathcal A, \sigma) \mapsto \D(\mathcal A[T_\sigma]).
\end{align*}
Moreover, for every morphism $f\colon (\mathcal A, \sigma) \to (\mathcal B, \sigma)$ in $\sigobj{\AnRing_{(-)}}$ the following is true:
\begin{lemenum}
	\item The induced symmetric monoidal functors $f^*\colon \D(\mathcal A)^\sigma \to \D(\mathcal B)^\sigma$ and $f^*\colon \D(\mathcal A[T_\sigma]) \to \D(\mathcal B[T_\sigma])$ act as $- \tensor_{\mathcal A} \mathcal B$ on the underlying modules. Moreover, the former functor is the restriction of the latter.

	\item The functor $f^*\colon \D(\mathcal A[T_\sigma]) \to \D(\mathcal B[T_\sigma])$ admits a right adjoint $f_*\colon \D(\mathcal B[T_\sigma]) \to \D(\mathcal A[T_\sigma])$ which acts as $f_*\colon \D(\mathcal B) \to \D(\mathcal A)$ on the underlying modules. If the base-change morphism $\sigma^* f_* \to f_* \sigma^*$ is an isomorphism of functors $\D(\mathcal B) \to \D(\mathcal A)$ (e.g. if $(\mathcal A, \sigma)$ and $(\mathcal B, \sigma)$ are perfect) then the aforementioned functor restricts to a functor $f_*\colon \D(\mathcal B)^\sigma \to \D(\mathcal A)^\sigma$ which is right adjoint to $f^*$.
\end{lemenum}
\end{lemma}
\begin{proof}
We first construct the functor $(\mathcal A, \sigma) \mapsto \D(\mathcal A[T_\sigma])$. Note that there is a natural functor $D\colon \sigobj{\AnRing_{(-)}} \cprod *_\N \to \infcatinf^\tensor$ coming from the functor $\sigobj{\AnRing_{(-)}} \to \sigobj{\infcatinf^\tensor} = \Fun(*_\N, \infcatinf^\tensor)$. Let $\delta\colon \D_\sigma(-) \to \sigobj{\AnRing_{(-)}} \cprod *_\N \cprod \opComm$ denote the coCartesian family of symmetric monoidal $\infty$-categories classifying the functor $D$. We define a simplicial set $\D(-[T_\sigma])$ over $S := \sigobj{\AnRing_{(-)}} \cprod \opComm$ by
\begin{align*}
	\Map_S(K, \D(-[T_\sigma])) = \Map_{S \cprod *_\N}(K \cprod *_\N, \D_\sigma(-))
\end{align*}
for all simplicial sets $K$ over $S$. Using \cite[Theorem B.4.2]{lurie-higher-algebra} one checks that the induced map $\delta'\colon \D(-[T_\sigma]) \to S$ is a categorical fibration, and in fact a coCartesian fibration of symmetric monoidal $\infty$-categories -- we leave the details to the reader as we do not really need this result later on; see the proofs of \cref{rslt:functoriality-of-composition-monoidal-structure} and \cite[Proposition 2.2.6.20]{lurie-higher-algebra} for a reference. We obtain the desired functor $(\mathcal A, \sigma) \mapsto \D(\mathcal A[T_\sigma])$. The functor $(\mathcal A, \sigma) \mapsto \D(\mathcal A)^\sigma$ can now be defined by passing to full subcategories. However, there is also a more direct way of constructing this functor: It is the composition
\begin{align*}
	\sigobj{\AnRing_{(V,\mm)}} \to \sigobj{\infcatinf^\tensor} = \Fun(*_\N, \infcatinf^\tensor) \to \Fun(*, \infcatinf^\tensor) = \infcatinf^\tensor,
\end{align*}
where the last functor is a right Kan extension.

It remains to prove (i) and (ii), so let $f\colon (\mathcal A, \sigma) \to (\mathcal B, \sigma)$ be given. Part (i) follows immediately from the construction. For part (ii) the claim about $f_*$ on $\D(-[T_\sigma])$ follows from \cref{rslt:weak-sigma-modules-equiv-modules-over-A-T-sigma}. If $\sigma^* f_* \isoto f_* \sigma^*$ is an isomorphism of functors $\D(\mathcal B) \to \D(\mathcal A)$ then clearly $f_*\colon \D(\mathcal B[T_\sigma]) \to \D(\mathcal A[T_\sigma])$ preserves $\sigma$-modules and hence restricts to $\D(\mathcal B)^\sigma \to \D(\mathcal A)^\sigma$.
\end{proof}

If $(\mathcal A, \sigma)$ is a perfect analytic $\sigma$-ring then there is an easy way to transform any lax $\sigma$-module over $\mathcal A$ into a $\sigma$-module over $\mathcal A$:

\begin{lemma} \label{rslt:right-adjoint-from-weak-sigma-mod-to-sigma-mod}
Let $(\mathcal A, \sigma)$ be a perfect analytic $\sigma$-ring over $(V,\mm)$. Then the inclusion $\D(\mathcal A)^\sigma \injto \D(\mathcal A[T_\sigma])$ admits a right adjoint
\begin{align*}
	(-)^=\colon \D(\mathcal A[T_\sigma]) \to \D(\mathcal A)^\sigma, \qquad M \mapsto M^=,
\end{align*}
which acts on the underlying $\mathcal A$-modules as
\begin{align*}
	M \mapsto \varprojlim(\dots \xto{\sigma^{3*} \sigma_M} \sigma^{3*} M \xto{\sigma^{2*}\sigma_M} \sigma^{2*} M \xto{\sigma^*\sigma_M} \sigma^* M \xto{\sigma_M} M).
\end{align*}
\end{lemma}
\begin{proof}
Since $(\mathcal A, \sigma)$ is perfect, the functor $\sigma^*\colon \D(\mathcal A) \isoto \D(\mathcal A)$ is an equivalence. We can thus extend the associated $\sigma$-object in $\infcatinf$ to a functor $*_\Z \to \infcatinf$ (by the usual Kan extension technique). Denote the classifying coCartesian fibration by $\D_{\sigma^\pm}(\mathcal A) \to *_\Z$, so that $\D_{\sigma^\pm}(\mathcal A)$ contains $\D_\sigma(\mathcal A)$. Note that every element in $\Fun_{*_\Z}(*_\Z, \D_{\sigma^\pm}(\mathcal A))$ maps all morphisms in $*_\Z$ to coCartesian morphisms in $\D_{\sigma^\pm}(\mathcal A)$ because every morphism in $*_\Z$ is an isomorphism. By a right Kan extension along $*_\N \to *_\Z$ we deduce that $\D(\mathcal A)^\sigma = \Fun_{*_\Z}(*_\Z, \D_{\sigma^\pm}(\mathcal A))$, that the inclusion $\D(\mathcal A)^\sigma \injto \D(\mathcal A[T_\sigma])$ is given by restriction along $*_\N \to *_\Z$ and that its right adjoint (given by the right Kan extension) has the expected description.
\end{proof}

\begin{remark} \label{rslt:sigma-modules-equiv-modules-over-A-T-sigma-pm-if-A-perfect}
Let $(\mathcal A, \sigma)$ be a perfect analytic $\sigma$-ring over $(V,\mm)$. By the proof of \cref{rslt:right-adjoint-from-weak-sigma-mod-to-sigma-mod} we can identify $\D(\mathcal A)^\sigma = \Fun_{*_\Z}(*_\Z, \D_{\sigma^\pm}(\mathcal A))$. Then by the same argument as in \cref{rslt:weak-sigma-modules-equiv-modules-over-A-T-sigma} we obtain an associative $\mathcal A$-algebra $\mathcal A[T_\sigma^\pm]$ such that $\D(\mathcal A)^\sigma = \D(\mathcal A[T_\sigma^\pm])$. Using this picture, the inclusion $\D(\mathcal A[T_\sigma^\pm]) = \D(\mathcal A)^\sigma \injto \D(\mathcal A[T_\sigma])$ is just the forgetful functor along the map $\mathcal A[T_\sigma] \to \mathcal A[T_\sigma^\pm]$ and its right adjoint is given by $M \mapsto \IHom_{\mathcal A[T_\sigma]}(\mathcal A[T_\sigma^\pm], M)$. From this description one also easily recovers the explicit computation of this right adjoint in \cref{rslt:right-adjoint-from-weak-sigma-mod-to-sigma-mod}.
\end{remark}

Finally, we discuss the functor $M \mapsto M^\sigma$ which associates to every $\sigma$-module $M$ its $\sigma$-invariants $M^\sigma$. Unless $\sigma$ is the identity on $\mathcal A$, we cannot view $M^\sigma$ as an $\mathcal A$-module anymore. This makes the following definition a bit subtle:

\begin{definition} \label{def:sigma-invariants-of-sigma-modules}
Let $(\mathcal A, \sigma)$ be an analytic $\sigma$-ring over $(V,\mm)$ such that $\sigma \isom \id$.
\begin{defenum}
	\item Since $\sigma$ is trivial on $\mathcal A$, we have $\D_\sigma(\mathcal A)^\tensor = *_\N \cprod \D(\mathcal A)^\tensor$ and thus
	\begin{align*}
		\D(\mathcal A[T_\sigma]) = \Fun_{*_\N}(*_\N, \D_\sigma(\mathcal A)) = \Fun(*_\N, \D_\sigma(\mathcal A)).
	\end{align*}
	Restriction along the functor $*_\N \to *$ defines a natural symmetric monoidal functor
	\begin{align*}
		- \tensor_{\mathcal A} (\mathcal A, \sigma)\colon \D(\mathcal A) \to \D(\mathcal A[T_\sigma]), \qquad M \mapsto (M, \id_M).
	\end{align*}

	\item The right Kan extension along $*_\N \to *$ produces a right adjoint
	\begin{align*}
		(-)^\sigma\colon \D(\mathcal A[T_\sigma]) \to \D(\mathcal A), \qquad M \mapsto M^\sigma = \equalizer MM\id\sigma.
	\end{align*}
	of the functor $- \tensor_{\mathcal A} (\mathcal A, \sigma)$ (where the computation via equalizer follows from  \cref{rslt:sigma-objects-equivalent-to-object-plus-endo}).

	\item Let $f\colon (\mathcal A, \sigma) \to (\mathcal B, \sigma)$ be a morphism in $\sigobj{\AnRing_{(-)}}$. We write
	\begin{align*}
		- \tensor_{\mathcal A} (\mathcal B, \sigma)\colon \D(\mathcal A) \to \D(\mathcal B[T_\sigma])
	\end{align*}
	for the symmetric monoidal functor obtained as the composition of $- \tensor_{\mathcal A} (\mathcal A, \sigma)$ and the base-change $\D(\mathcal A[T_\sigma]) \to \D(\mathcal B[T_\sigma])$. If there is no room for confusion we often write $- \tensor_{\mathcal A} \mathcal B$ for $- \tensor_{\mathcal A} (\mathcal B, \sigma)$.

	Whenever $f$ is clear from context we denote
	\begin{align*}
		(-)^\sigma\colon \D(\mathcal B[T_\sigma]) \to \D(\mathcal A), \qquad M \mapsto M^\sigma
	\end{align*}
	the right adjoint of $- \tensor_{\mathcal A} (\mathcal B, \sigma)$, given as the composition of the forgetful functor $\D(\mathcal B[T_\sigma]) \to \D(\mathcal A[T_\sigma])$ and the functor $(-)^\sigma\colon \D(\mathcal A[T_\sigma]) \to \D(\mathcal A)$.
\end{defenum}
\end{definition}

\begin{remark}
In the situation of \cref{def:sigma-invariants-of-sigma-modules}, note that $\mathcal A[T_\sigma] = \mathcal A[T]$ is just the polynomial algebra over $\mathcal A$ (with the induced analytic ring structure). The functor $- \tensor_{\mathcal A} (\mathcal A, \sigma)$ can thus be seen as the forgetful functor $\D(\mathcal A) \to \D(\mathcal A[T])$ along the map $\mathcal A[T] \to \mathcal A$ mapping $T \mapsto 1$. Its right adjoint $(-)^\sigma$ can be described as the functor $\IHom_{\mathcal A[T]}(\mathcal A, -)$.
\end{remark}

\begin{remark}
In the situation of \cref{def:sigma-invariants-of-sigma-modules} the functor $- \tensor_{\mathcal A} (\mathcal A, \sigma)\colon \D(\mathcal A) \to \D(\mathcal A[T_\sigma])$ factors over $\D(\mathcal A)^\sigma$. It follows that if $(\mathcal B, \sigma)$ is perfect then $(-)^\sigma\colon \D(\mathcal B[T_\sigma]) \to \D(\mathcal A)$ factors over $(-)^=\colon \D(\mathcal B[T_\sigma]) \to \D(\mathcal B)^\sigma$.
\end{remark}

\begin{remark} \label{rslt:trivial-sigma-structure-is-functorial}
By the usual techniques involving coCartesian fibrations the functor $- \tensor_{\mathcal A} (\mathcal A, \sigma)\colon \D(\mathcal A) \to \D(\mathcal A[T_\sigma])$ can be made compatible with pullback functors on the source and the target in a functorial way. We leave the details to the reader.
\end{remark}

\clearpage
\section{\texorpdfstring{$\ri_X^{+a}/\pi$}{O\textunderscore X\textasciicircum+a/pi}-Modules} \label{sec:ri-pi}

In this section we develop the 6-functor formalism for quasicoherent sheaves of $\ri_X^{+a}/\pi$-modules on small v-stacks $X$, thereby providing $p$-adic analogs of most of the $\ell$-adic results in \cite{etale-cohomology-of-diamonds}.

The section is structured as follows. In \cref{sec:ri-pi.descent-on-affperfd} we apply the descent formalism from \cref{sec:andesc} to prove various descent results about the $\infty$-category $\Dqcohri(A^+/\pi)$ associated to an affinoid perfectoid space $X = \Spa(A, A^+)$ with pseudouniformizer $\pi \in A^+$. We then apply these results in \cref{sec:ri-pi.def-on-vstack} to construct the symmetric monoidal $\infty$-category $\Dqcohri(X, \Lambda)$ of quasicoherent $\Lambda^a$-modules on a small v-stack $X$, where $\Lambda$ is a system of integral torsion coefficients (e.g. $\Lambda = \ri^+_{X^\sharp}/\pi$ for an untilt $X^\sharp$ of $X$ and a pseudouniformizer $\pi$ on $X^\sharp$). In \cref{sec:ri-pi.bd-and-disc} we study the full subcategories of bounded and of discrete objects in $\Dqcohri(X, \Lambda)$ and relate the latter $\infty$-category to classical sheaves. In \cref{sec:ri-pi.torsors} we introduce representation theory in the context of solid modules and show that quasicoherent $\Lambda^a$-modules on a classifying stack are the same as smooth representations. The \cref{sec:ri-pi.p-bounded} is devoted to single out a nice class of morphisms of small v-stacks, which have ``bounded $p$-cohomological dimension'' in an appropriate sense. These ``$p$-bounded'' morphisms allow us to define lower and upper shriek functors, thereby completing the 6-functor formalism -- this is done in \cref{sec:ri-pi.6-functor}. Afterwards, \cref{sec:ri-pi.compact-dualizable} studies compact, dualizable and perfect objects in $\Dqcohri(X,\Lambda)$. This will in particular provide a good understanding of invertible sheaves, which are needed to define $p$-cohomologically smooth maps in \cref{sec:ri-pi.smooth}. Afterwards, in \cref{sec:ri-pi.phi-mod} we study the relation between $\Fld_p$-sheaves and $\Dqcohri(X, \Lambda)$; the main result of that section is a $p$-torsion Riemann-Hilbert correspondence on small v-stacks. Finally, \cref{sec:ri-pi.poincare} shows that smooth maps of analytic adic spaces over $\Q_p$ are $p$-cohomologically smooth, from which we deduce $p$-adic Poincaré duality in rigid-analytic geometry.

\subsection{Descent on Affinoid Perfectoid Spaces} \label{sec:ri-pi.descent-on-affperfd}

In this subsection we apply the results of \cref{sec:andesc} to the setting of affinoid perfectoid spaces $X = \Spa(A, A^+)$ in order to prove various descent results for the associated $\infty$-category $\Dqcohri(A^+/\pi)$. While being rather technical in nature, the results in this subsection lie at the heart of the thesis.

In the following we fix a compatible choice of pseudo-uniformizer $\pi$ on every affinoid perfectoid space. More precisely we work in the following category:

\begin{definition}
Let $\Perfd$ denote the category of affinoid perfectoid spaces (as in \cite[Definition 8.1]{etale-cohomology-of-diamonds}) and let
\begin{align*}
	\AffPerfd_\pi := \{ (X, \pi) \setst \text{$X = \Spa(A, A^+)$ affinoid perfectoid, $\pi \in A^+$ pseudouniformizer} \}
\end{align*}
be the category of affinoid perfectoid spaces with fixed pseudouniformizer (where morphisms are required to send the chosen pseudouniformizers to each other via pullback).
\end{definition}

In the following, we will often omit the pseudouniformizer from the notation and simply write $X \in \AffPerfd_\pi$.
\begin{definition}
Let $X = \Spa(A, A^+) \in \AffPerfd_\pi$. We denote $\mm_A \subset A^+$ the ideal generated by the pseudouniformizers (i.e. topologically nilpotent units) in $A$. Then $(A^+, \mm_A)$ is an $\omega_1$-compact almost setup (cf. \cref{rslt:m-countably-generated-implies-omega-1-compact}), so we get a good theory of analytic rings over $(A^+, \mm_A)$, in particular solid rings as in \cref{def:solid-almost-ring}. As a main example we obtain the analytic ring $(A^+/\pi)^a_\solid$ and the associated $\infty$-category
\begin{align*}
	\Dqcohri(A^+/\pi)
\end{align*}
of $(A^+/\pi)^a_\solid$-modules. If $\Spa(A, A^+) \to \Spa(B, B^+)$ is a map of affinoid perfectoid spaces then the associated morphism of almost setups $(A^+, \mm_A) \to (B^+, \mm_B)$ is strict, so by \cref{rslt:functoriality-of-analytic-rings-along-strict-morphism} analytic rings over $(B^+, \mm_B)$ are the same as analytic rings over $(A^+, \mm_A)$ with a map from $B$. In particular there is no ambiguity of the notion of ``almost''.
\end{definition}

We now come to the central object of study in this subsection, the $\infty$-category $\DqcohriX X$ of ``quasicoherent solid almost $\ri^+_X/\pi$-modules'' on an affinoid perfectoid space $X$. Quasi-coherence means that a sheaf $\mathcal M \in \DqcohriX X$ should (quasi-pro-étale) locally, say on $Y = \Spa(B, B^+) \in X_\proet$ be identifiable with an object in $\Dqcohri(B^+/\pi)$. This suggests to define $\DqcohriX X$ by ``gluing'' the $\infty$-categories $\Dqcohri(B^+/\pi)$, as follows.

\begin{definition} \label{def:Dqcohri-on-AffPerfd}
We define the functor
\begin{align*}
	(\AffPerfd_\pi)^\opp \to \infcatinf, \qquad X \mapsto \Dqcohri(\ri^+_X/\pi)
\end{align*}
to be the sheafification (see \cref{rslt:sheafification}) of the presheaf $X = \Spa(A, A^+) \mapsto \Dqcohri(A^+/\pi)$ on the big quasi-pro-étale site on $\AffPerfd_\pi$.
\end{definition}

\begin{remark}
There is a minor set-theoretic issue in \cref{def:Dqcohri-on-AffPerfd}: Since $\AffPerfd_\pi$ is a large category, there is no sheafification of presheaves on $\AffPerfd_\pi$ in general. One way to circumvent this issue is by filtering $\AffPerfd_\pi$ by cardinal bounds $\kappa$ (i.e. considering the subsite of $\kappa$-small affinoid perfectoid spaces), for each of which we do have sheafification. Then by \cref{rslt:explicit-description-of-DqcohriX} the $\infty$-category $\DqcohriX X$ does not depend on $\kappa$, so we get a well-defined sheaf $X \mapsto \DqcohriX X$ on all of $\AffPerfd_\pi$.
\end{remark}

Our first goal is to derive a more explicit description of the $\infty$-category $\Dqcohri(\ri^+_X/\pi)$. This can be achieved by explicitly computing this category in the case that $X$ is totally disconnected. In fact, in this case no sheafification is required, as the following results show. To shorten notation, let us introduce the following terminology.

\begin{definition} \label{def:AffPerfd-mod-pi-property}
Let $f\colon Y = \Spa(B, B^+) \to X = \Spa(A, A^+)$ be a morphism in $\AffPerfd_\pi$ and let $\mathcal P$ be a property of analytic rings. Then we say that \emph{$f$ satisfies $\mathcal P$ mod $\pi$} if the morphism $(A^+/\pi)_\solid^a \to (B^+/\pi)_\solid^a$ satisfies $\mathcal P$. However, instead of saying ``mod-$\pi$ quasicoherent sheaves descend along $f$'' we will say ``mod-$\pi$ modules descend along $f$''.
\end{definition}

We begin with a very general base change result for affinoid perfectoid spaces, which will be applied frequently.

\begin{lemma} \label{rslt:base-change-for-affinoid-perfectoid}
Let
\begin{center}\begin{tikzcd}
	Y' \arrow[r] \arrow[d] & Y \arrow[d]\\
	X' \arrow[r] & X
\end{tikzcd}\end{center}
be a cartesian diagram in $\AffPerfd_\pi$. Write $X = \Spa(A, A^+)$, $X' = \Spa(A', A'^+)$, $Y = \Spa(B, B^+)$ and $Y' = \Spa(B', B'^+)$. Then the natural morphism of analytic rings
\begin{align*}
	(A'^+/\pi)_\solid^a \tensor_{(A^+/\pi)_\solid^a} (B^+/\pi)_\solid^a \isoto (B'^+/\pi)_\solid^a.
\end{align*}
is an isomorphism. In particular, base-change along the following diagram holds:
\begin{center}\begin{tikzcd}
	\Dqcohri(B'^+/\pi) \arrow[d] &&&& \arrow[llll,swap,"- \tensor_{(B^+/\pi)^a_\solid} (B'^+/\pi)_\solid"] \arrow[d] \Dqcohri(B^+/\pi)\\
	\Dqcohri(A'^+/\pi) &&&& \arrow[llll,swap,"- \tensor_{(A^+/\pi)^a_\solid} (A'^+/\pi)_\solid"] \Dqcohri(A^+/\pi)
\end{tikzcd}\end{center}
\end{lemma}
\begin{proof}
The claim about base-change follows from \cref{rslt:solid-base-change-for-discrete-rings} once we have shown the claimed isomorphism of analytic rings. To prove this isomorphism, first note that $B'^+$ is the completed integral closure of the image of $\pi_0(A'^+ \tensor_{A^+} B^+)$ in $\pi_0(A \tensor_A B)$ (see \cite[Proposition 6.18]{scholze-perfectoid-spaces}), from which it follows immediately that $B'^+/\pi$ is integral over $\pi_0(A'^+/\pi \tensor_{A^+/\pi} B^+/\pi)$. Thus by \cref{rslt:colimits-of-solid-discrete-almost-rings} the claimed isomorphism of analytic rings reduces to showing that the morphism
\begin{align}
	A'^{+a}/\pi \tensor_{A^{+a}/\pi} B^{+a}/\pi \isoto B'^{+a}/\pi \label{eq:base-change-for-affinoid-perfectoid-helper}
\end{align}
of discrete almost rings is an isomorphism (be aware that the tensor product on the left is derived).

Let us first consider the case that everything is in characteristic $p$. Then $A^+$, $A'^+$ and $B^+$ are perfect rings, hence $A'^+ \tensor_{A^+} B^+$ is concentrated in degree $0$ by \cite[Lemma 3.16]{bhatt-scholze-witt}. Moreover, as $A^+$, $A'^+$ and $B^+$ have no $\pi$-torsion we have $A'^+/\pi = A'^+ \tensor_A A/\pi$ (and similarly for $B^+$). We deduce
\begin{align*}
	A'^+/\pi \tensor_{A^+/\pi} B^+/\pi &= (A'^+ \tensor_{A^+} A^+/\pi) \tensor_{A^+/\pi} (B^+ \tensor_{A^+} A^+/\pi)\\
	&= (A'^+ \tensor_{A^+} B^+) \tensor_{A^+} A^+/\pi.
\end{align*}
Now $A'^+ \tensor_{A^+} B^+$ is a perfect ring and therefore has almost zero $\pi$-torsion. It follows that $A'^{+a}/\pi \tensor_{A^{+a}/\pi} B^{+a}/\pi$ is concentrated in degree $0$.

Still in characteristic $p$, the proof of the isomorphism \cref{eq:base-change-for-affinoid-perfectoid-helper} now reduces to showing that the map $\psi\colon \pi_0(A'^+/\pi \tensor_{A^+/\pi} B^+/\pi) \to B'^+/\pi$ is an almost isomorphism. To see that $\psi$ is almost injective we note that the kernel of the map $A'^+ \tensor_{A^+} B^+ \to A' \tensor_A B$ consists only of $\pi$-torsion (because $A' \tensor_A B$ is obtained from $A'^+ \tensor_{A^+} B^+$ by localizing at $\pi$) which is almost zero, as seen above. It is clear that $\psi$ is almost surjective. This finishes the proof in the case of characteristic $p$.

We now treat the general case (without assumption on the characteristic). If $\pi \divides p$ then there is a pseudouniformizer $\pi^\flat \in A^{\flat+}$ such that $A^+/\pi \isom A^{\flat+}/\pi^\flat$ (and similarly for $A'$ and $B$); this reduces everything to characteristic $p$, which was handled above. For general $\pi$ choose a pseudouniformizer $\varpi \in A^+$ such that $\varpi \divides p$ and choose $m \ge 1$ such that $\pi \divides \varpi^m$. Now for all $n \ge 1$ the exact sequence $0 \to A^+/\varpi \to A^+/\varpi^n \to A^+/\varpi^{n-1} \to 0$ (where the first map is multiplication by $\pi^{n-1}$) induces a cofiber sequence
\begin{align*}
	A'^+/\varpi \tensor_{A^+/\varpi} B^+/\varpi \to A'^+/\varpi^n \tensor_{A^+/\varpi^n} B^+/\varpi^n \to A'^+/\varpi^{n-1} \tensor_{A^+/\varpi^{n-1}} B^+/\varpi^{n-1}
\end{align*}
in $\D(A^+)$. One easily deduces \cref{eq:base-change-for-affinoid-perfectoid-helper} by induction on $n$.
\end{proof}

Let us now prove mod-$\pi$ descent over totally disconnected spaces. Due to the very simple nature of the connected components of a totally disconnected space, we get a strong general descent result (see \cref{rslt:v-cover-of-tot-disc-is-mod-pi-descendable} below).

\begin{lemma} \label{rslt:tors-free-ftype-over-tot-disc-implies-fppf-fibers}
Let $X = \Spa(A, A^+) \in \AffPerfd_\pi$ be totally disconnected and let $B^+$ be a $\pi$-torsion free $A^+$-algebra of finite type. Then for every connected component $x = \Spec (A^+/\pi)_x$ of $\Spec A^+/\pi$, the induced map $(A^+/\pi)_x \to B^+ \tensor_{A^+} (A^+/\pi)_x$ is flat and finitely presented.
\end{lemma}
\begin{proof}
Note that $\pi_0(\Spec A^+) = \pi_0(X) = \pi_0(\Spec A^+/\pi)$. In particular, the given connected component $x$ of $\Spec A^+/\pi$ corresponds to a connected component $\Spec A^+_x \subset \Spec A^+$ and we have $A^+_x/\pi = (A^+/\pi)_x$. Let $B^+_x := B^+ \tensor_{A^+} A^+_x$. Then $B^+_x$ is still $\pi$-torsion free and $B^+_x/\pi = B^+ \tensor_{A^+} (A^+/\pi)_x$, so we have to show that $A^+_x/\pi \to B^+_x/\pi$ is flat and finitely presented.

In the following all rings and constructions are considered in the classical sense, i.e. we ignore potential higher homotopies. Let $\hat A^+_x$ and $\hat B^+_x$ be the $\pi$-adic completions of $A^+_x$ and $B^+_x$, and let $B'^+_x \subset \hat B^+_x$ be the image of $B^+_x \tensor_{A^+_x} \hat A^+_x \to \hat B^+_x$. Then taking the chain of $A^+_x$-algebra maps $B^+_x \to B'^+_x \to \hat B^+_x$ modulo $\pi$ yields $B^+_x/\pi \to B'^+_x/\pi \to \hat B^+_x/\pi = B^+_x/\pi$ and the composition is the identity, hence $B^+_x/\pi \to B'^+/\pi$ is injective. Taking $B^+_x \to B^+_x \tensor_{A^+_x} \hat A^+_x \surjto B'^+_x$ modulo $\pi$ yields $B^+_x/\pi \isoto (B^+_x \tensor_{A^+_x} \hat A^+_x)/\pi \surjto B'^+/\pi$ (the first isomorphism follows from $\hat A^+_x/\pi = A^+_x/\pi$), so that $B^+_x/\pi \to B'^+_x/\pi$ is surjective. In total we deduce that $B^+_x/\pi = B'^+_x/\pi$, so that it is enough to show that $B'^+_x/\pi$ is flat and finitely presented over $A^+_x/\pi$.

We claim that even $B'^+_x$ is flat and finitely presented over $\hat A^+_x$. Indeed, note that $\hat A^+_x = K^+$ where $x \in \pi_0(X)$ is of the form $\Spa(K, K^+)$ for some perfectoid field $K$ and an open and bounded valuation subring $K^+$. Clearly $\hat B^+_x$ is $\pi$-torsion free so that the same is true for $B'^+_x$, hence $B'^+_x$ is flat over the valuation ring $\hat A^+_x$. Moreover, since $B'^+_x$ is of finite type over $\hat A^+_x$ it follows from \cite[Lemma 053E]{stacks-project} that $B'^+_x$ is finitely presented over $\hat A^+_x$.
\end{proof}

\begin{lemma} \label{rslt:v-cover-of-tot-disc-is-mod-pi-descendable}
Let $f\colon Y \to X$ be a v-cover in $\AffPerfd_\pi$ with $X$ totally disconnected. Then $f$ is mod-$\pi$ weakly fs-descendable of index $\le 4$. In particular mod-$\pi$ modules descend along $f$.
\end{lemma}
\begin{proof}
Let $X = \Spa(A, A^+)$ and $Y = \Spa(B, B^+)$. Write $B^+ = \varinjlim_{i\in I} B^+_i$, where $B^+_i \subset B^+$ are $A^+$-subalgebras of finite type. By \cref{rslt:tors-free-ftype-over-tot-disc-implies-fppf-fibers,rslt:fppf-cover-of-rings-is-descendable,rslt:fs-descendable-can-be-checked-on-connected-components} $(A^+/\pi)_\solid \to (B^+_i/\pi)_\solid$ is fs-descendable of index $\le 2$ for all $i \in I$. Thus the claim follows from \cref{rslt:filtered-colim-of-bounded-pi-0-desc-is-weakly-pi-0-desc} (use also \cref{rslt:descendability-stable-under-almost-localization} to deduce the statement for \emph{almost} rings).
\end{proof}

With the descent result \cref{rslt:v-cover-of-tot-disc-is-mod-pi-descendable} we can now compute $\DqcohriX X$ for totally disconnected $X$. As a corollary we obtain an explicit description of $\DqcohriX X$ for general affinoid perfectoid spaces $X$.

\begin{proposition} \label{rslt:computation-of-Dqcohri-in-td-case}
Suppose $X = \Spa(A, A^+) \in \AffPerfd_\pi$ admits a quasi-pro-étale map to a totally disconnected space. Then there is a natural equivalence
\begin{align*}
	\DqcohriX X = \Dqcohri(A^+/\pi).
\end{align*}
\end{proposition}
\begin{proof}
Let $\mathcal C$ be the pro-étale site of strictly totally disconnected spaces. Note that $\mathcal C$ admits fiber products because if $Y \to X$ is a pro-étale map and $X$ is strictly totally disconnected then $Y$ is also strictly totally disconnected, see \cite[Lemma 7.19]{etale-cohomology-of-diamonds}. It follows immediately from \cref{rslt:v-cover-of-tot-disc-is-mod-pi-descendable} (and \cref{rslt:base-change-for-affinoid-perfectoid}) that the functor $\mathcal C^\opp \to \infcatinf$, $\Spa(A, A^+) \mapsto \Dqcohri(A^+/\pi)$ is a sheaf of $\infty$-categories on $\mathcal C$. Now, from the sheafiness of $\Dqcohri(A^+/\pi)$ on $\mathcal C$ it follows that $\DqcohriX X = \Dqcohri(A^+/\pi)$ for every $X = \Spa(A, A^+)$ in $\mathcal C$, i.e. for every strictly totally disconnected $X$ (see \cref{rslt:sheaves-on-basis-equiv-sheaves-on-whole-site}).

Finally, let $X = \Spa(A, A^+) \in \AffPerfd_\pi$ such that there is a quasi-pro-étale map $X \to Z$ to some totally disconnected space $Z$. Choose a pro-étale cover $Y \to Z$ such that $Y$ is strictly totally disconnected (cf. \cite[Lemma 7.18]{etale-cohomology-of-diamonds}) and let $Y_\bullet \to Z$ be the associated Čech nerve. Let $Y_\bullet' := Y_\bullet \cprod_Z X$, which is a Čech cover of $X$ by strictly totally disconnected spaces $Y_n' = \Spa(B'^n, B'^{n+})$. Thus $\DqcohriX{Y'_n} = \Dqcohri(B'^{n+}/\pi)$ and in particular $\DqcohriX X = \varprojlim_{n\in\Delta} \Dqcohri(B'^{n+}/\pi)$. By \cref{rslt:v-cover-of-tot-disc-is-mod-pi-descendable,rslt:stability-of-pi-0-weakly-descendable-maps} (and \cref{rslt:base-change-for-affinoid-perfectoid}) we get that $(A^+/\pi)^a_\solid \to (B'^{0+}/\pi)^a_\solid$ is weakly fs-descendable and hence (using again \cref{rslt:base-change-for-affinoid-perfectoid}) we get $\DqcohriX X = \Dqcohri(A^+/\pi)$, as desired.
\end{proof}

\begin{corollary} \label{rslt:explicit-description-of-DqcohriX}
Let $X \in \AffPerfd_\pi$. Choose any quasi-pro-étale cover $Y \to X$ by a totally disconnected space $Y$, let $Y_\bullet$ be the associated Čech nerve, and write $Y_n = \Spa(B^n, B^{n+})$ for all $n$. Then there is a natural equivalence
\begin{align*}
	\DqcohriX X = \varprojlim_{n \in \Delta} \Dqcohri(B^{n+}/\pi).
\end{align*}
In particular $\DqcohriX X$ is a stable $\infty$-category.
\end{corollary}
\begin{proof}
The claim about stability follows from the first claim using \cite[Theorem 1.1.4.4]{lurie-higher-algebra}. The first claim follows directly from the definition of sheaves and \cref{rslt:computation-of-Dqcohri-in-td-case}, noting that all $Y_n$ are quasi-pro-étale over the totally disconnected space $Y$.
\end{proof}

Our next goal is to extend \cref{rslt:computation-of-Dqcohri-in-td-case} to a more general class of affinoid perfectoid spaces $X$, namely those which are ``weakly of perfectly finite type'' (see \cref{def:perfectly-finite-type-in-AffPerf} below) over some totally disconnected space. This result, see \cref{rslt:compute-Dqcohri-for-fin-type-over-tot-disc} below, is crucial for computing $\DqcohriX X$ on rigid spaces $X$ later on. The proof of \cref{rslt:compute-Dqcohri-for-fin-type-over-tot-disc} proceeds in two steps, captured by \cref{rslt:compute-Dqcohri-for-subset-of-Z'-over-Z} and \cref{rslt:fin-type-over-K-implies-mod-pi-desc}.

We start with \cref{rslt:compute-Dqcohri-for-subset-of-Z'-over-Z}, which essentially computes $\DqcohriX X$ for a somewhat exotic class of affinoid perfectoid spaces: Those $X$ which come as a relative compactification of a map of totally disconnected spaces of finite $\dimtrg$ (see \cite[Definition 21.5]{etale-cohomology-of-diamonds}). While this class of affinoid perfectoid spaces may not be of much interest in applications, it plays an important role in setting up the theory.

\begin{proposition} \label{rslt:compute-Dqcohri-for-subset-of-Z'-over-Z}
Let $f\colon Z' \to Z$ be a map of totally disconnected spaces in $\AffPerfd_\pi$ with $\dimtrg f < \infty$. Then there is a pro-étale cover $Y \to \overline Z'^{/Z}$ by a totally disconnected space $Y$ such that this cover is mod-$\pi$ weakly fs-descendable of index bounded by a constant which only depends on $\dimtrg f$.
\end{proposition}
\begin{proof}
Let $X := \overline{Z'}^{/Z}$. By \cite[Lemma 7.13]{etale-cohomology-of-diamonds} there is a pro-étale cover $Y = \varprojlim_{i\in I} Y_i \surjto X$ such that $Y$ is totally disconnected and each $Y_i$ is a disjoint union of rational open subsets of $X$. By \cref{rslt:filtered-colim-of-bounded-pi-0-desc-is-weakly-pi-0-desc} it is enough to show that there is some $d \ge 0$ such that for every $i \in I$, the map $Y_i \to X$ is mod-$\pi$ fs-descendable of index $\le d$. By definition of fs-descendability (and using \cref{rslt:fs-descendable-can-be-checked-on-connected-components} and $\pi_0(X) = \pi_0(\Spec A^+/\pi)$) we can reduce to the case that $Z$ and $Z'$ are connected. In the following we will show that $d := 2 c(\dimtrg f + 1)$ works, where $c(n)$ is the constant from \cref{rslt:steady-covering-is-descendable}.

We are now in the situation that $Z = \Spa(K, K^+)$ and $Z' = \Spa(K', K'^+)$ for some perfectoid fields $K$ and $K'$ and some open and bounded valuation subrings $K^+ \subset K$ and $K'^+ \subset K'$, so in particular $X = \overline Z^{/Z} = \Spa(K', A^+)$, where $A^+$ is the integral closure of $K^+ + \mm_{K'}$ in $K'$. Moreover, we are given a map $U = Y_i \to X$ such that $U = \bigdunion_{j=1}^n U_j$ for some cover $X = \bigunion_{j=1}^n U_j$ by rational open subsets $U_j \subset X$ and we want to show that the map $U \to X$ is mod-$\pi$ fs-descendable of index $\le d$. Thus we are roughly in the situation of \cref{rslt:descendability-of-open-cover-of-ZR-space}, but we cannot apply this result directly because we are working with mod-$\pi$ coefficients. In the following we adapt the proof of \cref{rslt:descendability-of-open-cover-of-ZR-space} to the present setting. We also make implicit use of \cref{rslt:descendability-stable-under-almost-localization} to get the same statements for \emph{almost} rings.

The space $X = \Spa(K', A^+)$ consists of all open valuation rings $V \subset \ri_{K'}$ which contain $K^+$. Modding out by the maximal ideal $\mm_{K'}$ of topologically nilpotent elements in $K'$ we obtain a canonical homeomorphism $\abs{\Spa(K', K^+)} \isoto \abs{\Spa(k', k^+)}$, where $k' = \ri_{K'}/\mm_{K'}$ and $k^+ = A^+/\mm_{K'} = K^+/\mm_K$. Now $k'$ is a discrete field and $k^+$ is a valuation ring such that $k'$ has transcendence degree $\le \dimtrg f$ over the fraction field $k = \ri_K/\mm_K$ of $k^+$, so we are exactly in the situation of \cref{rslt:descendability-of-open-cover-of-ZR-space}. By the proof of that result there exists an integral proper scheme $\overline S$ over $k^+$ which is covered by $d+1$ affines $\overline S = \bigunion_{l=1}^{d+1} \overline S_l$, $\overline S_l = \Spec \overline R_l$ and a morphism of discrete adic spaces $\Spa(k', k^+) \to \overline S$ such that the open subsets $\abs{U_j} \subset \abs{\Spa(k', k^+)}$ come via base-change from open subsets $\overline W_j \subset \overline S$.

As is true for any affinoid perfectoid space, the inclusion $A^+ \injto \ri_{K'}$ is an almost isomorphism. This implies
\begin{align*}
	(A^+/\pi)^a_\solid =  ((\ri_{K'}/\pi)^a, A^+/\pi)_\solid \end{align*}
The projection map $\rho\colon (\ri_{K'}/\pi, A^+/\pi) \to (k', k^+)$ of discrete adic rings induces a homeomorphism $\abs X \isom \abs{\Spa(k', k^+)} \isoto \abs{\Spa(\ri_{K'}/\pi, A^+/\pi)}$ because we are quotienting by nilpotent elements. Thus the covering $(U_j)_j$ of $X$ induces a qcqs open covering $\overline U_i$ of $\Spa(\ri_{K'}/\pi, A^+/\pi)$ and we have to show that the induced map $\bigdunion_{j=1}^n \overline U_j^\an \to \Spa(\ri_{K'}/\pi, A^+/\pi)^\an$ is descendable of index $\le d$.

For $l = 1, \dots, d+1$ we have compatible injections $\overline R_l \subset k'$ given by the dominant point $x_{\overline S} \in \overline S(k')$ (see the construction in the proof of \cref{rslt:descendability-of-open-cover-of-ZR-space}). Let $R_l \subset \ri_{K'}/\pi$ be the preimage of $R_l$ under the projection $\ri_{K'}/\pi \surjto k'$. The glueing datum for $(\overline S_l)_l$ lifts to a glueing datum for $(S_l = \Spec R_l)_l$ (because the former glueing datum is given by inclusions of subrings of $k'$), so we get a scheme $S = \bigunion_l S_l$ with $\abs S = \abs{\overline S}$. The morphism $\Spa(k', k^+) \to \overline S$ lifts to a morphism $\Spa(\ri_{K'}/\pi, A^+/\pi) \to S$ which is the same map on the underlying topological spaces. But then the affine covering $(\overline U_j)_j$ of $\Spa(\ri_{K'}/\pi, A^+/\pi)$ comes via base-change from an affine covering $(W_j)_j$ of $S$, so the claim follows from \cref{rslt:descendable-stable-under-base-change} and \cref{rslt:fppf-descendability-for-schemes} (as in the proof of \cref{rslt:descendability-of-open-cover-of-ZR-space}).
\end{proof}

\begin{remark} \label{rslt:Y-over-Z'-Z-has-mod-pi-Tor-dim-1}
It follows from the proof of \cref{rslt:compute-Dqcohri-for-subset-of-Z'-over-Z} that in the situation of that result the constructed map $Y \to \overline{Z'}^{/Z}$ has mod-$\pi$ Tor dimension $\le 1$. Namely, this reduces to showing the same for each $Y_i \to \overline{Z'}^{/Z}$, which by the proof is mod-$\pi$ a base-change of an open affine covering $(W_j)_j$ of some scheme $S$, hence the claim follows from \cref{rslt:open-immersion-of-affine-schemes-has-Tor-dim-1} (use the same arguments as in the proof of that result to reduce to the case that $S$ is affine).
\end{remark}

We now come to the second main step in the proof of \cref{rslt:compute-Dqcohri-for-fin-type-over-tot-disc}: In ``small'' cases, we will show a somewhat automatic descendability result using a global dimension argument (see \cref{rslt:fin-type-over-K-implies-mod-pi-desc} below).

\begin{definition} \label{def:perfectly-finite-type-in-AffPerf}
A morphism $\Spa(B, B^+) \to \Spa(A, A^+)$ in $\AffPerfd$ is called \emph{of weakly perfectly finite type} if for some $n \ge 0$ there is a surjective open continuous $A$-algebra map $A\langle T_1^{1/p^\infty}, \dots, T_n^{1/p^\infty}\rangle \surjto B$.
\end{definition}

\begin{lemma} \label{rslt:polynomial-over-valuation-ring-has-finite-gldim}
Let $V$ be a valuation ring of finite (Krull) dimension and let $n \ge 0$. Then the polynomial ring $V[T_1, \dots, T_n]$ has global dimension $\le n + 2$.
\end{lemma}
\begin{proof}
The case $n = 0$ is found in \cite[Theorem 5.1]{datta-vanishing-results-of-valuation-rings}. The general case then follows from \cite[Theorem 4.3.7]{weibel-homological-algebra}.
\end{proof}

\begin{corollary} \label{rslt:perfect-polynomial-ring-over-valuation-ring-has-finite-gldim}
Let $V$ be a valuation ring of finite (Krull) dimension and let $n \ge 0$. Then the perfect polynomial ring $V[T_1^{1/p^\infty}, \dots, T_n^{1/p^\infty}]$ has global dimension $\le n + 3$.
\end{corollary}
\begin{proof}
This is essentially the same proof as in \cite[Proposition 11.31]{bhatt-scholze-witt}. Namely, for all $i = 0, 1, \dots, \infty$ let $R_i := V[T_1^{1/p^i}, \dots, T_n^{1/p^i}]$. Given an $R_\infty$-module $M$ denote $M_i := M \tensor_{R_i} R_\infty$. Since each $R_i$ has global dimension $\le n + 2$ by \cref{rslt:polynomial-over-valuation-ring-has-finite-gldim} and $R_i \to R_\infty$ is flat it follows that all $M_i$ have projective dimension $\le n + 3$ (as $R_\infty$-modules). But also $M = \varinjlim_i M_i$, which implies that $M$ has projective dimension $\le n + 3$ (write this colimit as a cofiber of direct sums).
\end{proof}

\begin{proposition} \label{rslt:fin-type-over-K-implies-mod-pi-desc}
Let $Z \in \AffPerfd_\pi$ be totally disconnected and let $X \to Z$ be a map in $\AffPerfd_\pi$ which is of weakly perfectly finite type. Then for any v-cover $Y \surjto X$ the map $\overline Y^{/X} \to X$ is mod-$\pi$ fs-descendable of index bounded by a constant which only depends on the number of perfect generators for $X$ over $Z$.
\end{proposition}
\begin{proof}
We claim that the index of fs-descendability is $\le g(n + 3)$, where $g(-)$ is the constant from \cref{rslt:descendability-for-fin-global-dimension} and $n$ is the number of perfect generators needed for $X \to Z$. By definition of fs-descendability (and \cref{rslt:fs-descendable-can-be-checked-on-connected-components} together with $\pi_0(Z) = \pi_0(\ri^+_Z(Z)/\pi)$) we can reduce to the case that $Z$ is connected, i.e. $Z = \Spa(K, K^+)$ for some perfectoid field $K$ and an open and bounded valuation subring $K^+ \subset K$. Write $X = \Spa(A, A^+)$ and $Y = \Spa(B, B^+)$ so that $\overline Y^{/X} = \Spa(B, B'^+)$, where $B'^+$ is the integral closure of $A^+ + \mm_B$ in $B$. In particular, $B'^+/\pi$ is the integral closure of $A^+/\pi$ in $B^+/\pi$ (using that $\mm_B/\pi$ consists only of nilpotent and hence integral elements), so we have $(B'^+/\pi)^a_\solid = (B^{+a}/\pi, A^+/\pi)_\solid$. It is now enough to show that the map $(A^+/\pi)^a_\solid \to (B^{+a}/\pi, A^+/\pi)_\solid$ is descendable of index $\le g(n + 3)$. This is essentially an application of \cref{rslt:descendability-for-fin-global-dimension}, but the details are a bit subtle.

First note that by \cref{rslt:scheme-classical-desc-iff-canonical-compact-desc} we only need to show that the map $A^{+a}/\pi \to B^{+a}/\pi$ is discretely descendable of index $\le g(n+3)$, which is a property that only depends on the almost rings $A^{+a}/\pi$ and $B^{+a}/\pi$ (and not on $A^+/\pi$ and $B^+/\pi$). For every $k \ge 1$ let $\Spa(B^k, B'^{k+})$ be the $k$-fold fiber product of $\overline Y^{/X}$ over $X$. Then we can view $B'^{\bullet+,a}$ as a cosimplicial object in $\D(A^{+a})$. By \cite[Proposition 8.8]{etale-cohomology-of-diamonds} the natural map $A^{+a} \to \Tot(B'^{\bullet+,a})$ is an isomorphism in $\D(A^{+a})$. Now let $A' := K\langle T_1^{1/p^\infty}, \dots, T_n^{1/p^\infty}\rangle$, so that $A'^\circ = \ri_K \langle T_1^{1/p^\infty}, \dots, T_n^{1/p^\infty}\rangle$. By assumption there exists an open surjective continuous $K$-algebra map $f\colon A' \surjto A$. We can w.l.o.g. assume that $K^+ = \ri_K$ and that $A^+$ is the integral closure of $f(A'^\circ) \subset A$, as these changes leave the associated almost rings invariant. We further denote $A'^\circ_0 := \ri_K [T_1^{1/p^\infty}, \dots, T_n^{1/p^\infty}] \subset A'^\circ$. Applying the forgetful functor $\D(A^{+a}) \to \D(A'^{\circ a}_0)$ to the isomorphism $A^{+a} \isoto \Tot(B'^{\bullet+,a})$ yields the same isomorphism in $\D(A'^{\circ a}_0)$ (where we now compute the totalization in $\D(A'^{\circ a}_0)$).

By \cref{rslt:perfect-polynomial-ring-over-valuation-ring-has-finite-gldim} $A'^{\circ a}_0$ has global dimension $\le n + 3$ and hence the proof of \cref{rslt:descendability-for-fin-global-dimension} shows that the map $A^{+a} \to \Tot_{2n+7}(B'^{\bullet+,a})$ splits in $\D(A'^{+a}_0)$. Applying the functor $- \tensor_{A'^+_0} A'^+_0/\pi$ and using that $A'^\circ_0/\pi = A'^\circ/\pi$ we deduce that the map $A^{+a}/\pi \to \Tot_{2n+7}(B'^{\bullet+,a}/\pi)$ splits in $\D(A'^{\circ a}/\pi)$. In order to get this splitting also in $\D(A^{+a}/\pi)$ it is enough to see that
\begin{align}
	A^{+a}/\pi \tensor_{A'^{+a}/\pi} A^{+a}/\pi = A^{+a}/\pi, \qquad B'^{\bullet+,a}/\pi \tensor_{A'^{+a}/\pi} A^{+a}/\pi = B'^{\bullet+,a}/\pi, \label{eq:fin-type-over-K-mod-pi-descent-helper}
\end{align}
as then we can simply apply the functor $- \tensor_{A'^{+a}/\pi} A^{+a}/\pi$. Assuming this holds, we can finish by argueing in the same way as in \cref{rslt:descendability-for-fin-global-dimension}.

Thus it only remains to prove \cref{eq:fin-type-over-K-mod-pi-descent-helper}. Denote $X' := \Spa(A', A'^\circ)$, so that we have a map $X \to X'$. In fact this map is a Zariski closed immersion, so we have $X \cprod_{X'} X = X$ and $\overline Y^{/X} \cprod_{X'} X = \overline Y^{/X}$ (for example, using \cite[Remark 7.9]{etale-cohomology-of-diamonds} we can write $X$ as an intersection of rational open subsets of $X'$ and then use that fiber products are limits and hence commute with such an intersection). But then \cref{rslt:base-change-for-affinoid-perfectoid} immediately implies \cref{eq:fin-type-over-K-mod-pi-descent-helper}, as desired.
\end{proof}

We are finally in the position to prove our first main result, i.e. to compute $\DqcohriX X$ for ``small'' spaces $X \in \AffPerfd_\pi$:

\begin{theorem} \label{rslt:compute-Dqcohri-for-fin-type-over-tot-disc}
Suppose $X = \Spa(A, A^+) \in \AffPerfd_\pi$ is of weakly perfectly finite type over some totally disconnected space. Then there is a natural equivalence
\begin{align*}
	\DqcohriX X = \Dqcohri(A^+/\pi).
\end{align*}
\end{theorem}
\begin{proof}
Let $Z \in \AffPerfd_\pi$ be totally disconnected such that there is a map $X \to Z$ of weakly perfectly finite type. Let $n$ be the number of perfect generators needed for the morphism $X \to Z$. Choose any pro-étale map $Z' \surjto X$ with $Z'$ strictly totally disconnected and let $Z'_\bullet \to X$ denote the associated Čech nerve. Then all $Z'_k$ are (strictly) totally disconnected and we have $\dimtrg(Z' \to Z) \le n$ for all $k$. It follows from \cref{rslt:compute-Dqcohri-for-subset-of-Z'-over-Z} that each $\overline{Z'_k}^{/Z}$ admits a mod-$\pi$ weakly fs-descendable pro-étale cover by some totally disconnected space. The same is then true for each $Y_k := \overline{Z'_k}^{/X}$ because $Y_k$ is an intersection of qcqs open subsets of $\overline{Z'_k}^{/Z}$ and mod-$\pi$ weakly fs-descendable maps are preserved under base-change (use also that an intersection of qcqs open subsets of a totally disconnected space is still totally disconnected, see \cite[Lemma 7.6]{etale-cohomology-of-diamonds}). It thus follows from \cref{rslt:explicit-description-of-DqcohriX} that, writing $Y_k = \Spa(B^k, B^{k+})$, we have $\DqcohriX{Y_k} = \Dqcohri(B^+/\pi)$. Note furthermore that $Y_k$ is the $k$-fold fiber product of $Y_0$ over $X$. Since $\DqcohriX{(-)}$ is a pro-étale sheaf it follows that
\begin{align*}
	\DqcohriX X = \varprojlim_{k\in\Delta} \DqcohriX{Y_k} = \varprojlim_{k\in\Delta} \Dqcohri(B^{+k}/\pi).
\end{align*}
This means that we need to show that mod-$\pi$ modules descend along the map $Y_0 = \overline Z'^{/X} \to X$. But this follows from \cref{rslt:fin-type-over-K-implies-mod-pi-desc}.
\end{proof}

\begin{remark} \label{rslt:compute-Dqcohri-for-qproet-over-fin-type-over-tot-disc}
The proof of \cref{rslt:compute-Dqcohri-for-fin-type-over-tot-disc} can be generalized slightly: Suppose $X' = \Spa(A', A'^+) \in \AffPerfd_\pi$ is quasi-pro-étale over some $X = \Spa(A, A^+) \in \AffPerfd_\pi$ which is as in \cref{rslt:compute-Dqcohri-for-fin-type-over-tot-disc}; then we still have $\DqcohriX{X'} = \Dqcohri(A'^+/\pi)$. Namely, in the proof of \cref{rslt:compute-Dqcohri-for-fin-type-over-tot-disc} we can replace $Y_k$ by $Y'_k := Y_k \cprod_X X'$. Then $Y'_k = \Spa(B'^k, B'^{+k})$ is quasi-pro-étale over the totally disconnected space $Y_k$, so we still have $\DqcohriX{Y'_k} = \Dqcohri(B'^{+k}/\pi)$. Also, $Y'_0 \to X'$ is the base-change of $Y_0 \to X$ and hence still mod-$\pi$ fs-descendable (using \cref{rslt:stability-of-pi-0-weakly-descendable-maps}). An even more general version of \cref{rslt:compute-Dqcohri-for-fin-type-over-tot-disc} is found in \cref{rslt:compute-Dqcohri-for-p-bounded-over-tot-disc} below.
\end{remark}

Having computed $\DqcohriX X$ on a large class of affinoid perfectoid spaces $X$, we will now show that this $\infty$-category satisfies v-hyperdescent, i.e. the functor $X \mapsto \DqcohriX X$ is a hypercomplete v-sheaf. This will allow us to get a good theory of quasicoherent solid $\ri^{+a}/\pi$-sheaves on any small v-stack later on and it also provides a good amount of flexibility in proofs.

As a first step we will employ our theory on valuation rings (see \cref{sec:andesc.valrings}) to prove a very useful result on Tor dimensions, see \cref{rslt:mod-pi-Tor-dim-1-over-tot-disc-space} below. This result will be crucial for proving v-descent.

\begin{lemma} \label{rslt:any-map-in-AffPerf-is-pro-wh}
Any map $Y \to X$ in $\AffPerfd_\pi$ can be written as a cofiltered limit $Y = \varprojlim Y_i \to X$ such that all $Y_i \to X$ are of weakly perfectly finite type. Moreover, writing $Y = \Spa(B, B^+)$ and $Y_i = \Spa(B^i, B^{i+})$ we can arrange that the induced maps $B^i \to B$ are injective.
\end{lemma}
\begin{proof}
By \cite[Corollary 3.20]{etale-cohomology-of-diamonds} we can reduce to the case of characteristic $p$. More precisely, to see that weakly perfectly finite-type morphisms correspond to each other, note that a map of perfectoid spaces is of weakly perfectly finite type if and only if it factors as a closed immersion followed by a compactified relative perfectoid ball, and both operations are preserved under tilts and untilts (for Zariski closed immersions see e.g. \cite[Theorem 5.8]{etale-cohomology-of-diamonds}).

Now let $X = \Spa(A, A^+)$ be of characteristic $p$. For every finite subset $J \subset B^+$ let $B^J \subset B$ be the $A$-subalgebra given as the image of $A\langle T_j^{1/p^\infty} \setst j \in J \rangle \to B$ mapping $T_j$ to $j$. We equip $B^J$ with the quotient topology. Then $B^J$ is a perfect and complete Tate ring, hence uniform by \cite[Proposition 3.5]{etale-cohomology-of-diamonds} and in particular perfectoid. Let $B^{J+} := B^{J\circ} \isect B^+$, where the intersection is taken inside $B$. Then $B^{J+} \subset B^{J\circ}$ is open and integrally closed and hence $Y_J := \Spa(B^J, B^{J+})$ is an affinoid perfectoid space. Now let $I$ be the cofiltered set of finite subsets $J \subset B^+$, so that we get a cofiltered system $(Y_i)_{i\in I}$ of affinoid perfectoid spaces. Clearly each $Y_i \to X$ is of weakly perfectly finite type and each $B^i \to B$ is injective.

It remains to show that $Y \to \varprojlim_i Y_i$ is an isomorphism. This boils down to showing that the map $\varinjlim_i B^{i+}/\pi^n \to B^+/\pi^n$ is an isomorphism for all $n \ge 1$. It is clearly surjective. To show that it is injective, suppose we have $f_1 - f_2 = \pi^n g$ for some $f_1, f_2 \in B^{J+}$ and some $g \in B^+$; then $g \in B^{J'+}$ for any $J' \supset J$ with $g \in J'$, so that $f_1 - f_2 \in \pi^n B^{J'+}$.
\end{proof}

\begin{proposition} \label{rslt:mod-pi-Tor-dim-1-over-tot-disc-space}
Let $Y = \Spa(B, B^+) \to X = \Spa(A, A^+)$ be a map in $\AffPerfd_\pi$ and assume that $X$ is totally disconnected. Then $(A^+/\pi)^a_\solid \to (B^+/\pi)^a_\solid$ has Tor dimension $\le 1$.
\end{proposition}
\begin{proof}
Since Tor dimension is preserved under filtered colimits, by \cref{rslt:any-map-in-AffPerf-is-pro-wh} we reduce to the case that $Y \to X$ is of weakly perfectly finite type. Then by \cref{rslt:compute-Dqcohri-for-fin-type-over-tot-disc} we have $\DqcohriX Y = \Dqcohri(B^+/\pi)$. In other words, choosing any pro-étale Čech covering $\tilde Y_\bullet \to Y$ by strictly totally disconnected spaces $\tilde Y_n = \Spa(\tilde B^n, \tilde B^{n+})$, we get for every $M \in \Dqcohrile0{A^+/\pi}$ that
\begin{align*}
	M \tensor_{(A^+/\pi)^a_\solid} (B^+/\pi)^a_\solid = \varprojlim_{n\in\Delta} M \tensor_{(A^+/\pi)^a_\solid} (\tilde B^{n+}/\pi)^a_\solid.
\end{align*}
Since the totalization on the right preserves $\Dqcohrile1(A^+/\pi)$ it is enough to show that $(A^+/\pi)^a_\solid \to (\tilde B^{n+}/\pi)^a_\solid$ has Tor dimension $\le 1$ for all $n$. In other words, we can replace $Y$ by $\tilde Y_n$ and hence assume that $Y$ is totally disconnected.

By \cref{rslt:compare-catsldmod-with-sheaves-on-pi-0} it is enough to show that for every connected component $y = \Spec B^+_y/\pi$ of $\Spec B^+/\pi$, $(A^+/\pi)^a_\solid \to (B^+_y/\pi)^a_\solid$ has Tor dimension $\le 1$. Since $\pi_0(\Spec B^+/\pi) = \pi_0(Y)$ we are therefore reduced to the case that $Y$ is connected, i.e. $Y = \Spa(K', K'^+)$ for some perfectoid field $K'$ and some open and bounded valuation subring $K'^+ \subset K'$. Then the map $Y \to X$ factors over some connected component $x$ of $X$ and since the map $x \to X$ is mod-$\pi$ flat (again by \cref{rslt:compare-catsldmod-with-sheaves-on-pi-0}) we can reduce to the case that $X = x$, i.e. $X$ is connected and hence of the form $X = \Spa(K, K^+)$ for some perfectoid field $K$ and an open and bounded valuation subring $K^+$. Then $- \tensor_{(A^+/\pi)^a_\solid} (B^+/\pi)^a_\solid = - \tensor_{(K^+)^a_\solid} (K'^+)^a_\solid$ (where we treat $K^+$ and $K'^+$ as discrete rings and use \cref{rslt:solid-base-change-for-discrete-rings}), so it is enough to show that $(K^+)^a_\solid \to (K'^+)^a_\solid$ has Tor dimension $\le 1$. But this follows immediately from \cref{rslt:flat-map-of-val-rings-has-sld-Tor-dim-1}.
\end{proof}

\begin{remark}
One should see \cref{rslt:mod-pi-Tor-dim-1-over-tot-disc-space} as a variant of \cite[Proposition 7.23]{etale-cohomology-of-diamonds}: Instead of flatness we only get Tor dimension $\le 1$, but now it holds for the associated \emph{analytic} rings, which is a much stronger statement. See also \cref{rslt:mod-pi-Tor-dim-for-fin-type-over-tot-disc} for a generalization (where we allow more general spaces $X$).
\end{remark}

One of the main challenges of proving v-descent is that we have very little control over spaces of the form $Z = Y \cprod_X Y$ for (strictly) totally disconnected spaces $X$ and $Y$. In particular we cannot compute the associated category $\DqcohriX Z$, which is probably not given by $\Dqcohri(\ri^+_Z(Z)/\pi)$ in general. However, using \cref{rslt:mod-pi-Tor-dim-1-over-tot-disc-space} we can work around that issue, as follows.

\begin{lemma} \label{rslt:general-etale-descent-of-Dqcohri}
Let $Y \to X$ be an étale cover in $\AffPerfd_\pi$ with associated Čech cover $Y_\bullet \to X$ and write $X = \Spa(A, A^+)$ and $Y_n = \Spa(B^n, B^{n+})$. Then for every $M \in \Dqcohri(A^+/\pi)$ the natural map
\begin{align*}
	M \isoto \Tot(M \tensor_{(A^+/\pi)^a_\solid} (B^{\bullet+}/\pi)^a_\solid)
\end{align*}
is an isomorphism.
\end{lemma}
\begin{proof}
First consider the case of characteristic $p$. Let $K$ be the $\pi$-adic completion of $\Fld_p((\pi))(\pi^{1/p^\infty})$. By \cref{rslt:any-map-in-AffPerf-is-pro-wh} we can write $X$ as a cofiltered limit $X = \varprojlim_{i\in I} X_i$ where all $X_i$ are weakly of perfectly finite type over $\Spa(K, K^\circ)$. Then by \cite[Lemma 6.4.(iv)]{etale-cohomology-of-diamonds} there is some $i \in I$ and an étale cover $Y_i \to X_i$ in $\AffPerfd_\pi$ such that $Y = Y_i \cprod_{X_i} X$. Using \cref{rslt:base-change-for-affinoid-perfectoid} we easily reduce to the case $X = X_i$ and $Y = Y_i$. But then \cref{rslt:compute-Dqcohri-for-fin-type-over-tot-disc} implies $\Dqcohri(A^+/\pi) = \DqcohriX X$ and $\Dqcohri(B^{n+}/\pi) = \DqcohriX{Y_n}$, so the claim follows from étale descent of the sheaf $\DqcohriX{(-)}$.

Now consider the general case (i.e. not in characteristic $p$). If $\pi \divides p$ then we are reduced to characteristic $p$ (as in the proof of \cref{rslt:base-change-for-affinoid-perfectoid}). In general choose a pseudouniformizer $\varpi \in A$ with $\varpi \divides p$ and choose $m \ge 1$ such that $\pi \divides \varpi^m$. Then
\begin{align*}
	M \tensor_{(A^+/\pi)^a_\solid} (B^{n+}/\pi)^a_\solid = M \tensor_{(A^+/\varpi^m)^a_\solid} (B^{n+}/\varpi^m)^a_\solid
\end{align*}
so we reduce to the case $\pi = \varpi^m$. But then we can argue by induction by tensoring with the exact sequence $0 \to A^+/\varpi^{k-1} \to A^+/\varpi^k \to A^+/\varpi \to 0$ (using that $\Tot$ preserves cofiber sequences).
\end{proof}

\begin{lemma} \label{rslt:base-change-over-tot-disc-equiv}
Let $Y \to X$ be a map in $\AffPerfd_\pi$ such that $X$ is totally disconnected. Write $X = \Spa(A, A^+)$ and $Y = \Spa(B, B^+)$ and let $\mathcal C_1 \subset \Dqcohri(B^+/\pi)$ and $\mathcal C_2 \subset \DqcohriX Y$ be the essential images of the base-change functors $\Dqcohri(A^+/\pi) \to \Dqcohri(B^+/\pi)$ and $\DqcohriX X \to \DqcohriX Y$. Then the natural functor $\Dqcohri(B^+/\pi) \to \DqcohriX Y$ restricts to an equivalence $\mathcal C_1 \isoto \mathcal C_2$.
\end{lemma}
\begin{proof}
We have $\DqcohriX X = \Dqcohri(A^+/\pi)$ by \cref{rslt:computation-of-Dqcohri-in-td-case}, hence the functor $\DqcohriX X \to \DqcohriX Y$ factors over the functor $\Dqcohri(A^+/\pi) \to \Dqcohri(B^+/\pi)$. In particular the functor $\mathcal C_1 \to \mathcal C_2$ exists and is essentially surjective.

To show full faithfulness we note that the natural functor $\widetilde{(-)}\colon \Dqcohri(B^+/\pi) \to \DqcohriX Y$ has a right adjoint $\Gamma(Y, -)$. It can be explicitly described as follows: Choose any Čech cover $\tilde Y_\bullet \to Y$ of $Y$ by totally disconnected spaces $\tilde Y_n = \Spa(\tilde B^n, \tilde B^{n+})$. Then given $\mathcal M = M^\bullet \in \DqcohriX Y = \varprojlim_{n\in\Delta} \Dqcohri(\tilde B^{n+}/\pi)$ we have $\Gamma(Y, \mathcal M) = \Tot(M^\bullet)$ (where we view $M^\bullet$ as a cosimplicial object of $\Dqcohri(B^+/\pi)$ via the forgetful functors). It is therefore enough to show that the unit of the adjunction $M \to \Gamma(Y, \widetilde M)$ is an isomorphism for all $M \in \mathcal C_1$.

Now fix $M \in \mathcal C_1$, i.e. $M = M_0 \tensor_{(A^+/\pi)^a_\solid} (B^+/\pi)^a_\solid$ for some $M_0 \in \Dqcohri(A^+/\pi)$. Let $M^\bullet = M \tensor_{(B^+/\pi)^a_\solid} (\tilde B^{\bullet+}/\pi)^a_\solid$. We want to show that the map $M \to \Tot(M^\bullet)$ is an isomorphism. Note that we can write $M^\bullet = M_0 \tensor_{(A^+/\pi)^a_\solid} (\tilde B^{\bullet+}/\pi)^a_\solid$. By the finite Tor dimension of $- \tensor_{(A^+/\pi)^a_\solid} (\tilde B^{n+}/\pi)^a_\solid$ (see \cref{rslt:mod-pi-Tor-dim-1-over-tot-disc-space}) these base-change functors commute with Postnikov limits in $M_0$. As $\Tot(-)$ commutes with arbitrary limits, we can replace $M_0$ by $\tau_{\le k} M_0$ and hence reduce to the case $M_0 \in \Dqcohrile{k}(A^+/\pi)$. Without loss of generality we can assume $k = 0$.

By \cref{rslt:mod-pi-Tor-dim-1-over-tot-disc-space} we have $M^n \in \Dqcohrile1(\tilde B^{n+}/\pi)$ for all $n$ (this is a crucial application of our results on valuation rings), so that $\pi_k \Tot(M^\bullet) = \pi_k \Tot_n(M^\bullet)$ for all $k > -n+1$ (see \cite[Proposition 1.2.4.5.(5)]{lurie-higher-algebra}). On the other hand, we can write $\tilde Y_0 = \varprojlim_{i\in I} \tilde Y^i_0$ for étale covers $\tilde Y^i_0 \to Y$. Let $\tilde Y^i_\bullet \to Y$ be the associated Čech covers, write $\tilde Y^i_n = \Spa(\tilde B^n_i, \tilde B^{n+}_i)$ and let $M_i^\bullet = M \tensor_{(B^+/\pi)^a_\solid} (\tilde B^{\bullet+}_i/\pi)^a_\solid$. By \cref{rslt:general-etale-descent-of-Dqcohri} we know that $M = \Tot(M^\bullet_i)$ for all $i$. By the same reasoning as above we still have $\pi_k \Tot(M^\bullet_i) = \pi_k \Tot_n(M^\bullet_i)$ for $k > -n+1$. But $\Tot_n(-)$ is a finite limit and hence commutes with arbitrary colimits in a stable $\infty$-category, so we get
\begin{align*}
	\pi_k M &= \varinjlim_i \pi_k \Tot(M_i^\bullet) = \varinjlim_i \pi_k \Tot_n(M_i^\bullet) = \pi_k \Tot_n(\varinjlim M_i^\bullet) = \pi_k \Tot_n(M^\bullet)\\
	&= \pi_k \Tot(M^\bullet)
\end{align*}
for $k > -n+1$. Since $n$ can be chosen arbitrarily, we get $\pi_k M = \pi_k \Tot(M^\bullet)$ for all $k$ and consequently that $M \to \Tot(M^\bullet)$ is an isomorphism, as desired.
\end{proof}

\begin{corollary} \label{rslt:v-hyperdescent-over-td-can-be-checked-on-modules}
Let $Y_\bullet$ be a simplicial object in $\AffPerfd_\pi$ and write $Y_n = \Spa(B^n, B^{n+})$. If $Y_0$ is totally disconnected then the natural functor
\begin{align*}
	\varprojlim_{n\in\Delta} \Dqcohri(B^{n+}/\pi) \isoto \varprojlim_{n\in\Delta} \DqcohriX{Y_n}
\end{align*}
is an equivalence.
\end{corollary}
\begin{proof}
We have $\DqcohriX{Y_0} = \Dqcohri(B^{0+}/\pi)$ by \cref{rslt:computation-of-Dqcohri-in-td-case}. For every $n \ge 0$ let $\mathcal C_1^n \subset \Dqcohri(B^{n+}/\pi)$ and $\mathcal C_2^n \subset \DqcohriX{Y_n}$ be the essential images of the base-change functors $\Dqcohri(B^{0+}/\pi) \to \Dqcohri(B^{n+}/\pi)$ and $\DqcohriX{Y_0} \to \DqcohriX{Y_n}$, respectively. Then
\begin{align*}
	\varprojlim_{n\in\Delta} \Dqcohri(B^{n+}/\pi) = \varprojlim_{n\in\Delta} \mathcal C_1^n \isoto \varprojlim_{n\in\Delta} \mathcal C_2^n = \varprojlim_{n\in\Delta} \DqcohriX{Y_n},
\end{align*}
where the equivalence in the middle follows from \cref{rslt:base-change-over-tot-disc-equiv}.
\end{proof}

\begin{proposition} \label{rslt:Dqcohri-is-v-sheaf}
The presheaf $X \mapsto \DqcohriX X$ on $\AffPerfd_\pi$ is a v-sheaf.
\end{proposition}
\begin{proof}
Let $Y_\bullet \to X$ be a v-Čech cover in $\AffPerfd_\pi$. We have to show that the natural functor $\DqcohriX X \to \varprojlim_{n\in\Delta} \DqcohriX{Y_n}$ is an equivalence. By (quasi-)pro-étale descent (see \cref{rslt:check-sheafiness-on-finer-topology}) we can assume that $X$ and $Y_0$ are totally disconnected. Writing $X = \Spa(A, A^+)$ and $Y_n = \Spa(B^n, B^{n+})$ we now have $\DqcohriX X = \Dqcohri(A^+/\pi)$ by \cref{rslt:computation-of-Dqcohri-in-td-case} and $\varprojlim_{n\in\Delta} \DqcohriX{Y_n} = \varprojlim_{n\in\Delta} \Dqcohri(B^{n+}/\pi)$ by \cref{rslt:v-hyperdescent-over-td-can-be-checked-on-modules}. Thus the claim reduces to checking that mod-$\pi$ modules descend along $Y_0 \to X$, which is \cref{rslt:v-cover-of-tot-disc-is-mod-pi-descendable}.
\end{proof}

We now want to show that $X \mapsto \DqcohriX X$ is even a \emph{hypercomplete} v-sheaf, i.e. it satisfies descent along arbitrary v-hypercovers. Again the Tor dimension result \cref{rslt:mod-pi-Tor-dim-1-over-tot-disc-space} comes in very handy:

\begin{proposition} \label{rslt:v-hypercover-implies-fully-faithful-on-Dqcohri}
Let $Y_\bullet \to X$ be a v-hypercover in $\AffPerfd_\pi$. Then the natural functor $\DqcohriX X \injto \varprojlim_{n\in\Delta} \DqcohriX{Y_n}$ is fully faithful.
\end{proposition}
\begin{proof}
The described functor has a right adjoint $\varprojlim_{n\in\Delta} \DqcohriX{Y_n} \to \DqcohriX X$ and it is enough to show that the unit of the adjunction is isomorphic to the identity. By pro-étale descent we can assume that $X$ is totally disconnected (see \cref{rslt:sheaves-on-basis-equiv-sheaves-on-whole-site}). In this case, writing $X = \Spa(A, A^+)$ and $Y_n = \Spa(B^n, B^{n+})$, we have $\DqcohriX X = \Dqcohri(A^+/\pi)$ and by \cref{rslt:v-hyperdescent-over-td-can-be-checked-on-modules} the functor $\varprojlim_{n\in\Delta} \Dqcohri(B^{n+}/\pi) \to \varprojlim_{n\in\Delta} \DqcohriX{Y_n}$ restricts to an equivalence of the essential images of the natural functors $\Dqcohri(A^+/\pi) \to \varprojlim_{n\in\Delta} \Dqcohri(B^{n+}/\pi)$ and $\DqcohriX X \to \varprojlim_{n\in\Delta} \DqcohriX{Y_n}$. Hence it is enough to show that the functor $\Dqcohri(A^+/\pi) \to \varprojlim_{n\in\Delta} \Dqcohri(B^{n+}/\pi)$ is fully faithful.

By the same adjunction argument as before, we are reduced to showing the following: Let $M \in \Dqcohri(A^+/\pi)$ be given and let $M^\bullet = M \tensor_{(A^+/\pi)^a_\solid} (B^{n+}/\pi)^a_\solid$. Then the natural map $M \to \Tot(M^\bullet)$ is an isomorphism. As in the proof of \cref{rslt:base-change-over-tot-disc-equiv} we can assume that $M \in \Dqcohrile0(A^+/\pi)$, in which case $M^n \in \Dqcohrile1(B^{n+}/\pi)$ by \cref{rslt:mod-pi-Tor-dim-1-over-tot-disc-space} and hence $\pi_k \Tot(M^\bullet) = \pi_k \Tot_n(M^\bullet)$ for $k > -n+1$.

Now, for fixed $n_0 \ge 0$ let $\tilde Y_\bullet = \Spa(\tilde B^{\bullet}, \tilde B^{\bullet+}) := \cosk_{n_0} Y_\bullet$ be the $n_0$-th coskeleton of $Y_\bullet \to X$ and let $\tilde M^\bullet = M \tensor_{(A^+/\pi)^a_\solid} (\tilde B^{\bullet+}/\pi)^a_\solid$ be the associated cosimplicial object. Then $\tilde Y_m = Y_m$ for $m \le n_0$ so that $\tilde M^m = M^m$. Moreover, by \cref{rslt:Dqcohri-is-v-sheaf} and \cite[Lemma 6.5.3.9]{lurie-higher-topos-theory} we know that $\DqcohriX X = \varprojlim_{n\in\Delta} \DqcohriX{\tilde Y_n}$, so in particular $M = \Tot(\tilde M^\bullet)$. Thus for all $k > -n_0+1$ we deduce
\begin{align*}
	\pi_k \Tot(M^\bullet) = \pi_k \Tot_{n_0}(M^\bullet) = \pi_k \Tot_{n_0}(\tilde M^\bullet) = \pi_k \Tot(\tilde M^\bullet) = \pi_k M.
\end{align*}
Therefore, for any $k \in \Z$ we get $\pi_k \Tot(M^\bullet) = \pi_k(M)$ by choosing $n_0$ large enough. In particular the map $M \to \Tot(M^\bullet)$ is an isomorphism, as desired.
\end{proof}

Roughly, \cref{rslt:v-hypercover-implies-fully-faithful-on-Dqcohri} states that for every $X \in \AffPerfd_\pi$ and any $M \in \DqcohriX X$, $M$ defines a ``hypercomplete sheaf'' on $X_\vsite$. It is a formal consequence of this fact (and the sheafiness of $\DqcohriX X$) that $\DqcohriX X$ is itself hypercomplete:

\begin{theorem} \label{rslt:v-hyperdescent-for-Dqcohri}
The presheaf $X \mapsto \DqcohriX X$ on $\AffPerfd_\pi$ is a hypercomplete v-sheaf.
\end{theorem}
\begin{proof}
This follows from \cref{rslt:Dqcohri-is-v-sheaf,rslt:v-hypercover-implies-fully-faithful-on-Dqcohri} using the formal result \cref{rslt:fully-faithful-on-hypercovers-implies-hypercomplete}.
\end{proof}

\subsection{Basic Definitions on v-Stacks} \label{sec:ri-pi.def-on-vstack}

Having developed a good theory of quasicoherent $\ri^{+a}_X/\pi$-modules on affinoid perfectoid spaces, we will now extend these definitions to small v-stacks. Most of this is rather formal, but a little care is required to properly define the sheaf $\ri^{+a}_X/\pi$ on a general small v-stack. First recall the following definition of structure sheaves (cf. \cite[Definition 2.1]{mann-werner-simpson}):

\begin{definition}
Let $X$ be a small v-stack.
\begin{defenum}
	\item An \emph{untilt} $X^\sharp$ of $X$ is a map $X \to \Spa\Z_p$. A map of untilted small v-stacks $Y^\sharp \to X^\sharp$ is a map of small v-stacks $Y \to X$ over $\Spa\Z_p$. Note that $X$ comes equipped with a canonical untilt $X^\flat$ via $X \to * \to \Spa\Z_p$.

	For any untilt $X^\sharp$ of $X$ and any map $Y \to X$ of small v-stacks, $Y$ can uniquely be equipped with an untilt $Y^\sharp$ such that $Y \to X$ becomes a map $Y^\sharp \to X^\sharp$ of untilted small v-stacks. Moreover, if $X$ is a perfectoid space then untilts $X^\sharp$ of $X$ in the above sense correspond to untilts of $X$ as a perfectoid space.

	\item Fix an untilt $X^\sharp$ of $X$. Then for every map $Y \to X$ from a perfectoid space $Y$, let $Y^\sharp$ be the induced untilt. We define the sheaves
	\begin{align*}
		\ri_{X^\sharp}, \ri_{X^\sharp}^+,
	\end{align*}
	on $X_\vsite$ to be the unique sheaves which restrict to $\ri_{Y^\sharp}$ resp. $\ri^+_{Y^\sharp}$ for all perfectoid $Y \to X$.
\end{defenum}
\end{definition}

Some elementary properties of the structure sheaves $\ri_{X^\sharp}$ and $\ri_{X^\sharp}^+$ on an untilted small v-stack $X^\sharp$ are studied in \cite[\S2]{mann-werner-simpson} (there, these sheaves are denoted $\check\ri_X$ and $\check\ri_X^+$). Note that if $X = Z^\diamond$ for an analytic adic space $Z$, then the unique map $Z \to \Spa \Z_p$ induces a map $X \to \Spa\Z_p$, i.e. an untilt $X^\sharp$ and hence also a structure sheaf $\ri_{X^\sharp}$. This structure sheaf on $X$ restricts to the structure sheaf on $Z$ if $Z$ is a seminormal (e.g. smooth) rigid space over some non-archimedean field in characteristic $0$ (cf. \cite[Proposition 2.8]{mann-werner-simpson}).

In order to define the sheaf $\ri^{+a}_{X^\sharp}/\pi$ on an untilted small v-stack $X^\sharp$ we make the following definition.

\begin{definition}
Let $X^\sharp$ be an untilted small v-stack. A \emph{pseudouniformizer} on $X^\sharp$ is a classical sheaf of ideals $\pi \subset \ri^+_{X^\sharp}$ such that for all affinoid perfectoid spaces $Y = \Spa(B, B^+) \to X$, $\pi(Y) \subset B^{\sharp+}$ is a principal ideal generated by a topologically nilpotent unit in $B$. We say that $\pi$ is \emph{globally generated} if there is a global section $a \in H^0(X, \ri^+_{X^\sharp})$ such that $\pi = a \ri^+_{X^\sharp}$. In the case that $\pi$ is globally generated (e.g. $X$ is an affinoid perfectoid space), we will often confuse the ideal $\pi$ with any element that generates it.
\end{definition}

\begin{remark}
Pseudouniformizers do not necessarily exist on a small v-stack. For example, $X = *$ does not admit a pseudouniformizer. Note that in most practical situations the untilt and pseudouniformizer on $X$ come from an untilt and pseudouniformizer on some base field $K$ over which $X$ is defined. \end{remark}

While we are usually interested in the category of ``quasicoherent $\ri^+_{X^\sharp}/\pi$-modules'' on an untilted small v-stack $X^\sharp$, it can be useful to allow more general coefficient sheaves. Defining these coefficient sheaves is a bit subtle; we start by introducing almost setups on $X$.

\begin{remark}
The following definitions of almost setups and integral torsion coefficients are rather subtle and for the most part are not actually needed (but they help to make our notation much cleaner). The reader is invited to skip these definitions and always assume $\Lambda = \ri^+_{X^\sharp}/\pi$ for an untilt $X^\sharp$ and a pseudouniformizer $\pi$. In this case the reader should jump straight to \cref{rslt:def-of-qcoh-Lambda-modules}. Note that we will sometimes make the assumption ``$X \in X_\vsite^\Lambda$''. In the case that $\Lambda = \ri^+_{X^\sharp}/\pi$ this amounts to saying that there are ``enough global pseudouniformizers'' on $X^\sharp$ in an appropriate sense. This condition is always satisfied if $X$ admits a map to some affinoid perfectoid space.
\end{remark}

\begin{definition} \label{def:almost-setups-on-v-stacks}
Let $X$ be a small v-stack.
\begin{defenum}
	\item An \emph{almost setup} on $X$ is a triple $(\Lambda, X_\vsite^\Lambda, \mm_\Lambda)$, where $\Lambda$ is a classical sheaf of rings on $X_\vsite$, $X_\vsite^\Lambda \subset X_\vsite$ is a covering sieve and $\mm_\Lambda \subset \restrict\Lambda{X_\vsite^\Lambda}$ is a presheaf of ideals on $X_\vsite^\Lambda$ with the following properties:
	\begin{enumerate}[(i)]
		\item For every $Y \in X_\vsite^\Lambda$ the pair $(\Lambda(Y), \mm_\Lambda(Y))$ is an almost setup (see \cref{def:almost-setups}).

		\item For every map $Y' \to Y$ in $X_\vsite^\Lambda$ the induced morphism $\Lambda(Y) \to \Lambda(Y')$ provides a strict morphism of almost setups with respect to $\mm_\Lambda(Y)$ and $\mm_\Lambda(Y')$.
	\end{enumerate}
	We will usually omit the sieve $X_\vsite^\Lambda$ and the ideal sheaf $\mm_\Lambda$ from the notation and simply speak of an almost setup $\Lambda$.

	\item A \emph{morphism of almost setups} $\varphi\colon (\Lambda, X_\vsite^\Lambda, \mm_\Lambda) \to (\Lambda', X_\vsite^{\Lambda'}, \mm_{\Lambda'})$ on $X$ is a morphism of classical sheaves of rings $\varphi\colon \Lambda \to \Lambda'$ such that for every $Y \in X_\vsite^\Lambda \isect X_\vsite^{\Lambda'}$ the morphism $\varphi(Y)\colon \Lambda(Y) \to \Lambda'(Y)$ provides a morphism of almost setups with respect to $\mm_\Lambda(Y)$ and $\mm_{\Lambda'}(Y)$.

	The morphism $\varphi$ is called \emph{strict} if for all $Y \in  X_\vsite^\Lambda \isect X_\vsite^{\Lambda'}$ the morphism $\varphi(Y)$ is a strict morphism of almost setups.

	\item Let $\kappa$ be a regular cardinal and $\Lambda$ an almost setup on $X$. We say that $\Lambda$ is \emph{$\kappa$-compact} if for every $Y \in X_\vsite^\Lambda$ the almost setup $(\Lambda(Y), \mm_\Lambda(Y))$ is $\kappa$-compact.
\end{defenum}
\end{definition}

Of course, given an untilted small v-stack $X^\sharp$ with pseudouniformizer $\pi$, there should be an almost setup on $X$ given by the sheaf $\ri^+_{X^\sharp}/\pi$. If $X$ is not quasi-compact then the definition of this almost setup is a bit sublte, however. It can be done as follows:

\begin{definition}
Let $X^\sharp$ be an untilted small v-stack with pseudouniformizer $\pi$.
\begin{defenum}
	\item We say that a pseudouniformizer $\varpi$ on $(X^\sharp, \pi)$ is \emph{bounded} if $\pi \subset \varpi$ and there is an integer $n > 0$ with $\varpi^n \divides \pi$.

	\item We say that $(X^\sharp, \pi)$ has \emph{enough pseudouniformizers} if $\pi$ is globally generated and for every globally generated bounded pseudouniformizer $\varpi$ on $X^\sharp$ there is a globally generated bounded pseudouniformizer $\varpi'$ on $X^\sharp$ with $\varpi'^2 \divides \varpi$. We denote $\mm_{X^\sharp,\pi} \subset \ri^+_X(X)$ the ideal generated by all globally generated bounded pseudouniformizers.
\end{defenum}
\end{definition}

\begin{remark}
If $X^\sharp = \Spa(A^\sharp, A^{\sharp+}) \in \AffPerfd_\pi$ is an affinoid perfectoid space then $\mm_{X^\sharp,\pi} = A^{\sharp\circ\circ}$ is independent of the pseudouniformizer $\pi$. Clearly $X^\sharp$ has enough pseudouniformizers (and by \cref{rslt:enough-pseudouniformizers-along-map-of-v-stacks} the same is true for any $(Y^\sharp,\pi) \in \vStackspi$ which admits a map to $(X^\sharp, \pi)$).

If $X^\sharp$ is a quasi-compact small v-stack then still $\mm_{X^\sharp,\pi}$ is independent of $\pi$, as one checks by considering a cover $Y^\sharp \surjto X^\sharp$ by some affinoid perfectoid space $Y^\sharp$. In general however, $\mm_{X^\sharp,\pi}$ does depend on $\pi$. \end{remark}

\begin{lemma} \label{rslt:properties-of-enough-pseudouniformizers}
Let $X^\sharp$ be an untilted small v-stack with pseudouniformizer $\pi$. Assume that $(X^\sharp, \pi)$ has enough pseudouniformizers.
\begin{lemenum}
	\item The pair $(H^0(X, \ri^+_{X^\sharp}), \mm_{X^\sharp,\pi})$ is an $\omega_1$-compact almost setup.

	\item \label{rslt:enough-pseudouniformizers-along-map-of-v-stacks} Let $f\colon Y \to X$ be a map of small v-stacks. Then $(Y^\sharp,\pi)$ has enough pseudouniformizers and the induced map $(H^0(X, \ri^+_{X^\sharp}), \mm_{X^\sharp,\pi}) \to (H^0(Y, \ri^+_{Y^\sharp}), \mm_{Y^\sharp,\pi})$ is a strict morphism of almost setups.
\end{lemenum}
\end{lemma}
\begin{proof}
We first prove (i). To shorten notation let us denote $A = H^0(X, \ri_{X^\sharp})$ and $A^+ = H^0(X, \ri^+_{X^\sharp})$. It follows immediately from the definition that $\mm_{X^\sharp,\pi}^2 = \mm_{X^\sharp,\pi}$. We claim that $\mm_{X^\sharp,\pi}$ is a countable filtered colimit of copies of $A^+$; then in particular $\mm_{X^\sharp,\pi}$ is flat and $\omega_1$-compact, which implies (i). To prove the claim, pick any globally generated bounded pseudouniformizer $\varpi$ on $X^\sharp$. Then $\varpi$ is a unit in $A$ (because this is true locally), so in particular it is not a zero-divisor, hence the principal ideal $(\varpi) \subset \mm_{X^\sharp,\pi}$ is a free $A^+$-module. If $\varpi'$ is another globally generated bounded pseudouniformizer on $X^\sharp$ then there is some integer $n > 0$ such that $\varpi' \divides \varpi^n$. By definition of ``enough pseudouniformizers'' we can find some globally generated bounded pseudouniformizer $\varpi''$ such that $\varpi''^n \divides \varpi'$. Then $\varpi''^n \divides \varpi^n$, hence $(\varpi \varpi''^{-1})^n \in A^+$. As $A^+$ is integrally closed in $A$ (this can be checked locally and holds by definition on affinoid perfectoid spaces) we deduce $\varpi \varpi''^{-1} \in A^+$, hence $\varpi'' \divides \varpi$. All in all we deduce that the ideal $(\varpi, \varpi') \subset \mm_{X^\sharp,\pi}$ is contained in the principal ideal $(\varpi'') \subset \mm_{X^\sharp,\pi}$. By induction the same is true for any finitely generated subideal of $\mm_{X^\sharp,\pi}$, hence $\mm_{X^\sharp,\pi}$ is the union of its principal subideals $(\pi) \subset \mm_{X^\sharp,\pi}$ for bounded pseudouniformizers $\pi$. By a similar argument as above, a countable such union is enough: For each $n > 0$ pick a globally generated bounded pseudouniformizer $\pi_n$ such that $\pi_n^n \divides \pi$. Then every globally generated bounded pseudouniformizer is divisible by some $\pi_n$.

We now prove (ii) so let $f\colon (Y^\sharp,\pi) \to (X^\sharp,\pi)$ be given. In addition to the above notation $A$ and $A^+$ we denote $B = H^0(Y, \ri_{Y^\sharp})$ and $B^+ = H^0(Y, \ri^+_{Y^\sharp})$. Let $\varpi$ be a globally generated bounded pseudouniformizer on $Y^\sharp$. Choose any $n > 0$ such that $\pi \divides \varpi^n$ and a globally generated bounded pseudouniformizer $\varpi'$ on $X^\sharp$ with $\varpi'^n \divides \pi$. Then $(f^* \varpi')^n \divides \varpi^n$, hence by the same argument as in the proof of (i) we deduce $f^*\varpi' \divides \varpi$. It follows that $(Y^\sharp,\pi)$ has enough pseudouniformizers (pick $\varpi''$ on $X^\sharp$ such that $\varpi''^2 \divides \varpi'$; then $(f^*\varpi'')^2 \divides \varpi$) and that $\mm_{Y^\sharp,\pi}$ is the ideal generated by $f^* \mm_{X^\sharp,\pi}$, as desired.
\end{proof}

\begin{definition}
Let $X^\sharp$ be an untilted small v-stack with pseudouniformizer $\pi$.
\begin{defenum}
	\item Let $X_\vsite^{\sharp,\pi} \subset X_\vsite$ be the full subcategory spanned by the small v-stacks with enough pseudouniformizers.

	\item For every $Y \in X_\vsite^{\sharp,\pi}$, let $\mm_{\ri^+_{X^\sharp}/\pi}(Y) \subset (\ri^+_{X^\sharp}/\pi)(Y)$ be the ideal generated by the image of $\mm_{Y^\sharp,\pi}$.
\end{defenum}
By \cref{rslt:properties-of-enough-pseudouniformizers} the triple
\begin{align*}
	(\ri^+_{X^\sharp}/\pi, X_\vsite^{\sharp,\pi}, \mm_{\ri^+_{X^\sharp}/\pi})
\end{align*}
is an $\omega_1$-compact almost setup on $X$. In the following we will usually abbreviate it by $\ri^+_{X^\sharp}/\pi$.
\end{definition}

We have defined a canonical almost setup on every untilted small v-stack $X^\sharp$ with pseudouniformizer $\pi$. It is now straightforward to introduce the following slightly more general setting:

\begin{definition} \label{def:intergal-torsion-coefficients}
Let $X$ be a small v-stack.
\begin{defenum}
	\item \label{def:intergal-torsion-coefficients.1} A \emph{system of integral torsion coefficients on $X$} is an almost setup $\Lambda$ on $X$ with the following properties:
	\begin{enumerate}[(i)]
		\item For all affinoid perfectoid spaces $Y \in X_\vsite$ we have $Y \in X_\vsite^\Lambda$ and $H^i(Y, \Lambda)^a = 0$ for $i > 0$.

		\item There is a strict morphism $\Lambda_0 := \ri^+_{X^\sharp}/\pi \to \Lambda$ of almost setups for some untilt $X^\sharp$ of $X$ and some pseudouniformizer $\pi$ on $X^\sharp$ such that for all maps $Y' \to Y$ of affinoid perfectoid spaces in $X_\vsite$ the induced morphism of analytic rings
		\begin{align*}
			\Lambda(Y)^a_\solid \tensor_{\Lambda_0(Y)^a_\solid} \Lambda_0(Y')^a_\solid \isoto \Lambda(Y')^a_\solid
		\end{align*}
		is an isomorphism.
	\end{enumerate}
	A morphism $\Lambda \to \Lambda'$ of integral torsion coefficients on $X$ is a strict morphism of almost setups such that for all maps $Y' \to Y$ of affinoid perfectoid spaces in $X_\vsite$ the induced morphism of analytic rings
	\begin{align*}
		\Lambda'(Y)^a_\solid \tensor_{\Lambda(Y)^a_\solid} \Lambda(Y')^a_\solid \isoto \Lambda'(Y')^a_\solid
	\end{align*}
	is an isomorphism.

	\item \label{def:vStacksCoeff} We denote by $\vStacksCoeff$ the $2$-category of pairs $(X, \Lambda)$, where $X$ is a small v-stack and $\Lambda$ is a system of integral torsion coefficients on $X$. A morphism $(Y, \Lambda_Y) \to (X, \Lambda_X)$ is a morphism $Y \to X$ of small v-stacks such that $\Lambda_Y = \restrict{\Lambda_X}Y$. We will often omit $\Lambda$ from the notation and simply write $X \in \vStacksCoeff$ with an implicit choice of $\Lambda$.

	We denote by $\AffPerfCoeff \subset \vStacksCoeff$ the full subcategory spanned by the pairs $(X, \Lambda)$, where $X$ is an affinoid perfectoid space of characteristic $p$.
\end{defenum}
\end{definition}

\begin{remark}
In more concrete terms, the isomorphism of analytic rings in \cref{def:intergal-torsion-coefficients.1} boils down to the claim that the map of classical rings
\begin{align*}
	\Lambda'(Y) \tensor_{\Lambda(Y)} \Lambda(Y') \to \Lambda'(Y')
\end{align*}
is integral and an almost isomorphism (see \cref{rslt:colimits-of-solid-discrete-almost-rings}).
\end{remark}

Having established the general setup which we want to work in, it is time to construct the main object of interest: the $\infty$-category $\Dqcohri(X, \Lambda)$ of ``quasicoherent solid almost $\Lambda$-modules'' on $X$, for every $(X, \Lambda) \in \vStacksCoeff$. In order to construct this $\infty$-category we need to generalize the descent results from \cref{sec:ri-pi.descent-on-affperfd} to arbitrary intergal torsion coefficients $\Lambda$. This is mostly straightforward, but we need to handle the following subtlety:

\begin{lemma} \label{rslt:A+-mod-pi-equiv-ri-+-mod-pi}
Let $X = \Spa(A, A^+) \in \AffPerfd_\pi$. Then the natural morphism of rings $A^+/\pi \to (\ri^+_X/\pi)(X)$ induces an isomorphism of analytic rings
\begin{align*}
	(A^+/\pi)^a_\solid \isoto ((\ri^+_X/\pi)(X))^a_\solid.
\end{align*}
\end{lemma}
\begin{proof}
Note first that by \cite[Proposition 2.13]{mann-werner-simpson} it does not matter whether we compute the quotient $\ri^+_X/\pi$ on $X_\vsite$ or on $X_\proet$ (or on $X_\et$); moreover if $X$ is strictly totally disconnected then $(\ri^+_X/\pi)(X) = A^+/\pi$, so the claim is evidently true in this case. In general we know by \cite[Proposition 8.8]{etale-cohomology-of-diamonds} that $A^+/\pi \to (\ri^+_X/\pi)(X)$ is an almost isomorphism, from which it follows easily that $(A^+/\pi)^a_\solid \to ((\ri^+_X/\pi)(X))^a_\solid$ is a localization, i.e. the forgetful functor along this map is fully faithful. To prove the claimed isomorphism, it is therefore enough to show that the pullback functor $- \tensor_{(A^+/\pi)^a_\solid} ((\ri^+_X/\pi)(X))^a_\solid$ is conservative.

Suppose first that $\pi \divides p$. Then w.l.o.g. $X$ has characteristic $p$. In particular there is a canonical map $X \to Z = \Spa(K, \ri_K)$, where $K$ is the completed perfection of $\Fld((\pi))$. By \cref{rslt:any-map-in-AffPerf-is-pro-wh} we can write $X \to Z$ as a filtered colimit $X = \varprojlim_i X_i \to Z$ such that all $X_i = \Spa(A_i, A^+_i)$ are weakly of perfectly finite over $Z$. Then $A^+/\pi = \varinjlim_i A^+_i/\pi$ and we claim that the natural map
\begin{align*}
	\varinjlim_i (\ri^+_{X_i}/\pi)(X_i) \isoto (\ri^+_X/\pi)(X)
\end{align*}
is also an isomorphism. To see this, note that $\ri^+_X/\pi$ is the étale sheafification of the presheaf $\Spa(A', A'^+) \mapsto A'^+/\pi$. Thus the claimed isomorphism follows from the fact that every étale cover of $X$ comes via base-change from an étale cover of some $X_i$ (see \cite[Proposition 6.4.(iv)]{etale-cohomology-of-diamonds}), using that the presheaf $\Spa(A', A'^+) \mapsto A'^+/\pi$ satisfies the desired colimit property.

Using the above isomorphism and the fact that $(-)^a_\solid$ commutes with (non-empty) colimits by \cref{rslt:colimits-of-solid-discrete-rings}, we can reduce to the case that $X = X_i$ for some $i$, i.e. $X$ is weakly of perfectly finite type over the totally disconnected space $Z$. Pick any pro-étale cover $Y \surjto X$ by some strictly totally disconnected space $Y = \Spa(B, B^+)$. By \cref{rslt:compute-Dqcohri-for-fin-type-over-tot-disc} we have $\DqcohriX X = \Dqcohri(A^+/\pi)$ and $\DqcohriX Y = \Dqcohri(B^+/\pi)$, which implies that the pullback functor $- \tensor_{(A^+/\pi)^a_\solid} (B^+/\pi)^a_\solid$ is conservative. But $B^+/\pi = (\ri^+_Y/\pi)(Y)$ (by the first paragraph of the proof), hence $A^+/\pi \to B^+/\pi$ factors over $A^+/\pi \to (\ri^+_X/\pi)(X)$. This implies that the pullback functor $- \tensor_{(A^+/\pi)^a_\solid} ((\ri^+_X/\pi)(X))^a_\solid$ is conservative, as desired. We have finished the proof in the case $\pi \divides p$.

Now let $\pi$ be general. Choose a pseudouniformizer $\varpi \in A^+$ and an integer $n \ge 1$ such that $\varpi \divides p$, $\varpi^n \divides \pi$ and $\pi \divides \varpi^{n+1}$. We claim that the pullback $- \tensor_{(A^+/\pi)^a_\solid} (A^+/\varpi)^a_\solid$ is conservative. Namely, suppose we have some $M \in \Dqcohri(A^+/\pi)$ such that $M \tensor_{(A^+/\pi)^a_\solid} (A^+/\varpi)^a_\solid = 0$. Tensoring $M$ with the exact sequence $0 \to A^+/\frac\pi\varpi \to A^+/\pi \to A^+/\varpi \to 0$ shows that $M \tensor A^+/\frac\pi\varpi \isoto M$. Repeating this procedure with $\frac\pi{\varpi^m}$ in place of $\pi$ inductively shows that $M \isom M \tensor A^+/\frac\pi{\varpi^m}$ for all $m \le n$, in particular for $m = n$. On the other hand we have $\frac\pi{\varpi^n} \divides \varpi$, hence
\begin{align*}
	M \isom M \tensor A^+/\frac\pi{\varpi^n} = (M \tensor A^+/\varpi) \tensor A^+/\frac\pi{\varpi^n} = 0,
\end{align*}
as desired. Now consider the following commutative diagram of analytic rings:
\begin{center}\begin{tikzcd}
	(A^+/\pi)^a_\solid \arrow[r] \arrow[d] & ((\ri^+_X/\pi)(X))^a_\solid \arrow[d]\\
	(A^+/\varpi)^a_\solid \arrow[r] & ((\ri^+_X/\varpi)(X))^a_\solid
\end{tikzcd}\end{center}
By what we have shown above we know that the pullbacks along the lower horizontal map and the left vertical map are conservative. It follows that the pullback along their composition is conservative and in particular the pullback along the upper horizontal map is conservative, as desired.
\end{proof}

We can now restate the main results of \cref{sec:ri-pi.descent-on-affperfd} in terms of an arbitrary system of integral torsion coefficients $\Lambda$.

\begin{proposition} \label{rslt:basic-properties-of-integral-torsion-coeffs}
\begin{propenum}
	\item \label{rslt:base-change-for-int-tor-coeffs-on-aff-perf} Let $Y \to X \from X'$ be a diagram in $\AffPerfCoeff$. Then the natural map of analytic rings
	\begin{align*}
		\Lambda(Y)^a_\solid \tensor_{\Lambda(X)^a_\solid} \Lambda(X')^a_\solid \isoto \Lambda(Y \cprod_X X')^a_\solid
	\end{align*}
	is an isomorphism.

	\item \label{rslt:Tor-dim-and-desc-for-int-tor-coeffs-over-tot-disc} Let $Y \to X$ be a map in $\AffPerfCoeff$ and assume that $X$ is totally disconnected. Then the map $\Lambda(X)^a_\solid \to \Lambda(Y)^a_\solid$ has Tor dimension $\le 1$. If $f$ is a v-cover then this map is weakly fs-descendable of index $\le 4$.

	\item \label{rslt:Tor-dim-and-desc-for-int-tor-coeffs-on-Z'-over-Z} Let $f\colon Z' \to Z$ be a map of totally disconnected spaces in $\AffPerfCoeff$ with $\dimtrg f < \infty$. Then there is a pro-étale cover $Y \to \overline Z'^{/Z}$ by a totally disconnected space $Y$ such that $\Lambda(\overline Z'^{/Z})^a_\solid \to \Lambda(Y)^a_\solid$ has Tor dimension $\le 1$ and is weakly fs-descendable of index bounded by a constant which only depends on $\dimtrg f$.

	\item \label{rslt:fin-type-over-tot-disc-implies-Lambda-desc} Let $Z \in \AffPerfCoeff$ be totally disconnected and let $X \to Z$ be a map of affinoid perfectoid spaces which is of weakly perfectly finite type. Then for every v-cover $Y \surjto X$ the map $\Lambda(\overline Y^{/X})^a_\solid \to \Lambda(X)^a_\solid$ is fs-descendable of index bounded by a constant which only depends on the number of perfect generators for $X$ over $Z$.
\end{propenum}
\end{proposition}
\begin{proof}
In both (i) and (ii), fix a morphism $\ri^+_{X^\sharp}/\pi \to \Lambda$ of integral torsion coefficients for some untilt $X^\sharp$ of $X$ and some pseudouniformizer $\pi$ on $X^\sharp$. Then (i) reduces immediately to the case $\Lambda = \ri^+_{X^\sharp}/\pi$ (by definition of morphisms of integral torsion coefficients), which by \cref{rslt:A+-mod-pi-equiv-ri-+-mod-pi} reduces to \cref{rslt:base-change-for-affinoid-perfectoid}. Similarly, using \cref{rslt:A+-mod-pi-equiv-ri-+-mod-pi}, the first part of (ii) reduces to \cref{rslt:mod-pi-Tor-dim-1-over-tot-disc-space}, while the second part of (ii) reduces to \cref{rslt:v-cover-of-tot-disc-is-mod-pi-descendable} (using that weakly fs-descendable maps are stable under base-change, see \cref{rslt:stability-of-pi-0-weakly-descendable-maps}).

Part (iii) follows in a similar vein as (i) and (ii) using \cref{rslt:compute-Dqcohri-for-subset-of-Z'-over-Z,rslt:Y-over-Z'-Z-has-mod-pi-Tor-dim-1}. Similarly, part (iv) follows from \cref{rslt:fin-type-over-K-implies-mod-pi-desc}.
\end{proof}

\begin{proposition} \label{rslt:def-of-qcoh-Lambda-modules}
There is a unique hypercomplete v-sheaf
\begin{align*}
	(\vStacksCoeff)^\opp \to \infcatinf, \qquad (X, \Lambda) \mapsto \Dqcohri(X, \Lambda)
\end{align*}
of $\infty$-categories on $\vStacksCoeff$ which on all those $(X, \Lambda) \in \AffPerfCoeff$ with $X$ being of weakly perfectly finite type over some totally disconnected space takes the form
\begin{align*}
	\Dqcohri(X, \Lambda) = \Dqcohri(\Lambda(X)).
\end{align*}
\end{proposition}
\begin{proof}
Uniqueness of the sheaf is clear (see \cref{rslt:sheaves-on-basis-equiv-sheaves-on-whole-site}) so we only need to show existence, for which we can restrict to $\AffPerfCoeff$. Now let us temporarily denote by $(X, \Lambda) \mapsto \Dqcohri(X, \Lambda)$ the quasi-pro-étale sheafification of the presheaf $(X, \Lambda) \mapsto \Dqcohri(\Lambda(X))$ on $\AffPerfCoeff$ (to avoid set-theoretic issues, we can temporarily work with sheaves of $\infty$-categories in a bigger universe). Using \cref{rslt:basic-properties-of-integral-torsion-coeffs} we deduce by the same argument as in \cref{rslt:computation-of-Dqcohri-in-td-case} that $\Dqcohri(X, \Lambda) = \Dqcohri(\Lambda(X))$ whenever $X$ admits a quasi-pro-étale map to some totally disconnected space. It then follows by the same argument as in \cref{rslt:compute-Dqcohri-for-fin-type-over-tot-disc} (using again \cref{rslt:basic-properties-of-integral-torsion-coeffs}) that $\Dqcohri(X, \Lambda) = \Dqcohri(\Lambda(X))$ whenever $X$ is weakly of perfectly finite type over some totally disconnected space.

To finish the proof it remains to show that the quasi-pro-étale sheaf $(X, \Lambda) \mapsto \Dqcohri(X, \Lambda)$ is a hypercomplete v-sheaf. To see this, first note that \cref{rslt:general-etale-descent-of-Dqcohri} generalizes to $\Lambda$-coefficients by the same arguments as in the proof of \cref{rslt:basic-properties-of-integral-torsion-coeffs}. Using this and \cref{rslt:basic-properties-of-integral-torsion-coeffs} we can generalize \cref{rslt:base-change-over-tot-disc-equiv,rslt:v-hyperdescent-over-td-can-be-checked-on-modules,rslt:Dqcohri-is-v-sheaf} to $\Lambda$-coefficients; in particular $(X, \Lambda) \mapsto \Dqcohri(X, \Lambda)$ is a v-sheaf. Finally we note that \cref{rslt:v-hypercover-implies-fully-faithful-on-Dqcohri} also generalizes to $\Lambda$-coefficients: One can either deduce the generalized version from the existing version by the arguments in the proof of \cref{rslt:basic-properties-of-integral-torsion-coeffs} or one can repeat the argument. Thus \cref{rslt:fully-faithful-on-hypercovers-implies-hypercomplete} implies that $(X, \Lambda) \mapsto \Dqcohri(X, \Lambda)$ is hypercomplete.
\end{proof}

\begin{definition}
Let $X$ be a small v-stack with integral torsion coefficients $\Lambda$. The $\infty$-category $\Dqcohri(X, \Lambda)$ defined by \cref{rslt:def-of-qcoh-Lambda-modules} is called the $\infty$-category of \emph{quasicoherent solid almost $\Lambda$-modules on $X$}. In the case that $\Lambda = \ri^+_{X^\sharp}/\pi$ for some untilt $X^\sharp$ of $X$ and some pseudouniformizer $\pi$ on $X^\sharp$ we abbreviate $\DqcohriX{X^\sharp} = \Dqcohri(X, \ri^+_{X^\sharp}/\pi)$.
\end{definition}

\begin{remark} \label{rslt:explicit-description-of-DqcohriX-on-vstack}
More concretely we can describe the $\infty$-category $\Dqcohri(X, \Lambda)$ as follows: Choose any v-hypercover $Y_\bullet \to X$ by perfectoid spaces $Y_n$ such that each $Y_n = \bigdunion_{i\in I_n} Y_{n,i}$ is a disjoint union of totally disconnected spaces $Y_{n,i}$. Then
\begin{align*}
	\Dqcohri(X, \Lambda) = \varprojlim_{n\in\Delta} \prod_{i\in I_n} \Dqcohri(\Lambda(Y_{n,i})).
\end{align*}
In particular $\Dqcohri(X, \Lambda)$ is stable and admits a symmetric monoidal structure.
\end{remark}

In order to circumvent set-theoretic issues later on, let us also introduce a variant $\Dqcohri(X,\Lambda)_\kappa$ of $\Dqcohri(X,\Lambda)$ for every solid cutoff cardinal (see \cref{def:solid-cutoff-cardinal}) $\kappa$. The main advantage of $\Dqcohri(X,\Lambda)_\kappa$ over $\Dqcohri(X,\Lambda)$ is that the former is presentable (see \cref{rslt:properties-of-kappa-condensed-ri+-modules} below). Moreover, in the case $\kappa = \omega$ the resulting $\infty$-category $\Dqcohri(X,\Lambda)_\omega$ of \emph{discrete} quasicoherent $\Lambda$-modules enjoys a more explicit description in terms of actual sheaves on $X$, as we will discuss in \cref{sec:ri-pi.bd-and-disc}.

\begin{definition}
Given $X \in \vStacksCoeff$ and a solid cutoff cardinal $\kappa$ we denote
\begin{align*}
	\Dqcohri(X,\Lambda)_\kappa \subset \Dqcohri(X,\Lambda)
\end{align*}
the full subcategory of those $\mathcal M \in \Dqcohri(X,\Lambda)$ whose restriction to every totally disconnected space $Y \in X_\vsite$ lies in $\Dqcohri(\Lambda(Y))_\kappa \subset \Dqcohri(\Lambda(X))$. We call $\Dqcohri(X,\Lambda)_\kappa$ the $\infty$-category of \emph{$\kappa$-condensed sheaves} in $\Dqcohri(X,\Lambda)$. Moreover, the objects of $\Dqcohri(X,\Lambda)_\omega$ are called the \emph{discrete sheaves} in $\Dqcohri(X,\Lambda)$.
\end{definition}

\begin{remark}
The reader is invited to skip the following \cref{rslt:def-of-Dqcohri-kappa}, as it is only a subtle juggle of cardinalities and not very enlightening. If we had chosen a fixed solid cutoff cardinal $\kappa$ and worked with $\Dqcohri(X,\Lambda)_\kappa$ instead of $\Dqcohri(X,\Lambda)$ throughout, then \cref{rslt:def-of-Dqcohri-kappa} would not be needed.
\end{remark}

\begin{lemma} \label{rslt:def-of-Dqcohri-kappa}
Let $\kappa$ be a solid cutoff cardinal.
\begin{lemenum}
	\item The functor
	\begin{align*}
		(\vStacksCoeff)^\opp \to \infcatinf, \qquad (X, \Lambda) \mapsto \Dqcohri(X, \Lambda)_\kappa
	\end{align*}
	is a hypercomplete v-sheaf of $\infty$-categories on $\vStacksCoeff$ which on all those $(X, \Lambda) \in \AffPerfCoeff$ with $X$ being of weakly perfectly finite type over some totally disconnected space takes the form
	\begin{align*}
		\Dqcohri(X, \Lambda)_\kappa = \Dqcohri(\Lambda(X))_\kappa.
	\end{align*}

	\item \label{rslt:properties-of-kappa-condensed-ri+-modules} Let $X \in \vStacksCoeff$. Then $\Dqcohri(X,\Lambda)_\kappa$ is presentable and it is stable under all small colimits and the symmetric monoidal structure on $\Dqcohri(X,\Lambda)$. Moreover, for every cardinal $\lambda$ there is a cofinal class of pairs $\kappa_0 \le \kappa_1$ of solid cutoff cardinals such that every $\lambda$-small limit of objects in $\Dqcohri(X,\Lambda)_{\kappa_0}$ lies in $\Dqcohri(X,\Lambda)_{\kappa_1}$.
\end{lemenum}
\end{lemma}
\begin{proof}
Since $\Dqcohri(-,\Lambda)$ is a hypercomplete v-sheaf we can formally reduce (i) to the case of affinoid perfectoid spaces, where the claim boils down to showing that all the descent results from \cref{sec:ri-pi.descent-on-affperfd} also hold for $\Dqcohri(A^+/\pi)_\kappa$ in place of $\Dqcohri(A^+/\pi)$. For $\kappa > \omega$ this follows from the fact that the forgetful functor along any map $(A^+/\pi)^a_\solid \to (B^+/\pi)^a_\solid$ preserves $\Dqcohri(-)_\kappa \subset \Dqcohri(-)$ (see \cref{rslt:properties-of-solid-cutoff-cardinals}) and that $\Dqcohri(-)_\kappa$ is stable under totalizations in $\Dqcohri(-)$. For $\kappa = \omega$ we use \cref{rslt:descendable-implies-discretely-descendable} to see that if $(A^+/\pi)^a_\solid \to (B^+/\pi)^a_\solid$ is descendable then it implies in particular descent of discrete modules, so that all the descent results from \cref{sec:ri-pi.descent-on-affperfd} can easily be adapted to the case of discrete modules (in fact, many of the arguments can be simplified a lot if we disregard non-discrete modules). This proves (i).

We now prove (ii), so let $X \in \vStacksCoeff$ be given. Pick a hypercover $Y_\bullet \to X$ as in \cref{rslt:explicit-description-of-DqcohriX-on-vstack}, so that each $Y_n$ is of the form $Y_n = \bigdunion_{i\in I_n} Y_{n,i}$ with all $Y_{n,i}$ totally disconnected. Then
\begin{align*}
	\Dqcohri(X,\Lambda)_\kappa = \varprojlim_{n\in\Delta} \prod_{i\in I_n} \Dqcohri(\Lambda(Y_{n,i}))_\kappa,
\end{align*}
so since presentable $\infty$-categories are stable under limits (see \cite[Proposition 5.5.3.13]{lurie-higher-topos-theory}) it follows that $\Dqcohri(X,\Lambda)_\kappa$ is presentable. It is clear that $\Dqcohri(X,\Lambda)_\kappa$ is stable under all small colimits and the symmetric monoidal structure in $\Dqcohri(X,\Lambda)$. It remains to prove the claim about limits, so let the cardinal $\lambda$ be given. There is a cofinal class of solid cutoff cardinals $\kappa_0$ such that $\Dqcohri(\Lambda(Y_{n,i}))_{\kappa_0}$ is stable under $\lambda$-small limits in $\Dqcohri(\Lambda(Y_{n,i}))$; indeed, by \cref{rslt:properties-of-solid-cutoff-cardinals} the requirement on $\kappa_0$ is that $\D^a(\Lambda(Y_{n,i}))_{\kappa_0} \subset \D^a(\Lambda(Y_{n,i}))$ is stable under $\lambda$-small limits, so by \cref{rslt:condensed-objects-in-presentable-monoidal-cat} there is a cofinal class of such $\kappa_0$, which by the proof of \cite[Lemma 4.1]{etale-cohomology-of-diamonds} can be chosen to be solid cutoff cardinals. Now fix such a $\kappa_0$. We denote by $\Dqcohri(Y_\bullet,\Lambda) \to \Delta$ the coCartesian fibration classifying the functor $Y_n \mapsto \Dqcohri(Y_n,\Lambda)$. Let $\mathcal C$ be the $\infty$-category of sections $\Delta \to \Dqcohri(Y_\bullet,\Lambda)$, so that $\Dqcohri(X,\Lambda)$ can be identified with the full subcategory of $\mathcal C$ spanned by those sections which map every edge in $\Delta$ to a coCartesian edge in $\Dqcohri(Y_\bullet,\Lambda)$. For every solid cutoff cardinal $\kappa'$ we similarly get $\Dqcohri(Y_\bullet,\Lambda)_{\kappa'}$ and $\mathcal C_{\kappa'}$. By \cite[Proposition 5.5.3.17]{lurie-higher-topos-theory} $\mathcal C_{\kappa'}$ is presentable, so by the adjoint functor theorem the inclusion $\Dqcohri(X,\Lambda)_{\kappa'} \injto \mathcal C_{\kappa'}$ has a right adjoint $r_{\kappa'}$. Hence a limit in $\Dqcohri(X,\Lambda)_{\kappa'}$ is computed as a limit in $\mathcal C_{\kappa'}$ (which is computed pointwise) followed by $r_{\kappa'}$. Altogether this means that we only need to show that there is a solid cutoff cardinal $\kappa_1 \ge \kappa_0$ such that for every $\mathcal M \in \mathcal C_{\kappa_0}$ and every solid cutoff cardinal $\kappa' \ge \kappa_1$ have $r_{\kappa'}(\mathcal M) \in \Dqcohri(X,\Lambda)_{\kappa_1}$.

To prove the claim on $r_{\kappa'}$ we first show that this functor preserves certain filtered colimits. Choose an uncountable regular cardinal $\tau'$ such that $\sum_n \abs I_n < \tau'$ and pick another regular cardinal $\tau''$ such that $\tau'' \gg \tau'$. Since $\tau'$ is uncountable, each $\Dqcohri(\Lambda(Y_{n,i}))_{\kappa'}$ is $\tau'$-compactly generated, so by \cref{rslt:tau-small-limits-in-PrL-tau} $\Dqcohri(X,\Lambda)_{\kappa'}$ is $\tau'$-compactly generated with $\mathcal P \in \Dqcohri(X,\Lambda)_{\kappa'}$ being $\tau'$-compact if and only if each component $\mathcal P_{n,i} \in \Dqcohri(\Lambda(Y_{n,i}))$ is $\tau'$-compact. One checks that such a $\mathcal P$ is also $\tau'$-compact in $\mathcal C_{\kappa'}$ (e.g. by the proof of \cite[Proposition 5.4.7.11]{lurie-higher-topos-theory} one can inductively construct $\mathcal C_{\kappa'}$ from certain $\tau'$-small limits, for which we can apply \cref{rslt:tau-small-limits-in-PrL-tau} to compute the $\tau'$-compact objects; use \cite[Proposition 5.3.5.15]{lurie-higher-topos-theory} to handle the occuring term of $\Fun(\Delta^1, \mathcal D)$). By the proof of \cite[Proposition 5.4.7.7]{lurie-higher-topos-theory} this implies that $r_{\kappa'}$ preserves all $\tau''$-filtered colimits.

By the proof of \cref{rslt:characterize-condensed-objects-via-limit-property} we can choose a regular cardinal $\tau \ge \kappa_0$ such that for every $M \in \Dqcohri(\Lambda(Y_{n,i}))_{\kappa_0}$ and every $\tau$-cofiltered limit $S = \varprojlim_i S_i$ of profinite sets we have $M(S) = \varinjlim_i M(S_i)$. By \cref{rslt:characterize-condensed-objects-via-limit-property} we can then also choose a solid cutoff cardinal $\kappa_1 \ge \tau$ such that every $M \in \Dqcohri(\Lambda(Y_{n,i}))$ satisfying $M(S) = \varinjlim_i M(S_i)$ for all $\tau$-cofiltered limits $S = \varprojlim_i S_i$ is necessarily $\kappa_1$-condensed. We can assume $\tau \ge \tau''$, so that $r_{\kappa'}$ preserves $\tau$-filtered colimits for every solid cutoff cardinal $\kappa' \ge \kappa_1$. Now let $\mathcal M \in \mathcal C_{\kappa_0}$ and $\kappa' \ge \kappa_1$ be given. We claim that $r_{\kappa'}(\mathcal M) \in \Dqcohri(X,\Lambda)_{\kappa_1}$. Note first that for every $S = \varprojlim_i S_i$ as above (with all $S$ and $S_i$ being $\kappa'$-small) we have $\IHom(\Lambda^a_\solid[S], \mathcal M) = \varinjlim_i \IHom(\Lambda^a_\solid[S_i], \mathcal M)$, where the $\IHom$ is computed in $\mathcal C_{\kappa'}$ (i.e. componentwise in each $\Dqcohri(\Lambda(Y_{n,i}))_{\kappa'}$). By applying $r_{\kappa'}$ and commuting right adjoints we obtain $\IHom(\Lambda^a_\solid[S], r_{\kappa'}(\mathcal M)) = \varinjlim_i \IHom(\Lambda^a_\solid[S_i], r_{\kappa'}(\mathcal M))$ in $\Dqcohri(X,\Lambda)_{\kappa'}$. A similar property holds for $r_{\kappa_1}(\mathcal M)$ and moreover the natural morphism $r_{\kappa_1}(\mathcal M) \to r_{\kappa'}(\mathcal M)$ induces isomorphisms $\IHom(\Lambda^a_\solid[S], r_{\kappa_1}(\mathcal M))_{\kappa_1} = \IHom(\Lambda^a_\solid[S], r_{\kappa_1}(\mathcal M))_{\kappa_1}$ for all $\kappa_1$-small profinite sets $S$ (here we use implicitly that the componentwise functor $(-)_{\kappa_1}$ on $\mathcal C_{\kappa'}$ preserves coCartesian sections, which by the usual reductions boils down to the claim $M_{\kappa_1} \tensor_{A_\solid} A[T]_\solid = (M \tensor_{A_\solid} A[T]_\solid)_{\kappa_1}$ for all finite-type classical $\Z$-algebras $A$ and all $M \in \D_\solid(A)_{\kappa'}$; this can be checked on generators, where it follows from the fact that $- \tensor_{A_\solid} A[T]_\solid$ preserves limits on $\kappa_1$-condensed objects, see \cref{rslt:finite-type-solid-pullback-preserves-limits}). Since every $\kappa'$-small profinite set can be written as a $\tau$-cofiltered limit of $\kappa_1$-small profinite sets we deduce from the above colimit property that
\begin{align*}
	\IHom(\Lambda^a_\solid[S], r_{\kappa_1}(\mathcal M))_{\kappa_1} = \IHom(\Lambda^a_\solid[S], r_{\kappa_1}(\mathcal M))_{\kappa_1}
\end{align*}
for all $\kappa'$-small profinite sets $S$. Both sides can be computed componentwise in each $\Dqcohri(\Lambda(Y_{n,i})$ (as before this essentially follows from the fact that $- \tensor_{A_\solid} A[T]_\solid$ preserves limits), so we deduce $r_{\kappa_1}(\mathcal M) = r_{\kappa'}(\mathcal M)$, as desired.
\end{proof}

\begin{definition}
Let $X \in \vStacksCoeff$ and let $\kappa$ be a solid cutoff cardinal. Since $\Dqcohri(X,\Lambda)_\kappa$ is presentable by \cref{rslt:properties-of-kappa-condensed-ri+-modules}, by \cite[Remark 5.5.2.10]{lurie-higher-topos-theory} the inclusion $\Dqcohri(X,\Lambda)_\kappa \injto \Dqcohri(X,\Lambda)$ admits a right adjoint which we denote by
\begin{align*}
	(-)_\kappa\colon \Dqcohri(X,\Lambda) \to \Dqcohri(X,\Lambda)_\kappa, \qquad \mathcal M \mapsto \mathcal M_\kappa.
\end{align*}
In the case $\kappa = \omega$ we call the functor $(-)_\omega$ the \emph{discretization}.
\end{definition}

As a first application of $\Dqcohri(X,\Lambda)_\kappa \subset \Dqcohri(X,\Lambda)$ let us construct two right adjoint functors: the internal hom and the pushforward.

\begin{corollary} \label{rslt:existence-of-IHom-and-pushforward-on-Dqcohri}
\begin{corenum}
	\item \label{rslt:existence-of-IHom-on-Dqcohri} For every $X \in \vStacksCoeff$ the $\infty$-category $\Dqcohri(X,\Lambda)$ is stable, admits all small limits and colimits, and comes equipped with a natural closed symmetric monoidal structure.

	\item \label{rslt:existence-of-pushforward-on-Dqcohri} For every morphism $f\colon Y \to X$ in $\vStacksCoeff$ the restriction functor $f^*\colon \Dqcohri(X,\Lambda) \to \Dqcohri(Y,\Lambda)$ is symmetric monoidal and preserves all small colimits. It admits a right adjoint $f_*$.
\end{corenum}
\end{corollary}
\begin{proof}
We can construct $\Dqcohri(-,\Lambda)$ as a sheaf of symmetric monoidal $\infty$-categories (using that $\Dqcohri(\Lambda(-))$ is a functor to symmetric monoidal $\infty$-categories), which immediately produces a symmetric monoidal structure on $\Dqcohri(X,\Lambda)$ such that all pullback functors are symmetric monoidal. It follows from \cref{rslt:properties-of-kappa-condensed-ri+-modules} that $\Dqcohri(X,\Lambda)$ has all small limits and colimits.

The internal hom and pushforward can be constructed as follows: Using notation from the proof of \cref{rslt:properties-of-kappa-condensed-ri+-modules} it is easy to construct both functors on the $\infty$-categories $\mathcal C_\kappa$ of not necessarily coCartesian sections by an explicit componentwise definition. Then the actual functors are obtained by composing the just constructed functors with $r_\kappa\colon \mathcal C_\kappa \to \Dqcohri(X,\Lambda)_\kappa$. By the proof of \cref{rslt:properties-of-kappa-condensed-ri+-modules} the image of any fixed $\mathcal M \in \mathcal C_\kappa$ under $r_\kappa$ is independent of $\kappa$ for large enough $\kappa$, which allows us to construct the desired functors on all of $\Dqcohri(-,\Lambda)$.
\end{proof}

Using \cref{rslt:existence-of-IHom-and-pushforward-on-Dqcohri} we immediately get the first 4 functors of the desired 6-functor formalism:

\begin{definition} \label{def:4-functors}
Let $X$ be a small v-stack with integral torsion coefficients $\Lambda$.
\begin{defenum}
	\item We denote by
	\begin{align*}
		- \tensor -\colon \Dqcohri(X,\Lambda) \cprod \Dqcohri(X,\Lambda) \to \Dqcohri(X,\Lambda)
	\end{align*}
	the functor induced by the symmetric monoidal structure on $\Dqcohri(X, \Lambda)$ as constructed in \cref{rslt:existence-of-IHom-on-Dqcohri}. It preserves all small colimits in both arguments.

	\item By \cref{rslt:existence-of-IHom-on-Dqcohri} the symmetric monoidal structure $\tensor$ is closed. We denote
	\begin{align*}
		\IHom(-,-)\colon \Dqcohri(X,\Lambda)^\opp \cprod \Dqcohri(X,\Lambda) \to \Dqcohri(X,\Lambda)
	\end{align*}
	the internal hom functor.

	\item For a map $f\colon Y \to X$ of small v-stacks, the \emph{pullback functor}
	\begin{align*}
		f^*\colon \Dqcohri(X, \Lambda) \to \Dqcohri(Y, \Lambda)
	\end{align*}
	is defined to be the ``restriction'' morphism inherent in the sheaf $Y \to \Dqcohri(Y, \Lambda)$ for $Y \in X_\vsite$. By \cref{rslt:existence-of-pushforward-on-Dqcohri} $f^*$ is symmetric monoidal and preserves all small colimits.

	\item For a map $f\colon Y \to X$ of small v-stacks, the \emph{pushforward functor}
	\begin{align*}
		f_*\colon \Dqcohri(Y, \Lambda) \to \Dqcohri(X, \Lambda)
	\end{align*}
	is defined to be the right adjoint of $f^*$, as constructed in \cref{rslt:existence-of-pushforward-on-Dqcohri}.
\end{defenum}
\end{definition}

A natural question to ask is in what sense the objects $\mathcal M \in \Dqcohri(X, \Lambda)$ can be regarded as actual sheaves on $X$. Note that if $X \in X_\vsite^\Lambda$ (equivalently $X_\vsite^\Lambda = X_\vsite$) then $\Dqcohri(X, \Lambda)$ is the sheafification of the presheaf $Y \mapsto \Dqcohri(\Lambda(Y))$ on $X_\vsite$. In particular we obtain the following two functors:

\begin{definition}
Let $(X, \Lambda) \in \vStacksCoeff$ such that $X \in X_\vsite^\Lambda$.
\begin{defenum}
	\item We denote
	\begin{align*}
		\widetilde{(-)}\colon \Dqcohri(\Lambda(X)) \to \Dqcohri(X, \Lambda), \qquad M \mapsto \widetilde M
	\end{align*}
	the functor obtained on global sections from the morphism of the presheaf $\Dqcohri(\Lambda(-))$ to its sheafification, as explained above.

	\item We denote
	\begin{align*}
		\Gamma(X, -)\colon \Dqcohri(X, \Lambda) \to \Dqcohri(\Lambda(X)), \qquad \mathcal M \mapsto \Gamma(X, \mathcal M)
	\end{align*}
	a right adjoint of $\widetilde{(-)}$. (This right adjoint can for example be constructed explicitly using \cref{rslt:explicit-description-of-DqcohriX-on-vstack}.)
\end{defenum}
\end{definition}

If we have a morphism of integral torsion coefficients on a small v-stack, then there are associated base-change and forgetful functors:

\begin{lemma} \label{rslt:Dqcohri-change-of-coefficients}
Let $X$ be a small v-stack and $\alpha\colon \Lambda \to \Lambda'$ a morphism of integral torsion coefficients on $X$. Then there is a natural adjoint pair of functors
\begin{align*}
	- \tensor_{\Lambda^a_\solid} \Lambda'^a_\solid\colon \ \Dqcohri(X, \Lambda) \rightleftarrows \Dqcohri(X, \Lambda') \ \noloc \mathrm{forget}
\end{align*}
They satisfy the following properties:
\begin{lemenum}
	\item Let $f\colon Y \to X$ be any map of small v-stacks. Then both of the above functors commute with $f^*$ and the forgetful functor commutes with $f_*$. Moreover, the forgetful functor is conservative and preserves all small colimits.

	\item For every $Y \in X_\vsite^\Lambda \isect X_\vsite^{\Lambda'}$ there are natural commuting diagrams
	\begin{center}
		\begin{tikzcd}[column sep = large]
			\Dqcohri(\Lambda(Y)) \arrow[r,"\widetilde{(-)}"] \arrow[d,swap,"- \tensor_{\Lambda(Y)^a_\solid} \Lambda'(Y)^a_\solid"]  & \Dqcohri(Y, \Lambda)  \arrow[d,"- \tensor_{\Lambda^a_\solid} \Lambda'^a_\solid"] \\
			\Dqcohri(\Lambda'(Y)) \arrow[r,"\widetilde{(-)}"] & \Dqcohri(Y, \Lambda')
		\end{tikzcd}
										\qquad
		\begin{tikzcd}[column sep = large]
			\Dqcohri(\Lambda(Y))  & \Dqcohri(Y, \Lambda) \arrow[l,swap,"{\Gamma(Y, -)}"] \\
			\Dqcohri(\Lambda'(Y)) \arrow[u,"\mathrm{forget}"] & \Dqcohri(Y, \Lambda') \arrow[l,swap,"{\Gamma(Y, -)}"] \arrow[u,"\mathrm{forget}"]
		\end{tikzcd}
	\end{center}
\end{lemenum}
\end{lemma}
\begin{proof}
The base-change functor defines a morphism of sheaves $\Dqcohri(-, \Lambda) \to \Dqcohri(-, \Lambda')$ on the full subcategory of $X_\vsite$ consisting of the totally disconnected spaces, hence by \cref{rslt:sheaves-on-basis-equiv-sheaves-on-whole-site} extends to a functor $- \tensor_{\Lambda^a_\solid} \Lambda^a_\solid$. Clearly this functor satisfies (i) and (ii) and preserves all small colimits. It thus admits a right adjoint (e.g. this follows by the adjoint functor theorem on $\Dqcohri(-,-)_\kappa$ for all solid cutoff cardinals $\kappa$, because this $\infty$-category is presentable by \cref{rslt:def-of-Dqcohri-kappa}, and one easily sees that the resulting functor is independent of $\kappa$). Using base-change (see \cref{rslt:solid-base-change-for-discrete-rings}) one sees that this forgetful functor commutes with $f^*$ (use the explicit description in \cref{rslt:explicit-description-of-DqcohriX-on-vstack} and note that by the base-change, the forgetful functor preserves cocartesian edges). It is clear that the forgetful functor commutes with $f_*$ (pass to left adjoints). Finally, to see that the forgetful functor is conservative and preserves colimits we can argue locally and hence assume that $X$ is a totally disconnected perfectoid space. Then the claim is clear.
\end{proof}

\subsection{Bounded and Discrete Objects} \label{sec:ri-pi.bd-and-disc}

There is in general no natural $t$-structure on $\Dqcohri(X, \Lambda)$ because the maps $\Lambda(X)^a_\solid \to \Lambda(Y)^a_\solid$ are not flat, even if $X$ and $Y$ are totally disconnected and $\Lambda = \ri^+_{X^\sharp}/\pi$. In this subsection we discuss two ways to circumvent this deficit:
\begin{enumerate}[(a)]
	\item The maps $\Lambda(X)^a_\solid \to \Lambda(Y)^a_\solid$ have Tor dimension $\le 1$ if $X$ is totally disconnected (see \cref{rslt:Tor-dim-and-desc-for-int-tor-coeffs-over-tot-disc}). This allows us to speak of (left/right) \emph{bounded} objects in $\Dqcohri(-, \Lambda)$.

	\item The maps $\Lambda(X)^a \to \Lambda(Y)^a$ are flat (for $X$ totally disconnected). This allows us to define a natural $t$-structure on the full subcategory $\D^a(-, \Lambda)_\omega \subset \Dqcohri(-, \Lambda)$ of discrete objects.
\end{enumerate}
We begin with a discussion of bounded objects. They are defined as follows (also cf. \cref{def:almost-bounded-objects-in-D-A}):

\begin{definition}
Let $X$ be a small v-stack with integral torsion coefficients $\Lambda$.
\begin{defenum}
	\item \label{def:Dqcohri-bounded-subcategories} We define
	\begin{align*}
		\Dqcohrip(X, \Lambda), \Dqcohrim(X, \Lambda), \Dqcohrib(X, \Lambda) \subset \Dqcohri(X, \Lambda),
	\end{align*}
	to be the full subcategories of those objects which when restricted to any totally disconnected space $Y \in X_\vsite$ lie in the respective full subcategory of $\Dqcohri(\Lambda(Y))$. The objects of these full subcategories are called the \emph{left bounded}, \emph{right bounded} and \emph{bounded} sheaves, respectively.

	\item Suppose $X \in X_\vsite^\Lambda$. Then $\Dqcohri(X, \Lambda)$ is $\Lambda(X)^a$-enriched (see \cref{def:V-a-enriched-categories}) with
	\begin{align*}
		\Hom(\mathcal M, \mathcal N)^a = \tau_{\ge0}\Gamma(X, \IHom(\mathcal M, \mathcal N))_\omega
	\end{align*}
	for all $\mathcal M, \mathcal N \in \Dqcohri(X, \Lambda)$. We say that $\mathcal M \in \Dqcohri(X, \Lambda)$ is \emph{weakly almost left bounded} (resp. \emph{weakly almost right bounded}, resp. \emph{weakly almost bounded}) if for all $\varepsilon \in \mm_\Lambda(X)$ there is a left bounded (resp. right bounded, resp. bounded) sheaf $\mathcal M_\varepsilon \in \Dqcohri(X, \Lambda)$ such that $\mathcal M$ is an $\varepsilon$-retract of $\mathcal M_\varepsilon$. We denote the respective full subcategories by
	\begin{align*}
		\Dqcohriwap(X, \Lambda), \Dqcohriwam(X, \Lambda), \Dqcohriwab(X, \Lambda) \subset \Dqcohri(X, \Lambda).
	\end{align*}
\end{defenum}
\end{definition}

Of course one expects that containment in any of the subcategories of \cref{def:Dqcohri-bounded-subcategories} can be checked on any v-cover, and that on affinoid perfectoid spaces $X$ of weakly perfectly finite type over some totally disconnected space one has $\Dqcohriq(X, \Lambda) = \Dqcohriq(\Lambda(X))$. All of this is true, as the following results show.

\begin{lemma} \label{rslt:v-cover-of-tot-disc-implies-equiv-of-boundedness}
Let $Y \surjto X$ be a v-cover of totally disconnected spaces in $\AffPerfCoeff$. Then for any $? \in \{ +, -, b \}$ and $M \in \Dqcohri(\Lambda(X))$ we have
\begin{align*}
	M \in \Dqcohriq(\Lambda(X)) \quad \iff \quad M \tensor_{\Lambda(X)^a_\solid} \Lambda(Y)^a_\solid \in \Dqcohriq(\Lambda(Y)).
\end{align*}
\end{lemma}
\begin{proof}
The implication from left to right is clear for $? = -$ and follows from \cref{rslt:Tor-dim-and-desc-for-int-tor-coeffs-over-tot-disc} for $? = +$. We now prove the implication from right to left, first in the case $? = +$. Let $Y_\bullet \to X$ be the Čech nerve of $Y \to X$ and denote $M^n = M \tensor_{\Lambda(X)^a_\solid} \Lambda(Y_n)^a_\solid$. Then $M^n \in \Dqcohriq(\Lambda(Y_n))$ for all $n$; in fact $(M^n)_n$ is uniformly bounded to the left by \cref{rslt:Tor-dim-and-desc-for-int-tor-coeffs-over-tot-disc}. Moreover by \cref{rslt:Tor-dim-and-desc-for-int-tor-coeffs-over-tot-disc} we have $M = \Tot(M^\bullet)$. Thus the claim for $? = +$ follows immediately, as totalizations preserve $\D_{\le 0}$. Finally, to prove the implication from right to left in the case $? = -$ we note that $\Lambda(X)^a_\solid \to \Lambda(Y)^a_\solid$ is weakly fs-descendable of index $\le 4$ by \cref{rslt:Tor-dim-and-desc-for-int-tor-coeffs-over-tot-disc}, so the claim follows immediately from \cref{rslt:weakly-fs-descendable-implies-right-boundedness-descends}.
\end{proof}

\begin{proposition}
\begin{propenum}
	\item \label{rslt:bounded-Dqcohri-is-sheaf} For $? \in \{ +, -, b \}$ the presheaf $(X,\Lambda) \mapsto \Dqcohriq(X, \Lambda)$ is a hypercomplete v-sheaf on $\vStacksCoeff$ with
	\begin{align*}
		\Dqcohriq(X, \Lambda) = \Dqcohriq(\Lambda(X))
	\end{align*}
	for every totally disconnected affinoid space $X$. Moreover, $\Dqcohriq$ is preserved under any pullback and $\Dqcohrip$ is preserved under any pushforward.

	\item \label{rslt:right-bounded-Dqcohri-global-sections-criterion} Let $X \in \AffPerfCoeff$ and assume that there is a quasi-pro-étale cover $Y \to X$ by some totally disconnected space $Y$ such that $Y \to X$ is a finite composition of covers in $\AffPerfCoeff$ which on $\Lambda(-)^a_\solid$ are weakly fs-descendable morphisms of finite index. Then there are natural equivalences
	\begin{align*}
		\Dqcohri(X, \Lambda) = \Dqcohri(\Lambda(X)), \qquad \Dqcohrim(X, \Lambda) = \Dqcohrim(\Lambda(X)).
	\end{align*}

	\item \label{rslt:left-bounded-Dqcohri-global-sections-criterion} Let $X \in \AffPerfCoeff$ such that $\widetilde{(-)}\colon \Dqcohri(\Lambda(X)) \isoto \Dqcohri(X, \Lambda)$ is an equivalence. If there is a v-cover $Y \to X$ by some totally disconnected $Y$ such that the map $\Lambda(X)^a_\solid \to \Lambda(Y)^a_\solid$ has finite Tor dimension then
	\begin{align*}
		\Dqcohrip(X, \Lambda) = \Dqcohrip(\Lambda(X)).
	\end{align*}
\end{propenum}
\end{proposition}
\begin{proof}
To prove (i) we have to see that $\Dqcohriq(X, \Lambda) = \varprojlim_{n\in\Delta} \Dqcohriq(Y_n, \Lambda)$ for any v-hypercover $Y_\bullet \to X$ in $\vStacksCoeff$. Since we know v-hyperdescent for $\Dqcohri$ (by \cref{rslt:def-of-qcoh-Lambda-modules}) it is clear that the limit on the right-hand side identifies with a full subcategory $\mathcal C_X$ of $\Dqcohri(X, \Lambda)$, namely the subcategory of those objects in $\Dqcohri(X, \Lambda)$ whose pullback to $Y_n$ lies in $\Dqcohriq(Y_n, \Lambda)$ for all $n \ge 0$. It is thus clear that $\Dqcohriq(X, \Lambda) \subset \mathcal C_X$. To get ``$\supset$'' we reduce easily to the case that $X$ and all $Y_n$ are totally disconnected. Then $\Dqcohri(X, \Lambda) = \Dqcohri(\Lambda(X))$ and $\Dqcohri(Y_n, \Lambda) = \Dqcohri(\Lambda(Y_n))$, so the claim follows immediately from \cref{rslt:v-cover-of-tot-disc-implies-equiv-of-boundedness}. We also deduce from \cref{rslt:v-cover-of-tot-disc-implies-equiv-of-boundedness} that $\Dqcohriq(X, \Lambda) = \Dqcohriq(\Lambda(X))$ for totally disconnected $X$. Finally the claim about pullbacks is clear and the claim about pushforwards follows from the fact that totalizations preserve $\D_{\le0}$.

To prove (ii) let $X$ and $Y$ be given as in the claim. We argue by induction on the number of intermediate covers composing $Y \to X$. We can thus assume that $Y \to X$ factors as $Y \to X' \to X$ such that $X' \to X$ is a quasi-pro-étale cover with $\Lambda(X)^a_\solid \to \Lambda(X')^a_\solid$ a filtered colimit of fs-descendable maps of bounded index, $Y \to X'$ is a composition of such maps, and $\Dqcohri(X', \Lambda) = \Dqcohri(\Lambda(X'))$, $\Dqcohrim(X', \Lambda) = \Dqcohrim(\Lambda(X'))$. By \cref{rslt:Tor-dim-and-desc-for-int-tor-coeffs-over-tot-disc} we can assume that $Y$ is \emph{strictly} totally disconnected. Let $X'_\bullet \to X$ be the Čech nerve of $X' \to X$. Then each $Y_n := Y \cprod_{X'} X'_n$ is (strictly) totally disconnected and $Y_n \surjto X'_n$ is still a nice composition as in the claim, so we have
\begin{align*}
	\Dqcohri(X'_n, \Lambda) = \Dqcohri(\Lambda(X'_n)), \qquad \Dqcohrim(X'_n, \Lambda) = \Dqcohrim(\Lambda(X'_n)).
\end{align*}
In particular, since $\Lambda(X')^a_\solid \to \Lambda(X)^a_\solid$ is weakly fs-descendable we have
\begin{align*}
	\Dqcohri(X, \Lambda) = \varprojlim_{n\in\Delta} \Dqcohri(X'_n, \Lambda) = \varprojlim_{n\in\Delta} \Dqcohri(\Lambda(X'_n)) = \Dqcohri(\Lambda(X))
\end{align*}
by \cref{rslt:pi-0-weakly-descendable-implies-descent}. Since $\Lambda(X')^a_\solid \to \Lambda(X)^a_\solid$ is even weakly fs-descendable \emph{of finite index}, by \cref{rslt:weakly-fs-descendable-implies-right-boundedness-descends} the above equivalence also holds for right-bounded objects. This finishes the proof of (ii).

We now prove (iii) so let $X$ with $\Dqcohri(X, \Lambda) = \Dqcohri(\Lambda(X))$ be given and assume that there is some totally disconnected $Y$ such that $\Lambda(X)^a_\solid \to \Lambda(Y)^a_\solid$ has finite Tor dimension. Let $Y_\bullet \to X$ be any v-hypercover extending $Y \to X$ such that all $Y_n$ are totally disconnected. By (i) we have $\Dqcohrip(X, \Lambda) = \varprojlim_{n\in\Delta} \Dqcohrip(Y_n, \Lambda)$ and so the claim $\Dqcohrip(X, \Lambda) = \Dqcohrip(\Lambda(X))$ follows easily from the assumption on Tor dimension and the fact that totalizations preserve $\D_{\le0}$.
\end{proof}

\begin{lemma} \label{rslt:compute-bounded-Dqcohri-for-qproet-over-Z'-Z}
Let $f\colon Z' \to Z$ be a map of totally disconnected spaces in $\AffPerfCoeff$ with $\dimtrg f < \infty$ and let $X \to \overline{Z'}^{/Z}$ be any quasi-pro-étale map. Then $X$ satisfies the assumptions of \cref{rslt:right-bounded-Dqcohri-global-sections-criterion,rslt:left-bounded-Dqcohri-global-sections-criterion}. In particular the equivalence $\Dqcohri(X, \Lambda) = \Dqcohri(\Lambda(X))$ induces
\begin{align*}
	\Dqcohriq(X, \Lambda) = \Dqcohriq(\Lambda(X))
\end{align*}
for $? \in \{ +, -, b \}$.
\end{lemma}
\begin{proof}
This follows easily from \cref{rslt:Tor-dim-and-desc-for-int-tor-coeffs-on-Z'-over-Z}.
\end{proof}

\begin{proposition} \label{rslt:mod-pi-Tor-dim-for-fin-type-over-tot-disc}
Let $Y \to X$ be a map in $\AffPerfCoeff$ such that $X$ is weakly of perfectly finite type over some totally disconnected space $Z$. Then the map $\Lambda(X)^a_\solid \to \Lambda(Y)^a_\solid$ has Tor dimension $\le n + 1$, where $n$ is the number of perfect generators for $X$ over $Z$.
\end{proposition}
\begin{proof}
Pick a morphism $\ri^+_{Z^\sharp}/\pi \to \Lambda$ for some untilt $Z^\sharp$ of $Z$ and some pseudouniformizer $\pi$ on $Z^\sharp$. By definition of morphisms of integral torsion coefficients we immediately reduce to the case $\Lambda = \ri^+_{Z^\sharp}/\pi$. Let us now change notation slightly by writing $X, Y, Z \in \AffPerfd_\pi$ for $X^\sharp$, $Y^\sharp$ and $Z^\sharp$.

As in the proof of \cref{rslt:mod-pi-Tor-dim-1-over-tot-disc-space} we can reduce the claim to the case that $Y$ is totally disconnected, then we can further reduce (using \cref{rslt:compare-catsldmod-with-sheaves-on-pi-0}) to the case that $X$, $Y$ and $Z$ are connected. Writing $X \injto X' \to Z$, where $X'$ is the compactified relative perfectoid ball over $Z$ (of dimension $n$) and $X \injto X'$ is a closed immersion, the induced map $Y \to X'$ is still pro-étale and $Y = Y \cprod_{X'} X$. This allows us to further reduce to the case $X = X'$. We are now in the setting $Z = \Spa(K, K^+)$ and $Y = \Spa(L, L^+)$ for some perfectoid fields $K$ and $L$ and open and bounded valuation subrings $K^+ \subset K$ and $L^+ \subset L$. Moreover, $X = \Spa(A, A^+)$ with
\begin{align*}
	A = K\langle T_1^{1/p^\infty}, \dots, T_n^{1/p^\infty}\rangle, \qquad A^+ = K^+ + A^{\circ\circ}.
\end{align*}
Note that there are obvious maps $g$ and $h$ which make the following diagram commutative:
\begin{center}\begin{tikzcd}
	Y \arrow[r,"g"] \arrow[dr] & Y \cprod_Z X \arrow[d,"h"]\\
	& X
\end{tikzcd}\end{center}
Now $h$ is a base-change of the map $Y \to Z$, which has mod-$\pi$ Tor dimension $\le 1$ by \cref{rslt:mod-pi-Tor-dim-1-over-tot-disc-space}. Thus it is enough to show that $g$ has Tor dimension $\le n$. Note that
\begin{align*}
	(\ri^+_Y(Y)/\pi)^a_\solid &= (L^+/\pi)^a_\solid = (\ri_L^a/\pi, L^+/\pi)_\solid,\\
	(\ri^+_{Y \cprod_Z X}(Y \cprod_Z X)/\pi)^a_\solid &= (\ri_L^a/\pi[T_1^{1/p^\infty}, \dots, T_n^{1/p^\infty}], L^+/\pi)_\solid,
\end{align*}
so the claim reduces to showing that the map $\ri_L^a/\pi[T_1^{1/p^\infty}, \dots, T_n^{1/p^\infty}] \to \ri_L^{+a}/\pi$ of almost rings has Tor dimension $\le n$. But this map is obtained by successively replacing a coordinate $T_i$ with some $a_i \in \ri_L^{+a}/\pi$ and each of these operations has Tor dimension $1$.
\end{proof}

\begin{corollary} \label{rslt:computation-of-bounded-Dqcohri-if-fin-type-over-tot-disc}
Suppose $X \in \AffPerfCoeff$ is weakly of perfectly finite type over some totally disconnected space. Then $X$ satisfies the assumptions of \cref{rslt:right-bounded-Dqcohri-global-sections-criterion,rslt:left-bounded-Dqcohri-global-sections-criterion}. In particular the equivalence $\Dqcohri(X, \Lambda) = \Dqcohri(\Lambda(X))$ induces
\begin{align*}
	\Dqcohriq(X, \Lambda) = \Dqcohriq(\Lambda(X))
\end{align*}
for $? \in \{ +, -, b \}$.
\end{corollary}
\begin{proof}
This follows from \cref{rslt:fin-type-over-tot-disc-implies-Lambda-desc,rslt:compute-bounded-Dqcohri-for-qproet-over-Z'-Z} (for \cref{rslt:right-bounded-Dqcohri-global-sections-criterion}) and from \cref{rslt:mod-pi-Tor-dim-for-fin-type-over-tot-disc} (for \cref{rslt:left-bounded-Dqcohri-global-sections-criterion}).
\end{proof}

Having defined (left/right) bounded objects in $\Dqcohri(X, \Lambda)$, we now study some of their properties. One of the main advantages of working with left bounded complexes is that they satisfy unconditional base-change, similar to the $\ell$-adic case (cf. \cite[Proposition 17.6]{etale-cohomology-of-diamonds}):

\begin{proposition} \label{rslt:base-change-for-bounded-Dqcohri}
Let
\begin{center}\begin{tikzcd}
	Y' \arrow[r,"g'"] \arrow[d,"f'"] & Y \arrow[d,"f"]\\
	X' \arrow[r,"g"] & X
\end{tikzcd}\end{center}
be a cartesian diagram in $\vStacksCoeff$ and assume that $f$ is qcqs. Then the natural base-change morphism
\begin{align*}
	g^* f_* \isoto f'_* g'^*
\end{align*}
is an equivalence of functors from $\Dqcohrip(Y, \Lambda)$ to $\Dqcohrip(X', \Lambda)$.

If $X \in X_\vsite^\Lambda$, then the above isomorphism extends to weakly almost left-bounded sheaves, i.e. to an isomorphism of functors $\Dqcohriwap(Y, \Lambda) \to \Dqcohriwap(X', \Lambda)$.
\end{proposition}
\begin{proof}
Formal arguments using v-hyperdescent reduce the claim to the case that $X$ and $X'$ are totally disconnected perfectoid spaces. In this case we have $\Dqcohri(X, \Lambda) = \Dqcohri(\Lambda(X))$ and $\Dqcohri(X', \Lambda) = \Dqcohri(\Lambda(X'))$. If $Y$ is also a totally disconnected perfectoid space then we can use \cref{rslt:base-change-over-tot-disc-equiv} to reduce the claim to \cref{rslt:base-change-for-int-tor-coeffs-on-aff-perf}. For general $Y$, since $f$ is qcqs, we can choose a v-hypercover $Y_\bullet \to Y$ by totally disconnected spaces $Y_n$. Then given any $\mathcal M = (M^n)_n \in \Dqcohrip(Y, \Lambda) = \varprojlim_{n\in\Delta} \Dqcohrip(\Lambda(Y_n))$ we have $f_* \mathcal M = \Tot(M^\bullet)$. But by \cref{rslt:Tor-dim-and-desc-for-int-tor-coeffs-over-tot-disc} the cosimplicial object $M^\bullet$ in $\Dqcohri(X, \Lambda)$ is uniformly bounded to the left and $g^*$ has finite Tor dimension; both together imply that $g^*$ commutes with the totalization, which reduces the claim to the case that $Y$ is totally disconnected. This proves the first part of the claim.

Now assume that $X \in X_\vsite^\Lambda$ and let $\mathcal M \in \Dqcohriwap(Y, \Lambda)$ be weakly almost left-bounded. For every $\varepsilon \in \mm_\Lambda(X)$ let $\mathcal M_\varepsilon$ be a left-bounded sheaf on $Y$ such that $\mathcal M$ is an $\varepsilon$-retract of $\mathcal M_\varepsilon$. Then we get a diagram
\begin{center}\begin{tikzcd}
	g^* f_* \mathcal M \arrow[r] \arrow[d,shift left] & f'_* g'^* \mathcal M \arrow[d,shift left]\\
	g^* f_* \mathcal M_\varepsilon \arrow[r,"\sim"] \arrow[u,shift left] & f'_* g'^* \mathcal M_\varepsilon \arrow[u,shift left]
\end{tikzcd}\end{center}
of sheaves on $X'$ (where either the upwards or the downwards vertical maps form a commutative diagram), where the bottom map is an isomorphism by the left-bounded case shown above and the composition of each downward map with its corresponding upward map is multiplication by $\varepsilon$. Pull this diagram back to a cover of $X'$ by totally disconnected spaces. Then a diagram chase shows that on all homotopy groups, the kernel and cokernel of the top horizontal map are killed by $\varepsilon$. As this is true for all $\varepsilon$, the top horizontal map is an isomorphism, as desired.
\end{proof}

The category $\Dqcohri(X, \Lambda)$ on a small v-stack $X$ with integral torsion coefficients $\Lambda$ behaves a lot like the derived category of a certain category of ``\emph{étale} sheaves of $\Lambda^a_\solid$-modules''. In the following we will strengthen this intuition by proving an analogue of \cite[Proposition 14.9]{etale-cohomology-of-diamonds}. In the proof we will need the following argument, which was also used in the proof of \cite[Lemma 11.23]{etale-cohomology-of-diamonds}:

\begin{lemma} \label{rslt:cover-diagram-of-vstacks-by-tot-disc-spaces}
Let $I$ be a partially ordered set and $(X_i)_{i\in I}$ diagram of qcqs small v-stacks. Then there is a cofinal subset $I' \subset I$ and a diagram $(Y_{i,\bullet} \to X_i)_{i \in I'}$ of v-hypercovers with the following properties:
\begin{lemenum}
	\item All $Y_{i,n}$ and $Y_n := \varprojlim_i Y_{i,n}$ are totally disconnected perfectoid spaces.
	\item $Y_\bullet \to X := \varprojlim_i X_i$ is a v-hypercover.
\end{lemenum}
\end{lemma}
\begin{proof}
We can w.l.o.g. assume that $I$ contains a final object $0 \in I$. Let $\mathcal Z$ denote the set of all pairs $(J, (Y_{j,\bullet} \to X_j)_{j\in J})$ where $0 \in J \subset I$ and $(Y_{j,\bullet} \to X_j)_j$ is a diagram of hypercovers such that all $Y_{j,n}$ are totally disconnected perfectoid spaces. We put a partial order on $\mathcal Z$ by saying $(J, (Y_{j,\bullet} \to X_j)_{j\in J}) \le (J', (Y'_{j,\bullet} \to X_j)_{j\in J'})$ if $J \subset J'$ and $(Y_{j,\bullet} \to X_j) = (Y'_{j,\bullet} \to X_j)$ for all $j \in J$. It is clear that every totally ordered subset of $\mathcal Z$ has an upper bound (take the union of the $J$'s), so by Zorn's lemma there is some maximal element $(I', (Y_{i,\bullet} \to X_i)_{i\in I}) \in \mathcal Z$.

For every non-empty $0 \in J \subset I'$ denote $Y_{J,\bullet} := \varprojlim_{j\in J} Y_{j,\bullet}$ and $X_J := \varprojlim_{j\in J} X_j$. Then $X_J$ and and all $Y_{J,n}$ are qcqs and $Y_{J,\bullet} \to X_j$ is a hypercover. Indeed, by decomposing the limit over $J$ into simpler pieces, this reduces to the case of fiber products and to infinite fiber products over $X_0$ resp. $Y_{0,\bullet}$, i.e. to cofiltered limits; using \cite[Lemma 12.11]{etale-cohomology-of-diamonds}, the case of fiber products is handled by \cite[Proposition 12.10]{etale-cohomology-of-diamonds} and cofiltered limits are handled by \cite[Lemma 12.17]{etale-cohomology-of-diamonds}.

Now suppose that $I'$ is not cofinal in $I$. Then there is some $i_0 \in I \setminus I'$ which is not smaller than any element of $I'$. Let $J := \{ i_0 \} \union I'$ and $J' := \{ i \in I' \setst i \le i_0 \}$. Choose a hypercover $Y_{i_0,\bullet} \to X_{i_0}$ refining the hypercover $Y_{J,\bullet} \cprod_{X_J} X_{i_0}$ and such that all $Y_{i_0,n}$ are totally disconnected. This gives a pair $(J, (Y_{j,\bullet} \to X_j)_j) \in \mathcal Z$, contradicting the maximality of $I'$. The only thing that is left to show is that all $Y_n = \varprojlim_{i\in I'} Y_{i,n}$ are totally disconnected. By definition this boils down to showing that every qcqs cover $Y_n = \bigunion_{k=1}^m U_k$ splits; but this cover comes via base-change from a similar cover of some $Y_{i,n}$, which splits as $Y_{i,n}$ is totally disconnected.
\end{proof}

We can now prove the promised analog of \cite[Proposition 14.9]{etale-cohomology-of-diamonds}, showing that $\Lambda^a_\solid$-modules behave a lot like étale sheaves. There is also a version for cohomology, closer to loc. cit., in a more restricted setting, see \cref{rslt:colimit-of-cohomology-qproet-unbounded-version} below.

\begin{proposition} \label{rslt:colim-of-pushforward-pullback}
Let $(X_i)_{i\in I}$ be a cofiltered diagram in $\vStacksCoeff$ with qcqs transition maps and let $X := \varprojlim_i X_i$ be the limit. Assume that $I$ has a final object $0 \in I$ and denote $f_i\colon X_i \to X_0$, $f\colon X \to X_0$ the induced maps. Then for $\mathcal M \in \Dqcohri(X_0, \Lambda)$ the natural map
\begin{align*}
	\varinjlim_i f_{i*} f_i^* \mathcal M \isoto f_* f^* \mathcal M
\end{align*}
is an isomorphism in the following cases:
\begin{propenum}
	\item $\mathcal M$ is left-bounded, i.e. $\mathcal M \in \Dqcohrip(X_0, \Lambda)$.
	\item $X_0 \in X_{0\vsite}^\Lambda$ and $\mathcal M$ is weakly almost left-bounded, i.e. $\mathcal M \in \Dqcohriwap(X_0, \Lambda)$.
	\item All transition maps in the diagram $(X_i)_i$ are quasi-pro-étale.
\end{propenum}
\end{proposition}
\begin{proof}
By \cite[Lemma 12.17]{etale-cohomology-of-diamonds} the map $f$ is qcqs. Thus in any of (i), (ii) and (iii), $f_*$ satisfies base-change along any map $X'_0 \to X_0$ (in (i) and (ii) this is \cref{rslt:base-change-for-bounded-Dqcohri} and in (iii) $f$ is quasi-pro-étale, so base-change is easy, cf. \cref{rslt:base-change-for-Dqcohri-along-qcqs-p-bounded-map}). We can therefore replace $X_0$ by $X'_0$ and $X_i$ by $X_i' := X_0' \cprod_{X_0} X_i$ to assume that $X_0$ is totally disconnected. In case (iii) the claim follows easily from \cref{rslt:computation-of-Dqcohri-in-td-case}.

It remains to prove (i) and (ii). We first prove (i), so assume that $\mathcal M \in \Dqcohri(\Lambda(X_0))$ is left-bounded. Choose a diagram $(Y_{i,\bullet} \surjto X_i)_{i\in I}$ of hypercovers as in \cref{rslt:cover-diagram-of-vstacks-by-tot-disc-spaces}. In particular all $Y_{i,n}$ and $Y_n$ are totally disconnected, so we have
\begin{align*}
	\Dqcohri(X_i, \Lambda) = \varprojlim_{n\in\Delta} \Dqcohri(\Lambda(Y_{i,n})), \qquad \Dqcohri(X, \Lambda) = \varprojlim_{n\in\Delta} \Dqcohri(\Lambda(Y_n))
\end{align*}
Thus the claim boils down to the claim that the natural map
\begin{align*}
	\varinjlim_i \Tot(\mathcal M \tensor_{\Lambda(X_0)^a_\solid} \Lambda(Y_{i,\bullet})^a_\solid) \isoto \Tot(\mathcal M \tensor_{\Lambda(X_0)^a_\solid} \Lambda(Y_\bullet)^a_\solid)
\end{align*}
is an isomorphism of $\Lambda(X_0)^a_\solid$-modules. But all the involved modules are uniformly bounded to the left (by \cref{rslt:Tor-dim-and-desc-for-int-tor-coeffs-over-tot-disc}), so the totalizations are computed via some spectral sequences over $\catsldmoda{\Lambda(X_0)}$. In particular, we can pull the filtered colimit into the left totalization, which immediately shows the desired isomorphism.

We now prove (ii), so assume that $X_0 \in X_{0\vsite}^\Lambda$ and $\mathcal M \in \Dqcohriwap(X_0, \Lambda)$. For every $\varepsilon \in \mm_\Lambda(X_0)$ let $\mathcal M_\varepsilon$ be a left-bounded sheaf on $X_0$ such that $\mathcal M$ is an $\varepsilon$-retract of $\mathcal M_\varepsilon$. We get a diagram (where both the version with the upward maps and the one with the downward maps commutes)
\begin{center}\begin{tikzcd}
	\varinjlim_i f_{i*} f_i^* \mathcal M \arrow[r] \arrow[d,shift left] & f_* f^* \mathcal M \arrow[d,shift left]\\
	\varinjlim_i f_{i*} f_i^* \mathcal M_\varepsilon \arrow[r,"\sim"] \arrow[u,shift left] & f_* f^* \mathcal M_\varepsilon \arrow[u,shift left]
\end{tikzcd}\end{center}
of $\Lambda(X_0)^a_\solid$-modules, where the bottom horizontal map is an isomorphism by case (i) and the composition of each downward map with its corresponding upward map is multiplication by $\varepsilon$. As in the proof of \cref{rslt:base-change-for-bounded-Dqcohri}, a diagram chase shows that on all homotopy groups, the kernel and cokernel of the top horizontal map are killed by $\varepsilon$, proving the claim.
\end{proof}

We now come to the study of discrete objects in $\Dqcohri(X, \Lambda)$, i.e. the full subcategory $\Dqcohri(X, \Lambda)_\omega$. Let us start by collecting some basic properties:

\begin{lemma}
Let $X \in \vStacksCoeff$.
\begin{lemenum}
	\item The stable $\infty$-category $\Dqcohri(X, \Lambda)_\omega$ comes naturally equipped with a left-complete $t$-structure which is compatible with $\Dqcohrip$ and $\Dqcohrim$. For every map $Y \to X$ of small v-stacks the pullback functor $f^*\colon \Dqcohri(X, \Lambda)_\omega \to \Dqcohri(Y, \Lambda)_\omega$ is $t$-exact.

	\item \label{rslt:qcqs-pushforward-preserves-left-bounded-discrete-sheaves} Let $f\colon Y \to X$ be a qcqs map of small v-stacks. Then $f_*\colon \Dqcohrip(Y, \Lambda) \to \Dqcohrip(X, \Lambda)$ preserves discrete sheaves and hence restricts to a functor $\Dqcohrip(Y, \Lambda)_\omega \to \Dqcohrip(X, \Lambda)_\omega$.
\end{lemenum}
\end{lemma}
\begin{proof}
Part (i) follows easily from the fact that for every map $\Spa(A, A^+) \to \Spa(B, B^+)$ of totally disconnected spaces and every pseudouniformizer $\pi \in A^+$ the map $A^+/\pi \to B^+/\pi$ is flat (e.g. define $\Dqcohrige0(X, \Lambda)_\omega$ to consist of those $\mathcal M \in \Dqcohri(X, \Lambda)_\omega$ whose pullback to any totally disconnected space $U = \Spa(A, A^+)$ lies in $\D^a_{\ge0}(\Lambda(U))$. It is clear that pullbacks are $t$-exact with respect to this $t$-structure.

To prove (ii) we can immediately reduce to the case that $X$ is a totally disconnected space. Then choose a hypercover of $Y$ by totally disconnected spaces and compute $f_*$ as the totalization along this hypercover. The claim follows by observing that $\Dqcohri(\Lambda(X))_\omega \subset \Dqcohri(\Lambda(X))$ is stable under uniformly left-bounded totalizations (as these can be computed using a spectral sequence).
\end{proof}

\begin{remark}
If the map $f\colon Y \to X$ in \cref{rslt:qcqs-pushforward-preserves-left-bounded-discrete-sheaves} is $p$-bounded (see \cref{def:p-bounded} below) then the conclusion of the result even holds for unbounded complexes. This follows by the usual Postnikov tower argument, using that $f_*$ has finite cohomological dimension.
\end{remark}

It turns out that the stable $\infty$-category $\Dqcohri(X, \Lambda)_\omega$ admits a more explicit description in terms of actual sheaves on $X$, at least if $X$ is a locally spatial diamond. In order to provide this explicit description, we need to work with sheaves of almost modules. They can be defined as follows:

\begin{definition}
Let $X$ be a locally spatial diamond and $\Lambda$ an almost setup on $X$ with $X \in X_\vsite^\Lambda$.
\begin{defenum}
	\item We denote $\Shv(X_\et, \Lambda)$ the derived $\infty$-category of classical sheaves of $\Lambda$-modules on $X_\et$. We similarly define $\Shv(X_\qproet, \Lambda)$.

	\item A \emph{classical étale sheaf of $\Lambda^a$-modules on $X$} is a sheaf $\mathcal M$ on $X_\et$ with values in the abelian category $\D(\Lambda(X)^a)^\heartsuit_\omega$ of classical $\Lambda(X)^a$-modules together with a morphism $\Lambda^a \cprod \mathcal M \to \mathcal M$ which gives every $\mathcal M(U)$ the structure of a $\Lambda(U)^a$-module. We denote the derived $\infty$-category of classical étale sheaves of $\Lambda^a$-modules on $X$ by
	\begin{align*}
		\Shv(X_\et, \Lambda^a).
	\end{align*}
	Given a map $f\colon Y \to X$ from another locally spatial diamond $Y$ there is a natural pair of adjoint functors
	\begin{align*}
		f^*\colon \Shv(X_\et, \Lambda^a) \rightleftarrows \Shv(Y_\et, \Lambda^a) \noloc f_*,
	\end{align*}
	see \cite[Sections 06YV, 07A5]{stacks-project}. Similar definitions apply to $X_\qproet$ in place of $X_\et$.
\end{defenum}
\end{definition}

As in \cref{sec:andesc.almmath} we can view $\Shv(X, \Lambda^a)$ as a localization of $\Shv(X, \Lambda)$. More precisely, we have the following result:

\begin{lemma} \label{rslt:general-properties-of-almost-localization-for-sheaves-on-diamond}
Let $X$ be a locally spatial diamond and $\Lambda$ an almost setup on $X$ with $X \in X_\vsite^\Lambda$.
\begin{lemenum}
	\item There is a $t$-exact almost localization functor
	\begin{align*}
		(-)^a\colon \Shv(X_\et, \Lambda) \to \Shv(X_\et, \Lambda^a), \qquad \mathcal M \mapsto \mathcal M^a
	\end{align*}
	which can be computed as $\mathcal M^a(U) = \mathcal M(U)^a$. It preserves all small limits and colimits and commutes with pullbacks and pushforwards along maps of locally spatial diamonds.

	\item The localization functor $(-)^a$ admits a left $t$-exact fully faithful right adjoint
	\begin{align*}
		(-)_*\colon \Shv(X_\et, \Lambda^a) \to \Shv(X_\et, \Lambda), \qquad \mathcal M \mapsto \mathcal M_*
	\end{align*}
	which can be computed as $\mathcal M_*(U) = \mathcal M(U)_*$. It commutes with pushforwards along maps of locally spatial diamonds.

	\item The localization functor $(-)^a$ admits a $t$-exact fully faithful left adjoint
	\begin{align*}
		(-)_!\colon \Shv(X_\et, \Lambda^a) \to \Shv(X_\et, \Lambda), \qquad \mathcal M \mapsto \mathcal M_!
	\end{align*}
	which can be computed as $\mathcal M_!(U) = \mathcal M(U)_!$. It commutes with pullbacks along maps of locally spatial diamonds.
\end{lemenum}
\end{lemma}
\begin{proof}
The functor $(-)^a$ is easily constructed on the hearts by the formula $\mathcal M^a(U) = \mathcal M(U)^a$. One immediately checks that this functor is exact and hence extends to a $t$-exact functor on the full derived $\infty$-categories (with the same formula). One checks similarly (using the properties in \cref{rslt:properties-of-almost-V-modules}) that the functors $(-)_*$ and $(-)_!$ defined by the given formulas preserve the sheaf property and are thus well-defined. They provide the desired adjoints and by the explicit formulas it is clear that they are fully faithful.

One sees immediately that the functor $(-)^a$ commutes with pullbacks (e.g. by explicit computation using K-flat resolutions). Since $(-)_!$ is $t$-exact, the functor $(-)^a$ preserves K-injective complexes and hence it is easy to see that $(-)^a$ commutes with pushforwards. The other claims follow formally from adjunctions.
\end{proof}

We now come to the main definition in the study of discrete objects in $\Dqcohri(X, \Lambda)$, the so-called overconvergent sheaves:

\begin{definition}
Let $X \in \vStacksCoeff$ be a locally spatial diamond. A sheaf $\mathcal M \in \Shv(X_\et, \Lambda^a)$ is called \emph{overconvergent} if it satisfies the following property: Let $x = \Spa(C, C^+) \to X$ be a quasi-pro-étale map from an algebraically closed perfectoid field $C$ with open and bounded valuation subring $C^+$ and let $x^\circ\colon \Spa(C, C^\circ) \to X$ be the induced map; then the natural morphism of stalks $\mathcal M_x \isoto \mathcal M_{x^\circ}$ is an isomorphism. We denote
\begin{align*}
	\Shv(X_\et, \Lambda^a)^\oc \subset \Shv(X_\et, \Lambda^a)
\end{align*}
the full subcategory of overconvergent sheaves.
\end{definition}

The following results show that overconvergent sheaves behave nicely and that they indeed provide a description of $\Dqcohri(X, \Lambda)_\omega$. In particular this shows that our theory of quasicoherent $\Lambda^a$-modules can be used to compute cohomology in classical rigid-analytic geometry (see \cref{rslt:discrete-pushforward-computes-classical-cohomology} below).

\begin{lemma} \label{rslt:basic-properties-of-overconvergent-sheaves}
Let $X \in \vStacksCoeff$ be a locally spatial diamond and assume that $X \in X_\vsite^\Lambda$.
\begin{lemenum}
	\item A sheaf $\mathcal M \in \Shv(X_\et, \Lambda^a)$ is overconvergent if and only if all the classical sheaves $\pi_n \mathcal M$ are overconvergent. In particular $\Shv(X_\et, \Lambda^a)^\oc$ comes naturally equipped with a $t$-structure.

	\item If $X = \Spa(A, A^+)$ is a strictly totally disconnected perfectoid space then there is a natural equivalence of $\infty$-categories $\Shv(X_\et, \Lambda^a)^\oc = \D(\Lambda(X)^a)_\omega$.

	\item \label{rslt:qcqs-pushforward-preserves-overconvergent-sheaves} For every map $f\colon Y \to X$ of locally spatial diamonds, the pullback $f^*\colon \Shv(X_\et, \Lambda^a) \to \Shv(Y_\et, \Lambda^a)$ preserves overconvergent sheaves. If $f$ is qcqs then the pushforward $f_*$ also preserves overconvergent sheaves.
\end{lemenum}
\end{lemma}
\begin{proof}
Part (i) follows immediately from the fact that the stalk functor $\mathcal M \mapsto \mathcal M_x$ is $t$-exact and that isomorphisms can be checked on cohomology. We now prove (ii) so assume that $X = \Spa(A, A^+)$ is a strictly totally disconnected perfectoid space. Then there is a natural adjunction
\begin{align*}
	\widetilde{(-)}\colon \Dqcohri(\Lambda(X)^a)_\omega \rightleftarrows \Shv(X_\et, \Lambda^a) \noloc \Gamma(X, -).
\end{align*}
Here $\Gamma(X, -)$ is just the evaluation of a sheaf at $X$ and $\widetilde{(-)}$ exists for example by the adjoint functor theorem. We can equivalently view $\widetilde{(-)}$ as the pullback along the map of ringed sites $(X_\et, \Lambda^a) \to (*, \Lambda(X)^a)$. In fact $\widetilde{(-)}$ is the left derived functor of its restriction to the heart, where it is the left adjoint to the evaluation-at-$X$ functor and is thus computed by $\widetilde M(U) = M \tensor_{\Lambda(X)^a} \Lambda(U)^a$. Note that $\Lambda(U)^a$ is flat over $\Lambda(X)^a$ (because so is $\ri^+_X(U)/\pi$ over $\ri^+_X(X)/\pi$) which shows that $\widetilde{(-)}$ is $t$-exact. Since $X$ is strictly totally disconnected, $\Gamma(X, -)$ is also $t$-exact. Hence $\Gamma(X, \widetilde M) = M$ for every $M \in \Dqcohri(\Lambda(X)^a)$ (by $t$-exactness this can be checked on the hearts, where it is clear). This implies that $\widetilde{(-)}$ is fully faithful.

The functor $\widetilde{(-)}$ clearly factors over $\Shv(X_\et, \Lambda^a)^\oc$, so to finish the proof of (ii) it only remains to see that every overconvergent $\Lambda^a$-module lies in the essential image of $\widetilde{(-)}$. Let $\mathcal M \in \Shv(X_\et, \Lambda^a)^\oc$ and $M := \Gamma(X, \mathcal M)$. We need to show that the natural map $\widetilde M \to \mathcal M$ is an isomorphism. This can be checked on stalks, so fix any $x \in X$; we need to show that the map $\widetilde M_x \to \mathcal M_x$ is an isomorphism. Since both $\widetilde M$ and $\mathcal M$ are overconvergent, we can reduce to the case that $x$ is the closed point of its connected component in $X$. Then $x$ is the intersection of all clopen neighbourhoods which contain it, so we reduce to showing that for every clopen subset $U \subset X$ the natural map $\widetilde M(U) \to \mathcal M(U)$ is an isomorphism. But since $U$ is clopen we have $X = U \dunion V$ for some clopen $V \subset X$ and hence $\widetilde M(U) \cprod \widetilde M(V) = M = \mathcal M(U) \cprod \mathcal M(V)$. The claim follows easily.

We now prove (iii), so let $f\colon Y \to X$ be given. Fix any $\mathcal M \in \Shv(X_\et, \Lambda^a)^\oc$; we need to show that $f^* \mathcal M$ is overconvergent. Let $y = \Spa(C', C'^+) \to Y$ be a quasi-pro-étale map from a geometric point $y$. Then the composition $y \to Y \to X$ factors over a quasi-pro-étale map $x = \Spa(C, C^+) \to X$. The stalk of $f^* \mathcal M$ at $y$ is the pullback of $f^* \mathcal M$ along $y \to Y$ and hence equals the pullback of $\mathcal M$ along $y \to x \to X$. This allows us to reduce to the case $X = \Spa(C, C^+)$ and $Y = \Spa(C', C'^+)$, where the claim is easy ($\mathcal M$ is a constant sheaf).

Now assume that $f$ is qcqs and let $\mathcal N \in \Shv(Y_\et, \Lambda^a)^\oc$ be given. We want to show that $f_* \mathcal N \in \Shv(X_\et, \Lambda^a)$ is overconvergent. By taking a K-injective resolution of $\mathcal N$, using (i) and noting that stalks are $t$-exact, the claim can be checked in the hearts, i.e. we can assume that $\mathcal N$ is concentrated in degree $0$ and we only need to see that $\mathcal M := \pi_0 f_* \mathcal N$ is overconvergent. By abuse of notation we work solely in the abelian setting from now on (e.g. we write $\mathcal M(X)$ for $\pi_0 \mathcal M(X)$). The overconvergence of $\mathcal M$ can be checked on an open cover of $X$, so we can assume that $X$ and $Y$ are qcqs. Pick a quasi-pro-étale map $x = \Spa(C, C^+) = \varprojlim_i U_i \to X$, where all $U_i \to X$ are qcqs étale. Then
\begin{align*}
	\mathcal M_x = \varinjlim_i \mathcal M(U_i) = \varinjlim_i \mathcal N(U_i \cprod_X Y).
\end{align*}
Since all $U_i \cprod_X Y$ are qcqs, the above colimit is equal to $(i^* \mathcal N)(Y_x)$, where $i\colon Y_x := Y \cprod_X x \to Y$ is the natural map (cf. \cite[Proposition 14.9]{etale-cohomology-of-diamonds}). Using the first part of (iii) we can thus assume $X = x = \Spa(C, C^+)$ from now on. Let $X^\circ := \Spa(C, C^\circ)$ and $Y^\circ := Y \cprod_X X^\circ$. Then $Y^\circ \subset Y$ is an open subspace and we need to show that the natural map $\mathcal N(Y) \to \mathcal N(Y^\circ)$ is an isomorphism. In fact we claim that the natural map $\mathcal N \to j_*j^* \mathcal N$ is an isomorphism, where $j\colon Y^\circ \injto Y$ is the inclusion. This can be checked on stalks, which reduces the claim to the case $Y = \Spa(C', C'^+)$. But then the claim follows immediately from the definition of overconvergence.
\end{proof}

\begin{proposition} \label{rslt:discrete-sheaves-equiv-overconvergent-sheaves}
Let $X \in \vStacksCoeff$ be a locally spatial diamond and assume that $X \in X_\vsite^\Lambda$. Then there is a natural equivalence of $\infty$-categories
\begin{align*}
	\Dqcohri(X, \Lambda)_\omega = \Shv\cplt(X_\et, \Lambda^a)^\oc,
\end{align*}
where $\Shv\cplt$ denotes the left completion of $\Shv$. This equivalence is compatible with pullbacks along maps of locally spatial diamonds.
\end{proposition}
\begin{proof}
As in the proof of \cite[Proposition 14.10]{etale-cohomology-of-diamonds} one shows that the pullback functor along $X_\qproet \to X_\et$ induces a $t$-exact fully faithful embedding $\Shv^+(X_\et, \Lambda^a) \injto \Shv^+(X_\qproet, \Lambda^a)$ (where $\Shv^+$ denotes left-bounded objects); alternatively one can deduce this result from loc. cit. by almostification (as in \cref{rslt:general-properties-of-almost-localization-for-sheaves-on-diamond}). The functor $X_\qproet^\opp \to \infcatinf$, $U \mapsto \Shv^+(U_\qproet, \Lambda^a)$ (obtained from \cref{rslt:D+-is-a-functor}) is a sheaf of $\infty$-categories: By \cref{rslt:derived-descent-from-abelian-descent} this reduces to showing descent of the hearts, which is immediately clear because quasi-pro-étale sheaves satisfy quasi-pro-étale descent. It follows from \cite[Theorem 14.12.(i)]{etale-cohomology-of-diamonds} that the induced functor $X_\qproet^\opp \to \infcatinf$, $U \mapsto \Shv^+(U_\et, \Lambda^a)$ is a sheaf of $\infty$-categories. We claim that also the functor
\begin{align*}
	X_\qproet^\opp \to \infcatinf, \qquad U \mapsto \Shv^+(U_\et, \Lambda^a)^\oc
\end{align*}
is a sheaf of $\infty$-categories. This amounts to showing the following: Let $f\colon V \surjto U$ be a cover in $X_\qproet$ and let $\mathcal M \in \Shv^+(U_\et, \Lambda^a)$ such that $f^* \mathcal M$ is overconvergent; then $\mathcal M$ is overconvergent. But this is straightforward to see from the definitions.

We now have the two sheaves $\Shv^+((-)_\et, \Lambda^a)^\oc$ and $\Dqcohrip(-, \Lambda)_\omega$ on $X_\qproet$. In order to show that they agree, it is thus enough to do so on a basis (see \cref{rslt:sheaves-on-basis-equiv-sheaves-on-whole-site}), e.g. on the subsite of strictly totally disconnected spaces in $X_\qproet$. But then the claim follows from \cref{rslt:basic-properties-of-overconvergent-sheaves}. By taking left completions we arrive at the claimed equivalence of $\infty$-categories.
\end{proof}

\begin{example} \label{rslt:discrete-pushforward-computes-classical-cohomology}
Let $C$ be an algebraically closed non-archimedean field over $\Q_p$ and let $f\colon X \to \Spa(C, \ri_C)$ be a qcqs rigid-analytic variety. Let $\mathcal M$ be a classical sheaf of $\ri^+_X/p$-modules on $X_\et$ such that $\mathcal M^a$ is overconvergent (e.g. $\mathcal M$ is locally free). Then by \cref{rslt:discrete-sheaves-equiv-overconvergent-sheaves} we can naturally view $\mathcal M^a$ as an object of $\Dqcohri(X^\diamond, \ri^+_X/p)$. Moreover by \cref{rslt:qcqs-pushforward-preserves-overconvergent-sheaves,rslt:qcqs-pushforward-preserves-left-bounded-discrete-sheaves,rslt:general-properties-of-almost-localization-for-sheaves-on-diamond} the pushforward $f_* \mathcal M^a \in \Dqcohri(\ri_C/p)$ computes the classical almost sheaf cohomology of $\mathcal M$, i.e.
\begin{align*}
	\pi_{-n} (f_* \mathcal M^a) = H^n(X_\et, \mathcal M)^a
\end{align*}
for all $n \in \Z$.
\end{example}

\subsection{Torsors and Representations} \label{sec:ri-pi.torsors}

This subsection is an aside on representation theory in the setting of solid (almost) modules. We aim to apply this theory to better understand the category of quasicoherent $\ri^{+a}_X/\pi$-modules on torsors, which is in particular important to get a handle on finite-type maps of analytic adic spaces later on.

We will start with the basic definition of $G$-representations on analytic rings $\mathcal A$. Here $G$ will be a condensed group, i.e. a group object in the category of condensed sets. In order to define $G$-representations in a general setting we would need to work with condensed $\infty$-categories (or at least $\Cond(\Ani)$-enriched $\infty$-categories), which seems to be rather subtle. Instead we chose to pursue a much more ad-hoc approach where we derive everything we need from a static setting, where we can use very explicit constructions. In particular our theory will only work for static analytic rings $\mathcal A$, i.e. those where $\underline{\mathcal A}$ is a static ring.

In the following we fix an almost setup $(V,\mm)$. Without further ado, let us define $G$-representations over analytic rings over $(V,\mm)$:

\begin{definition}
Let $G$ be a condensed group and $A$ a static $(V,\mm)$-algebra.
\begin{defenum}
	\item A \emph{$G$-action on $A$} is a $G$-action on $A_{**}$, i.e. a map of condensed sets $\rho\colon G \cprod A_{**} \to A_{**}$ satisfying the obvious compatibilities. Given such a $G$-action, we obtain a map of $(V,\mm)$-algebras $\rho_g\colon A \to A$ for every $g \in G$.

	\item Let $\rho$ be a $G$-action on $A$. We define the $A$-algebra $A[G]_\rho$ to be the ring on $A[G]$ given by the multiplication law $[g] \times [g'] = [gg']$ and $(ga) \times [g] = [g] \times a$ for $g \in G$ and $a \in A$. More precisely, this defines a classical associative ring structure on $A_{**}(S)[G(S)]$ for every extremally disconnected set $S$, hence an associative $A_{**}$-algebra structure on $A_{**}[G]$ by sheafification and thus also an associative $A$-algebra structure on $A[G]$ via almost localization.

	\item Let $\mathcal A$ be a static analytic ring over $(V,\mm)$ (i.e. $\underline{\mathcal A}$ is static). A \emph{$G$-action on $\mathcal A$} is a $G$-action on $\underline{\mathcal A}$ such that for each $g \in G$ the induced map $\rho_g\colon \underline{\mathcal A} \to \underline{\mathcal A}$ is a map of analytic rings $\mathcal A \to \mathcal A$. Given such a $G$-action $\rho$, we let $\mathcal A[G]_\rho$ denote the associative $\mathcal A$-algebra given by the induced associative analytic ring structure on $\underline{\mathcal A}[G]_\rho$ over $\mathcal A$.
\end{defenum}
\end{definition}

\begin{example}
Let $G$ be a profinite group and let $\rho$ be the trivial $G$-action on $\Z$. Then $\rho$ is also a $G$-action on $\Z_\solid$ and we have $\Z_\solid[G]_\rho = \varprojlim_{H \subset G} \Z[G/H]$, where $H$ ranges over the compact open subgroups of $G$. Note that in general $\mathcal A[G]_\rho$ may not be static.
\end{example}

\begin{definition}
Let $G$ be a condensed group acting on a static analytic ring $\mathcal A$ over $(V,\mm)$. A \emph{semilinear $G$-action} on an $\mathcal A$-module $M \in \D(\mathcal A)$ is a $\mathcal A[G]_\rho$-module structure on $M$. We denote
\begin{align*}
	\D(\mathcal A)^G := \D(\mathcal A[G]_\rho)
\end{align*}
the $\infty$-category of $\mathcal A$-modules equipped with a semilinear $G$-action.
\end{definition}

\begin{remark}
In our applications it will always be the case that $\mathcal A[G]_\rho[S]$ is static for all extremally disconnected sets $S$. In this case $\D(\mathcal A)^G$ is the derived $\infty$-category of its heart. This heart has an explicit description: An object in the heart is a static $\mathcal A$-module $M$ together with a map $G \cprod \pi_0 M_* \to \pi_0 M_*$ of condensed sets which satisfies the obvious compatibilities to make it a semilinear $G$-action.
\end{remark}

Having introduced $G$-representations for a condensed group $G$, we now define the $G$-invariance functor. On cohomology groups this recovers the usual group cohomology.

\begin{definition} \label{def:group-cohomology}
Let $G$ be a condensed group and let $f\colon \mathcal A \to \mathcal B$ be a $G$-equivariant map of static analytic rings in $\AnRing_{(-)}$ equipped with $G$-actions $\rho$. Assume that the $G$-action on $A$ is trivial.
\begin{defenum}
	\item There is a natural morphism $\mathcal A[G]_\rho \to \mathcal A$ of $\mathcal A$-algebras given by $[g] \mapsto 1$ (more precisely this defines morphisms $\underline{\mathcal A}_{**}(S)[G(S)] \to \underline{\mathcal A}_{**}(S)$ which by sheafification produce a morphism $\underline{\mathcal A}_{**}[G]_\rho \to \underline{\mathcal A}_{**}$ and then by almost localization and $\mathcal A$-induction we get the desired morphism $\mathcal A[G]_\rho \to \mathcal A$). The forgetful functor along this morphism is denoted
	\begin{align*}
		- \tensor_{\mathcal A} (\mathcal A, G)\colon \D(\mathcal A) \to \D(\mathcal A)^G.
	\end{align*}

	\item The functor from (a) admits a right adjoint which we denote by
	\begin{align*}
		(-)^G\colon \D(\mathcal A)^G \to \D(\mathcal A), \qquad M \mapsto M^G.
	\end{align*}
	It can be computed by $M^G = \IHom_{\mathcal A[G]_\rho}(\mathcal A, M)$.

	\item The map $f$ induces a base-change functor $f^*\colon \D(\mathcal A)^G \to \D(\mathcal B)^G$. Combining this with the functor from (a) we get a functor
	\begin{align*}
		- \tensor_{\mathcal A} (\mathcal B, G)\colon \D(\mathcal A) \to \D(\mathcal B)^G.
	\end{align*}
	If $f$ is clear from context then we denote the right adjoint of this functor by
	\begin{align*}
		(-)^G\colon \D(\mathcal B)^G \to \D(\mathcal A), \qquad M \mapsto M^G
	\end{align*}
	It can be computed as $M^G = \IHom_{\mathcal A[G]_\rho}(\mathcal A, f_* M)$ for $M \in \D(\mathcal B)^G$. We will sometimes denote $\Gamma(G, M) = M^G$ and $H^k(G, M) = \pi_{-k}(M^G)$ and call it the \emph{$G$-cohomology of $M$}.
\end{defenum}
\end{definition}

In the solid setting, the above definition of group cohomology recovers the classical definition of continuous group cohomology. More precisely, we have the following explicit computation:

\begin{proposition} \label{rslt:compute-group-cohom-via-Hom-from-G}
In the situation of \cref{def:group-cohomology} assume that $G$ is locally profinite and that $\mathcal A = (A, A^+)_\solid$ and $\mathcal B = (B, B^+)_\solid$ are static discrete Huber pairs over $(V,\mm)$. Let $M \in \D_\solid(B, B^+)^G$ be static. Let $\mathcal C^\bullet$ be the cochain complex in $\D_\solid(A, A^+)^\heartsuit$ which is zero for $n < 0$ and for $n \ge 0$ has
\begin{align*}
	\mathcal C^n = \IHom_G(G^{n+1}, M),
\end{align*}
the space of $G$-equivariant maps from $G^{n+1}$ to $M$ (with the diagonal $G$-action on $G^{n+1}$), and whose $n$-th boundary map is given by
\begin{align*}
	d_n(f)(g_0, \dots, g_{n+1})) = \sum_{i=0}^n (-1)^i f(g_0, \dots, \hat g_i, \dots, g_{n+1}), \qquad \text{for $f \in \mathcal C^n$}.
\end{align*}
Then $\mathcal C^\bullet$ computes the group cohomology of $M$ in $\D_\solid(A, A^+)$, i.e.
\begin{align*}
	H^n(G, M) = H^n \mathcal C^\bullet, \qquad \text{for $n \ge 0$}.
\end{align*}
\end{proposition}
\begin{proof}
For every $n \ge 0$, note that the diagonal action of $G$ on $G^{n+1}$ induces an $(A, A^+)_\solid[G]_\rho$-module structure on $(A, A^+)_\solid[G^{n+1}]$, making it a projective $(A, A^+)_\solid[G]_\rho$-module (projectivity follows from the fact that it is isomorphic to $(A, A^+)_\solid[G]_\rho[G^n]$). There is a canonical projective resolution
\begin{align*}
	\dots \to (A, A^+)_\solid[G^{n+1}] \to \dots \to (A, A^+)_\solid[G^2] \to (A, A^+)_\solid[G] \to A \to 0
\end{align*}
in the heart of $\D_\solid(A, A^+)^G$. Now $M^G = \IHom_{(A, A^+)_\solid[G]_\rho}(A, M)$, so from the above projective resolution of $A$ we obtain the desired result.
\end{proof}

If we apply our theory of quasicoherent $\ri^+_X/\pi$-modules to a $G$-torsor we naturally end up with a certain category of $G$-representations. It follows from the ``étale nature'' of this category that it will only contain particularly nice $G$-representations, namely those whose $G$-action is continuous for the discrete topology, classically called smooth representations. The precise definitions are as follows.

\begin{definition}
Let $(A, A^+)$ be a static discrete Huber pair over $(V,\mm)$, let $S$ be a profinite set and let $M \in \D_\solid(A, A^+)$.
\begin{defenum}
	\item We denote
	\begin{align*}
		\IHom(S, M) := \IHom_A(A[S], M) \in \D_\solid(A, A^+)
	\end{align*}
	and call its objects the \emph{continuous maps from $S$ to $M$}.

	\item We denote
	\begin{align*}
		\smIHom(S, M) := \IHom(S, A) \tensor_{(A, A^+)_\solid} M \in \D_\solid(A, A^+)
	\end{align*}
	and call its objects the \emph{locally constant maps from $S$ to $M$}.

					\end{defenum}
\end{definition}

\begin{lemma} \label{rslt:computation-of-smIHom}
Let $(A, A^+)$ be a static discrete Huber pair over $(V,\mm)$, let $M \in \D_\solid(A, A^+)^\heartsuit$ and let $S$ be a profinite set. Then:
\begin{lemenum}
	\item Both $\IHom(S, M)$ and $\smIHom(S, M)$ are static. If $M$ is discrete then they are both discrete.

	\item For every profinite set $T$ we have a natural isomorphism
	\begin{align*}
		\smIHom(S, M)(T) = \IHom(S, M(T))
	\end{align*}
	of classical $A$-modules.

	\item The natural map
	\begin{align*}
		\smIHom(S, M) \injto \IHom(S, M)
	\end{align*}
	is injective. If $M$ is discrete then it is an isomorphism.
\end{lemenum}
\end{lemma}
\begin{proof}
Letting $S = \varprojlim_i S_i$ with finite $S_i$ we get
\begin{align*}
	\smIHom(S, M) = \IHom(S, A) \tensor_{(A, A^+)_\solid} M = \varinjlim_i \IHom(S_i, A) \tensor_{(A, A^+)_\solid} M = \varinjlim_i M^{S_i}.
\end{align*}
This immediately implies that $\smIHom(S, M)$ is static. Moreover we have $\IHom(S, M) = \IHom_A((A, A^+)_\solid[S], M)$, which is static because $(A, A^+)_\solid[S]$ is projective in $\D_\solid(A, A^+)^\heartsuit$. If $M$ is discrete then so is $\smIHom(S, M)$ by the above computation, and in this case we also have $\IHom(S, M) = \smIHom(S, M)$ (e.g. because $M$ is nuclear by \cref{rslt:almost-solid-discrete-equiv-nuclear}). This finishes the proof of (i).

Part (ii) follows immediately from the above computation of $\smIHom(S, M)$. Part (iii) also follows from that computation by using that the maps $\IHom(S_i, M) \injto \IHom(S, M)$ are clearly injective.
\end{proof}

\begin{definition} \label{def:smooth-group-rep}
Let $G$ be a locally profinite group acting on a static discrete Huber pair $(A, A^+)$ over $(V,\mm)$.
\begin{defenum}
	\item For every static $M \in \D_\solid(A, A^+)^\heartsuit$ we define $\smIHom(G, M) \in \D_\solid(A, A^+)$ via
	\begin{align*}
		\smIHom(G, M)(S) := \IHom(G, M(S))
	\end{align*}
	for all profinite sets $S$. Given any open compact subgroup $H \subset G$ we note by \cref{rslt:computation-of-smIHom} that $\smIHom(G, M) = \prod_{G/H} \smIHom(H, M)$, which is indeed a (static) $(A, A^+)_\solid$-module. We also deduce that the natural map $\smIHom(G, M) \injto \IHom(G, M)$ is injective.

	\item A static $G$-representation $M \in \D_\solid(A, A^+)^G$ is called \emph{smooth} if the induced map $M \to \IHom(G, M)$ factors over $\smIHom(G, M)$.
\end{defenum}
\end{definition}

\begin{lemma} \label{rslt:smooth-reps-are-Grothendieck-abelian}
Let $G$ be a locally profinite group acting on a static discrete Huber pair $(A, A^+)$ over $(V,\mm)$. For every strong limit cardinal $\kappa$, the full subcategory $\mathcal C_\kappa \subset \D_\solid(A, A^+)_\kappa^G$ spanned by the smooth static $G$-representations is a Grothendieck abelian category.
\end{lemma}
\begin{proof}
The only non-trivial part is the existence of generators. They can be given by $(A, A^+)_\solid[S][G/H]$ for compact open subgroups $H \subset G$ and $\kappa$-small profinite sets $S$, equipped with the obvious $G$-action. That these elements indeed form a family of generators follows from the fact that every element in $M(S)$ for any smooth $G$-representation $M$ is fixed by some $H$.
\end{proof}

\begin{definition}
Let $G$ be a locally profinite group acting on a static discrete Huber pair $(A, A^+)$ over $(V,\mm)$. We denote
\begin{align*}
	\Drepsldsm G{(A, A^+)} := \varinjlim_\kappa \D(\mathcal C_\kappa)
\end{align*}
with $\mathcal C_\kappa$ as in \cref{rslt:smooth-reps-are-Grothendieck-abelian}. The elements of $\Drepsldsm G{(A, A^+)}$ are called the \emph{smooth $G$-representations on $(A, A^+)_\solid$-modules}.
\end{definition}

\begin{remark} \label{rslt:any-rep-on-discrete-module-is-smooth}
By \cref{rslt:computation-of-smIHom}, if $M$ is a discrete $(A, A^+)_\solid$-module then $\smIHom(G, M) = \IHom(G, M)$. In particular every $G$-action on $M$ is automatically smooth. It follows that the $\infty$-category of smooth $G$-representations on discrete $A$-modules embeds fully faithfully into $\Drepsldsm G{(A, A^+)}$.
\end{remark}

We now study the $\infty$-category of smooth $G$-representations. Our first result is the following characterization of smooth $G$-representations as ``families of discrete $G$-representations'':

\begin{lemma} \label{rslt:smooth-G-rep-is-family-of-discrete}
Let $G$ be a locally profinite group acting on a static discrete Huber pair $(A, A^+)$ over $(V,\mm)$ and let $M$ be a static $(A, A^+)_\solid$-module. A smooth semilinear $G$-representation on $M$ can be equivalently described by a functor $\rho$ from $*_\proet$ to the category of continuous semilinear $G$-actions on discrete $A$-modules such that for every $S \in *_\proet$, the underlying $A$-module of $\rho(S)$ is $M(S)$ and $\rho(S)$ is $A(S)$-semilinear.
\end{lemma}
\begin{proof}
By definition, a smooth $G$-action on $M$ is given by an $A$-linear map $\rho'\colon M \to \smIHom(G, M)$ satisfying certain axioms. Such a map is specified by a functor which associates to every $S \in *_\proet$ a map of discrete $A$-modules $\rho'(S)\colon M(S) \to \smIHom(G, M)(S) = \IHom(G, M(S))$. The group action and $A$-semilinearity axioms for $\rho'$ translate to group action and $A$-similinearity axioms for $\rho'(S)$ (and vice versa).
\end{proof}

Similar to the setting of general representations, we can introduce group cohomology for smooth representations. This \emph{smooth} group cohomology is in general not equal to the group cohomology defined in \cref{def:group-cohomology} via the functor $\Drepsldsm G{(A, A^+)} \to \D_\solid(A, A^+)^G$. In particular, the latter functor is \emph{not} fully faithful in general. We will provide an explicit computation of smooth group cohomology, which sheds more light on the difference between continuous group cohomology and smooth group cohomology, see \cref{rmk:difference-of-smooth-group-cohom-and-group-cohom} below.

\begin{definition} \label{def:group-cohom-of-smooth-rep}
Let $G$ be a locally profinite group acting on a static discrete Huber pair $(B, B^+)$ over $(V,\mm)$. Let $f\colon (A, A^+) \to (B, B^+)$ be a $G$-equivariant map of static discrete Huber pairs over $(V,\mm)$, where $G$ acts trivially on $(A, A^+)$. We define
\begin{align*}
	\Gamma_\sm(G, -)\colon \Drepsldsm G{(B, B^+)} \to \D_\solid(A, A^+)
\end{align*}
as the right derived functor of the functor $M \mapsto \IHom(B, M)$, where $B$ denotes the $G$-representation on the $B$-module $B$ induced by the given $G$-action on the ring $B$ and the $\IHom$ denotes the $\D_\solid(A, A^+)^\heartsuit$-enriched $\Hom$ in $\Drepsldsm G{(B, B^+)}^\heartsuit$.
\end{definition}

\begin{lemma} \label{rslt:smooth-group-cohom-is-computed-on-each-S-separately}
In the situation of \cref{def:group-cohom-of-smooth-rep} let $M \in \Drepsldsm G{(B, B^+)}$. Then for every profinite set $S$ we have a natural isomorphism
\begin{align*}
	\Gamma_\sm(G, M)(S) = \Gamma(G, M(S))
\end{align*}
in $\D(A(S))$, i.e. smooth group cohomology can be computed on each $S$ separately.
\end{lemma}

\begin{remark}
By \cref{rslt:smooth-group-cohom-is-computed-on-each-S-separately} the computation of the smooth group cohomology $\Gamma_\sm(G, M)$ is reduced to computing the continuous group cohomologies $\Gamma(G, M(S))$ of the discrete $G$-representations $M(S)$.
\end{remark}

\begin{proof}[Proof of \cref{rslt:smooth-group-cohom-is-computed-on-each-S-separately}]
Let $r\colon \D_\solid(B, B^+)^\heartsuit \to \Drepsldsm G{(B, B^+)}^\heartsuit$ be the right adjoint of the forgetful functor. For every $M \in \Drepsldsm G{(B, B^+)}^\heartsuit$ and every $n \ge 0$ we denote $r^n(M) \in \Drepsldsm G{(B, B^+)}^\heartsuit$ the application of $n$ times $r$ to $M$ (with implicit forgetful functors in between). By adjunction there is a natural map $M \to r(M)$, and in general there are $n+1$ natural maps $r^n(M) \to r^{n+1}(M)$ for $n \ge 0$. Taking the alternating sum over these maps gives a complex
\begin{align}
	0 \to M \to r(M) \to r^2(M) \to \dots \label{eq:canonical-resolution-of-smooth-rep}
\end{align}
in $\Drepsldsm G{(B, B^+)}$. We claim that this complex is exact. One first checks that
\begin{align*}
	r(M) = \smIHom(G, M)
\end{align*}
with $G$ acting by restricting the natural $G$-action on $\IHom(G, M)$ (use \cref{rslt:smooth-G-rep-is-family-of-discrete} to reduce everything to the case of discrete modules). It follows from the explicit description of $r(M)$ that the map $M \to r(M)$ splits $B$-linearly (but not $G$-equivariantly). From this one easily shows that the complex \cref{eq:canonical-resolution-of-smooth-rep} is split exact in $\D_\solid(B, B^+)$ and in particular exact in $\Drepsldsm G{(B, B^+)}$. We have thus constructed a functorial resolution of $M$.

We now claim that (for $M$ still lying in the heart)
\begin{align*}
	\Gamma_\sm(G, M) = [H^0_\sm(G, r(M)) \to H^0_\sm(G, r^2(M)) \to \dots],
\end{align*}
for which we need to check that $H^k(\Gamma_\sm(G, r^n(M))) = 0$ for all $k, n > 0$. Clearly we can assume $n = 1$, and we claim that in fact $\Gamma_\sm(G, r(M)) = M$. To see this, note first that $r$ is exact (by the above explicit description) and hence extends to a $t$-exact functor $r\colon \D_\solid(B, B^+) \to \Drepsldsm G{(B, B^+)}$, which is still right adjoint to the forgetful functor (now on the derived level). Let $(-)^\triv\colon \D_\solid(A, A^+) \to \Drepsldsm G{(B, B^+)}$ be a left adjoint of $\Gamma_\sm(G, -)$ (which exists by the adjoint functor theorem). Then $\Gamma_\sm(G, r(-))$ is right adjoint to the composition of $(-)^\triv$ with the forgetful functor, so it is enough to show that the latter composition is equal to $- \tensor_{(A, A^+)_\solid} (B, B^+)_\solid$. This can be checked on compact projective generators, in which case the claim reduces to the abelian level, where it is easily verified.

We have now managed to explicitly compute $\Gamma_\sm(G, M)$ for $M \in \Drepsldsm G{(B, B^+)}^\heartsuit$. Next we claim that for any profinite set $S$ we have
\begin{align*}
	\Gamma(G, M(S)) = [H^0(G, r(M)(S)) \to H^0(G, r^2(M)(S)) \to \dots],
\end{align*}
for which we need to show that $H^k(G, r^n(M)(S)) = 0$ for all $k, n > 0$. Note that $r(M)(S) = r(M(S))$ (where we view $M(S)$ as a discrete $B(S)$-module), so we can reduce to the case $n = 1$. Denoting $N := M(S)$ (equipped with the inherited $G$-action), the claim now is
\begin{align*}
	H^k(G, \IHom(G, N)) = 0 \qquad \text{for $k > 0$.}
\end{align*}
But this follows easily using Hom-tensor adjunction:
\begin{align*}
	\Gamma(G, \IHom(G, N)) &= \Hom_{(A, A^+)_\solid[G]}(A, \IHom_A((A, A^+)_\solid[G], N))\\
	&= \Hom_{(A, A^+)_\solid[G]}((A, A^+)_\solid[G] \tensor_{(A, A^+)_\solid} A, N)\\
	&= \Hom_{(A, A^+)_\solid[G]}((A, A^+)_\solid[G], \IHom_A(A, N))\\
	&= N.
\end{align*}
(Alternatively use an explicit computation via \cref{rslt:compute-group-cohom-via-Hom-from-G}.) All in all we see that the claimed $\Gamma_\sm(G, M)(S) = \Gamma(G, M(S))$ for $M \in \Drepsldsm G{(B, B^+)}$ reduces to
\begin{align*}
	H^0_\sm(G, r^n(M))(S) = H^0(G, r^n(M)(S))
\end{align*}
The right-hand side are just the $G$-invariants (in the classical sense) of the $G$-action on $r^n(M)(S)$. By \cref{rslt:smooth-G-rep-is-family-of-discrete} the same follows for the left-hand side. This finishes the proof in the case that $M$ is static.

Now let $M \in \Drepsldsm G{(B, B^+)}$ be general. We can find a K-injective complex $I^\bullet$ which is isomorphic to $M$. Then $\Gamma_\sm(G, M)(S)$ is computed by the complex $H^0_\sm(G, I^\bullet)(S)$. On the other hand by what we have shown above each $I(S)$ is $\Gamma$-acyclic. Hence $\Gamma(G, M(S))$ is computed by the complex $H^0(G, I^\bullet(S))$. As before, it follows easily from \cref{rslt:smooth-G-rep-is-family-of-discrete} that both complexes are equal.
\end{proof}

\begin{corollary} \label{rslt:compute-smooth-group-cohom-via-smIHoms}
In the situation of \cref{def:group-cohom-of-smooth-rep} let $M \in \Drepsldsm G{(B, B^+)}^\heartsuit$. Let $\mathcal C_\sm^\bullet$ be the cochain complex in $\D_\solid(A, A^+)$ which is zero for $n < 0$ and for $n \ge 0$ has
\begin{align*}
	\mathcal C_\sm^n = \IHom^\sm_G(G^{n+1}, M),
\end{align*}
the space of $G$-equivariant locally constant maps from $G^{n+1}$ to $M$ (with the diagonal $G$-action on $G^{n+1}$), and whose $n$-th boundary map is given by
\begin{align*}
	d_n(f)(g_0, \dots, g_{n+1})) = \sum_{i=0}^n (-1)^i f(g_0, \dots, \hat g_i, \dots, g_{n+1}), \qquad \text{for $f \in \mathcal C^n_\sm$}.
\end{align*}
Then $\mathcal C_\sm^\bullet$ computes the smooth group cohomology of $M$, i.e.
\begin{align*}
	H_\sm^n(G, M) = H^n \mathcal C_\sm^\bullet, \qquad \text{for $n \ge 0$}.
\end{align*}
\end{corollary}
\begin{proof}
Evaluate both sides on all profinite sets $S$ and use \cref{rslt:smooth-group-cohom-is-computed-on-each-S-separately,rslt:computation-of-smIHom,rslt:compute-group-cohom-via-Hom-from-G}.
\end{proof}

\begin{remark} \label{rmk:difference-of-smooth-group-cohom-and-group-cohom}
By comparing \cref{rslt:compute-group-cohom-via-Hom-from-G} and \cref{rslt:compute-smooth-group-cohom-via-smIHoms} we can describe the difference between (continuous) $G$-representations and smooth $G$-representations (and their cohomologies) as follows. A $G$-representation on an $(A, A^+)_\solid$-module $M$ is just a continuous representation on $M$. On the other hand, a smooth $G$-representation on $M$ is a $G$-action on $M$ which is even continuous for the \emph{discrete} topology on $M$. Now $\Gamma(G, M)$ computes the classical continuous group cohomology of $M$ while $\Gamma_\sm(G, M)$ computes the group cohomology of the \emph{discrete} $G$-representation on $M$, but afterwards equips each $H^k_\sm(G, M)$ with a topology coming from the topology on $M$. From this picture it is clear that if $M$ is discrete then a $G$-representation on $M$ is the same as a smooth $G$-representation on $M$ (and $\Gamma(G, M) = \Gamma_\sm(G, M)$), but for general $M$ the two notions differ.
\end{remark}

\begin{lemma} \label{rslt:smooth-rep-underlying-module-is-colim-of-H-invariants}
In the situation of \cref{def:group-cohom-of-smooth-rep} let $M \in \Drepsldsm G{(B, B^+)}$. Then there is a natural isomorphism
\begin{align*}
	M = \varinjlim_{H \subset G} \Gamma_\sm(H, M)
\end{align*}
in $\D_\solid(A, A^+)$. Here $H$ ranges over all compact open subgroups of $G$ and $\Gamma_\sm(H, M)$ is computed by viewing $M$ as a complex of smooth $H$-representations via the forgetful functor.
\end{lemma}
\begin{proof}
One checks that there is a compatible system of canonical maps $\Gamma(H, M) \to M$; we need to show that the resulting map $\varinjlim_H \Gamma(H, M) \to M$ is an isomorphism. By \cite[Theorem 079P]{stacks-project} there is a K-injective complex $I^\bullet$ in $\Drepsldsm G{(B, B^+)}$ and an isomorphism $M \isom I^\bullet$. Then $\Gamma_\sm(G, M)$ is computed by the complex $H^0_\sm(G, I^\bullet)$ (see \cite[Lemma 070L]{stacks-project}). For every compact open subgroup $H \subset G$, $I^\bullet$ is also K-injective as a complex of smooth $H$-representations, hence $\Gamma_\sm(H, M)$ is computed by the complex $H^0_\sm(H, I^\bullet)$. This reduces the claim to the case that $M = I$ is an (injective) object in $\Drepsldsm G{(B, B^+)}$, in which case we want to show that
\begin{align*}
	I = \varinjlim_{H \subset G} I^H.
\end{align*}
This can be checked on each profinite set $S$ separately, which by \cref{rslt:smooth-group-cohom-is-computed-on-each-S-separately} reduces the claim to the case that $I$ is a discrete $B$-module with continuous $G$-action. But now the claim is clear because every $x \in I$ has an open and closed stabilizer $H \subset G$ (this works also in the almost setting, e.g. by writing $I = (I_!)^a$).
\end{proof}

In general the $\infty$-category $\Drepsldsm G{(A, A^+)}$ is not left complete. However, if we put some finiteness conditions on $G$ and $A$ then this left-completeness follows. The relevant definition is as follows. Note that our definition of the $\ell$-cohomological dimension is a priori different from the classical definition, since the classical one only takes $G$-representations on \emph{discrete} $\Fld_\ell$-modules into account. However, in practice both definitions agree, as \cref{rslt:cd-l-coincides-with-classical-version} shows.

\begin{definition}
Let $G$ be a condensed group acting on a static analytic ring $\mathcal A$. We say that $G$ has \emph{cohomological dimension $\le n$ over $\mathcal A$} if for all static $M \in \D(\mathcal A)^G$ and all $i > n$ we have $H^i(G, M) = 0$ in $\D(\Z)$. We denote $\cd_{\mathcal A} G$ the cohomological dimension of $G$ over $\mathcal A$.

Given a prime $\ell$, we abbreviate $\cd_\ell \, G := \cd_{\Fld_{\ell\solid}} G$ and call it the \emph{$\ell$-cohomological dimension} of $G$.
\end{definition}

\begin{lemma} \label{rslt:group-cd-simple-bounds}
Let $G$ be a locally profinite group and $\mathcal A$ a static analytic ring.
\begin{lemenum}
	\item Suppose $G$ also acts on a static analytic ring $\mathcal B$ and there is a $G$-equivariant map $\mathcal A \to \mathcal B$ of analytic rings. Then $\cd_{\mathcal B} G \le \cd_{\mathcal A} G$.
	\item For every open subgroup $H \subset G$ we have $\cd_{\mathcal A} H \le \cd_{\mathcal A} G$.
\end{lemenum}
\end{lemma}
\begin{proof}
Part (i) follows from the fact that the group cohomology functor $\D(\mathcal B)^G \to \D(\Z)$ factors over the forgetful functor $\D(\mathcal B)^G \to \D(\mathcal A)^G$. To prove (ii), let $\Ind^G_H\colon \D(\mathcal A)^H \to \D(\mathcal A)^G$ be the right adjoint of the forgetful functor. Since the forgetful functor preserves projective objects in the hearts, $\Ind^G_H$ is $t$-exact. It follows from the definitions that for any $M \in \D(\mathcal A)^H$ we have $\Gamma(H, M) = \Gamma(G, \Ind^G_H M)$, which implies (ii).
\end{proof}

\begin{proposition} \label{rslt:cd-l-coincides-with-classical-version}
Let $G$ be a profinite group and $\ell$ a prime. Assume that for every static $G$-representation $M$ on a finite $\Fld_\ell$-vector space and all $i \ge 0$, $H^i(G, M)$ is finite. Then $\cd_\ell \, G \in [0, \infty]$ is the smallest number such that for all $i > \cd_\ell \, G$ and all static finite discrete $M \in \D(\Fld_\ell)_\omega^G$ we have $H^i(G, M) = 0$.
\end{proposition}
\begin{proof}
The following argument is based on ideas of a workshop on condensed mathematics held in La Tourette in Fall 2021.

Suppose we have some integer $n \ge 0$ such that for all static finite discrete $M \in \D(\Fld_\ell)_\omega^G$ and all $i > n$ we have $H^i(G, M) = 0$. Let $\mathcal C \subset \D(\Fld_{\ell\solid})^G$ be the full subcategory of static objects $M$ with $H^i(G, M) = 0$ for $i > n$. Let $\mathcal C_0 \subset \D(\Fld_{\ell\solid})^G$ be the full subcategory of static finite discrete objects, so that $\mathcal C_0 \subset \mathcal C$ by assumption. We claim that the canonical functor
\begin{align*}
	\Ind(\Pro(\mathcal C_0)) \isoto (\D(\Fld_{\ell\solid})^G)^\heartsuit
\end{align*}
is essentially surjective (it is even an equivalence). This follows in the same way as in the proof of \cite[Theorem 2.9]{scholze-analytic-spaces} using that the compact projective generators $\Fld_{\ell\solid}[G][S]$ of the heart of $\D(\Fld_{\ell\solid})^G$ lie in $\Pro(\mathcal C_0)$. From the explicit computation in \cref{rslt:compute-group-cohom-via-Hom-from-G} it follows easily that $H^i(G, -)$ commutes with filtered colimits, hence $\mathcal C$ is stable under filtered colimits. It is therefore enough to show that $\mathcal C$ contains the image of $\Pro(\mathcal C_0)$, i.e. that it contains every cofiltered limit of objects in $\mathcal C_0$. Suppose $M = \varprojlim_i M_i$ is of this form. Then in fact $M = \varprojlim_i M_i$ holds true in the derived category (i.e. $R^i\lim_i M_i = 0$ for $i > 0$): Since the forgetful functor $\Drepsld G{\Fld_\ell} \to \D_\solid(\Fld_\ell)$ is conservative and preserves limits, we can check the claim on underlying $\Fld_\ell$-vector spaces. Then in the derived category we have
\begin{align*}
	\varprojlim_i M_i = \varprojlim_i \IHom_{\Fld_\ell}(\IHom_{\Fld_\ell}(M_i, \Fld_\ell), \Fld_\ell) = \IHom_{\Fld_\ell}(\varinjlim_i \IHom_{\Fld_\ell}(M_i, \Fld_\ell)).
\end{align*}
But $N := \varinjlim_i \IHom_{\Fld_\ell}(M_i, \Fld_\ell)$ is discrete, hence $\IHom_{\Fld_\ell}(N, \Fld_\ell)$ can be computed in the derived category $\D(\Fld_\ell)$ of discrete $\Fld_\ell$-vectors spaces. But then $\Fld_\ell$ is injective and hence the $\IHom$ is concentrated in degree $0$, as claimed. It follows that $\Gamma(G, M) = \varprojlim_i \Gamma(G, M_i)$. By assumption each $H^k(G, M_i)$ is finite, from which we deduce $R^n \varprojlim_i H^k(G, M_i) = 0$ for $n > 0$ (by the same argument as above). Since also $H^k(G, M_i) = 0$ for $k > n$ it follows that $H^k(G, M) = 0$ for $k > n$, as desired.
\end{proof}

\begin{lemma} \label{rslt:group-cd-bounds-smooth-cd}
Let $G$ be a locally profinite group acting on a static discrete Huber pair $(A, A^+)$. Then for every $i > \cd_{(A, A^+)_\solid} G$ and every static $M \in \Drepsldsm G{(A, A^+)}$ we have $H^i_\sm(G, M) = 0$.
\end{lemma}
\begin{proof}
This follows directly from \cref{rslt:smooth-group-cohom-is-computed-on-each-S-separately} using that we have $\cd_{(A(S), A^+)_\solid} G \le \cd_{(A, A^+)_\solid} G$ by \cref{rslt:group-cd-simple-bounds}.
\end{proof}

We now come to the promised left-completeness result of the $\infty$-category of smooth representations. The relevant notion is the following.

\begin{definition} \label{def:locally-finite-cohom-dim}
Let $G$ be a locally profinite group acting on a static analytic ring $\mathcal A$. We say that $G$ has \emph{locally finite cohomological dimension over $\mathcal A$} if there is some open subgroup $H \subset G$ such that $\cd_{\mathcal A} H < \infty$.

Given a prime $\ell$, we say that $G$ has \emph{locally finite $\ell$-cohomological dimension} if $G$ has locally finite cohomological dimension over $\Fld_{\ell\solid}$.
\end{definition}

\begin{lemma} \label{rslt:D-smooth-rep-is-left-complete-if-G-has-fin-cohom-dim}
Let $G$ be a locally profinite group acting on a static discrete Huber pair $(A, A^+)$ over $(V,\mm)$. If $G$ has locally finite cohomological dimension over $(A_{**\omega}, A^+)$ then $\Drepsldsm G{(A, A^+)}$ is left complete.
\end{lemma}
\begin{proof}
First assume that $(V,\mm) = (\Z, \Z)$. Given $M \in \Drepsldsm G{(A, A^+)}$, we need to show that the natural map $M \to \varprojlim_n \tau_{\le n} M$ is an isomorphism. The forgetful functor $\Drepsldsm G{(A, A^+)} \to \D(\Z)$ is conservative, so it is enough to check the desired isomorphism after applying this forgetful functor. Then by \cref{rslt:smooth-rep-underlying-module-is-colim-of-H-invariants} we are reduced to showing that
\begin{align*}
	\Gamma_\sm(H, M) \to \varprojlim_n \Gamma_\sm(H, \tau_{\le n} M)
\end{align*}
is an isomorphism for a basis of open subgroups $H \subset G$. By \cref{rslt:group-cd-simple-bounds,rslt:group-cd-bounds-smooth-cd} we can find some $d \ge 0$ and a basis of $H$'s such that each functor $\Gamma_\sm(H,-)$ has cohomological dimension $\le d$. This implies that for all $k \in \Z$ we have
\begin{align*}
	\tau_{\le k} \varprojlim_n \Gamma_\sm(H, \tau_{\le n} M) = \tau_{\le k} \Gamma_\sm(H, \tau_{\le k+d} M) = \tau_{\le k} \Gamma_\sm(H, M).
\end{align*}
Using that $\D(\Z)$ is left complete, we arrive at the desired result.

Now let $(V,\mm)$ be general. Then there is a $t$-exact almost localization functor
\begin{align*}
	(-)^a\colon \Drepsldsm G{(A_{**\omega}, A^+)} \to \Drepsldsm G{(A, A^+)}
\end{align*}
which has a fully faithful left adjoint $(-)_!$ (these functors act as $(-)^a$ and $(-)_!$ on the underlying modules). In particular $(-)^a$ is essentially surjective and preserves small limits, including Postnikov limits. Thus the left-completeness of $\Drepsldsm G{(A, A^+)}$ follows from the left-completeness of $\Drepsldsm G{(A_{**\omega}, A^+)}$ shown above.
\end{proof}

We are finally in the position to apply the above representation theory to the theory of $\ri^+_X/\pi$-modules. In fact, if $X = [\Spa(C, \ri_C)/G]$ is the classifying stack of some locally profinite group $G$ over an algebraically closed field $C$ then one would expect $\DqcohriX X$ to be some kind of representation category. This is true:

\begin{lemma} \label{rslt:compute-Dqcohri-for-classifying-stack}
Let $X = \Spa(A, A^+) \in \AffPerfd_\pi$ be totally disconnected, let $G$ be a locally profinite group acting on $X$ and let $X/G$ be the (stack) quotient.
\begin{lemenum}
	\item There is a natural equivalence of $\infty$-categories
	\begin{align*}
		\Dqcohrip(\ri^+_{X/G}/\pi) = \Drepsldasmp G{A^+/\pi}.
	\end{align*}

	\item Assume that $G$ has locally finite $p$-cohomological dimension. Then the equivalence in (i) extends to an equivalence
	\begin{align*}
		\DqcohriX{X/G} = \Drepsldasm G{A^+/\pi}.
	\end{align*}

	\item Let $f\colon X \to Y = \Spa(B, B^+)$ be a $G$-equivariant map in $\AffPerfd_\pi$, where $G$ acts trivially on $Y$ and $Y$ is totally disconnected. Then under the equivalences in (i) and (ii) the pushforward $f_*$ computes smooth group cohomology.
\end{lemenum}
\end{lemma}

\begin{remark}
We prove a more general version of \cref{rslt:compute-Dqcohri-for-classifying-stack} below (see \cref{rslt:compute-Dqcohri-on-stack-quotient-of-p-bounded-space}).
\end{remark}

\begin{proof}
Let $Y_\bullet \to X/G$ be the Čech nerve of the projection $X \to X/G$. Then $Y_n \isom X \cprod G^n$, so all $Y_n$ are disjoint unions of totally disconnected spaces and the pullback on $\Dqcohri$ along any map $Y_n \to Y_m$ with $n \ge m$ is $t$-exact. It follows that the assumptions of \cref{rslt:derived-descent-from-abelian-descent} are satisfied and hence that $\Dqcohrip(\ri^+_{X/G}/\pi) = \D^+(\mathcal A)$ for a certain abelian category $\mathcal A$. Using \cref{rslt:descent-computation-of-abelian-limit} we can compute $\mathcal A$ to be the category of pairs $(M, \alpha)$, where $M \in \D^a_\solid(A^+/\pi)^\heartsuit$ and $\alpha$ is a $\IHom(G, A^+/\pi)^a$-linear automorphism of $\smIHom(G, M)$ (built out of a $G/H$-family of $\IHom(H, A^+/\pi)^a$-linear automorphisms of $\smIHom(H, M)$ for any open compact subgroup $H \subset G$) satisfying a certain cocycle condition. By adjunction, $\alpha$ can equivalently be described as an $A^{+a}/\pi$-linear map $M \to \smIHom(G, M)$ and one easily checks that the cocycle condition on $\alpha$ is precisely the group representation cocycle condition. Thus $\mathcal A$ is the category of smooth $G$-representations on $(A^+/\pi)^a_\solid$-modules, proving (i).

To prove (ii) let us assume that $G$ has locally finite $p$-cohomological dimension. By \cref{rslt:derived-descent-from-abelian-descent} it is enough to show that then $\Drepsldasm G{A^+/\pi}$ is left complete. By \cref{rslt:D-smooth-rep-is-left-complete-if-G-has-fin-cohom-dim} it is enough to show that $G$ has locally finite cohomological dimension over $A^+/\pi$. If $\pi \divides p$ then this follows immediately from the fact that $G$ has locally finite $p$-cohomological dimension. For general $\pi$ we have $\pi \divides \varpi^n$ for some $\varpi \divides p$ and some $n > 0$; we can then argue by induction on $n$.

It remains to prove (iii), so let $f$ be given. Then under the equivalences in (i) and (ii) we have $A^{+a}/\pi = f^* (B^{+a}/\pi)$ (where the left-hand side denotes the $G$-representation on $A^{+a}/\pi$ given by the $G$-action on $X$). Thus it follows easily from the adjunction of $f^*$ and $f_*$ that $f_*$ corresponds to $\Gamma_\sm(G, -)$.
\end{proof}

\begin{corollary} \label{rslt:classifying-stack-over-tot-disc-pushforward-bounded}
Let $X \in \AffPerfd_\pi$ be totally disconnected, let $G$ be a locally profinite group, let $X/G$ be the stack quotient (with respect to the trivial action of $G$) and let $f\colon X/G \to X$ denote the projection. If $G$ has finite $p$-cohomological dimension then $f_*\colon \DqcohriX{X/G} \to \DqcohriX X$ preserves $\Dqcohrim$.
\end{corollary}
\begin{proof}
This follows immediately from \cref{rslt:compute-Dqcohri-for-classifying-stack} and \cref{rslt:group-cd-bounds-smooth-cd} (using the argument in the proof of \cref{rslt:compute-Dqcohri-for-classifying-stack} to show that $G$ has finite cohomological dimension over $\ri^+_X(X)/\pi$).
\end{proof}

\subsection{Cohomologically Bounded Morphisms} \label{sec:ri-pi.p-bounded}

In preparation of the definition of the lower-shriek functor, the goal of the present section is to improve the base-change result \cref{rslt:base-change-for-bounded-Dqcohri} by extending it from $\Dqcohrip$ to $\Dqcohri$ under a stronger assumption on $f\colon Y \to X$, which we call $p$-boundedness. We will then provide a large class of $p$-bounded maps and show that they enjoy excellent properties.

To motivate the definition of $p$-bounded maps it helps to look at the $\ell$-adic case. There, the requirement for $f$ to satisfy unbounded base-change is that $f$ has finite $\ell$-cohomological dimension, i.e. that the pushforward $f_*$ has finite cohomological dimension on $\ell$-torsion sheaves. We therefore need to introduce a notion of finite $p$-cohomological dimension for $f$, which should roughly state that the pushforward $f_*\colon \Dqcohri(Y, \Lambda) \to \Dqcohri(X, \Lambda)$ has finite cohomological dimension. However, as there is in general no $t$-structure on $\Dqcohri(X, \Lambda)$, we cannot really define the cohomological dimension of $f_*$. The non-existence of the $t$-structure stems from the fact that for a map $\Spa(B, B^+) \to \Spa(A, A^+)$ of totally disconnected perfectoid spaces with pseudouniformizer $\pi$ the map $(A^+/\pi)^a_\solid \to (B^+/\pi)^a_\solid$ of analytic almost rings is not flat. On the other hand, the map $A^+/\pi \to B^+/\pi$ of classical rings \emph{is} flat, as is consequently the map $(A^+/\pi, \Z)^a_\solid \to (B^+/\pi, \Z)^a_\solid$ of analytic almost rings. This suggests to glue $\Dqcohri(\Lambda(X), \Z)$ in order to get a stable $\infty$-category $\Dqcohri(X, (\Lambda, \Z))$ closely related to $\Dqcohri(X, \Lambda)$, but with a $t$-structure. This glueing is possible by the following result.

\begin{lemma}
There is a unique hypercomplete v-sheaf
\begin{align*}
	(\vStacksCoeff)^\opp \to \infcatinf, \qquad (X, \Lambda) \mapsto \Dqcohri(X, (\Lambda, \Z))
\end{align*}
of $\infty$-categories on $\vStacksCoeff$ which on all those $(X, \Lambda) \in \AffPerfCoeff$ with $X$ being weakly of perfectly finite type over some totally disconnected space takes the form
\begin{align*}
	\Dqcohri(X, (\Lambda, \Z)) = \Dqcohri(\Lambda(X), \Z).
\end{align*}
It has the following properties:
\begin{lemenum}
	\item There is a natural left-complete $t$-structure on $\Dqcohri(X, (\Lambda, \Z))$ and all pullback functors are $t$-exact.

	\item \label{rslt:bounded-DqcohriZ-is-sheaf} The full subcategories
	\begin{align*}
		\Dqcohrip(X, (\Lambda, \Z)), \Dqcohrim(X, (\Lambda, \Z)), \Dqcohrib(X, (\Lambda, \Z)) \subset \Dqcohri(X, (\Lambda, \Z))
	\end{align*}
	of (left/right) bounded objects form themselves hypercomplete v-sheaves of $\infty$-categories.

	\item \label{rslt:DqcohriZ-invariant-under-compactification} If $X \in \vStacksCoeff$ is separated and $\Lambda$ extends to $\overline X$ then there is a natural equivalence
	\begin{align*}
		\Dqcohri(X, (\Lambda, \Z)) = \Dqcohri(\overline X, (\Lambda, \Z)).
	\end{align*}
\end{lemenum}
\end{lemma}
\begin{proof}
One checks easily using \cref{rslt:discrete-desc-equiv-abs-compactification-desc,rslt:descendable-implies-discretely-descendable} that all the descent results used in \cref{rslt:def-of-qcoh-Lambda-modules} extend from $\Dqcohri(\Lambda(X))$ to $\Dqcohri(\Lambda(X), \Z)$. This proves the first part of the claim. Part (i) follows from the fact that for any map $f\colon Y \to X$ of affinoid perfectoid spaces in $\vStacksCoeff$ with $X$ totally disconnected, the map $(\Lambda(X), \Z)^a_\solid \to (\Lambda(Y), \Z)^a_\solid$ is flat -- indeed, base-change along this map amounts to $- \tensor_{\Lambda(X)^a} \Lambda(Y)^a$ which is exact by \cite[Proposition 7.23]{etale-cohomology-of-diamonds}. Part (ii) follows by the same argument as in \cref{rslt:bounded-Dqcohri-is-sheaf}, noting that \cref{rslt:v-cover-of-tot-disc-implies-equiv-of-boundedness} adapts immediately to the new setting.

It remains to prove part (iii), which reduces immediately to the case that $X$ is a totally disconnected perfectoid space. Then $(\Lambda(X), \Z)^a_\solid = (\Lambda(\overline X), \Z)^a_\solid$, so we only need to show $\Dqcohri(\overline X, (\Lambda, \Z)) = \Dqcohri(\Lambda(\overline X), \Z)$. This reduces to showing that for some pro-étale cover $Y \surjto \overline X$ by a strictly totally disconnected space $Y$, the map $(\Lambda(\overline X), \Z)^a_\solid \to (\Lambda(Y), \Z)^a_\solid$ is weakly fs-descendable. Let $Y_0 := Y \cprod_{\overline X} X$. Then $(\Lambda(Y_0), \Z)^a_\solid = (\Lambda(Y), \Z)^a_\solid$ (e.g. by \cref{rslt:base-change-for-affinoid-perfectoid}), so it is enough to show that $(\Lambda(X), \Z)^a_\solid \to (\Lambda(Y_0), \Z)^a_\solid$ is weakly fs-descendable. This is implied by the fact that $\Lambda(X)^a_\solid \to \Lambda(Y_0)^a_\solid$ is weakly fs-descendable, which in turn follows from \cref{rslt:Tor-dim-and-desc-for-int-tor-coeffs-over-tot-disc}.
\end{proof}

\begin{definition}
For every map $f\colon Y \to X$ in $\vStacksCoeff$ we denote
\begin{align*}
	f^*\colon \Dqcohri(X, (\Lambda, \Z)) \rightleftarrows \Dqcohri(Y, (\Lambda, \Z)) \noloc f_*
\end{align*}
the pullback $f^*$, i.e. the ``restriction map'' of the sheaf $\Dqcohri(-, (\Lambda, \Z))$, and the pushforward $f_*$, i.e. the right adjoint of $f^*$ (this can be constructed similarly as in the proof of \cref{rslt:existence-of-pushforward-on-Dqcohri}). In the special case that $\Lambda = \ri^+_{X^\sharp}/\pi$ for some untilt $X^\sharp$ of $X$ and some pseudouniformizer $\pi$ on $X^\sharp$ we abbreviate
\begin{align*}
	\DqcohriXZ{X^\sharp} := \Dqcohri(X, (\ri^+_{X^\sharp}/\pi, \Z)).
\end{align*}
\end{definition}

Before we can come to the definition of $p$-bounded maps, we need to make one more definition. The following definition captures the property of a map $f\colon Y \to X$ of small v-stacks to have locally finite $\dimtrg$ (we need the extra definition because $\dimtrg f$ is not defined for general maps of small v-stacks).

\begin{definition}
A $0$-truncated map $f\colon Y \to X$ of small v-stacks has \emph{locally bounded dimension} if for every diamond $X'$ with a map $X' \to X$, there is a v-cover $(Z_i' \to Y')_{i\in I}$ of $Y' := Y \cprod_X X'$ by diamonds $Z_i'$ such that each map $f_i\colon Z_i' \to X'$ has $\dimtrg f_i < \infty$. We say that $f$ has \emph{bounded dimension} if $f$ has locally bounded dimension and is quasi-compact.
\end{definition}

\begin{remark}
If a $0$-truncated map $f\colon Y \to X$ of small v-stacks has locally bounded dimension, then for every diamond $X'$ with a map $X' \to X$ and every map $Z' \to Y' := Y \cprod_X X'$ from a diamond $Z'$, if $Z' \to X'$ has $\dimtrg < \infty$ then $Z' \cprod_{Y'} Z' \to X$ has locally bounded dimension. Indeed, it is clearly true for $Z' \cprod_{X'} Z' \to X'$ and since $f$ is $0$-truncated the map $Z' \cprod_{Y'} Z' \injto Z' \cprod_{X'} Z'$ is injective (this is why we assume $f$ to be $0$-truncated). In particular there is a hypercover $Z'_\bullet \to Y$ by diamonds $Z'_n$ such that all $Z'_n \to X'$ have locally $\dimtrg < \infty$. If one wants to drop the $0$-truncatedness hypothesis on $f$ then this hypercover property should be the definition of ``locally bounded dimension''.
\end{remark}

We can finally come to the main definition of the present subsection: that of $p$-bounded maps. Following the explanations at the beginning of the subsection we roughly define a map $f$ to be $p$-bounded if $f$ has universally finite $p$-cohomological dimension. The precise definition is as follows.

\begin{definition}
Let $f\colon Y \to X$ be a map of small v-stacks.
\begin{defenum}
	\item Assume that $X$ is an affinoid perfectoid space. We say that $f$ has \emph{$p$-cohomological dimension $\le d$} if for some pseudouniformizer $\pi$ on $X^\flat$ the functor
	\begin{align*}
		f_*\colon \DqcohriXZ{Y^\flat} \to \DqcohriXZ{X^\flat}
	\end{align*}
	has cohomological dimension $\le d$, i.e. sends $\D_{\ge0}$ to $\D_{\ge -d}$.

	\item \label{def:p-bounded} A map $f\colon Y \to X$ of small v-stacks is called \emph{$p$-bounded} if $f$ has locally bounded dimension and is locally separated and satisfies the following property: For every strictly totally disconnected space $X'$ with a map $X' \to X$ and pullback $f'\colon Y' := Y \cprod_X X' \to X'$ there is a v-cover of $Y'$ by quasiseparated subsheaves $Z'_i \subset Y'$ such that every $Z'_i \to X'$ has finite $p$-cohomological dimension.

	We say that $f$ is \emph{$p$-bounded of rank $\le d$} if in the above definition, the $Z_i'$ can be chosen disjoint and such that all $Z_i' \to X'$ have $p$-cohomological dimension $\le d$.
\end{defenum}
\end{definition}

\begin{remark}
In the definition of $p$-bounded maps, one can weaken the requirement ``locally separated'' slightly by only requiring that after pullback to any separated v-sheaf, $Y$ can be covered by subsheaves all of which are separated over $X$. However, we saw little use of this more general version of $p$-boundedness and chose to not make things more complicated than necessary.
\end{remark}

We begin the study of $p$-bounded maps with the following few results, which establish some of the technical properties we need in order to prove our main results on $p$-boundedness afterwards.

\begin{lemma}
\begin{lemenum}
	\item \label{rslt:p-bounded-implies-finite-cohom-dimension-on-DqcohriZ} Suppose $f\colon Y \to X$ is a qcqs $p$-bounded map in $\vStacksCoeff$. Then the functor
	\begin{align*}
		f_*\colon \Dqcohri(Y, (\Lambda, \Z)) \to \Dqcohri(X, (\Lambda, \Z))
	\end{align*}
	has finite cohomological dimension and preserves all small colimits.

	\item \label{rslt:p-bounded-implies-base-change-on-DqcohriZ} Let
	\begin{center}\begin{tikzcd}
		Y' \arrow[r,"g'"] \arrow[d,"f'"] & Y \arrow[d,"f"]\\
		X' \arrow[r,"g"] & X
	\end{tikzcd}\end{center}
	be a cartesian diagram in $\vStacksCoeff$ and assume that $f$ is $p$-bounded and qcqs. Then the natural base-change morphism
	\begin{align*}
		g^* f_* \isoto f'_* g'^*
	\end{align*}
	is an equivalence of functors from $\Dqcohri(Y, (\Lambda, \Z))$ to $\Dqcohri(X', (\Lambda, \Z))$.
\end{lemenum}
\end{lemma}
\begin{proof}
We first prove (i) in the case that $X$ is a strictly totally disconnected space. We first note that for any small v-stack $g\colon Z \to X$, if $g$ has finite $p$-cohomological dimension, then so does $g'\colon Z' \to X$ for any quasicompact injection $Z' \injto Z$. Namely, this follows from the fact that pushforward along $Z' \injto Z$ is $t$-exact, which itself follows easily from the fact that $Z' \injto Z$ is quasi-pro-étale by \cite[Corollary 10.6]{etale-cohomology-of-diamonds}.

Now by definition of $p$-boundedness there is some cover of $Y$ by quasiseparated subsheaves $Z_i \subset Y$, $i \in I$ such that all the maps $Z_i \to X$ have finite $p$-cohomological dimension. For each $i$ choose a cover of $Z_i$ by quasicompact v-sheaves $Z'_{ij}$ (e.g. totally disconnected spaces) and denote $Z_{ij} \subset Z_i$ the image of $Z'_{ij} \to Z_i$. Then the $Z_{ij}$ form a cover of $Y$ by qcqs subsheaves and by the previous paragraph each $Z_{ij} \to X$ has finite $p$-cohomological dimension. Thus, replacing the cover $(Z_i)_i$ by $(Z_{ij})_{ij}$ we can from now on assume that all $Z_i$ are qcqs. For every finite non-empty subset $J \subset I$ denote $Z_J := \bigisect_{i \in J} Z_i$; then by the discussion in the previous paragraph also $Z_J \to X$ has finite $p$-cohomological dimension. Moreover, we have
\begin{align*}
	\Dqcohri(Y, (\Lambda, \Z)) = \varprojlim_{J \subset I} \Dqcohri(Z_J, (\Lambda, \Z))
\end{align*}
and the pushforward $f_*\colon \Dqcohri(Y, (\Lambda, \Z)) \to \Dqcohri(X, (\Lambda, \Z))$ is computed as the limit of the functors $f_{J*}\colon \Dqcohri(Z_J, (\Lambda, \Z)) \to \Dqcohri(X, (\Lambda, \Z))$. Since $Y$ is quasicompact we can assume that $I$ is finite, hence to prove that $f_*$ has finite cohomological dimension and preserves small colimits it is enough to prove that each $f_{J*}$ has these properties. Therefore we can from now on assume that $f$ has $p$-cohomological dimension $\le d$ for some $d$.

If we have a map $\Lambda' \to \Lambda$ of integral torsion coefficients on $X$ and the claim holds for $\Lambda'$ then it also holds for $\Lambda$: This follows from the fact that the forgetful functor $\Dqcohri(-, (\Lambda, \Z)) \to \Dqcohri(-, (\Lambda', \Z))$ is $t$-exact and conservative, preserves colimits and commutes with pushforward (as in \cref{rslt:Dqcohri-change-of-coefficients}). We can therefore reduce to the case $\Lambda = \ri^+_{X^\sharp}/\pi$ for some untilt $X^\sharp$ of $X$ and some pseudouniformizer $\pi$ on $X^\sharp$. By assumption there is a pseudouniformizer $\varpi$ on $X^\flat$ such that $f^\flat_*\colon \Dqcohri(\ri^+_{Y^\flat}/\varpi, \Z) \to \Dqcohri(\ri^+_{X^\flat}/\varpi, \Z)$ has cohomological dimension $\le d$. By the standard arguments (using short exact sequences of the form $0 \to A^{\sharp+}/\pi^{n-1} \to A^{\sharp+}/\pi^n \to A^{\sharp+}/\pi \to 0$, cf. the proof of \cref{rslt:general-etale-descent-of-Dqcohri}) we easily deduce that also $f_*\colon \DqcohriXZ{Y^\sharp} \to \DqcohriXZ{X^\sharp}$ has cohomological dimension $\le d$, as desired. It remains to show that $f_*$ preserves small colimits. Thus suppose we are given a small family $(\mathcal M_i)_{i\in I}$ of objects in $\DqcohriXZ{Y^\sharp}$; we need to show that $f_* \bigdsum_i \mathcal M_i = \bigdsum_i f_* \mathcal M_i$. As $f_*$ has finite cohomological dimension, by looking at Postnikov towers we easily reduce the claim to the case that all $\mathcal M_i$ are uniformly bounded to the left. Pick any hypercover $Y_\bullet \to Y$ by totally disconnected spaces $Y_n = \Spa(B^n, B^{n+})$. Then $f_* \mathcal M = \varprojlim_{n\in\Delta} f_{n*} M^n$ for any object $\mathcal M = (M^n)_n \in \DqcohriXZ{Y^\sharp} = \varprojlim_{n\in\Delta} \DqcohriXZ{Y^\sharp_n}$, where $f_n\colon Y_n \to X$ is the projection. For all $n \ge 0$ we have $\DqcohriXZ{Y^\sharp_n} = \Dqcohri(B^{n+}/\pi, \Z)$ and since the $\mathcal M_i$ are uniformly bounded to the left, also $\bigdsum_i \mathcal M_i$ is bounded to the left. Thus $f_* \bigdsum_i \mathcal M_i$ is computed via a spectral sequence from the $f_{n*} \bigdsum_i M_i^n$, so it is enough to show that each $f_{n*}$ preserves direct sums. But $f_{n*}$ is just a forgetful functor. This finishes the proof of (i) in the case that $X$ is strictly totally disconnected.

We now prove (ii), so let everything be given as in the claim. We can formally reduce to the case that $X$ and $X'$ are strictly totally disconnected perfectoid spaces. Then by what we have already proved for part (i), $f_*\colon \Dqcohri(Y, (\Lambda, \Z)) \to \Dqcohri(X, (\Lambda, \Z))$ has finite cohomological dimension, i.e. there is some $d \ge 0$ such that $f_*$ maps $\Dqcohrige0$ to $\Dqcohrige{-d}$. Given any $\mathcal M \in \Dqcohri(Y, (\Lambda, \Z))$ we have $\mathcal M = \varprojlim_k \tau_{\le k} \mathcal M$ by left-completeness. Then $f_* \mathcal M = \varprojlim_k f_* \tau_{\le k} \mathcal M$ and by the cohomological boundedness of $f_*$ we note that the homotopy groups of $f_* \tau_{\le k} \mathcal M$ and $f_* \tau_{\le k-1} \mathcal M$ differ at most at the places $k$ to $k-d$; this implies that $g^*$ commutes with the limit. The same proof as in \cref{rslt:base-change-for-bounded-Dqcohri} shows that the desired base-change isomorphism holds on left-bounded sheaves; altogether we get
\begin{align*}
	g^* f_* \mathcal M = \varprojlim_k g^* f_* \tau_{\le k} \mathcal M = \varprojlim_k f'_* g'^* \tau_{\le k} \mathcal M = f'_* \varprojlim_k g'^* \tau_{\le k} \mathcal M = f'_* g'^* \mathcal M,
\end{align*}
as desired.

We can now finish the proof of (i) in the case of general $X$. Namely, by the base-change result (ii) and the sheafiness of $\Dqcohrim(-, (\Lambda, \Z))$ (see \cref{rslt:bounded-DqcohriZ-is-sheaf}) the claim about cohomological dimension can be checked after passing to any v-cover of $X$. Choosing such a v-cover consisting of strictly totally disconnected spaces, the claim follows. Similarly, the claim about preservation of colimits reduces by the base-change (ii) to the case that $X$ is strictly totally disconnected, which was handled above.
\end{proof}

\begin{lemma} \label{rslt:stronger-property-of-p-boundedness}
Let $Y \to X$ be a $p$-bounded map in $\vStacksCoeff$ and assume that $X$ is a quasicompact separated small v-sheaf. Then $Y$ is a v-sheaf and there is a cover of $Y$ by qcqs subsheaves $Z_i \subset Y$ such that each $g_i\colon Z_i \to X$ is $p$-bounded, separated, of bounded dimension, and the pushforward $g_{i*}\colon \Dqcohri(Z_i, (\Lambda, \Z)) \to \Dqcohri(X, (\Lambda, \Z))$ has finite cohomological dimension. The same is then true for every finite intersection of the $Z_i$'s.
\end{lemma}
\begin{proof}
By definition of $p$-bounded maps, $Y \to X$ is locally separated, so $Y$ is a v-sheaf (instead of a general v-stack) and there is a cover of $Y$ by open subspaces $U_j \subset Y$ such that each $U_j \to X$ is separated. For every $j$, choose a v-cover of $U_j$ by quasicompact small v-sheaves $Z'_{jk}$ (e.g. totally disconnected perfectoid spaces) and for every $k$ let $Z_{jk} \subset U_j$ denote the image of $Z'_{jk} \to U_j$. Then every $Z_{jk}$ is quasicompact and every map $Z_{jk} \to X$ is separated. After relabeling we get a cover of $Y$ by qcqs subsheaves $Z_i \subset Y$ such that each $Z_i \to X$ is separated. As $Y \to X$ has locally bounded dimension, every $Z_i \to X$ has bounded dimension. In view of \cref{rslt:p-bounded-implies-finite-cohom-dimension-on-DqcohriZ} the claim about finite cohomological dimension of the pushforward reduces to showing that each map $Z_i \to X$ is $p$-bounded. Thus, pick any cover $X' \surjto X$ by some strictly totally disconnected space $X'$ and denote $Y'$ and $Z'_i$ their base-change to $X'$. Then by $p$-boundedness of $Y \to X$ there is a cover of $Y'$ by quasiseparated subsheaves $W_l \subset Y'$ such that each $W_l \to X$ has finite $p$-cohomological dimension. Then the family $Z''_{il} \subset Z'_i \isect W_l \subset Z'_i$ forms a cover of $Z'_i$ by quasiseparated subsheaves and since every inclusion $Z''_{il} \subset W_l$ is quasicompact, the first part of the proof of \cref{rslt:p-bounded-implies-finite-cohom-dimension-on-DqcohriZ} shows that every map $Z''_{il} \to X$ has finite $p$-cohomological dimension. This proves that $Z_i$ is $p$-bounded.

Finally, we need to verify that every finite intersection of the $Z_i$'s satisfies the same properties. It is clear that this intersection is separated over $X$ and of bounded dimension. By the same argument as before it also follows that this intersection is $p$-bounded, as desired.
\end{proof}

\begin{lemma} \label{rslt:comparison-of-Dqcohri-and-DqcohriZ}
Let $f\colon Y \to X$ be a proper $p$-bounded map in $\vStacksCoeff$ and assume that $X$ is a totally disconnected perfectoid space. Then there are natural fully faithful functors
\begin{align*}
	\Dqcohri(X, \Lambda) \injto \Dqcohri(X, (\Lambda, \Z)), \qquad \Dqcohri(Y, \Lambda) \injto \Dqcohri(Y, (\Lambda, \Z))
\end{align*}
which preserve all small colimits and commute with $f^*$ and $f_*$. Moreover, the $t$-structure on $\Dqcohri(Y, (\Lambda, \Z))$ restricts to a left-complete $t$-structure on $\Dqcohri(Y, \Lambda)$ which is compatible with the full subcategories $\Dqcohriq(Y, \Lambda)$ for $? \in \{ +, -, b \}$.
\end{lemma}
\begin{proof}
The first functor is just the inclusion $\Dqcohri(\Lambda(X)) \injto \Dqcohri(\Lambda(X), \Z)$, which clearly satisfies all of the claimed properties. To construct the second functor, pick any hypercover $Y_\bullet \to Y$ by totally disconnected spaces with $\dimtrg(Y_n/X) < \infty$ and let $Y'_\bullet := \overline Y^{/X}_\bullet$. Then because $f$ is proper, $Y'_\bullet \to Y$ is still a hypercover and by \cref{rslt:compute-bounded-Dqcohri-for-qproet-over-Z'-Z} we have $\Dqcohri(Y'_n, \Lambda) = \Dqcohri(\Lambda(Y'_n)) = \Dqcohri(\Lambda(Y_n), \Lambda(X))$. Thus we have
\begin{align*}
	\Dqcohri(Y, \Lambda) &= \varprojlim_{n\in\Delta} \Dqcohri(\Lambda(Y_n), \Lambda(X)),\\
	\Dqcohri(Y, (\Lambda, \Z)) &= \varprojlim_{n\in\Delta} \Dqcohri(\Lambda(Y_n), \Z).
\end{align*}
Now for each $n$ there is a natural inclusion $\Dqcohri(\Lambda(Y_n), \Lambda(X)) \subset \Dqcohri(\Lambda(Y_n), \Z)$, which is compatible with all the transition functors in the above limit diagrams (these transition functors are just $- \tensor_{\Lambda(Y_n)^a} \Lambda(Y_m)^a$ in both cases). We therefore obtain the desired fully faithful functor $\Dqcohri(Y, \Lambda) \injto \Dqcohri(Y, (\Lambda, \Z))$. It is clear that this functor preserves all small colimits, as these are computed componentwise and the above inclusion of $\infty$-categories on each $Y_n$ preserves all small colimits. It is also clear that the constructed inclusion commutes with $f^*$ (which is computed as $- \tensor_{\Lambda(X)^a} \Lambda(Y_\bullet)^a$ in both cases) and with $f_*$ (which is computed as a totalization in both cases). The $t$-structure on $\Dqcohri(Y, \Lambda)$ is the one inherited from the $t$-structures on $\Dqcohri(\Lambda(Y_n), \Z)$ along the above limit diagram (using that all transition functors in the diagram are $t$-exact). Its compatibility with $\Dqcohriq(Y, \Lambda)$ follows directly from $\Dqcohriq(Y'_n, \Lambda) = \Dqcohriq(\Lambda(Y_n), \Lambda(X))$ (see \cref{rslt:compute-bounded-Dqcohri-for-qproet-over-Z'-Z}).

Regarding naturality: With a bit more effort it should be possible to show that the above construction is independent of the chosen hypercover $Y_\bullet \to Y$. A more elegant approach is to consider the subsite $\mathcal C \subset \vStacksCoeff$ spanned by all those spaces which are proper of bounded dimension over $X$ and compare the two sheaves $\Dqcohri(-, \Lambda)$ and $\Dqcohri(-, (\Lambda, \Z))$ using that they are the sheafifications of the presheaves $\Dqcohri(\Lambda(-))$ and $\Dqcohri(\Lambda(-), \Z)$. In any case, we do not need the naturality of the construction in the following, so we contend ourselves with the more explicit hypercover construction used above.
\end{proof}

Working directly with the definition of $p$-boundedness (\cref{def:p-bounded}) is a bit cumbersome, but in practice one should usually be able to deduce $p$-boundedness of a given map from simple ``building blocks'' using various stability properties. Two of these building blocks are quasi-pro-étale maps (see \cref{rslt:pro-etale-maps-are-p-bounded} below) and maps which are ``of finite type'' (see \cref{rslt:finite-type-implies-p-bounded} below). These two cases and various stability properties will be presented in the following.

\begin{lemma} \label{rslt:stability-of-p-bounded-maps}
\begin{lemenum}
	\item \label{rslt:p-boundedness-is-analytically-local} The property of being $p$-bounded is analytically local on source and target.

	\item \label{rslt:p-boundedness-is-v-local-on-target} Let $f\colon Y \to X$ be a map of small v-stacks which is locally separated and has locally bounded dimension. If there is a v-cover $\tilde X \to X$ such that the base-change $\tilde f\colon \tilde Y \to \tilde X$ is $p$-bounded, then $f$ is $p$-bounded.

	\item \label{rslt:p-bounded-maps-stable-under-base-change-composition} $p$-bounded maps are stable under composition and base-change.

	\item \label{rslt:p-boundedness-can-be-checked-after-compactification} Let $f\colon Y \to X$ be a separated map of small v-stacks. Then $f$ is $p$-bounded if and only if $\overline f^{/X}\colon \overline Y^{/X} \to X$ is $p$-bounded. If moreover $X$ is a separated v-sheaf then $f$ is $p$-bounded if and only if $\overline f\colon \overline Y \to \overline X$ is $p$-bounded.

	\item \label{rslt:pro-etale-maps-are-p-bounded} Every quasi-pro-étale map and every injection\footnote{A map $f\colon Y \to X$ of small v-stacks is \emph{injective} if for all $S \in \Perf$ the induced functor $Y(S) \to X(S)$ is fully faithful and injective on isomorphism classes. This property is stable under any base-change, and if $X$ is a v-\emph{sheaf} then $Y$ is also a v-sheaf and $f$ is injective as a map of v-sheaves.} of small v-stacks is $p$-bounded.

	\item \label{rslt:p-bounded-G-torsor-criterion} Let $Y \to X$ be a map of small v-stacks and $G$ a profinite group acting on $Y$ over $X$. If $Y \to X$ is $p$-bounded, $Y/G \to X$ is locally separated and $G$ has finite $p$-cohomological dimension then $Y/G \to X$ is $p$-bounded.

	\item \label{rslt:p-bounded-maps-satisfy-2-out-of-3} Let $f\colon Y \to X$ and $g\colon Z \to Y$ be maps of small v-stacks. If $f \comp g$ is $p$-bounded and $f$ is $0$-truncated (e.g. $p$-bounded), then $g$ is $p$-bounded.

	\item \label{rslt:p-boundedness-can-be-checked-on-fibers} Let $f\colon Y \to X$ be a map of small v-stacks which is locally separated and has locally bounded dimension. If there is some integer $d$ such that for every geometric point $x = \Spa(C, C^+) \to X$ the fiber $f_x\colon Y_x \to x$ is $p$-bounded of rank $\le d$, then $f$ is $p$-bounded.

	\item \label{rslt:p-boundedness-can-be-checked-on-connected-components-of-source} Let $f\colon Y \to X$ be a map of small v-stacks which is locally separated and representable in locally spatial diamonds with $\dimtrg f < \infty$. If there is some integer $d$ such that for every connected component $y \subset Y$ the map $y \to X$ is $p$-bounded of rank $\le d$, then $f$ is $p$-bounded.
\end{lemenum}
\end{lemma}
\begin{proof}
We start with the proof of (v): Assume first that $f\colon Y \to X$ is quasi-pro-étale. Then $f$ is locally separated by definition and clearly has locally bounded dimension. Moreover, after pullback to any strictly totally disconnected space, $Y$ admits a cover by qcqs open subspaces $U \subset Y$, which are then strictly totally disconnected as well (as they are pro-étale over a strictly totally disconnected space), hence $U \to X$ has $p$-cohomological dimension $0$. This proves that $f$ is $p$-bounded. Now assume that $f\colon Y \injto X$ is an injection of small v-stacks. Then clearly $f$ is (locally) separated and has locally bounded dimension (it is even representable in diamonds by \cite[Proposition 11.10]{etale-cohomology-of-diamonds}). After pullback to any strictly totally disconnected space, we can cover $Y$ by subsheaves $Y_i$ such that each $Y_i \to X$ is quasi-pro-étale (see \cite[Proposition 10.5]{etale-cohomology-of-diamonds}) and hence $p$-bounded. This implies that $Y \to X$ is also $p$-bounded.

We now prove (iii). It is clear that $p$-bounded maps are stable under any base-change. To handle composition, let $f\colon Y \to X$ and $g\colon Z \to Y$ be $p$-bounded maps of small v-stacks. Pick any strictly totally disconnected space $X'$ with a map $X' \to X$ and denote $f'\colon Y' \to X'$ and $g'\colon Z' \to Y'$ the base-changes. By \cref{rslt:stronger-property-of-p-boundedness} there is a v-cover of $Y'$ by qcqs substacks $Y'_j \subset Y'$ such that each $Y'_j \to X'$ is separated, $p$-bounded and has finite $p$-cohomological dimension. We can replace $Y'$ by $Y'_j$ to reduce to the case that $Y'$ is a separated and quasi-compact v-sheaf and the map $Y' \to X'$ has finite $p$-cohomological dimension. By applying \cref{rslt:stronger-property-of-p-boundedness} to $g$ we can find a cover of $Z'$ by qcqs subsheaves $Z'_i$ such that each map $Z'_i \to Y$ is separated and has finite $p$-cohomological dimension. But then every map $Z'_i \to X'$ has finite $p$-cohomological dimension, proving that $f \comp g$ is indeed $p$-bounded.

It is now easy to see (i): It is clear from the definition that a map of small v-stacks is $p$-bounded if it is so on an open cover of the source. Conversely, by (iii) and (v) the property of being $p$-bounded is preserved by passing to any open subset of the source; this proves that $p$-boundedness is analytically local on the source. To show that it is analytically local on the target, let $f\colon Y \to X$ be a map of small v-stacks and assume that there is an open cover $X = \bigunion_{i\in I} U_i$ by open substacks $U_i \subset X$ such that each $f_i\colon V_i := Y \cprod_X U_i \to U_i$ is $p$-bounded. By (iii) and (v) each $V_i \to X$ is $p$-bounded, so since $p$-boundedness is local on the source and the $V_i$'s form an open cover of $Y$, it follows that $f$ is $p$-bounded.

We now prove (ii), so assume that $f\colon Y \to X$ and $\tilde f\colon \tilde Y \to \tilde X$ are as in the claim. We can assume that $\tilde X$ is a totally disconnected space and by (i) we can further assume that $f$ is separated. It is then enough to show that $f$ is $p$-bounded after base-change to every totally disconnected space, so we can assume that $X$ is a separated v-sheaf, in which case the same is true for $Y$. Now cover $Y$ by qcqs subsheaves $Y_i \injto Y$ (e.g. let the $Y_i$'s be the images of a cover $(Y_i \to Y)_i$ by totally disconnected spaces $Y_i$). It follows from the definition of $p$-boundedness that it is enough to show that each map $Y_i \to X$ is $p$-bounded, so we can reduce to the case that $f$ is qcqs (use also (iii) and (v) to see that $\tilde f$ can still be assumed to be $p$-bounded). Then $\tilde f$ is qcqs, so by \cref{rslt:p-bounded-implies-finite-cohom-dimension-on-DqcohriZ} $\tilde f$ has finite $p$-cohomological dimension. By base-change for left-bounded complexes in $\Dqcohri(-, (\Lambda, \Z))$ (which follows in the same way as in \cref{rslt:base-change-for-bounded-Dqcohri}) and sheafiness of $\Dqcohrim(-, (\Lambda, \Z))$ (see \cref{rslt:bounded-DqcohriZ-is-sheaf}) we deduce that $f$ has finite $p$-cohomological dimension. As this argument holds after any base-change of $f$ to some totally disconnected space, we deduce that $f$ is $p$-bounded.

To prove (iv), we note that the property of having bounded dimension is preserved under taking compactifications, and that the condition of finite $p$-cohomological dimension in the definition of $p$-boundedness is invariant under taking compactifications by \cref{rslt:DqcohriZ-invariant-under-compactification}.

We now prove (vi), so let $X$, $Y$ and $G$ be given as in the claim. The situation is stable under base-change, so by (ii) we can assume that $X$ is a totally disconnected space. Passing to an open cover of $Y/G$ we can furthermore assume that $Y/G \to X$ is separated. Since $Y \to Y/G$ is separated, also $Y \to X$ is separated. As in the proof of (ii) we can construct a v-cover $(Z_i \injto Y/G)_i$ by quasicompact subsheaves $Z_i \subset Y/G$, which reduces the claim to each $Z_i$, hence allows us to w.l.o.g. assume that $f\colon Y/G \to X$ is quasicompact. Note that $f$ factors as $Y/G \to X/G \to X$, where $X/G$ is the (stack) quotient of $X$ by the trivial $G$-action. Now $Y \to X$ is $p$-bounded and qcqs, so by (ii) the same is true for $Y/G \to X/G$, hence by \cref{rslt:p-bounded-implies-finite-cohom-dimension-on-DqcohriZ} the map $Y/G \to X/G$ has finite $p$-cohomological dimension. On the other hand, $\Dqcohri(X/G, (\Lambda, \Z))$ is equivalent to the $\infty$-category of smooth $G$-representations on $(\Lambda(X), \Z)^a_\solid$-modules (by the same argument as in \cref{rslt:compute-Dqcohri-for-classifying-stack} applied to $\Dqcohri(-, (\Lambda, \Z))$ instead of $\Dqcohri(-, \Lambda)$) and the pushforward along $X/G \to X$ is computed by smooth $G$-cohomology. Thus by the same argument as in \cref{rslt:classifying-stack-over-tot-disc-pushforward-bounded} we deduce that the map $X/G \to X$ has finite $p$-cohomological dimension. Altogether we deduce that $f\colon Y/G \to X$ has finite $p$-cohomological dimension. As the same argument works after base-change along any map $X' \to X$ for some strictly totally disconnected $X'$, we deduce that $f$ is $p$-bounded.

We now prove (vii), so let $f\colon Y \to X$ and $g\colon Z \to Y$ be such that $f \comp g$ is $p$-bounded and $f$ is $0$-truncated. Then $Z \to Y$ factors as $Z \injto Z \cprod_X Y \to Y$, where the first map is the graph. But $Z \injto Z \cprod_X Y$ is injective, hence $p$-bounded by (v), and $Z \cprod_X Y \to Y$ is the base-change of $f \comp g$, hence also $p$-bounded; we are thus finished by (iii).

We now prove (viii), so let $f$ be given as in the claim. By the usual reductions as in the proof of (ii) we can assume that $X$ is strictly totally disconnected and $f$ is separated and quasicompact, in which case we want to show that $f$ has finite $p$-cohomological dimension. Thus, for some pseudouniformizer $\pi$ on $X^\flat$ we need to see that the functor $f_*\colon \DqcohriXZ{Y^\flat} \to \DqcohriXZ{X^\flat}$ has finite cohomological dimension. Let $\mathcal M \in \DqcohriXZ{Y^\flat}^\heartsuit$ be given. We claim that $f_* \mathcal M$ lies in $\D_{\ge-d}$. Indeed, this can be checked after applying $g_x^*\colon \DqcohriXZ{X^\flat} \to \DqcohriXZ{x^\flat}$ for every connected component $x = \Spa(C, C^+)$ of $X$, because the family of functors $g_x^*$ is $t$-exact and conservative (see \cref{rslt:compare-catsldmod-with-sheaves-on-pi-0}). By left-bounded base-change this reduces the claim to showing that every $f_x\colon Y_x \to x$ has $p$-cohomological dimension $\le d$, which follows immediately from the assumption that $f_x$ is $p$-bounded of rank $\le d$.

It remains to prove (ix), so let $f$ be given as in the claim. By the usual reductions as in the proof of (ii) we can assume that $X$ is a strictly totally disconnected perfectoid space and $f$ is separated and quasicompact. Define the topological space $T$ as the image of $\abs Y \to \abs X \cprod_{\pi_0(X)} \pi_0(Y)$. By assumption $Y$ is spatial, so by \cite[Proposition 12.14.(iv)]{etale-cohomology-of-diamonds} the map $\abs Y \to \abs X$ is spectral and generalizing. It follows that $\abs Y \to \abs X \cprod_{\pi_0(X)} \pi_0(Y)$ is also a spectral and generalizing map of spectral spaces and hence that $T$ is a generalizing and pro-constructible (see \cite[Lemma 2.3]{etale-cohomology-of-diamonds}) subset of $\abs X \cprod_{\pi_0(X)} \pi_0(T)$. Thus by \cite[Corollary 7.22]{etale-cohomology-of-diamonds} $T$ corresponds to some affinoid perfectoid space $X'$ which is pro-étale over $X$, and $Y \to X$ factors as $Y \to X' \to X$. By (ii) and (v) it is enough to show that $Y \to X'$ is $p$-bounded, so we can replace $X$ by $X'$ and hence assume that the map $\pi_0(Y) \isoto \pi_0(X)$ is bijective. But then the fiber of $f$ over each connected component of $X$ is a connected component of $Y$, so we are done by (viii).
\end{proof}

Having established some general stability properies of $p$-bounded maps, we now turn our focus on providing some explicit examples. The first example is that of a map of perfectoid spaces which is locally weakly of finite type, as follows.

\begin{definition} \label{def:perfectly-finite-type-in-Perf}
Let $f\colon Y \to X$ be a map of perfectoid spaces. We say that $f$ is \emph{locally weakly of perfectly finite type} if for every $y \in Y$ there are open neighborhoods $V \subset Y$ of $y$ and $U \subset X$ of $f(y)$ such that $U = \Spa(A, A^+)$ and $V = \Spa(B, B^+)$ are affinoid perfectoid, $f(V) \subset U$ and the induced map $\restrict fV\colon V \to U$ is weakly of perfectly finite type in the sense of \cref{def:perfectly-finite-type-in-AffPerf}.
\end{definition}

\begin{remark}
It is unclear to us whether a map $f\colon \Spa(B, B^+) \to \Spa(A, A^+)$ which is locally weakly of perfectly finite type (in the sense of \cref{def:perfectly-finite-type-in-Perf}) is automatically weakly of perfectly finite type (in the sense of \cref{def:perfectly-finite-type-in-AffPerf}). To avoid confusion, we refrain from saying that a map of perfectoid spaces is weakly of perfectly finite type if it is quasicompact and locally weakly of perfectly finite type.
\end{remark}

As remarked in the beginning of the proof of \cref{rslt:any-map-in-AffPerf-is-pro-wh}, the property of being (locally) of weakly perfectly finite type is invariant under passing to tilts and untilts and can hence be seen as an inherent property of the associated diamonds without any room for ambiguity.

\begin{lemma} \label{rslt:perf-fin-type-implies-p-bounded}
Let $f\colon Y \to X$ be a map of perfectoid spaces which is locally weakly of perfectly finite type. Then $f$ is $p$-bounded.

If $X = \Spa(A, A^+)$ and $Y = \Spa(B, B^+)$ are affinoid and $f$ is weakly of perfectly finite type then $f$ is $p$-bounded of rank bounded by a constant which only depends on the number of perfect generators for $B$ over $A$.
\end{lemma}
\begin{proof}
Since $p$-boundedness is analytically local on source and target (see \cref{rslt:p-boundedness-is-analytically-local}) we can assume that $X = \Spa(A, A^+)$ and $Y = \Spa(B, B^+)$ are affinoid perfectoid and $f$ is weakly of perfectly finite type, so we are reduced to showing the second claim. To verify the requirement for $p$-boundedness on $f$, we can now reduce to the case that $X$ is (strictly) totally disconnected, in which case it is enough to show that $f$ has $p$-cohomological dimension $\le c(n)$, where $c(n)$ is a constant which only depends on the number $n$ of perfect generators needed for $B$ over $A$. Fix any pseudouniformizer $\pi$ on $X^\flat$. We have $\DqcohriXZ{X^\flat} = \Dqcohri(A^{\flat+}/\pi, \Z)$ and $\DqcohriXZ{Y^\flat} = \Dqcohri(B^{\flat+}/\pi, \Z)$, so that $f_*$ is just the forgetful functor. It is tempting to guess that $f_*$ is $t$-exact, but the natural $t$-structure on $\DqcohriXZ{Y^\flat}$ may differ from the canonical $t$-structure on $\Dqcohri(B^{\flat+}/\pi, \Z)$, so we need a more elaborate argument, as follows.

Fix any $M \in \DqcohriXZ{Y^\flat}^\heartsuit$. Now pick a quasi-pro-étale cover $Z = \Spa(C, C^+) \surjto Y$ by some totally disconnected space $Z$ with projection $g\colon Z \to Y$. Then $\DqcohriXZ{Z^\flat} = \Dqcohri(C^{\flat+}/\pi, \Z)$ and by the definition of the $t$-structure on $\DqcohriXZ{Y^\flat}$ we have $g^* M \in \Dqcohrige0(C^{\flat+}/\pi, \Z)$. But by \cref{rslt:fin-type-over-K-implies-mod-pi-desc} the map $(B^{\flat+}/\pi, \Z)^a_\solid \to (C^{\flat+}/\pi, \Z)^a_\solid$ is fs-descendable of index bounded by some constant depending only on $n$, so by \cref{rslt:weakly-fs-descendable-implies-right-boundedness-descends}, there is a constant $b$ only depending on $n$ such that $M \in \Dqcohrige{b}(B^{\flat+}/\pi, \Z)$. Thus also $f_*M \in \Dqcohrige{b}(A^{\flat+}/\pi, \Z)$, as desired.
\end{proof}

The above study of $p$-bounded maps culminates in the following result, which produces the most important class of $p$-bounded maps: those which are (weakly) of finite type in the adic sense. This includes all maps between rigid-analytic varieties over some fixed base-field, which will allow us to construct all 6-functors in that setting without having to impose any further conditions.

\begin{proposition} \label{rslt:finite-type-implies-p-bounded}
Let $f\colon Y \to X$ be a map of analytic adic spaces over $\Spa\Z_p$ which is locally weakly of finite type. Then $f^\diamond\colon Y^\diamond \to X^\diamond$ is $p$-bounded.
\end{proposition}
\begin{proof}
Since $p$-boundedness is analytically local on source and target (see \cref{rslt:p-boundedness-is-analytically-local}) we can reduce to the case that $X$ and $Y$ are affinoid. Then $f$ factors as the composition of a Zariski closed immersion and a compactified relative ball. Since $p$-boundedness is stable under composition (see \cref{rslt:p-bounded-maps-stable-under-base-change-composition}) it is enough to consider the cases of Zariski closed immersions and compactified relative balls separately. On the other hand, Zariski closed immersions are quasi-pro-étale after applying $(-)^\diamond$ (by a similar argument as in \cite[Remark 7.9]{etale-cohomology-of-diamonds} or via \cite[Corollary 10.6]{etale-cohomology-of-diamonds}) and hence $p$-bounded by \cref{rslt:pro-etale-maps-are-p-bounded}. We are thus reduced to the case that $f$ is a compactified relative ball and by invoking the stability of $p$-boundedness under compositions we can even assume that $f$ has relative dimension $1$. Covering the ball by tori we can then reduce to the case that $f$ is a compactified relative $1$-dimensional torus. By \cref{rslt:p-boundedness-can-be-checked-on-fibers} it is enough to show that the geometric fibers of $f$ are strictly $p$-bounded of bounded rank, so we will from now on assume that $X = \Spa(C, C^+)$ is a geometric point. To simplify notation, we can reduce a bit further: Since $p$-boundedness is preserved under compactifications (see \cref{rslt:p-boundedness-can-be-checked-after-compactification} and note that the proof also works for $p$-bounded maps of given rank) we can assume $C^+ = \ri_C$ and that $f$ is the (non-compactified) torus over $X$.

We are now in the situation that $X = \Spa(C, \ri_C)$ and $Y = \Spa(C\langle T^{\pm1} \rangle, \ri_C\langle T^{\pm1} \rangle)$, where $C$ is an algebraically closed perfectoid field. If $C$ has characteristic $p$ then $Y^\diamond$ is perfectoid (the perfection of $Y$) and $f^\diamond$ is (weakly) of perfectly finite type, so we are done by \cref{rslt:perf-fin-type-implies-p-bounded}. If $C$ has characteristic $0$ then there is a $\Z_p$-torsor $\tilde Y = \Spa(C\langle T^{\pm1/p^\infty} \rangle, \ri_C \langle T^{\pm1/p^\infty} \rangle) \to Y$, where $\tilde Y$ is perfectoid and $\tilde Y \to X$ is (weakly) of perfectly finite type and hence strictly $p$-bounded (of bounded rank). As $\Z_p$ has $p$-cohomological dimension $1$, \cref{rslt:p-bounded-G-torsor-criterion} implies that $Y = \tilde Y/\Z_p \to X$ is strictly $p$-bounded of bounded rank (the claim about the rank can easily be deduced from the proof of the reference).
\end{proof}

We have produced a good amount of examples of $p$-bounded maps. There will be a few more examples at the end of this subsection (see \cref{rslt:p-bounded-on-aff-perfd-is-fs-desc-local-on-source,rslt:map-of-tot-disc-spaces-with-fin-dimtrg-is-p-bounded} below), but they require a better understanding of the general properties enjoyed by $p$-bounded maps, which will be our next goal to study. We start by showing that the pushforward along qcqs $p$-bounded maps preserves colimits and satisfies arbitrary base-change, which will be important to define the functors $f_!$ and $f^!$ later on.

\begin{proposition}
\begin{propenum}
	\item \label{rslt:qcqs-p-bounded-pushforward-preserves-colimits-and-Dqcohrim} Suppose $f\colon Y \to X$ is a qcqs $p$-bounded map in $\vStacksCoeff$. Then the functor
	\begin{align*}
		f_*\colon \Dqcohri(Y, \Lambda) \to \Dqcohri(X, \Lambda)
	\end{align*}
	preserves right-boundedness and all small colimits.

	\item \label{rslt:base-change-for-Dqcohri-along-qcqs-p-bounded-map} Let
	\begin{center}\begin{tikzcd}
		Y' \arrow[r,"g'"] \arrow[d,"f'"] & Y \arrow[d,"f"]\\
		X' \arrow[r,"g"] & X
	\end{tikzcd}\end{center}
	be a cartesian diagram in $\vStacksCoeff$ and assume that $f$ is $p$-bounded and qcqs. Then the natural base-change morphism
	\begin{align*}
		g^* f_* \isoto f'_* g'^*
	\end{align*}
	is an equivalence of functors from $\Dqcohri(Y, \Lambda)$ to $\Dqcohri(X', \Lambda)$.
\end{propenum}
\end{proposition}
\begin{proof}
We first prove (i) in the case that $X$ is a totally disconnected perfectoid space. By \cref{rslt:stronger-property-of-p-boundedness} there is a (w.l.o.g. finite) cover of $Y$ by qcqs subsheaves $Z_i \subset Y$ such that every finite intersection of the $Z_i$'s is $p$-bounded and separated over $X$. Thus by the same finite limit argument as in the proof of \cref{rslt:p-bounded-implies-finite-cohom-dimension-on-DqcohriZ} we can reduce to the case that $f$ itself is separated. Now factor $f$ as the composition of $j\colon Y \injto \overline Y^{/X}$ and $\overline f^{/X}\colon \overline Y^{/X} \to X$. It is enough to show that both $j_*$ and $\overline f^{/X}_*$ preserves right-boundedness and all small colimits. For the latter functor this follows immediately from \cref{rslt:p-bounded-implies-finite-cohom-dimension-on-DqcohriZ,rslt:comparison-of-Dqcohri-and-DqcohriZ}. To get the same for $j_*$, it is enough to show that $j_*$ is quasi-pro-étale (then (ii) is easy to show for $j$ in place of $f$, so by the sheafiness of $\Dqcohrim$ we reduce to the case that both source and target of $j$ are strictly totally disconnected, where the claim is obvious). To see that $j_*$ is indeed quasi-pro-étale, note first that $j$ is injective as a map of v-stacks (by \cite[Proposition 10.10]{etale-cohomology-of-diamonds} and the explicit definition of $\overline Y^{/X}$ in \cite[Proposition 18.6]{etale-cohomology-of-diamonds}). Moreover, $j$ is quasicompact because $Y$ is quasicompact and $\overline Y^{/X}$ is (quasi-)separated. But now \cite[Corollary 10.6]{etale-cohomology-of-diamonds} shows that $j$ is quasi-pro-étale. This finishes the proof of (i) in the case that $X$ is totally disconnected.

We now prove (ii), for which we can apply similar reductions as in (i): We can formally reduce to the case that $X$ and $X'$ are totally disconnected and then further reduce to the case that $f$ is separated. We can then factor $f$ into a quasi-pro-étale map $j$ and the relative compactification $\overline f^{/X}$. Base-change along $j$ is easy, so we are reduced to showing base-change along $\overline f^{/X}$. Thus we can from now on assume that $f$ is proper. Then \cref{rslt:comparison-of-Dqcohri-and-DqcohriZ} induces a left-complete $t$-structure on $\Dqcohri(Y, \Lambda)$ under which $f_*$ has finite cohomological dimension by \cref{rslt:p-bounded-implies-finite-cohom-dimension-on-DqcohriZ}. The rest of the proof of (ii) is now the same as in the proof of \cref{rslt:p-bounded-implies-base-change-on-DqcohriZ}.

Using (ii) we can deduce (i) in general, in the same way as in the proof of \cref{rslt:p-bounded-implies-finite-cohom-dimension-on-DqcohriZ}.
\end{proof}

In practice one is very often in the situation that one wants to study a small v-stack $X$ which admits a $p$-bounded map to some totally disconnected space; this is for example the case for rigid-analytic varieties over some algebraically closed field by \cref{rslt:finite-type-implies-p-bounded}. Small v-stacks $X$ which are of this form inherit some of the nice properties of totally disconnected spaces with regards to the behavior of $\Dqcohri(X, \Lambda)$ and hence deserve an extra name:

\begin{definition}
A small v-stack $X \in \vStacksCoeff$ is called \emph{$p$-bounded} if $X$ admits a $p$-bounded map (in $\vStacksCoeff$) to some totally disconnected perfectoid space.
\end{definition}

An important property of qcqs $p$-bounded small v-stacks is that their cohomology preserves all small colimits:

\begin{proposition} \label{rslt:global-sections-on-qcqs-p-bounded-v-stack-preserve-colimits}
Let $X \in \vStacksCoeff$ be qcqs and $p$-bounded. Then $X \in X_\vsite^\Lambda$ and the functor
\begin{align*}
	\Gamma(X, -)\colon \Dqcohri(X, \Lambda) \to \Dqcohri(\Lambda(X))
\end{align*}
preserves all small colimits.
\end{proposition}
\begin{proof}
By definition there is some $p$-bounded map $f\colon X \to Z$ to some totally disconnected space $Z$ such that the integral torsion coefficients $\Lambda$ extend to $Z$. Then $Z \in Z_\vsite^\Lambda$ and hence $X \in Z_\vsite^\Lambda$, so in particular $X \in X_\vsite^\Lambda$. Moreover, we have a commuting diagram
\begin{center}\begin{tikzcd}
	\Dqcohri(X, \Lambda) \arrow[r,"{\Gamma(X,-)}"] \arrow[d,"f_*"] & \Dqcohri(\Lambda(X)) \arrow[d,"\mathrm{forget}"]\\
	\Dqcohri(Z, \Lambda) \arrow[r,"{\Gamma(Z,-)}"] & \Dqcohri(\Lambda(Z))
\end{tikzcd}\end{center}
Since the forgetful functor on the right is conservative and preserves all small colimits, it is enough to show that the functor $\Gamma(Z, -) \comp f_*$ preserves all small colimits. But $\Gamma(Z, -)$ is an isomorphism and $f_*$ preserves colimits by \cref{rslt:qcqs-p-bounded-pushforward-preserves-colimits-and-Dqcohrim}.
\end{proof}

As a corollary we obtain the following variant of \cref{rslt:colim-of-pushforward-pullback}, which is closer to \cite[Proposition 14.9]{etale-cohomology-of-diamonds}:

\begin{corollary} \label{rslt:colimit-of-cohomology-qproet-unbounded-version}
Let $X \in \vStacksCoeff$ be a $p$-bounded diamond. Then for every $\mathcal M \in \Dqcohri(X, \Lambda)$ and every cofiltered inverse system $(U_i)_{i\in I}$ of qcqs diamonds in $X_\qproet$ with limit $U$, the natural morphism
\begin{align*}
	\varinjlim_i \Gamma(U_i, \mathcal M) \isoto \Gamma(U, \mathcal M)
\end{align*}
is an isomorphism in $\Dqcohri(\Lambda(X))$.
\end{corollary}
\begin{proof}
W.l.o.g. $I$ contains a final object $0 \in I$. Replacing $X$ by $U_0$ we can assume that $X$ is qcqs. Now combine \cref{rslt:colim-of-pushforward-pullback} and \cref{rslt:global-sections-on-qcqs-p-bounded-v-stack-preserve-colimits}.
\end{proof}

\begin{corollary} \label{rslt:etale-cohomology-is-conservative-on-p-bounded-diam}
Let $X \in \vStacksCoeff$ be a $p$-bounded locally spatial diamond. Then the family of functors $\Gamma(U, -)\colon \Dqcohri(X, \Lambda) \to \Dqcohri(\Lambda(X))$ for qcqs $U \in X_\et$ is conservative.
\end{corollary}
\begin{proof}
By \cite[Proposition 11.24]{etale-cohomology-of-diamonds} we can find a cover $(U_i \to X)_i$ by totally disconnected spaces $U_i$ such that each $U_i = \varprojlim_j U_{ij}$ is a cofiltered limit of qcqs $U_{ij} \in X_\et$. Then the family of functors $\Gamma(U_i, -)$ is conservative, and by \cref{rslt:colimit-of-cohomology-qproet-unbounded-version} the same follows for the family of functors $\Gamma(U_{ij}, -)$.
\end{proof}

In \cref{rslt:compute-Dqcohri-for-fin-type-over-tot-disc} we showed that if an affinoid perfectoid space $X = \Spa(A, A^+)$ is weakly of perfectly finite type over some totally disconnected space then $\DqcohriX X = \Dqcohri(A^+/\pi)$. This conclusion holds more generally if we only assume that $X$ is \emph{$p$-bounded} over some totally disconnected space, as we will show now. This allows more flexibility in proofs (although in practice \cref{rslt:compute-Dqcohri-for-qproet-over-fin-type-over-tot-disc} is often enough). As a preparation for the proof we need the following result.

\begin{lemma} \label{rslt:A-solid-S-is-flat-on-discrete-modules}
Let $A$ be a classical ring and $S$ a profinite set. Then the functor
\begin{align*}
	- \tensor_{A_\solid} A_\solid[S]\colon \D(A)_\omega \to \D_\solid(A)
\end{align*}
is $t$-exact.
\end{lemma}
\begin{proof}
It is enough to check $t$-exactness after composition with the forgetful functor $\D_\solid(A) \to \D_\solid(\Z)$. As such, the functor $- \tensor_{A_\solid} A_\solid[S]$ can be written as
\begin{align*}
	- \tensor_{A_\solid} A_\solid[S] = \varinjlim_i - \tensor_{A_{i\solid}} A_{i\solid}[S],
\end{align*}
where $(A_i)_i$ is a filtered diagram of finite-type $\Z$-algebras such that $A = \varinjlim_i A_i$. We can therefore reduce to the case that $A = A_i$ for some $i$, i.e. that $A$ is of finite type over $\Z$. Then $A_\solid[S] = \prod_I A$ for some set $I$ (see \cref{rslt:A-solid-is-analytic-ring}). We need to show that for every classical $A$-module $M$ the solid $A$-module $M \tensor_{A_\solid} \prod_I A$ is concentrated in degree $0$. By writing $M$ as a filtered colimit of finite $A$-modules, we can reduce to the case that $M$ is finite. Then $M$ admits a projective resolution by finite free $A$-modules. Clearly this resolution remains exact after tensoring with $\prod_I A$, as desired.
\end{proof}

\begin{theorem} \label{rslt:compute-Dqcohri-for-p-bounded-over-tot-disc}
If $X \in \AffPerfCoeff$ is $p$-bounded then there is a natural equivalence
\begin{align*}
	\Dqcohri(X, \Lambda) = \Dqcohri(\Lambda(X)).
\end{align*}
\end{theorem}
\begin{proof}
We claim that the natural pair of adjoint functors
\begin{align*}
	\widetilde{(-)}\colon \Dqcohri(\Lambda(X)) \leftrightarrows \Dqcohri(X, \Lambda)\noloc \Gamma(X, -)
\end{align*}
forms the desired equivalence. This reduces to the following two claims:
\begin{enumerate}[(a)]
	\item $\widetilde{(-)}$ is fully faithful, i.e. the unit of the adjunction $M \isoto \Gamma(X, \tilde M)$ is an equivalence for all $M \in \Dqcohri(\Lambda(X))$.
	\item $\Gamma(X, -)$ is conservative.
\end{enumerate}
To prove (a) we note that $\Gamma(X, -)$ preserves all small colimits by \cref{rslt:global-sections-on-qcqs-p-bounded-v-stack-preserve-colimits}, hence both sides of the desired isomorphism $M \isoto \Gamma(X, \tilde M)$ commute with colimits in $M$. We can thus reduce to the case $M = \Lambda(X)^a_\solid[S]$ for some profinite set $S$. In this case we can compute $\Gamma(X, \tilde M)$ as follows. Let $f\colon X \to Z$ be a $p$-bounded map with $Z$ totally disconnected and let $Y_\bullet \to X$ be a pro-étale Čech cover by totally disconnected spaces $Y_n$. Since $f$ has bounded dimension we have $\dimtrg(Y_n/X) < \infty$ for all $n$. Note that $Y'_n := \overline{Y_n}^{/X} \injto \overline{Y_n}^{/Z}$ is quasi-pro-étale, so by \cref{rslt:compute-bounded-Dqcohri-for-qproet-over-Z'-Z} we have $\Dqcohri(Y'_n, \Lambda) = \Dqcohri(\Lambda(Y'_n)) = \Dqcohri(\Lambda(Y_n), \Lambda(X))$. Since $Y'_\bullet \to X$ is still a Čech cover we have
\begin{align*}
	\Dqcohri(X, \Lambda) = \varprojlim_{n\in\Delta} \Dqcohri(Y'_n, \Lambda) = \varprojlim_{n\in\Delta} \Dqcohri(\Lambda(Y_n), \Lambda(X)).
\end{align*}
Thus for $M = \Lambda(X)^a_\solid[S]$ we can compute $\Gamma(X, \tilde M)$ with the complex
\begin{align*}
	M \tensor_{\Lambda(X)^a_\solid} \Lambda(Y_0)^a \to M \tensor_{\Lambda(X)^a_\solid} \Lambda(Y_1)^a \to M \tensor_{\Lambda(X)^a_\solid} \Lambda(Y_2)^a \to \dots
\end{align*}
which we want to be quasi-isomorphic to $M$. By \cref{rslt:A-solid-S-is-flat-on-discrete-modules} we can reduce to the case $M = \Lambda(X)^a$, in which case the above complex becomes the complex $\Lambda(Y_0)^a \to \Lambda(Y_1)^a \to \dots$, which computes Čech cohomology of the sheaf $\Lambda^a$ on $X_\vsite$. But by definition of integral torsion coefficients we have $H^i(X, \Lambda)^a = 0$ for $i > 0$, proving that this complex is indeed quasi-isomorphic to $\Lambda(X)^a$. This finishes the proof of (a).

We now prove (b), i.e. that $\Gamma(X, -)$ is conservative. Thus suppose we are given a morphism $\mathcal M \to \mathcal N$ in $\Dqcohri(X, \Lambda)$ such that the induced morphism $\Gamma(X, \mathcal M) \to \Gamma(X, \mathcal N)$ is an isomorphism. Pick any cover $U = \varprojlim U_i \to X$ with $U$ totally disconnected and all $U_i \to X$ étale. Then $\Gamma(U, -)$ is conservative by \cref{rslt:def-of-qcoh-Lambda-modules} and hence the collection of functors $\Gamma(U_i, -)$, $i \in I$ is conservative by \cref{rslt:colimit-of-cohomology-qproet-unbounded-version}. We thus need to show that $\Gamma(U_i, \mathcal M) \to \Gamma(U_i, \mathcal N)$ is an isomorphism for every $i \in I$. From now on fix some $i \in I$. By \cref{rslt:any-map-in-AffPerf-is-pro-wh} we can write $X$ as a cofiltered limit $X = \varprojlim_{j\in J} Z_j$ with $Z_j$ weakly of perfectly finite type over $Z$. Then by \cite[Proposition 6.4.(iv)]{etale-cohomology-of-diamonds} we have $U_i = V \cprod_{Z_j} X$ for some $j \in J$ and some étale $g\colon V \to Z_j$. But by \cref{rslt:def-of-qcoh-Lambda-modules} we have $\Dqcohri(Z_j, \Lambda) = \Dqcohri(\Lambda(Z_j))$, hence the pushforward $f_{j*}\colon \Dqcohri(X, \Lambda) \to \Dqcohri(Z_j, \Lambda)$ along $f_j\colon X \to Z_j$ factors over $\Gamma(X, -)$. Thus $f_{j*} \mathcal M \isoto f_{j*} \mathcal N$ is an isomorphism, and by base-change (see \cref{rslt:base-change-for-Dqcohri-along-qcqs-p-bounded-map}) we have $\Gamma(U_i, -) = g^* f_{j*}$, so that $\Gamma(U_i, \mathcal M) \isoto \Gamma(U_i, \mathcal N)$ is an isomorphism, as desired.
\end{proof}

\begin{warning}
Suppose $X \in \AffPerfCoeff$ is $p$-bounded. Then by \cref{rslt:compute-Dqcohri-for-p-bounded-over-tot-disc} we have $\Dqcohri(X, \Lambda) = \Dqcohri(\Lambda(X))$, so in particular $\Dqcohri(X, \Lambda)$ inherits a $t$-structure from $\Dqcohri(\Lambda(X))$. However, this $t$-structure may in general not be compatible with the subcategories $\Dqcohrip$ and $\Dqcohrim$ defined in \cref{def:Dqcohri-bounded-subcategories}!
\end{warning}

As a corollary of \cref{rslt:compute-Dqcohri-for-p-bounded-over-tot-disc} we obtain an even more general class of small v-stacks $X$ on which we can explicitly compute $\Dqcohri(X, \Lambda)$, see in particular \cref{rslt:explicit-formula-for-Dqcohri-on-torus} below.

\begin{corollary} \label{rslt:compute-Dqcohri-on-stack-quotient-of-p-bounded-space}
Let $X \in \AffPerfCoeff$ be $p$-bounded, let $G$ be a locally profinite group acting on $X$ and let $X/G$ denote the (stack) quotient. Assume that $\Lambda$ extends to $X/G$ and that $G$ has locally finite $p$-cohomological dimension (see \cref{def:locally-finite-cohom-dim}).
\begin{corenum}
	\item There is a natural equivalence $\Dqcohri(X/G, \Lambda) = \Drepsldasm G{\Lambda(X)}$, where the action of $G$ on $\Lambda(X)$ is induced by the existence of $\Lambda$ on $X/G$.

	\item Let $Y \in \AffPerfCoeff$ be also $p$-bounded, let $X \to Y$ be a $G$-equivariant map in $\AffPerfCoeff$ (with $G$ acting trivially on $Y$) and let $f\colon X/G \to Y$ be the induced map. Then
	\begin{align*}
		f_*\colon \Dqcohri(X/G, \Lambda) \to \Dqcohri(Y, \Lambda)
	\end{align*}
	coincides with smooth group cohomology $\Gamma_\sm(G, -)$ via the equivalence in (i).
\end{corenum}
\end{corollary}
\begin{proof}
Using \cref{rslt:compute-Dqcohri-for-p-bounded-over-tot-disc} we can apply the same argument as in \cref{rslt:compute-Dqcohri-for-classifying-stack} (note that all $X \cprod G^n$ are $p$-bounded, as they are pro-étale over $X$).
\end{proof}

\begin{example} \label{rslt:explicit-formula-for-Dqcohri-on-torus}
Let $K$ be an algebraically closed non-archimedean field of characteristic $0$ and let $X = \Spa(K\langle T^\pm\rangle, \ri_K\langle T^\pm\rangle)$ be the torus over $\Spa(K, \ri_K)$. Then for every pseudouniformizer $\pi \in K$ we have
\begin{align*}
	\DqcohriX X = \Drepsldasm{\Z_p}{\ri_K/\pi[T^{\pm 1/p^\infty}]},
\end{align*}
where $\Z_p$ acts on $\ri_K/\pi[T^{\pm 1/p^\infty}]$ via $\gamma \times T^{1/p^n} := \zeta_{p^n}^\gamma T^{1/p^n}$, where $(\zeta_{p^n})_n$ is a compatible system of $p^n$-th roots of unity in $K$ (cf. \cite[Lemma 5.5]{rigid-p-adic-hodge}). Moreover, the pushforward
\begin{align*}
	f_*\colon \DqcohriX X \to \Dqcohri(\ri_K/\pi)
\end{align*}
is computed by smooth group cohomology $\Gamma_\sm(\Z_p, -)$.
\end{example}

As promised above, we now present more $p$-boundedness criteria. They rely on \cref{rslt:compute-Dqcohri-for-p-bounded-over-tot-disc} and thus needed to be placed after that result.

\begin{lemma} \label{rslt:p-bounded-on-aff-perfd-is-fs-desc-local-on-source}
Let $Z \to Y \to X$ be maps of affinoid perfectoid spaces with the following properties:
\begin{enumerate}[(a)]
	\item The maps $Z \to X$ and $Z \to Y$ are $p$-bounded of some rank $\le r$.

	\item The map $Z \to Y$ is a v-cover and there is a pseudouniformizer $\pi$ on $X^\flat$ such that $Z \to Y$ is mod-$\pi$ weakly fs-descendable of some index $d$.
\end{enumerate}
Then $Y \to X$ is $p$-bounded of rank bounded by a constant which only depends on $d$ and $r$.
\end{lemma}
\begin{proof}
We can assume that $X$ is strictly totally disconnected, in which case we need to show that the map $Y \to X$ has $p$-cohomological dimension $\le c(d, r)$, where $c(d, r)$ is a constant only depending on $d$ and $r$. Denote $Z_\bullet \to Y$ the Čech nerve of $Z \to Y$ and let $Z_n = \Spa(C^n, C^{n+})$. We implicitly trat all affinoid perfectoid spaces as tilted to avoid the superscript $\flat$ everywhere. Fix a pseudouniformizer $\pi$ on $X$ as in (b). By (a) and \cref{rslt:p-bounded-maps-stable-under-base-change-composition} the maps $Z_n \to X$ are $p$-bounded, hence by \cref{rslt:compute-Dqcohri-for-p-bounded-over-tot-disc} we have $\DqcohriXZ{Z_n} = \Dqcohri(C^{n+}/\pi, \Z)$ for all $n$. By (b) the map $Z \to Y$ is mod-$\pi$ weakly fs-descendable (of index $\le d$), hence also $\DqcohriXZ Y = \Dqcohri(B^+/\pi, \Z)$ (here we implicitly use \cref{rslt:discrete-desc-equiv-abs-compactification-desc,rslt:descendable-implies-discretely-descendable} to see that $(B^+/\pi, \Z)^a_\solid \to (C^+/\pi, \Z)^a_\solid$ is weakly fs-descendable). In particular, the pushforward along $f\colon Y \to X$ is just a forgetful functor. However, as in the proof of \cref{rslt:perf-fin-type-implies-p-bounded} the $t$-structure on $\DqcohriXZ Y$ may differ from the canonical $t$-structure on $\Dqcohri(B^+/\pi, \Z)$, so the pushforward needs not be $t$-exact. We can argue as follows: Pick some $M \in \Dqcohrige0(\ri^+_Y/\pi, \Z)$. Then the pullback $g^* M$ of $M$ along $g\colon Z \to Y$ lies in $\Dqcohrige0(\ri^+_Z/\pi, \Z)$. Since $h\colon Z \to X$ has $p$-cohomological dimension $\le r$, the pushforward $h_* (g^* M)$ lies in $\Dqcohrige{-r}(A^+/\pi, \Z)$. But $h_*\colon \Dqcohri(C^{0+}/\pi, \Z) \to \Dqcohri(A^+/\pi, \Z)$ is just the forgetful functor, which implies $g^* M \in \Dqcohrige{-r}(C^{0+}/\pi, \Z)$. Finally, since $(B^+/\pi, \Z)^a_\solid \to (C^{0+}/\pi, \Z)^a_\solid$ is weakly fs-descendable of index $\le d$, \cref{rslt:weakly-fs-descendable-implies-right-boundedness-descends} implies that $M \in \Dqcohrige{c(d,r)}(B^+/\pi, \Z)$, where $c(r,d) = -r + b(d)$. In particular $f_* M \in \Dqcohrige{c(d,r)}(A^+/\pi, \Z)$, hence $f$ has $p$-cohomological dimension $\le c(d,r)$, as desired.
\end{proof}

\begin{lemma} \label{rslt:map-of-tot-disc-spaces-with-fin-dimtrg-is-p-bounded}
Let $f\colon Y \to X$ be a map of totally disconnected spaces with $\dimtrg f < \infty$. Then $f$ is $p$-bounded.
\end{lemma}
\begin{proof}
We claim that if $Y$ is connected then $f$ is $p$-bounded of rank bounded by a constant $t(d)$ only depending on $d = \dimtrg f$. By \cref{rslt:p-boundedness-can-be-checked-on-connected-components-of-source} we can now reduce to the case that $Y = \Spa(K', K'^+)$ is connected. Then $f$ factors over a connected component of $X$, so we can assume that $X = \Spa(K, K^+)$ is also connected. Assume first that $\trc(K'/K) = d$. Then $Y \to X$ is a composition of a pro-étale map and a $d$-dimensional compactified perfectoid ball. The latter map is $p$-bounded of bounded rank by \cref{rslt:perf-fin-type-implies-p-bounded} and one deduces easily that the same is true after composition with a pro-étale map. This finishes the proof in the case $\trc(K'/K) \le d$.

In general, recall that by definition of $\dimtrg$ there is some field extension $K''/K'$ such that $\trc(K''/K) = d$. Let $Z := \Spa(K'', K''^+)$, where $K''^+$ is an extension of the valuation ring $K'^+$. By the previous paragraph, both $Z \to X$ and $Z \to Y$ are $p$-bounded of rank bounded by a constant which only depends on $d$. Moreover, $Z \to Y$ is mod-$\pi$ fs-descendable of index $\le 4$ by \cref{rslt:v-cover-of-tot-disc-is-mod-pi-descendable}. Thus we can apply \cref{rslt:p-bounded-on-aff-perfd-is-fs-desc-local-on-source} to conclude that $Y \to X$ is $p$-bounded of bounded rank.
\end{proof}

\subsection{The 6-Functor Formalism} \label{sec:ri-pi.6-functor}

We are finally in the position to construct the 6-functor formalism for quasicoherent $\ri^{+a}_X/\pi$-modules. In \cref{def:4-functors} we have already introduced the four functors $f_*$, $f^*$, $- \tensor -$ and $\IHom(-, -)$. It remains to construct the functor $f_!$ and its right adjoint $f^!$ and to verify all the expected compatibilities. As we are working in an $\infty$-categorical setting, these ``compatibilities'' involve infinitely many coherences of higher homotopies, which are hard to check by hand. Fortunately, all of the higher coherences fall in place automatically once we check everything on the $1$-categorical level: This is shown in \cref{sec:infcat.sixfun}, following previous work by Liu-Zheng \cite{enhanced-six-operations}. Thus our strategy will be to check everything on the $1$-categorical level by hand and then use the machinery of loc. cit. to get an elegant formulation of the 6-functor formalism including all higher coherences ``for free''.

Given a ``nice'' qcqs map $f\colon Y \to X$ in $\vStacksCoeff$, the functor $f_!\colon \Dqcohri(Y, \Lambda) \to \Dqcohri(X, \Lambda)$ will be defined as
\begin{align*}
	f_! := g_* \comp j_!
\end{align*}
for any composition $f = g \comp j$, where $j\colon Y \to Z$ is étale and $g\colon Z \to X$ is proper. Here $j_!$ will be the left adjoint of $j^*$, which will be our first goal to construct. This is the bridge between the lower-shriek functor for quasicoherent sheaves on diamonds and quasicoherent sheaves on schemes, so unsurprisingly it relies on the following result:

\begin{lemma} \label{rslt:tot-disc-space-mod-pi-isom-to-scheme}
Let $X = \Spa(A, A^+)$ be a totally disconnected perfectoid space, $\pi \in A$ a pseudouniformizer and $\tilde X := \Spec A^+/\pi$. Then there is a natural isomorphism
\begin{align*}
	(\abs{X}, \ri^+_X/\pi) = (\abs{\tilde X}, \ri_{\tilde X})
\end{align*}
of locally ringed spaces.
\end{lemma}
\begin{proof}
There is a uniquely determined map $\psi\colon (\abs{X}, \ri^+_X/\pi) \to (\abs{\tilde X}, \ri_{\tilde X})$ of locally ringed spaces induced by the identity map on global sections (see \cite[Lemma 01I1]{stacks-project}). We will show that $\psi$ is an isomorphism.

If $X = \Spa(K, K^+)$ is connected, then the claim is easy: $\abs X$ is parametrized by valuation rings $K^+ \subseteq V \subseteq \ri_K$ and $\abs{\tilde X}$ is parametrized by valuation rings $K^+/\mm_K \subseteq \overline V \subseteq \ri_K/\mm_K$; one checks easily that both sets of valuation rings agree via taking images/preimages under $\ri_K \surjto \ri_K/\mm_K$. In the general case, note that $\pi_0(X) = \pi_0(\tilde X)$ because both sides are determined by the idempotents in $A$ (which automatically lie in $A^+$). Together with the connected case we already deduce that $\psi$ is bijective on the topological spaces.

Consider a distinguished open subset $\tilde U = D(\overline f) \subset \tilde X$ for some $\overline f \in A^+/\pi$. Choose any lift $f \in A^+$ of $\overline f$. Then one checks easily that $\tilde U = \{ \abs f \ge 1 \} \subset X$ (this reduces to the case that $X$ is connected, where it is easy), so in particular $\tilde U$ is open in $X$. Thus to prove that $\psi$ is a homeomorphism $\abs X \isoto \abs{\tilde X}$, we need to see that every open subset of $\abs X$ is a union of $D(\overline f)$'s. Given $U \subset X$ open, let $x \in U$ be any point and let $c \in \pi_0(X)$ be the connected component of $x$. Then since $c = \Spa(K, K^+)$ for a perfectoid field $K$, there is some $\overline{f_c} \in K^+/\pi$ such for any lift $f_c \in K^+$ we have $U \isect c = \{ \abs{f_c} \ge 1 \}$. Let $\pi\colon X \to \pi_0(X)$ denote the projection. Then $K^+/\pi$ is the colimit of $\ri^+_X(\pi^{-1} W)/\pi$ for compact open neighbourhoods $W$ of $c$; in particular we can extend $\overline{f_c}$ to an element $\overline{f_W} \in \ri^+_X(\pi^{-1} W)/\pi$ for some such $W$. Pick any lift $f_W \in \ri^+_X(\pi^{-1}W)$ of $\overline{f_W}$. Then
\begin{align*}
	\bigisect_{c\in W' \subseteq W} (X \setminus U \isect \{ \abs{f_W} \ge 1 \}) \isect \pi^{-1} W' = \emptyset.
\end{align*}
Here $W'$ ranges through all compact open neighbourhoods of $c \in \pi_0(X)$ contained in $W$. Each term in the above intersection is pro-constructible and in particular closed in the constructible topology on $X$. It follows that some term must already be empty (otherwise the union of the complements would be an open cover in the constructible topology without finite subcover), i.e. there is some $c \in W' \subset W$ such that
\begin{align*}
	(X \setminus U \isect \{ \abs{f_W} \ge 1 \}) \isect \pi^{-1} W' = \emptyset
\end{align*}
Now replace $f_W$ by the function $f \in A^+$ which is equal to $f_W$ on $W'$ and equal to $0$ outside of $W'$. Then
\begin{align*}
	\{ \abs f \ge 1 \} = \{ \abs{f_W} \ge 1 \} \isect \pi^{-1}W' \subset U,
\end{align*}
i.e. if $\overline f \in A^+/\pi$ is the reduction of $f$ mod $\pi$ then $D(\overline f) \subset U$. Repeating the procedure for all $x \in U$ we deduce that $U$ is a union of $D(\overline f)$'s, as desired. This finally proves that $\abs X = \abs{\tilde X}$ as topological spaces.

It remains to show that the morphism of sheaves $\ri_{\tilde X} \to \ri^+_X/\pi$ induced by $\psi$ is an isomorphism. This can be checked on stalks, so we are back in the situtation that $X = \Spa(K, K^+)$ is connected, where the claim is easy.
\end{proof}

We are now prepared to prove the existence of $j_!$. The following result is analogous to \cite[Proposition 19.1]{etale-cohomology-of-diamonds}:

\begin{lemma} \label{rslt:existence-of-lower-shriek-for-etale-maps}
Let $j\colon U \to X$ be an étale map in $\vStacksCoeff$. Then $j^*\colon \Dqcohri(X, \Lambda) \to \Dqcohri(U, \Lambda)$ admits a left adjoint
\begin{align*}
	j_!\colon \Dqcohri(U, \Lambda) \to \Dqcohri(X, \Lambda).
\end{align*}
Moreover, for every map $g\colon X' \to X$ of small v-stacks with base-change $g'\colon U' := U \cprod_X X' \to X'$, $j'\colon U' \to U$, the natural transformation
\begin{align*}
	j'_! g'^* \isoto g^* j_!
\end{align*}
of functors $\Dqcohri(U, \Lambda) \to \Dqcohri(X', \Lambda)$ is an equivalence.
\end{lemma}
\begin{proof}
First assume that $X$ is strictly totally disconnected. We note that in this case the claimed base-change for $j_!$ is automatic (assuming the existence of $j_!$), by the following argument. Passing to right-adjoints, we need to show that the natural transformation $j^* g_* \isoto g'_* j'^*$ is an isomorphism. First assume that $j$ is quasicompact and separated (in particular $U$ is strictly totally disconnected). Choosing any v-hypercover $X'_\bullet \to X'$ such that each $X'_n$ is a disjoint union of totally disconnected spaces and using that $j^*$ preserves limits (by the existence of the left adjoint $j_!$) we can reduce to the case that $X'$ is strictly totally disconnected (hence also $U'$ is strictly totally disconnected). Then the base-change follows from \cref{rslt:base-change-for-int-tor-coeffs-on-aff-perf}. Now assume that $j$ is only separated (not necessarily quasicompact). Then $U$ is still a perfectoid space, so we can write $U = \bigunion_{i\in I} U_i$ as a filtered union of quasicompact open subspaces $U_i \subset U$. Then $\Dqcohri(U, \Lambda) = \varprojlim_i \Dqcohri(U_i, \Lambda)$ and all $U_i$ are strictly totally disconnected, so the claimed base-change follows easily from the case of qcqs $j$ (apply it both to the maps $U_i \to X$ and the maps $U_j \injto U_i$). Finally, if $j$ is general then the open subsets of $U$ which are separated over $X$ form a basis of $U$; we can therefore argue as before.

Still in the case of $X$ being strictly totally disconnected, we now prove the existence of $j_!$. First assume that $j$ is quasicompact and separated. Then $U$ is a disjoint union of quasicompact open subsets of $X$, so the claim reduces to the case that $j\colon U \injto X$ is a quasicompact open subset. Thus by \cite[Lemmas 7.5, 7.6]{etale-cohomology-of-diamonds} $U$ is also totally disconnected and in particular affinoid. Pick a morphism of integral torsion coefficients $\ri^+_{X^\sharp}/\pi \to \Lambda$ for some untilt $X^\sharp = \Spa(A, A^+)$ of $X$ and some pseudouniformizer $\pi \in A$. Denote $U^\sharp = \Spa(B, B^+)$. By \cref{rslt:tot-disc-space-mod-pi-isom-to-scheme}, $\Spec B^+/\pi \injto \Spec A^+/\pi$ is an open immersion of schemes, hence for $\Lambda(U)' := \Lambda(X) \tensor_{A^+/\pi} B^+/\pi$ also $j'\colon \Spec \Lambda(U)' \injto \Spec \Lambda(X)$ is an open immersion of schemes. By \cref{rslt:lower-shriek-commutes-with-almost-localization} it follows that the pullback functor $j'^{a*}\colon \Dqcohri(\Lambda(X)) \to \Dqcohri(\Lambda(U)')$ admits a left adjoint $j'^a_!\colon \Dqcohri(\Lambda(U)') \to \Dqcohri(\Lambda(X))$. But by the definition of morphisms of integral torsion coefficients we have $\Lambda(U)'^a_\solid = \Lambda(U)^a_\solid$, so by \cref{rslt:def-of-qcoh-Lambda-modules} we have $\Dqcohri(X, \Lambda) = \Dqcohri(\Lambda(X))$ and $\Dqcohri(U, \Lambda) = \Dqcohri(\Lambda(U))$; all in all this shows that $j_! = j'^a_!$ is the desired left-adjoint. Also note that it follows from \cref{rslt:scheme-6-functor-cutoff-cardinals} that for any solid cutoff cardinal $\kappa$ we have $j_!(\Dqcohri(U, \Lambda)_\kappa) \subset \Dqcohri(X, \Lambda)_\kappa$.

Now assume that $j$ is not necessarily quasicompact anymore (but still separated). Let $\mathcal I_j$ be the category of quasicompact open subsets of $U$. Since $U$ is a perfectoid space, $\mathcal I_j$ is a basis of $\abs U$ and hence $\Dqcohri(U, \Lambda)_\kappa = \varprojlim_{V \in \mathcal I_j} \Dqcohri(V, \Lambda)_\kappa$, where the transition maps in the diagram on the right are the pullbacks. By the case of quasicompact $j$ all of these pullback functors preserve all small limits, which easily implies that limits in $\Dqcohri(U, \Lambda)_\kappa$ are computed componentwise with respect to the above limit (because the componentwise product preserves cartesian objects). Thus it follows from the above case of quasicompact $j$ that $j^*\colon \Dqcohri(X, \Lambda)_\kappa \to \Dqcohri(U, \Lambda)_\kappa$ preserves small limits and thus admits a left adjoint $j_{!\kappa}$ by the adjoint functor theorem. Using base-change and the case of quasicompact separated $j$, one deduces that for any cutoff cardinal $\kappa' \ge \kappa$ we have $j_{!\kappa'}(\Dqcohri(U, \Lambda)_\kappa) \subset \Dqcohri(X, \Lambda)_\kappa$, i.e. $j_{!\kappa}$ does not depend on $\kappa$. This shows the existence of $j_!$.

Now if $j$ is general (but $X$ still strictly totally disconnected) we can argue as in the separated case by arguing over the category $\mathcal J_j$ of open subsets of $U$ which are separated over $X$. This finishes the proof of all claims in the case that $X$ is strictly totally disconnected.

Finally, if $X$ is general we can argue as in \cite[Proposition 19.1]{etale-cohomology-of-diamonds}: Choose any v-hypercover $X_\bullet \to X$ such that each $X_n$ is a disjoint union of strictly totally disconnected spaces, and let $U_\bullet := X_\bullet \cprod_X U$, with $j_\bullet\colon U_\bullet \to X_\bullet$ the associated map of simplicial objects. Then $\Dqcohri(X, \Lambda)$ consists of the coCartesian sections $\Delta \to \Dqcohri(X_\bullet, \Lambda)$ (and similarly for $U$ in place of $X$). By the case of strictly totally disconnected $X$, the componentwise $j_{\bullet!}\colon \Dqcohri(U_\bullet, \Lambda) \to \Dqcohri(X_\bullet, \Lambda)$ preserves coCartesian morphisms (i.e. each $j_{n!}$ satisfies base-change) and hence restricts to a functor $j_!\colon \Dqcohri(U, \Lambda) \to \Dqcohri(X, \Lambda)$. Clearly $j_{\bullet!}$ is left adjoint to the componentwise $j_\bullet^*$, hence the same is true for their restrictions to cartesian objects, i.e. $j_!$ is left adjoint to $j_*$. Base-change for $j_!$ is formal.
\end{proof}

Let us also record the following lemma, which is of little interest by itself, but turns out to be very useful in some proofs below.

\begin{lemma} \label{rslt:etale-lower-shriek-of-generator-is-bounded}
Let $j\colon U \to X$ be a qcqs étale map in $\vStacksCoeff$. Then for every profinite set $S$ the sheaf $j_! \Lambda^a_\solid[S] \in \Dqcohri(X, \Lambda)$ is bounded.
\end{lemma}
\begin{proof}
Boundedness can be checked after pullback to any cover of $X$, so by base-change for $j_!$ we reduce to the case that $X$ is strictly totally disconnected. Then $U$ admits a finite open cover by separated perfectoid spaces, so we can assume that $U$ is separated. In this case $U$ is a finite disjoint union of quasicompact open subsets of $X$, so we can further reduce to the case that $U \subset X$ is a quasicompact open subset. Then as in the proof of \cref{rslt:existence-of-lower-shriek-for-etale-maps} we have $\Dqcohri(X, \Lambda) = \Dqcohri(\Lambda(X))$ and $\Dqcohri(U, \Lambda) = \Dqcohri(\Lambda(U)')$ for some classical ring $\Lambda(U)'$ such that $\Spec \Lambda(U)' \injto \Spec \Lambda(X)$ is an open immersion of schemes. By \cref{rslt:lower-shriek-commutes-with-almost-localization} $j_!$ commutes with the almost localization functor. But as a functor $\D_\solid(\Lambda(U)') \to \D_\solid(\Lambda(X))$ on concrete (i.e. non-almost) objects it preserves compact objects. This proves that $j_!\Lambda^a_\solid[S]$ is the almost localization of a compact object in $\D_\solid(X)$ and in particular bounded.
\end{proof}

Next up we show that the functor ``$f_! = g_* \comp j_!$'' satisfies the projection formula. As we have no good definition of $f_!$ yet (i.e. it is not clear yet that the given description is independent of the representation $f = g \comp j$) we will formulate the projection formula for $j_!$ and $g_*$ separately, which is the content of the next two results.

\begin{lemma} \label{rslt:projection-formula-for-open-immersion}
Let $j\colon U \to X$ be an étale map in $\vStacksCoeff$. Then for all $\mathcal M \in \Dqcohri(X, \Lambda)$ and $\mathcal N \in \Dqcohri(U, \Lambda)$ the natural morphism
\begin{align*}
	j_! (\mathcal N \tensor j^* \mathcal M) \isoto (j_! \mathcal N) \tensor \mathcal M
\end{align*}
is an isomorphism.
\end{lemma}
\begin{proof}
By base-change we can reduce to the case that $X = \Spa(A, A^+)$ is a strictly totally disconnected perfectoid space. In particular $U$ is a perfectoid space, so we can write $U = \bigunion_{i\in I} U_i$ as a filtered union of qcqs open subspaces. Then $\Dqcohri(U, \Lambda) = \varprojlim_i \Dqcohri(U_i, \Lambda)$ and $j_!$ is computed as $j_! = \varinjlim_i j_{i!}$, where $j_i$ is the map $j_i\colon U_i \to X$ (this follows formally from the adjunction of $j_!$ and $j^*$). As $\tensor$ preserves colimits in each argument, we can reduce to the case $U = U_i$ for some $i$, i.e. w.l.o.g. $U$ is qcqs. By a similar reduction we can further reduce to the case that $U$ is separated. Then $U$ is a finite disjoint union of quasicompact open subset of $X$, so we can even further reduce to the case that $U \subset X$ is an open subset. In particular $U = \Spa(B, B^+)$ is also totally disconnected and $\Dqcohri(X, \Lambda) = \Dqcohri(\Lambda(X))$ and $\Dqcohri(U, \Lambda) = \Dqcohri(\Lambda(U))$. By \cref{rslt:tot-disc-space-mod-pi-isom-to-scheme} the map $\Spec B^+/\pi \injto \Spec A^+/\pi$ is an open immersion of schemes, hence also $\Lambda(X)^a_\solid \to \Lambda(U)^a_\solid$ comes from an open immersion of schemes. By \cref{rslt:lower-shriek-commutes-with-almost-localization} the claim can therefore be reduced to the projection formula for an open immersion of schemes, which is \cref{rslt:scheme-6-functor-projection-formula}.
\end{proof}

\begin{lemma} \label{rslt:projection-formula-for-proper-p-bounded}
Let $f\colon Y \to X$ be a proper $p$-bounded map in $\vStacksCoeff$. Then for all $\mathcal M \in \Dqcohri(X, \Lambda)$ and $\mathcal N \in \Dqcohri(Y, \Lambda)$ the natural morphism
\begin{align*}
	(f_* \mathcal N) \tensor \mathcal M \isoto f_* (\mathcal N \tensor f^* \mathcal M)
\end{align*}
is an isomorphism.
\end{lemma}
\begin{proof}
By base-change for $f_*$ (see \cref{rslt:base-change-for-Dqcohri-along-qcqs-p-bounded-map}) we can assume that $X$ is a totally disconnected perfectoid space, in particular $\Dqcohri(X, \Lambda) = \Dqcohri(\Lambda(X))$. By \cref{rslt:qcqs-p-bounded-pushforward-preserves-colimits-and-Dqcohrim} both sides of the claimed isomorphism commute with filtered colimits in $\mathcal M$, so we can assume $\mathcal M = \Lambda(X)^a_\solid[S]$. By \cref{rslt:comparison-of-Dqcohri-and-DqcohriZ,rslt:p-bounded-implies-finite-cohom-dimension-on-DqcohriZ} there is a natural left complete $t$-structure on $\Dqcohri(Y, \Lambda)$ under which $f_*$ has finite cohomological dimension; thus both sides of the claimed isomorphism commute with Postnikov limits in $\mathcal N$, so that we can w.l.o.g. assume that $\mathcal N$ is left-bounded, say $\mathcal N \in \D_{\le0}$. Now pick a hypercover $Y_\bullet \to Y$ by totally disconnected spaces $Y_n$ of finite $\dimtrg$ over $X$, so that as in the proof of \cref{rslt:comparison-of-Dqcohri-and-DqcohriZ} we have
\begin{align*}
	\Dqcohri(Y, \Lambda) &= \varprojlim_{n\in\Delta} \Dqcohri(\Lambda(Y_n), \Lambda(X)),\\
	\Dqcohri(Y, (\Lambda, \Z)) &= \varprojlim_{n\in\Delta} \Dqcohri(\Lambda(Y_n), \Z).
\end{align*}
Then $\mathcal N \in \Dqcohrile0(Y, \Lambda)$ is representated by a cosimplicial object $N^\bullet$ in $\Dqcohri(\Lambda(X))$ such that $N^n \in \Dqcohrile0(\Lambda(Y_n), \Lambda(X))$ for all $n$. We deduce
\begin{align*}
	f_* (\mathcal N \tensor f^* \mathcal M) &= \Tot(N^\bullet \tensor_{(\Lambda(Y_\bullet), \Lambda(X))^a_\solid} (\Lambda(Y_\bullet), \Lambda(X))^a_\solid[S])\\
	&= \Tot\left(\varinjlim_{A \subset \Lambda(X)} N^\bullet \tensor_{(\Lambda(X), A)^a_\solid} (\Lambda(X), A)^a_\solid[S]\right),
\intertext{where $A$ ranges through all subrings of $\Lambda(X)$ which are of finite type over $\Z$. For each fixed $A$, the associated cosimplicial object on the right of the above equation naturally represents an object in $\Dqcohri(Y, (\Lambda, \Z))$ (this follows from the fact that the transition functors in the above limit representation for $\Dqcohri(Y, (\Lambda, \Z))$ are of the form $- \tensor_{\Lambda(Y_n)^a} \Lambda(Y_m)^a = - \tensor_{(\Lambda(Y_n), A)^a_\solid} (\Lambda(Y_m), A)^a_\solid$). Since $f_*$ commutes with filtered colimits by \cref{rslt:p-bounded-implies-finite-cohom-dimension-on-DqcohriZ}, the above expression becomes}
	&= \varinjlim_{A \subset \Lambda(X)} \Tot(N^\bullet \tensor_{(\Lambda(X), A)^a_\solid} (\Lambda(X), A)^a_\solid[S]).
\intertext{Now by \cref{rslt:finite-type-A+-implies-finite-flat-dimension-of-generators} the functor $- \tensor_{(\Lambda(X), A)^a_\solid} (\Lambda(X), A)^a_\solid[S]$ is left-bounded by some constant (depending on $A$) and thus commutes with totalizations of uniformly left-bounded cosimplicial objects. Therefore the above expression becomes}
	&= \varinjlim_{A \subset \Lambda(X)} \Tot(N^\bullet) \tensor_{(\Lambda(X), A)^a_\solid} (\Lambda(X), A)^a_\solid[S]\\
	&= \Tot(N^\bullet) \tensor_{\Lambda(X)^a_\solid} \Lambda(X)^a_\solid[S]\\
	&= (f_* \mathcal N) \tensor \mathcal M,
\end{align*}
as desired.
\end{proof}

Before we proceed with the definition of $f_!$, let us verify a different kind of projection formula for $f_!$, namely one with respect to a change of integral torsion coefficients $\Lambda \to \Lambda'$. Again, we do this separately for the étale and the proper case:

\begin{lemma} \label{rslt:coeff-projection-formula-for-etale-maps}
Let $j\colon U \to X$ be an étale map of small v-stacks and let $\Lambda \to \Lambda'$ be a morphism of integral torsion coefficients on $X$. Then for all $\mathcal N \in \Dqcohri(U, \Lambda)$ the natural morphism
\begin{align*}
	j_! (\mathcal N \tensor_{\Lambda^a_\solid} \Lambda'^a_\solid) \isoto (j_! \mathcal N) \tensor_{\Lambda^a_\solid} \Lambda'^a_\solid
\end{align*}
is an isomorphism.
\end{lemma}
\begin{proof}
By passing to right-adjoints, the claim reduces to showing that $j^*$ commutes with the forgetful functor $\Dqcohri(-, \Lambda') \to \Dqcohri(-, \Lambda)$. This is shown in \cref{rslt:Dqcohri-change-of-coefficients}.
\end{proof}

\begin{lemma} \label{rslt:coeff-projection-formula-for-proper-p-bounded}
Let $f\colon Y \to X$ be a proper $p$-bounded map of small v-stacks and let $\Lambda \to \Lambda'$ be a morphism of integral torsion coefficients on $X$. Then for all $\mathcal N \in \Dqcohri(Y, \Lambda)$ the natural morphism
\begin{align*}
	(f_* \mathcal N) \tensor_{\Lambda^a_\solid} \Lambda'^a_\solid \isoto f_* (\mathcal N \tensor_{\Lambda^a_\solid} \Lambda'^a_\solid)
\end{align*}
is an isomorphism.
\end{lemma}
\begin{proof}
By \cref{rslt:Dqcohri-change-of-coefficients,rslt:base-change-for-Dqcohri-along-qcqs-p-bounded-map} both sides of the claim satisfy base-change so we can reduce to the case that $X$ is totally disconnected. Note that the morphism $\Lambda \to \Lambda'$ is essentially given by its action $\Lambda(X) \to \Lambda'(X)$ on global sections: If the latter map is an isomorphism then also $\Lambda(Z)^a_\solid \isoto \Lambda'(Z)^a_\solid$ is an isomorphism for all $Z \in X_\vsite$. We can now factor $\Lambda(X) \to \Lambda'(X)$ as a composition $\Lambda(X) \to \Lambda''(X) \to \Lambda'(X)$, where $\Lambda''(X)$ is a polynomial algebra over $\Lambda(X)$ (with possibly infinitely many indeterminants) and $\Lambda''(X) \to \Lambda'(X)$ is finite (or even surjective). One can then define a new system of integral torsion coefficients $\Lambda''$ on $X$ by setting $\Lambda''(Z) := \Lambda(Z) \tensor_{\Lambda(X)} \Lambda''(X)$ together with the obvious maps $\Lambda \to \Lambda'' \to \Lambda'$. It is now enough to prove the claimed isomorphism for each of these two maps individually. For the second map this is easy: After applying the forgetful functors along $\Lambda'' \to \Lambda'$ (which is conservative), $- \tensor_{\Lambda^a_\solid} \Lambda'^a_\solid$ simply becomes $- \tensor \Lambda'^a$ in $\Dqcohri(-, \Lambda'')$, so the claimed isomorphism is a special case of the projection formula \cref{rslt:projection-formula-for-proper-p-bounded}.

It remains to handle the first map, i.e. from now on we can assume that $\Lambda'$ is a polynomial algebra over $\Lambda$ (with possibly infinitely many indeterminants). As in the proof of \cref{rslt:comparison-of-Dqcohri-and-DqcohriZ} we choose a hypercover $Y_\bullet \to Y$ by totally disconnected spaces which have finite $\dimtrg$ over $X$ in order to get
\begin{align*}
	\Dqcohri(Y, \Lambda) = \varprojlim_{n\in\Delta} \Dqcohri(\Lambda(Y_n), \Lambda(X))
\end{align*}
(and the same for $\Lambda'$ in place of $\Lambda$). This induces a left-complete $t$-structure on $\Dqcohri(Y, \Lambda)$ and both sides of the claimed isomorphism commute with Postnikov limits in $\mathcal N$, so we can assume that $\mathcal N$ is left-bounded. Thus $\mathcal N$ is represented by a uniformly left-bounded cosimplicial object $N^\bullet$ in $\Dqcohri(\Lambda(X))$ such that each $N^n$ comes from an object in $\Dqcohri(\Lambda(Y_n), \Lambda(X))$. We get
\begin{align*}
	(f_* \mathcal N) \tensor_{\Lambda^a_\solid} \Lambda'^a_\solid &= \Tot(N^\bullet) \tensor_{\Lambda(X)^a_\solid} \Lambda'(X)^a_\solid.
\intertext{It follows from \cref{rslt:finite-type-solid-pullback-preserves-limits} (and base-change) that $- \tensor_{\Lambda(X)^a_\solid} \Lambda'(X)^a_\solid$ is $t$-exact and in particular commutes with totalizations of uniformly left-bounded cosimplicial objects. Thus the above expression becomes}
	&= \Tot(N^\bullet \tensor_{\Lambda(X)^a_\solid} \Lambda'(X)^a_\solid)\\
	&= f_*(\mathcal N \tensor_{\Lambda^a_\solid} \Lambda'^a_\solid),
\end{align*}
as desired.
\end{proof}

A priori it is not clear that the above definition ``$f_! = g_* \comp j_!$'' is independent of the chosen representation $f = g \comp j$ (this problem can be avoided by using the canonical compactification) and that the so-defined $(-)_!$ is functorial. Both claims boil down to the following statement, which is analogous to \cite[Proposition 19.2]{etale-cohomology-of-diamonds}:

\begin{lemma} \label{rslt:etale-lower-shriek-compatible-with-pushforward}
Let
\begin{center}\begin{tikzcd}
	U' \arrow[r,"j'"] \arrow[d,"f'"] & X' \arrow[d,"f"]\\
	U \arrow[r,"j"] & X
\end{tikzcd}\end{center}
be a cartesian diagram in $\vStacksCoeff$. Assume that $j$ is an open immersion and $f$ is $p$-bounded and proper. Then the natural transformation
\begin{align*}
	j_! f'_* \isoto f_* j'_!
\end{align*}
of functors $\Dqcohri(U', \Lambda) \to \Dqcohri(X, \Lambda)$ is an equivalence.
\end{lemma}
\begin{proof}
This follows from the projection formula (see \cref{rslt:projection-formula-for-open-immersion,rslt:projection-formula-for-proper-p-bounded}) and base-change (see \cref{rslt:existence-of-lower-shriek-for-etale-maps,rslt:base-change-for-Dqcohri-along-qcqs-p-bounded-map}): Since $j$ is an open immersion, it follows easily from base-change for $j_!$ that $j_!$ and $j'_!$ are fully faithful. Thus, given any $\mathcal M \in \Dqcohri(U', \Lambda)$ we can write $\mathcal M = j'^* \mathcal N$ for some $\mathcal N \in \Dqcohri(X', \Lambda)$. Then
\begin{align*}
	&f_* j'_! \mathcal M = f_* j'_! (j^* \mathcal N \tensor \Lambda^a) = f_* (\mathcal N \tensor j'_! \Lambda^a) = f_* (\mathcal N \tensor f^* j_! \Lambda^a) = (f_* \mathcal N) \tensor j_!\Lambda^a =\\&\qquad= j_! j^* f_* \mathcal N = j_! f'_* \mathcal M,
\end{align*}
as desired.
\end{proof}

We now have developed all the necessary results to finally introduce the 6-functor formalism for $\DqcohriX X$. In order to properly introduce the 6-functor formalism, we need a precise definition of which morphisms $f\colon Y \to X$ will get an associated functor $f_!$:

\begin{definition} \label{def:bdcs-maps}
We say that a map $f\colon Y \to X$ of small v-stacks is \emph{bdcs}\footnote{The name comes from ``$p$-BoundeD, Compactifiable and Spatial''.} if it $p$-bounded, locally compactifiable (i.e. there is an analytic cover of $Y$ on which $f$ is compactifiable) and representable in locally spatial diamonds.
\end{definition}

\begin{lemma} \label{rslt:stability-of-bdcs-maps}
\begin{lemenum}
	\item The property of being bdcs is analytically local on both source and target.
	\item Bdcs maps are stable under composition and base-change.
	\item Every étale map is bdcs.
	\item Let $f\colon Y \to X$ and $g\colon Z \to Y$ be maps of small v-stacks. If $f$ and $f \comp g$ are bdcs, then so is $g$.
\end{lemenum}
\end{lemma}
\begin{proof}
Note that all of the claimed properties are true for ``$p$-bounded'' in place of ``bdcs'' by \cref{rslt:stability-of-p-bounded-maps}. We deduce (i), (ii) and (iii). To prove (iv) let $f$ and $g$ be given as in the claim. By passing to open covers of $Z$ and $Y$ we can assume that $g$ and $f \comp g$ are compactifiable. Then $g$ is compactifiable by \cite[Proposition 22.3]{etale-cohomology-of-diamonds}. It remains to show that $g$ is representable in locally spatial diamonds, so pick any v-cover $Y' \surjto Y$ by a locally spatial diamond $Y'$. Then we have a cartesian diagram
\begin{center}\begin{tikzcd}
	Z \cprod_Y Y' \arrow[r] \arrow[d] & Z \cprod_X Y' \arrow[d]\\
	Y \arrow[r] & Y \cprod_X Y
\end{tikzcd}\end{center}
of small v-stacks. Since $f \comp g$ is bdcs, $Z \cprod_X Y'$ is a locally separated locally spatial diamond. Since $g$ is separated, the diagonal $Y \to Y \cprod_X Y$ is a closed immersion and in particular a quasicompact injection. It follows from \cite[Proposition 11.20]{etale-cohomology-of-diamonds} that $Z \cprod_Y Y'$ is a locally spatial diamond, as desired.
\end{proof}

\begin{example}
Every map of rigid-analytic varieties over a fixed base field induces a bdcs map on the associated diamonds. Namely, by \cref{rslt:finite-type-implies-p-bounded} it is $p$-bounded and by \cref{rslt:stability-of-bdcs-maps} we can reduce to the case of affinoid spaces, where the map is automatically compactifiable.
\end{example}

We now have everything at hand to construct the desired 6-functor formalism for quasi-coherent $\ri^{+a}_X/\pi$-sheaves. By making use of our theory of abstract 6-functor formalisms (see \cref{sec:infcat.sixfun}) this construction is rather straightforward. In the following, recall the definition of the $\infty$-operad $\Corr(\mathcal C)_{E,all}$ and of 6-functor formalisms, see \cref{def:category-of-correspondences,def:correspondence-category-operad,def:abstract-6-functor-formalism}.

\begin{theorem} \label{rslt:main-6-functor-formalism}
There is a 6-functor formalism
\begin{align*}
	\Dqcohri\colon \Corr(\vStacksCoeff)_{bdcs,all} \to \infcatinf
\end{align*}
with the following properties:
\begin{thmenum}
	\item Restricted to the symmetric monoidal subcategory $(\vStacksCoeff)^\opp$ (equipped with the coproduct monoidal structure), $\Dqcohri$ coincides with the functor constructed in \cref{rslt:def-of-qcoh-Lambda-modules}.

	\item For every bdcs morphism $f\colon Y \to X$ in $\vStacksCoeff$, the functor
	\begin{align*}
		f_! := \Dqcohri([Y \xfrom{\id} Y \xto{f} X])\colon \Dqcohri(Y, \Lambda) \to \Dqcohri(X, \Lambda)
	\end{align*}
	preserves all small colimits. If $f = j$ is étale then $j_!$ is left adjoint to $j^*$ and if $f$ is proper then $f_! = f_*$.
\end{thmenum}
\end{theorem}
\begin{proof}
Fix an uncountable solid cutoff cardinal $\kappa$. We will first construct the 6-functor formalism $\Dqcohri(-)_\kappa$, where the main advantage is that all $\Dqcohri(X, \Lambda)_\kappa$ are presentable. At the end we will take the colimit over $\kappa$ to construct the desired 6-functor formalism $\Dqcohri$.

Let $P$ be the class of proper $p$-bounded maps and $I$ the class of quasicompact open immersions in $\vStacksCoeff$. Moreover, let $bdcqc$ denote the class of maps in $\vStacksCoeff$ which are $p$-bounded, compactifiable and quasicompact. Using \cref{rslt:p-bounded-maps-satisfy-2-out-of-3} it follows easily that the pair $I, P \subset bdcqc$ is a suitable decomposition of $bdcqc$ (see \cref{def:sixfun-suitable-decomposition-of-E}). We can thus apply \cref{rslt:construct-6-functor-formalism-out-of-I-P}: Condition (a) of that result is satisfied by \cref{rslt:base-change-for-Dqcohri-along-qcqs-p-bounded-map,rslt:projection-formula-for-proper-p-bounded}, condition (b) is satisfied by \cref{rslt:existence-of-lower-shriek-for-etale-maps,rslt:projection-formula-for-open-immersion} and condition (c) is satisfied by \cref{rslt:etale-lower-shriek-compatible-with-pushforward}. We obtain a 6-functor formalism
\begin{align*}
	\D_{qc,\kappa}\colon \Corr(\vStacksCoeff)_{bdcqc,all} \to \infcatinf, \qquad (X, \Lambda) \mapsto \Dqcohri(X,\Lambda)_\kappa.
\end{align*}
We will now extend $\D_{qc,\kappa}$ from $bdcqc$ to $bdcs$ using the extension results in \cref{sec:infcat.sixfun}. This will be performed in two steps:
\begin{enumerate}[1.]
	\item In the first step, we first restrict $\D_{qc,\kappa}$ to $\Corr(\vStacksCoeff)_{bdcsqc,all}$, where $bdcsqc \subset bdcqc$ is the subset of those edges which are representable in locally spatial diamonds. We now wish to extend $\D_{qc,\kappa}$ from $bdcsqc$ to the class of edges $bdcss$ consisting of those bdcs morphisms which are separated. By \cref{rslt:extend-6-functors-locally-on-target} (and \cref{rslt:sheaves-on-basis-equiv-sheaves-on-whole-site}) this extension can be performed on the full subcategory $\mathcal C_1 \subset \vStacksCoeff$ consisting of the separated locally spatial diamonds, i.e. we need to extend the $\infty$-operad map
	\begin{align*}
		\D_{qc,\kappa}\colon \Corr(\mathcal C_1)_{bdcsqc,all} \to \infcatinf
	\end{align*}
	to an $\infty$-operad map
	\begin{align*}
		\D_{s,\kappa}\colon \Corr(\mathcal C_1)_{bdcss,all} \to \infcatinf.
	\end{align*}
	We first apply \cref{rslt:extend-6-functor-to-disjoint-unions-on-source}, which allows us to extend $\D_{qc,\kappa}$ from $bdcsqc$ to the collection $E_1$ of edges of the form $\bigdunion_i Y_i \to X$, where each $Y_i \to X$ lies in $bdcsqc$; let us denote the new 6-functor formalism by $\D'_{qc,\kappa}$. We now apply \cref{rslt:extend-6-functors-locally-on-source} to extend $\D = \D'_{qc,\kappa}$ from $E_1$ to $E'_1 := bdcss$. Here we use the collection $S_1 \subset E_1$ of edges of the form $\bigdunion_i U_i \to X$ for covers $X = \bigunion_i U_i$ by quasicompact open immersions $U_i \injto X$. Then for $j \in S_1$ we have $j^! = j^*$, hence condition (b) of \cref{rslt:extend-6-functors-locally-on-source} follows from the sheafiness of $\D(-)_\kappa$. Condition (c) amounts to saying that every separated bdcs map $Y \to X$ of separated locally spatial diamonds admits a cover $Y = \bigunion_i V_i$ by quasicompact open immersions $V_i \injto Y$ such that each map $V_i \to X$ is quasicompact (it is automatically compactifiable by \cite[Proposition 22.3.(v)]{etale-cohomology-of-diamonds}); but this is easily satisfied, e.g. pick the $V_i$ to be any open cover of $Y$ by quasicompact open subsets (then the maps $V_i \injto Y$ and $V_i \to X$ are quasicompact because both $X$ and $Y$ are separated). Finally, condition (d) follows easily from the fact that all the spaces in $\mathcal C_1$ are separated. This finishes the construction of the 6-functor formalism
	\begin{align*}
		\D_{s,\kappa}\colon \Corr(\vStacksCoeff)_{bdcss,all} \to \infcatinf
	\end{align*}
	(where we implicitly used \cref{rslt:extend-6-functors-locally-on-target} to extend from $\mathcal C_1$ to $\vStacksCoeff$).

	\item In the second extension step we extend $\D_{s,\kappa}$ to the desired $\infty$-operad map
	\begin{align*}
		\D_\kappa\colon \Corr(\vStacksCoeff)_{bdcs,all} \to \infcatinf.
	\end{align*}
	This extension is similar to the previous one, albeit somewhat simpler: We can perform the extension directly on $\mathcal C_2 = \vStacksCoeff$ by applying \cref{rslt:extend-6-functors-locally-on-source} to $E_2 = bdcss$ and $E_2' = bdcs$ with $S_2 \subset E_2$ being the collection of all open immersions.
\end{enumerate}
The above construction of $\D_\kappa$ is functorial in $\kappa$ (for the extension results this is clear by construction; for the initial construction of $\D_{qc,\kappa}$ we note that the same construction works verbatim for $\Dqcohri(-)$ in place of $\Dqcohri(-)_\kappa$, which easily implies the desired functoriality), so that we can define
\begin{align*}
	\Dqcohri := \varinjlim_\kappa \D_\kappa\colon \Corr(\vStacksCoeff)_{bdcs,all} \to \infcatinf.
\end{align*}
By construction $\Dqcohri$ is a pre-6-functor formalism on $\vStacksCoeff$, which in particular assigns a functor $f_!\colon \Dqcohri(Y,\Lambda) \to \Dqcohri(X,\Lambda)$ to every bdcs map $f\colon Y \to X$ in $\vStacksCoeff$. To finish the proof, it remains to verify the following two properties of $\Dqcohri$:
\begin{enumerate}[(a)]
	\item For every étale map $j\colon U \to X$, $j_!$ is left adjoint to $j^*$. To see this, we first apply \cref{rslt:construct-6-functor-formalism-out-of-I-P} to the case that $E = I = et$ is the collection of étale maps in $\vStacksCoeff$ and $P$ consists only of degenerate edges; then conditions (b) and (c) are vacuous and condition (a) is satisfied by \cref{rslt:existence-of-lower-shriek-for-etale-maps,rslt:projection-formula-for-open-immersion}. We thus obtain a 6-functor formalism
	\begin{align*}
		\D_\et\colon \Corr(\vStacksCoeff)_{et,all} \to \infcatinf, \qquad (X,\Lambda) \mapsto \Dqcohri(X,\Lambda).
	\end{align*}
	We need to show that $\D_\et$ is equivalent to the restriction of $\Dqcohri$ to $\Corr(\vStacksCoeff)_{et,all}$. By the uniqueness of the extension results \cref{rslt:extend-6-functors-locally-on-source,rslt:extend-6-functor-to-disjoint-unions-on-source,rslt:extend-6-functors-locally-on-target} we can show this equivalence on the full subcategory $\mathcal C \subset \vStacksCoeff$ consisting of locally spatial diamonds and we can then restrict to the subset $etsqc \subset et$ of separated quasicompact étale maps. We can now further reduce to the full subcategory $\mathcal C' \subset \mathcal C$ consisting of strictly totally disconnected spaces. But note that every map of strictly totally disconnected spaces which lies in $etsqc$ is of the form $\bigdunion_{i=1}^n U_i \to X$ for quasicompact open immersions $U_i \injto X$, so we can further replace $etsqc$ by the collection of quasicompact open immersions. But in this case $\D_\et$ and $\Dqcohri$ agree by construction.

	\item For every bdcs map $f\colon Y \to X$ in $\vStacksCoeff$, the functor $f_!$ admits a right adjoint $f^!\colon \Dqcohri(X,\Lambda) \to \Dqcohri(Y,\Lambda)$. To see this, we first note that for a fixed uncountable solid cutoff cardinal $\kappa$ the functor $f_!$ maps $\Dqcohri(Y,\Lambda)_\kappa$ to $\Dqcohri(X,\Lambda)_\kappa$ and preserves all small colimits, so that it admits a right adjoint $f^!_\kappa\colon \Dqcohri(X,\Lambda)_\kappa \to \Dqcohri(Y,\Lambda)_\kappa$. Thus in order to prove the claim it is enough to show that $f^!_\kappa$ does not depend on $\kappa$ for large $\kappa$ (depending on $f$), i.e. if $\kappa'$ is a large enough cutoff cardinal then there is a cutoff cardinal $\kappa'' \ge \kappa'$ such that for a cofinal class of cutoff cardinal $\kappa \ge \kappa''$, $f^!_\kappa$ maps $\Dqcohri(-, \Lambda)_{\kappa'}$ to $\Dqcohri(-, \Lambda)_{\kappa''}$.

	Pick any hypercover $X_\bullet \to X$ by locally spatial diamonds $X_\bullet$ and let $f_\bullet\colon Y_\bullet \to X_\bullet$ be the pullback of $f$. In the same way as in the proof of \cref{rslt:existence-of-pushforward-on-Dqcohri} we can reduce to showing that $f_{n,\kappa}^!$ is independent of $\kappa$ for large $\kappa$ and all $n$, i.e. we can from now on assume that $X$ is a locally spatial diamond.
		Then $Y$ is also a locally spatial diamond (by definition of bdcs maps) and admits a cover by open subsets $U \subset Y$ which are quasicompact and compactifiable over $X$. Then $\Dqcohri(Y, \Lambda)_\kappa = \varprojlim_U \Dqcohri(U, \Lambda)_\kappa$, where $U$ ranges through the aforementioned open subsets of $Y$, and since $j^!_\kappa = j^*$ for open immersions $j$ the functor $f^!_\kappa$ is computed componentwise with respect to this limit of $\infty$-categories. We can therefore reduce to the case $U = Y$ for some $U$, i.e. we can from now on assume that $f$ is quasicompact and separated. By passing to a hypercover of $X$ consisting of disjoint unions of totally disconnected spaces, we can further assume that $X$ is totally disconnected. Then $Y$ is a spatial diamond, so we can pick a pro-étale cover $U = \varprojlim_i U_i \to Y$ such that $U$ is totally disconnected and all $U_i \to Y$ are quasicompact separated étale (see \cite[Proposition 11.24]{etale-cohomology-of-diamonds}). Denote $g\colon U \to Y$ and $g_i\colon U_i \to Y$ the projections. It is enough to show that $g^* f^!_\kappa$ is independent of $\kappa$ (i.e. preserves $\Dqcohri(-, \Lambda)_{\kappa'}$ with $\kappa' \le \kappa$ as above). Since $f_* g_*\colon \Dqcohri(\Lambda(U)) \to \Dqcohri(\Lambda(X))$ is just the forgetful functor, it is then enough to verify that $f_* g_* g^* f^!_\kappa$ is independent of $\kappa$. By \cref{rslt:colim-of-pushforward-pullback} we have $g_* g^* = \varinjlim_i g_{i*} g_i^*$ and by \cref{rslt:qcqs-p-bounded-pushforward-preserves-colimits-and-Dqcohrim} $f_*$ preserves this colimit, so since $\Dqcohri(-, \Lambda)_\kappa \subset \Dqcohri(-, \Lambda)$ is stable under all small colimits, we reduce to showing that $f_* g_{i*} g_i^* f^!_\kappa$ is independent of $\kappa$. Replacing $f$ by $f \comp g_i$ we are now reduced to showing that $f_* f^!_\kappa$ is independent of $\kappa$. It follows formally from the projection formula for $f_!$ (which holds by \cref{rslt:abstract-6-functor-formalism-properties}) and adjunctions that $f_* f^!_\kappa = \IHom(f_! \Lambda^a, -)_\kappa$. But $f_!\Lambda^a$ is just some fixed object in $\Dqcohri(\Lambda(X))$, hence $\IHom(f_! \Lambda^a, -)_\kappa$ does indeed not depend on $\kappa$ for large enough $\kappa$ (see the proof of the existence of $\IHom$ in \cref{rslt:condensed-objects-in-presentable-monoidal-cat}).
\end{enumerate}
This finishes the construction of the 6-functor formalism $\Dqcohri$ with all the desired properties.
\end{proof}

From \cref{rslt:main-6-functor-formalism} we can extract the following definition of the shriek functors, thereby completing our collection of six functors:

\begin{definition} \label{def:shriek-functors}
Let $f\colon Y \to X$ be a bdcs map in $\vStacksCoeff$.
\begin{defenum}
	\item We define $f_!\colon \Dqcohri(Y, \Lambda) \to \Dqcohri(X, \Lambda)$ to be the functor $f_! := \Dqcohri(Y \xfrom{\id} Y \xto{f} X)$, where $\Dqcohri$ is as in \cref{rslt:main-6-functor-formalism}.

	\item We define $f^!\colon \Dqcohri(X, \Lambda) \to \Dqcohri(Y, \Lambda)$ to be the right adjoint of $f_!$, which exists by \cref{rslt:main-6-functor-formalism}.
\end{defenum}
\end{definition}

\begin{remark}
The construction of the functor $f_!$ in \cref{rslt:main-6-functor-formalism} is not very explicit, so we provide a more direct description:
\begin{enumerate}[1.]
	\item We first define $f_! := (\overline f^{/X})_* \comp j_!$ for every quasicompact compactifiable $p$-bounded map $f\colon Y \to X$, which automatically factors as the composition of the open immersion $j\colon Y \injto \overline Y^{/X}$ and the proper $p$-bounded map $\overline f^{/X}\colon \overline Y^{/X} \to X$.

	\item Suppose that $f\colon Y \to X$ is a bdcs map of locally spatial diamonds. Let $\mathcal I$ be the category of open subsets $V \subset Y$ which are quasicompact and compactifiable over $X$. Then for every solid cutoff cardinal $\kappa$ (see \cref{def:solid-cutoff-cardinal}) we have $\Dqcohri(Y, \Lambda)_\kappa = \varinjlim_{V \in \mathcal I} \Dqcohri(V, \Lambda)_\kappa$ in $\catPrL$, where the transition maps on the right are given by $j_!$ for open immersions $j\colon V \injto V'$. We thus need to define $f_!$ to be the unique morphism in $\catPrL$ which restricts to the functor $g_!$ defined above for every composition $g\colon V \injto Y \to X$ with $V \in \mathcal I$.

	\item Suppose that $f\colon Y \to X$ is a bdcs map of small v-stacks. Choose a hypercover $X_\bullet \to X$ such that all $X_n$ are locally spatial diamonds, and let $f_\bullet\colon Y_\bullet \to X_\bullet$ be the base-change of $f$. Then $f_{\bullet!}\colon \Dqcohri(Y_\bullet, \Lambda) \to \Dqcohri(X_\bullet, \Lambda)$ preserves coCartesian edges and hence restricts to a functor $f_!\colon \Dqcohri(Y, \Lambda) \to \Dqcohri(X, \Lambda)$. Note that $f_!$ must necessarily be defined this way if we want it to satisfy base-change.
\end{enumerate}
Note that instead of constructing the 6-functor formalism in \cref{rslt:main-6-functor-formalism} one could carry out the above 3 steps directly, cf. \cite[\S22]{etale-cohomology-of-diamonds}. However, this way we do not get all the higher coherences we expect $f_!$ to satisfy.
\end{remark}

For the convenience of the reader we extract the standard 1-categorical properties of the 6-functor formalism in \cref{rslt:main-6-functor-formalism}:
\begin{corollary} \label{rslt:main-6-functor-formalism-properties}
The adjoint pairs of functors $(\tensor, \IHom)$, $(f^*, f_*)$, $(f_!, f^!)$ satisfy the following properties:
\begin{corenum}
	\item (Functoriality) The assignments $f \mapsto f^*$, $f \mapsto f_*$, $f \mapsto f_!$ and $f \mapsto f^!$ are functorial in $f$, i.e. they define functors $(\vStacksCoeff)^\opp \to \infcatinf$, $\vStacksCoeff \to \infcatinf$, $(\vStacksCoeff)_{bdcs} \to \infcatinf$ and $((\vStacksCoeff)_{bdcs})^\opp \to \infcatinf$.

	\item (Special Cases) If $j\colon U \to X$ is an étale map in $\vStacksCoeff$ then $j^! = j^*$. If $f\colon Y \to X$ is a proper bdcs map in $\vStacksCoeff$ then $f_! = f_*$.

	\item \label{rslt:main-6-functor-formalism-projection-formula} (Projection Formula) Let $f\colon Y \to X$ be a bdcs map in $\vStacksCoeff$. Then for all $\mathcal M \in \Dqcohri(X, \Lambda)$ and $\mathcal N \in \Dqcohri(Y, \Lambda)$ there is a natural isomorphism
	\begin{align*}
		f_!(\mathcal N \tensor f^* \mathcal M) = (f_! \mathcal N) \tensor \mathcal M.
	\end{align*}

	\item \label{rslt:main-6-functor-formalism-proper-base-change} (Proper Base-Change) Let
	\begin{center}\begin{tikzcd}
		Y' \arrow[r,"g'"] \arrow[d,"f'"] & Y \arrow[d,"f"]\\
		X' \arrow[r,"g"] & X
	\end{tikzcd}\end{center}
	be a cartesian diagram in $\vStacksCoeff$ such that $f$ is bdcs. Then there is a natural equivalence
	\begin{align*}
		g^* f_! = f'_! g'^*
	\end{align*}
	of functors $\Dqcohri(Y, \Lambda) \to \Dqcohri(X', \Lambda)$.
\end{corenum}
\end{corollary}
\begin{proof}
Every 6-functor formalism satisfies these properties (see \cref{rslt:abstract-6-functor-formalism-properties}), so in particular they are satisfied by the 6-functor formalism of \cref{rslt:main-6-functor-formalism}.
\end{proof}

The 6-functor formalism from \cref{rslt:main-6-functor-formalism} is functorial in the chosen integral torsion coefficients $\Lambda$, in the following sense:

\begin{proposition} \label{rslt:coeff-projection-formula}
Let $f\colon Y \to X$ be a bdcs map of small v-stacks and let $\Lambda \to \Lambda'$ be a morphism of integral torsion coefficients on $X$. Then for all $\mathcal N \in \Dqcohri(Y, \Lambda)$ there is a natural isomorphism
\begin{align*}
	f_! (\mathcal N \tensor_{\Lambda^a_\solid} \Lambda'^a_\solid) = (f_! \mathcal N) \tensor_{\Lambda^a_\solid} \Lambda'^a_\solid.
\end{align*}
\end{proposition}
\begin{proof}
Let $\vStacksCoeffVar$ be the category of pairs $(X,\Lambda) \in \vStacksCoeff$ but where morphisms $(Y,\Lambda') \to (X,\Lambda)$ consist of a morphism $f\colon Y \to X$ of small v-stacks together with a morphism $f^*\Lambda \to \Lambda'$ of integral torsion coefficients. We let $bdcs$ denote the collection of morphisms $f\colon (Y,\Lambda') \to (X,\Lambda)$ in $\vStacksCoeffVar$ such that $f\colon Y \to X$ is bdcs and $f^*\Lambda \isoto \Lambda'$ is an isomorphism. Then the same arguments as in \cref{rslt:main-6-functor-formalism} allow us to construct a 6-functor formalism
\begin{align*}
	\Dqcohri\colon \Corr(\vStacksCoeffVar)_{bdcs,all} \to \infcatinf, \qquad (X,\Lambda) \mapsto \Dqcohri(X,\Lambda).
\end{align*}
In fact the proof works verbatim, except that in the first step (i.e. in the application of \cref{rslt:construct-6-functor-formalism-out-of-I-P}) we need to additionally use \cref{rslt:coeff-projection-formula-for-etale-maps,rslt:coeff-projection-formula-for-proper-p-bounded}. Having this 6-functor formalism on $\vStacksCoeffVar$ at hand, the desired projection formula is just proper base-change along the diagram
\begin{center}\begin{tikzcd}
	(Y,\Lambda') \arrow[r] \arrow[d] & (Y,\Lambda) \arrow[d]\\
	(X,\Lambda') \arrow[r] & (X,\Lambda)
\end{tikzcd}\end{center}
in $\vStacksCoeffVar$, which holds by \cref{rslt:abstract-6-functor-formalism-properties}.
\end{proof}

We also record the following formal consequences of the 6-functor formalism, showing some more compatibilities of the six functors.

\begin{corollary}
Let $f\colon Y \to X$ be a bdcs map in $\vStacksCoeff$.
\begin{corenum}
	\item \label{rslt:IHom-adjunction-for-shriek-functors} For all $\mathcal M \in \Dqcohri(X, \Lambda)$ and $\mathcal N \in \Dqcohri(Y, \Lambda)$ there is a natural isomorphism
	\begin{align*}
		\IHom(f_! \mathcal N, \mathcal M) = f_* \IHom(\mathcal N, f^! \mathcal M).
	\end{align*}

	\item \label{rslt:upper-shriek-of-IHom-formula} For all $\mathcal M, \mathcal N \in \Dqcohri(X, \Lambda)$ there is a natural isomorphism
	\begin{align*}
		f^! \IHom(\mathcal N, \mathcal M) \isom \IHom(f^* \mathcal N, f^! \mathcal M).
	\end{align*}
\end{corenum}
\end{corollary}
\begin{proof}
Both properties follow formally from adjunctions and the projection formula, cf. \cite[Proposition 23.3]{etale-cohomology-of-diamonds}.
\end{proof}

\subsection{Compact, Perfect and Dualizable Sheaves} \label{sec:ri-pi.compact-dualizable}

We now focus our attention on some important special types of objects in $\Dqcohri(X, \Lambda)$, namely (almost) compact, perfect and dualizable ones.

Almost compact objects can be defined and studied in the general context of $V^a$-enriched $\infty$-categories for some almost setup $(V, \mm)$ (see \cref{def:almost-compact-object} and \cref{rslt:almost-compact-generators-implies-generated}). We introduce the following terminology:

\begin{definition}
Let $(V, \mm)$ be an almost setup, $\Lambda$ a discrete $V$-algebra and $P \in \D^a_\solid(\Lambda)$.
\begin{defenum}
	\item $P$ is called \emph{almost compact} if it is so as an object of the $V^a$-enriched $\infty$-category $\D^a_\solid(\Lambda)$, equivalently if the functor $\IHom(P,-)$ preserves all small colimits. We denote $\D^a_\solid(\Lambda)^{a\omega} \subset \D^a_\solid(\Lambda)$ the full subcategory of almost compact objects.

	\item $P$ is called \emph{almost perfect} if for every $\varepsilon \in \mm$ there exists a perfect discrete $\Lambda$-module $P_\varepsilon \in \D(\Lambda)_\omega$ such that $P \approx_\varepsilon P^a_\varepsilon$.

	\item $P$ is called \emph{weakly almost perfect} if for every $\varepsilon \in \mm$ there is a perfect discrete $\Lambda$-module $P_\varepsilon \in \D(\Lambda)_\omega$ such that $P$ is an $\varepsilon$-retract of $P^a_\varepsilon$. We denote $\D^a_\solid(\Lambda)_\waperf \subset \D^a_\solid(\Lambda)$ the full subcategory of weakly almost perfect objects.
\end{defenum}
\end{definition}

From the general study of almost compact objects in $V^a$-enriched $\infty$-categories we immediately obtain the following result, showing that almost compact objects in $\D^a_\solid(\Lambda)$ can be approximated by actual compact objects, similar to weakly almost perfect objects:

\begin{lemma} \label{rslt:approximate-almost-compact-by-compact}
Let $(V, \mm)$ be an almost setup, $\Lambda$ a discrete $V$-algebra and $P \in \D^a_\solid(\Lambda)$. Then the following are equivalent:
\begin{lemenum}
	\item $P$ is almost compact.
	\item For every $\varepsilon \in \mm$ there is a compact object $P_\varepsilon \in \D_\solid(\Lambda)$ such that $P$ is an $\varepsilon$-retract of $P^a_\varepsilon$.
\end{lemenum}
Moreover, if $P$ is discrete (and almost compact) then the $P_\varepsilon$ can be chosen discrete as well.
\end{lemma}
\begin{proof}
This follows directly from \cref{rslt:characterization-of-almost-compact-objects-in-general-C}.
\end{proof}

Apart from compact and perfect objects, the third notion we want to study is that of dualizable objects in $\D^a_\solid(\Lambda)$. Dualizability can be defined in a much more general context, as follows:

\begin{definition}
Let $\mathcal C$ be a symmetric monoidal $\infty$-category. An object $P \in \mathcal C$ is called \emph{dualizable} if there are an object $P^*$, called the \emph{dual} of $P$, and morphisms
\begin{align*}
	\ev_P\colon P^* \tensor P \to 1, \qquad i_P\colon 1 \to P \tensor P^*,
\end{align*}
called the \emph{evaluation map} and \emph{coevaluation map} respectively, such that there are homotopy coherent diagrams
\begin{center}
	\begin{tikzcd}
		P \arrow[dr,"\id"] \arrow[r,"i_P \tensor \id"] & P \tensor P^* \tensor P \arrow[d,"\id \tensor \ev_P"] \\
		& P
	\end{tikzcd}
	\qquad
	\begin{tikzcd}
		P^* \arrow[dr,"\id"] \arrow[r,"\id \tensor i_P"] & P^* \tensor P \tensor P^* \arrow[d,"\ev_P \tensor \id"] \\
		& P^*
	\end{tikzcd}
\end{center}
\end{definition}

Let us recall a few equivalent definitions for dualizable objects in case the symmetric monoidal $\infty$-category $\mathcal C$ is closed:

\begin{lemma} \label{rslt:equivalent-defs-for-dualizable}
Let $\mathcal C$ be a closed symmetric monoidal $\infty$-category, pick any $P \in \mathcal C$ and let $P^* := \IHom(P, -)$. Then the following are equivalent:
\begin{lemenum}
	\item $P$ is dualizable.
	\item The natural map $P \tensor P^* \isoto \IHom(P, P)$ is an isomorphism.
	\item For every $X \in \mathcal C$ the natural map $X \tensor P^* \isoto \IHom(P, X)$ is an isomorphism.
	\item The functor $- \tensor P\colon \mathcal C \to \mathcal C$ is left adjoint to the functor $- \tensor P^*\colon \mathcal C \to \mathcal C$.
\end{lemenum}
In this case $P^*$ is the dual of $P$ and $P^*$ is dualizable with dual $P$.
\end{lemma}
\begin{proof}
Assume first that $\mathcal C$ is a $1$-category. It is clear that (iii) implies (ii). To check that (i) implies (iv) we assume that $P$ is dualizable with dual $P'$. Then one checks by hand that the evaluation and coevaluation maps provide an adjunction of the functors $- \tensor P$ and $- \tensor P'$. This easily implies $P' = P^*$, proving (iv). To see that (iv) implies (iii), we simply note that under the assumption of (iv) we have for all $X, Y \in \mathcal C$
\begin{align*}
	\Hom(Y, X \tensor P^*) = \Hom(Y \tensor P, X) = \Hom(Y, \IHom(P, X)).
\end{align*}
We now prove that (ii) implies (i), so assume that $P$ satisfies (ii). There is a natural evaluation map $\ev_P\colon P^* \tensor P \to 1$ and a coevaluation map $i_P\colon 1 \to \IHom(P, P) = P \tensor P^*$ and one checks easily by hand that they satisfy the required commutative diagrams. This finishes the proof in the case that $\mathcal C$ is a $1$-category.

Now let $\mathcal C$ be general. Note that (i), (ii) and (iii) only depend on the homotopy category of $\mathcal C$, which is a symmetric monoidal $1$-category. Thus their equivalence follows from the $1$-categorical case discussed above. Moreover, it is easy to see that (iii) and (iv) are equivalent.
\end{proof}

The following result relates all types of modules introduced above, showing that they are essentially all the same:

\begin{proposition} \label{rslt:equiv-of-compact-perfect-dualizable-over-Lambda}
Let $(V, \mm)$ be an almost setup, $\Lambda$ a discrete $V$-algebra and $P \in \D^a_\solid(\Lambda)$. Then the following are equivalent:
\begin{propenum}
	\item $P$ is almost compact and discrete.
		\item $P$ is weakly almost perfect.
	\item $P$ is dualizable.
\end{propenum}
If this is the case then $P$ is weakly almost bounded.
\end{proposition}
\begin{proof}
In the $\infty$-category $\D(\Lambda)_\omega$ of discrete $\Lambda$-modules, the notions of compact and perfect objects are equivalent, see e.g. \cite[Proposition 07LT]{stacks-project}. It thus follows from \cref{rslt:approximate-almost-compact-by-compact} that (i) implies (ii). We now prove that (ii) implies (i), so assume that $P$ is weakly almost perfect. By \cref{rslt:approximate-almost-compact-by-compact} $P$ is almost compact, so it only remains to see that $P$ is discrete. Let $P_\omega$ be the discretization of $P$ (i.e. the image under the right-adjoint of $\D^a(\Lambda)_\omega \injto \D^a_\solid(\Lambda)$), with adjunction map $\delta\colon P_\omega \to P$. We claim that $\delta$ is an isomorphism. Fix some $\varepsilon \in \mm$ and choose a perfect object $P_\varepsilon \in \D(\Lambda)_\omega$ and maps $f_\varepsilon\colon P \to P^a_\varepsilon$, $g_\varepsilon\colon P^a_\varepsilon \to P$ such that $g_\varepsilon f_\varepsilon \isom \varepsilon \id$. Then $P^a_\varepsilon$ is discrete, i.e. the associated map $(P^a_\varepsilon)_\omega \isoto P^a_\varepsilon$ from its discretization is an isomorphism. A quick diagram chase (similar to \cref{rslt:characterization-of-almost-compact-objects-in-general-C}) reveals that on homotopy groups both kernel and cokernel of $\delta$ are killed by $\varepsilon$. As this is true for every $\varepsilon$, we deduce that $\delta$ is an isomorphism. This finishes the proof that (i) and (ii) are equivalent.

We now prove that (iii) implies (i), so assume that $P$ is dualizable. Then by \cref{rslt:equivalent-defs-for-dualizable} we have $\IHom(P, -) = P^* \tensor -$, which clearly preserves all small colimits; hence $P$ is almost compact. To show that $P$ is discrete, we can equivalently show that $P$ is nuclear by \cref{rslt:solid-discrete-equiv-nuclear}. But by the adjunction of $P \tensor -$ and $P^* \tensor -$ (see \cref{rslt:equivalent-defs-for-dualizable}) we have for all profinite sets $S$
\begin{align*}
	&P(S) = \Hom(\Lambda_\solid^a[S], P) = \Hom(\Lambda_\solid^a[S] \tensor P^*, \Lambda) = \Hom(P^*, \Lambda^a_\solid[S]^\vee) =\\&\qquad= \Hom(\Lambda, \Lambda^a_\solid[S]^\vee \tensor P) = (\Lambda^a_\solid[S]^\vee \tensor P)(*),
\end{align*}
as desired (this does not exactly show nuclearity of $P$, but is enough to apply the proof of \cref{rslt:solid-discrete-equiv-nuclear}). We now show that (i) implies (iii), so assume that $P$ is almost compact and discrete. By \cref{rslt:equivalent-defs-for-dualizable} we need to see that the natural map $\IHom(P, \Lambda^a) \tensor P \isoto \IHom(P, P)$ is an isomorphism. For this it is enough to show that for every discrete $M \in \D^a(\Lambda)_\omega$ we the map $\IHom(P, \Lambda^a) \tensor Q \isoto \IHom(P, Q)$ is an isomorphism. As $Q$ is discrete we can write it as a colimit of copies of $\Lambda^a$. But by almost compactness of $P$, both sides of the claimed isomorphism commute with colimits in $Q$, so w.l.o.g. $Q = \Lambda^a$. In this case the claimed isomorphism is obvious.

Finally, it follows immediately from \cref{rslt:approximate-almost-compact-by-compact} that if $P$ is almost compact then it is weakly almost bounded (because every compact object in $\D_\solid(\Lambda)$ is bounded). This proves the last part of the claim.
\end{proof}

Let us now extend the above definitions to the geometric setting, i.e. to objects in $\Dqcohri(X, \Lambda)$ for $(X, \Lambda) \in \vStacksCoeff$. As this $\infty$-category is symmetric monoidal, we have a notion of dualizable objects inside it. However, we need an extra definition for (almost) compact and perfect objects:

\begin{definition}
Let $(X, \Lambda) \in \vStacksCoeff$.
\begin{defenum}
	\item Suppose $X \in X_\vsite^\Lambda$. Then $\Dqcohri(X, \Lambda)$ is naturally $\Lambda(X)^a$-enriched via the symmetric monoidal functor $\widetilde{(-)}\colon \D^a_{\ge0}(\Lambda(X))_\omega \to \Dqcohri(X, \Lambda)$. More concretely, under this enrichment we have
	\begin{align*}
		\Hom(\mathcal M, \mathcal N)^a = \tau_{\ge0} \Gamma(X, \IHom(\mathcal M, \mathcal N))_\omega \in \D^a_{\ge0}(\Lambda(X))_\omega
	\end{align*}
	for any two objects $\mathcal M, \mathcal N \in \Dqcohri(X, \Lambda)$. An object $\mathcal P \in \Dqcohri(X, \Lambda)$ is called \emph{almost compact} if it is almost compact with respect to the enrichment just discussed, i.e. if the functor $\Hom(\mathcal P, -)^a$ preserves all small colimits. We denote $\Dqcohri(X, \Lambda)^{a\omega} \subset \Dqcohri(X, \Lambda)$ the full subcategory of almost compact objects.

	\item \label{def:weakly-almost-perfect-module-on-v-stack} An object $\mathcal P \in \Dqcohri(X, \Lambda)$ is called \emph{weakly almost perfect} if for every map $f\colon Y \to X$ from a totally disconnected space $Y$ the pullback $f^* \mathcal P \in \Dqcohri(\Lambda(Y))$ is weakly almost perfect. We denote $\Dqcohri(X, \Lambda)_\waperf \subset \Dqcohri(X, \Lambda)$ the full subcategory of weakly almost perfect objects.
\end{defenum}
\end{definition}

We start with a few basic results on almost compact and weakly almost perfect objects in $\Dqcohri(X, \Lambda)$.

\begin{lemma}
Let $X \in \vStacksCoeff$ with $X \in X_\vsite^\Lambda$.
\begin{lemenum}
	\item \label{rslt:almost-hom-adjunction-for-*} For every map $f\colon Y \to X$ of small v-stacks and all $\mathcal M \in \Dqcohri(X, \Lambda)$ and $\mathcal N \in \Dqcohri(Y, \Lambda)$ there is a natural equivalence
	\begin{align*}
		\Hom(f^* \mathcal M, \mathcal N)^a = \Hom(\mathcal M, f_* \mathcal N)^a
	\end{align*}
	in $\D^a_{\ge0}(\Lambda(X))_\omega$.

	\item \label{rslt:almost-hom-adjunction-for-etale-!} For every bdcs map $f\colon Y \to X$ and all $\mathcal M \in \Dqcohri(Y, \Lambda)$ and $\mathcal N \in \Dqcohri(X, \Lambda)$ there is a natural equivalence
	\begin{align*}
		\Hom(f_! \mathcal M, \mathcal N)^a = \Hom(\mathcal M, f^! \mathcal N)^a
	\end{align*}
	in $\D^a_{\ge0}(\Lambda(X))_\omega$.
\end{lemenum}
\end{lemma}
\begin{proof}
We first prove (i). Using that $\Gamma(Y, -) = \Gamma(X, f_*(-))$ (up to a forgetful functor) one easily constructs a natural map from right to left of the claimed identity. Using Yoneda and adjunctions it is easy to check that this map is an isomorphism.

Part (ii) follows in a similar way from adjunctions (and the projection formula, see \cref{rslt:main-6-functor-formalism-projection-formula}).
\end{proof}

\begin{corollary}
Let $X \in \vStacksCoeff$ with $X \in X_\vsite^\Lambda$.
\begin{corenum}
	\item \label{rslt:pullback-of-almost-compact-along-qcqs-pbd-is-almost-compact} Let $f\colon Y \to X$ be a qcqs $p$-bounded map. If $\mathcal P \in \Dqcohri(X, \Lambda)$ is almost compact then so is $f^*\mathcal P \in \Dqcohri(Y, \Lambda)$.

	\item \label{rslt:etale-lower-shriek-of-almost-compact-is-almost-compact} Let $j\colon U \to X$ be an étale map. If $\mathcal P \in \Dqcohri(U, \Lambda)$ is almost compact then so is $j_! \mathcal P \in \Dqcohri(X, \Lambda)$.
\end{corenum}
\end{corollary}
\begin{proof}
Part (i) follows easily from \cref{rslt:almost-hom-adjunction-for-*}: By \cref{rslt:qcqs-p-bounded-pushforward-preserves-colimits-and-Dqcohrim} the pushforward $f_*$ preserves all small colimits, so if $\mathcal P$ is almost compact then $\Hom(f^* \mathcal P, -)^a = \Hom(\mathcal P, f_*(-))^a$ preserves all small colimits, as desired. In a similar vein, part (ii) follows directly from \cref{rslt:almost-hom-adjunction-for-etale-!}
\end{proof}

\begin{lemma} \label{rslt:geometric-properties-of-dualizable-and-perfect-sheaves}
Let $(X, \Lambda) \in \vStacksCoeff$.
\begin{lemenum}
	\item \label{rslt:weakly-alm-perf-sheaves-satisfy-v-descent} The assignment $(X, \Lambda) \mapsto \Dqcohri(X, \Lambda)_\waperf$ defines a hypercomplete v-sheaf on $\vStacksCoeff$.

	\item \label{rslt:weakly-alm-perf-equiv-dualizable} An object $\mathcal P \in \Dqcohri(X, \Lambda)$ is weakly almost perfect if and only if it is dualizable. In this case, for every map $f\colon Y \to X$ of small v-stacks and every $\mathcal M \in \Dqcohri(X, \Lambda)$ there is a natural isomorphism
	\begin{align*}
		f^* \IHom(\mathcal P, \mathcal M) = \IHom(f^* \mathcal P, f^* \mathcal M).
	\end{align*}

	\item \label{rstl:discrete-compact-implies-perfect-on-loc-spatial-diam} Suppose $X$ is a locally spatial diamond. Then every discrete almost compact object in $\Dqcohri(X, \Lambda)$ is weakly almost perfect.
\end{lemenum}
\end{lemma}
\begin{proof}
Suppose $\mathcal P \in \Dqcohri(X, \Lambda)$ is dualizable and let $f\colon Y \to X$ be any map of small v-stacks. Then $f^*\colon \Dqcohri(X, \Lambda) \to \Dqcohri(Y, \Lambda)$ is symmetric monoidal, hence it follows directly from the definition that $f^* \mathcal P \in \Dqcohri(Y, \Lambda)$ is also dualizable with dual $(f^* P)^* = f^* P^*$. In particular, for every $\mathcal M \in \Dqcohri(X, \Lambda)$ we have by \cref{rslt:equivalent-defs-for-dualizable}
\begin{align*}
	f^* \IHom(\mathcal P, \mathcal M) = f^* (\mathcal P^* \tensor \mathcal M) = (f^* \mathcal P)^* \tensor f^* \mathcal M = \IHom(f^* \mathcal P, f^* \mathcal M).
\end{align*}
This proves the second part of (ii). Now drop the assumption that $\mathcal P$ is dualizable and instead assume that $f^* \mathcal P$ is dualizable and $f$ is a cover; we will deduce that $\mathcal P$ is dualizable. Let $Y_\bullet \to X$ be a hypercover with $Y_0 = Y$ and let $f_n\colon Y_n \to X$ denote the projection. Then $\Dqcohri(X, \Lambda) = \varprojlim_{n\in\Delta} \Dqcohri(Y_n, \Lambda)$ and for every $\mathcal M \in \Dqcohri(X, \Lambda)$ we have
\begin{align*}
	\IHom(\mathcal P, \mathcal M) = \varprojlim_{n\in\Delta} f_{n*} \IHom(f_n^*\mathcal P, f_n^* \mathcal M).
\end{align*}
On the other hand, by what we have shown above the object $\IHom(f_\bullet^* \mathcal P, f_\bullet^* \mathcal M) \in \Dqcohri(Y_\bullet, \Lambda)$ is coCartesian, which implies $f_\bullet^* \IHom(\mathcal P, \mathcal M) = \IHom(f_\bullet^* \mathcal P, f_\bullet^* \mathcal M)$. This in particular implies that we have $f^*\IHom(\mathcal P, \mathcal M) = \IHom(f^* \mathcal P, f^* \mathcal M)$. Now it is easy to see that $\mathcal P$ is dualizable: By \cref{rslt:equivalent-defs-for-dualizable} we need to check that the natural map $P \tensor \IHom(P, \Lambda^a) \isoto \IHom(P, P)$ is an isomorphism. Since $f$ is a cover, it is enough to show this isomorphism after applying $f^*$. But by the identity on $\IHom$ which we have just shown, this turns into $f^* P \tensor \IHom(f^* P, \Lambda^a) \isoto \IHom(f^* P, f^* P)$, which is true by \cref{rslt:equivalent-defs-for-dualizable}.

Altogether we obtain from the above that the full subcategory of $\Dqcohri(X, \Lambda)$ consisting of the dualizable objects satisfies v-hyperdescent. In particular if an object $\mathcal P \in \Dqcohri(X, \Lambda)$ is weakly almost perfect then on a cover by totally disconnected spaces it is weakly almost perfect, hence by \cref{rslt:equiv-of-compact-perfect-dualizable-over-Lambda} dualizable. By descent, $\mathcal P$ itself is dualizable. Conversely if $\mathcal P$ is dualizable then locally on totally disconnected spaces it is almost weakly perfect (by \cref{rslt:equiv-of-compact-perfect-dualizable-over-Lambda}), so that $\mathcal P$ is weakly almost perfect by definition. This finishes the proof of (ii). Now (i) follows because we have just shown that dualizable objects satisfy hypercomplete v-descent.

It remains to show (iii), so assume that $X$ is a locally spatial diamond and $\mathcal P \in \Dqcohri(X, \Lambda)$ is almost compact and discrete. Then by \cite[Proposition 11.24]{etale-cohomology-of-diamonds} we can find a cover $(f_i\colon Y_i \to X)_i$ of $X$ by totally disconnected spaces $Y_i$ such that all $f_i$ are qcqs and quasi-pro-étale (in particular $p$-bounded). Thus \cref{rslt:pullback-of-almost-compact-along-qcqs-pbd-is-almost-compact} implies that each $f_i^* \mathcal P$ is almost compact (and discrete), so by \cref{rslt:equiv-of-compact-perfect-dualizable-over-Lambda} it is weakly almost perfect. By (i) $\mathcal P$ is weakly almost perfect.
\end{proof}

For almost compact objects in $\Dqcohri(X, \Lambda)$ to become really useful, one would like them to generate all of $\Dqcohri(X, \Lambda)$ under colimits. At least if $X$ is a $p$-bounded locally spatial diamond, this is indeed the case:

\begin{proposition}
Let $X$ be a $p$-bounded locally spatial diamond.
\begin{propenum}
	\item For every étale map $j\colon U \to X$ with $U$ qcqs and every profinite set $S$, the object
	\begin{align*}
		\Lambda^a_\solid[S][U] := j_! \Lambda^a_\solid[S] \in \Dqcohri(X, \Lambda)
	\end{align*}
	is almost compact.

	\item \label{rslt:almost-compact-objects-generate-Dqcohri-on-loc-spat-diam} Every $\mathcal M \in \Dqcohri(X, \Lambda)$ is a filtered colimit of copies of $\Lambda^a_\solid[S][U][n]$ for varying profinite sets $S$, étale $U \to X$ with qcqs $U$ and integers $n$.

	\item \label{rslt:compact-implies-bounded-on-p-bd-loc-spat-diam} Suppose that $X$ is quasiseparated. Then every almost compact object in $\Dqcohri(X, \Lambda)$ is weakly almost bounded.
\end{propenum}
\end{proposition}
\begin{proof}
We first prove (i), so let $U$ and $S$ be given. By \cref{rslt:etale-lower-shriek-of-almost-compact-is-almost-compact} it is enough to show that $\Lambda^a[S] \in \Dqcohri(U, \Lambda)$ is almost compact. We have
\begin{align*}
	\Hom(\Lambda^a_\solid[S], -)^a = \IHom_{\Lambda(U)^a}(\Lambda(U)^a_\solid[S], \Gamma(U, -)),
\end{align*}
so the almost compactness of $\Lambda^a_\solid[S]$ reduces to the almost compactness of $\Lambda(U)^a_\solid[S]$ and the fact that $\Gamma(U, -)$ preserves small colimits by \cref{rslt:global-sections-on-qcqs-p-bounded-v-stack-preserve-colimits}.

To prove (ii), note that by \cref{rslt:etale-cohomology-is-conservative-on-p-bounded-diam} the family of functors
\begin{align*}
	\Gamma(U, -) = \Hom(\Lambda^a_\solid[U], -)^a\colon \Dqcohri(X, \Lambda) \to \Dqcohri(\Lambda(X))
\end{align*}
is conservative. Thus (ii) follows immediately from \cref{rslt:almost-compact-generators-implies-generated}.

We now prove (iii), so assume that $X$ is quasiseparated and let $\mathcal P \in \Dqcohri(X, \Lambda)$ be almost compact. Fix any $\varepsilon \in \mm_\Lambda(X)$. By \cref{rslt:characterization-of-almost-compact-objects-in-general-C} there is some object $\mathcal P_\varepsilon \in \Dqcohri(X, \Lambda)$ which is obtained via finite (co)limits and retracts from $\Lambda^a_\solid[S][U][n]$'s such that $\mathcal P$ is an $\varepsilon$-retract of $\mathcal P_\varepsilon$. It is thus enough to show that $\mathcal P_\varepsilon$ is bounded, which reduces immediately to showing that every $\Lambda^a_\solid[S][U]$ is bounded. But since $U$ is quasicompact and $X$ is quasiseparated, the map $U \to X$ is quasicompact, so the claim follows from \cref{rslt:etale-lower-shriek-of-generator-is-bounded}.
\end{proof}

It follows easily from the above results that \cref{rslt:equiv-of-compact-perfect-dualizable-over-Lambda} holds verbatim on $p$-bounded spatial diamonds:

\begin{proposition} \label{rslt:equiv-of-compact-perfect-dualizable-on-p-bd-spat-diam}
Let $X \in \vStacksCoeff$ be a $p$-bounded spatial diamond and let $\mathcal P \in \Dqcohri(X, \Lambda)$. Then the following are equivalent:
\begin{propenum}
	\item $\mathcal P$ is almost compact and discrete.
	\item $\mathcal P$ is weakly almost perfect.
	\item $\mathcal P$ is dualizable.
\end{propenum}
If this is the case then $\mathcal P$ is weakly almost bounded.
\end{proposition}
\begin{proof}
The equivalence of (ii) and (iii) is found in \cref{rslt:weakly-alm-perf-equiv-dualizable}. By \cref{rstl:discrete-compact-implies-perfect-on-loc-spatial-diam} (i) implies (ii). It remains to show that (iii) implies (i), so assume that $\mathcal P$ is dualizable. Then it is discrete by \cref{rslt:equiv-of-compact-perfect-dualizable-over-Lambda} (because discreteness can be checked locally, e.g. on totally disconnected spaces) and by \cref{rslt:equivalent-defs-for-dualizable} we have
\begin{align*}
	\Hom(\mathcal P, -)^a = \Gamma(X, \IHom(\mathcal P, -)) = \Gamma(X, - \tensor \mathcal P^*).
\end{align*}
But $\Gamma(X, -)$ preserves small colimits by \cref{rslt:global-sections-on-qcqs-p-bounded-v-stack-preserve-colimits}, hence $\mathcal P$ is almost compact. The last part of the claim follows from \cref{rslt:compact-implies-bounded-on-p-bd-loc-spat-diam}.
\end{proof}

We now turn our focus on affinoid perfectoid spaces. In this case we have the following description of almost compact and (weakly) almost perfect objects:

\begin{proposition}
Let $X$ be an affinoid perfectoid space with integral torsion coefficients $\Lambda$. Then there is a natural fully faithful inclusion
\begin{align*}
	\Dqcohri(\Lambda(X))_\waperf \subset \Dqcohri(X, \Lambda)_\waperf.
\end{align*}
If $X$ is $p$-bounded then there are natural equivalence of $\infty$-categories
\begin{align*}
	\Dqcohri(X, \Lambda)^{a\omega} &= \Dqcohri(\Lambda(X))^{a\omega},\\
	\Dqcohri(X, \Lambda)_\waperf &= \Dqcohri(\Lambda(X))_\waperf.
\end{align*}
\end{proposition}
\begin{proof}
We first prove the second claim, so assume that $X$ is $p$-bounded. Then by \cref{rslt:compute-Dqcohri-for-p-bounded-over-tot-disc} we have $\Dqcohri(X, \Lambda) = \Dqcohri(\Lambda(X))$, which immediately implies the claim about almost compact objects. The claim about weakly almost perfect objects follows from the fact that on both sides of the equivalence the weakly almost perfect objects are precisely the dualizable objects (see \cref{rslt:equiv-of-compact-perfect-dualizable-over-Lambda,rslt:weakly-alm-perf-equiv-dualizable}).

We now prove the first claim, so let $X$ be general. We need to show that the functor
\begin{align*}
	\widetilde{(-)}\colon \Dqcohri(\Lambda(X))_\waperf \injto \Dqcohri(X, \Lambda)
\end{align*}
is fully faithful. This amounts to showing that for every weakly almost perfect object $P \in \Dqcohri(\Lambda(X))_\waperf$ the natural morphism $P \isoto \Gamma(X, \widetilde P)$ is an isomorphism. If $P = \Lambda(X)^a$ then this amounts to saying that the cohomology of $\Lambda^a$ on $X$ vanishes, which is the case by definition of integral torsion coefficients. Both sides of the claimed isomorphism commute with finite (co)limits and retracts, hence we deduce that it is true in the case that $P = P'^a$ for some perfect object $P' \in \D(\Lambda(X))_\omega$. Now let $P$ be an arbitrary almost perfect object and fix any $\varepsilon \in \mm_\Lambda(X)$. Then there is a perfect object $P_\varepsilon \in \D(\Lambda(X))_\omega$ such that $P$ is an $\varepsilon$-retract of $P^a_\varepsilon$. We obtain a diagram
\begin{center}\begin{tikzcd}
	P \arrow[r] \arrow[d,shift left] & \Gamma(X, \widetilde P) \arrow[d,shift left]\\
	P_\varepsilon^a \arrow[r,"\sim"] \arrow[u,shift left] & \Gamma(X, \widetilde{P_\varepsilon^a}) \arrow[u,shift left]
\end{tikzcd}\end{center}
Like in previous proofs, a quick diagram chase shows that on homotopy groups both kernel and cokernel of the upper horizontal map are killed by $\varepsilon$. Hence the upper horizontal map is an isomorphism, as desired.
\end{proof}

A special case of dualizable objects are invertible objects in a symmetric monoidal $\infty$-category. They are of particular interest in the definition of cohomologically smooth maps later on.

\begin{definition}
Let $\mathcal C$ be a symmetric monoidal $\infty$-category. An object $L \in \mathcal C$ is called \emph{invertible} if there is an object $L' \in \mathcal C$ such that $L \tensor L' \isom 1_{\mathcal C}$.
\end{definition}

As invertible objects are in particular dualizable, all of the above results apply to them. In particular, it follows easily that they satisfy v-descent:

\begin{lemma} \label{rslt:invertible-is-v-local}
Let $f\colon Y \surjto X$ be a v-cover in $\vStacksCoeff$ and $\mathcal L \in \Dqcohri(X, \Lambda)$. Then $\mathcal L$ is invertible if and only if $f^* \mathcal L$ is.
\end{lemma}
\begin{proof}
The ``only if'' part is clear because $f^*$ is symmetric monoidal. We now prove the ``if'' part, so assume that $f^* \mathcal L$ is invertible. Then it is in particular dualizable, hence by \cref{rslt:geometric-properties-of-dualizable-and-perfect-sheaves} $\mathcal L$ is dualizable and $(f^* \mathcal L)^* = f^* \mathcal L^*$. But as $f^* \mathcal L$ is invertible, the natural map $f^*(\mathcal L \tensor \mathcal L^*) = f^* \mathcal L \tensor (f^* \mathcal L)^* \isoto \Lambda^a$ is an isomorphism, so since $f^*$ is conservative, also $\mathcal L \tensor \mathcal L^* \isoto \Lambda^a$ is an isomorphism, as desired.
\end{proof}

The v-descent result \cref{rslt:invertible-is-v-local} for invertible sheaves is enough for all our purposes, but one may still wonder what precise structure invertible sheaves in $\Dqcohri(X, \Lambda)$ possess. For general $\Lambda$ it is hard to say much, e.g. on totally disconnected spaces the $\infty$-category $\Dqcohri(X, \Lambda) = \Dqcohri(\Lambda(X))$ is a rather general category of modules. One may still hope to understand invertible sheaves in $\Dqcohri(\ri^+_{X^\sharp}/\pi)$ for an untilt $X^\sharp$ with pseudouniformizer $\pi$, e.g. one may ask if these sheaves are (étale/pro-étale/v) locally concentrated in a single degree and free. However, invertible sheaves turn out to be quite subtle due to the almost mathematics involved (this is already the case for almost modules over rings, cf. the remark before \cite[Definition 4.4.32]{almost-ring-theory}). In the following we present some partial results on invertible sheaves in $\Dqcohri(\ri^+_{X^\sharp}/\pi)$, which are independent of the rest of the thesis, but might be of use in some applications.

We start by analyzing invertible sheaves on a point. Even in this case it is not true that they are automatically isomorphic to a shift of the structure sheaf. For example, if $K$ is a perfectoid field with pseudouniformizer $\pi$, then for every $\gamma \in \R_{>0}$ we can form the ideal $I_\gamma = \{ x \in \ri_K \setst v(x) \ge \gamma \}$, where $v$ is the valuation on $K$. If $\gamma$ does not lie in the image of $v$ then $I_\gamma$ is non-trivial, even $I^a_\gamma \not\isom \ri^a_K$. But it is still invertible in $\D^a_\solid(\ri_K)$; in particular $I^a_\gamma/\pi I^a_\gamma$ is a non-trivial invertible object in $\D^a_\solid(\ri_K/\pi)$. We believe that all invertible objects in $\D^a_\solid(\ri_K/\pi)$ are of this form (up to shift), but we are not able to show this. In fact, it is not even clear to us why every invertible object in $\D^a_\solid(\ri_K/\pi)$ is bounded. On the other hand, assuming boundedness, we get quite close to the conjectured form $I^a_\gamma/\pi I^a_\gamma$. This is the content of our first result on invertible sheaves (see \cref{rslt:strictly-invertible-on-point-equiv-approx-trivial} below). As a preparation we first need a few lemmas.

\begin{lemma} \label{rslt:ri-K-mod-p-is-almost-injective-over-itself}
Let $K$ be a perfectoid field with pseudouniformizer $\pi$. Then $\ri^a_K/\pi$ is almost injective as a discrete module over itself, i.e. the functor
\begin{align*}
	\IHom(- , \ri^a_K/\pi)(*)\colon \D(\ri^a_K/\pi)_\omega^\opp \to \D(\ri^a_K/\pi)_\omega
\end{align*}
is $t$-exact.
\end{lemma}
\begin{proof}
We have to show that for every almost injection $N \to M$ of classical $\ri_K/\pi$-modules the induced map $\pi_0\IHom(M^a, \ri_K^a/\pi)(*) \to \pi_0\IHom(N^a, \ri_K^a/\pi)(*)$ is surjective. By replacing $M$ by $\pi_0 (M^a)_*$ and $N$ by $\pi_0 (N^a)_*$ we can assume that $N \injto M$ is an actual injection. By \cref{rslt:compute-almost-Hom-over-V} we have $\pi_0\IHom(M^a, \ri_K^a/\pi)(*) = \pi_0\IHom(M \tensor \mm_K, \ri_K/\pi)(*)^a$ (and similarly for $N$), so we need to show that the natural map
\begin{align*}
	\pi_0\IHom(M \tensor \mm_K, \ri_K/\pi)(*) \to \pi_0\IHom(N \tensor \mm_K, \ri_K/\pi)(*)
\end{align*}
of classical $\ri_K/\pi$-modules is almost surjective. Thus pick any $\ri_K/\pi$-linear map $f\colon N \tensor \mm_K \to \ri_K/\pi$ and any $\varepsilon \in \mm_K$. We need to extend $\varepsilon f$ from $N \tensor \mm_K$ to $M \tensor \mm_K$. By Zorn's Lemma we reduce to the case $M = N + (\ri_K/\pi) x$ for some $x \in M$. Let $\mm_x \subset \mm_K$ be the ideal of those $m \in \mm_K$ such that $m \tensor x \in \mm_K \tensor N$. We seek to define $g(m \tensor x)$ for all $m \in \mm_K$ such that $g(m \tensor x) = \varepsilon f(m \tensor x)$ if $m \in \mm_x$.

In the case $\mm_x = \mm_K$ there is nothing to do, so we can assume $\mm_x \subsetneq \mm_K$. Fix any $m_0 \in \mm_K \setminus \mm_x$ such that $\varepsilon m_0 \in \mm_x$. Then
\begin{align*}
	\frac\pi{m_0} f(\varepsilon m_0 \tensor x) = f(\varepsilon \pi \tensor x) = 0,
\end{align*}
hence $f(\varepsilon m_0 \tensor x) = m_0 y$ for some $y \in \ri_K/\pi$. We can thus define $g(m \tensor x) = my$ for all $m \in \mm_K$, which clearly satisfies the required compatibility with $\varepsilon f$.
\end{proof}

\begin{lemma} \label{rslt:tot-disc-cohom-of-almost-compact-object-is-almost-fin-pres}
Let $X = \Spa(A, A^+)$ be a totally disconnected perfectoid space and $P \in \D^a(A^+/\pi)_\omega$ a discrete and almost compact object. Then for all $n \in \Z$, $\pi_n(P)$ is an almost finitely presented $A^{+a}/\pi$-module.
\end{lemma}
\begin{proof}
Fix $n \in \Z$. Given $\varepsilon \in \mm_A$, we use \cref{rslt:approximate-almost-compact-by-compact} to choose a discrete compact object $P_\varepsilon \in \D(A^+/\pi)_\omega$ such that $P$ is an $\varepsilon$-retract of $P^a_\varepsilon$. In particular $\pi_n(P)$ is an $\varepsilon$-retract of $\pi_n(P_\varepsilon)^a$. As $P_\varepsilon$ is perfect and $A^+/\pi$ is coherent (see \cite[Proposition 10.5.2]{scholze-berkeley-lectures}), $\pi_n(P_\varepsilon)$ is a finitely presented $A^+/\pi$-module, hence the functor $\pi_0\IHom(\pi_n(P_\varepsilon)^a, -)\colon \catsldmoda{A^+/\pi} \to \catsldmoda{A^+/\pi}$ preserves all filtered colimits. By the same diagram chase argument as in the proof of \cref{rslt:characterization-of-almost-compact-objects-in-general-C} we deduce that the functor $\pi_0\IHom(\pi_n(P), -)$ also preserves filtered colimits. By \cite[Proposition 2.3.15.(ii)]{almost-ring-theory} it follows that $\pi_n(P)$ is almost finitely presented, as desired.
\end{proof}

\begin{lemma} \label{rslt:strictly-invertible-on-point-equiv-approx-trivial}
Let $K$ be a perfectoid field with an open and bounded valuation subring $K^+$ and let $L \in \Dqcohri(K^+/\pi)$. Then the following are equivalent:
\begin{lemenum}
	\item $L$ is invertible and bounded.
	\item There is an integer $n$ such that $L \approx K^{+a}/\pi[n]$.
\end{lemenum}
\end{lemma}
\begin{proof}
We first prove that (i) implies (ii), so assume that $L$ is invertible and bounded. By \cref{rslt:equiv-of-compact-perfect-dualizable-over-Lambda} $L$ is discrete, and since the dual $L^\vee$ is also invertible, $L^\vee$ is also discrete. Thus \cref{rslt:ri-K-mod-p-is-almost-injective-over-itself} implies that $L^\vee$ is bounded. Let $n$ and $n'$ be the lowest homological degrees such that $\pi_n(L) \ne 0$ and $\pi_{n'}(L^\vee) \ne 0$. Then $\pi_{n+n'}(L \tensor L^\vee) = \pi_0(\pi_n(L) \tensor \pi_{n'}(L^\vee))$. By \cref{rslt:tot-disc-cohom-of-almost-compact-object-is-almost-fin-pres} both $\pi_n(L)$ and $\pi_{n'}(L^\vee)$ are almost finitely presented and hence by \cite[Theorem 2.5]{rigid-p-adic-hodge} (note that w.l.o.g. $K^+ = \ri_K$ because $K^{+a} \isom \ri^a_K$ and hence $\D^a(K^+/\pi) = \D^a(\ri_K/\pi)$) they are $\approx$ to a direct sum of copies of $K^{+a}/I_\gamma$ for certain non-negative real numbers $\gamma$, where $I_\gamma \subset K^+$ is the ideal generated by all elements of valuation at least $\gamma$. Clearly $\pi_0(K^{+a}/I_\gamma \tensor K^{+a}/I_{\gamma'}) = K^{+a}/I_{\min(\gamma,\gamma')}$, which is only zero if one of the two factors is zero. Thus $\pi_0(\pi_n(L) \tensor \pi_{n'}(L^\vee))$ is nonzero and hence we must have $n + n' = 0$ and $\pi_0(\pi_n(L) \tensor \pi_{n'}(L^\vee)) = K^{+a}/\pi$. Clearly $L^\vee$ vanishes in homological degrees above $-n$, hence it must be concentrated in the single degree $n' = -n$. Repeating the argument for $L^\vee$ in place of $L$ shows that $L = L^{\vee\vee}$ is also concentrated in a single degree. Finally, by again considering tensor products of $K^{+a}/I_\gamma$ as above we deduce that we must have $\pi_n(L) \approx \pi_{n'}(L^\vee) \approx K^{+a}/\pi$.

We now prove that (ii) implies (i), so assume that $L \approx K^{+a}/\pi[n]$ for some $n$. Clearly this implies that $L$ is bounded, so we only need to show that it is invertible. W.l.o.g. $n = 0$. Pick any $\varepsilon \in \mm_K$ and any $f_\varepsilon\colon L \rightleftarrows K^{+a}/\pi\noloc g_\varepsilon$ realising $L \approx_\varepsilon K^{+a}/\pi$. Then $f_\varepsilon^\vee$ and $g_\varepsilon^\vee$ realise $L^\vee \approx_\varepsilon (K^{+a}/\pi)^\vee = K^{+a}/\pi$ and the following diagrams (with either only the upward or only the downward maps) commute:
\begin{center}\begin{tikzcd}
	L \tensor L^\vee \arrow[r] \arrow[d,"f_\varepsilon \tensor g_\varepsilon^\vee",shift left] & K^{+a}/\pi \arrow[d,"\varepsilon",shift left]\\
	K^{+a}/\pi \tensor K^{+a}/\pi \arrow[u,"g_\varepsilon \tensor f_\varepsilon^\vee",shift left] \arrow[r,"\sim"] & (K^{+a}/\pi)^\vee \arrow[u,"\varepsilon",shift left]
\end{tikzcd}\end{center}
From the fact that the bottom horizontal map is an isomorphism, it follows by diagram chase that both the kernel and cokernel of the top horizontal map are killed by $\varepsilon$. As $\varepsilon$ was arbitrary, the top map is an isomorphism, as desired.
\end{proof}

Having partially understood the case of points, we can now extend our study of invertible sheaves to more general spaces. We obtain the following result, showing that invertible sheaves behave roughly as expected, at least if one assumes them to be bounded (which we conjecture to be true automatically):

\begin{proposition} \label{rslt:strictly-invertible-sheaf-is-locally-concentrated-in-single-degree}
Let $X$ be an untilted small v-stack with pseudouniformizer $\pi$ and let $\mathcal L \in \DqcohriX X$ be invertible and bounded. Then $\mathcal L$ is discrete and there is a disjoint open cover $X = \bigdunion_{n\in\Z} U_n$ such that for all $n$, $\restrict{\mathcal L}{U_n} \in \D^a(\ri^+_U/\pi)_\omega$ is concentrated in homological degree $n$ and flat.
\end{proposition}
\begin{proof}
Discreteness of $\mathcal L$ can be checked on fibers and thus follows from \cref{rslt:strictly-invertible-on-point-equiv-approx-trivial}. Now assume that $X = \Spa(A, A^+)$ is totally disconnected. For every $x \in \pi_0(X)$ let $n(x)$ be the unique integer such that $\pi_{n(x)}(\mathcal L_x)$ is non-zero (see \cref{rslt:strictly-invertible-on-point-equiv-approx-trivial}). We claim that the map $\alpha\colon \pi_0(X) \to \Z$, $x \mapsto n(x)$ is continuous. Clearly $\mathcal L$ is almost compact hence by \cref{rslt:tot-disc-cohom-of-almost-compact-object-is-almost-fin-pres} all homology groups $\pi_n(\mathcal L)$ are almost finitely presented. Fix $n$, pick any $\varepsilon \in \mm_A$ and choose some finitely presented $A^+/\pi$-module $\mathcal L_\varepsilon$ such that $\pi_n(\mathcal L) \approx_\varepsilon \mathcal L_\varepsilon$. Then if for some $x \in \pi_0(X)$, $\pi_n((\mathcal L_\varepsilon)_x)$ is killed by $\varepsilon$, then the same is true for all points in some open neighbourhood $U \subset \pi_0(X)$ of $x$. In other words, $n(x) \ne n$ implies that $n(y) \ne n$ for all $y$ in some neighbourhood $U$ of $x$. This implies that $\alpha^{-1}(\Z \setminus \{ n \}) \subset \pi_0(X)$ is open. As the image of $\alpha$ is finite (by boundedness of $\mathcal L$), we deduce the continuity of $\alpha$. Thus letting $U_n := \alpha^{-1}(n)$ we get the desired open decomposition $X = \bigdunion_n U_n$ such that each $\restrict{\mathcal L}{U_n}$ is concentrated in degree $n$, say $\restrict{\mathcal L}{U_n} = L_n[n]$ for some $L_n \in \D^a(\ri^+_U/\pi)^\heartsuit_\omega$. The same then follows for $\mathcal L^\vee$, from which we deduce easily that in fact $L_n$ is invertible as an object in the heart. This implies that $L_n$ is flat, finishing the proof of the claim in the case that $X$ is totally disconnected.

Now let $X$ be general. It is still true that for every point $x\colon \Spa(K, K^+) \to X$, $x^*\mathcal L$ is concentrated in a single degree $n(x)$ which is independent of representative $x$ in the equivalence class of the induced point in $\abs X$. This defines a map $\alpha\colon \abs X \to \Z$. To check that $\alpha$ is continuous, pick any cover $Y \surjto X$, where $Y = \bigdunion_i Y_i$ is a disjoint union of totally disconnected perfectoid spaces $Y_i$. Then $\abs X$ is (by definition) equipped with the quotient topology from $\abs Y$, so to show that $\alpha$ is continuous, it is enough to show that the composition $\abs Y \to \abs X \to \Z$ is continuous. This then reduces to the restricted maps $Y_i \to \Z$, which were seen in the previous paragraph to be continuous.
\end{proof}

\subsection{Cohomologically Smooth Morphisms} \label{sec:ri-pi.smooth}

Since there is no general notion of smooth maps of diamonds, the 6-functor formalism \cref{rslt:main-6-functor-formalism} so far lacks a version of Poincaré duality. In this section we fix this issue by providing a definition of ``$p$-cohomologically smooth'' maps of small v-stacks, similar to the $\ell$-adic analog in \cite[\S23]{etale-cohomology-of-diamonds}. We then derive some basic properties of $p$-cohomologically smooth maps and in particular show that they are stable under the expected operations. Apart from étale maps, we will not provide any concrete examples of $p$-cohomologically smooth maps however; this will be postponed to \cref{sec:ri-pi.poincare}.

Here is the definition of $p$-cohomologically smooth maps:

\begin{definition} \label{def:p-cohomologically-smooth-map}
A bdcs map $f\colon Y \to X$ of small v-stacks is \emph{$p$-cohomologically smooth} if for every map $X' \to X$ from a totally disconnected perfectoid space $X' = \Spa(A', A'^+)$ there is a pseudouniformizer $\pi \in A'^+$ such that the map $f'\colon Y' := Y \cprod_X X' \to X'$ satisfies the following property: The natural morphism
\begin{align*}
	f'^! (\ri^{+a}_{X'}/\pi) \tensor f'^* \isoto f'^!
\end{align*}
is an isomorphism of functors $\DqcohriX{X'} \to \DqcohriX{Y'}$ and $f'^! (\ri^{+a}_{X'}/\pi)$ is invertible.
\end{definition}

As a first step we show that $p$-cohomologically smooth maps satisfy Poincaré duality on any base $X$ (not just totally disconnected spaces) and with respect to any system of integral torsion coefficients $\Lambda$ on $X$. This will be the content of \cref{rslt:p-cohom-smoothness-description-of-upper-shriek} below. In order to allow arbitrary systems of coefficients, we need the following compatibility result (analogous to \cite[Remark 23.2]{etale-cohomology-of-diamonds}).

\begin{lemma} \label{rslt:upper-shriek-commutes-with-forgetful-functor-and-change-of-rings}
Let $f\colon Y \to X$ be a bdcs map of small v-stacks and let $\Lambda \to \Lambda'$ be a morphism of integral torsion coefficients on $X$. Then $f^!$ commutes with the forgetful functor along $\Lambda \to \Lambda'$, i.e. there is a commuting diagram
\begin{center}\begin{tikzcd}
	\Dqcohri(X, \Lambda') \arrow[r,"f^!"] \arrow[d,swap,"\mathrm{forget}"] & \Dqcohri(Y, \Lambda') \arrow[d,"\mathrm{forget}"]\\
	\Dqcohri(X, \Lambda) \arrow[r,"f^!"] & \Dqcohri(Y, \Lambda)
\end{tikzcd}\end{center}
\end{lemma}
\begin{proof}
By passing to left adjoints, the claim reduces to \cref{rslt:coeff-projection-formula}.
\end{proof}

\begin{lemma} \label{rslt:discrete-projection-formula-along-hom-of-coeff}
Let $X$ be a small v-stack and $\Lambda \to \Lambda'$ a morphism of integral torsion coefficients on $X$. Then for all $\mathcal N \in \Dqcohri(X, \Lambda')$ and all discrete $\mathcal M \in \Dqcohri(X, \Lambda)_\omega$ the natural morphism
\begin{align*}
	\mathcal N \tensor_{\Lambda^a_\solid} \mathcal M \isoto \mathcal N \tensor_{\Lambda'^a_\solid} (\mathcal M \tensor_{\Lambda^a_\solid} \Lambda'^a_\solid)
\end{align*}
is an isomorphism in $\Dqcohri(X, \Lambda)$. Moreover, we have $\mathcal M \tensor_{\Lambda^a_\solid} \Lambda'^a_\solid = \mathcal M \tensor_{\Lambda^a_\solid} \Lambda'^a$, where on the right-hand side we view $\Lambda'^a$ as an algebra object in $\Dqcohri(X, \Lambda)$.
\end{lemma}
\begin{proof}
As both sides of the claimed isomorphism commute with pullbacks (see \cref{rslt:Dqcohri-change-of-coefficients}), we can w.l.o.g. assume that $X$ is a totally disconnected perfectoid space. In this case we have $\Dqcohri(X, \Lambda) = \Dqcohri(\Lambda(X))$ and since $\mathcal M$ is discrete, it is a colimit of copies of $\Lambda^a$. Since both sides of the claimed isomorphism commute with colimits in $\mathcal M$ we can reduce to the case $\mathcal M = \Lambda^a$. But then there is nothing to prove.
\end{proof}

\begin{proposition} \label{rslt:p-cohom-smoothness-identities}
\begin{propenum}
	\item \label{rslt:p-cohom-smoothness-description-of-upper-shriek} Let $f\colon Y \to X$ be a bdcs and $p$-cohomologically smooth map in $\vStacksCoeff$. Then the natural morphism
	\begin{align*}
		f^! \Lambda^a \tensor f^* \isoto f^!
	\end{align*}
	of functors $\Dqcohri(X, \Lambda) \to \Dqcohri(Y, \Lambda)$ is an isomorphism and $f^! \Lambda^a$ is invertible.

	\item \label{rslt:p-cohom-smoothness-base-change-for-upper-shriek} Let
	\begin{center}\begin{tikzcd}
		Y' \arrow[r,"g'"] \arrow[d,"f'"] & Y \arrow[d,"f"]\\
		X' \arrow[r,"g"] & X
	\end{tikzcd}\end{center}
	be a cartesian square in $\vStacksCoeff$ and assume that $f$ is bdcs and $p$-cohomologically smooth. Then the natural morphism
	\begin{align*}
		g'^* f^! \isoto f'^! g^*
	\end{align*}
	is an isomorphism of functors $\Dqcohri(X, \Lambda) \to \Dqcohri(Y', \Lambda)$.
\end{propenum}
\end{proposition}
\begin{proof}
We first prove (i) in the case that $X$ is totally disconnected. Note that (i) reduces to the following two claims:
\begin{enumerate}[(a)]
	\item If $\Lambda \to \Lambda'$ is a morphism of integral torsion coefficients on $X$ and (i) holds for $\Lambda$ then it also holds for $\Lambda'$.

	\item If there is some untilt $X^\sharp$ of $X$ and some pseudouniformizer $\pi$ on $X^\sharp$ such that (i) holds for $\Lambda = \ri^+_{X^\sharp}/\pi$, then (i) also holds for $\Lambda' = \ri^+_{X^\sharp}/\pi^2$.
\end{enumerate}
We first prove (a), so let $\Lambda \to \Lambda'$ be given such that (i) holds for $\Lambda$. In the following we use subscripts $\Lambda$ and $\Lambda'$ for all functors to distinguish their versions on $\Dqcohri(-, \Lambda)$ and $\Dqcohri(-, \Lambda')$. It is enough to show the claimed isomorphism of functors after applying the forgetful functor along $\Lambda \to \Lambda'$ because the forgetful functor is conservative by \cref{rslt:Dqcohri-change-of-coefficients}. By \cref{rslt:upper-shriek-commutes-with-forgetful-functor-and-change-of-rings,rslt:Dqcohri-change-of-coefficients} we compute, for every $\mathcal M \in \Dqcohri(X, \Lambda')$,
\begin{align*}
	f_{\Lambda'}^! \mathcal M = f_\Lambda^! \mathcal M = f_\Lambda^* \mathcal M \tensor_\Lambda f_\Lambda^! \Lambda^a = f_{\Lambda'}^* \mathcal M \tensor_\Lambda f_\Lambda^! \Lambda^a.
\end{align*}
Plugging in $\mathcal M = \Lambda'^a$ we obtain $f_{\Lambda'}^! \Lambda'^a = \Lambda'^a \tensor_\Lambda f_\Lambda^! \Lambda^a$. Thus by \cref{rslt:discrete-projection-formula-along-hom-of-coeff} we have $f_{\Lambda'}^! \Lambda'^a = f_\Lambda^! \Lambda^a \tensor_{\Lambda^a_\solid} \Lambda'^a_\solid$, which shows that $f_{\Lambda'}^! \Lambda'^a$ is invertible. Moreover, \cref{rslt:discrete-projection-formula-along-hom-of-coeff} implies
\begin{align*}
	f_{\Lambda'}^* \mathcal M \tensor_{\Lambda'} f_{\Lambda'}^! \Lambda'^a = f_{\Lambda'}^* \mathcal M \tensor_{\Lambda'} (f_\Lambda^! \Lambda^a \tensor_{\Lambda^a_\solid} \Lambda'^a_\solid) = f_{\Lambda'}^* \mathcal M \tensor_\Lambda f_\Lambda^! \Lambda^a.
\end{align*}
Comparing this to the above computation of $f_{\Lambda'}^! \mathcal M$ proves the desired isomorphism of functors. This finishes the proof of claim (a).

We now prove (b), so fix the untilt $X^\sharp = \Spa(A, A^+)$ and the pseudouniformizer $\pi \in A^+$ and assume that (i) holds for $\Lambda = \ri^+_{X^\sharp}/\pi$. To avoid confusion let us denote $f^!_\pi$ and $f^!_{\pi^2}$ the two versions of upper shriek we are dealing with. Given any $M \in \Dqcohri(A^+/\pi^2)$ we write $M_\pi := M \tensor_{A^{+a}/\pi^2} A^{+a}/\pi$. Then there is a triangle $M_\pi \to M \to M_\pi$ in $\Dqcohri(A^+/\pi^2)$, and as all our functors are exact, the proof of $f_{\pi^2}^! (A^{+a}/\pi^2) \tensor f_{\pi^2}^* M \isoto f_{\pi^2}^! M$ reduces to the case of $M_\pi$, i.e. that $M$ comes via forgetful functor from $\Dqcohri(A^+/\pi)$. In this case we get from \cref{rslt:upper-shriek-commutes-with-forgetful-functor-and-change-of-rings}
\begin{align*}
	f^!_{\pi^2} M = f^!_\pi M = f^!_\pi(A^{+a}/\pi) \tensor_\pi f^*_\pi M.
\end{align*}
By \cref{rslt:discrete-projection-formula-along-hom-of-coeff} it is therefore enough to show that $f^!_\pi(A^{+a}/\pi) = f^!_{\pi^2}(A^{+a}/\pi^2) \tensor_{\pi^2} \ri^{+a}_{Y^\sharp}/\pi$, i.e. we are reduced to the case $M = A^{+a}/\pi$. By assumption $f^!_\pi$ preserves all colimits. Using the triangle $M_\pi \to M \to M_\pi$ from above (and \cref{rslt:upper-shriek-commutes-with-forgetful-functor-and-change-of-rings}) we deduce that $f^!_{\pi^2}$ preserves all colimits. But $A^{+a}/\pi$ is a colimit of copies of $A^{+a}/\pi^2$ (inside $\Dqcohri(A^+/\pi^2)$, hence as both sides of the claimed identity preserve colimits, we reduce to the case $M = A^{+a}/\pi^2$; then there is nothing to prove. From now on we can drop the subscripts $\pi$ and $\pi^2$ for $f^!$.

To finish the proof of (b) it remains to show that $\mathcal L := f^!(A^{+a}/\pi^2)$ is invertible in $\Dqcohri(\ri^+_{Y^\sharp}/\pi^2)$. We claim that the functor
\begin{align*}
	F := \IHom_{\pi^2}(\mathcal L, -)\colon \Dqcohri(\ri^+_{Y^\sharp}/\pi^2) \to \Dqcohri(\ri^+_{Y^\sharp}/\pi^2)
\end{align*}
preserves all colimits. By plugging in the triangles $M_\pi \to M \to M_\pi$ used above we see that it is enough to show that $F$ precomposed with the forgetful functor $\DqcohriX{Y^\sharp} \to \Dqcohri(\ri^+_{Y^\sharp}/\pi^2)$ preserves colimits. But this composition of functors is computed by $\IHom_\pi(\mathcal L_\pi, -)$ and as $\mathcal L_\pi = f^!(A^+/\pi)$ is invertible, this functor preserves colimits. By writing $\ri^{+a}_{Y^\sharp}/\pi$ as a colimit of copies of $\ri^{+a}_{Y^\sharp}/\pi^2$ (in $\Dqcohri(\ri^+_{Y^\sharp}/\pi^2)$) we deduce
\begin{align*}
	(\mathcal L^\vee)_\pi &= \IHom_{\pi^2}(\mathcal L, \ri^{+a}_{Y^\sharp}/\pi^2) \tensor_{\pi^2} \ri^{+a}_{Y^\sharp}/\pi = \IHom_{\pi^2}(\mathcal L, \ri^{+a}_{Y^\sharp}/\pi) = \IHom_\pi(\mathcal L_\pi, \ri^{+a}_{Y^\sharp}/\pi)\\
	&= (\mathcal L_\pi)^\vee.
\end{align*}
Now in order to show that $\mathcal L$ is invertible we need to check that the natural map $\mathcal L \tensor_{\pi^2} \mathcal L^\vee \to \ri^{+a}_{Y^\sharp}/\pi^2$ is an isomorphism. Using the triangles $M_\pi \to M \to M_\pi$ from above it is enough to show this isomorphism after applying $(-)_\pi$. But by the computation $(\mathcal L^\vee)_\pi = (\mathcal L_\pi)^\vee$ this turns into $\mathcal L_\pi \tensor_\pi (\mathcal L_\pi)^\vee \to \ri^{+a}_{Y^\sharp}/\pi$, which is an isomorphism because $\mathcal L_\pi$ is invertible in $\DqcohriX{Y^\sharp}$. This finishes the proof of claim (b) and thus also the proof of (i) in the case that $X$ is totally disconnected.

We now prove (ii). First assume that $X$ and $X'$ are totally disconnected perfectoid spaces. By what we have shown in (i) the claimed isomorphism of functors reduces to showing that the natural map
\begin{align*}
	\alpha\colon g'^* f^! \Lambda^a \to f'^! \Lambda^a
\end{align*}
of invertible sheaves on $Y'$ is an isomorphism. This can be checked locally, so we can assume that $f$ is $p$-bounded and qcqs; in particular $Y$ is a $p$-bounded spatial diamond. Note that $g_*$ is conservative because it is just a forgetful functor. Similarly $g'_*$ is conservative on weakly almost left-bounded sheaves: By base-change (see \cref{rslt:base-change-for-bounded-Dqcohri}) we can check this after pullback to any cover of $Y$. Choose a $p$-bounded cover $\tilde Y \to Y$ by a totally disconnected perfectoid space. Then $\tilde Y' := \tilde Y \cprod_Y Y' = \tilde Y \cprod_X X'$ is a $p$-bounded affinoid perfectoid space (the map $\tilde Y' \to X'$ is $p$-bounded), hence by \cref{rslt:compute-Dqcohri-for-p-bounded-over-tot-disc} the functor $\tilde g'_*$ is a forgetful functor and thus conservative; here $\tilde g'\colon \tilde Y' \to \tilde Y$ denotes the base-change of $g'$. Since invertible sheaves are weakly almost left-bounded (see \cref{rslt:equiv-of-compact-perfect-dualizable-on-p-bd-spat-diam}), it is now enough to show that
\begin{align*}
	g'_*\alpha\colon g'_* g'^* f^! \Lambda^a \to g'_* f'^! \Lambda^a
\end{align*}
is an isomorphism. Let us denote $\mathcal L := f^! \Lambda^a$. By using the identity $g'_* f'^! = f^! g_*$ (which follows from proper base-change, see \cref{rslt:main-6-functor-formalism-proper-base-change}, by passing to right adjoints) and (i) we can rewrite $g'_* \alpha$ as
\begin{align*}
	g'_*\alpha\colon g'_* g'^* \mathcal L \to f^! g_* \Lambda^a = \mathcal L \tensor f^* g_* \Lambda^a.
\end{align*}
Applying weakly almost left bounded base-change again, we are reduced to showing that the natural map
\begin{align*}
	g'_*\alpha\colon g'_* g'^* \mathcal L \to \mathcal L \tensor g'_* \Lambda^a
\end{align*}
is an isomorphism. This can be checked on a v-cover, so pick $p$-bounded affinoid perfectoid covers $\tilde Y$ and $\tilde Y'$ as above and let $\tilde{\mathcal L}$ and $\tilde g'$ be the pullbacks of $\mathcal L$ and $g'$ to $\tilde Y$. Applying weakly almost left-bounded base-change again we need to see that the natural map
\begin{align*}
	\tilde g'_* \tilde g'^* \tilde{\mathcal L} \to \tilde{\mathcal L} \tensor \tilde g'_* \Lambda^a
\end{align*}
is an isomorphism. By \cref{rslt:compute-Dqcohri-for-p-bounded-over-tot-disc} we have $\Dqcohri(\tilde Y', \Lambda) = \Dqcohri(\Lambda(\tilde Y'))$ and $\Dqcohri(\tilde Y, \Lambda) = \Dqcohri(\Lambda(\tilde Y))$. Denoting $L = \tilde{\mathcal L} \in \Dqcohri(\Lambda(\tilde Y))$, the above map translates to
\begin{align*}
	L \tensor_{\Lambda(\tilde Y)^a_\solid} \Lambda(\tilde Y')^a_\solid \to L \tensor_{\Lambda(\tilde Y)^a_\solid} \Lambda(\tilde Y')^a.
\end{align*}
As $L$ and $\Lambda(\tilde Y')^a$ are discrete, this is an isomorphism; it is in fact isomorphic to the simple tensor product $L \tensor_{\Lambda(\tilde Y)^a} \Lambda(\tilde Y')^a$ of almost modules -- no solidification is needed. This finishes the proof of (ii) in the case that $X$ and $X'$ are totally disconnected.

It is now easy to deduce (ii) in general: Choose a hypercover $X_\bullet \to X$ by disjoint unions of totally disconnected spaces and let $Y_\bullet \to Y$ be the base-change of this hypercover to $Y$ with map $f_\bullet\colon Y_\bullet \to X_\bullet$ of hypercovers. By the above shown case of (ii) we deduce that the componentwise computed $f_\bullet^!$ preserves coCartesian edges and hence computes $f^!$ (cf. the proof of \cite[Proposition 23.12]{etale-cohomology-of-diamonds}). In particular, the base-change (ii) holds for the map $g\colon X_0 \to X$, i.e. for any hypercover by a disjoint union of totally disconnected spaces. For general $g\colon X' \to X$ choose a cover $g_0\colon X'_0 \to X_0$ of $g$ by a map of disjoint unions of totally disconnected spaces. Then we know that (ii) holds for $g_0$ and for the covers $X'_0 \surjto X'$ and $X_0 \surjto X$, from which it follows that it also holds for $g$.

It remains to prove (i) for general $X$. The claimed isomorphism of functors can be checked on any v-cover and by (ii) still has the same form on that cover, so the isomorphism is reduced to the case that $X$ is totally disconnected, which was shown above. By \cref{rslt:invertible-is-v-local} the same applies to the claim that $f^! \Lambda^a$ is invertible.
\end{proof}

Having established the fundamental properties of $p$-cohomologically smooth maps, we now show that they enjoy stability under the expected operations.

\begin{lemma} \label{rslt:stability-of-p-cohom-smooth}
\begin{lemenum}
	\item The condition of being bdcs and $p$-cohomologically smooth is analytically local on both source and target.

	\item Among bdcs maps, the condition of being $p$-cohomologically smooth is v-local on the target and $p$-cohomologically smooth local on the source.

	\item Bdcs and $p$-cohomologically smooth maps are stable under composition and base-change.

	\item Every étale map is bdcs and $p$-cohomologically smooth.
\end{lemenum}
\end{lemma}
\begin{proof}
For all of the claims above the required stability properties for bdcs maps are shown in \cref{rslt:stability-of-bdcs-maps}, so we only need to verify the $p$-cohomological smoothness.

We first prove part (iv). If $f\colon Y \to X$ is étale then $f^! = f^*$ (after pullback to any totally disconnected space and choosing any pseudouniformizer) which immediately implies that $f$ satisfies the smoothness condition \cref{def:p-cohomologically-smooth-map}.

We now prove part (iii). It is clear that $p$-cohomologically smooth maps are stable under base-change. Stability under composition follows formally from \cref{rslt:p-cohom-smoothness-identities}.

Part (i) is a special case of (ii), so it remains to prove (ii). We first show that $p$-cohomological smoothness is v-local on the target. Thus assume we are given a bdcs map $f\colon Y \to X$ and a v-cover $g\colon X' \surjto X$ such that the base-changed map $f'\colon Y' \to X'$ is $p$-cohomologically smooth. We want to show that $f$ is $p$-cohomologically smooth. By definition this reduces to the case that $X$ is totally disconnected. Now extend $g$ to a hypercover $X'_\bullet \to X$ with $X'_0 = X'$ and let $f'_\bullet\colon Y'_\bullet \to X'_\bullet$ be the base-change. Then by \cref{rslt:p-cohom-smoothness-base-change-for-upper-shriek} the componentwisely defined $f'^!_\bullet$ preserves coCartesian edges and hence computes $f^!$ (as in the proof of \cref{rslt:p-cohom-smoothness-base-change-for-upper-shriek}). In particular $f^!$ commutes with pullback to $X'$, so the $p$-cohomological smoothness of $f'$ implies that $f$ is $p$-cohomologically smooth as well (use also that invertibility is v-local, \cref{rslt:invertible-is-v-local}).

We now show that $p$-cohomological smoothness is $p$-cohomologically smooth on the source. Thus suppose we are given bdcs maps $g\colon Z \to Y$ and $f\colon Y \to X$ such that $g$ and $f \comp g$ are $p$-cohomologically smooth and $g$ is surjective. We want to show that $f$ is $p$-cohomologically smooth, for which we can assume that $X$ is totally disconnected; pick any pseudouniformizer $\pi$ on $X^\flat$ and denote $\Lambda = \ri^+_{X^\flat}/\pi$. Note that $g^!$ is of the form $\mathcal L_g \tensor g^*$ for an invertible sheaf $\mathcal L_g$ and hence conservative. Thus we can check that the natural map $f^! \Lambda^a \tensor f^* \isoto f^!$ is an isomorphism after applying $g^!$. But
\begin{align*}
	g^! f^! = (fg)^! = (fg)^!\Lambda^a \tensor (fg)^* = \mathcal L_g \tensor g^* f^! \Lambda^a \tensor g^*f^* = g^! (f^!\Lambda^a \tensor f^*),
\end{align*}
as desired. Similarly, we can verify that $f^!\Lambda^a$ is invertible after applying $g^! = \mathcal L_g \tensor g^*$, but $g^! f^!\Lambda^a$ is invertible because $f \comp g$ is $p$-cohomologically smooth.
\end{proof}

Next up we verify the expected smooth base-change results and compatibility with $\IHom$. The following results and proofs are completely analogous to \cite[Proposition 23.16]{etale-cohomology-of-diamonds} and \cite[Proposition 23.17]{etale-cohomology-of-diamonds}.

\begin{proposition}
Let
\begin{center}\begin{tikzcd}
	Y' \arrow[r,"g'"] \arrow[d,"f'"] & Y \arrow[d,"f"]\\
	X' \arrow[r,"g"] & X
\end{tikzcd}\end{center}
be a cartesian square in $\vStacksCoeff$.
\begin{propenum}
	\item Assume that $g$ is bdcs. Then the natural morphism
	\begin{align*}
		g^! f_* \isoto f'_* g'^!
	\end{align*}
	is an isomorphism of functors $\Dqcohri(Y, \Lambda) \to \Dqcohri(X', \Lambda)$.

	\item Assume that $g$ is bdcs and $p$-cohomologically smooth. Then the natural morphism
	\begin{align*}
		g^* f_* \isoto f'_* g'^*
	\end{align*}
	is an isomorphism of functors $\Dqcohri(Y, \Lambda) \to \Dqcohri(X', \Lambda)$.

	\item \label{rslt:smooth-upper-shriek-base-changes} Assume that $f$ is bdcs and that either $f$ or $g$ is bdcs and $p$-cohomologically smooth. Then the natural morphism
	\begin{align*}
		g'^* f^! \isoto f'^! g^*
	\end{align*}
	is an isomorphism of functors $\Dqcohri(X, \Lambda) \to \Dqcohri(Y', \Lambda)$.
\end{propenum}
\end{proposition}
\begin{proof}
Part (i) follows from proper base-change (see \cref{rslt:main-6-functor-formalism-proper-base-change}) by passing to right adjoints. Then part (ii) follows from (i) and \cref{rslt:p-cohom-smoothness-description-of-upper-shriek} by noting that $g^*$ and $g^!$ only differ by a strictly invertible sheaf.

It remains to prove (iii). In the case that $f$ is $p$-cohomologically smooth, this is just \cref{rslt:p-cohom-smoothness-base-change-for-upper-shriek}. Now assume instead that $g$ is bdcs and $p$-cohomologically smooth (hence the same is true for $g'$). We can check the claimed isomorphism after tensoring with the invertible sheaf $\mathcal L_{g'} := g'^! \Lambda^a$. Denoting also $\mathcal L_g := g^! \Lambda^a$, the left-hand side becomes
\begin{align*}
	\mathcal L_{g'} \tensor g'^* f^! = g'^! f^! = f'^! g^! = f'^!(\mathcal L_g \tensor g^*).
\end{align*}
We claim that the natural map $f'^* \mathcal L_g \tensor f'^! g^* \isoto f'^!(\mathcal L_g \tensor g^*)$ is an isomorphism. This follows from Yoneda by using the projection formula for $f'_!$ and that $\mathcal L_g$ is invertible. Thus, the above identity continues to
\begin{align*}
	\mathcal L_{g'} \tensor g'^* f^! = f'^* \mathcal L_g \tensor f'^! g^* = \mathcal L_{g'} \tensor f'^! g^*.
\end{align*}
Here in the last step we used that $f'^* g^! = g'^! f^*$ which holds by the first case of (iii).
\end{proof}

\begin{proposition}
Let $f\colon Y \to X$ be a bdcs and $p$-cohomologically smooth map in $\vStacksCoeff$. Then for all $\mathcal M, \mathcal N \in \Dqcohri(X, \Lambda)$ the natural morphism
\begin{align*}
	f^* \IHom(\mathcal M, \mathcal N) \isoto \IHom(f^* \mathcal M, f^* \mathcal N)
\end{align*}
is an isomorphism.
\end{proposition}
\begin{proof}
This follows formally from \cref{rslt:upper-shriek-of-IHom-formula} using that $f^! = f^* \tensor \mathcal L$ for some invertible sheaf $\mathcal L$.
\end{proof}

\subsection{\texorpdfstring{$\varphi$}{Phi}-Modules and Local Systems} \label{sec:ri-pi.phi-mod}

We now equip the $\ri^+_X/\pi$-modules studied in the previous subsections with $\varphi$-actions, where $\varphi\colon \ri^+_X/\pi \to \ri^+_X/\pi$ is the Frobenius (assuming that $\pi \divides p$). The resulting $\infty$-category $\DqcohriX X^\varphi$ of almost $\varphi$-modules over $\ri^+_X/\pi$ has a strong connection to the $\infty$-category of overconvergent étale $\Fld_p$-sheaves on the small v-stack $X$. In fact the latter embeds fully faithfully into the former and induces an equivalence of dualizable objects (see \cref{rslt:global-Riemann-Hilbert} below); this is a version of a $p$-torsion Riemann-Hilbert correspondence. This correspondence can be seen as an analog of the Riemann-Hilbert correspondence discussed in \cite{riemann-hilbert-mod-p}. The authors of loc. cit. also work on a Riemann-Hilbert correspondence in rigid geometry, which will be similar to ours but on the coherent side extends the sheaves to a formal model (rather then working with $\ri^+_X/\pi$-modules on the generic fiber). Our Riemann-Hilbert correspondence has many interesting consequences; among others, we can deduce a general version of the Primitive Comparison Theorem, see \cref{rslt:primitive-comparison} below.

We start by introducing $\varphi$-modules over analytic rings. They can be defined as a special case of the $\sigma$-modules studied in \cref{sec:andesc.adic-sigobj}.

\begin{definition}
Let $(V,\mm)$ be an almost setup such that $p = 0$ in $V$ and let $\mathcal A$ be an analytic ring over $(V,\mm)$.
\begin{defenum}
	\item The Frobenius $\varphi\colon \mathcal A \to \mathcal A$ induces the structure of an analytic $\sigma$-ring $(\mathcal A, \varphi)$ over $(V,\mm)$. This construction is functorial in $\mathcal A$ as follows: Let $\AnRing^p_{(-)} \subset \AnRing_{(-)}$ be the full subcategory spanned by the analytic rings lying over an almost setup $(V,\mm)$ with $p = 0$ in $V$. The forgetful functor $\sigobj{\AnRing^p_{(-)}} \to \AnRing^p_{(-)}$ induces an equivalence of $\AnRing^p_{(-)}$ with the full subcategory of $\sigobj{\AnRing^p_{(-)}}$ spanned by the objects $(\mathcal A, \varphi)$. Inverting this equivalence yields a fully faithful functor
	\begin{align*}
		\AnRing^p_{(-)} \injto \sigobj{\AnRing^p_{(-)}}, \qquad \mathcal A \mapsto (\mathcal A, \varphi).
	\end{align*}
	We say that $\mathcal A$ is \emph{perfect} if $(\mathcal A, \varphi)$ is perfect. Since the induced morphism $\varphi_V\colon (V,\mm) \to (V,\mm)$ is always strict (by \cite[Proposition 2.1.7.(ii)]{almost-ring-theory}), $\mathcal A$ is perfect if and only if the map $\varphi\colon \mathcal A \to \varphi_{V*} \mathcal A$ is an isomorphism.

	\item A \emph{$\varphi$-module over $\mathcal A$} is an object of $\D(\mathcal A)^\varphi$. A \emph{lax $\varphi$-module over $\mathcal A$} is an object of $\D(\mathcal A[T_\varphi])$.
\end{defenum}
\end{definition}

We first study $\varphi$-modules in the local setting, i.e. on affinoid perfectoid spaces. The main result of this investigation is a local Riemann-Hilbert correspondence (see \cref{rslt:local-Riemann-Hilbert}), which lies at the heart of the proof of the global Riemann-Hilbert correspondence. In the local Riemann-Hilbert correspondence we want to relate the $\infty$-category $\D(\Fld_p(X))_\omega$ of discrete $\Fld_p(X)$-modules (where $\Fld_p(X)$ is the ring of continuous maps $\abs X \to \Fld_p$) to the $\infty$-category $\Dqcohri(A^+/\pi)^\varphi$ of $\varphi$-modules over $(A^+/\pi)^a_\solid$, for any strictly totally disconnected space $X = \Spa(A, A^+)$ of characteristic $p$ and any pseudouniformizer $\pi \in A^+$. The general strategy is as follows:
\begin{enumerate}[(1)]
	\item Given a $\varphi$-module $M \in \Dqcohri(A^+/\pi)^\varphi$, we can use the $\varphi$-structure on $M$ to ``lift'' $M$ to an object in $\Dqcohri(A^+)^\varphi$. Roughly this operation takes $M$ to $\varprojlim_{\varphi_M} M$. For example, if $M = A^{+a}/\pi$ then the associated $\varphi$-module in $\Dqcohri(A^+)^\varphi$ is $A^{+a}$.

	\item The functor $\D(\Fld_p(X))_\omega \to \Dqcohri(A^+/\pi)$ is given by $- \tensor (A^{+a}/\pi)$ and it admits the right adjoint $(-)^\varphi_\omega$ (cf. \cref{def:sigma-invariants-of-sigma-modules}). By the previous step $(-)^\varphi$ factors over $\Dqcohri(A^+)^\varphi$, which makes the situation more pleasent as $A^+$ is perfect. However, the functor $(-)^\varphi\colon \Dqcohri(A^+)^\varphi \to \D(\Fld_p(X))$ is still a bit more subtle than one might expect, because it entails the functor $(-)_*\colon \Dqcohri(A^+) \to \D_\solid(A^+)$, passing from almost modules to honest modules. The functor $(-)_*$ does not preserves colimits, so it is a priori not clear that $(-)^\varphi$ does. Here the following observation comes to the rescue: We have $M^\varphi = (M \tensor_{A^+} A)^\varphi$. This is helpful because $\pi$ is invertible in $A$ so that the almost world over $A$ collapses to the non-almost world.

	\item The previous step will eventually allow us to reduce the local Riemann-Hilbert correspondence to an understanding of the relation between $\D(\Fld_p(X))_\omega$ and $\D(A_\omega)^\varphi_\omega$. After reducing to connected components (so that $A = C$ is an algebraically closed field) we end up in the classical situation of comparing $\varphi$-modules on $C$-vector spaces to $\Fld_p$-vector spaces.
\end{enumerate}
In the following we carry out the above agenda. Step (1) is handled by the following results:

\begin{lemma} \label{rslt:phi-modules-over-A+-mod-pi-equiv-pi-complete-phi-modules}
Let $\Spa(A, A^+)$ be an affinoid perfectoid space of characteristic $p$ and let $\pi \in A^+$ be a pseudouniformizer. Then the functor
\begin{align*}
	- \tensor_{(A^+)^a_\solid} (A^+/\pi)^a_\solid\colon \D^a_\solid(A^+)^\varphi_{\hat\pi} \isoto \D^a_\solid(A^+/\pi)^\varphi
\end{align*}
is an equivalence of $\infty$-categories. Here $\D^a_\solid(A^+)^\varphi_{\hat\pi} \subset \D^a_\solid(A^+)^\varphi$ denotes the full subcategory of those $\varphi$-modules whose underlying $A^+$-module is $\pi$-adically complete.
\end{lemma}
\begin{proof}
We first claim that $\D^a_\solid(A^+)^\varphi_{\hat\pi} = \varprojlim_n \D^a_\solid(A^+/\pi^n)^\varphi$ via the natural functor. Namely, for lax $\varphi$-modules this can be checked on underlying modules (where all pushforward functors behave as expected) and thus follows from \cref{rslt:aidc-complete-modules-equiv-limit-over-A-mod-I-n}. It then remains to check that if $M \in \D((A^+)^a_\solid[T_\varphi])$ satisfies $M \tensor_{(A^+)^a_\solid} (A^+/\pi^n)^a_\solid \in \D^a_\solid(A^+/\pi^n)^\varphi$ then $M \in \D^a_\solid(A^+)^\varphi$. But this follows easily from $M = \varprojlim_n M/\pi^n$.

It is now enough to show that all the transition functors $\D^a_\solid(A^+/\pi^{n+1})^\varphi \to \D^a_\solid(A^+/\pi^n)^\varphi$ are equivalences. In fact it is enough to check that the base-change functor $\D^a_\solid(A^+/\pi^p)^\varphi \to \D^a_\solid(A^+/\pi)^\varphi$ is an equivalence (then replace $\pi$ by $\pi^p$ and repeat). Using that the Frobenius on $A^+/\pi^p$ factors as $A^+/\pi^p \to A^+/\pi \xto{\varphi} A^+/\pi^p$, this follows immediately from \cref{rslt:alpha-beta-sigma-invariants}.
\end{proof}

\begin{lemma} \label{rslt:dualizable-pi-complete-phi-module-is-bounded}
Let $X = \Spa(A, A^+)$ be an affinoid perfectoid space of characteristic $p$ with pseudouniformizer $\pi \in A^+$. Then every dualizable object in $\D^a_\solid(A^+)^\varphi_{\hat\pi}$ (equipped with the $\pi$-completed tensor product) is bounded and dualizable as an object in $\D^a_\solid(A^+)^\varphi$.
\end{lemma}
\begin{proof}
Let $P \in \D^a_\solid(A^+)^\varphi_{\hat\pi}$ be dualizable. We only need to show that $P$ is bounded; then it is automatically dualizable in $\D^a_\solid(A^+)^\varphi$ because its dual $P^*$ is then also bounded and hence by \cref{rslt:solid-tensor-product-preserves-adic-complete} the $\pi$-completed tensor product $P \hat\tensor P^*$ coincides with the solid tensor product $P \tensor P^*$ in $\D^a_\solid(A^+)^\varphi$.

Since $P$ is dualizable, so is $P/\pi \in \D^a_\solid(A^+/\pi)$. By \cref{rslt:equiv-of-compact-perfect-dualizable-over-Lambda} there is some $N \ge 0$ such that for all integers $n$ with $\abs n \ge N$, $\pi_n(P/\pi)$ is killed by $\pi^{1/p}$. By considering the long exact homology sequence associated to the exact triangle $P \xto\pi P \to P/\pi$ in $\D^a_\solid(A^+)$ it follows that for $\abs n > N$ both kernel and cokernel of the map $\pi_n(P) \xto\pi \pi_n(P)$ are annihilated by $\pi^{1/p}$. We claim that this implies that this map is actually an isomorphism. Namely, let $P_n := \pi_0(\pi_n(P)_*)$, which is a static $\varphi$-module over $A^+_\solid$ such that $P_n^a = \pi_n(P)$. Choose any $x \in P_n$ (more precisely, any $x \in P_n(S)$ for any profinite set $S$) such that $\pi x = 0$. Then $\pi^{1/p} x = 0$ by the above observation, hence $\pi \varphi_P(x) = \varphi_P(\pi^{1/p} x) = 0$. But then $\varphi_P(x)$ lies in the kernel of multiplication by $\pi$ and thus $\pi^{1/p} \varphi_P(x) = 0$. Applying $\varphi_P^{-1}$ we deduce $\pi^{1/p^2} x = 0$. Thus $\pi \varphi_P^2(x) = 0$, which implies $\pi^{1/p} \varphi_P^2(x) = 0$ and therefore $\pi^{1/p^3} x = 0$. Repeating this argument we deduce that $\varepsilon x = 0$ for all $\varepsilon \in \mm_A$, hence $x = 0$. This shows that multiplication by $\pi$ is injective on $P_n$. We can now similarly show that it is surjective: Given any $y \in P_n$, we can find some $x_0 \in P_n$ such that $\pi x_0 = \pi^{1/p} y$. Then $\pi^{1/p^2} \varphi_P^{-1}(y) = \pi^{1/p} \varphi_P^{-1}(x_0)$, hence there is some $x_1 \in P_n$ such that $\pi x_1 = \pi^{1/p^2} \varphi_P^{-1}(y)$. It follows that $\pi^p \varphi_P(x_1) = \pi^{1/p} y$ and since multiplication by $\pi$ is injective on $P_n$ we deduce $y = \pi \cdot (\pi^{p - 1 - 1/p} \varphi_P(x_1))$, i.e. $y$ lies in the image of multiplication by $\pi$.

We have now shown that for $\abs n \gg 0$ the map $\pi_n(P) \xto\pi \pi_n(P)$ is an isomorphism. But by \cref{rslt:M-is-I-adically-complete-iff-all-pi-n-M-are-so} $\pi_n(P)$ is $\pi$-adically complete, so by definition we have $0 = \varprojlim(\dots \xto\pi \pi_n(P) \xto\pi \pi_n(P) \xto\pi \pi_n(P))$. Thus $\pi_n(P) = 0$, as desired.
\end{proof}

We now come to step (2) of the above proof strategy. The following result shows the promised $M^\varphi = (M \tensor_{A^+} A)^\varphi$, while the subsequent result deduces that the functor $(-)^\varphi\colon \Dqcohri(A^+/\pi)^\varphi \to \D_\solid(\Fld_p(X))$ preserves small colimits.

\begin{lemma} \label{rslt:phi-invariants-on-A+-equiv-phi-invariants-of-tensor-A}
Let $X = \Spa(A, A^+)$ be an affinoid perfectoid space of characteristic $p$ and let $R := \Fld_p(X)$ (i.e. the ring of continuous maps $\abs X \to \Fld_p$). Then for all $M \in \D^a_\solid(A^+)^\varphi$ the natural morphism
\begin{align*}
	M^\varphi \isoto (M \tensor_{(A^+)^a_\solid} (A, A^+)_\solid)^\varphi
\end{align*}
is an isomorphism in $\D_\solid(R)$.
\end{lemma}
\begin{proof}
Both sides commute with the fully faithful inclusion from the solid theory to all topological modules, so it is enough to prove the claimed isomorphism in the latter setting. Moreover, both sides can be computed on sections on extremally disconnected sets $S$, so we reduce to showing the claim in the discrete setting: Given $M \in \D(A^{+a}_\omega)_\omega$ we need to see that the natural map $M^\varphi \isoto (M \tensor_{A^{+a}_\omega} A_\omega)^\varphi$ is an isomorphism in $\D(R)_\omega$. From now on we drop the subscript $\omega$ and simply assume that everything is discrete.

Note that $\D(A^{+a})^\varphi$ has a left-complete $t$-structure given by the $t$-structure on $\D(A^{+a})$. In fact, by \cref{rslt:sigma-modules-equiv-modules-over-A-T-sigma-pm-if-A-perfect} it can be identified with $\D(A^{+a}[T_\varphi^{\pm}])$, where $A^{+a}[T_\varphi^\pm]$ is the associative $A^{+a}$-algebra given by the rule $T_\varphi a = \varphi(a) T_\varphi$. Both of the described functors commute with the limit $M = \varprojlim_n \tau_{\le n} M$, so we can assume that $M$ is left-bounded. By the above description of $\D(A^{+a})^\varphi$ we note that this $\infty$-category is the derived $\infty$-category of its heart, which itself is a Grothendieck abelian category. In particular we can resolve $M$ by injective objects to reduce to the case that $M$ itself lies in the heart and is injective. In this case the $R$-module $M^\varphi$ is also concentrated in degree $0$ (the functor $(-)^\varphi$ is a derived functor); it can be described explicitly by
\begin{align*}
	M^\varphi = \{ x \in M_* \setst \varphi_M(x) = x \},
\end{align*}
where $\varphi_M\colon M \to \varphi_* M$ is the $\varphi$-linear map belonging to the $\varphi$-module $M$. We now observe that $(M \tensor_{A^{+a}} A)^\varphi$ is also concentrated in degree $0$: It is computed as $(M_* \tensor_{A^+} A^+[1/\pi])^\varphi$ for any pseudouniformizer $\pi \in A^+$, so the claim follows from \cite[Lemma 7.1.2]{riemann-hilbert-mod-p}. We get
\begin{align*}
	(M \tensor_{A^{+a}} A)^\varphi = \{ a \tensor x \in A \tensor_{A^+} M_* \setst a^p \tensor \varphi_M(x) = a \tensor x \}.
\end{align*}
It remains to check that the map of classical $R$-modules $\alpha\colon M^\varphi \to (M \tensor_{A^{+a}} A)^\varphi$ is an isomorphism, which amounts to showing that it is both injective and surjective. Let us also denote $\alpha'\colon M_* \to (M_* \tensor_{A^+} A)$ the natural map (so that $\alpha$ is the restriction of $\alpha'$) and $\varphi_M'$ the Frobenius action on $M \tensor_{A^+} A$.

We first show that $\alpha$ is injective. Suppose we have $x \in M^\varphi$ such that $\alpha(x) = 0$. Then there is some pseudouniformizer $\pi \in A^+$ such that $\pi x = 0$ (in $M_*$). But then also $0 = \varphi_M^{-1}(\pi x) = \pi^{1/p} x$. We inductively deduce that $\varepsilon x = 0$ for all $\varepsilon \in \mm_A$. But this implies $x = 0$ in $M_* = \IHom_{A^+}(\mm_A, M_*)$, as desired.

It remains to show that $\alpha$ is surjective, so suppose we are given some $y \in (M \tensor_{A^+} A)^\varphi$. Then there is some pseudouniformizer $\pi \in A^+$ and some $x' \in M_*$ such that $y = \pi^{-1} \alpha'(x')$. We deduce $\alpha'(\varphi_M(x') - \pi^{p-1} x') = \varphi_{M'}(\pi y) - \pi^p y = \pi^p(\varphi_{M'}(y) - y) = 0$, hence there is some pseudouniformizer $\varpi \in A^+$ such that $\pi \divides \varpi$ and $\frac{\varpi^p}{\pi^p} (\varphi_M(x') - \pi^{p-1} x') = 0$. Let us rewrite the latter condition as
\begin{align*}
	\frac\varpi\pi x' = \varphi_M^{-1}(\frac{\varpi^p}\pi x').
\end{align*}
We now define an element $x \in M_* = \IHom_{A^+}(\mm_A, M_*)$ by requiring
\begin{align*}
	x(\varpi_n) = \varphi_M^{-n}(\frac\varpi\pi x')
\end{align*}
for all $n \ge 0$, where $\varpi_n := \varpi^{1/p^n}$ For this to define an element of $M_*$ we need to verify the following compatibility: For all $n \ge 0$ we have
\begin{align*}
	&\frac{\varpi_n}{\varpi_{n+1}} x(\varpi_{n+1}) = \frac{\varpi_n}{\varpi_{n+1}} \varphi_M^{-n-1}(\frac\varpi\pi x') = \varphi_M^{-n-1}(\frac{\varpi^p}{\varpi} \frac\varpi\pi x') = \varphi_M^{-n}(\varphi_M^{-1}(\frac{\varpi^p}\pi x')) =\\&\qquad= \varphi_M^{-n}(\frac\varpi\pi x') = x(\varpi_n).
\end{align*}
It is clear that $x \in M_*$ is $\varphi_M$-invariant: By definition $\varphi_M$ acts on $x$ by $(\varphi_M(x))(\varepsilon) = \varphi_M(x(\varepsilon^{1/p}))$ for all $\varepsilon \in \mm_A$; thus for $\varepsilon = \varpi_n$ we have $(\varphi_M(x))(\varepsilon) = x(\varepsilon)$ and for all other $\varepsilon$ it follows from that case. One also checks immediately that $\alpha(x) = y$, as desired.
\end{proof}

\begin{lemma} \label{rslt:phi-invariants-on-aff-perf-space-preserve-colim}
Let $X = \Spa(A, A^+)$ be an affinoid perfectoid space of characteristic $p$ and let $R := \Fld_p(X)$. Then the functor
\begin{align*}
	(-)^\varphi\colon \D^a_\solid(A^+/\pi)^\varphi \to \D_\solid(R), \qquad M \mapsto M^\varphi
\end{align*}
preserves all small colimits.
\end{lemma}
\begin{proof}
It is clear that $(-)^\varphi\colon \D_\solid(A^+/\pi)^\varphi \to \D_\solid(R)$ preserves all small colimits, since it is given by an equalizer and hence a finite limit. If we work with almost coefficients instead then it is less clear that $(-)^\varphi$ preserves small colimits because it factors as $\D^a_\solid(A^+/\pi)^\varphi \xto{(-)_*} \D_\solid(A^+/\pi)^\varphi \xto{(-)^\varphi} \D_\solid(R)$ and the functor $M \mapsto M_*$ does not preserve all small colimits. To circumvent this issue we make use of \cref{rslt:phi-invariants-on-A+-equiv-phi-invariants-of-tensor-A} which allows us to move from the almost world to the non-almost world in a colimit-preserving way.

More concretely, we can argue as follows: Let $(M_i)_i$ be a small family of objects in $\D^a_\solid(A^+/\pi)^\varphi$. We need to show that $(\bigdsum_i M_i)^\varphi = \bigdsum_i M_i^\varphi$. By \cref{rslt:phi-modules-over-A+-mod-pi-equiv-pi-complete-phi-modules} we can lift $M_i$ to $M'_i \in \D^a_\solid(A^+)^\varphi$ and since $(-)^\varphi\colon \D^a_\solid(A^+) \to \D_\solid(R)$ preserves all small colimits by \cref{rslt:phi-invariants-on-A+-equiv-phi-invariants-of-tensor-A}, the claimed identity of $\varphi$-invariants is equivalent to the identity $(\bigdsum_i M'_i)^\varphi = (\widehat\bigdsum_i M'_i)^\varphi$ (where $\widehat\bigdsum$ denotes the $\pi$-completed direct sum). Let us denote $N := \fib(\bigdsum_i M'_i \to \widehat\bigdsum_i M'_i) \in \D^a_\solid(A^+)^\varphi$. Then the claimed identity on $\varphi$-invariants is equivalent to $N^\varphi = 0$. Now reverse the above equivalences in the non-almost world, with $M_{i*} \in \D_\solid(A^+/\pi)^\varphi$ in place of $M_i$, $M'_{i*} \in \D_\solid(A^+)^\varphi$ in place of $M'_i$ and $N' := \fib(\bigdsum_i M'_{i*} \to \widehat\bigdsum_i M'_{i*})$ in place of $N$ (note that \cref{rslt:phi-modules-over-A+-mod-pi-equiv-pi-complete-phi-modules} holds in the non-almost setting by the same proof). From the fact that $(-)^\varphi\colon \D_\solid(A^+/\pi)^\varphi \to \D_\solid(R)$ preserves all small colimits we can therefore deduce that $N'^\varphi = 0$. There is a natural almost isomorphism $N' \to N_*$. Moreover, by the proof of \cref{rstl:characterizations-of-I-adically-complete-modules} the underlying $(A^+)_\solid$-module of $N_*$ is $\varprojlim(\dots \xto\pi (\bigdsum_i M_i)_* \xto\pi (\bigdsum_i M_i)_*)$ and a similar description holds for $N'$. In particular multiplication by $\pi$ is an isomorphism on $N_*$ and on $N'$ so that the almost isomorphism $N' \to N_*$ is actually an isomorphism. It follows that $N^\varphi = (N_*)^\varphi = N'^\varphi = 0$, as desired.
\end{proof}

We have now proved step (2) of the above agenda. In order for it to be very useful we still need a better understanding of the functor $- \tensor_{A^+} A$ in order to reduce everything to $\varphi$-modules over $A$ (as in step (3)). The following result comes in very handy:

\begin{lemma} \label{rslt:tensor-compact-A+-phi-module-with-A-is-conservative}
Let $\Spa(A, A^+)$ be an affinoid perfectoid space of characteristic $p$. Then the base-change functor
\begin{align*}
	- \tensor_{(A^+)^a_\solid} A_\solid\colon \D^a_\solid(A^+)^{a\omega,\varphi} \to \D_\solid(A, A^+)^\varphi
\end{align*}
is conservative. Here $\D^a_\solid(A^+)^{a\omega,\varphi} \subset \D^a_\solid(A^+)^\varphi$ denotes the full subcategory of those $\varphi$-modules whose underlying $(A^+)^a_\solid$-module is almost compact.
\end{lemma}
\begin{proof}
Fix a pseudouniformizer $\pi \in A^+$ and any $P \in \D^a_\solid(A^+)^{a\omega,\varphi}$ such that $P \tensor_{(A^+)^a_\solid} (A, A^+)_\solid = 0$; we need to see that $P = 0$. The condition $P \tensor_{(A^+)^a_\solid} (A, A^+)_\solid = 0$ means that $\varinjlim(P \xto\pi P \xto\pi P \xto\pi \dots) = 0$. Since $P$ is almost compact this implies that
\begin{align*}
	\varinjlim_\pi \pi_0 \IHom(P, P)_\omega = \pi_0 \IHom(P, \varinjlim_\pi P)_\omega = 0,
\end{align*}
from which we deduce that $\varinjlim_\pi \pi_0 \Hom(P, P)$ is almost zero (where we view $\Hom$ as enriched over $A^+_\omega$). This means that for every $f \in \pi_0 \Hom(P, P)$ there is some $n \ge 0$ such that $\pi^n f = 0$. Applying this to $f = \id$ we deduce that the map $P \xto{\pi^n} P$ is zero. Via the isomorphism $\varphi^* P \isoto P$ we deduce that multiplication by $\pi^{n/p}$ is also zero on $P$. Inductively it follows that multiplication by any $\varepsilon \in \mm_A$ is zero on $P$, hence $P = 0$.
\end{proof}

Before we can come to the proof of the local Riemann-Hilbert correspondence, we observe that the whole situation commutes with passage to connected components:

\begin{lemma} \label{rslt:aff-perf-phi-module-functors-commute-with-stalks}
Let $X = \Spa(A, A^+)$ be a totally disconnected space of characteristic $p$ with pseudouniformizer $\pi$, $x = \Spa(K, K^+) \in \pi_0(X)$ a connected component and $R := \Fld_p(X)$. The map $x \to X$ induces a map $R \to \Fld_p(x) = \Fld_p$ and the following diagrams commute
\begin{center}
	\begin{tikzcd}[column sep = huge]
		\D_\solid(\Fld_p) \arrow[r,"- \tensor (K^+/\pi)^a_\solid"] & \D^a_\solid(K^+/\pi)^\varphi\\
		\D_\solid(R) \arrow[u,"- \tensor \Fld_p"] \arrow[r, "- \tensor (A^+/\pi)^a_\solid"] & \D^a_\solid(A^+/\pi)^\varphi \arrow[u,swap,"- \tensor K^{+a}/\pi"]
	\end{tikzcd}
	\qquad
	\begin{tikzcd}[column sep = large]
		\D_\solid(\Fld_p) & \D^a_\solid(K^+/\pi)^\varphi \arrow[l,swap,"(-)^\varphi"]\\
		\D_\solid(R) \arrow[u,"- \tensor \Fld_p"] & \D^a_\solid(A^+/\pi)^\varphi \arrow[u,swap,"- \tensor K^{+a}/\pi"] \arrow[l,swap,"(-)^\varphi"]
	\end{tikzcd}
\end{center}
\end{lemma}
\begin{proof}
The commutativity of the left diagram is clear (both paths compute $- \tensor_{R_\solid} (K^+/\pi)^a_\solid$, cf. \cref{rslt:trivial-sigma-structure-is-functorial}). It remains to prove that the right diagram commutes, so let $M \in \D^a_\solid(A^+/\pi)^\varphi$ be given. By the commutativity of the left diagram we get a natural morphism $M^\varphi \tensor \Fld_p \to (M \tensor K^{+a}/\pi)^\varphi$ which we claim to be an isomorphism. Now write $x = \bigisect_i U_i$ for a filtered system $(U_i)_i$ of clopen neighborhoods $x \in U_i \subset X$. Write $U_i = \Spa(A_i, A^+_i)$ and $R_i = \Fld_p(U_i)$. Then $A^+_i/\pi$ is a retract of $A^+/\pi$ and $R_i$ is a retract of $R$ in a compatible way, which easily implies that the map $M^\varphi \tensor R_i \isoto (M \tensor A^{+a}/\pi)^\varphi$ is an isomorphism for all $i$. Note that $\Fld_p = \varinjlim_i R_i$ and $K^+/\pi = \varinjlim_i A^+_i/\pi$, so since $(-)^\varphi$ preserves all small colimits (see \cref{rslt:phi-invariants-on-aff-perf-space-preserve-colim}) we can pull the colimits out of $M^\varphi \tensor \Fld_p \to (M \tensor K^{+a}/\pi)^\varphi$ to conclude that this map is an isomorphism.
\end{proof}

Finally, we are in the position to prove the promised local Riemann-Hilbert correspondence. The combination of the above results leads to the following:

\begin{proposition} \label{rslt:local-Riemann-Hilbert}
Let $X = \Spa(A, A^+)$ be a strictly totally disconnected space of characteristic $p$ with pseudouniformizer $\pi \in A^+$ and denote $R := \Fld_p(X)$. Then the symmetric monoidal functor
\begin{align*}
	- \tensor_R (A^{+a}/\pi)\colon \D(R)_\omega \injto \D^a_\solid(A^+/\pi)^\varphi
\end{align*}
is fully faithful and induces an equivalence of the full subcategories of dualizable objects on both sides.
\end{proposition}
\begin{proof}
The given functor has the right adjoint $(-)^\varphi\colon \D^a_\solid(A^+/\pi)^\varphi \to \D(R)_\omega$ (where we implicitly pass to the discretization). Thus in order to prove the claimed full faithfulness we need to show that for every discrete $R$-module $M \in \D(R)_\omega$ the natural map $M \isoto (M \tensor_R A^{+a}/\pi)^\varphi$ is an isomorphism. This can be checked on stalks on $\Spec R$. Now $R$ is the ring of continuous maps from the profinite set $\pi_0(X)$ to $\Fld_p$, so the the topological space $\Spec R$ can naturally be identified with $\pi_0(X)$ (e.g. write $\pi_0(X) = \varprojlim S_i$ for finite $S_i$, so that $\Spec R = \varprojlim_i \Spec \Fld_p(S_i) = \varprojlim_i S_i = \pi_0(X)$). By \cref{rslt:aff-perf-phi-module-functors-commute-with-stalks} $(M \tensor_R A^{+a}/\pi)^\varphi$ commutes with passage to stalks, so we can reduce to the case that $X$ is connected and thus of the form $X = \Spa(C, C^+)$ for some algebraically closed field $C$ and $R = \Fld_p$. Writing $M \in \D(\Fld_p)_\omega$ as a colimit of copies of $\Fld_p$ and using that $(-)^\varphi$ preserves colimits (see \cref{rslt:phi-invariants-on-aff-perf-space-preserve-colim}) we reduce to the case $M = \Fld_p$. Then $(M \tensor (C^{+a}/\pi))^\varphi = (C^{+a}/\pi)^\varphi$, which we need to show to be isomorphic to $\Fld_p$. These $\varphi$-invariants can also be computed after applying the right adjoint of $\D^a_\solid(C^+)^\varphi \to \D^a_\solid(C^+/\pi)^\varphi$, which by \cref{rslt:phi-modules-over-A+-mod-pi-equiv-pi-complete-phi-modules} maps $C^{+a}/\pi$ to $C^{+a}$. By \cref{rslt:phi-invariants-on-A+-equiv-phi-invariants-of-tensor-A} we have $(C^{+a})^\varphi = C^\varphi$. We thus need to check that $C^\varphi = \Fld_p$. We obviously have $\pi_0 C^\varphi = \Fld_p$ (these are just the $\varphi$-invariant elements in $C$), so it remains to verify that $\pi_1 C^\varphi = 0$, i.e. that the map $\varphi - 1\colon C \to C$ is surjective. This follows immediately from the fact that $C$ is algebraically closed. We have thus finished the proof that $- \tensor_R (A^{+a}/\pi)$ is fully faithful.

We now prove the second part of the claim, i.e. that $- \tensor_R (A^{+a}/\pi)$ induces an equivalence of dualizable objects on both sides. By the full faithfulness of this functor it only remains to check that its essential image contains all dualizable objects of $\D^a_\solid(A^+/\pi)^\varphi$. Fix a dualizable object $P \in \D^a_\solid(A^+/\pi)^\varphi$. We need to see that the natural map $P^\varphi \tensor_R (A^{+a}/\pi) \to P$ is an isomorphism. This can be checked after pullback to the connected components of $\Spec A^+/\pi$ (see \cref{rslt:compare-catsldmod-with-sheaves-on-pi-0}) which correspond to $\pi_0(X)$. By \cref{rslt:aff-perf-phi-module-functors-commute-with-stalks} we can thus reduce to the case that $X = \Spa(C, C^+)$ is connected (and hence $R = \Fld_p$), as in the proof of full faithfulness. By \cref{rslt:phi-modules-over-A+-mod-pi-equiv-pi-complete-phi-modules} we can identify $P$ with a dualizable object $P' \in \D^a_\solid(C^+)^\varphi_{\hat\pi}$ (where the tensor product is the $\pi$-completed tensor product in $\D^a_\solid(C^+)^\varphi$). Then by \cref{rslt:dualizable-pi-complete-phi-module-is-bounded} $P'$ is bounded and dualizable as an object in $\D^a_\solid(C^+)^\varphi$ (where the tensor product is the solid tensor product). Moreover, it follows immediately from the definitions that $P'^\varphi = P^\varphi$, so it is now enough to show that the natural map $P'^\varphi \tensor_R C^{+a} \to P'$ is an isomorphism. By \cref{rslt:tensor-compact-A+-phi-module-with-A-is-conservative} this boils down to showing that $P'^\varphi$ is dualizable (and hence $P'^\varphi \tensor_R C^{+a}$ is almost compact) and that for $P'' := P' \tensor_{C^{+a}} C \in \D_\solid(C, C^+)^\varphi$ the natural map $P''^\varphi \tensor_R C \to P''$ is an isomorphism. Clearly $P''$ is dualizable in $\D_\solid(C, C^+)^\varphi$ and by \cref{rslt:phi-invariants-on-A+-equiv-phi-invariants-of-tensor-A} we have $P'^\varphi = P''^\varphi$, so we are reduced to showing the following claim: The adjoint pair of functors
\begin{align*}
	- \tensor_R C\colon \D(R)_\omega \rightleftarrows \D_\solid(C, C^+)^\varphi \noloc (-)^\varphi
\end{align*}
restricts to quasi-inverse equivalences on the full subcategories of dualizable objects. By \cite[Theorems 5.9, 5.50]{andreychev-condensed-huber-pairs} the symmetric monoidal $\infty$-category of dualizable objects in $\D_\solid(C, C^+)$ is naturally equivalent to the symmetric monoidal $\infty$-category of dualizable objects in $\D(C_\omega)_\omega$ (and hence the same follows for $\varphi$-modules in the respective $\infty$-categories by \cref{rslt:sigma-module-dualizable-iff-underlying-module-is-so}), so we can replace $\D_\solid(C, C^+)^\varphi$ by $\D(C_\omega)^\varphi_\omega$ in the above claim. Given any dualizable $Q \in \D(C_\omega)^\varphi_\omega$ then the underlying $C_\omega$-module is perfect and hence a finite direct sum of copies of $C_\omega[n]$ for varying $n \in \Z$. We need to verify that such a decomposition also exists as $\varphi$-module. This boils down to showing that every $\varphi$-module structure on a finite-dimensional $C_\omega$-vector space admits a basis of $\varphi$-fixed vectors. This is a classical result, see e.g. \cite[Exposé XXII, Proposition 1.1]{sga-7}; alternatively, one can deduce the equivalence of dualizable objects in $\D(R)_\omega$ and $\D(C_\omega)_\omega^\varphi$ directly from \cite[Corollary 12.1.7, Theorem 12.4.1]{riemann-hilbert-mod-p}.
\end{proof}

Having understood the Riemann-Hilbert correspondence in the local setting, we now want to globalize it in order to get a similar result on every small v-stack $X$. This is a formal gluing procedure, but in order for the obtained result to be really useful we need to better understand the involved $\infty$-categories of sheaves on $X$. We start by introducing $\varphi$-modules over $\ri^+_X/\pi$.

\begin{lemma} \label{rslt:Frobenius-is-morphism-of-integral-torsion-coefficients}
Let $X^\sharp$ be an untilted small v-stack with pseudouniformizer $\pi$ such that $p \in \pi$. Then the Frobenius $\varphi\colon \ri^+_{X^\sharp}/\pi \to \ri^+_{X^\sharp}/\pi$ is a morphism of integral torsion coefficients and in particular induces a symmetric monoidal functor
\begin{align*}
	\varphi^*\colon \DqcohriX{X^\sharp} \to \DqcohriX{X^\sharp}.
\end{align*}
\end{lemma}
\begin{proof}
By \cite[Proposition 2.1.7.(ii)]{almost-ring-theory} for every almost setup $(V,\mm)$ with $p = 0$ in $V$ the Frobenius $\varphi\colon V \to V$ defines a strict morphism of almost setups $(V,\mm) \to (V,\mm)$. Thus using \cref{rslt:A+-mod-pi-equiv-ri-+-mod-pi} it only remains to show that for every map $\Spa(B, B^+) \to \Spa(A, A^+)$ of affinoid perfectoid spaces with pseudouniformizer $\pi \in A^+$ such that $\pi \divides p$ the diagram
\begin{center}\begin{tikzcd}
	(A^+/\pi)^a_\solid \arrow[d,"\varphi"] \arrow[r] & (B^+/\pi)^a_\solid \arrow[d,"\varphi"]\\
	(A^+/\pi)^a_\solid \arrow[r] & (B^+/\pi)^a_\solid
\end{tikzcd}\end{center}
is a pushout square of analytic rings. This question only depends on the tilts of $\Spa(A, A^+)$ and $\Spa(B, B^+)$, so we can assume that both are of characteristic $p$. Then $\varphi\colon (A^+/\pi)^a_\solid \to (A^+/\pi)^a_\solid$ factors as the projection $(A^+/\pi)^a_\solid \to (A^+/\pi^{1/p})^a_\solid$ followed by an isomorphism $(A^+/\pi^{1/p})^a_\solid \isoto (A^+/\pi)^a_\solid$ (and similarly for $B^+$ in place of $A^+$). Clearly the isomorphisms form a pushout square, so it remains to check that
\begin{align*}
	(A^+/\pi^{1/p})^a_\solid \tensor_{(A^+/\pi)^a_\solid} (B^+/\pi)^a_\solid = (B^+/\pi^{1/p})^a_\solid.
\end{align*}
By \cref{rslt:colimits-of-solid-discrete-almost-rings} this follows from the observation that $A^+/\pi^{1/p} \tensor_{A^+/\pi} B^+/\pi = B^+/\pi^{1/p}$ (in the derived sense, as always).
\end{proof}

\begin{definition}
\begin{defenum}
	\item We denote $\vStackspip \subset \vStacksCoeff$ the full subcategory of those pairs $(X, \Lambda)$ where $\Lambda = \ri^+_{X^\sharp}/\pi$ for some untilt $X^\sharp$ of $X$ and some pseudouniformizer $\pi$ on $X^\sharp$ such that $p \in \pi$. We denote the objects of $\vStackspip$ by $(X^\sharp, \pi)$ or often simply by $X$ (with the untilt and $\pi$ being implicit).

	\item \label{def:phi-module-in-DqcohriX} Let $X \in \vStackspip$. A \emph{$\varphi$-module over $\ri^+_X/\pi$} is an object $\mathcal M \in \DqcohriX X$ together with an isomorphism $\varphi^* \mathcal M \isoto \mathcal M$ (here $\varphi^*$ comes from \cref{rslt:Frobenius-is-morphism-of-integral-torsion-coefficients}). We denote
	\begin{align*}
		\DqcohriX X^\varphi := \equalizer{\DqcohriX X}{\DqcohriX X}\id{\varphi^*}
	\end{align*}
	the symmetric monoidal $\infty$-category of $\varphi$-modules over $\ri^+_X/\pi$.
\end{defenum}
\end{definition}

\begin{remark} \label{rslt:phi-modules-form-hypercomplete-v-sheaf}
By commuting limits with limits one sees immediately that the assignment $X \mapsto \DqcohriX X^\varphi$ defines a hypercomplete v-sheaf of symmetric monoidal $\infty$-categories on $\vStackspip$. In fact this sheaf is just the equalizer of the two morphism $\id, \varphi^*\colon \DqcohriX{(-)} \to \DqcohriX{(-)}$ of sheaves of symmetric monoidal $\infty$-categories.
\end{remark}

Before we continue, let us investigate how the six functors from our 6-functor formalism extend to $\varphi$-modules. In fact, this works as nicely as one could hope:

\begin{lemma} \label{rslt:compatibility-of-phi-6-functors-with-forgetful-functor}
Let $f\colon Y \to X$ be a morphism in $\vStackspip$.
\begin{lemenum}
	\item The forgetful functor $\DqcohriX X^\varphi \to \DqcohriX X$ is symmetric monoidal, conservative and preserves all small limits, colimits and internal homs.

	\item There is a natural symmetric monoidal pullback functor $f^*\colon \DqcohriX X^\varphi \to \DqcohriX Y^\varphi$ which acts as $f^*$ on the underlying $\ri^+_X/\pi$-modules.

	\item \label{rslt:pushforward-of-phi-modules} The functor from (ii) admits a right adjoint $f_*\colon \DqcohriX Y^\varphi \to \DqcohriX X^\varphi$ which acts as $f_*$ on the underlying $\ri^+_Y/\pi$-modules.

	\item If $f$ is bdcs then there is a natural functor $f_!\colon \DqcohriX Y^\varphi \to \DqcohriX X^\varphi$ which acts as $f_!$ on the underlying $\ri^+_Y/\pi$-modules.

	\item \label{rslt:upper-shriek-of-phi-modules} If $f$ is bdcs then the functor from (iv) admits a right adjoint $f^!\colon \DqcohriX X^\varphi \to \DqcohriX Y^\varphi$ which acts as $f^!$ on the underlying $\ri^+_X/\pi$-modules.
\end{lemenum}
\end{lemma}
\begin{proof}
We start the proof with some general observations. Let $Y_\bullet \to X$ be a hypercover of $X$ such that each $Y_n = \bigdunion_{i\in I_n} Y_{n,i}$ is a disjoint union of totally disconnected spaces $Y_{n,i} = \Spa(B_{n,i}, B^+_{n,i})$. Let $\DqcohriX{Y_\bullet} \to \Delta$ be the coCartesian fibration classifying the functor $\Delta \to \infcatinf$, $n \mapsto \DqcohriX{Y_n}$ and let $\mathcal C$ be the $\infty$-category of sections $\Delta \to \DqcohriX{Y_\bullet}$. Then $\DqcohriX X$ is the full subcategory of $\mathcal C$ spanned by the coCartesian sections. Let $r\colon \mathcal C \to \DqcohriX X$ be a right adjoint to the inclusion (it exists by the proof of \cref{rslt:properties-of-kappa-condensed-ri+-modules}). We can similarly construct $\DqcohriX{Y_\bullet}^\varphi$ and the $\infty$-category $\mathcal C^\varphi$ of sections $\Delta \to \DqcohriX{Y_\bullet}^\varphi$, together with the right adjoint $r^\varphi\colon \mathcal C^\varphi \to \DqcohriX X^\varphi$ of the inclusion. Here we can identify $\mathcal C^\varphi$ with the equalizer of the two maps $\id, \varphi^*\colon \mathcal C \to \mathcal C$, so in particular there is a natural forgetful functor $\mathcal C^\varphi \to \mathcal C$. We claim that the diagram
\begin{center}\begin{tikzcd}
	\mathcal C^\varphi \arrow[r,"r^\varphi"] \arrow[d] & \DqcohriX X^\varphi \arrow[d]\\
	\mathcal C \arrow[r,"r"] & \DqcohriX X
\end{tikzcd}\end{center}
commutes, where the vertical maps are the forgetful functors. To see this, we first observe that each forgetful functor $\Dqcohri(B_{n,i}^+/\pi)^\varphi \to \Dqcohri(B_{n,i}^+/\pi)$ admits a left adjoint $\ell$; namely, this essentially boils down to showing that this functor preserve small limits, which follows easily from the comparison with $\pi$-adically complete $B_{n,i}^{+a}$-modules (see \cref{rslt:phi-modules-over-A+-mod-pi-equiv-pi-complete-phi-modules}) on which the forgetful functor evidently preserves limits (see \cref{rslt:properties-of-sigma-module-forgetful-functors}). Thus the commutativity of the above diagram can be checked after passing to left adjoints, which reduces the claim to checking that the functors $\ell$ preserve coCartesian edges in $\DqcohriX{Y_\bullet}$. This boils down to showing that for a map $g\colon \Spa(A', A'^+) \to \Spa(A, A^+)$ of totally disconnected spaces in $\vStackspip$ the diagram
\begin{center}\begin{tikzcd}
	\Dqcohri(A^+/\pi)^\varphi \arrow[r,"g^*"] & \Dqcohri(A'^+/\pi)^\varphi\\
	\Dqcohri(A^+/\pi) \arrow[u,"\ell"] \arrow[r,"g^*"] & \Dqcohri(A'^+/\pi) \arrow[u,"\ell"]
\end{tikzcd}\end{center}
commutes. Passing to right adjoint in this diagram, the claim reduces to showing that the forgetful functor from $\varphi$-modules to the underlying modules commutes with the forgetful functor $\Dqcohri(A'^+/\pi) \to \Dqcohri(A^+/\pi)$. This in turn follows from the fact that the latter forgetful functor commutes with $\varphi^*$ by \cref{rslt:Frobenius-is-morphism-of-integral-torsion-coefficients}. This finally proves that the first diagram above commutes.

We will now prove the actual claims of the lemma. In part (i) everything is obvious except that the forgetful functor preserves small limits and internal homs. To check these two claims, we note that limits and internal homs in $\DqcohriX X^\varphi$ can be computed by first computing them in $\mathcal C^\varphi$ and then applying $r^\varphi$. Now since in $\mathcal C^\varphi$ both limits and internal homs are computed componentwise (in each $\DqcohriX{Y_{n,i}}^\varphi$) and the above diagram commutes, we can reduce the claim to each $Y_{n,i}$, i.e. we can assume that $X = \Spa(A, A^+)$ is totally disconnected. Then by \cref{rslt:phi-modules-over-A+-mod-pi-equiv-pi-complete-phi-modules} we reduce the claim to showing that the forgetful functor $\Dqcohri(A^+)^\varphi_{\hat\pi} \to \Dqcohri(A^+)_{\hat\pi}$ and the base-change functor $\Dqcohri(A^+) \to \Dqcohri(A^+/\pi)$ preserve small limits and internal homs. For the former functor this is true since $A^+$ is perfect (so that $\varphi^*$ is an isomorphism) and for the latter functor it follows because the base-change functor is computed by a simple cofiber sequence $M \mapsto \cofib(M \xto\pi M)$.

Part (ii) is obvious from the definitions (cf. \cref{rslt:phi-modules-form-hypercomplete-v-sheaf}). We now prove (iii). Note first that there is a right adjoint of $f^*$ on the $\varphi$-module $\infty$-categories (this can be constructed as in \cref{rslt:existence-of-pushforward-on-Dqcohri}; up to set-theoretic nonsense this is just the adjoint functor theorem). As in the proof of (i) we can make use of the commuting diagram involving $r$ and $r^\varphi$ above in order to reduce to the case that $X$ is totally disconnected. By choosing a hypercover of $Y$ by disjoint unions of totally disconnected spaces we see that $f_*$ is computed as a combination of forgetful functors $\Dqcohri(B^+/\pi)^\varphi \to \Dqcohri(A^+/\pi)^\varphi$ along maps $\Spa(A, A^+) \to \Spa(B, B^+)$ and a fixed limit (built out of products and a totalization). Both of these operations commute with the forgetful functor from $\varphi$-modules to the underlying modules (as shown in the proof of the commuting diagram above). Thus $f_*$ commutes with the forgetful functor.

We now prove (iv), so assume that $f$ is bdcs. Then $f_!$ on $\varphi$-modules can be constructed in the same way as it was constructed in \cref{sec:ri-pi.6-functor}: Essentially we let $f_! = g_* \comp j_!$ for a (local) factorization of $f$ into a proper map $g$ and an étale map $j$. By (iii) it only remains to show that $j_!$ exists on $\varphi$-modules and commutes with the forgetful functor. But this follows immediately from the fact that $j_!$ commutes with $\varphi^*$ by \cref{rslt:coeff-projection-formula-for-etale-maps} (namely, apply the usual coCartesian section argument to the equalizer in the definition of $\varphi$-modules).

It remains to prove (v), so again assume that $f$ is bdcs. It follows in the same way as in the proof of \cref{rslt:main-6-functor-formalism} that the functor $f_!$ admits a right adjoint $f^!$ on $\varphi$-modules; we have to show that $f^!$ commutes with the forgetful functor from $\varphi$-modules to the underlying modules. By passing to a hypercover of $X$ by disjoint unions of totally disconnected spaces and replacing $Y$ by the pullback hypercover, we see that $f^!$ can be computed componentwise on the hypercover plus an application of the functor $r$ constructed above but with respect to the hypercover of $Y$ and not of $X$. Thus the commuting diagram at the beginning of the proof allows us to reduce the claim to the case that $X = \Spa(A, A^+)$ is totally disconnected. Then since $f$ is bdcs we know that $Y$ is a locally spatial diamond and locally compactifiable over $X$, so we can pass to an open cover of $Y$ in order to reduce to the case that $f$ is quasicompact and compactifiable. Now choose a pro-étale cover $g\colon U = \varprojlim_i U_i \to Y$ such that $U$ is totally disconnected and all $g_i\colon U_i \to Y$ are separated quasicompact étale (see \cite[Proposition 11.24]{etale-cohomology-of-diamonds}). Then $g^*$ is conservative, so it is enough to show that $g^* f^!$ commutes with the forgetful functor from $\varphi$-modules to modules. Moreover, $f_* g_*$ is a forgetful functor and hence conservative, so we reduce to showing that $f_* g_* g^* f^!$ commutes with the forgetful functor from $\varphi$-modules to modules. By \cref{rslt:colim-of-pushforward-pullback} we have $g_* g^* = \varinjlim_i g_{i*} g_i^*$ and by \cref{rslt:qcqs-p-bounded-pushforward-preserves-colimits-and-Dqcohrim} $f_*$ preserves this colimit, so we reduce to showing that $f_* g_{i*} g_i^* f^!$ commutes with the forgetful functor from $\varphi$-modules to modules. By replacing $f$ by $f \comp g_i$ we are left with showing that $f_* f^!$ commutes with the forgetful functor from $\varphi$-modules to modules. By factoring $f$ into an open immersion and a $p$-bounded proper map we see that $f_!$ satisfies the projection formula (even for the $\varphi$-module version). It then follows formally from adjunctions that $f_* f^! = \IHom(f_! \ri^+_Y/\pi, -)$. But we know by (i) and (iv) that the right-hand side commutes with the forgetful functor from $\varphi$-modules to modules, as desired.
\end{proof}

\begin{proposition} \label{rslt:6-functor-formalism-for-phi-modules}
The six functors $\tensor$, $\IHom$, $f^*$, $f_*$, $f_!$ and $f^!$ extend naturally to $\varphi$-modules over $\vStackspip$. These $\varphi$-module versions also form a 6-functor formalism.
\end{proposition}
\begin{proof}
The 6-functor formalism can be constructed in the same way as in \cref{rslt:main-6-functor-formalism}. Here \cref{rslt:compatibility-of-phi-6-functors-with-forgetful-functor} guarantees that the $\varphi$-module version of each of the six functors coincides with the original version after applying the forgetful functor from $\varphi$-modules to modules.
\end{proof}

\begin{remark}
In view of \cref{rslt:6-functor-formalism-for-phi-modules} and the Riemann-Hilbert correspondence (see \cref{rslt:global-Riemann-Hilbert} below) one may view $\DqcohriX X^\varphi$ as the central $\infty$-category of this thesis (instead of the version $\DqcohriX X$ without $\varphi$-module structure), at least if $\pi \divides p$.
\end{remark}

We now single out a special subcategory of $\DqcohriX X^\varphi$: the perfect sheaves.\footnote{In view of the Riemann-Hilbert correspondence we decided to call the dualizable objects in $\DqcohriX X^\varphi$ ``perfect'' and not ``weakly almost perfect''.} Roughly, they are the ones that are locally given by an almost finite amount of data:

\begin{definition}
Let $X \in \vStackspip$. A $\varphi$-module $\mathcal P \in \DqcohriX X^\varphi$ is called \emph{perfect} if the underlying object in $\DqcohriX X$ is weakly almost perfect (see \cref{def:weakly-almost-perfect-module-on-v-stack}). We denote
\begin{align*}
	\DqcohriX X^\varphi_\perf \subset \DqcohriX X^\varphi
\end{align*}
the full subcategory spanned by the perfect $\varphi$-modules.
\end{definition}

\begin{lemma}
\begin{lemenum}
	\item \label{rslt:perfect-phi-modules-satisfy-v-hyperdescent} The assignment $X \mapsto \DqcohriX X^\varphi_\perf$ defines a hypercomplete v-sheaf on $\vStackspip$.

	\item Let $X \in \vStackspip$. Then $\mathcal P \in \DqcohriX X^\varphi$ is perfect if and only if it is dualizable.
\end{lemenum}
\end{lemma}
\begin{proof}
Part (i) follows from \cref{rslt:weakly-alm-perf-sheaves-satisfy-v-descent}. Part (ii) follows from \cref{rslt:weakly-alm-perf-equiv-dualizable} by the same argument as in \cref{rslt:sigma-module-dualizable-iff-underlying-module-is-so}.
\end{proof}

With the above definitions and results we understand one side of the Riemann-Hilbert correspondence: the $\varphi$-modules. We now introduce the other side of the correspondence, i.e. the $\Fld_p$-sheaves. The right category here is the $\infty$-category of overconvergent étale $\Fld_p$-sheaves. This contains all $\Fld_p$-local systems and all Zariski constructible sheaves and thus most of the interesting examples of $\Fld_p$-sheaves. The definition is as follows.

\begin{definition} \label{def:overconvergent-Fp-sheaves}
Let $X$ be a small v-stack. A sheaf $\mathcal F \in \D_\et(X, \Fld_p)$ is called \emph{overconvergent} if for all maps $i\colon \Spa(C, C^+) \to X$ from a geometric point $\Spa(C, C^+)$ and induced maps $i'\colon \Spa(C, C^\circ) \to X$, $j\colon \Spa(C, C^\circ) \injto \Spa(C, C^+)$, the natural map $i^* \mathcal F \isoto j_* i'^* \mathcal F$ is an isomorphism. We denote
\begin{align*}
	\D_\et(X, \Fld_p)^\oc \subset \D_\et(X, \Fld_p)
\end{align*}
the full subcategory spanned by the overconvergent sheaves.
\end{definition}

We get the following basic results on overconvergent $\Fld_p$-sheaves. Notably they imply that $\D_\et(-, \Fld_p)^\oc$ is a hypercomplete v-sheaf which on strictly totally disconnected spaces is of the form $\D(\Fld_p(X))_\omega$, i.e. forms the left-hand side of the local Riemann-Hilbert correspondence from \cref{rslt:local-Riemann-Hilbert}.

\begin{lemma}
\begin{lemenum}
	\item \label{rslt:overconvergent-F-p-sheaves-form-v-hypersheaf} The assignment $X \mapsto \D_\et(X, \Fld_p)^\oc$ defines a hypercomplete v-sheaf of symmetric monoidal $\infty$-categories on the site of small v-stacks.

	\item A sheaf $\mathcal F \in \D_\et(X, \Fld_p)$ is overconvergent if and only if all $\pi_n(\mathcal F)$ are overconvergent.

	\item \label{rslt:overconvergent-F-p-sheaves-on-std-space-are-modules} Let $X$ be a strictly totally disconnected perfectoid space. Then there is a natural equivalence $\D_\et(X, \Fld_p)^\oc = \D(\Fld_p(X))_\omega$.

	\item Let $f\colon Y \to X$ be a qcqs map of locally spatial diamonds. Then the étale pushforward $f_*\colon \D^+_\et(Y, \Fld_p) \to \D^+_\et(X, \Fld_p)$ preserves overconvergent sheaves and hence restricts to a functor $f_*\colon \D^+_\et(Y, \Fld_p)^\oc \to \D^+_\et(X, \Fld_p)^\oc$.
\end{lemenum}
\end{lemma}
\begin{proof}
To prove (i) we first note that the assignment $X \mapsto \D_\et(X, \Fld_p)$ is a hypercomplete v-sheaf of $\infty$-categories: By \cite[Theorem 14.12]{etale-cohomology-of-diamonds} this reduces to the claim that the assignment $X \mapsto \D(X_\vsite, \Fld_p)$ is a hypercomplete v-sheaf, which is formal (the required functoriality of the assignment can be constructed as in \cite[\S2]{enhanced-six-operations}; then v-hyperdescent is formal, see \cite[Proposition 17.3]{etale-cohomology-of-diamonds}). Thus in order to prove (i) we only have to verify the following: Given a v-cover $f\colon Y \surjto X$ of small v-stacks and some $\mathcal F \in \D_\et(X, \Fld_p)$ such that $f^* \mathcal F$ is overconvergent, then $\mathcal F$ is overconvergent. We can immediately reduce this claim to the case that $X = \Spa(C, C^+)$ and $Y = \Spa(C', C'^+)$ are geometric points. In this case the étale sites of $X$ and $Y$ are equivalent to the analytic sites, and $\abs X$ and $\abs Y$ consist of a totally ordered chain of points under specialization. By assumption we have $(f^* \mathcal F)(V) = (f^* \mathcal F)(Y)$ for every open subset $V \subset Y$. One checks that $(f^* \mathcal F)(V) = \mathcal F(f(V))$, hence $\mathcal F(U) = \mathcal F(X)$ for all open subsets $U \subset X$, as desired. The functoriality of the symmetric monoidal structure on $\D_\et(X, \Fld_p)^\oc$ can either be constructed directly or deduced from (iii) and \cref{rslt:sheaves-on-basis-equiv-sheaves-on-whole-site}.

Part (ii) follows easily from the fact that all the functors involved in the definition of overconvergent sheaves are $t$-exact and that isomorphisms can be checked on $\pi_n$.

We now prove (iii), so assume that $X$ is a strictly totally disconnected perfectoid space. Then the étale site of $X$ is equivalent to the analytic site of $X$ and we have $\D_\et(X, \Fld_p) = \D(X_\et, \Fld_p) = \D(\abs X, \Fld_p)$. Let $\rho\colon \abs X \to \pi_0(X)$ be the projection. Then $\rho^*\colon \D(\pi_0(X), \Fld_p) \injto \D(\abs X, \Fld_p)$ is fully faithful: This amounts to the claim that $\id \isoto \rho_* \rho^*$ is an isomorphism, which can be checked on stalks and hence in the case that $X = \Spa(C, C^+)$ is connected; in this case it is obvious. A sheaf $\mathcal F \in \D(\abs X, \Fld_p)$ lies in the essential image of $\rho^*$ if and only if the map $\rho^* \rho_* \mathcal F \to \mathcal F$ is an isomorphism. By passing to connected components of $X$ one sees that this condition is satisfied precisely if $\mathcal F$ is overconvergent. We therefore obtain a natural equivalence $\rho^*\colon \D(\pi_0(X), \Fld_p) \isoto \D_\et(X, \Fld_p)^\oc$ of symmetric monoidal $\infty$-categories. Now the topological space $\pi_0(X)$ admits a basis by clopen subsets, each of which is a retract of $\pi_0(X)$. It follows easily that the natural functor $\D(\Fld_p(X))_\omega \to \D(\pi_0(X), \Fld_p)$ (given by pullback along $\pi_0(X) \to *$) is an equivalence of symmetric monoidal $\infty$-categories.

Finally, part (iv) can be shown in the same way as \cref{rslt:qcqs-pushforward-preserves-overconvergent-sheaves}, noting that the overconvergence condition can be checked on \emph{quasi-pro-étale} maps $\Spa(C, C^+) \to X$ by the arguments in (i).
\end{proof}

Similarly as for the $\varphi$-module side, we also introduce perfect objects in $\D_\et(X, \Fld_p)^\oc$ and prove some basic descent results.

\begin{definition}
Let $X$ be a small v-stack. An overconvergent sheaf $\mathcal F \in \D_\et(X, \Fld_p)^\oc$ is called \emph{perfect} if for every map $f\colon Y \to X$ from a strictly totally disconnected space $Y$ the pullback $f^* \mathcal F \in \D_\et(Y, \Fld_p)^\oc = \D(\Fld_p(Y))$ is perfect. We denote
\begin{align*}
	\D_\et(X, \Fld_p)^\oc_\perf \subset \D_\et(X, \Fld_p)^\oc
\end{align*}
the full subcategory of perfect sheaves.
\end{definition}

\begin{lemma}
\begin{lemenum}
	\item \label{rslt:perfect-overconvergent-F-p-sheaves-satisfy-v-hyperdescent} The assignment $X \mapsto \D_\et(X, \Fld_p)^\oc_\perf$ is a hypercomplete v-sheaf of symmetric monoidal $\infty$-categories.

	\item Let $X$ be a small v-stack. Then an object $\mathcal F \in \D_\et(X, \Fld_p)^\oc$ is perfect if and only if it is dualizable.
\end{lemenum}
\end{lemma}
\begin{proof}
As in the proof of \cref{rslt:geometric-properties-of-dualizable-and-perfect-sheaves} the full subcategory of dualizable objects in $\D_\et(X, \Fld_p)^\oc$ satisfies v-hyperdescent. Thus both (i) and (ii) reduce to the following claim: If $X$ is strictly totally disconnected then $\mathcal F \in \D_\et(X, \Fld_p)^\oc = \D(\Fld_p(X))_\omega$ is dualizable if and only if it is perfect. This follows from \cref{rslt:equiv-of-compact-perfect-dualizable-over-Lambda}, noting that a weakly almost perfect object in $\D(\Fld_p(X))_\omega$ is a retract of a perfect object and hence itself perfect.
\end{proof}

We have now introduced all the necessary $\infty$-categories of sheaves appearing in the Riemann-Hilbert correspondence. As a last step before stating the correspondence we need to introduce the Riemann-Hilbert functor:

\begin{definition}
Let $X \in \vStackspip$.
\begin{defenum}
	\item We denote
	\begin{align*}
		- \tensor \ri^{+a}_X/\pi\colon \D_\et(X, \Fld_p)^\oc \to \DqcohriX X^\varphi
	\end{align*}
	the symmetric monoidal functor which is obtained as follows: If $X = \Spa(A, A^+)$ is a strictly totally disconnected perfectoid space (w.l.o.g. of characteristic $p$) then it is the functor $- \tensor_R A^{+a}/\pi\colon \D(R)_\omega \to \Dqcohri(A^+/\pi)^\varphi$ using the identification $\D_\et(X, \Fld_p)^\oc = \D(R)$ for $R = \Fld_p(X)$ from \cref{rslt:overconvergent-F-p-sheaves-on-std-space-are-modules}. In general we define $- \tensor \ri^{+a}_X/\pi$ by gluing using \cref{rslt:overconvergent-F-p-sheaves-form-v-hypersheaf} (see \cref{rslt:sheaves-on-basis-equiv-sheaves-on-whole-site,rslt:trivial-sigma-structure-is-functorial}).

	\item By the adjoint functor theorem (see \cite[Remark 5.5.2.10]{lurie-higher-topos-theory}) the functor in (i) admits a right adjoint which we denote by
	\begin{align*}
		(-)^\varphi\colon \DqcohriX X^\varphi \to \D_\et(X, \Fld_p)^\oc.
	\end{align*}
\end{defenum}
\end{definition}

\begin{remark}
Note that the functor $- \tensor \ri^{+a}_X/\pi\colon \D_\et(X, \Fld_p)^\oc \to \DqcohriX X^\varphi$ factors over the full subcategory $\DqcohriX X^\varphi_\omega \subset \DqcohriX X^\varphi$ of \emph{discrete} $\varphi$-modules.
\end{remark}

We are finally in the position to state and prove the $p$-torsion Riemann-Hilbert correspondence on small v-stacks:

\begin{theorem} \label{rslt:global-Riemann-Hilbert}
Let $X \in \vStackspip$. Then the functor
\begin{align*}
	- \tensor \ri^{+a}_X/\pi\colon \D_\et(X, \Fld_p)^\oc \injto \DqcohriX X^\varphi
\end{align*}
is fully faithful and induces an equivalence of symmetric monoidal $\infty$-categories
\begin{align*}
	\D_\et(X, \Fld_p)^\oc_\perf \isoto \DqcohriX X^\varphi_\perf.
\end{align*}
\end{theorem}
\begin{proof}
By v-hyperdescent (see \cref{rslt:phi-modules-form-hypercomplete-v-sheaf}, \cref{rslt:overconvergent-F-p-sheaves-form-v-hypersheaf}, \cref{rslt:perfect-phi-modules-satisfy-v-hyperdescent} and \cref{rslt:perfect-overconvergent-F-p-sheaves-satisfy-v-hyperdescent}) the claim can be checked v-locally on $X$, so we can assume that $X$ is a strictly totally disconnected perfectoid space. Then the claim follows from \cref{rslt:local-Riemann-Hilbert}.
\end{proof}

The Riemann-Hilbert correspondence allows us to relate $\Fld_p$-sheaves and $\Fld_p$-cohomology to our theory of quasicoherent $\ri^{+a}_X/\pi$-modules. More concretely we deduce the following general version of the Primitive Comparison Theorem. Note that the relevant condition on $f\colon Y \to X$ (i.e. that $f^!$ preserves small colimits) can be checked \emph{locally on $Y$}. It is in particular satisfied by smooth maps of analytic adic spaces over $\Q_p$ by \cref{rslt:adic-smooth-maps-over-Qp-are-p-cohom-smooth} below, so that we recover \cite[Theorem 5.1]{rigid-p-adic-hodge}.

\begin{corollary} \label{rslt:primitive-comparison}
Let $f\colon Y \to X$ be a proper bdcs map of $p$-bounded locally spatial diamonds in $\vStackspip$. Assume that $f^!\colon \DqcohriX X \to \DqcohriX Y$ preserves all small colimits. Then:
\begin{corenum}
	\item The pushforward $f_*\colon \D^+_\et(Y, \Fld_p)^\oc \to \D^+_\et(X, \Fld_p)^\oc$ preserves perfect objects.

	\item For every $\mathcal F \in \D_\et(Y, \Fld_p)^\oc_\perf$ the natural morphism
	\begin{align*}
		f_* \mathcal F \tensor \ri^{+a}_X/\pi \isoto f_*(\mathcal F \tensor \ri^{+a}_Y/\pi)
	\end{align*}
	is an isomorphism in $\DqcohriX X$.
\end{corenum}
\end{corollary}
\begin{proof}
Both claims can be checked locally on $X$, so we can assume that $X$ and $Y$ are qcqs. Since $f^!$ preserves all small colimits, it follows from \cref{rslt:almost-hom-adjunction-for-etale-!} that $f_* = f_!\colon \DqcohriX Y \to \DqcohriX X$ preserves almost compact objects and it follows from \cref{rslt:qcqs-pushforward-preserves-left-bounded-discrete-sheaves} that this functor also preserves left-bounded discrete objects. It thus follows from \cref{rslt:equiv-of-compact-perfect-dualizable-on-p-bd-spat-diam} that $f_*$ preserves left-bounded weakly almost perfect objects. Combining this with \cref{rslt:pushforward-of-phi-modules} we deduce that $f_*\colon \DqcohriX Y^\varphi \to \DqcohriX X^\varphi$ preserves left-bounded perfect $\varphi$-modules. By \cref{rslt:global-Riemann-Hilbert} every perfect $\varphi$-module is bounded, hence $f_*$ preserves perfect $\varphi$-modules.

Now let $\mathcal F \in \D_\et(Y, \Fld_p)^\oc_\perf$ be given. Then by \cref{rslt:global-Riemann-Hilbert} we have $\mathcal F = (\mathcal F \tensor \ri^{+a}_Y/\pi)^\varphi$. It follows formally from adjunctions that $(-)^\varphi$ commutes with pushforwards, so we deduce
\begin{align*}
	(f_* (\mathcal F \tensor \ri^{+a}_Y/\pi))^\varphi = f_* (\mathcal F \tensor \ri^{+a}_Y/\pi)^\varphi = f_* \mathcal F
\end{align*}
On the other hand, by what we have shown above, $\mathcal M' := f_* (\mathcal F \tensor \ri^{+a}_Y/\pi)$ is perfect and thus $\mathcal M'^\varphi \tensor \ri^{+a}_X/\pi = \mathcal M'$ by \cref{rslt:global-Riemann-Hilbert}. Therefore (ii) follows from the above computation by applying $- \tensor \ri^{+a}_X/\pi$ and (i) follows from \cref{rslt:global-Riemann-Hilbert} because $f_* \mathcal F = \mathcal M'^\varphi$.
\end{proof}

Another corollary of the Riemann-Hilbert correspondence is that the invertible sheaves appearing in the definition of $p$-cohomologically smooth morphisms are rather simple (which is a priori not clear because invertible objects in the almost setting are quite subtle). In fact, they are given by invertible $\Fld_p$-sheaves and are in particular étale locally free:

\begin{corollary} \label{rslt:dualizing-complex-automatically-given-by-F-p-locsys}
Let $f\colon Y \to X$ be a bdcs and $p$-cohomologically smooth map of locally spatial diamonds in $\vStackspip$ and let $\mathcal L := f^! \ri^{+a}_X/\pi$. Then there is a unique invertible object $\mathcal L_0 \in \D_\et(Y, \Fld_p)$ such that
\begin{align*}
	\mathcal L = \mathcal L_0 \tensor_{\Fld_p} \ri^{+a}_Y/\pi.
\end{align*}
In particular, $\mathcal L$ is étale locally isomorphic to $\ri^{+a}_Y/\pi[n]$ for some integer $n$.
\end{corollary}
\begin{proof}
By \cref{rslt:p-cohom-smoothness-description-of-upper-shriek} $\mathcal L$ is an invertible object of $\DqcohriX Y$. By \cref{rslt:upper-shriek-of-phi-modules} it comes naturally equipped with a $\varphi$-action and can thus be seen as an invertible object of $\DqcohriX Y^\varphi$. In particular it is a perfect $\varphi$-module and hence by \cref{rslt:global-Riemann-Hilbert} corresponds to a perfect object $\mathcal L_0 \in \D_\et(Y, \Fld_p)$. Since \cref{rslt:global-Riemann-Hilbert} is an equivalence of symmetric monoidal $\infty$-categories, $\mathcal L_0$ is necessarily invertible.
\end{proof}

We expect many more applications of the Riemann-Hilbert correspondence to the rigid-analytic world, which we aim to investigate in subsequent papers. We end this subsection with the following neat little observation:

\begin{example}
Let $X = \Spa \Q_p$, let $\tilde K := \Q_p(\mu_{p^\infty})$ and denote $\tilde X := \Spa \tilde K$. Then by \cref{rslt:global-Riemann-Hilbert,rslt:compute-Dqcohri-on-stack-quotient-of-p-bounded-space} the category of $\Gal(\overline{\Q_p}/\Q_p)$-representations on finite discrete $\Fld_p$-vector spaces is equivalent to the category of finite free $\ri^a_{\tilde K}/p$-modules equipped with a compatible pair of a $\varphi$-action and a $\Gamma := \Z_p^\cprod$-action. This is a version of the classical $(\varphi, \Gamma)$-correspondence.
\end{example}

\subsection{Poincaré Duality for Rigid-Analytic Varieties} \label{sec:ri-pi.poincare}

We now apply the formalism developed in the previous subsections to the rigid-analytic setting in order to prove a version of Poincaré duality for $\Fld_p$-cohomology. By definition of $p$-cohomologically smooth morphisms this essentially amounts to saying that if $f\colon Y \to X$ is a smooth map of analytic adic spaces over $\Q_p$ then the associated map of diamonds $Y^\diamond \to X^\diamond$ is $p$-cohomologically smooth. By the stability properties of $p$-cohomological smoothness (see \cref{rslt:stability-of-p-cohom-smooth}) we can reduce this question to the case that $X$ is totally disconnected and $Y$ is a torus over $X$ -- this is where the 6-functor formalism really simplifies the computation, as it allows us to localize on $Y$. In order to handle the local case we perform a computation involving formal models similar to Faltings' original work \cite{faltings-p-adic-hodge}.

The local computation relies on a good understanding of the interplay of the theory of quasicoherent $\ri^{+a}_X/\pi$-modules with quasicoherent sheaves on a formal model. To motivate the following definitions, recall that for any (qcqs) rigid-analytic variety $X$ over a non-archemedian field $K$, there is a formal model $\mathfrak X$ of $X$ and a morphism
\begin{align*}
	(X, \ri^+_X) \to (\mathfrak X, \ri_{\mathfrak X})
\end{align*}
of locally ringed topological spaces. Given a pseudouniformizer $\pi \in K$ we denote $\mathfrak X_\pi := V(\pi) \subset \mathfrak X$. Then the above morphism also induces a morphism
\begin{align*}
	t\colon (X_\proet, \ri^+_X/\pi) \to (X, \ri^+_X/\pi) \to (\mathfrak X_\pi, \ri_{\mathfrak X_\pi})
\end{align*}
of locally ringed spaces (cf. \cref{def:morphism-t-to-formal-model} below). Now $\mathfrak X_\pi$ is a scheme and hence admits a good theory of solid quasicoherent sheaves $\D_\solid(\ri_{\mathfrak X_\pi})$. Similarly, on $X^\diamond_\proet$ we defined the category $\DqcohriX{X^\diamond}$ of quasicoherent $\ri^{+a}/\pi$-modules. One would thus expect there to be a pair of adjoint functors
\begin{align*}
	t^*\colon \D_\solid(\ri_{\mathfrak X_\pi}) \rightleftarrows \DqcohriX{X^\diamond}\noloc t_*.
\end{align*}
We will now construct this adjunction and study its properties. Of course we want the formation of $t^*$ and $t_*$ to be functorial in $t$ (in a suitable sense) and there is no reason for us to restrict to formal models -- any map of ringed spaces as above works. This leads to the following definition:

\begin{definition}
Let $S = \Spa(R, R^+) \in \AffPerfd_\pi$ be a totally disconnected perfectoid space with pseudouniformizer $\pi$. Let $\StkSch_S$ be the following $2$-category:
\begin{enumerate}[(i)]
	\item The objects of $\StkSch_S$ are given by the disjoint union of small v-stacks over $S$ and classical schemes over $R^+/\pi$.

	\item $1$-morphisms $f\colon Y \to X$ of two objects $X, Y \in \StkSch_S$ are given as follows: If both $X$ and $Y$ are small v-stacks or both are schemes then the morphisms $Y \to X$ are the ones in the respective category (over $S$ or over $\Spec R^+/\pi$, respectively). If $Y$ is a scheme and $X$ is not, then there is no morphism $Y \to X$. Finally, if $X$ is a scheme and $Y$ is a small v-stack, then morphisms $Y \to X$ are the morphisms of locally ringed spaces
	\begin{align*}
		(\abs Y, \ri^+_Y/\pi) \to (X, \ri_X).
	\end{align*}
	over $(\abs S, \ri^+_S/\pi) = \Spec R^+/\pi$ (cf. \cref{rslt:tot-disc-space-mod-pi-isom-to-scheme}).

	\item $2$-morphisms in $\StkSch_S$ are given as follows: If all associated objects are small v-stacks then we use the definition of $2$-morphisms of small v-stacks. If one of the associated objects is a scheme then there are no non-trivial $2$-morphisms.
\end{enumerate}
\end{definition}

Our aim is to define natural functors $t_*$ and $t^*$ associated to every morphism $t\colon Y \to X$ in $\StkSch_S$. In order to do so, it is useful to first do this on the full subcategory where only those small v-stacks are allowed which are totally disconnected perfectoid spaces and afterwards extend to all small v-stacks. For the second step it is handy to have a site on $\StkSch_S$ which ``restricts'' to the sites of small v-stacks and discrete adic spaces respectively:

\begin{definition} \label{def:site-on-StkSch}
Let $S = \Spa(R, R^+) \in \AffPerfd_\pi$ be totally disconnected. We define a site on $\StkSch_S$ by declaring a sieve $\mathcal C_X$ over an object $X \in \StkSch_S$ to be a covering sieve in the following cases:
\begin{enumerate}[(i)]
	\item $X$ is a small v-stack and there are small v-stacks $U_i \in \mathcal C_X$ such that the family $(U_i \to X)_i$ defines a v-cover.
	\item $X$ is a scheme and there are schemes $U_i \in \mathcal C_X$ such that the family $(U_i \to X)_i$ is an open cover.
\end{enumerate}
Note that for every map $Y \to X$ in $\StkSch_S$ and every open immersion $U \injto X$ the fiber product $Y \cprod_X U$ exists and is an open subset of $Y$, namely the preimage of $U$ under $\abs Y \to \abs X$. This implies that the above definition of the site on $\StkSch_S$ is valid.
\end{definition}

We are finally in the position to define the functors $t_*$ and $t^*$ promised above. With the category $\StkSch_S$ at hand, this is essentially formal:

\begin{proposition} \label{rslt:D-solid-functor-on-StkSch}
Let $S = \Spa(R, R^+) \in \AffPerfd_\pi$ be totally disconnected. There is a unique (up to isomorphism) sheaf of $\infty$-categories
\begin{align*}
	\mathcal F\colon \StkSch_S^\opp \to \infcatinf
\end{align*}
with the following properties:
\begin{propenum}
	\item Restricted to the category of small v-stacks, $\mathcal F$ is the sheaf $(X, \pi) \mapsto \DqcohriX X$.
	\item Restricted to the category of schemes, $\mathcal F$ is the sheaf $X \mapsto \D^a_\solid(\ri_X)$.
	\item Let $X = \Spa(A, A^+)$ be a totally disconnected space over $S$, let $\tilde X = \Spec A^+/\pi$ and let $\psi\colon (\abs{X}, \ri^+_X/\pi) \isoto (\abs{\tilde X}, \ri_{\tilde X})$ be the isomorphism of ringed spaces from \cref{rslt:tot-disc-space-mod-pi-isom-to-scheme}. Then
	\begin{align*}
		\mathcal F(\psi)\colon \D^a_\solid(A^+/\pi) \isoto \Dqcohri(A^+/\pi)
	\end{align*}
	is the identity.
\end{propenum}
\end{proposition}
\begin{proof}
Let $\StkSch'_S \subset \StkSch_S$ be the full subcategory which contains all (classical) schemes but only those small v-stacks which are totally disconnected perfectoid spaces. Then the functor $\mathcal F\colon (\StkSch'_S)^\opp \to \infcatinf$ is easily constructed using the functoriality of the construction $X \mapsto \D^a_\solid(X)$ for schemes (cf. \cref{rslt:functoriality-of-sheaves-on-analytic-spaces-over-AlmSetup}). It is clear that $\mathcal F$ is uniquely determined (on $\StkSch'_S$) by the conditions (i), (ii) and (iii). But $\StkSch'_S$ is a basis of $\StkSch_S$ with respect to the site \cref{def:site-on-StkSch}, so by \cref{rslt:sheaves-on-basis-equiv-sheaves-on-whole-site}, $\mathcal F$ can uniquely be extended to $\StkSch_S$.
\end{proof}

\begin{definition}
Let $S = \Spa(R, R^+) \in \AffPerfd_\pi$ be totally disconnected and let $t\colon Y \to X$ be a morphism in $\StkSch_S$.
\begin{defenum}
	\item We denote $t^* := \mathcal F(t)$, where $\mathcal F$ is defined as in \cref{rslt:D-solid-functor-on-StkSch}. In particular, if $Y$ is a small v-stack and $X$ is a scheme then we obtain the functor
	\begin{align*}
		t^*\colon \D^a_\solid(\ri_X) \to \DqcohriX Y.
	\end{align*}
	The existence of $\mathcal F$ makes the construction $t \mapsto t^*$ functorial, in particular we have $t_1^* \comp t_2^* = (t_2 \comp t_1)^*$ for all composable morphisms $t_1$ and $t_2$. It follows easily from the definitions that $t^*$ preserves all small colimits.

	\item We denote $t_*$ a right adjoint of $t^*$. This is automatically functorial, so in particular we have $t_{1*} \comp t_{2*} = (t_1 \comp t_2)_*$ for all composable morphisms $t_1$ and $t_2$. If $Y$ is a small v-stack and $X$ is a scheme then we obtain the functor
	\begin{align*}
		t_*\colon \DqcohriX Y \to \D^a_\solid(\ri_X).
	\end{align*}
\end{defenum}
\end{definition}

\begin{remark}
The functors $t_*$ and $t^*$ can be described more explicitly as follows. Let $t\colon Y \to X$ be a morphism in $\StkSch_S$ with $X$ a scheme and $Y$ a small v-stack.
\begin{itemize}
 	\item Given any $\mathcal M \in \D^a_\solid(\ri_X)$, $t^* \mathcal M \in \DqcohriX Y$ is the object whose pullback to every totally disconnected $Z = \Spa(C, C^+) \to Y$ is $t_Z^* \mathcal M$ where $t_Z$ is the induced map of schemes $\Spec C^+/\pi \to X$.

 	\item To describe $t_*$, pick any hypercover $Z_\bullet \to Y$ such that each $Z_n = \bigdunion_{i\in I_n} Z_{n,i}$ is a disjoint union of (untilted) totally disocnnected spaces $Z_{n,i} = \Spa(C^{n,i}, C^{n,i+})$. Let $t_{n,i}\colon \Spec C^{n,i+}/\pi \to X$ denote the induced map of schemes for all $n$ and $i$. Then, given any
 	\begin{align*}
 		\mathcal N = (N^{n,i})_{n,i} \in \DqcohriX Y = \varprojlim_{n\in\Delta} \prod_{i\in I_n} \Dqcohri(C^{n,i+}/\pi)
 	\end{align*}
 	(cf. \cref{rslt:explicit-description-of-DqcohriX-on-vstack}) we have
 	\begin{align*}
 		t_* \mathcal N = \varprojlim_{n\in \Delta} \prod_{i\in I_n} t_{n,i*} N^{n,i} \in \D_\solid(\ri_X).
 	\end{align*}
\end{itemize}
\end{remark}

Both for schemes (see \cref{rslt:scheme-6-functor,rslt:lower-shriek-commutes-with-almost-localization}) and for small v-stacks (see \cref{sec:ri-pi.6-functor}) we have a six-functor formalism and in particular a good notion of the functors $f_!$ and $f^!$ for suitable morphisms $f$. In order to reduce the computation of these functors from rigid spaces to their formal models we need to show some compatibilities between the scheme 6-functor formalism and the small v-stack 6-functor formalism. We start with some general results:

\begin{lemma} \label{rslt:open-immersion-lower-shriek-base-change-in-StkSch}
Let $S = \Spa(R, R^+) \in \AffPerfd_\pi$ be totally disconnected and let
\begin{center}\begin{tikzcd}
	U' \arrow[d,"j'",hook] \arrow[r,"g'"] & U \arrow[d,"j",hook]\\
	X' \arrow[r,"g"] & X
\end{tikzcd}\end{center}
be a cartesian diagram in $\StkSch_S$ with $j$ being an open immersion. Then the natural transformation
\begin{align*}
	j'_! g'^* \isoto g^* j_!
\end{align*}
of functors $\mathcal F(U) \to \mathcal F(X')$ (with $\mathcal F$ as in \cref{rslt:D-solid-functor-on-StkSch}) is an isomorphism.
\end{lemma}
\begin{proof}
If all involved spaces are schemes or all involved spaces are small v-stacks then the claim was already proven (for the latter case see \cref{rslt:existence-of-lower-shriek-for-etale-maps}), thus we can assume that $X$ and $U$ are schemes and $X'$ and $U'$ are untilted small v-stacks. Since, as mentioned, $j_!$ satisfies arbitrary base-change in schemes and small v-stacks, we can easily reduce to the case that $X = \Spec A$ is affine and $X' = \Spa(A', A'^+)$ is totally disconnected. In fact, we can even reduce to the case $A = A'^+/\pi$. As in the proof of \cref{rslt:etale-lower-shriek-compatible-with-pushforward}, $j_!$ is the colimit of the $j_{i!}$'s for any open cover $U = \bigunion_i U_i$ with maps $j_i\colon U_i \to X$, so since $g^*$ and $g'^*$ preserve colimits we can reduce to the case that $U = U_i$ is a distinguished open subset. Now $g^*$ and $g'^*$ are the identity functors, so the claim is clear.
\end{proof}

\begin{corollary} \label{rslt:base-change-in-StkSch-along-open-immersions}
Let $S = \Spa(R, R^+) \in \AffPerfd_\pi$ be totally disconnected and let
\begin{center}\begin{tikzcd}
	U' \arrow[d,"f'"] \arrow[r,"j'",hook] & X' \arrow[d,"f"]\\
	U \arrow[r,"j",hook] & X
\end{tikzcd}\end{center}
be a cartesian diagram in $\StkSch_S$ with $j$ being an open immersion. Then the natural transformation
\begin{align*}
	j^* f_* \isoto f'_* j'^*
\end{align*}
of functors $\mathcal F(X') \to \mathcal F(U)$ (with $\mathcal F$ as in \cref{rslt:D-solid-functor-on-StkSch}) is an isomorphism.
\end{corollary}
\begin{proof}
Pass to right adjoints in \cref{rslt:open-immersion-lower-shriek-base-change-in-StkSch}.
\end{proof}

Having studied a rather general version of the functors $t_*$ and $t^*$, we will now turn our focus on our main application for this theory: formal models. Let us first make precise what we mean by the morphism $t$ in this scenario. We do not want to deal with subtleties arising from the non-noetherian setting over the totally disconnected space $S$ so we contend ourselves with a somewhat ad-hoc approach by only using \emph{affine} formal models. The following definitions and results should extend to all suitable formal models (or even ``almost adic spaces'') over $S$ in a straightforward way.

\begin{definition} \label{def:generic-and-special-fiber-of-formal-scheme}
Let $S = \Spa(R, R^+)$ be a totally disconnected perfectoid space and let $\mathfrak X = \Spf A^+$ be an affine formal scheme with an adic map to $\mathfrak S := \Spf R^+$ (i.e. $A^+$ is a classically $\pi$-adically complete $R^+$-algebra).
\begin{defenum}
	\item Let $\mathfrak X^\diamond$ be the untilted small v-sheaf associated to $\mathfrak X$ (see \cite[\S18.1]{scholze-berkeley-lectures}). The \emph{generic fiber} $X$ of $\mathfrak X$ is the small v-sheaf
	\begin{align*}
		X = \mathfrak X^\diamond \cprod_{\mathfrak S^\diamond} S^\diamond.
	\end{align*}

	\item For a pseudouniformizer $\pi \in R$ we denote
	\begin{align*}
		\mathfrak X_\pi := V(\pi) \subset \mathfrak X.
	\end{align*}
	As $\mathfrak X \to \mathfrak S$ is adic by assumption, $\mathfrak X_\pi$ is a (classical) scheme over $\Spec R^+/\pi$.
\end{defenum}
\end{definition}

\begin{lemma} \label{rslt:map-of-ringed-spaces-from-formal-scheme-diamond}
Let $\mathfrak X$ be an affine formal scheme over $\Spf \Z_p$ and let $\mathfrak X^\diamond$ be the associated untilted small v-sheaf (see \cite[\S18.1]{scholze-berkeley-lectures}). Then there is a natural surjective morphism of locally ringed spaces
\begin{align*}
	(\abs{\mathfrak X^\diamond}, \ri^+_{X^\diamond}) \to (\abs{\mathfrak X}, \ri_{\mathfrak X}).
\end{align*}
\end{lemma}
\begin{proof}
By \cite[\S18.2]{scholze-berkeley-lectures} there is a natural continuous surjective map $f\colon \abs{\mathfrak X^\diamond} \to \abs{\mathfrak X}$, so it only remains to find a suitable map $\ri_{\mathfrak X} \to f_* \ri^+_{\mathfrak X^\diamond}$ of sheaves on $\mathfrak X$. Choose any v-hypercover $Y_\bullet \to \mathfrak X^\diamond$ by perfectoid spaces $Y_n$ and denote $g_n\colon Y_n \to \mathfrak X^\diamond$ the projection. Then
\begin{align*}
	f_* \ri^+_{\mathfrak X^\diamond} = \eq(f_* g_{0*} \ri^+_{Y_0} \rightrightarrows f_* g_{1*} \ri^+_{Y_1}).
\end{align*}
Hence it is enough to construct compatible maps $\ri_{\mathfrak X} \to f_* g_{n*}\ri^+_{Y_n}$. But by definition of $\mathfrak X^\diamond$ the map $g_n$ is given by a map $g'_n\colon Y_n \to \mathfrak X$ of pre-adic spaces and clearly $g'_n = f \comp g_n$ on underlying topological spaces. This produces the desired morphism.
\end{proof}

\begin{definition} \label{def:morphism-t-to-formal-model}
Let $S = \Spa(R, R^+) \in \AffPerfd_\pi$ be totally disconnected and let $\mathfrak X$ be an adic affine formal scheme over $\mathfrak S = \Spf R^+$ with generic fiber $X$. We define the morphism
\begin{align*}
	t_{\mathfrak X}\colon (X, \pi) \to \mathfrak X_\pi
\end{align*}
in $\StkSch_S$ as follows: Using \cref{rslt:map-of-ringed-spaces-from-formal-scheme-diamond} we get a map of locally ringed spaces $(\abs X, \ri^+_X) \to (\abs{\mathfrak X^\diamond}, \ri^+_{\mathfrak X^\diamond}) \to (\abs{\mathfrak X}, \ri_{\mathfrak X})$. Modding out by $\pi$ and noting that $\abs{\mathfrak X} = \abs{\mathfrak X_\pi}$ we arrive at the desired map
\begin{align*}
	(\abs X, \ri^+_X/\pi) \to (\abs{\mathfrak X_\pi}, \ri_{\mathfrak X_\pi})
\end{align*}
of locally ringed spaces. In the following we will often write $t$ for $t_{\mathfrak X}$.
\end{definition}

If the generic fiber of the formal scheme $\mathfrak X$ (as defined in \cref{def:generic-and-special-fiber-of-formal-scheme}) is $p$-bounded over the base then the morphism $t_{\mathfrak X}$ associated to $\mathfrak X$ as in \cref{def:morphism-t-to-formal-model} behaves a lot like a $p$-bounded proper morphism. This will be made more precise in the following results.

\begin{setup} \label{stp:p-bounded-formal-scheme-over-tot-disc}
We will frequently make use of the following setup: $S = \Spa(R, R^+) \in \AffPerfd_\pi$ is a totally disconnected perfectoid space with pseudouniformizer $\pi$, $\mathfrak X = \Spf A^+$ is an adic affine formal scheme over $\Spf R^+$ and $X = \Spa(A, A^+)^\diamond$ is the generic fiber of $\mathfrak X$. We furthermore assume that $X \to S$ is $p$-bounded.
\end{setup}

\begin{lemma} \label{rslt:compute-Dqcohri-for-rel-compactification-of-tot-disc-over-rigid-space}
Let $S$, $\mathfrak X$ and $X$ be as in \cref{stp:p-bounded-formal-scheme-over-tot-disc}. Then for every totally disconnected space $Y = \Spa(B, B^+)$ and map $Y \to X$ with $\dimtrg(Y/X) < \infty$ we have
\begin{align*}
	\DqcohriX{\overline Y^{/X}} = \Dqcohri(B^{+a}/\pi, A^+/\pi),
\end{align*}
compatibly with $\Dqcohrip$, $\Dqcohrim$ and $\Dqcohrib$ on both sides.
\end{lemma}
\begin{proof}
We have $\overline Y^{/X} = \Spa(B, (A^+ + \mm_B)')$, where $(-)'$ denotes the integral closure in $B$. In particular $\overline Y^{/X}$ is an affinoid perfectoid space and
\begin{align*}
	((A^+ + \mm_B)'/\pi)^a_\solid = (B^{+a}/\pi, A^+/\pi)_\solid.
\end{align*}
Note that $\overline Y^{/X} \to \overline Y^{/S}$ is pro-étale, so the claim follows from \cref{rslt:compute-bounded-Dqcohri-for-qproet-over-Z'-Z}.
\end{proof}

\begin{proposition}
Let $S$, $\mathfrak X$ and $X$ be as in \cref{stp:p-bounded-formal-scheme-over-tot-disc}.
\begin{propenum}
	\item $t_*\colon \DqcohriX X \to \D^a_\solid(\ri_{\mathfrak X_\pi})$ preserves all small colimits.

	\item \label{rslt:projection-formula-for-t} $t_*$ satisfies the projection formula, i.e. for all $\mathcal M \in \D^a_\solid(\ri_{\mathfrak X_\pi})$ and $\mathcal N \in \DqcohriX X$ the natural morphism
	\begin{align*}
		(t_* \mathcal N) \tensor \mathcal M \isoto t_* (\mathcal N \tensor t^* \mathcal M)
	\end{align*}
	is an isomorphism.
\end{propenum}
\end{proposition}
\begin{proof}
To prove (i), note that the pushforward $f_*\colon \D_\solid(\ri_{\mathfrak X_\pi}) \to \D_\solid(\ri_K/\pi)$ is just a forgetful functor (since $\mathfrak X$ is affine) and in particular it is conservative and preserves small colimits. It is thus enough to show that the composed functor $f_* t_*\colon \DqcohriX X \to \D^a_\solid(\ri_K/\pi)$ preserves small colimits. Looking at the commutative diagram
\begin{center}\begin{tikzcd}
	X \arrow[r,"t"] \arrow[d,"f'"] & \mathfrak X_\pi \arrow[d,"f"]\\
	\Spa(K, \ri_K) \arrow[r,"t'"] & \Spec(\ri_K/\pi)
\end{tikzcd}\end{center}
in $\StkSch_S$ we note that it is enough to show that $t'_* f'_*$ preserves small colimits. But $t'_*$ is the identity functor and $f'$ is $p$-bounded by assumption so that $f'_*$ preserves small colimits by \cref{rslt:qcqs-p-bounded-pushforward-preserves-colimits-and-Dqcohrim}.

To prove (ii) we argue as in \cref{rslt:projection-formula-for-proper-p-bounded}: Since $t_*$ commutes with colimits by (i), both sides of the claimed isomorphism commute with colimits in $\mathcal M$, so we can assume that $\mathcal M = (A^+/\pi)^a_\solid[S]$ for some profinite set $S$. By a similar argument as in (i) we see that the functor $t_*\colon \DqcohriXZ X \to \D^a_\solid(A^+/\pi, \Z)$ (which can be constructed in the same way as $t_*$ by adapting \cref{rslt:D-solid-functor-on-StkSch}) also preserves small colimits and has finite cohomological dimension with respect to the canonical $t$-structures on both sides. It follows that both sides of the claimed isomorphism commute with Postnikov limits in $\mathcal N$, so that we can assume that $\mathcal N$ is left-bounded. Choose any hypercover $Y_\bullet \to X$ by totally disconnected spaces $Y_n = \Spa(B^n, B^{n+})$. Then the relative compactifications of $Y_\bullet$ over $X$ still form a hypercover, so by \cref{rslt:compute-Dqcohri-for-rel-compactification-of-tot-disc-over-rigid-space} we have
\begin{align*}
	\DqcohriX X &= \varprojlim_{n\in\Delta} \Dqcohri(B^{n+,a}/\pi, A^+/\pi),\\
	\DqcohriXZ X &= \varprojlim_{n\in\Delta} \Dqcohri(B^{n+,a}/\pi, \Z).
\end{align*}
Thus $\mathcal N$ is represented by a cosimplicial object $N^\bullet$ in $\Dqcohri(A^+/\pi)$ with $N^n \in \Dqcohri(B^{n+,a}/\pi, A^+/\pi)$ for all $n$, and $t_*$ computes the totalization. Now the same computation as in \cref{rslt:projection-formula-for-proper-p-bounded} implies the claim.
\end{proof}

The main purpose of introducing the morphism $t_{\mathfrak X}$ for formal schemes is that it often allows us to reduce the computation of the functor $f_!$ on quasicoherent $\ri^{+a}_X/\pi$-modules on small v-stacks to the computation of the functor $f_!$ in the setting of schemes. This should be true quite generally, but for simplicity we only show it in the case we care about:

\begin{proposition} \label{rslt:compute-lower-shriek-via-pushforward-to-formal-model}
Let $S$, $\mathfrak X$ and $X$ be as in \cref{stp:p-bounded-formal-scheme-over-tot-disc} and assume that $A^+ = R^+\langle T^\pm \rangle$, i.e. $X$ is the relative torus over $S$. Then there is a natural isomorphism
\begin{align*}
	f_! \isom \mathfrak f_{\pi !} \comp t_{\mathfrak X *}
\end{align*}
of functors $\DqcohriX X \to \Dqcohri(R^+/\pi)$.
\end{proposition}
\begin{proof}
We let $\mathfrak Y$ be the relative projective line over $\Spf R^+$, glued out of two copies of $\Spf R^+\langle T \rangle$ along their intersection. Then $\mathfrak Y_\pi := V(\pi) \subset \mathfrak Y$ is the relative projective line $\mathfrak g_\pi\colon \mathfrak Y_\pi \to \Spec R^+/\pi$ and the generic fiber of $\mathfrak Y$ is the relative projective line $g\colon Y \to S$. We get a canonical map $t = t_{\mathfrak Y}\colon Y \to \mathfrak Y_\pi$ in $\StkSch_S$ (even though we have defined $t$ only for \emph{affine} formal schemes, it is easy to define it for $\mathfrak Y$ as well). We have a natural open immersion $\mathfrak j_\pi\colon \mathfrak X_\pi \injto \mathfrak Y_\pi$, whose pullback along $Y \to \mathfrak Y_\pi$ is the natural open immersion $j\colon X \injto Y$. We get the following diagram of $\infty$-categories:
\begin{center}\begin{tikzcd}
	\DqcohriX X \arrow[r,"j_!"] \arrow[d,"t_*"] & \DqcohriX Y \arrow[r,"g_*"] \arrow[d,"t_*"] & \Dqcohri(R^+/\pi)\\
	\Dqcohri(\ri_{\mathfrak X_\pi}) \arrow[r,swap,"\mathfrak j_{\pi!}"] & \Dqcohri(\ri_{\mathfrak Y_\pi}) \arrow[ur,swap,"\mathfrak g_{\pi*}"]
\end{tikzcd}\end{center}
The composition of the upper horizontal maps computes $f_!$, while the composition of the bottom horizontal map and the diagonal map computes $\mathfrak f_{\pi!}$. It is clear that the triangle on the right commutes, so in order to finish the claim it is enough to show that the left square also commutes. In other words we need to see that the natural morphism $\mathfrak j_{\pi!} \comp t \isoto t \comp j_!$ is an isomorphism. But this is a formal consequence of the projection formulas for $j_!$, $\mathfrak j_{\pi!}$ and $t_*$ (see \cref{rslt:projection-formula-for-open-immersion,rslt:projection-formula-for-t}) and base-change for all these maps along open immersions (see \cref{rslt:open-immersion-lower-shriek-base-change-in-StkSch,rslt:base-change-in-StkSch-along-open-immersions}), cf. the proof of \cref{rslt:etale-lower-shriek-compatible-with-pushforward}.
\end{proof}

\begin{remark}
As mentioned above, \cref{rslt:compute-lower-shriek-via-pushforward-to-formal-model} holds in much greater generality, i.e. for arbitrary formal models over $S$ or more generally for certain (almost) adic spaces over $S$. However, stating such a result requires a good understanding of (almost) adic spaces over $\Spa(R^+/\pi)$, which we avoided here for simplicity.
\end{remark}

With \cref{rslt:compute-lower-shriek-via-pushforward-to-formal-model} at hand we finally have a good way of computing the lower shriek functor on the torus by passage to its canonical formal model. Thus there is nothing in the way of performing Faltings' computation in our setting now:

\begin{lemma} \label{rslt:Koszul-duality-for-smoothness-proof}
Let $X = \Spa(A, A^+)$ be a totally disconnected perfectoid space over $\Q_p(\mu_{p^\infty})$. Let $\mathfrak Y := \Spa A^+\langle T^{\pm1} \rangle$ be the relative torus, with generic fiber $Y \to X$. Then for all $N \in \Dqcohri(A^+/\pi[T^{\pm1}])$ there is a functorial isomorphism
\begin{align*}
	t_*(\ri^{+a}_Y/\pi) \tensor N \isom \IHom(t_*(\ri^{+a}_Y/\pi), N[-1])
\end{align*}
in $\Dqcohri(A^+/\pi[T^{\pm1}])$.
\end{lemma}
\begin{proof}
For ease of notation denote $B^{(+)}_n := A^{(+)}\langle T^{\pm1/p^n}\rangle$ and $\tilde B^{(+)} := A^{(+)}\langle T^{\pm1/p^\infty}\rangle$. Let us further denote $\Delta := \Hom(\Z[1/p]/\Z, \mu_{p^\infty}(\Q_p(\mu_{p^\infty})) = \Z_p(1)$. Then $\Delta$ acts continuously on $\tilde B^+$ by
\begin{align*}
	\gamma \times T^n = \gamma(n \bmod \Z) \, T^n \qquad \text{for $\gamma \in \Delta$ and $n \in \Z[1/p]$}
\end{align*}
and this action induces a $\Delta$-torsor $\Spa(\tilde B, \tilde B^+) \to Y$. By \cref{rslt:compute-Dqcohri-on-stack-quotient-of-p-bounded-space}, $\DqcohriX Y$ is equivalent to the $\infty$-category of smooth $\Delta$-representations on $(\tilde B^+/\pi)^a_\solid$-modules and $\Gamma(Y, -)$ computes smooth $\Delta$-cohomology. Note that $\Gamma(Y, -)$ viewed as a functor to $\Dqcohri(A^+/\pi)$ is the composition of $t_*$ and the forgetful functor $\Dqcohri(B_0^+/\pi) \to \Dqcohri(A^+/\pi)$, hence $t_*$ is computed by smooth $\Delta$-cohomology. By \cref{rslt:compute-smooth-group-cohom-via-smIHoms} smooth group cohomology agrees with classical continuous group cohomology on discrete objects. Since $\Delta \isom \Z_p$, it is well-known how to compute the cohomology:
\begin{align*}
	t_*(\ri^{+a}_Y/\pi) &= \Gamma(\Delta, \tilde B^{+a}/\pi)\\
	&\isom [\tilde B^{+a}/\pi \xto{1 - \gamma} \tilde B^{+a}/\pi]\\
	&= \bigdsum_{n\in \Z[\frac1p]/\Z} [B^+_0/\pi \xto{1 - \gamma(n)} B^+_0/\pi]
\end{align*}
where $\gamma$ is any fixed generator of $\Delta$ and the complex on the right-hand side sits in cohomological degrees $0$ and $1$.
We now obtain a natural morphism
\begin{align*}
	t_*(\ri^{+a}_Y/\pi) \tensor t_*(\ri^{+a}_Y/\pi) \to t_*(\ri^{+a}_Y/\pi) \to H^1(t_*(\ri^{+a}_Y/\pi))[-1] \surjto B^+_0/\pi[-1](-1),
\end{align*}
where the last morphism is the projection to the summand for $n = 0$.
After tensoring with $N$, this induces a morphism
\begin{align*}
	\varphi\colon t_*(\ri^{+a}_Y/\pi) \tensor N \to \IHom(t_*(\ri^{+a}_Y/\pi), N[-1](-1))
\end{align*}
in $\Dqcohri(B^+_0/\pi)$, which we claim to be an isomorphism. Assuming for a moment that $N = B^{+a}_0/\pi$, one checks that for each $n \in \Z[1/p]/\Z$, $\varphi$ maps the $n$-part of $t_*(\ri^{+a}_Y/\pi)$ isomorphically to the dual of the $-n$-part of $t_*(\ri^{+a}_Y/\pi)$. Thus $\varphi$ is essentially given by the natural map
\begin{align*}
	\bigdsum_{n\in \Z[\frac1p]/\Z} [B^+_0/\pi \xto{1 - \gamma(n)} B^+_0/\pi] \to \prod_{n\in \Z[\frac1p]/\Z} [B^+_0/\pi \xto{1 - \gamma(n)} B^+_0/\pi].
\end{align*}
Note that the cohomology groups of the $n$-part in the above sum and product is killed by pseudouniformizers which get arbitrarily close to $1$ for increasing powers of $p$ in the denominator of $n$. This implies that the above map from sum to product is indeed an (almost) isomorphism, proving that $\varphi$ is an isomorphism.

For general $N$ it is now enough to show that we can pull $N$ out of the $\IHom$. But by the same argument as before, this $\IHom$ is computed as a direct sum of $\IHom$'s for all $n$-parts, and $N$ can obviously be pulled out of each of these $\IHom$'s.
\end{proof}

We finally arrive at the first main result of the present subsection, which is also one of the main applications of our 6-functor formalism:

\begin{theorem} \label{rslt:adic-smooth-maps-over-Qp-are-p-cohom-smooth}
Let $f\colon Y \to X$ be a smooth map of analytic adic spaces over $\Q_p$. Then the induced map $f^\diamond\colon Y^\diamond \to X^\diamond$ of diamonds is $p$-cohomologically smooth.
\end{theorem}
\begin{proof}
In the following we will drop the diamond superscript and simply denote by $X$, $Y$ and $f$ their counterparts in the world of diamonds. By \cref{rslt:finite-type-implies-p-bounded} $f$ is $p$-bounded, hence bdcs. We now show that $f$ is $p$-cohomologically smooth. By \cref{rslt:stability-of-p-cohom-smooth} this reduces to the case that $X = \Spa(A, A^+)$ is a strictly totally disconnected perfectoid space and then, as $f$ is locally étale over a relative torus, to the case that $f$ is a relative torus. Finally, by stability of $p$-cohomological smoothness under compositions we reduce to the case that $Y = \Spa(A\langle T^{\pm1} \rangle, A^+\langle T^{\pm1} \rangle)$ is a relative torus of dimension $1$. Now pick a pseudouniformizer $\pi$ on $X$. Given any $M \in \Dqcohri(A^+/\pi)$ we will show that the natural morphism
\begin{align*}
	f^!(A^{+a}/\pi) \tensor f^* M \to f^! M
\end{align*}
is an isomorphism in $\DqcohriX Y$ and that $f^!(A^{+a}/\pi) \isom \ri^{+a}_Y/\pi[2]$. Note that $f$ is the generic fiber of the map $\mathfrak f\colon \mathfrak Y = \Spf A^+\langle T^{\pm1} \rangle \to \Spf A^+$ of formal schemes. Thus by \cref{rslt:compute-lower-shriek-via-pushforward-to-formal-model} we get $f_! = \mathfrak f_{\pi!} \comp t_*$ (with notation as in that result). By adjunctions it follows that providing a map $\ri^{+a}_Y/\pi[2] \to f^!(A^{+a}/\pi)$ is equivalent to providing a map $t_* \ri^{+a}_Y/\pi[2] \to \mathfrak f_\pi^!(A^{+a}/\pi) \isom A^{+a}/\pi[1]$ (where the computation of $\mathfrak f_\pi^!$ follows from \cref{rslt:scheme-6-functor-poincare-duality,rslt:lower-shriek-commutes-with-almost-localization}). Such a map exists by the explicit computations in \cref{rslt:compute-lower-shriek-via-pushforward-to-formal-model} (it is given by a projection map on the first homology groups). It now remains to show that the induced map
\begin{align*}
	f^* M[2] \to f^!(A^{+a}/\pi) \tensor f^* M \to f^! M
\end{align*}
is an isomorphism.

Let $\tilde Y = \varprojlim_n Y_n \to Y$ be the standard pro-étale tower, where $Y_n$ is obtained by adjoining the $p^n$-th root of $T$ (the maps $Y_n \to Y$ are étale because we are working over $\Q_p$). As $\tilde Y \surjto Y$ is surjective, it is enough to verify the above isomorphism after pullback to $\tilde Y$. Then, since $\tilde Y$ is a $p$-bounded affinoid perfectoid space, by \cref{rslt:compute-Dqcohri-for-p-bounded-over-tot-disc} it is enough to verify the isomorphism after applying $\Gamma(\tilde Y, -)$. By \cref{rslt:colimit-of-cohomology-qproet-unbounded-version} it is even enough to verify the isomorphism after applying $\Gamma(Y_n, -)$ for all $n$. Each $Y_n$ is itself a relative torus of dimension $1$ over $X$, so the claim boils down to showing the isomorphism after applying $\Gamma(Y, -)$, i.e. we need to see that the morphism
\begin{align*}
	\Gamma(Y, f^* M[2]) \to \Gamma(Y, f^! M)
\end{align*}
is an isomorphism in $\Dqcohri(A^+/\pi)$. Let us compute the right-hand side. Note that
\begin{align*}
	\Gamma(Y, f^! M) = \Hom(\ri^{+a}_Y/\pi, f^! M)^a = \Hom(f_!(\ri^{+a}_Y/\pi), M)^a
\end{align*}
By \cref{rslt:compute-lower-shriek-via-pushforward-to-formal-model} we get
\begin{align*}
	\Gamma(Y, f^! M) = \Hom(\mathfrak f_{\pi!} t_* (\ri^{+a}_Y/\pi), M)^a = \IHom(t_*(\ri^{+a}_Y/\pi), \mathfrak f_\pi^! M).
\end{align*}
By Poincaré duality in the setting of schemes (see \cref{rslt:scheme-6-functor-poincare-duality,rslt:lower-shriek-commutes-with-almost-localization}) we have $\mathfrak f_\pi^! M \isom \mathfrak f_\pi^* M[1]$. Thus altogether we obtain
\begin{align*}
	\Gamma(Y, f^! M) \isom \IHom(t_*(\ri^{+a}_Y/\pi), \mathfrak f_\pi^* M)[1].
\end{align*}
We now compute $\Gamma(Y, f^* M[2])$: By the projection formula for $t_*$ (see \cref{rslt:projection-formula-for-t}) we have
\begin{align*}
	\Gamma(Y, f^*M) = t_* f^* M = t_*(\ri^{+a}_Y/\pi) \tensor \mathfrak f_\pi^* M.
\end{align*}
Thus the desired isomorphism $\Gamma(Y, f^*M[2]) \isom \Gamma(Y, f^! M)$ follows from \cref{rslt:Koszul-duality-for-smoothness-proof}.
\end{proof}

In order to deduce Poincaré duality for $\Fld_p$-cohomology on smooth rigid-analytic varieties from \cref{rslt:adic-smooth-maps-over-Qp-are-p-cohom-smooth} it remains to compute the ``dualizing complex'' $f^!(\ri^{+a}_X/p)$. More specifically we need to show that if $f$ has pure dimension $d$ then $f^!(\ri^{+a}_X/p) = \ri^{+a}_Y/p(d)[2d]$. The proof will rely on the following version of \cref{rslt:p-cohom-smoothness-identities}:

\begin{proposition} \label{rslt:p-cohom-smoothness-identities-for-immersion}
Let $X, Y \to S$ be bdcs and $p$-cohomologically smooth maps of locally spatial diamonds in $\vStackspip$ and let $f\colon Y \to X$ be an $S$-morphism.
\begin{propenum}
	\item \label{rslt:p-cohom-smoothness-description-of-upper-shriek-for-immersion} For all invertible $\mathcal L \in \DqcohriX X^\varphi$ the natural morphism
	\begin{align*}
		f^!(\ri^{+a}_X/\pi) \tensor f^* \mathcal L \isoto f^! \mathcal L
	\end{align*}
	is an isomorphism and $f^!(\ri^{+a}_X/\pi)$ is invertible.

	\item \label{rslt:p-cohom-smoothness-base-change-for-immersion} Let $S' \to S$ be a map of locally spatial diamonds and let $f'\colon Y' \to X'$ denote the base-change of $f$ along $S' \to S$, so that the following diagram is cartesian:
	\begin{center}\begin{tikzcd}
		Y' \arrow[r,"g'"] \arrow[d,"f'"] & Y \arrow[d,"f"]\\
		X' \arrow[r,"g"] & X
	\end{tikzcd}\end{center}
	Then for all invertible $\mathcal L \in \DqcohriX X^\varphi$ the natural morphism
	\begin{align*}
		g'^* f^! \mathcal L \isoto f'^! g^* \mathcal L
	\end{align*}
	is an isomorphism.
\end{propenum}
\end{proposition}
\begin{proof}
We first prove (i). To show that $f^!(\ri^{+a}_X/\pi)$ is invertible, let us denote $h_X\colon X \to S$ and $h_Y\colon Y \to S$ the two projections. Then by \cref{rslt:dualizing-complex-automatically-given-by-F-p-locsys} there is an étale cover $(U_i \to X)_{i\in I}$ such that the restriction of $h_X^!(\ri^{+a}_S/\pi)$ to each $U_i$ is free up to a shift. Since the condition of being invertible is étale local we can replace $X$ by $U_i$ and $Y$ by $Y \cprod_X U_i$ in order to reduce to the case that $h_X^!(\ri^{+a}_S/\pi) \isom \ri^{+a}_X/\pi[n]$ for some $n \in \Z$. Then
\begin{align*}
	f^!(\ri^{+a}_X/\pi) = f^!(\ri^{+a}_X/\pi[n])[-n] \isom f^! h_X^!(\ri^{+a}_S/\pi)[-n] = h_Y^!(\ri^{+a}_S/\pi)[-n],
\end{align*}
which is invertible because $h_Y$ is $p$-cohomologically smooth. Now let $\mathcal L \in \DqcohriX X^\varphi$ be a given invertible sheaf; we need to show that the map $f^!(\ri^{+a}_X/\pi) \tensor f^* \mathcal L \isoto f^! \mathcal L$ is an isomorphism. By \cref{rslt:global-Riemann-Hilbert} $\mathcal L$ comes from an invertible $\Fld_p$-sheaf on $X$ and is in particular étale locally free up to shift. Since the claimed isomorphism is compatible with étale pullback we can pass to this étale cover of $X$ (and the induced étale cover of $Y$) to reduce to the case that $\mathcal L = \ri^{+a}_X/\pi$. But then the claimed isomorphism is clear.

We now prove (ii), so let $S'$, $X'$, $Y'$ and $\mathcal L$ be given as in the claim, and denote $h_X\colon X \to S$ and $h_Y\colon Y \to S$ the projections, with pullback $h_{X'}\colon X' \to S'$ and $h_{Y'}\colon Y' \to S'$. As in the proof of (i) we can pass to an étale cover of $X$ in order to assume that $\mathcal L$ and $h_X^!(\ri^{+a}_S/\pi)$ are free up to shift. After applying a shift to $\mathcal L$ (which is harmless) we reduce to the case that $\mathcal L = h_X^!(\ri^{+a}_S/\pi)$. Now consider the following diagram of pullback squares:
\begin{center}\begin{tikzcd}
	Y' \arrow[r,"g'"] \arrow[d,"f'"] \arrow[dd,bend right,"h_{Y'}",swap] & Y \arrow[d,swap,"f"] \arrow[dd,bend left,"h_Y"]\\
	X' \arrow[r,"g"] \arrow[d,"h_{X'}"] & X \arrow[d,swap,"h_X"]\\
	S' \arrow[r,"g_S"] & S
\end{tikzcd}\end{center}
Since $h_X$ and $h_Y$ are $p$-cohomologically smooth we can apply \cref{rslt:p-cohom-smoothness-base-change-for-upper-shriek} in order to get
\begin{align*}
	&g'^* f^! \mathcal L = g'^* f^! h_X^! (\ri^{+a}_S/\pi) = g'^* h_Y^! (\ri^{+a}_S/\pi) = h_{Y'}^! g_S^* (\ri^{+a}_S/\pi) = f'^! h_{X'}^! g_S^* (\ri^{+a}_S/\pi) =\\&\qquad= f'^! g^* h_X^! (\ri^{+a}_S/\pi) = f'^! g^* \mathcal L,
\end{align*}
as desired.
\end{proof}

The following lemma provides another situation in which the conclusion of \cref{rslt:smooth-upper-shriek-base-changes} holds. This computation will be applied in a very specific situation in \cref{rslt:dualizing-complex-for-smooth-maps-over-Qp}.

\begin{lemma} \label{rslt:upper-shriek-base-change-for-smooth-var-and-split-surjection}
Let $C$ be an algebraically closed non-archimedean field over $\Q_p$ and let
\begin{center}\begin{tikzcd}
	Y' \arrow[r,"g'"] \arrow[d,"f'"] & Y \arrow[d,"f"]\\
	X' \arrow[r,"g"] & X
\end{tikzcd}\end{center}
be a cartesian diagram of connected smooth rigid-analytic varieties over $C$ such that $f$ is split surjective and $\dim X - \dim X' = \dim Y - \dim Y'$. Then for every invertible $\mathcal L \in \Dqcohri(\ri^+_X/p)^\varphi$ the natural morphisms
\begin{align*}
	g'^* f^! \mathcal L &\isoto f'^! g^* \mathcal L,\\
	f'^* g^! \mathcal L &\isoto g'^! f^* \mathcal L
\end{align*}
are isomorphisms.
\end{lemma}
\begin{proof}
We first prove the second isomorphism. By the proof of \cref{rslt:adic-smooth-maps-over-Qp-are-p-cohom-smooth} we know that for a smooth rigid-analytic variety $Z$ over $C$ the map $Z \to \Spa C$ is $p$-cohomologically smooth and the $p$-adic dualizing complex is concentrated in cohomological degree $-2 \dim Z$. By \cref{rslt:p-cohom-smoothness-description-of-upper-shriek-for-immersion} both $f'^* g^! \mathcal L$ and $g'^! f^* \mathcal L$ are invertible and from the proof of that statement and the assumption on dimensions it follows that they are concentrated in the same cohomological degree $\dim X - \dim X'$. Employing also \cref{rslt:global-Riemann-Hilbert} we find that the map $f'^* g^! \mathcal L \to g'^! f^* \mathcal L$ is essentially a map of $\Fld_p$-local systems on $Y'$. Since $Y'$ is connected this map of local systems is automatically an isomorphism as soon as it is non-zero. Assume for a contradiction that this map is zero; then by adjunctions and proper base-change also the natural map $f^* g_! g^! \mathcal L \to f^* \mathcal L$ is zero. Since $f$ is split surjective it follows that the map $g_! g^! \mathcal L \to \mathcal L$ is zero. By adjunction this means that the identity $g^! \mathcal L \to g^! \mathcal L$ is zero, i.e. that $g^! \mathcal L$ is zero. On the other hand we know by \cref{rslt:p-cohom-smoothness-description-of-upper-shriek-for-immersion} that $g^! \mathcal L$ is invertible; contradiction! This finishes the proof of the second isomorphism.

It remains to prove the first isomorphism. But this follows from the second isomorphism in the same way as in the proof of \cref{rslt:smooth-upper-shriek-base-changes} using \cref{rslt:p-cohom-smoothness-description-of-upper-shriek-for-immersion}.
\end{proof}

Finally we are in the position to compute the dualizing complex $f^!(\ri^{+a}_X/p)$ of a smooth map $f\colon Y \to X$ of analytic adic spaces over $\Q_p$. The following proof is somewhat technical and may be a bit hard to grasp on a first read, so we refer the reader to the following \cref{rmk:proof-strategy-for-dualizing-complex} for an overview of the proof strategy.

\begin{theorem} \label{rslt:dualizing-complex-for-smooth-maps-over-Qp}
Let $f\colon Y \to X$ be a smooth map of analytic adic spaces over $\Q_p$. If $f$ has pure dimension $d$ then there is a natural isomorphism
\begin{align*}
	f^{\diamond!}(\ri^{+a}_X/p) = \ri^{+a}_Y/p(d)[2d].
\end{align*}
\end{theorem}
\begin{proof}
To improve the readability of the proof we make the following simplifications to the notation: Firstly, we view every analytic adic space over $\Spa\Q_p$ implicitly as a locally spatial diamond. Secondly, we use \cref{rslt:global-Riemann-Hilbert} to identify invertible $\ri^{+a}_S/p$-sheaves with invertible $\Fld_p$-sheaves on every such locally spatial diamond $S$. In particular we write expressions like $f^! \Fld_p$ for $f^!(\ri^{+a}_X/p)$. This notation is justified because in all relevant situations the functor $f^!$ restricts to a functor on invertible $\Fld_p$-sheaves by \cref{rslt:p-cohom-smoothness-description-of-upper-shriek-for-immersion}.

Let $\mathcal C$ be the category whose objects are smooth maps $f\colon Y \to X$ of analytic adic spaces over $\Q_p$ of pure dimension $d$ and whose morphisms are commuting squares
\begin{center}\begin{tikzcd}
	Y' \arrow[r,"g'"] \arrow[d,"f'"] & Y \arrow[d,"f"]\\
	X' \arrow[r,"g"] & X
\end{tikzcd}\end{center}
such that the induced map $Y' \to X' \cprod_X Y$ is étale. By (the proof of) \cref{rslt:adic-smooth-maps-over-Qp-are-p-cohom-smooth} we know that for every $f\colon Y \to X$ in $\mathcal C$ the ``dualizing complex'' $\omega_f := f^! \Fld_p$ is given by an invertible $\Fld_p$-sheaf concentrated in cohomological degree $-2d$. Using \cref{rslt:p-cohom-smoothness-base-change-for-upper-shriek} we observe that the assignment
\begin{align*}
	\omega\colon \mathcal C^\opp \to \D(\Fld_p)^\heartsuit, \qquad [f\colon Y \to X] \mapsto H^{-2d}(Y, \omega_f)
\end{align*}
defines a functor. The goal is to show that this functor is naturally isomorphic to the functor
\begin{align*}
	\omega'\colon \mathcal C^\opp \to \D(\Fld_p)^\heartsuit, \qquad [f\colon Y \to X] \mapsto H^0(Y, \Fld_p(d)).
\end{align*}
Let $\mathcal C_t \subset \mathcal C$ be the full subcategory of those maps $f\colon Y \to X$ such that $f$ admits an étale map to a $d$-dimensional torus over $X$. Since every $f\colon Y \to X$ in $\mathcal C$ admits an étale covering $(Y_i \to Y)_i$ such that all the induced maps $Y_i \to X$ lie in $\mathcal C_t$, both $\omega$ and $\omega'$ are determined by their restriction to $\mathcal C_t$. In other words, it is enough to show that $\omega$ and $\omega'$ are naturally isomorphic on $\mathcal C_t$. Now let $\mathcal C_{tb} \subset \mathcal C_t$ be the full subcategory of those maps $f\colon Y \to X$ where $X$ is a strictly totally disconnected perfectoid space; it is enough to show that the restrictions of $\omega$ and $\omega'$ to $\mathcal C_{tb}$ are naturally isomorphic.

Let $\mathcal C_{tp} \subset \mathcal C_{tb}$ be the full subcategory of those maps $f\colon Y \to X$ where $X = \Spa(C, C^+)$ for some algebraically closed extension $C$ of $\Q_p$ and some open and bounded valuation subring $C^+ \subset C$. We claim that in order to show that $\omega$ and $\omega'$ are naturally isomorphic on $\mathcal C_{tb}$, it is enough to show that this is true on $\mathcal C_{tp}$. To see this, assume that the latter condition is satisfied. Then given any $f\colon Y \to X$ in $\mathcal C_{tb}$ we need to find an isomorphism $\omega(f) \isoto \omega'(f)$ which is functorial in $f$ and compatible with the given isomorphism of $\omega$ and $\omega'$ on $\mathcal C_{tp}$. By looking at stalks it becomes clear that there is at most one isomorphism $\omega(f) \isoto \omega'(f)$ which satisfies the latter condition, and should it exist then automatically satisfies the former condition. It thus remains to show that for fixed $f$ there is an isomorphism $\omega(f) \isoto \omega'(f)$ compatible with the isomorphism of $\omega$ and $\omega'$ on $\mathcal C_{tp}$. First assume that $Y = T^d_X = T^d \cprod X$ is the $d$-dimensional torus over $X$. Note that the assignment $S \mapsto \omega(T^d_S \to S)$ (for perfectoid $S$ over $\Q_p$) defines an invertible $\Fld_p$-sheaf on $\Spa\Q_p$ (here we use that by the proof of \cref{rslt:adic-smooth-maps-over-Qp-are-p-cohom-smooth} the dualizing complex on the relative torus is free), which by the isomorphism of $\omega$ and $\omega'$ on $\mathcal C_{tp}$ is isomorphic to the invertible $\Fld_p$-sheaf $S \mapsto \omega'(T^d_S \to S)$. In particular this provides the desired isomorphism $\omega(T^d_X \to X) \isoto \omega'(T^d_X \to X)$. This isomorphism is even $H^0(X, \Fld_p)$-linear (instead of merely being $\Fld_p$-linear) and since $H^0(X, \Fld_p) = H^0(T^d_X, \Fld_p)$ we deduce that we get an isomorphism of free sheaves $\omega_f[-2d] \isoto \Fld_p(d)$ on $T_X^d$. In particular, if now $f\colon Y \to X$ in $\mathcal C_{tb}$ is general, i.e. admits an étale map $Y \to T^d_X$, then we get an induced isomorphism $\omega(Y \to X) \isoto \omega'(Y \to X)$ which is still compatible with the isomorphism of $\omega$ and $\omega'$ on $\mathcal C_{tp}$, as desired.

We have reduced the proof to checking that $\omega$ and $\omega'$ are naturally isomorphic on $\mathcal C_{tp}$. For simplicity we can make two more adjustments: Firstly, from now on there is no point in restricting to those $f\colon Y \to X$ where $Y$ is étale over a relative torus anymore. Secondly, we can even assume that $X = \Spa(C, \ri_C)$ (i.e. no general $C^+$ is required) by a similar argument as in the first part of the proof of \cite[Theorem 5.1]{rigid-p-adic-hodge} (alternatively argue directly with tori as in the previous paragraph).

Altogether we have reduced the proof to the following claim: Let $\mathcal C' \subset \mathcal C$ be the full subcategory of those maps $f\colon Y \to X$ such that $X = \Spa(C, \ri_C)$ for some algebraically closed extension $C$ of $\Q_p$; then the functors $\omega$ and $\omega'$ are naturally isomorphic on $\mathcal C'$. In particular we have managed to put ourselves into the realm of classical rigid-analytic geometry. To construct the desired isomorphism, we perform further transformations of the problem at hand. The following argument was suggested by Scholze.

Fix some $f\colon Y \to X$ in $\mathcal C'$ and consider the diagonal $\Delta_f\colon Y \injto Y \cprod_X Y$ and the two projections $g_1, g_2\colon Y \cprod_X Y \to Y$. Using \cref{rslt:p-cohom-smoothness-description-of-upper-shriek-for-immersion} (which applies to $\Delta_f$, $g_1$ and $g_2$) and \cref{rslt:p-cohom-smoothness-base-change-for-upper-shriek} we compute
\begin{align*}
	&\omega_f = \Delta_f^! g_1^! \omega_f = \Delta_f^! (g_1^* \omega_f \tensor g_1^!\Fld_p) = \Delta_f^! (g_1^* \omega_f \tensor g_2^* \omega_f) = \Delta_f^* (g_1^* \omega_f \tensor g_2^* \omega_f) \tensor \Delta_f^!\Fld_p =\\&\qquad= \omega_f \tensor \omega_f \tensor \Delta_f^!\Fld_p.
\end{align*}
It follows that there is a natural isomorphism $\omega_f = (\Delta_f^! \Fld_p)^\vee$, so it is enough to compute the right-hand side. Note that $\Delta_f$ is a local complete intersection in the sense of \cite[Definition 5.4]{guo-li-derived-de-rham} because $Y$ and $Y \cprod_X Y$ are smooth over $X$ (e.g. by \cite[Lemma IV.4.14]{fargues-scholze-geometrization} this map is locally a section of a relative ball). Since $\Delta_f$ is also a closed immersion, it follows from \cite[Proposition 5.3.(1)]{guo-li-derived-de-rham} that $\Delta_f$ is locally given by a regular closed immersion (in the sense of algebraic geometry) of affinoid $C$-algebras. By adapting the standard construction of the ``deformation to the normal cone'' from algebraic geometry to the rigid-analytic world (to handle blow-ups in the rigid-analytic world see e.g. \cite{schoutens-blow-up}) we construct the following diagram of smooth rigid-analytic varieties over $X = \Spa C$:
\begin{center}\begin{tikzcd}
	Y \cprod_X \mathbb P^1_X \arrow[r,hook,"i"] \arrow[d] & \tilde Y \arrow[dl]\\
	\mathbb P^1_X
\end{tikzcd}\end{center}
Here the left vertical map is the projection to the second factor and the map $i$ is a closed immersion such that for a classical point $x \in \mathbb P^1_X$ the fiber
\begin{align*}
	i_x\colon Y = Y \cprod_X \{ x \} \injto \tilde Y_x := \tilde Y \cprod_{\mathbb P^1_X} \{ x \}
\end{align*}
has the following form: If $x \ne \infty$ then $\tilde Y_x = Y \cprod_X Y$ and $i_x = \Delta_f$, but if $x = \infty$ then $\tilde Y_\infty$ is the normal cone of $Y$ in $Y \cprod_X Y$ (which is a vector bundle over $Y$) and $i_\infty$ is the $0$-section. More informally, $\tilde Y$ is a family of maps parametrized by $\mathbb P^1_X$ which deforms the map $\Delta_f$ into the zero-section map of the normal cone. We will now use this deformation in order to reduce the computation of $\Delta_f^!\Fld_p$ to the computation of $i_\infty^! \Fld_p$. To make this reduction we first observe the following diagram of pullback squares:
\begin{center}\begin{tikzcd}
	Y \arrow[d,hook,"i_x"] \arrow[r,hook,"j_x"] \arrow[dd,"h_2'",bend right,swap] & Y \cprod_X \mathbb P^1_X \arrow[d,hook,"i",swap] \arrow[dd,"h_2",bend left]\\
	\tilde Y_x \arrow[r,hook,"j'_x"] \arrow[d,"h_1'"] & \tilde Y \arrow[d,"h_1",swap]\\
	\{ x \} \arrow[r,hook,"j''_x"] & \mathbb P^1_X
\end{tikzcd}\end{center}
We claim that for every invertible $\Fld_p$-sheaf $\mathcal L$ on $\tilde Y$ the natural morphism $j_x^* i^! \mathcal L \isoto i_x^! j_x'^* \mathcal L$ is an isomorphism. As in the proof of \cref{rslt:p-cohom-smoothness-base-change-for-immersion} it is enough to show that a similar claim holds for the lower square and the outer square. For the outer square this follows immediately from \cref{rslt:smooth-upper-shriek-base-changes} because $h_2$ is smooth. For the lower square we can apply \cref{rslt:upper-shriek-base-change-for-smooth-var-and-split-surjection}: After passing to an open cover of $Y$ we can assume that $Y$ is connected, in which case all the spaces occurring in the lower square are connected (for $\tilde Y_x$ use the explicit description provided above); also note that $h_1$ is split surjective, with a section given by $\mathbb P^1_X \to Y \cprod_X \mathbb P^1_X \xto{i} \tilde Y$, where the first map is induced by the choice of any fixed $C$-point in $Y$. This proves the claim $j_x^* i^! \mathcal L = i_x^! j_x'^* \mathcal L$.

Let us now denote $\Omega_f := i^!\Fld_p$, which is an invertible sheaf on $Y \cprod_X \mathbb P^1_X$. By the argument in the previous paragraph we know that for $x \in \mathbb P^1_X$ with induced map $j_x\colon Y \injto Y \cprod_X \mathbb P^1_X$ we have
\begin{align*}
	j_x^* \Omega_f = \begin{cases} \Delta_f^!\Fld_p, & \text{if $x \ne \infty$},\\ i_\infty^!\Fld_p, & \text{if $x = \infty$.} \end{cases}.
\end{align*}
On the other hand we note that via the projection $h\colon Y \cprod_X \mathbb P^1_X \to Y$, invertible $\Fld_p$-sheaves on $Y$ and on $Y \cprod_X \mathbb P^1_X$ are identified (the equivalence being the functors $h^*$ and $h_*$) -- in fact this is true for any small v-stack $Y$ over $X$: One can check this after pullback to a cover of $Y$ by strictly totally disconnected spaces and then pass to stalks in order to reduce to the case that $Y$ is a geometric point; then the claim follows from the fact that $\mathbb P^1$ is simply connected. It follows that there is a natural isomorphism
\begin{align*}
	\Delta_f^!\Fld_p = j_0^* \Omega_f = j_0^* h^* h_* \Omega_f = h_* \Omega_f = j_\infty^* h^* h_* \Omega_f = j_\infty^* \Omega_f = i_\infty^!\Fld_p
\end{align*}
of invertible sheaves on $Y$. Altogether we obtain
\begin{align*}
	\omega_f = (\Delta_f^!\Fld_p)^\vee = (i_\infty^!\Fld_p)^\vee.
\end{align*}
In other words, we are left with computing $i_\infty^!\Fld_p$. Now $i_\infty\colon Y \injto \tilde Y_\infty$ is the zero-section of a rank-$d$ vector bundle over $Y$, so it is enough to show the following: Let $\mathcal C''$ be the category fibered in groupoids over the category of locally spatial diamonds in $\Q_p$ which for to every such diamond $S$ associates the groupoid of rank-$d$ vector bundles $E \to S$. We obtain a functor
\begin{align*}
	\chi\colon \mathcal C''^\opp \to \D(\Fld_p)^\heartsuit, \qquad [E \to S] \mapsto H^{2d}(S, s_E^!\Fld_p),
\end{align*}
where $s_E\colon S \injto E$ denotes the zero-section (the required functoriality for $\chi$ follows from \cref{rslt:p-cohom-smoothness-base-change-for-immersion} by observing that locally every vector bundle is of the form $S \cprod \mathbb A^d$). We need to see that $\chi$ is naturally isomorphic to the functor
\begin{align*}
	\chi'\colon \mathcal C''^\opp \to \D(\Fld_p)^\heartsuit, \qquad [E \to S] \mapsto H^0(S, \Fld_p(-d)).
\end{align*}
Note that $\mathcal C''$ corresponds to the classifying stack $\Spa\Q_p/\GL_d$, hence $\chi$ and $\chi'$ define $\Fld_p$-local systems on $\Spa\Q_p/\GL_d$. By the usual arguments involving the dual version of the Barr-Beck theorem (applied to the pullback-pushforward adjunction along the projection $\Spa\Q_p \surjto \Spa\Q_p/\GL_d$) it follows that an $\Fld_p$-local system $\mathcal L$ on $\Spa\Q_p/\GL_d$ is the same as an $\Fld_p$-local system $\mathcal L_0$ on $\Spa \Q_p$ together with a co-action $\mathcal L \to \mathcal L \tensor H^0(\GL_d, \Fld_p)$. Since $\GL_d$ is connected and thus $H^0(\GL_d, \Fld_p) = \Fld_p$ we deduce that the categories of $\Fld_p$-local systems on $\Spa\Q_p$ and $\Spa\Q_p/\GL_d$ are naturally isomorphic. In particular it is enough to find a natural isomorphism of the pullbacks of $\chi$ and $\chi'$ to $\Spa\Q_p$.

The map $\Spa\Q_p \to \Spa\Q_p/\GL_d$ is given by associating to every locally spatial diamond $S$ over $\Spa\Q_p$ the trivial vector bundle $S \cprod \mathbb A^d$. In particular the pullback $\chi_0$ of $\chi$ to $\Spa\Q_p$ is given by the functor $S \mapsto H^{2d}(S, s_{S \cprod \mathbb A^d}^!\Fld_p)$, i.e. we have $\chi_0 = s^!\Fld_p[2d]$, where $s\colon \Spa\Q_p \injto \mathbb A^d_{\Q_p}$ is the embedding of the origin. By choosing a suitable open embedding of the torus $T^d_{\Q_p} = \Spa \Q_p\langle T_1^\pm, \dots, T_d^\pm \rangle$ into affine space we get $\chi_0 = t^!\Fld_p[2d]$, where $t\colon \Spa\Q_p \injto T^d_{\Q_p}$ is the embedding of $1$. Let $f\colon T^d_{\Q_p} \to \Spa\Q_p$ be the projection. Then by \cref{rslt:p-cohom-smoothness-description-of-upper-shriek-for-immersion} we have
\begin{align*}
	\Fld_p = t^! f^! \Fld_p = t^!\Fld_p \tensor t^* f^!\Fld_p = \chi_0[2d] \tensor t^* f^!\Fld_p,
\end{align*}
which implies that
\begin{align*}
	\chi_0[2d] = (t^* f^!\Fld_p)^\vee.
\end{align*}
Therefore we have reduced the computation of $f^!\Fld_p$ for general smooth $f\colon Y \to X$ to the special case $f\colon T^d_{\Q_p} \to \Spa\Q_p$. This reduces further to the case $d = 1$. But the computation in the torus case was done in the proof of \cref{rslt:adic-smooth-maps-over-Qp-are-p-cohom-smooth} over algebraically closed fields and by carefully keeping track of Galois actions we also deduce the desired result over $\Q_p$.
\end{proof}

\begin{remark} \label{rmk:proof-strategy-for-dualizing-complex}
While the proof of \cref{rslt:dualizing-complex-for-smooth-maps-over-Qp} is somewhat lengthy and technical, it is surprisingly formal; in fact, the same proof strategy can be applied in different 6-functor formalisms and even in different geometric situations. For the convenience of the reader we summarize the main line of argument:
\begin{enumerate}[(a)]
	\item We first make several abstract reductions in order to reduce the computation of $f^!\Fld_p$ to the case that $X = \Spa C$ for some algebraically closed extension of $\Q_p$. This step may not be strictly necessary, but it allows us to work in the comfortable setting of classical rigid-analytic geometry.

	\item Denoting $\Delta_f\colon Y \injto Y \cprod_X Y$ the diagonal, a quick computation shows that $f^! \Fld_p = (\Delta_f^! \Fld_p)^\vee$, so it is enough to compute $\Delta_f^! \Fld_p$.

	\item By deforming the regular closed immersion $\Delta_f$ to its normal cone $i_\infty\colon Y \injto \tilde Y_\infty$ we deduce $\Delta_f^! \Fld_p = i_\infty^! \Fld_p$, so it is enough to compute $i_\infty^! \Fld_p$. Here we use the fact that every invertible $\Fld_p$-sheaf on $\mathbb P^1$ is trivial.

	\item We note that the normal cone $\tilde Y_\infty$ is a vector bundle over $Y$ and $i_\infty$ is the zero-section of that vector bundle. Therefore the computation of $i_\infty^! \Fld_p$ is a special case of computing $s_E^! \Fld_p$ along the zero section $s_E\colon S \injto E$ of a general diamond $S$ over $\Q_p$ with rank-$d$ vector bundle $E \to S$. The vector bundles on diamonds over $\Q_p$ are classified by the stack $\Spa\Q_p/\GL_d$ and one sees easily that the formula $[E \to S] \mapsto H^{2d}(S, s_E^! \Fld_p)$ defines an invertible $\Fld_p$-sheaf $\chi$ on this stack. It remains to compute $\chi$.

	\item Since $H^0(\GL_d, \Fld_p) = \Fld_p$, a standard argument shows that invertible $\Fld_p$-sheaves on $\Spa\Q_p$ and $\Spa\Q_p/\GL_d$ are naturally equivalent. It is thus enough to determine the pullback $\chi_0$ of $\chi$ to $\Spa\Q_p$. One sees immediately that $\chi_0 = s^! \Fld_p[2d]$, where $s\colon \Q_p \to \mathbb A^d_{\Q_p}$ is the embedding of the origin. The computation of this special case was essentially done in the proof of \cref{rslt:adic-smooth-maps-over-Qp-are-p-cohom-smooth}.
\end{enumerate}
\end{remark}

Combining \cref{rslt:dualizing-complex-for-smooth-maps-over-Qp} with the results of \cref{sec:ri-pi.phi-mod} we deduce the following application to $\Fld_p$-cohomology:

\begin{corollary} \label{rslt:relative-F-p-Poincare-duality}
Let $K$ be a non-archimedean field extension of $\Q_p$ and let $f\colon Y \to X$ be a proper smooth map of rigid-analytic varieties over $K$. Assume that $f$ has pure dimension $d$. Then for every $\Fld_p$-local system $\mathcal L$ on $Y$ there is a natural isomorphism
\begin{align*}
	f_* \mathcal L^\vee (d)[2d] = (f_* \mathcal L)^\vee
\end{align*}
of perfect $\Fld_p$-sheaves on $X$.
\end{corollary}
\begin{proof}
By \cref{rslt:adic-smooth-maps-over-Qp-are-p-cohom-smooth,rslt:dualizing-complex-for-smooth-maps-over-Qp} $f$ is $p$-cohomologically smooth and $f^! (\ri^{+a}_X/p) = (\ri^{+a}_Y/p)(d)[2d]$. Using the Primitive Comparison Theorem (see \cref{rslt:primitive-comparison}) and \cref{rslt:IHom-adjunction-for-shriek-functors} we deduce
\begin{align*}
	f_* \mathcal L^\vee (d)[2d] \tensor_{\Fld_p} \ri^{+a}_X/p &= f_*(\IHom(\mathcal L, \Fld_p(d)[2d]) \tensor_{\Fld_p} \ri^{+a}_Y/p)\\
	&= f_* \IHom(\mathcal L \tensor_{\Fld_p} \ri^{+a}_Y/p, f^!(\ri^{+a}_X/p))\\
	&= \IHom(f_*(\mathcal L \tensor_{\Fld_p} \ri^{+a}_Y/p), \ri^{+a}_X/p)\\
	&= \mathcal L^\vee \tensor_{\Fld_p} \ri^{+a}_X/p.
\end{align*}
The above computation is compatible with the natural $\varphi$-module structures, so we can pass to $\varphi$-invariants to conclude by \cref{rslt:global-Riemann-Hilbert}.
\end{proof}

\begin{corollary} \label{rslt:Fp-Poincare-duality-over-field}
Let $K$ be a non-archimedean field extension of $\Q_p$ and let $X$ be a proper smooth rigid-analytic variety of pure dimension $d$ over $K$. Then for all $i \in \Z$ there is a natural perfect pairing
\begin{align*}
	H^i_\et(X, \Fld_p) \tensor_{\Fld_p} H^{2d - i}_\et(X, \Fld_p) \to \Fld_p(-d).
\end{align*}
\end{corollary}
\begin{proof}
Apply \cref{rslt:relative-F-p-Poincare-duality} to $f\colon X \to \Spa K$ and $\mathcal L = \Fld_p$.
\end{proof}

\clearpage
\appendix
\section{Appendix: \texorpdfstring{$\infty$}{Infinity}-Categories} \label{sec:infcat}

In the appendix we gather several results about $\infty$-categories which we have not found in the literature. More concretely, the appendix is structured as follows. In \cref{sec:infcat.derived} we show functoriality of the derived $\infty$-category and compare derived descent with abelian descent. In \cref{sec:infcat.presentable} we prove some additional results on presentable $\infty$-categories and in particular deduce that filtered colimits commute with finite limits of $\infty$-categories. In \cref{sec:infcat.sheaves} we gather a few basic results about sheaves in the $\infty$-categorical setting, including very explicit descriptions for sheaves on nice sites. In \cref{sec:infcat.enriched} we study enriched $\infty$-categories and prove basic results like enriched Yoneda in that setting. In \cref{sec:infcat.sixfun} we introduce a general formalism for capturing 6-functor formalisms in the $\infty$-categorical setting and provide some tools to construct them.

\subsection{Derived \texorpdfstring{$\infty$}{Infinity}-Categories} \label{sec:infcat.derived}

We collect a few straightforward results about derived $\infty$-categories. In the following let $\catcat \subset \infcatinf$ denote the full subcategory spanned by the ($1$-)categories. The first result shows that the construction $\mathcal A \mapsto \D^+(\mathcal A)$ is functorial in an appropriate sense:

\begin{lemma} \label{rslt:D+-is-a-functor}
Let $\catcat^{AbI,ex} \subset \catcat$ be the subcategory of abelian categories with enough injectives and exact functors. Then the construction from \cite[Definition 1.3.2.7]{lurie-higher-algebra} defines a functor
\begin{align*}
	\D^+\colon \catcat^{AbI,ex} \to \infcatinf, \qquad \mathcal A \mapsto \D^+(\mathcal A).
\end{align*}
\end{lemma}
\begin{proof}
Let us denote by $\infcatinf^{st,tex}$ the $\infty$-category of stable $\infty$-categories equipped with $t$-structures whose heart has enough injectives, where morphisms are $t$-exact functors. We can easily construct the functor $(-)^\heartsuit\colon \infcatinf^{st,tex} \to \catcat^{AbI,ex}$, $\mathcal C \to \mathcal C^\heartsuit$. Now let $\mathcal X \subset \infcatinf^{st,tex}$ be the full subcategory spanned by those $\infty$-categories which are equivalent to $D^+(\mathcal A)$ for some $\mathcal A \in \catcat^{AbI,ex}$. By \cite[Theorem 1.3.3.2]{lurie-higher-algebra} the restriction of the functor $(-)^\heartsuit$ to $\mathcal X$ is fully faithful, hence an equivalence. But then taking any quasi-inverse produces the desired functor $\D^+$.
\end{proof}

We now prove an abstract descent result. It shows that derived descent along exact functors is equivalent to deriving the abelian descent.

\begin{proposition} \label{rslt:derived-descent-from-abelian-descent}
Let $\mathcal A^\bullet$ be a cosimplicial object of $\catcat$ with limit $\mathcal A := \varprojlim_{n\in\Delta} \mathcal A^n$ in $\catcat$. Assume that the following hold:
\begin{enumerate}[(a)]
    \item All $\mathcal A^n$ are Grothendieck abelian categories and all transition functors in the diagram $\mathcal A^\bullet$ are exact and preserve colimits.
    \item For every morphism $\alpha\colon [n] \to [m]$ in $\Delta$ let $\alpha^*\colon \mathcal A^n \to \mathcal A^m$ denote the associated functor and let $\alpha_*\colon \mathcal A^m \to \mathcal A^n$ be a right adjoint. Then $\alpha_*$ is exact.
	\item Letting $d_n^i\colon [n] \to [n+1]$ denote the $i$th coface map, then for all $\alpha\colon [n] \to [m]$ the following diagram satisfies base-change:
    \begin{center}\begin{tikzcd}
        \mathcal A^{m+1} \arrow[d,swap,"d^0_{m*}"] & \arrow[l,swap,"\alpha'^*"] \arrow[d,"d^0_{n*}"] \mathcal A^{n+1}\\
        \mathcal A^m & \arrow[l,"\alpha^*"] \mathcal A^n
    \end{tikzcd}\end{center}
    In other words the natural morphism of functors $\alpha'^* d^0_{n*} \isoto d^0_{m*} \alpha^*$ is an equivalence. (Here $\alpha'\colon [n+1] \to [m+1]$ is defined by $\alpha'(0) = 0$, $\alpha'(k+1) = \alpha(k)+1$ for $k \ge 0$.)
\end{enumerate}
Then:
\begin{propenum}
    \item \label{rslt:descent-computation-of-abelian-limit} The category $\mathcal A$ is abelian and admits the following description: An object of $\mathcal A$ is a pair $(X, \alpha)$, where $X$ is an object of $\mathcal A^0$ and $\alpha$ is an isomorphism $\alpha\colon d_0^{0*}X \isoto d_0^{1*}X$ in $\mathcal A^1$ such that the following diagram in $\mathcal A^2$ commutes:
    \begin{center}\begin{tikzcd}
        d_1^{0*} d_0^{0*} X \isom d_1^{1*} d_0^{0*} X \arrow[dr,swap,"d_1^{0*}\alpha"] \arrow[rr,"d_1^{1*}\alpha"] & & d_1^{1*} d_0^{1*} X \isom d_1^{2*} d_0^{1*} X \\
        & d_1^{0*} d_0^{1*} X \isom d_1^{2*} d_0^{0*} X \arrow[ur,swap,"d_1^{2*}\alpha"]
    \end{tikzcd}\end{center}
    A morphism $(X_1, \alpha_1) \to (X_2, \alpha_2)$ in $\mathcal A$ is a morphism $f\colon X_1 \to X_2$ such that $d_0^{1*}f \comp \alpha_1 = \alpha_2 \comp d_0^{0*}f$.

    \item Assume furthermore that $\mathcal A$ has enough injectives. Then there is a natural equivalence of $\infty$-categories
    \begin{align*}
        \D^+(\mathcal A) \isoto \varprojlim_{n\in\Delta} \D^+(\mathcal A^n),
    \end{align*}
    where the limit on the right is computed in $\infcatinf$ (use \cref{rslt:D+-is-a-functor} to produce the diagram $\D^+(\mathcal A^\bullet)$). Moreover, if $\mathcal A$ is Grothendieck and $\D(\mathcal A)$ and all $\D(\mathcal A^n)$ are left complete then the above equivalence extends to an equivalence
    \begin{align*}
    	\D(\mathcal A) \isoto \varprojlim_{n\in\Delta} \D(\mathcal A^n).
    \end{align*}
\end{propenum}
\end{proposition}
\begin{proof}
Part (i) can easily be checked by hand: Note first that the explicit description of $\mathcal A$ implies that $\mathcal A$ is abelian, so we only have to show that $\mathcal A$ has the claimed description. To make things easier we note that the full subcategory $\catcat \subset \infcatinf$ of ordinary categories is equivalent to a 2-category (this follows from \cite[Corollary 2.3.4.20]{lurie-higher-topos-theory} and \cite[Proposition 2.3.4.18]{lurie-higher-topos-theory}). It follows from the proof of \cite[Lemma 1.3.3.10]{lurie-higher-algebra} that $\varprojlim_{n\in\Delta} \mathcal A^n = \varprojlim_{n\in\Delta^{\le2}} \mathcal A^n$, i.e. the limit depends only on the first three terms. The rest is an explicit computation.

To prove (ii), note first that all transition functors in the \emph{augmented} cosimplicial category $\mathcal A \to \mathcal A^\bullet$ are exact and admit exact right adjoints. Indeed, we only need to check this for the functor $d_{-1}^{0*}\colon \mathcal A \to \mathcal A^0$. By the explicit description of $\mathcal A$ it is clear that this functor is exact. Its right adjoint $d_{-1*}^0\colon \mathcal A^0 \to \mathcal A$ can be computed as follows: For $X \in \mathcal A^0$ we have $d_{-1*}^0 X = (d_{0*}^0 d_0^{1*}X, \alpha)$, where $\alpha$ is the isomorphism
\begin{align*}
    d_0^{0*} (d_{0*}^0 d_0^{1*}X) \isom d_{1*}^0 d_1^{1*} d_0^{1*} X \isom d_{1*}^0 d_1^{2*} d_0^{1*} X \isom d_0^{1*} (d_{0*}^0 d_0^{1*} X).
\end{align*}
Note that in this chain of isomorphisms we used the base-change assumption of (c) and that some of the isomorphisms are not canonical (but can be chosen independently of $X$). It is straightforward to check that this construction does indeed provide a description of $d_{-1*}^0$. From the fact that $d_0^{1*}$ and $d_{0*}^0$ are exact one deduces easily that $d_{-1*}^0$ is exact. It also follows that the base change of (c) holds for $\alpha\colon \emptyset \to [n]$ as well: This immediately reduces to $n = 0$, where it follows directly from the explicit description of $d_{-1*}^0$.

Now apply the functor $\D^+$ from \cref{rslt:D+-is-a-functor} in order to get an augmented cosimplicial $\infty$-category $\D^+(\mathcal A) \to \D^+(\mathcal A^\bullet)$. By abuse of notation we still denote $\alpha^*\colon \D^+(\mathcal A^n) \to \D^+(\mathcal A^m)$ resp. $\alpha_*\colon \D^+(\mathcal A^m) \to \D^+(\mathcal A^n)$ the transition functor resp. the right adjoint of the transition functor associated to a map $\alpha\colon [n] \to [m]$ in $\Delta_+$. These are the induced functors from the corresponding functors in the augmented cosimplicial category $\mathcal A \to \mathcal A^\bullet$. Since the latter functors are all exact, the functors $\alpha^*$ and $\alpha_*$ can be computed directly on homology. It follows immediately that the base-change of (c) also holds in $\D^+(\mathcal A) \to \D^+(\mathcal A^\bullet)$ (even for $\alpha\colon \emptyset \to [n]$ by the previous paragraph).

To prove that the induced functor $\D^+(\mathcal A) \to \varprojlim_{n\in\Delta} \D^+(\mathcal A^n)$ is an equivalence we employ \cite[Corollary 4.7.5.3]{lurie-higher-algebra} (applied to the opposite categories). Condition (2) of that reference was shown above. It remains to prove that the functor $F := d_{-1}^{0*}\colon \D^+(\mathcal A) \to \D^+(\mathcal A^0)$ is conservative and preserves totalizations of $F$-split cosimplicial objects in $\D^+(\mathcal A)$. We can check that $F$ is conservative directly on homology and hence reduce to showing that the functor $d_{-1}^{0*}\mathcal A \to \mathcal A_0$ is conservative, which is clear by the explicit description of $\mathcal A$. The preservation of totalizations can be proved as in the proof of \cite[Proposition D.6.4.6]{lurie-spectral-algebraic-geometry}.

Finally the claim about unbounded derived categories is obtained by writing $\D(-) = \varprojlim_k \D_{\le k}(-)$ (by definition of left-completeness) and observing that $\D_{\le k}(A) \isoto \varprojlim_{n\in\Delta} \D_{\le k}(\mathcal A^n)$ is an equivalence by the $\D^+$-case. Passing to the limit over $k$ and commuting limits with limits produces the desired result.
\end{proof}

\subsection{Presentable \texorpdfstring{$\infty$}{Infinity}-Categories} \label{sec:infcat.presentable}

In this subsection we prove various results about presentable $\infty$-categories. Among others, we show that the $\infty$-category $\infcatinf$ of (small) $\infty$-categories is compactly generated and that this implies that for any regular cardinal $\kappa$, $\kappa$-small limits of $\infty$-categories commute with $\kappa$-filtered colimits. Moreover, we give criteria for detecting that an $\infty$-category is $\kappa$-compactly generated; most prominently we will show that if $\kappa$ is uncountable then a $\kappa$-small limit of $\kappa$-compactly generated $\infty$-categories is still $\kappa$-compactly generated.

We start with the following criterion for detecting $\kappa$-compactly generated $\infty$-categories.

\begin{lemma} \label{rslt:compact-generators-implies-compactly-generated}
Let $\mathcal C$ be an $\infty$-category with all small colimits, $\kappa$ a regular cardinal and $(X_i)_i$ a small family of $\kappa$-compact objects of $\mathcal C$ such that the family of functors $\Hom(X_i, -)\colon \mathcal C \to \Ani$ is conservative. Then:
\begin{lemenum}
	\item $\mathcal C$ is $\kappa$-compactly generated (in particular presentable). Every object of $\mathcal C$ is an iterated small colimit of $X_i$'s.
	\item Every $\kappa$-compact object of $\mathcal C$ is a retract of an iterated $\kappa$-small colimit of $X_i$'s.
\end{lemenum}
\end{lemma}
\begin{proof}
Let $\mathcal D(0) \subset \mathcal C$ be the full subcategory spanned by the objects $X_i$. For every ordinal $\alpha$ we define the full subcategory $\mathcal D(\alpha) \subset \mathcal C$ by transfinite induction: For limit ordinals we just take the union of the previous steps; for a successor we let $\mathcal D(\alpha + 1)$ be the full subcategory of $\mathcal C$ spanned by all $\kappa$-small colimits of objects in $\mathcal D(\alpha)$. Then each $\mathcal D(\alpha)$ consists of $\kappa$-compact objects of $\mathcal C$ (by transfinite induction) and $\mathcal D := \mathcal D(\kappa)$ is stable under all $\kappa$-small colimits. By \cite[Proposition 5.3.5.11]{lurie-higher-topos-theory} the inclusion $\mathcal D \injto \mathcal C$ induces a $\kappa$-accessible fully faithful functor
 \begin{align*}
 	\Ind_\kappa(\mathcal D) \xinjto{F} \mathcal C.
 \end{align*}
We claim that $F$ is an equivalence, which proves (i) and (ii) (for (ii) use also \cite[Lemma 5.4.2.4]{lurie-higher-topos-theory}). By \cite[Proposition 5.5.1.9, Remark 5.5.2.10]{lurie-higher-topos-theory} the functor $F$ admits a right adjoint $G$. It is enough to show that $G$ is conservative, so suppose we have a map $Y \to Y'$ in $\mathcal C$ such that $G(Y) \isoto G(Y')$ is an isomorphism. For every $X_i \in \mathcal D(0)$ we have
\begin{align*}
	\Hom(X_i, Y) = \Hom(F(X_i), Y) = \Hom(X_i, G(Y))
\end{align*}
(and similarly for $Y'$ in place of $Y$), hence the map
\begin{align*}
	\Hom(X_i, Y) \isoto \Hom(X_i, Y')
\end{align*}
is an isomorphism. Since the family of functors $\Hom(X_i, -)$ is conservative, this implies that $Y \isoto Y'$ is an isomorphism, as desired.
\end{proof}

Our next goal is to show that for uncountable cardinals $\kappa$, $\kappa$-small limits of $\kappa$-compactly generated $\infty$-categories are still $\kappa$-compactly generated, see \cref{rslt:tau-small-limits-in-PrL-tau} below (these limits are still presentable by \cite[Proposition 5.5.3.13]{lurie-higher-topos-theory}, but we need a more refined version of that statement). We start with several technical results about filtered $\infty$-categories. Recall the definition of $\kappa$-cofinal functors from \cite[Definition 5.4.5.8]{lurie-higher-topos-theory}.

\begin{definition}
Let $\tau > \kappa$ be two regular cardinals. Let $\catfiltcofin\tau\kappa \subset \infcatinf$ be the subcategory whose objects are spanned by the $\tau$-filtered $\infty$-categories which admit $\tau$-small $\kappa$-filtered colimits and whose morphisms are spanned by the $\kappa$-cofinal functors which preserve $\tau$-small $\kappa$-filtered colimits.
\end{definition}

\begin{lemma} \label{rslt:catfiltcofin-stable-under-undercategories}
Let $\tau > \kappa$ be two regular cardinals, let $f\colon \mathcal J \to \mathcal I$ be a morphism in $\catfiltcofin\tau\kappa$ and let $p\colon K \to \mathcal J$ be a diagram where $K$ is $\kappa$-small. Then the induced functor $f_{p/}\colon \mathcal J_{p/} \to \mathcal I_{fp/}$ lies in $\catfiltcofin\tau\kappa$.
\end{lemma}
\begin{proof}
By \cite[Lemma 5.4.5.11]{lurie-higher-topos-theory} the $\infty$-categories $\mathcal J_{p/}$ and $\mathcal I_{fp/}$ are $\tau$-filtered. By \cite[Lemma 5.4.6.4]{lurie-higher-topos-theory} the functor $f_{p/}$ is $\kappa$-cofinal.

We now show that $\mathcal J_{p/}$ admits $\tau$-small $\kappa$-filtered colimits (the same argument then applies to $\mathcal I_{fp/}$). This is just a refinement of \cite[Lemma 5.4.5.14]{lurie-higher-topos-theory}: Let $\mathcal K$ be a $\tau$-small $\kappa$-filtered $\infty$-category and let $q_0\colon \mathcal K \to \mathcal J_{p/}$ be a diagram corresponding to a map $q\colon K \star \mathcal K \to \mathcal J$. Observe that $K \star \mathcal K$ is still $\tau$-small and $\kappa$-filtered, so that there is a colimit diagram $\overline q\colon (K \star \mathcal K)^\triangleright \to \mathcal J$. Then note that the diagram $\overline q$ can also be identified with a colimit diagram for $q_0$.

A similar argument as in the previous paragraph shows that $f_{p/}$ preserves $\tau$-small $\kappa$-filtered colimits, which finishes the proof.
\end{proof}

\begin{lemma} \label{rslt:catfiltcofin-stable-under-tau-small-limits}
Let $\tau > \kappa$ be two regular cardinals. The $\infty$-category $\catfiltcofin\tau\kappa$ admits all $\tau$-small limits and the inclusion $\catfiltcofin\tau\kappa \subset \infcatinf$ preserves $\tau$-small limits.
\end{lemma}
\begin{proof}
By the proof of \cite[Proposition 4.4.2.6]{lurie-higher-topos-theory} it is enough to prove the claim for $\tau$-small products and fiber products. For $\tau$-small products this is easy to see: Given a $\tau$-small family $(\mathcal I_i)_{i\in I}$ in $\catfiltcofin\tau\kappa$, their product $\mathcal I := \prod_i \mathcal I_i$ is again in $\catfiltcofin\tau\kappa$ (this follows easily from the fact that colimits in $\mathcal I$ are computed componentwise). One checks similarly that a functor $\mathcal J \to \mathcal I$ of $\infty$-categories is $\kappa$-cofinal and preserves $\tau$-small $\kappa$-filtered colimits if and only if all the projections $\mathcal J \to \mathcal I_i$ do (using \cref{rslt:catfiltcofin-stable-under-undercategories}, reduce the $\kappa$-cofinality proof to a proof of weak cofinality).

It remains to prove the claim for fiber products, so let us assume that we are given a homotopy Cartesian (with respect to the Joyal model structure) diagram
\begin{center}\begin{tikzcd}
	\mathcal J' \arrow[r,"q'"] \arrow[d,"p'"] & \mathcal J \arrow[d,"p"]\\
	\mathcal I' \arrow[r,"q"] & \mathcal I
\end{tikzcd}\end{center}
of $\infty$-categories such that $\mathcal J$, $\mathcal I$, $\mathcal I'$, $p$ and $q$ all lie in $\catfiltcofin\tau\kappa$. By \cite[Lemmas 5.4.5.5 and 5.4.6.5]{lurie-higher-topos-theory}, $\mathcal J'$, $p'$ and $q'$ all lie in $\catfiltcofin\tau\kappa$ as well. To finish the proof, we have to show the following: Given an object $\mathcal K \in \catfiltcofin\tau\kappa$ and a functor $f\colon \mathcal K \to \mathcal J'$ such that $p' \comp f$ and $q' \comp f$ both lie in $\catfiltcofin\tau\kappa$, then $f$ lies in $\catfiltcofin\tau\kappa$. Given such $f$, it follows from \cite[Lemma 5.4.5.5]{lurie-higher-topos-theory} that $f$ preserves $\tau$-small $\kappa$-filtered colimits. It remains to show that $f$ is $\kappa$-cofinal.

Let $K$ be a $\kappa$-small and weakly contractible simplicial set and $g\colon K \to \mathcal K$ a diagram. We have to show that the functor $\mathcal K_{g/} \to \mathcal J'_{fg/}$ induced by $f$ is weakly cofinal. By \cref{rslt:catfiltcofin-stable-under-undercategories} and \cite[Lemma 5.4.5.4]{lurie-higher-topos-theory} the induced diagram
\begin{center}\begin{tikzcd}
	\mathcal J'_{fg/} \arrow[r,"q'"] \arrow[d,"p'"] & \mathcal J_{p'fg/} \arrow[d,"p"]\\
	\mathcal I'_{q'fg/} \arrow[r,"q"] & \mathcal I_{pq'fg/}
\end{tikzcd}\end{center}
of $\infty$-categories is still homotopy cartesian and lies in $\catfiltcofin\tau\kappa$. This allows us to replace $\mathcal K$ by $\mathcal K_{g/}$, $\mathcal J'$ by $\mathcal J'_{fg/}$ etc. in order to reduce to the case $K = \emptyset$. In other words, we want to show that $f$ is weakly cofinal.

Let an object $j' \in \mathcal J'$ be given. We have to find some $k \in \mathcal K$ and a morphism $j' \to f(k)$ in $\mathcal J'$. Since the diagram $\mathcal J' \to \mathcal I' \cprod \mathcal J \from \mathcal K$ lies in $\catfiltcofin\tau\kappa$, it satisfies the hypothesis of \cite[Lemma 5.4.6.3]{lurie-higher-topos-theory}. The proof of loc. cit. then constructs a commutative diagram
\begin{center}\begin{tikzcd}
	A \arrow[r] \arrow[d, "\alpha"] & (\kappa) \arrow[d, "\gamma"] & A' \arrow[l] \arrow[d, "\beta"]\\
	\mathcal J' \arrow[r] & \mathcal I' \cprod \mathcal J & \mathcal K \arrow[l]
\end{tikzcd}\end{center}
where $(\kappa)$ is the set of ordinals less than $\kappa$, which splits as a disjoint union $(\kappa) = A \dunion A'$ of even and odd ordinals, respectively. The inductive nature of the proof of loc. cit. allows us to choose $\alpha(0) = j'$. Now let $\overline j'$ be a colimit of $\alpha$ and $k$ a colimit of $\beta$. Then the images of $\overline j'$ and $k$ in $\mathcal I' \cprod \mathcal J$ are equivalent (both are colimits of $\gamma$). We can assume that $p$ and $q$ are categorical fibrations and hence that $\mathcal J' = \mathcal I' \cprod_{\mathcal I} \mathcal J$. It follows that $f(k) \in \mathcal J'$ is equivalent to $\overline j'$. On the other hand there is clearly a map $j' \to \overline j'$. We therefore obtain a map $j' \to f(k)$, as desired.
\end{proof}

\begin{lemma} \label{rslt:projlim-of-cat-has-colimits}
Let $\mathcal C = \varprojlim_i \mathcal C_i$ be a limit of $\infty$-categories and let $K$ be a simplicial set. Assume that all $\mathcal C_i$ admit colimits indexed by $K$ and that all transition functors in the diagram $(\mathcal C_i)_i$ preserve colimits indexed by $K$. Then:
\begin{lemenum}
	\item $\mathcal C$ admits colimits indexed by $K$.
	\item A diagram $\overline p\colon K^\triangleright \to \mathcal C$ is a colimit diagram if and only if all the induced diagrams $K^\triangleright \to \mathcal C_i$ are colimit diagrams. In particular, all the functors $\mathcal C \to \mathcal C_i$ preserve colimits indexed by $K$.
\end{lemenum}
\end{lemma}
\begin{proof}
By the proof of \cite[Proposition 4.4.2.6]{lurie-higher-topos-theory} it is enough to prove the claim in the case that the limit $\varprojlim_i \mathcal C_i$ is a small product or a fiber product. The former case is clear because colimits in a product are computed pointwise. The latter case is \cite[Lemma 5.4.5.5]{lurie-higher-topos-theory}.
\end{proof}

For the following, recall the definition of $\catAcc\kappa$ in \cite[Definition 5.4.2.16]{lurie-higher-topos-theory} for any regular cardinal $\kappa$: This category is the category of $\kappa$-accessible $\infty$-categories whose functors preserve $\kappa$-filtered colimits and $\kappa$-compact objects. Also recall the definition of $\tau \gg \kappa$ for two regular cardinals $\kappa$ and $\tau$ (see \cite[Definition 5.4.2.8]{lurie-higher-topos-theory}).

\begin{proposition} \label{rslt:limit-of-Acc-kappa-is-in-Acc-tau}
Let $\kappa$ and $\tau$ be two regular cardinals such that $\tau > \kappa$ and $\tau \gg \kappa$. Let $(\mathcal C_i)_{i\in I} \in \catAcc\tau$ be a $\tau$-small diagram in $\catAcc\tau$. Assume that each $\mathcal C_i$ admits all small $\kappa$-filtered colimits and all transition functors in $(\mathcal C_i)_i$ preserves these colimits. Let $\mathcal C = \varprojlim_i \mathcal C_i$ be the limit in $\infcatinf$. Then $\mathcal C$ is also the limit of $(\mathcal C_i)_i$ in $\catAcc\tau$, i.e.
\begin{propenum}
	\item $\mathcal C$ lies in $\catAcc\tau$,
	\item given any $\mathcal D \in \catAcc\tau$ then a functor $\mathcal D \to \mathcal C$ lies in $\catAcc\tau$ if and only if all induced functors $\mathcal D \to \mathcal C_i$ lie in $\catAcc\tau$.
\end{propenum}
Moreover, we have a natural equivalence $\mathcal C^\tau = \varprojlim_i \mathcal C_i^\tau$.
\end{proposition}
\begin{proof}
We generalize the argument of \cite[Proposition 5.4.6.6]{lurie-higher-topos-theory}. By \cref{rslt:projlim-of-cat-has-colimits}, $\mathcal C$ admits small $\tau$-filtered colimits. Thus in order to show that $\mathcal C \in \catAcc\tau$ we only have to show that there is an essentially small full subcategory $\mathcal C' \subset \mathcal C$ such that $\mathcal C'$ consists of $\tau$-compact objects and every object in $\mathcal C$ can be written as a $\tau$-filtered colimit in $\mathcal C'$ (cf. \cite[Proposition 5.4.2.2.(3)]{lurie-higher-topos-theory}). We let $\mathcal C' := \varprojlim_i \mathcal C_i^\tau$ (the limit being taken in $\infcatinf$). Then we have to show the following:
\begin{itemize}
	\item The natural functor $\mathcal C' \to \mathcal C$ is fully faithful (so we will write $\mathcal C' \subset \mathcal C$) and $\mathcal C'$ consists of $\tau$-compact objects in $\mathcal C$. To see this, by the proof of \cite[Proposition 4.4.2.6]{lurie-higher-topos-theory} it is enough to check that both claims are stable under fiber products and $\tau$-small products. Full faithfulness is clear. Regarding $\tau$-compactness, for a $\tau$-small product this follows from \cite[Proposition 5.3.4.10]{lurie-higher-topos-theory} and for fiber products it follows from \cite[Lemma 5.4.5.7]{lurie-higher-topos-theory}.

	\item Every object in $\mathcal C$ is the colimit of a $\tau$-filtered diagram in $\mathcal C'$. To see this, let $X \in \mathcal C$ be given and for each $i \in I$ let $X_i \in \mathcal C_i$ be its image in $\mathcal C_i$. Let $\mathcal C'_{/X}$ be the full subcategory of $\mathcal C_{/X}$ consisting of those morphisms $Y \to X$ such that $Y \in \mathcal C'$. We similarly define $\mathcal C_{i,/X_i}^\tau$ for each $i \in I$.

	Note that $\mathcal C'_{/X}$ is a limit of $(\mathcal C_{i,/X_i}^\tau)_i$ in $\infcatinf$: This can again be checked separately for ($\tau$-)small products and fiber products; for ($\tau$-)small products this is clear and for fiber products this follows in the same way as in \cite[Lemma 5.4.5.4]{lurie-higher-topos-theory} (note that the argument of \cite[Lemma 5.4.5.3]{lurie-higher-topos-theory} also applies to overcategories, cf. \cite[Corollary 4.2.1.3]{lurie-higher-topos-theory}).

	From $\mathcal C_i \in \catAcc\tau$ it follows that $\mathcal C_{i,/X_i}^\tau$ is $\tau$-filtered. Moreover, by \cite[Corollary 5.3.4.15]{lurie-higher-topos-theory} the categories $\mathcal C_i^\tau$ are stable under $\tau$-small colimits in $\mathcal C_i$ so that from the assumptions on $\mathcal C_i$ it follows that $\mathcal C_i^\tau$ is admits all $\tau$-small $\kappa$-filtered colimits; the same is then true for $\mathcal C_{i,/X_i}^\tau$. All in all we see $\mathcal C_{i,/X_i}^\tau \in \catfiltcofin\tau\kappa$ for all $i$. By our assumption on the diagram $(\mathcal C_i)_i$, all transition functors in the diagram $(\mathcal C_{i,/X_i}^\tau)_i$ preserve $\tau$-small $\kappa$-filtered colimits. Moreover, from \cite[Lemma 5.4.6.1]{lurie-higher-topos-theory} it follows that all transition functors in the diagram $(\mathcal C_{i,/X_i}^\tau)_i$ are $\kappa$-cofinal. Altogether this implies that the diagram $(\mathcal C_{i,/X_i}^\tau)_i$ factors through $\catfiltcofin\tau\kappa$. From \cref{rslt:catfiltcofin-stable-under-tau-small-limits} we deduce that $\mathcal C'_{/X} \in \catfiltcofin\tau\kappa$ and that all the functors $\mathcal C'_{/X} \to \mathcal C_{i,/X_i}^\tau$ are $\kappa$-cofinal.

	There is an obvious map of diagrams
	\begin{center}\begin{tikzcd}
		(\mathcal C'_{/X})^\triangleright \arrow[r,"h'"] \arrow[d,"(p_i)_i"] & \mathcal C \arrow[d,"(q_i)_i"]\\
		((\mathcal C_{i,/X_i}^\tau)^\triangleright)_i \arrow[r,"(h_i)_i"] & (\mathcal C_i)_i
	\end{tikzcd}\end{center}
	Since each $\mathcal C_i$ lies in $\catAcc\tau$, all the functors $h_i$ are colimit diagrams. From the above discussion we obtain that the functors $p_i$ are $\kappa$-cofinal, so by \cite[Lemma 5.4.5.12]{lurie-higher-topos-theory} they are in particular cofinal. It follows that the functors $h_i \comp p_i = q_i \comp h'$ are colimit diagrams. Since $\mathcal C$ is a limit of $(\mathcal C_i)_i$, it follows using \cref{rslt:projlim-of-cat-has-colimits} that $h'$ is a colimit diagram. Thus $X$ is a $\tau$-filtered colimit of objects in $\mathcal C_i$, as desired.
\end{itemize}
We have finished the proof of (i). Part (ii) is now easy to see: Let $\mathcal D \in \catAcc\tau$, let $f\colon \mathcal D \to \mathcal C$ be a functor and let $f_i\colon \mathcal D \to \mathcal C_i$ be the induced functors; the only nontrivial parts of (ii) are the following:
\begin{itemize}
	\item If all $f_i$ preserve $\tau$-filtered colimits then so does $f$. As before, this problem can be decomposed into the cases of ($\tau$-)small products and fiber products of $\infty$-categories. For the latter case use \cite[Lemma 5.4.5.5]{lurie-higher-topos-theory}.

	\item If all $f_i$ preserve $\tau$-compact objects then so does $f$. Again the only interesting case is that of a fiber product, where we apply \cite[Lemma 5.4.5.7]{lurie-higher-topos-theory}.
\end{itemize}
It remains to show the last part of the claim, i.e. $\mathcal C' = \mathcal C^\tau$. By what we have shown above we know that $\mathcal C = \Ind_\tau(\mathcal C')$. By \cite[Lemma 5.4.2.4]{lurie-higher-topos-theory} all $\mathcal C_i^\tau$ are idempotent complete, hence by \cref{rslt:projlim-of-cat-has-colimits} the same is true for $\mathcal C'$. Thus, by again applying \cite[Lemma 5.4.2.4]{lurie-higher-topos-theory} we deduce $\mathcal C' = \mathcal C^\tau$.
\end{proof}

We now arrive at the promised result about $\kappa$-small limits of $\kappa$-compactly generated $\infty$-categories. Recall the definition of $\catPrLk\kappa$ in \cite[Notation 5.5.7.7]{lurie-higher-topos-theory} for any regular cardinal $\kappa$: It is the $\infty$-category of $\kappa$-compactly generated $\infty$-categories whose functors preserve all small colimits and $\kappa$-compact objects.

\begin{corollary} \label{rslt:tau-small-limits-in-PrL-tau}
Let $\tau$ be an uncountable regular cardinal. Then $\catPrLk\tau$ admits $\tau$-small limits and the forgetful functor $\catPrLk\tau \to \infcatinf$ preserves these limits. Moreover, for every $\tau$-small diagram $(\mathcal C_i)_i$ in $\catPrLk\tau$ we have $\mathcal C^\tau = \varprojlim_i \mathcal C_i^\tau$.
\end{corollary}
\begin{proof}
Let $(\mathcal C_i)_i$ be a $\tau$-small diagram in $\catPrLk\tau$ with limit $\mathcal C = \varprojlim_i \mathcal C_i$ in $\infcatinf$. Since $\tau$ is uncountable we have $\tau \gg \omega$ and $\tau > \omega$. Thus we can apply \cref{rslt:limit-of-Acc-kappa-is-in-Acc-tau} with $\kappa = \omega$ to see that $\mathcal C$ is $\tau$-accessible. By \cref{rslt:projlim-of-cat-has-colimits} it has all small colimits and hence is $\tau$-compactly generated. We also obtain the identity $\mathcal C^\tau = \varprojlim_i \mathcal C_i^\tau$. To finish the proof it remains to show that for any $\tau$-compactly generated $\infty$-category $\mathcal D$ and any functor $f\colon \mathcal D \to \mathcal C$, $f$ lies in $\catPrLk\tau$ if and only all induced functors $f_i\colon \mathcal D \to \mathcal C_i$ lie in $\catPrLk\tau$. By \cref{rslt:limit-of-Acc-kappa-is-in-Acc-tau} $f$ preserves $\tau$-compact objects if and only if all $f_i$ do. Moreover, $f$ preserves small colimits if and only if all $f_i$ do: This follows as in the proof of \cref{rslt:limit-of-Acc-kappa-is-in-Acc-tau} by using \cite[Lemma 5.4.5.5]{lurie-higher-topos-theory}.
\end{proof}

We now turn our focus to limits and colimits inside a presentable $\infty$-category. It is a classical fact that filtered colimits preserve finite limits in the category of sets. This generalizes to a similar statement in the $\infty$-category of anima by \cite[Proposition 5.3.3.3]{lurie-higher-topos-theory}. In the following we will prove that the same property also holds in the $\infty$-category of (small) $\infty$-categories. First we note that this property holds in any compactly generated $\infty$-category:

\begin{lemma} \label{rslt:filtered-colim-preserve-small-lim-in-compactly-generated-cat}
Let $\kappa$ be a regular cardinal and let $\mathcal C$ be a $\kappa$-compactly generated $\infty$-category. Then for every regular cardinal $\tau \gg \kappa$, $\tau$-filtered colimits commute with $\tau$-small limits in $\mathcal C$.

More precisely, for every $\tau$-filtered $\infty$-category $\mathcal I$ the colimit functor $\varinjlim_{\mathcal I}\colon \Fun(\mathcal I, \mathcal C) \to \mathcal C$ constructed in \cite[\S5.3.3]{lurie-higher-topos-theory} preserves $\tau$-small limits.
\end{lemma}
\begin{proof}
By \cite[Proposition 5.4.2.11]{lurie-higher-topos-theory} $\mathcal C$ is $\tau$-compactly generated, so we can w.l.o.g. assume $\tau = \kappa$. By \cite[Theorem 5.5.1.1.(4)]{lurie-higher-topos-theory} there is a small $\infty$-category $\mathcal D$ which admits $\kappa$-small colimits and an equivalence $\mathcal C \isom \Ind_\kappa(\mathcal D)$. Then by \cite[Corollary 5.3.5.4]{lurie-higher-topos-theory} the category $\Ind_\kappa(\mathcal D) \subset \mathcal P(\mathcal D) := \Fun(\mathcal D^\opp, \Ani)$ can be described as the full subcategory spanned by the functors which preserve $\kappa$-small limits. It follows immediately that $\Ind_\kappa(\mathcal D) \subset \mathcal P(\mathcal D)$ is stable under ($\kappa$-)small limits (as those are computed point-wise, see \cite[Corollary 5.1.2.3]{lurie-higher-topos-theory}) and by \cite[Proposition 5.3.5.3]{lurie-higher-topos-theory} it is also stable under $\kappa$-filtered colimits. It is thus enough to show that $\kappa$-filtered colimits commute with $\kappa$-small limits in $\mathcal P(\mathcal D)$. But since both colimits and limits are computed pointwise, this reduces immediately to showing that $\kappa$-filtered colimits commute with $\kappa$-small limits in $\Ani$. This is \cite[Proposition 5.3.3.3]{lurie-higher-topos-theory}.
\end{proof}

In order to get the desired commutation of filtered colimits and finite limits in $\infcatinf$, it now remains to show that $\infcatinf$ is compactly generated. To prove this, we will make use of the theory of localizations. Recall that for an $\infty$-category $\mathcal C$ and a subset $S$ of morphisms in $\mathcal C$, an object $X \in \mathcal C$ is called \emph{$S$-local} if for every morphism $s\colon Y \to Z$, composition with $s$ induces an isomorphism $\Hom_{\mathcal C}(Z, X) \isoto \Hom_{\mathcal C}(Y, X)$ of anima (cf. \cite[Definition 5.5.4.1]{lurie-higher-topos-theory}).

In the following, denote by $S^{-1}\mathcal C \subset \mathcal C$ the full subcategory spanned by $S$-local objects (be aware that in general this is not a localization of $\infty$-categories; however this is true if $S$ is spanned by a small subset and $\mathcal C$ is presentable, see \cite[Proposition 5.5.4.15]{lurie-higher-topos-theory}).

\begin{lemma} \label{rslt:localization-compactly-generated-criterion}
Let $\kappa$ be a regular cardinal, let $\mathcal C$ be a $\kappa$-compactly generated $\infty$-category and let $S$ be a small collection of morphisms between $\kappa$-compact objects of $\mathcal C$. Then $S^{-1}\mathcal C$ is $\kappa$-compactly generated.
\end{lemma}
\begin{proof}
By \cite[Proposition 5.5.4.15.(3)]{lurie-higher-topos-theory} the inclusion $S^{-1}\mathcal C \subset \mathcal C$ has a left adjoint $L\colon \mathcal C \to S^{-1}\mathcal C$, which is then automatically a localization functor. Thus by \cite[Corollary 5.5.7.3]{lurie-higher-topos-theory} it is enough to show that $S^{-1}\mathcal C \subset \mathcal C$ is stable under $\kappa$-filtered colimits. This follows as in the proof of \cite[Proposition 5.5.4.2.(3)]{lurie-higher-topos-theory}: Let $\overline p\colon K^\triangleright \to \mathcal C$ be a colimit diagram such that $K$ is small and $\kappa$-cofiltered and $p := \restrict pK$ factors through $S^{-1}\mathcal C$. Let $s\colon X \to Y$ be a morphism in $S$ and let $s'\colon F_X \to F_Y$ be the corresponding map of corepresentable functors $\mathcal C \to \Ani$ (where we let $\Ani$ denote the $\infty$-category of \emph{large} anima in this proof). By definition of $\kappa$-compactness, since both $X$ and $Y$ are $\kappa$-compact, the maps $\overline p_X = F_X \comp \overline p$ and $\overline p_Y = F_Y \comp \overline p$ are colimit diagrams in $\Ani$. The map $s'$ induces a map $\overline p_X \to \overline p_Y$ which is an equivalence when restricted to $K$ (since $\overline p(k)$ is $S$-local for $k \in K$) and hence must be an equivalence in general. Evaluation at $\infty \in K^\triangleright$ shows that $\overline p(\infty)$ is $\{ s \}$-local. Repeating this for all $s \in S$ proves $\overline p(\infty) \in S^{-1}\mathcal C$, as desired.
\end{proof}

\begin{proposition} \label{rslt:infcatinf-is-compactly-generated}
The $\infty$-category $\infcatinf$ of (small) $\infty$-categories is compactly generated.
\end{proposition}
\begin{proof}
Let $\mathbf A := \catset_\Delta^+$ be the simplicial model category of marked simplicial sets with the Cartesian model structure (cf. \cite[Corollar 3.1.4.4]{lurie-higher-topos-theory}). By \cite[Theorem 3.1.5.1.(A0)]{lurie-higher-topos-theory} this model category is equivalent to the category of simplicial sets with the Joyal model structure, so that $\infcatinf$ is equivalent to the simplicial nerve of the fibrant-cofibrant objects in $\mathbf A$, i.e. $\infcatinf \isom N(\mathbf A^\circ)$. Let $\mathbf B$ be the simplicial model category of complete Segal spaces, which has as underlying category the category $\Fun(\Delta^\opp, \catset_\Delta)$ of simplicial spaces; see \cite[Theorem 7.2]{rezk-segal-spaces}. By \cite[Theorem 4.12]{segal-spaces-vs-infty-categories} there is a Quillen equivalence of $\mathbf B$ with the Joyal model category. All in all we obtain a Quillen equivalence $\mathbf A \isom \mathbf B$ of model categories and one can check by explicit computation that the right Quillen adjoint $\mathbf A \to \mathbf B$ is a simplicial functor. We can thus apply \cite[Corollary A.3.1.12]{lurie-higher-topos-theory} to see that $N(\mathbf B^\circ) \isom N(\mathbf A^\circ) \isom \infcatinf$.

By the beginning of \cite[\S12]{rezk-segal-spaces} the model structure $\mathbf B$ is defined as a Bousfield localization of the Reedy model structure on $\Fun(\Delta^\opp, \catset_\Delta)$. By \cite[Example A.2.9.21]{lurie-higher-topos-theory} the Reedy model structure on $\Fun(\Delta^\opp, \catset_\Delta)$ coincides with the injective model structure. Thus the proof of \cite[Proposition A.3.7.6]{lurie-higher-topos-theory} shows that $\infcatinf$ is a localization of the category $\mathcal P(\Delta) = \Fun(\Delta^\opp, \Ani)$. Looking at the definition of $\mathbf B$, we see more precisely that $\infcatinf = S^{-1}\mathcal P(\Delta)$, where $S$ is the following collection of morphisms: Let $\catset \subset \Ani$ be the inclusion of (small) discrete spaces and let $c\colon \catset_\Delta \to \mathcal P(\Delta)$ be the map $\catset_\Delta = \Fun(\Delta^\opp, \catset) \subset \Fun(\Delta^\opp, \Ani) = \mathcal P(\Delta)$. Then $S$ is given as the collection of maps
\begin{align*}
	S = \{ c(\mathrm{Sp}^n) \injto c(\Delta^n) \}_{n\ge0} \union \{ c(N(\{ a \isotofrom b \})) \to c(*) \},
\end{align*}
where the \emph{spine} $\mathrm{Sp}^n \subset \Delta^n$ is defined to be the union of the edges $0 \to 1, 1 \to 2, \dots, n-1 \to n$ and $\{ a \isotofrom b \}$ denotes the groupoid with two objects and precisely one isomorphism between them. By the proof of \cite[Theorem 5.5.1.1.(5)]{lurie-higher-topos-theory} the category $\mathcal P(\Delta)$ is compactly generated. In order to finish the proof it is therefore enough to show that all objects occuring as target and source of morphisms in $S$ are compact; then we can apply \cref{rslt:localization-compactly-generated-criterion}.

It is now enough to prove the following claim: Given any finite simplicial set $K$ the object $c(K) \in \mathcal P(\Delta)$ is compact. To see this, choose $n \ge 0$ large enough so that $K$ has no non-degenerate simplices in dimension $>n$. Note that the restriction functor $i^*\colon \mathcal P(\Delta) \to \mathcal P(\Delta_{\le n})$ has a left adjoint $i_!\colon \mathcal P(\Delta_{\le n}) \to \mathcal P(\Delta)$ given by left Kan extension (cf. \cite[Proposition 4.3.2.17]{lurie-higher-topos-theory}). The composition $i_! i^*\colon \mathcal P(\Delta) \to \mathcal P(\Delta)$ is the $n$-skeleton functor. By the choice of $n$ we get $i_! i^* c(K) = c(K)$. On the other hand, $i^* c(K) \in \mathcal P(\Delta_{\le n})$ is compact by \cite[Proposition 5.3.4.13]{lurie-higher-topos-theory} and $i_!$ preserves compact objects because its right adjoint $i^*$ preserves (filtered) colimits (cf. \cite[Proposition 5.5.7.2.(1)]{lurie-higher-topos-theory}). Thus $c(K) \in \mathcal P(\Delta)$ is compact, as desired.
\end{proof}

\begin{corollary} \label{rslt:filtered-colim-preserve-finite-lim-in-infcatinf}
Let $\kappa$ be a regular cardinal. Then $\kappa$-filtered colimits preserve $\kappa$-small limits in $\infcatinf$.
\end{corollary}
\begin{proof}
Combine \cref{rslt:filtered-colim-preserve-small-lim-in-compactly-generated-cat} and \cref{rslt:infcatinf-is-compactly-generated} and note that $\kappa \gg \omega$ for every regular cardinal $\kappa$.
\end{proof}

\subsection{Sheaves of \texorpdfstring{$\infty$}{Infinity}-Categories} \label{sec:infcat.sheaves}

Throughout this thesis we work with sheaves valued in various $\infty$-categories, in particular sheaves of anima and sheaves of $\infty$-categories. Most of the relevant theory is available in the literature, but we collect it here for the convenience of the reader. Unlike in \cite{lurie-higher-topos-theory} we keep notation with respect to the $\infty$-site $\mathcal C$ and do not pass to the associated $\infty$-topos $\mathcal X = \catshs{\mathcal C}$, as this seems the most intuitive way to work with in the thesis.

We start be recalling the definition of a site. The modern definition of sites (or Grothendieck topologies) is based on the notion of covering sieves rather than explicit covering families, so we stick to that convention. Later we formulate a more explicit way of constructing sites, more in line with how sites are usually handled in algebraic geometry.

\begin{definition}
Let $\mathcal C$ be an $\infty$-category and $X \in \mathcal C$ an object.
\begin{defenum}
	\item A \emph{sieve} on $X$ is a full subcategory $\mathcal C^0_{/X} \subset \mathcal C_{/X}$ such that for every morphism $Z \to Y$ in $\mathcal C_{/X}$, if $Y$ belongs to $\mathcal C^0_{/X}$ then so does $Z$.

	\item Given a sieve $\mathcal C^0_{/X}$ on $X \in \mathcal C$ and a morphism $f\colon Y \to X$ in $\mathcal C$, we define the \emph{pullback} $f^* \mathcal C^0_{/X}$ to be the sieve on $Y$ spanned by those morphisms $\varphi\colon Z \to Y$ such that $f \comp \varphi\colon Z \to X$ lies in $\mathcal C^0_{/X}$.
\end{defenum}
\end{definition}

\begin{definition}
An \emph{$\infty$-site} is a small $\infty$-category $\mathcal C$ equipped with a Grothendieck topology, i.e. a collection of \emph{covering sieves} on $\mathcal C$ satisfying the following properties:
\begin{enumerate}[(i)]
	\item For every object $X \in \mathcal C$ the overcategory $\mathcal C_{/X}$ is a covering sieve on $X$.
	\item Let $f\colon Y \to X$ be a morphism in $\mathcal C$ and $\mathcal C^0_{/X}$ a covering sieve on $X$. Then $f^* \mathcal C^0_{/X}$ is a covering sieve on $Y$.
	\item Let $X \in \mathcal C$ be an object and let $\mathcal C^0_{/X}, \mathcal C^1_{/X} \subset \mathcal C_{/X}$ be two sieves on $X$. Suppose that $\mathcal C^0_{/X}$ is a covering sieve and that for every object $\varphi\colon Y \to X$ of $\mathcal C^0_{/X}$ the sieve $\varphi^* \mathcal C^1_{/X}$ is a covering sieve on $Y$. Then $\mathcal C^1_{/X}$ is a covering sieve on $X$.
\end{enumerate}
\end{definition}

With the definition of $\infty$-sites at hand, we can now define the associated $\infty$-category of sheaves and hypercomplete sheaves. For the following definition cf. \cite[\S1.3.1]{lurie-spectral-algebraic-geometry}.

\begin{definition}
Let $\mathcal C$ be an $\infty$-site and let $\mathcal D$ be any $\infty$-category which admits small limits. A \emph{$\mathcal D$-valued sheaf on $\mathcal C$} is a functor $\mathcal F \colon \mathcal C^\opp \to \mathcal D$ such that for every object $X \in \mathcal C$ and every covering sieve $\mathcal C^0_{/X}$ of $X$ the composite functor
\begin{align*}
	(\mathcal C^0_{/X})^\vartriangleright \subset (\mathcal C_{/X})^\vartriangleright \to \mathcal C \xto{\mathcal F^\opp} \mathcal D^\opp
\end{align*}
is a colimit diagram. We denote by
\begin{align*}
	\catsh{\mathcal C}{\mathcal D}
\end{align*}
the full subcategory of $\Fun(\mathcal C^\opp, \mathcal D)$ spanned by the $\mathcal D$-valued sheaves on $\mathcal C$.

In the case that $\mathcal D = \Ani$ is the $\infty$-category of anima we abbreviate $\catshs{\mathcal C} = \catsh{\mathcal C}{\Ani}$. We call $\mathcal X = \catshs{\mathcal C}$ the \emph{$\infty$-topos associated to $\mathcal C$}.
\end{definition}

Given an $\infty$-topos $\mathcal X$ and an $\infty$-category $\mathcal D$ with all small limits, there is a definition $\mathcal D$-valued sheaves on $\mathcal X$: Namely it is the $\infty$-category of limit preserving functors $\mathcal F\colon \mathcal X^\opp \to \mathcal D$ (cf. \cite[Notation 6.3.5.16]{lurie-higher-topos-theory}). A priori there might be a possible ambiguity in notation if $\mathcal X$ is the $\infty$-topos associated to an $\infty$-site $\mathcal C$, which is solved by the following:

\begin{lemma} \label{rslt:equivalence-of-sheaf-defs-on-site}
Let $\mathcal C$ be an $\infty$-site with associated $\infty$-topos $\mathcal X = \catshs{\mathcal C}$. Then for any $\infty$-category $\mathcal D$ with all small limits, there is a natural equivalence
\begin{align*}
	\catsh{\mathcal X}{\mathcal D} \isoto \catsh{\mathcal C}{\mathcal D}.
\end{align*}
Moreover, if $\mathcal D$ is presentable then $\catsh{\mathcal C}{\mathcal D}$ is presentable.
\end{lemma}
\begin{proof}
The first part of the claim is {\cite[Proposition 1.3.1.7]{lurie-spectral-algebraic-geometry}}. Now assume that $\mathcal D$ is presentable. We know that $\mathcal X$ is presentable as well, hence of the form $\mathcal X = \Ind_\kappa(\mathcal X_0)$ for some regular cardinal $\kappa$ and some small $\infty$-category $\mathcal X_0$ which admits $\kappa$-small colimits. Then
\begin{align*}
	\catsh{\mathcal C}{\mathcal D} = \catsh{\mathcal X}{\mathcal D} = \Fun'(\mathcal X_0^\opp, \mathcal D),
\end{align*}
where $\Fun'(\mathcal X_0^\opp, \mathcal D)$ denotes the full subcategory of $\Fun(\mathcal X_0, \mathcal D)$ spanned by those functors which preserve $\kappa$-small limits (see \cite[Proposition 5.5.1.9]{lurie-higher-topos-theory}). Now $\Fun(\mathcal X_0^\opp, \mathcal D)$ is presentable, so we need to show that $\Fun'(\mathcal X_0^\opp, \mathcal D)$ is a strongly reflective subcategory. For each limit diagram $f\colon K^\triangleleft \to \mathcal X_0^\opp$ let $\mathcal E(f) \subset \Fun(\mathcal X_0^\opp, \mathcal D)$ denote the full subcategory spanned by those functors $F\colon \mathcal X_0^\opp \to \mathcal D$ such that $F\comp f\colon K^\triangleleft \to \mathcal D$ is a limit. Then $\Fun'(\mathcal X_0^\opp, \mathcal D) = \bigisect_f \mathcal E(f)$, so by \cite[Lemma 5.5.4.18]{lurie-higher-topos-theory} it is enough to show that each $\mathcal E(f) \subset \Fun(\mathcal X_0^\opp, \mathcal D)$ is a strongly reflective subcategory. From now on we fix the limit diagram $f\colon K^\triangleleft \to \mathcal C$. Let $\mathcal A := \Fun(K^\triangleleft, \mathcal D)$ and let $\mathcal A_0 \subset \mathcal A$ denote the full subcategory of all limit diagrams. By \cite[Lemma 5.5.4.19]{lurie-higher-topos-theory} $\mathcal A_0$ is a strongly reflective subcategory of $\mathcal A$. Now $f$ induces a functor $f^*\colon \Fun(\mathcal X_0^\opp, \mathcal D) \to \mathcal A$ via composition with $f$ and $\mathcal E(f)$ is precisely the category of those functors $F\colon \mathcal X_0^\opp, \mathcal D$ such that $f^*F \in \mathcal A_0$. By \cite[Lemma 5.5.4.17]{lurie-higher-topos-theory} $\mathcal E(f)$ is strongly reflective, as desired.
\end{proof}

In the setting of $\infty$-categories, there is a stronger version of the sheaf property, leading to hypercomplete sheaves. For the following definition, recall the notion of hypercoverings in an $\infty$-topos in \cite[Definition 6.5.3.2]{lurie-higher-topos-theory}.

\begin{definition}
Let $\mathcal C$ be an $\infty$-site with associated $\infty$-topos $\mathcal X = \catshs{\mathcal C}$ and let $\mathcal D$ be an $\infty$-category which admits all small limits. A $\mathcal D$-valued sheaf $\mathcal F$ on $\mathcal C$ is called \emph{hypercomplete} if, viewed as a functor $\mathcal F\colon \mathcal X^\opp \to \mathcal D$, $\mathcal F$ satisfies the following property: For every object $X \in \mathcal X$ and every hypercovering $U_\bullet\to X$ in $\mathcal X_{/X}$, the associated morphism
\begin{align*}
	\mathcal F(X) \isoto \varprojlim_{n\in\Delta} \mathcal F(U_n)
\end{align*}
is an isomorphism. We denote by
\begin{align*}
	\cathsh{\mathcal C}{\mathcal D} \subset \catsh{\mathcal C}{\mathcal D}
\end{align*}
the full subcategory consisting of the hypercomplete sheaves.
\end{definition}

An analogous version of \cref{rslt:equivalence-of-sheaf-defs-on-site} holds also for hypercomplete sheaves, making it unambiguous:

\begin{lemma} \label{rslt:equivalence-of-hypercomplete-sheaf-defs}
Let $\mathcal C$ be an $\infty$-site with associated $\infty$-topos $\mathcal X = \catshs{\mathcal C}$ and its hypercompletion $\hat{\mathcal X}$ (see \cite[\S6.5.3]{lurie-higher-topos-theory}). Then for any $\infty$-category $\mathcal D$ with all small limits, there is a natural equivalence
\begin{align*}
	\catsh{\hat{\mathcal X}}{\mathcal D} \isoto \cathsh{\mathcal C}{\mathcal D}.
\end{align*}
Moreover, if $\mathcal D$ is presentable then $\cathsh{\mathcal C}{\mathcal D}$ is presentable.
\end{lemma}
\begin{proof}
By \cite[Corollary 6.5.3.13]{lurie-higher-topos-theory}, $\hat{\mathcal X}$ is the localization of $\mathcal X$ at the collection of morphisms $X \to \abs{U_\bullet}$ for all hypercoverings $U_\bullet \to X$ in $\mathcal X$. Thus the claimed equivalence follows from \cref{rslt:equivalence-of-sheaf-defs-on-site} and \cite[Proposition 5.5.4.20]{lurie-higher-topos-theory}.

If $\mathcal D$ is presentable then the presentability of $\cathsh{\mathcal C}{\mathcal D} = \catsh{\hat{\mathcal X}}{\mathcal D}$ follows in the same way as in the proof of \cref{rslt:equivalence-of-sheaf-defs-on-site} (the proof only used that $\mathcal X$ is presentable).
\end{proof}

One of the fundamental operations when working with sheaves is sheafification. Of course this also works in the $\infty$-categorical setting, as follows.

\begin{lemma} \label{rslt:sheafification}
Let $\mathcal C$ be an $\infty$-site and let $\mathcal D$ be a presentable $\infty$-category.
\begin{lemenum}
	\item The inclusion $\catsh{\mathcal C}{\mathcal D} \subset \Fun(\mathcal C^\opp, \mathcal D)$ admits a left adjoint
	\begin{align*}
		\Fun(\mathcal C^\opp, \mathcal D) \to \catsh{\mathcal C}{\mathcal D}, \qquad \mathcal F \mapsto \mathcal F^\sharp.
	\end{align*}

	\item The inclusion $\cathsh{\mathcal C}{\mathcal D} \subset \catsh{\mathcal C}{\mathcal D}$ admits a left adjoint
	\begin{align*}
		\catsh{\mathcal C}{\mathcal D} \to \cathsh{\mathcal C}{\mathcal D}, \qquad \mathcal F \mapsto \hat{\mathcal F}.
	\end{align*}
\end{lemenum}
\end{lemma}
\begin{proof}
To prove (i) we need to show that $\catsh{\mathcal C}{\mathcal D} \subset \Fun(\mathcal C^\opp, \mathcal D)$ is a strongly reflective subcategory. This can be shown in the same way as in the proof of \cref{rslt:equivalence-of-sheaf-defs-on-site}: $\catsh{\mathcal C}{\mathcal D}$ is defined to be the full subcategory consisting of those functors $\mathcal F\colon \mathcal C^\opp \to \mathcal D$ which transform a certain set of diagrams $f\colon K^\triangleleft \to \mathcal C^\opp$ into limit diagrams in $\mathcal D$, so the claim follows from \cite[Propositions 5.5.4.17, 5.5.4.18, 5.5.4.19]{lurie-higher-topos-theory}.

For part (ii), note that after the identifications of \cref{rslt:equivalence-of-hypercomplete-sheaf-defs} the inclusion $\cathsh{\mathcal C}{\mathcal D} \subset \catsh{\mathcal C}{\mathcal D}$ is induced by the right Kan extension $\Fun(\hat{\mathcal X}^\opp, \mathcal D) \to \Fun(\mathcal X^\opp, \mathcal D)$ along $\hat{\mathcal X} \subset \mathcal X$, whose left adjoint is given by the restriction functor (cf. \cite[Proposition 4.3.2.17]{lurie-higher-topos-theory}).
\end{proof}

\begin{definition}
Let $\mathcal C$ be a site and let $\mathcal D$ be a presentable $\infty$-category.
\begin{defenum}
	\item For every presheaf $\mathcal F \in \Fun(\mathcal C^\opp, \mathcal D)$ we call $\mathcal F^\sharp \in \catsh{\mathcal C}{\mathcal D}$ the \emph{sheafification} of $\mathcal F$.

	\item For every sheaf $\mathcal F \in \catsh{\mathcal C}{\mathcal D}$ we call $\hat{\mathcal F} \in \cathsh{\mathcal C}{\mathcal D}$ the \emph{hypercompletion} of $\mathcal F$.
\end{defenum}
\end{definition}

We also record the following property of covers (i.e. effective epimorphisms) in the $\infty$-category of sheaves on a site, which is well-known in the $1$-categorical case.

\begin{lemma} \label{rslt:cover-of-sheaves-characterized-by-sections}
Let $\mathcal C$ be an $\infty$-site and let $(f_i\colon \mathcal F_i \surjto \mathcal G)_{i\in I}$ be a small family of maps in $\catshs{\mathcal C}$. Then the following are equivalent:
\begin{lemenum}
	\item The map $\bigdunion_i \mathcal F_i \to \mathcal G$ is an effective epimorphism.
	\item For every $X \in \mathcal C$ and every section $s \in \pi_0 \mathcal G(X)$ there are a cover $(U_j \to X)_{j\in J}$ in $\mathcal C$, a map $\alpha\colon J \to I$ and sections $t_j \in \pi_0\mathcal F_{\alpha(j)}(U_j)$ such that $f_{\alpha(j)}(t_j) = \restrict s{U_j}$.
\end{lemenum}
\end{lemma}
\begin{proof}
We start with the proof that (i) implies (ii). First assume that $I = \{ 0 \}$ and denote the map $f_0\colon \mathcal F_0 \to \mathcal G$ by $f\colon \mathcal F \to \mathcal G$. By taking the pullback along the map $X \to \mathcal G$ induced by $s$ and replacing $\mathcal C$ by $\mathcal C_{/X}$ we can assume that $\mathcal G = *$ and $X = *$ are final objects of $\mathcal C$. Let $\mathcal F'$ denote the geometric realization of the Čech nerve of $f\colon \mathcal F \to *$ in the $\infty$-category of presheaves on $\mathcal C$. Then the induced map $\mathcal F' \to *$ is a monomorphism and becomes an isomorphism after sheafification. Thus by \cite[Lemma 6.2.2.16]{lurie-higher-topos-theory} the full subcategory of $\mathcal C$ consisting of those $U$ such that $\mathcal F'(U) \ne \emptyset$ is a covering sieve for $\mathcal C$. Then the claim follows by picking any covering $(U_j \to *)_j$ with all $U_j$ from that covering sieve and noting that for each $U_i$ the map $\pi_0\mathcal F(U_j) \to \pi_0\mathcal F'(U_j)$ is surjective because $\mathcal F(U_j) \to \mathcal F'(U_j)$ is an effective epimorphism.

Now let $I$ be general. By the previous case we reduce to the case that $\mathcal G = \bigdunion_i \mathcal F_i$. Letting $\mathcal G'$ denote the disjoint union of the $\mathcal F_i$ in the category of presheaves on $\mathcal C$, the induced map $\mathcal G' \injto \mathcal G$ is a monomorphism. Now argue as in the previous case.

By reversing the above argument we deduce that (ii) implies (i).
\end{proof}

Sheaves are locally defined objects by design, i.e. an equivalence of sheaves can be checked locally on the site. Moreover, a sheaf on a site can uniquely be constructed by providing its values on any ``basis'' of the site. More precisely, we have the following:

\begin{definition}
A \emph{basis} of an $\infty$-site $\mathcal C$ is a full subcategory $\mathcal B \subset \mathcal C$ such that for every object $X \in \mathcal C$ there is a collection of morphisms $(U_i \to X)_{i\in I}$ generating a covering sieve of $X$ such that all $U_i$ lie in $\mathcal B$.

If $\mathcal B$ is a basis of the $\infty$-site $\mathcal C$, then there is a unique Grothendieck topology on $\mathcal B$ such that a sieve $\mathcal B^0_{/X}$ is a covering sieve if and only if its image under $\mathcal B_{/X} \injto \mathcal C_{/X}$ generates a covering sieve of $\mathcal C$. We implicitly regard $\mathcal B$ as an $\infty$-site using this Grothendieck topology.
\end{definition}

\begin{proposition} \label{rslt:sheaves-on-basis-equiv-sheaves-on-whole-site}
Let $\mathcal C$ be an $\infty$-site with basis $\mathcal B$ and let $\mathcal D$ be an $\infty$-category which has all small limits.
\begin{propenum}
	\item Right Kan extension along $\mathcal B \injto \mathcal C$ defines an equivalence
	\begin{align*}
		\cathsh{\mathcal B}{\mathcal D} \isoto \cathsh{\mathcal C}{\mathcal D},
	\end{align*}
	whose inverse is given by restriction.

	\item Right Kan extension along $\mathcal B \injto \mathcal C$ defines a fully faithful functor
	\begin{align*}
		\catsh{\mathcal B}{\mathcal D} \injto \catsh{\mathcal C}{\mathcal D}.
	\end{align*}
	This functor is an equivalence in either of the following situations:
	\begin{enumerate}[(a)]
		\item $\mathcal B$ and $\mathcal C$ have all finite limits and are $n$-categories for some $n$.

		\item The inclusion $\mathcal B \injto \mathcal C$ preserves fiber products and for every $X \in \mathcal C$ there is a cover $(U_i \to X)_i$ in $\mathcal C$ such that all fiber products $U_{i_0} \cprod_X \dots \cprod_X U_{i_n}$ exist and belong to $\mathcal B$.
	\end{enumerate}
\end{propenum}
\end{proposition}
\begin{proof}
If $\mathcal D = \Ani$ is the $\infty$-category of anima then (i) is \cite[Corollary A.7]{aoki-tensor-triangulated-geometry} and (ii) is \cite[Proposition A.5]{aoki-tensor-triangulated-geometry}, \cite[Corollary A.8]{aoki-tensor-triangulated-geometry} and \cite[Lemma C.3]{hoyois-quadratic-refinement-of-trace-formula}. The proof of \cite[Proposition A.5]{aoki-tensor-triangulated-geometry} immediately generalizes to general $\mathcal D$, proving the first part of (ii). Moreover, if $\mathcal X_{\mathcal B} = \catshs{\mathcal B}$ and $\mathcal X_{\mathcal C} = \catshs{\mathcal C}$ are the associated topoi then $\hat{\mathcal X}_{\mathcal B} = \hat{\mathcal X}_{\mathcal C}$ and under the assumptions of (ii).(a) and (ii).(b) we even have $\mathcal X_{\mathcal B} = \mathcal X_{\mathcal C}$. Using \cref{rslt:equivalence-of-hypercomplete-sheaf-defs,rslt:equivalence-of-sheaf-defs-on-site} we deduce the equivalence claims for general $\mathcal D$ (to see that the equivalences are indeed given by right Kan extensions we need to verify that these Kan extensions produce sheaves; this can be done in a similar way as in the last part of the proof of \cite[Proposition B.6.6]{lurie-ultracategories}).
\end{proof}

In practice one often constructs ($\infty$-)sites from explicit covering families, as is the case for all $\infty$-sites occurring in this thesis. In this case also the sheaf property is more explicit, as follows.

\begin{definition} \label{def:explicit-covering-site}
An \emph{explicit covering site} is a small $\infty$-category $\mathcal C$ together with a collection of morphisms $S$ in $\mathcal C$ satisfying the following properties:
\begin{enumerate}[(i)]
	\item $S$ contains all equivalences. Moreover, if $f\colon Y \to X$ and $g\colon Z \to Y$ are two morphisms in $\mathcal C$ with composition $h \isom f \comp g$ then: if $f, g \in S$ then $h \in S$ and if $h \in S$ then $f \in S$. In particular, if two morphisms $f$ and $f'$ in $\mathcal C$ are homotopic, then $f \in S$ iff $f' \in S$.

	\item $\mathcal C$ admits pullbacks and $S$ is stable under pullbacks.

	\item $\mathcal C$ admits finite coproducts and $S$ is stable under finite coproducts.

	\item Finite coproducts in $\mathcal C$ are universal. That is, given a diagram $\bigdunion_{1 \le i \le n} X_i \to Y \from Y'$ in $\mathcal C$, the canonical map
	\begin{align*}
		\bigdunion_{1 \le i \le n} (X_i \cprod_Y Y') \isoto \big(\bigdunion_{1\le i \le n} X_i\big) \cprod_Y Y'
	\end{align*}
	is an isomorphism.

	\item Coproducts in $\mathcal C$ are disjoint. That is, for all $X, Y \in \mathcal C$ the fiber product $X \cprod_{X \dunion Y} Y$ is an initial object of $\mathcal C$.
\end{enumerate}
\end{definition}

\begin{definition}[{cf. \cite[Proposition A.3.2.1]{lurie-spectral-algebraic-geometry}}]
Let $(\mathcal C, S)$ be an explicit covering site. The \emph{associated $\infty$-site} is given by the Grothendieck topology on $\mathcal C$ where a sieve $\mathcal C^0_{/X}$ is a covering sieve if and only if it contains a finite collection of morphisms $(U_i \to X)_{1 \le i \le n}$ such that the induced map $\bigdunion_i U_i \to X$ is in $S$. In the following we will always implicitly view an explicit covering site as an $\infty$-site by this construction.
\end{definition}

\begin{proposition} \label{rslt:sheaves-on-explicit-covering-site}
Let $(\mathcal C, S)$ be an explicit covering site, $\mathcal D$ any $\infty$-category which has all small limits and $\mathcal F\colon \mathcal C^\opp \to \mathcal D$ a functor. Then $\mathcal F$ is a $\mathcal D$-valued sheaf on $\mathcal C$ if and only if it satisfies the following properties:
\begin{enumerate}[(a)]
	\item $\mathcal F$ preserves finite products,
	\item for every morphism $f\colon U_0 \to X$ in $S$ with associated Čech nerve $U_\bullet \to X$, the morphism
	\begin{align*}
		\mathcal F(X) \isoto \varprojlim_{n\in \Delta} \mathcal F(U_n)
	\end{align*}
	is an isomorphism in $\mathcal D$.
\end{enumerate}
\end{proposition}
\begin{proof}
See \cite[Proposition A.3.3.1]{lurie-spectral-algebraic-geometry}.
\end{proof}

There is a similar description of hypercomplete sheaves on an explicit covering site if we replace the Čech covers by hypercoverings, as follows.

\begin{definition}
Let $(\mathcal C, S)$ be an explicit covering site. A \emph{hypercovering} in $\mathcal C$ is an augmented simplicial object $X_\bullet \to X_{-1}$ in $\mathcal C$ such that for every $n \ge 0$ the canonical map
\begin{align*}
	X_n \to (\cosk_{n-1}(X_\bullet \to X_{-1}))_n
\end{align*}
belongs to $S$. Here $\cosk_n$ denotes the augmented coskeleton functor, defined as right Kan extension along the inclusion $\Delta^\opp_{+,\le n} \subset \Delta^\opp_+$.
\end{definition}

\begin{proposition} \label{rslt:hypercomplete-sheaves-on-explicit-covering-site}
Let $(\mathcal C, S)$ be an explicit covering site, $\mathcal D$ any $\infty$-category which has all small limits and $\mathcal F\colon \mathcal C^\opp \to \mathcal D$ a functor. Then $\mathcal F$ is a hypercomplete $\mathcal D$-valued sheaf on $\mathcal C$ if and only if it satisfies the following properties:
\begin{enumerate}[(a)]
	\item $\mathcal F$ preserves finite products,
	\item for every hypercovering $U_\bullet \to X$ in $\mathcal C$, the morphism
	\begin{align*}
		\mathcal F(X) \isoto \varprojlim_{n\in \Delta} \mathcal F(U_n)
	\end{align*}
	is an isomorphism in $\mathcal D$.
\end{enumerate}
\end{proposition}
\begin{proof}
This is \cite[Proposition A.5.7.2]{lurie-spectral-algebraic-geometry}, except that in the reference they work with \emph{semi}-simplicial hypercoverings and not simplicial ones. To get the same result for simplicial hypercoverings, we adapt the proof. Namely, everything except the last part remain unchanged. We are now in the situation that we are given an $\infty$-connective sheaf $\mathcal F \in \catshs{\mathcal C}$ satisfying condition (ii) above and a final object $X \in \mathcal C$ and we need to show that $\mathcal F(X)$ is non-empty. As in the proof in the reference, consider the right fibration $\widetilde{\mathcal C} \to \mathcal C$ classified by the functor $\mathcal F\colon \mathcal C^\opp \to \Ani$. As in the reference it is enough to construct a simplicial hypercovering $X_\bullet \to X$ in $\mathcal C$ together with a lifting $Y_\bullet\colon \Delta^\opp \to \widetilde{\mathcal C}$ of $X_\bullet$.

We construct $X_{\bullet\le n}$ and $Y_{\bullet\le n}$ by induction on $n$. Assume we have finished the construction up to $n-1$. Let $M_n(X_\bullet) = (\cosk_{n-1}(X_\bullet))_n$. For the induction step we can argue in the same way as in the reference to obtain an object $D \in \mathcal C$ with a map $D \to M_n(X_\bullet)$ in $S$, together with extensions $X'_{\bullet\le n}\colon (\Delta'_{\le n})^\opp \to \mathcal C$ and $Y'_{\bullet\le n}\colon (\Delta'_{\le n})^\opp \to \widetilde{\mathcal C}$, where $\Delta'_{\le n}$ is defined as
\begin{align*}
	\Delta'_{\le n} := \Delta_{s,\le n} \dunion_{\Delta_{s,\le n-1}} \Delta_{\le n-1}.
\end{align*}
Now define $X_{\bullet\le n}\colon \Delta_{\le n}^\opp \to \mathcal C$ by left Kan extension of $X'_{\bullet\le n}$ along $(\Delta'_{\le n})^\opp \injto \Delta_{\le n}^\opp$. This leaves $X_m$ untouched for $m < n$ and replaces $X'_n$ by
\begin{align*}
	X_n = \bigdunion_{[n] \surjto [m]} X'_m.
\end{align*}
In particular, since finite disjoint unions exist in $\mathcal C$, the left Kan extension exists. We can similarly Kan extend $Y'_{\bullet\le n}$ to a functor $Y_{\bullet\le n}\colon \Delta_{\le n}^\opp \to \widetilde{\mathcal C}$ lifting $X$ (using that $\mathcal F\colon \mathcal C^\opp \to \Ani$ preserves finite products). Now note that $M_n(X_{\bullet\le n}) = M_n(X_{\bullet\le n-1})$ is unchanged, and the map $X_n \to M_n(X_{\bullet\le n})$ factors over $X'_n \to M_n(X_{\bullet\le n})$. As the latter map lies in $S$ (by construction), so does the former, proving that $X_\bullet \to X_{-1}$ is indeed a hypercovering.
\end{proof}

From now on (and throughout the thesis) we will implicitly use \cref{rslt:sheaves-on-explicit-covering-site,rslt:hypercomplete-sheaves-on-explicit-covering-site} whenever we are working with an explicit covering site.

Given different sites on the same underlying category $\mathcal C$, the sheafiness of a presheaf with respect to one site helps for checking sheafiness with respect to the other site:

\begin{lemma} \label{rslt:check-sheafiness-on-finer-topology}
Let $\mathcal D$ be a presentable $\infty$-category and let $\mathcal C$ be small a category which is equipped with two explicit covering sites $(\mathcal C, S)$ and $(\mathcal C, S')$ such that $S \subset S'$. Given a functor $\mathcal F\colon \mathcal C^\opp \to \mathcal D$, the following are equivalent:
\begin{lemenum}
	\item $\mathcal F$ is a sheaf on $(\mathcal C, S')$.
	\item $\mathcal F$ is a sheaf on $(\mathcal C, S)$ and for every $f\colon Y \to X$ in $S'$ there is a commutative diagram
	\begin{center}\begin{tikzcd}
		Y' \arrow[r,"g'"] \arrow[d,"f'"] & Y \arrow[d,"f"]\\
		X' \arrow[r,"g"] & X
	\end{tikzcd}\end{center}
	(not necessarily cartesian) such that $g$ and $g'$ are in $S$ and given the Čech nerve $Y'_\bullet \to X'$ of $Y' \to X'$, the natural morphism
	\begin{align*}
		\mathcal F(X') \to \varprojlim_{n\in\Delta} \mathcal F(Y'_\bullet)
	\end{align*}
	is an isomorphism.
\end{lemenum}
\end{lemma}
\begin{proof}
It is clear that (i) implies (ii). For the converse, let us say that $\mathcal F$ descends along a morphism $Y \to X$ in $\mathcal C$ if, given the Čech nerve $Y_\bullet \to X$ of $Y \to X$, the map $\mathcal F(X) \to \varprojlim_{n\in\Delta} \mathcal F(Y_\bullet)$ is an isomorphism. We say that $\mathcal F$ descends universally along $Y \to X$ if $\mathcal F$ descends along $Y \cprod_Z X \to Z$ for any $Z \to X$ (cf. \cite[Definition 3.1.1]{enhanced-six-operations}).

Now assume that $\mathcal F$ satisfies (ii). Given a morphism $f\colon Y \to X$ in $S'$ we have to show that $\mathcal F$ descends along $f$. Choose $X'$, $Y'$, $g$, $g'$ and $f'$ as in (ii). By the $(\mathcal C, S)$-sheafiness of $\mathcal F$, $\mathcal F$ descends universally along $g$ and $g'$ and by (ii) $\mathcal F$ descends along $f'$. Thus by \cite[Lemma 3.1.2.(4)]{enhanced-six-operations} $\mathcal F$ descends along $g \comp f'$ and hence along $h := f \comp g'$. Consider the commutative diagram
\begin{center}\begin{tikzcd}
	Y' \arrow[r] \arrow[dr,"g'",swap] & Y \cprod_X Y' \arrow[r,"\pr{Y'}"] \arrow[d,"\pr{Y}"] & Y' \arrow[d,"h"]\\
	& Y \arrow[r,"f"] & X
\end{tikzcd}\end{center}
in $\mathcal C$. Since $\mathcal F$ descends universally along $g'$, it also descends universally along $\pr{Y}$ by \cite[Lemma 3.1.2.(3)]{enhanced-six-operations}. Moreover, $\mathcal F$ descends universally along $\pr{Y'}$ because $\pr{Y'}$ is a retract (precomposing $\pr{Y'}$ with the natural map $Y' \to Y \cprod_X Y'$ yields the identity on $Y'$), see \cite[Lemma 3.1.2.(1)]{enhanced-six-operations}. Thus by \cite[Lemma 3.1.2.(2)]{enhanced-six-operations}, $\mathcal F$ descends along $f$.
\end{proof}

Our next result concerns hypercompleteness of sheaves. We will provide a criterion for when a sheaf of $\infty$-categories on a site is automatically hypercomplete. We first need some preliminaries.

\begin{definition}
We borrow the following notation from \cite[\S4.7.2]{lurie-higher-algebra}.
\begin{defenum}
	\item \label{def:functor-rho-on-Delta} Let $\rho\colon \Delta_+ \to \Delta$ be the functor sending $[n] \mapsto [n+1]$ and a map $\alpha\colon [n] \to [m]$ to the map $\rho(\alpha)\colon [n+1] \to [m+1]$ with $\rho(\alpha)(0) = 0$ and $\rho(\alpha)(k) = \alpha(k-1)$ for $k \ge 1$.

	\item \label{def:functor-T-on-simplicial-objects} Let $\mathcal C$ be an $\infty$-category. Then the functor $\rho$ induces a functor $T\colon X_\bullet \mapsto (X_{\bullet+1} \to X_0)$ from the $\infty$-category $\Fun(\Delta^\opp, \mathcal C)$ of simplicial objects of $\mathcal C$ to the $\infty$-category $\Fun(\Delta_+^\opp, \mathcal C)$ of augmented simplicial objects of $\mathcal C$. There is a natural map $X_{\bullet+1} \to X_\bullet$ of simplicial objects.
\end{defenum}
\end{definition}

We will apply the operator $T$ from \cref{def:functor-T-on-simplicial-objects} to hypercoverings $Y_\bullet \to X$ in a site $\mathcal C$. In the case that $Y_\bullet \to X$ is a Čech covering, one sees easily that $TY_\bullet = (Y_{\bullet+1} \to Y_0)$ is the pullback $TY_\bullet = (Y_\bullet \to X) \cprod_X Y_0$. We should therefore regard $T$ as a generalization of this pullback functor (the pullback functor itself does not behave as nicely on hypercoverings).

\begin{lemma} \label{rslt:compute-cosk-of-TX-bullet}
Let $\mathcal C$ be an $\infty$-category which has all finite limits and let $X_\bullet$ be a simplicial object of $\mathcal C$. Then for all $n \ge 0$ there is a natural isomorphism
\begin{align*}
	(\cosk_n X_\bullet)_{n+1} = (\cosk_{n-1} TX_\bullet)_n \cprod_{(\cosk_{n-1} X_\bullet)_n} X_n
\end{align*}
in $\mathcal C$. Here $\cosk_{n-1} TX_\bullet$ denotes the coskeleton of the \emph{augmented} simplicial object $TX_\bullet$.
\end{lemma}
\begin{proof}
Recall that for any simplicial object $Z_\bullet$ in $\mathcal C$ we have
\begin{align*}
	(\cosk_n Z_\bullet)_{n+1} = \varprojlim_{k \in (\Delta_{/[n+1]})^\opp_{\le n}} Z_k
\end{align*}
(see e.g. \cite[Lemma 0183]{stacks-project}). Let us abbreviate $\mathcal D_n := (\Delta_{/[n+1]})^\opp_{\le n}$ and $\mathcal D_{n,+}$ for the same with $\Delta_+$ in place of $\Delta$. Now fix $n$. Let $\mathcal D' \subset \mathcal D_n$ be the full subcategory spanned by the maps $\alpha\colon [k] \to [n+1]$ with $\alpha(0) = 0$. Similarly let $\mathcal D'' \subset \mathcal D_n$ be the full subcategory spanned by the maps $\alpha\colon [k] \to [n+1]$ with $k < n$ and $\alpha(0) > 0$. Finally let $v \in \mathcal D_n$ be the map $v\colon [n] \to [n+1]$ with $v(k) = k+1$. Then $\mathcal D'$, $\mathcal D''$ and $\{ v \}$ form a partition of the objects of $\mathcal D_n$. Moreover, the only maps in $\mathcal D_n$ which are not contained in one of these subcategories go from $\mathcal D''$ to $\mathcal D'$ or from $\mathcal D''$ to $v$. We deduce
\begin{align*}
	(\cosk_n X_\bullet)_{n+1} = \varprojlim_{k \in \mathcal D^\opp_n} X_k = (\varprojlim_{k \in {\mathcal D'}^\opp} X_k) \cprod_{(\varprojlim_{k \in {\mathcal D''}^\opp} X_k)} X_v.
\end{align*}
Now note that there is a natural equivalence of categories $\psi''\colon \mathcal D_{n-1} \isoto \mathcal D''$ mapping $\alpha\colon [k] \to [n]$ to the map $\psi(\alpha)\colon [k] \to [n+1]$ with $\psi(\alpha)(l) = \alpha(l)+1$. Under this equivalence we get a natural isomorphism
\begin{align*}
	\varprojlim_{k \in {\mathcal D''}^\opp} X_k = (\cosk_{n-1} X_\bullet)_n
\end{align*}
Similarly there is a natural isomorphism $\psi'\colon \mathcal D_{n-1,+} \isoto \mathcal D'$ sending $\alpha\colon [k] \to [n]$ to the map $\psi'(\alpha)\colon [k+1] \to [n+1]$ with $\psi'(\alpha)(l) = 0$ for $l = 0$ and $\psi'(\alpha)(l) = \alpha(l-1)+1$ for $l > 0$. Note that $\psi'$ factors through the map $\rho$ from \cref{def:functor-rho-on-Delta}, so that we deduce
\begin{align*}
	\varprojlim_{k \in {\mathcal D'}^\opp} X_k = (\cosk_{n-1} TX_\bullet)_n.
\end{align*}
The claim follows.
\end{proof}

\begin{corollary} \label{rslt:T-of-hypercover-is-split-hypercover}
Let $(\mathcal C, S)$ be an explicit covering site and $Y_\bullet \to X$ a hypercovering in $\mathcal C$. Then $TY_\bullet = (Y_{\bullet+1} \to Y_0)$ is a hypercovering in $\mathcal C$ and split in the sense of \cite[Defitinion 4.7.2.2]{lurie-higher-algebra}.
\end{corollary}
\begin{proof}
It follows from \cite[Corollary 4.7.2.15]{lurie-higher-algebra} that $TY_\bullet$ is split. To see that it is a hypercovering, note that for every $n \ge -1$ the map $(TY_\bullet)_{n+1} = Y_{n+2} \to (\cosk_n TY_\bullet)_{n+1}$ factors by \cref{rslt:compute-cosk-of-TX-bullet} as
\begin{align*}
	Y_{n+2} \to (\cosk_{n+1} Y_\bullet)_{n+2} = (\cosk_n TY_\bullet)_{n+1} \cprod_{(\cosk_n Y_\bullet)_{n+1}} Y_{n+1} \to (\cosk_n TY_\bullet)_{n+1},
\end{align*}
where the second map is the projection to the first factor. Here the non-augmented coskeleta are computed in $\mathcal C_{/X}$. We want the composed map to be a cover in $\mathcal C$ (i.e. $\in S$), so it is enough to see that both maps are covers. For the first map this is clear by definition of hypercoverings. For the second map it follows by base-change from the fact that the map $Y_{n+1} \to (\cosk_n Y_\bullet)_{n+1}$ is a cover in $\mathcal C$ (again by definition of $Y_\bullet \to X$ being a hypercovering).
\end{proof}

We conclude this section with the following criterion for checking whether a sheaf of $\infty$-categories is hypercomplete. Note that $\infcatinf$ is presentable (e.g. by \cref{rslt:infcatinf-is-compactly-generated}), so all of the above results on sheaves apply to $\infcatinf$-valued sheaves.

\begin{proposition} \label{rslt:fully-faithful-on-hypercovers-implies-hypercomplete}
Let $(\mathcal C, S)$ be an explicit covering site and let $\mathcal F\colon \mathcal C^\opp \to \infcatinf$ be a sheaf of $\infty$-categories on $\mathcal C$. Assume that for every hypercovering $Y_\bullet \to X$ in $\mathcal C$, the natural functor
\begin{align*}
	\mathcal F(X) \injto \varprojlim_{n\in\Delta} \mathcal F(Y_n)
\end{align*}
is fully faithful. Then $\mathcal F$ is hypercomplete.
\end{proposition}

\begin{remark}
Intuitively, the full faithfulness condition in \cref{rslt:fully-faithful-on-hypercovers-implies-hypercomplete} means that for every $X \in \mathcal C$, every object $x \in \mathcal F(X)$ is a hypercomplete sheaf on $\mathcal C_{/X}$. In that spirit, \cref{rslt:fully-faithful-on-hypercovers-implies-hypercomplete} is similar to the proof of \cite[Theorem 11.2.(4)]{bhatt-scholze-witt} (on page 51), but we do not require the a priori existence of a nice map $\mathcal F \to \mathcal G$ to a hypercomplete sheaf $\mathcal G$.
\end{remark}

\begin{proof}[Proof of \cref{rslt:fully-faithful-on-hypercovers-implies-hypercomplete}]
Fix any hypercovering $Y_\bullet \to X$. Let $\tilde Y_{\bullet,\bullet}$ be the Čech nerve of the map $TY_\bullet \to (Y_\bullet \to X)$ of augmented simplicial objects. Let us also denote by $\tilde X_\bullet \to X$ the Čech nerve of $Y_0 \to X$. Then for fixed $k \ge 0$, $\tilde Y_{\bullet,k} \to \tilde X_k$ is a hypercovering in $\mathcal C$ (by \cref{rslt:T-of-hypercover-is-split-hypercover}) and for fixed $k \ge 0$, $\tilde Y_{n,\bullet} \to Y_n$ is a Čech cover in $\mathcal C$ (to see that the map $Y_{n+1} \to Y_n$ is a cover we can e.g. use \cref{rslt:compute-cosk-of-TX-bullet}). Since $\mathcal F$ is a sheaf on $\mathcal C$ we have
\begin{align*}
	\varprojlim_{n\in\Delta} \mathcal F(Y_n) = \varprojlim_{n\in\Delta} \varprojlim_{k\in\Delta} \mathcal F(\tilde Y_{n,k}) = \varprojlim_{k\in\Delta} \left(\varprojlim_{n\in\Delta} \mathcal F(\tilde Y_{n,k})\right).
\end{align*}
For every $k \ge 0$ let $\mathcal C^k \subset \varprojlim_{n\in\Delta} \mathcal F(\tilde Y_{n,k})$ be the essential image of the functor $\varprojlim_{n\in\Delta} \mathcal F(Y_n) \to \varprojlim_{n\in\Delta} \mathcal F(\tilde Y_{n,k})$. Then we have functors
\begin{align*}
	\varprojlim_{n\in\Delta} \mathcal F(Y_n) \to \varprojlim_{k\in\Delta} \mathcal C^k \injto \varprojlim_{k\in\Delta} \left(\varprojlim_{n\in\Delta} \mathcal F(\tilde Y_{n,k})\right),
\end{align*}
whose composition is the above equivalence and where the second functor is fully faithful. It follows that both functors are equivalences, i.e.
\begin{align*}
	\varprojlim_{n\in\Delta} \mathcal F(Y_n) = \varprojlim_{k\in\Delta} \mathcal C^k.
\end{align*}
By \cref{rslt:T-of-hypercover-is-split-hypercover} the hypercovering $(\tilde Y_{\bullet,0} \to \tilde X_0) = TY_\bullet$ is split in the sense of \cite[Defitinion 4.7.2.2]{lurie-higher-algebra}, from which it follows automatically that the functor
\begin{align*}
	\mathcal F(\tilde X_0) \isoto \varprojlim_{n\in\Delta} \mathcal F(\tilde Y_{n,0})
\end{align*}
is an equivalence (see \cite[Remark 4.7.2.3]{lurie-higher-algebra}).\footnote{This is saying that any descent datum for $Y_\bullet \to X$ becomes effective for $Y_{\bullet+1} \to Y_0$, i.e. after ``pullback'' from $X$ to $Y_0$.} It follows that for all $k \ge 0$, $\mathcal C^k$ is contained in the essential image of the functor $\mathcal F(\tilde X_k) \to \varprojlim_{n\in\Delta} \mathcal F(\tilde Y_{n,k})$ (for $k = 0$ this is clear by the above equivalence; then it follows also for $k > 0$). But by assumption on $\mathcal F$ the functor $\mathcal F(\tilde X_k) \injto \varprojlim_{n\in\Delta} \mathcal F(\tilde Y_{n,k})$ is fully faithful, so we can view $\mathcal C^k$ as a full subcategory of $\mathcal F(\tilde X_k)$. Then
\begin{align*}
	\mathcal F(X) \injto \varprojlim_{n\in\Delta} \mathcal F(Y_n) = \varprojlim_{k\in\Delta} \mathcal C^k \injto \varprojlim_{k\in\Delta} \mathcal F(\tilde X_k) = \mathcal F(X)
\end{align*}
and the composition is the identity. It follow that both fully faithful functors are equivalences, i.e. $\mathcal F(X) \isoto \varprojlim_{n\in\Delta} \mathcal F(Y_n)$ as desired.
\end{proof}

\subsection{Enriched \texorpdfstring{$\infty$}{Infinity}-Categories} \label{sec:infcat.enriched}

In the following we prove some basic facts about enriched $\infty$-categories. The goal is not to give a full introduction to enrichment in the $\infty$-categorical setting, but merely to provide enough machinery so that we can use them efficiently in the main part of the thesis.

We will make use of the operadic model \cite[\S4.2.1]{lurie-higher-algebra} for enriched $\infty$-categories, as it seems to provide the easiest access to enriched functors. Let us briefly recall the main definitions. In the following, we use the $\infty$-operad $\opLM$ from \cite[Definition 4.2.1.7]{lurie-higher-algebra}, whose objects over $\langle 1 \rangle \in \catFinAst$ are denoted $\mathfrak a$ and $\mathfrak m$.

\begin{definition}
A \emph{weak enrichment} of an $\infty$-category $\mathcal C$ over a monoidal $\infty$-category $\mathcal V^\tensor$ is a fibration of $\infty$-operads $\mathcal C^\tensor \to \opLM$ together with isomorphisms $\mathcal C_{\mathfrak a}^\tensor \isom \mathcal V^\tensor$ and $\mathcal C^\tensor_{\mathfrak m} \isom \mathcal C$.
\end{definition}

\begin{definition}
Let $q\colon \mathcal C^\tensor \to \opLM$ be a weak enrichment of the $\infty$-category $\mathcal C \isom \mathcal C^\tensor_{\mathfrak m}$ over the monoidal $\infty$-category $\mathcal V^\tensor \isom \mathcal C_{\mathfrak a}^\tensor$.
\begin{defenum}
	\item We say that $q$ exhibits $\mathcal C$ as \emph{(left-)tensored over $\mathcal V^\tensor$} if $q$ is a coCartesian fibration of $\infty$-operads.

	\item \label{def:enriched-infty-category} We say that $q$ exhibits $\mathcal C$ as \emph{enriched over $\mathcal V^\tensor$} if $q$ satisfies the following properties:
	\begin{enumerate}[(i)]
		\item For every sequence of objects $V_1, \dots, V_n \in \mathcal V$ and every pair of objects $X, Y \in \mathcal C$ we denote by $\Hom(\{ V_1, \dots, V_n \} \tensor X, Y)$ the summand of the mapping space $\Mul_{\mathcal C^\tensor}(\{ V_1, \dots, V_n, X \}, \{ Y \})$ corresponding to the linear ordering $1 \le \dots \le n$ on $\langle n \rangle$. We then require that the natural map
		\begin{align*}
			\Hom(\{ V_1 \tensor \dots \tensor V_n \} \tensor X, Y) \isoto \Hom(\{ V_1, \dots, V_n \} \tensor X, Y)
		\end{align*}
		is an isomorphism.

		\item For every pair of objects $X, Y \in \mathcal C$ there is an object $\IHom_{\mathcal C}(X, Y) \in \mathcal V$ together with a map $\alpha \in \Hom(\{ \IHom_{\mathcal C}(X, Y) \} \tensor X, Y)$ with the following universal property: For every object $V \in \mathcal V$, composition with $\alpha$ induces an isomorphism of spaces
		\begin{align*}
			\Hom(V, \IHom_{\mathcal C}(X, Y)) \isoto \Hom(\{ V \} \tensor X, Y).
		\end{align*}

		\item For every object $X \in \mathcal C$ there is a $q$-coCartesian morphism $\id_X \in \Hom(\{ 1_{\mathcal V} \} \tensor X, X)$ and composition with $\id_X$ yields an isomorphism of spaces
		\begin{align*}
			\Hom(X, X) \isoto \Hom(\{ 1_{\mathcal V} \} \tensor X, X).
		\end{align*}
	\end{enumerate}
\end{defenum}
\end{definition}

Note that if $\mathcal C$ is left-tensored over $\mathcal V^\tensor$ then there is a tensor product $\tensor\colon \mathcal V \cprod \mathcal C \to \mathcal C$. In this case condition (i) of \cref{def:enriched-infty-category} is automatically satisfied (both sides of the desired isomorphism identify with $\Hom(V_1 \tensor \dots \tensor V_n \tensor X, Y)$) and similarly for condition (iii). Moreover, condition (ii) of \cref{def:enriched-infty-category} is satisfied if for every $X \in \mathcal C$ the functor $- \tensor X\colon \mathcal V \to \mathcal C$ admits a right adjoint.

\begin{example} \label{ex:enriched-category-via-monoidal-functor}
Suppose $F\colon \mathcal V \to \mathcal C$ is a monoidal functor of monoidal $\infty$-categories. Then $F$ gives $\mathcal C$ the structure of a left-tensored $\infty$-category over $\mathcal V$ (e.g. view $F$ as a morphism of algebras in $\infcatinf$ so that $\mathcal C$ is in particular a module over $\mathcal V$, cf. \cite[Remark 2.4.2.6]{lurie-higher-algebra}). Here the tensor product is given by
\begin{align*}
	\mathcal V \cprod \mathcal C \to \mathcal C, \qquad (V, X) \mapsto F(V) \tensor_{\mathcal C} X.
\end{align*}
If for every $X \in \mathcal C$ the functor $- \tensor X\colon \mathcal V \to \mathcal X$ admits a right adjoint, then $\mathcal C$ is also enriched over $\mathcal V$.
\end{example}

Having defined a notion of enriched $\infty$-categories, it is time to introduce the categories of functors in that setting. The following definition is not found in \cite{lurie-higher-algebra}, but is closely related to the stricter definition of \emph{linear functors}, see \cite[Definition 4.6.2.7]{lurie-higher-algebra}.

\begin{definition} \label{def:enriched-functors}
Let $\mathcal V^\tensor$ be a monoidal $\infty$-category and let $p\colon \mathcal C^\tensor \to \opLM$ and $q\colon \mathcal D^\tensor \to \opLM$ be fibrations of $\infty$-categories which exhibit $\mathcal C := \mathcal C^\tensor_{\mathfrak m}$ and $\mathcal D := \mathcal D^\tensor_{\mathfrak m}$ as enriched over $\mathcal V^\tensor$. A \emph{$\mathcal V$-enriched functor} from $\mathcal C$ to $\mathcal D$ is a lax $\opLM$-monoidal functor $\mathcal C^\tensor \to \mathcal D^\tensor$ which is the identity on $\mathcal V^\tensor$, i.e. it is a functor $f\colon \mathcal C^\tensor \to \mathcal D^\tensor$ satisfying the following conditions:
\begin{enumerate}[(i)]
	\item $f$ is a functor over $\opLM$, i.e. $p = q \comp F$.
	\item $f$ carries inert morphisms in $\mathcal C^\tensor$ to inert morphisms in $\mathcal D^\tensor$.
	\item The restriction of $f$ to $\mathcal V^\tensor = \mathcal C^\tensor_{\mathfrak a} = \mathcal D^\tensor_{\mathfrak a}$ is the identity.
\end{enumerate}
We denote
\begin{align*}
	\Fun_{\mathcal V}(\mathcal C, \mathcal D) \subset \Fun_{\opLM}(\mathcal C^\tensor, \mathcal D^\tensor) \cprod_{\Fun_{\opAssoc}(\mathcal V^\tensor, \mathcal V^\tensor)} \{ \id \}
\end{align*}
the full subcategory spanned by the $\mathcal V$-enriched functors.
\end{definition}

In the setting of \cref{def:enriched-functors}, suppose that $f\colon \mathcal C \to \mathcal D$ is a $\mathcal V$-enriched functor. Then $f$ has an underlying functor $\mathcal C \to \mathcal D$ of $\infty$-categories (given by restricting $f$ to $\mathcal C \subset \mathcal C^\tensor$) and for every pair of objects $X, Y \in \mathcal C$, $f$ induces a natural morphism
\begin{align*}
	\IHom_{\mathcal C}(X, Y) \to \IHom_{\mathcal D}(f(X), f(Y))
\end{align*}
in $\mathcal V$. Namely, this morphism is by definition the same as an object in $\Hom(\{ \IHom_{\mathcal C}(X, Y) \} \tensor f(X), f(Y))$, for which we can take the image of the natural morphism under the map
\begin{align*}
	\Hom(\{ \IHom_{\mathcal C}(X, Y) \} \tensor X, Y) \to \Hom(\{ \IHom_{\mathcal C}(X, Y) \} \tensor f(X), f(Y))
\end{align*}
induced by $f$. This shows that a $\mathcal V$-enriched functor as defined above in \cref{def:enriched-functors} really contains the expected data. Note that if $\mathcal C$ and $\mathcal D$ are additionally tensored over $\mathcal V$, then the $\mathcal V$-enrichment of a functor $f\colon \mathcal C \to \mathcal D$ can equivalently be described by specifying for every $V \in \mathcal V$ and $X \in \mathcal C$ a natural morphism
\begin{align*}
	V \tensor f(X) \to f(V \tensor X)
\end{align*}
in $\mathcal D$ (plus higher coherences). Note that $f$ is $\mathcal V$-linear (as defined in \cite[Definition 4.6.2.7]{lurie-higher-algebra}) if the morphisms $V \tensor f(X) \isoto f(V \tensor X)$ are isomorphisms.

We are particularly interested in the $\infty$-category of enriched functors, so in the following we study some of its basic properties. We start by recording the expected properties of the forgetful functor:

\begin{lemma} \label{rslt:enriched-forgetful-functor-properties}
Let $\mathcal V$ be a monoidal $\infty$-category and let $\mathcal C$ and $\mathcal D$ be $\mathcal V$-enriched $\infty$-categories.
\begin{lemenum}
	\item \label{rslt:enriched-forgetful-functor-is-conservative} There is a natural forgetful functor
	\begin{align*}
		\Fun_{\mathcal V}(\mathcal C, \mathcal D) \to \Fun(\mathcal C, \mathcal D)
	\end{align*}
	which is conservative.

	\item \label{rslt:enriched-forgetful-functor-preserves-lim-colim} Suppose that $\mathcal D$ is even tensored over $\mathcal V$ and let $K$ be a simplicial set.
	\begin{enumerate}[(a)]
		\item Suppose that $\mathcal D$ admits $K$-indexed limits. Then the functor categories in (i) admit $K$-indexed limits and the forgetful functor preserves them.

		\item Suppose that $\mathcal D$ admits $K$-indexed colimits and the operation $\tensor\colon \mathcal V \cprod \mathcal D \to \mathcal D$ preserves $K$-indexed colimits in the second argument. Then the functor categories in (i) admit $K$-indexed colimits and the forgetful functor preserves them.
	\end{enumerate}
\end{lemenum}
\end{lemma}
\begin{proof}
Let $p\colon \mathcal C^\tensor \to \opLM$ and $q\colon \mathcal D^\tensor \to \opLM$ be the fibrations of $\infty$-operads exhibiting the enrichments of $\mathcal C$ and $\mathcal D$ over $\mathcal V$. The forgetful functor in (i) is given by restricting a $\mathcal V$-enriched functor $\mathcal C^\tensor \to \mathcal D^\tensor$ to $\mathcal C$. To show that the forgetful functor is conservative, let $f \to g$ be any morphism of $\mathcal V$-enriched functors $\mathcal C^\tensor \to \mathcal D^\tensor$ whose restriction to $\mathcal C$ is an equivalence. By definition both $f$ and $g$ restrict to the identity on $\mathcal V \subset \mathcal C^\tensor$ and the morphism $f \to g$ restricts to an isomorphism on $\mathcal V$. But any object in $\mathcal C^\tensor_{\langle 1 \rangle}$ lies either in $\mathcal V$ or in $\mathcal C$, hence the morphism $f \to g$ is an equivalence on $\mathcal C^\tensor_{\langle 1 \rangle}$. Using that $\mathcal C^\tensor_{\langle n \rangle} \isom \mathcal (C^\tensor_{\langle 1 \rangle})^n$ for every $n \ge 0$, compatibly with $f$ and $g$ (because they preserve inert morphisms), we deduce that $f \to g$ is an isomorphism on all of $\mathcal C^\tensor$, as desired. This finishes the proof of (i).

We now prove (ii), so assume that $q$ is a coCartesian fibration. Suppose first that $\mathcal D$ admits $K$-indexed limits and let $f\colon K \to \Fun_{\mathcal V}(\mathcal C, \mathcal D)$ be a diagram, which we can equivalently view as a functor $F\colon K \cprod \mathcal C^\tensor \to \mathcal D^\tensor$. We claim that $F$ admits a $q$-right Kan extension $\overline F\colon K^\triangleleft \cprod \mathcal C^\tensor \to \mathcal D^\tensor$. This boils down to checking that for every $X \in \mathcal C^\tensor$ the diagram $F_X\colon K \to \mathcal D^\tensor$ given by $F_X = F(-, X)$ admits a $q$-limit. Note that the whole diagram $F_X$ lies in the fiber $\mathcal D^\tensor_{\mathfrak x}$ over $\mathfrak x := p(X) \in \opLM$ and by \cite[Corollary 4.3.1.15]{lurie-higher-topos-theory} it is enough to show that the restricted diagram $F_X\colon K \to \mathcal D^\tensor_{\mathfrak x}$ admits a limit (this is then automatically a $q$-limit in $\mathcal D^\tensor$). By definition of a coCartesian fibration of $\infty$-operads (see \cite[Definition 2.1.2.13]{lurie-higher-algebra}) $\mathcal D^\tensor_{\mathfrak x}$ is isomorphic to a product of copies of $\mathcal D_{\mathfrak a} = \mathcal V$ and $\mathcal D_{\mathfrak m} = \mathcal D$, compatibly with the diagram $F$. We can therefore reduce to the cases $\mathfrak x = \mathfrak a$ (i.e. $X \in \mathcal V$) and $\mathfrak x = \mathfrak m$ (i.e. $X \in \mathcal C$). In the former case the diagram $F_X$ is trivial by definition of $\Fun_{\mathcal V}(\mathcal C, \mathcal D)$ and in the latter case the diagram $F_X$ is a $K$-indexed diagram in $\mathcal D$, which admits a limit by assumption. We have thus shown the existence of the $q$-right Kan extension $\overline F$, which we can also view as a functor $\overline f\colon K^\triangleleft \to \Fun(\mathcal C^\tensor, \mathcal D^\tensor)$. One checks easily using the above computations that $\overline f$ factors over $\Fun_{\mathcal V}(\mathcal C, \mathcal D)$. By \cite[Lemma 3.2.2.9]{lurie-higher-algebra} $\overline f$ is a $q'$-limit diagram, where $q'$ is the functor
\begin{align*}
	q'\colon \Fun(\mathcal C^\tensor, \mathcal D^\tensor) \to \Fun(\mathcal C^\tensor, \opLM)
\end{align*}
given by postcomposing with $q$. Note that $q'$ is an inner fibration by \cite[Corollary 2.3.2.5]{lurie-higher-topos-theory} and the restriction of $q'$ to $\Fun_{\mathcal V}(\mathcal C, \mathcal D)$ is constant with image $p$. By \cite[Proposition 4.3.1.5.(4)]{lurie-higher-topos-theory} the $q'$-limit diagram $\overline f$ restricts to a limit diagram in the fiber of $q'$ over $p$, proving the existence of the limit over $f$. By the explicit computation of $f$ above, the forgetful functor in (i) preserves the limit of $f$. This finishes the proof of (ii).(a).

We now prove (ii).(b), so assume that $\mathcal D$ admits $K$-indexed colimits and that the tensor product operations preserve $K$-indexed colimits in the second argument. Let $f\colon K \to \Fun_{\mathcal V}(\mathcal C, \mathcal D)$ be a diagram, which we can view as a functor $F\colon K \cprod \mathcal C^\tensor \to \mathcal D^\tensor$. We claim that $F$ admits a $q$-left Kan extension $\overline F\colon K^\triangleright \cprod \mathcal C^\tensor \to \mathcal D^\tensor$. This amounts to showing that for every $X \in \mathcal C$ the diagram $F_X\colon K \to \mathcal D^\tensor$ given by $F_X = F(-, X)$ admits a $q$-colimit. By \cite[Proposition 4.3.1.10]{lurie-higher-topos-theory} it is enough to show that the diagram $F_X\colon K \to \mathcal D^\tensor_{q(X)}$ has a colimit and for every morphism $\alpha\colon q(X) \to \mathfrak y$ in $\opLM$ the associated functor $\alpha_!\colon \mathcal D^\tensor_{q(X)} \to \mathcal D^\tensor_{\mathfrak y}$ preserves that colimit. The first part reduces to the existence of $K$-indexed colimits in $\mathcal D$ (as in the proof of (ii).(a) above) and is thus true by assumption. For the second part, let $\alpha$ be given. If $\alpha$ is inert then $\alpha_!$ is a projection of a product of two $\infty$-categories to one factor, so there is nothing to show. We can thus assume that $\alpha$ is active and we can then reduce to the case that $\mathfrak y = \mathfrak m$ (using that projections preserves colimits again). In this case $\mathcal D^\tensor_{q(X)} \isom \mathcal V^n \cprod \mathcal D$ for some $n$, $\mathcal D^\tensor_{\mathfrak y} = \mathcal D$ and $\alpha_!$ is the functor $(V_1, \dots, V_n, X) \mapsto (V_1 \tensor \dots \tensor V_n) \tensor X$. Since the diagram $F_X$ is constant on the $\mathcal V$-terms, showing that $\alpha_!$ preserves the colimit of $F_X$ follows immediately from the fact that the tensor product $\mathcal V \cprod \mathcal D \to \mathcal D$ preserves $K$-indexed colimits in the second argument. This shows the existence of the Kan extension $\overline F$. The rest of the argument is the same as in (ii).(a).
\end{proof}

The $\infty$-category $\infcatinf$ contains an object $*$ characterized by the property that $\Fun(*, \mathcal C) = \mathcal C$ for every $\infty$-category $\mathcal C$. This phenomenon generalizes to the enriched setting, as follows.

\begin{definition}
Let $\mathcal V$ be a monoidal $\infty$-category. Then $\mathcal V$ is tensored over itself and hence defines a coCartesian fibration $\mathcal V^{\tensor\tensor} \to \opLM$ of $\infty$-operads. We denote
\begin{align*}
	*^\tensor_{\mathcal V} \subset \mathcal V^{\tensor\tensor}
\end{align*}
the full suboperad spanned by the full subcategory $*_{\mathcal V} := \{ 1_{\mathcal V} \} \subset \mathcal V$, where $1_{\mathcal V}$ denotes the unit object of $\mathcal V$. Note that $*^\tensor_{\mathcal V}$ exhibits $*_{\mathcal V}$ as a $\mathcal V$-enriched $\infty$-category with $\IHom_{*_{\mathcal V}}(1_{\mathcal V}, 1_{\mathcal V}) = 1_{\mathcal V}$.
\end{definition}

\begin{lemma} \label{rslt:univ-prop-of-trivial-enriched-category}
Let $\mathcal V$ be a monoidal $\infty$-category and let $\mathcal C$ be a $\mathcal V$-enriched $\infty$-category. Then there is a natural equivalence of $\infty$-categories
\begin{align*}
	\Fun_{\mathcal V}(*_{\mathcal V}, \mathcal C) = \mathcal C.
\end{align*}
\end{lemma}
\begin{proof}
This is more subtle than one might think at first, because the $\infty$-operad $*^\tensor_{\mathcal V}$ contains $\mathcal V^\tensor$ (the $\infty$-operad exhibiting the monoidal structure on $V$) and is therefore rather non-trivial. However, Lurie's theory of free algebras (which are essentially Kan extensions in the $\infty$-operad setting) comes to the rescue. Namely, let $\opTriv$ denote the trivial $\infty$-operad (see \cite[Example 2.1.1.20]{lurie-higher-algebra}) and let $\mathcal T^\tensor := \mathcal V^\tensor \boxplus \opTriv$, where $\boxplus$ denotes the coproduct of $\infty$-operads constructed in \cite[Construction 2.2.3.3]{lurie-higher-algebra}. Then there is a natural map of $\infty$-operads $\mathcal T^\tensor \to *_{\mathcal V}^\tensor$ which is the identity on $\mathcal V^\tensor$ and sends $\langle 1 \rangle \in \opTriv$ to the object $1_{\mathcal V} \in *_{\mathcal V}$. This map identifies $\mathcal T^\tensor$ as a subcategory $\mathcal T^\tensor \subset *_{\mathcal V}^\tensor$. Let $p\colon \mathcal C^\tensor \to \opLM$ denote the fibration of $\infty$-operads exhibiting the enrichment of $\mathcal C$. We claim that there is an adjoint pair of functors
\begin{align*}
	F\colon \Alg'_{\mathcal T^\tensor/\opLM}(\mathcal C^\tensor) \rightleftarrows \Alg'_{*_{\mathcal V}^\tensor/\opLM}(\mathcal C^\tensor) \noloc G,
\end{align*}
where $G$ is the functor induced by restriction along $\mathcal T^\tensor \subset *_{\mathcal V}^\tensor$ and $\Alg'$ denotes the subcategory of those maps of $\infty$-operads which are the identity on $\mathcal V^\tensor$ and where do not allow non-trivial morphisms on the restriction to $\mathcal V^\tensor$. Here $F$ will be an ``operadic left Kan extension'' along $\mathcal T^\tensor \subset *_{\mathcal V}^\tensor$. More precisely, we employ \cite[Corollary 3.1.3.4]{lurie-higher-algebra}, so we have to show the following: Given an object $f \in \Alg'_{\mathcal T^\tensor/\opLM}(\mathcal C^\tensor)$ and an object $X \in (*_{\mathcal V}^\tensor)_{\langle 1 \rangle}$, the diagram $f\colon (\mathcal T^\tensor_\act)_{/X} \to \mathcal C^\tensor_\act$ can be extended to an operadic $p$-colimit diagram over $\opLM$. Here $(\mathcal T^\tensor_\act)_{/X}$ denotes the $\infty$-category of active morphisms $Y \to X$ in $*_{\mathcal V}^\tensor$, with morphisms in $(\mathcal T^\tensor_\act)_{/X}$ being the morphisms $Y \to Y'$ over $X$ which lie in $\mathcal T^\tensor_\act$. If $X \in \mathcal V^\tensor$ then $(\mathcal T^\tensor_\act)_{/X}$ has the final object $\id\colon X \to X$, so there is nothing to show. Thus from now on we can assume $X = 1_{\mathcal V} \in *_{\mathcal V}$; we will still denote it $X$ to distinguish it from the element $1_{\mathcal V} \in \mathcal V^\tensor$. Now the objects of $(\mathcal T^\tensor_\act)_{/1_{\mathcal V}}$ are the active morphisms $V_1 \dsum \dots V_n \dsum 1_{\mathcal V} \to \mathcal V$ in $*_{\mathcal V}^\tensor$ where $n \ge 0$ and $V_1, \dots, V_n \in \mathcal V$ (here $\dsum$ is defined as in \cite[Remark 2.1.1.15]{lurie-higher-algebra}). Let $\mathcal I' \subset (\mathcal T^\tensor_\act)_{/1_{\mathcal V}}$ denote the full subcategory spanned by only those morphisms with $n = 1$. We claim that $\mathcal I$ is cofinal in $(\mathcal T^\tensor_\act)_{/1_{\mathcal V}}$. First note that removing the morphisms with $n = 0$ is cofinal: Take $\alpha\colon X \to 1_{\mathcal V} \dsum X$ to be the morphism obtained from the coCartesian morphism $0 \to 1_{\mathcal V}$ in $V^\tensor$ by $- \dsum X$. Then morphisms $X \to X$ are the same as morphisms $1_{\mathcal V} \dsum X \to X$ lying over the fixed morphism in $\opLM$ corresponding to the natural odering of $\langle 2 \rangle$, and $\alpha$ provides maps between them, which easily shows that removing the morphisms $X \to X$ is indeed cofinal (e.g. use \cite[Theorem 4.1.3.1]{lurie-higher-topos-theory}). Now for $n \ge 2$, note that morphisms $V_1 \dsum \dots V_n \dsum X \to X$ over any given active morphism in $\opLM$ are naturally the same as morphisms $(V_1 \tensor \dots V_n) \dsum X \to X$ (over the obvious morphism in $\opLM$) with canonical interpolating morphisms between them, from which one deduces that $\mathcal I'$ is indeed cofinal in $(\mathcal T^\tensor_\act)_{/1_{\mathcal V}}$.

Let now $\mathcal I \subset \mathcal I'$ be the full subcategory consisting only of those morphisms $V \dsum X \to X$ which lie over the morphism in $\opLM$ that is the standard ordering on $\langle 2 \rangle$. Then $\mathcal I$ is cofinal in $\mathcal I'$ because there is an inert morphism $V \dsum X \to V \dsum X$ in $\mathcal T^\tensor$ interpolating between the two orderings on $\langle 2 \rangle$. It now remains to find an operadic $p$-colimit over the diagram $f\colon \mathcal I \to \mathcal C^\tensor$, $V \dsum X \mapsto V \dsum f(X)$. But note that $\mathcal I$ is in fact equivalent to $\mathcal V_{/1_{\mathcal V}}$ and hence admits the final object $1_{\mathcal V}$. We are thus reduced to showing that the diagram $\{ * \} \to \mathcal C^\tensor$, $* \mapsto 1_{\mathcal V} \dsum f(X)$ admits an operadic $p$-colimit in $\mathcal C$. By \cite[Lemma 3.1.1.14]{lurie-higher-algebra} this amounts to saying that there is a $p$-coCartesian active morphism $1_{\mathcal V} \dsum f(X) \to f(X)$ in $\mathcal C^\tensor$, which is true by definition of enriched $\infty$-categories.

We have now constructed the functor $F$. Note that since $\opTriv$-algebras in any $\infty$-operad are given by objects in the underlying $\infty$-category (see \cite[Remark 2.1.3.6]{lurie-higher-algebra}), we can naturally identify $\Alg'_{\mathcal T^\tensor/\opLM}(\mathcal C^\tensor)$ with $\mathcal C$. Note furthermore that $\Alg'_{*_{\mathcal V}^\tensor/\opLM}(\mathcal C^\tensor) = \Fun_{\mathcal V}(*_{\mathcal V}, \mathcal C)$ by definition. Under this identification, $G$ is the functor $\Fun_{\mathcal V}(*_{\mathcal V}, \mathcal C) \to \mathcal C$ given by evaluation at $1_{\mathcal V} \in *_{\mathcal V}$. This functor factors over the forgetful functor to $\Fun(*_{\mathcal V}, \mathcal C)$ which is conservative by \cref{rslt:enriched-forgetful-functor-is-conservative}; thus $G$ is conservative. Moreover, the above explicit computation of operadic colimits shows easily that the unit $\id \isoto GF$ is an isomorphism. Altogether this implies that $G$ and $F$ are quasi-inverse equivalences, as desired.
\end{proof}

Next up we show a weak version of the enriched Yoneda lemma. This allows a much better understanding of enriched functor categories in many situations. In the following, note that every closed monoidal $\infty$-category is automatically enriched over itself.

\begin{lemma} \label{rslt:weak-enriched-yoneda}
Let $\mathcal V$ be a closed monoidal $\infty$-category and let $\mathcal C$ be a $\mathcal V$-enriched $\infty$-category. Then for every $X \in \mathcal C$ the evaluation-at-$X$ functor
\begin{align*}
	\ev_X\colon \Fun_{\mathcal V}(\mathcal C, \mathcal V) \to \mathcal V, \qquad f \mapsto f(X)
\end{align*}
admits a left adjoint $\ell_X\colon \mathcal V \to \Fun_{\mathcal V}(\mathcal C, \mathcal V)$ such that for every $V \in \mathcal V$, the underlying functor of $\ell_X(V)$ is $\IHom_{\mathcal C}(X, -) \tensor V$. In particular, for every $f \in \Fun_{\mathcal V}(\mathcal C, \mathcal V)$ and every $V \in \mathcal V$ there is a natural isomorphism of anima
\begin{align*}
	\Hom(\IHom_{\mathcal C}(X, -) \tensor V, f) = \Hom(V, f(X)).
\end{align*}
\end{lemma}
\begin{proof}
Let $p\colon \mathcal C^\tensor \to \opLM$ and $q\colon \mathcal V^{\tensor\tensor} \to \opLM$ denote the fibrations of $\infty$-categories inducing the enrichments of $\mathcal C$ and $\mathcal V$ over $\mathcal V$. As in the proof of \cref{rslt:univ-prop-of-trivial-enriched-category} the object $X \in \mathcal C$ corresponds to an $\opLM$-operad map $i_X\colon \mathcal T^\tensor \to \mathcal C$ (where $\mathcal T^\tensor = \mathcal V \boxplus \opTriv$) and the functor $\ev_X$ can be written as the functor
\begin{align*}
	\ev_X\colon \Alg'_{\mathcal C^\tensor/\opLM}(\mathcal V^{\tensor\tensor}) \to \Alg'_{\mathcal T^\tensor/\opLM}(\mathcal V^{\tensor\tensor})
\end{align*}
given by precomposing with $i_X$. Here $\Alg'$ denotes the subcategory of those functors which are the identity on $\mathcal V^\tensor$ (the fiber over $\opAssoc \subset \opLM$) and with non-nontrivial morphisms on the restriction to $\mathcal V^\tensor$ (as in the proof of \cref{rslt:univ-prop-of-trivial-enriched-category}). We can therefore apply the operadic Kan extension construction \cite[Corollary 3.1.3.4]{lurie-higher-algebra}, which reduces the existence of the left adjoint $\ell_X$ to the following problem: Given an object $V \in \mathcal V$, viewed as an $\opLM$-operad map $i_V\colon \mathcal T^\tensor \to \mathcal V^{\tensor\tensor}$, and an object $Y \in \mathcal C$, the induced diagram $i_V\colon (\mathcal T^\tensor_\act)_{/Y} \to \mathcal V^{\tensor\tensor}$ has an operadic $p$-colimit. The objects of $(\mathcal T^\tensor_\act)_{/Y}$ can be identified with active morphisms $V_1 \dsum \dots \dsum V_n \dsum X \to Y$ in $\mathcal C^\tensor$, where $n \ge 0$ and $V_1, \dots, V_n \in \mathcal V$. As in the proof of \cref{rslt:univ-prop-of-trivial-enriched-category} one shows that the full subcategory $\mathcal I \subset (\mathcal T^\tensor_\act)_{/Y}$ spanned by objects of $\Hom(\{ V_1 \} \tensor X, Y)$ (where $V_1 \in \mathcal V$) is cofinal, hence we only need to compute an operadic $p$-colimit over $\mathcal I$. But now $\mathcal I$ can be identified with $\mathcal V_{/\IHom_{\mathcal C}(X, Y)}$, which has the final object $\IHom_{\mathcal C}(X, Y)$. We are thus reduced to finding an operadic $p$-colimit of the diagram $\{ * \} \to \mathcal V^{\tensor\tensor}$, $* \mapsto \IHom_{\mathcal C}(X, Y) \dsum V$ in $\mathcal V$. Employing \cite[Lemma 3.1.1.14]{lurie-higher-algebra} such an operaidc $p$-colimit is given by the $p$-coCartesian edge $\IHom_{\mathcal C}(X, Y) \dsum V \to \IHom_{\mathcal C}(X, Y) \tensor V$. This shows the existence of $\ell_X$ and we also immediately get the desired description of $\ell_X(V)$.
\end{proof}

\begin{corollary} \label{rslt:weak-enriched-yoneda-is-functorial}
Let $\mathcal V$ be a closed monoidal $\infty$-category and let $\mathcal C$ be a $\mathcal V$-enriched $\infty$-category. Then there is a natural functor
\begin{align*}
	\mathcal C^\opp \cprod \mathcal V \to \Fun_{\mathcal V}(\mathcal C, \mathcal V), \qquad (X, V) \mapsto \IHom(X, -) \tensor V.
\end{align*}
Restricting this functor to $\mathcal C^\opp \cprod \{ 1_{\mathcal V} \}$ provides a fully faithful embedding
\begin{align*}
	\mathcal C^\opp \injto \Fun_{\mathcal V}(\mathcal C, \mathcal V).
\end{align*}
\end{corollary}
\begin{proof}
There is a natural functor
\begin{align*}
	\alpha\colon \mathcal C \cprod \mathcal V^\opp \cprod \Fun_{\mathcal V}(\mathcal C, \mathcal V) \to \Ani, \qquad (X, V, f) \mapsto \Hom(V, f(X)).
\end{align*}
Let us abbreviate $\mathcal F := \Fun_{\mathcal V}(\mathcal C, \mathcal V)$. Then we can equivalently see $\alpha$ as a functor $\mathcal C \cprod \mathcal V^\opp \to \mathcal P(\mathcal F^\opp)$, where $\mathcal P(-) = \Fun((-)^\opp, \Ani)$ denotes the presheaf $\infty$-category. By \cref{rslt:weak-enriched-yoneda}, for every $(X, V) \in \mathcal C \cprod \mathcal V^\opp$ the object $\alpha(X, V) \in \mathcal P(\mathcal F^\opp)$ is representable by $\IHom(X, -) \tensor V \in \mathcal F^\opp$, i.e. it lies in the image of the Yoneda embedding $\mathcal F^\opp \injto \mathcal P(\mathcal F^\opp)$. Composing $\alpha$ with a quasi-inverse of the Yoneda embedding on its image we obtain a functor $\mathcal C \cprod \mathcal V^\opp \to \mathcal F^\opp$. Now apply $(-)^\opp$ to get the desired functor. The claim that $\mathcal C ^\opp \injto \Fun_{\mathcal V}(\mathcal C, \mathcal V)$ is fully faithful follows immediately from the universal property of $\IHom(X, -) \tensor V$.
\end{proof}

\begin{corollary} \label{rslt:enriched-presheaves-are-presentable}
Let $\mathcal V$ be a presentably monoidal $\infty$-category and $\mathcal C$ a small $\mathcal V$-enriched $\infty$-category. Then the $\infty$-category $\Fun_{\mathcal V}(\mathcal C, \mathcal V)$ is presentable and generated under small colimits by the objects $\IHom(X, -) \tensor V$ for $X \in \mathcal C$ and $V \in \mathcal V$.
\end{corollary}
\begin{proof}
By \cref{rslt:enriched-forgetful-functor-preserves-lim-colim} the $\infty$-category $\Fun_{\mathcal V}(\mathcal C, \mathcal V)$ admits all small colimits. Now let $\kappa$ be a regular cardinal such that $\mathcal V$ is $\kappa$-compactly generated. Then by \cref{rslt:enriched-forgetful-functor-is-conservative,rslt:weak-enriched-yoneda} the family of objects $\IHom(X, -) \tensor V \in \Fun_{\mathcal V}(\mathcal C, \mathcal V)$ for $X \in \mathcal C$ and $\kappa$-compact $V \in \mathcal V$ satisfies the assumptions of \cref{rslt:compact-generators-implies-compactly-generated}, hence $\Fun_{\mathcal V}(\mathcal C, \mathcal V)$ which implies the claim.
\end{proof}

We now want to study the behavior of the $\infty$-category $\Fun_{\mathcal V}(\mathcal C, \mathcal V)$ under a change of the monoidal $\infty$-category $\mathcal V$. In general it seems to be surprisingly subtle to ``forget'' a $\mathcal V'$-enrichment on an $\infty$-category $\mathcal C$ along a lax monoidal functor $\nu_*\colon \mathcal V' \to \mathcal V$, so we contend ourselves with the special case where $\nu_*$ admits a monoidal left adjoint. Here the following technical results comes in handy.

\begin{lemma} \label{rslt:automatic-lax-monoidal-structure-of-right-adjoint}
Let $\nu^*\colon \mathcal V \to \mathcal V'$ be a colimit-preserving (symmetric) monoidal functor of presentably (symmetric) monoidal $\infty$-categories. Then there is a lax (symmetric) monoidal functor $\nu_*\colon \mathcal V' \to \mathcal V$ which is right adjoint to $\nu^*$ and such that for every coCartesian fibration of $\infty$-operads $\mathcal O^\tensor \to \mathcal B^\tensor$, where $\mathcal B^\tensor \in \{ \opAssoc, \opComm \}$ (depending on whether the claim is about monoidal or symmetric monoidal structures) and $\mathcal O^\tensor$ is small, the functors
\begin{align*}
	(\nu^*)_*\colon \Alg_{\mathcal O/\mathcal B}(\mathcal V') \rightleftarrows \Alg_{\mathcal O/\mathcal B}(\mathcal V) \noloc (\nu_*)_*
\end{align*}
are adjoint.
\end{lemma}
\begin{proof}
The same proof as in \cite[Proposition A.5.11]{gepner-haugseng-enriched} applies.
\end{proof}

\begin{lemma} \label{rslt:forget-enrichment-properties}
Let $\nu^*\colon \mathcal V \to \mathcal V'$ be a colimit-preserving monoidal functor of presentably monoidal $\infty$-categories and let $\nu_*\colon \mathcal V' \to \mathcal V$ be a right adjoint of $\nu^*$. Let $\mathcal C$ be a $\mathcal V'$-enriched $\infty$-category.
\begin{lemenum}
	\item \label{rslt:forget-enrichment-along-presentable-V} $\mathcal C$ is naturally enriched over $\mathcal V$ with
	\begin{align*}
		\IHom_{\mathcal C, \mathcal V}(-, -) = \nu_* \comp \IHom_{\mathcal C, \mathcal V'}(-, -),
	\end{align*}
	where $\IHom_{\mathcal C, \mathcal V^{(')}}$ denotes the $\mathcal V^{(')}$-enriched $\Hom$ in $\mathcal C$.

	\item \label{rslt:connecting-map-for-forgetting-enrichment} Let $p'\colon \mathcal C^\tensor_{\mathcal V'} \to \opLM$ and $p\colon \mathcal C^\tensor_{\mathcal V} \to \opLM$ denote the fibrations of $\infty$-operads exhibiting the enrichment of $\mathcal C$ over $\mathcal V'$ and $\mathcal V$. Then there is a natural map of $\infty$-operads over $\opLM$,
	\begin{align*}
		(\nu_*, \id)\colon \mathcal C^\tensor_{\mathcal V'} \to \mathcal C^\tensor_{\mathcal V},
	\end{align*}
	which restricts to $\nu_*$ over $\mathfrak a$ and to the identity over $\mathfrak m$.
					\end{lemenum}
\end{lemma}

\begin{remark}
Our construction for \cref{rslt:forget-enrichment-properties} is rather ad-hoc, which is why we need the strong assumptions on $\mathcal V$, $\mathcal V'$ and $\nu_*$. As suggested by Lurie, there should be a much more general and elegant construction, working along any lax monoidal functor $\nu_*\colon \mathcal V' \to \mathcal V$: Choose $\mathcal C^\tensor_{\mathcal V}$ to be the pushout (in the $\infty$-category of $\infty$-operads) of $\mathcal C^\tensor_{\mathcal V'}$ and $\mathcal V^\tensor$ over $\mathcal V'^\tensor$. Using this approach, it is also easy to write down a forgetful functor $\Fun_{\mathcal V'}(\mathcal C, \mathcal D) \to \Fun_{\mathcal V}(\mathcal C, \mathcal D)$ for any $\mathcal V'$-enriched $\infty$-categories $\mathcal C$ and $\mathcal D$, which would remove all the assumptions in \cref{rslt:forgetful-functor-of-enriched-presheaves} below. On the other hand, computing colimits of $\infty$-operads is generally a non-trivial task, so we stick with our ad-hoc version.
\end{remark}

\begin{proof}[Proof of \cref{rslt:forget-enrichment-properties}]
We start with the proof of (i). Let $\mathcal V^\tensor$ and $\mathcal V'^\tensor$ denote the $\infty$-operads exhibiting the monoidal structures on $\mathcal V$ and $\mathcal V'$. Let $\mathcal D$ be any weakly $\mathcal V'$-enriched $\infty$-category, exhibited by a fibration $q'\colon \mathcal D^\tensor_{\mathcal V'} \to \opLM$, i.e. we have $(\mathcal D^\tensor_{\mathcal V'})_{\mathfrak m} = \mathcal D$ and $(\mathcal D^\tensor_{\mathcal V'})_{\mathfrak a} = \mathcal V'^\tensor$. We claim that then there is a natural fibration $q\colon \mathcal D^\tensor_{\mathcal V} \to \opLM$ with $(\mathcal D^\tensor_{\mathcal V})_{\mathfrak m} = \mathcal D$ and $(\mathcal D^\tensor_{\mathcal V})_{\mathfrak a} = \mathcal V^\tensor$ together with an $\opLM$-monoidal map of $\infty$-operads $(\nu^*, \id)\colon \mathcal D^\tensor_{\mathcal V} \to \mathcal D^\tensor_{\mathcal V'}$ which restricts to $\nu^*$ over $\mathfrak a$ and is the identity over $\mathfrak m$. Namely, first assume that $\mathcal D$ is tensored over $\mathcal V'$, i.e. $q'$ is a coCartesian fibration. Then we can view $\nu^*\colon \mathcal V \to \mathcal V'$ as a map of algebras in $\infcatinf$ and $\mathcal D$ as a $\mathcal V'$-module (see \cite[Remark 2.4.2.6]{lurie-higher-algebra}). By forgetting the $\mathcal V'$-module structure of $\mathcal D$ along $\nu^*$, we can also view $\mathcal D$ as a $\mathcal V$-module, which produces the desired (coCartesian) fibration $q$ and the map $(\nu^*, \id)$. Now if $\mathcal D$ is general, form the $\opLM$-monoidal envelope of $q'$ (see \cite[\S2.2.4]{lurie-higher-algebra}) and restrict to the full suboperad spanned by all objects over $\mathfrak m$ but only those in $\mathcal V'$ over $\mathfrak a$. This way we obtain a $\mathcal V'$-tensored $\infty$-category $\mathcal D'$ such that $\mathcal D \subset \mathcal D'$ compatibly with the weakly enriched structure. Now the desired construction of $q$ and $(\nu^*, \id)$ reduces easily to $\mathcal D'$, which was handled before.

Applying the above construction to the $\mathcal V'$-enriched $\infty$-category $\mathcal D = \mathcal C$ we obtain a weak enrichment of $\mathcal C$ over $\mathcal V$. One checks easily that this weak enrichment is actually an enrichment and that the claimed formula for the enriched $\Hom$'s holds; this proves (i).

We now prove (ii), so let $p'$ and $p$ be given as in the claim. Then by the construction in (i), there is an $\opLM$-monoidal map $(\nu^*, \id)\colon \mathcal C^\tensor_{\mathcal V} \to \mathcal C^\tensor_{\mathcal V'}$ which restricts to $\nu^*$ over $\mathfrak a$ and to the identity over $\mathfrak m$. Starting with this functor $(\nu^*, \id)$, the construction of the desired functor $(\nu_*, \id)$ can be done similar to the proof of \cite[Proposition A.5.11]{gepner-haugseng-enriched}, as follows. Let $\mathcal O^\tensor \to \opLM$ be a map of $\infty$-operads such that $\mathcal O_{\mathfrak m} = \mathcal C$. We denote
\begin{align*}
	\Alg'_{\mathcal O^\tensor/\opLM}(\mathcal C_{\mathcal V}^\tensor) \subset \Alg_{\mathcal O^\tensor/\opLM}(\mathcal C_{\mathcal V}^\tensor)
\end{align*}
the subcategory spanned by those maps $\mathcal O^\tensor \to \mathcal C_{\mathcal V}^\tensor$ which restrict to the identity on $\mathcal C$, and where we only allow morphisms of functors that are trivial over $\mathfrak m \in \opLM$. Now suppose that $\mathcal O_{\mathfrak a}$ is small. We claim that then $\Alg'_{\mathcal O^\tensor/\opLM}(\mathcal C_{\mathcal V}^\tensor)$ is presentable. Namely, let $\mathcal O^\tensor_{triv} := \mathcal O_{\mathfrak a,triv} \bsum \mathcal C_{triv}$, where on the right-hand side the subscript ``$triv$'' means the trivial $\infty$-operad with underlying $\infty$-category the given $\infty$-category. Then $\mathcal O^\tensor_{triv}$ has a natural map to $\opLM$ and we have
\begin{align*}
	\Alg'_{\mathcal O^\tensor_{triv}/\opLM}(\mathcal C_{\mathcal V}^\tensor) = \Fun(\mathcal O_{\mathfrak a}, \mathcal V),
\end{align*}
which is presentable. Moreover, it follows easily from the theory of operadic left Kan extensions (see \cite[Corollary 3.1.3.4]{lurie-higher-algebra}) that there is an adjoint pair of functors
\begin{align*}
	F\colon \Alg'_{\mathcal O^\tensor_{triv}/\opLM}(\mathcal C_{\mathcal V}^\tensor) \rightleftarrows \Alg'_{\mathcal O^\tensor/\opLM}(\mathcal C_{\mathcal V}^\tensor) \noloc G,
\end{align*}
where $G$ is the forgetful functor. By the same argument as in the proof of \cref{rslt:enriched-forgetful-functor-preserves-lim-colim} one shows that $\Alg'_{\mathcal O^\tensor/\opLM}(\mathcal C_{\mathcal V}^\tensor)$ admits all small colimits and that $G$ preserves small colimits. Moreover, $G$ is easily seen to be conservative (after all, $\mathcal O^\tensor_{triv}$ and $\mathcal O^\tensor$ have the same objects). Hence the adjunction of $F$ and $G$ is monadic by \cite[Theorem 4.7.3.5]{lurie-higher-algebra}, so we can apply \cite[Proposition A.5.9]{gepner-haugseng-enriched} to conclude that $\Alg'_{\mathcal O^\tensor/\opLM}(\mathcal C_{\mathcal V}^\tensor)$ is indeed presentable.

Back to the construction of $(\nu_*,\id)$, the adjoint functor theorem now provides us with a natural right adjoint (the fact that $(\nu^*,\id)_*$ preserves all colimits follows from the fact that $(\nu^*,\id)$ is $\opLM$-monoidal, i.e. sends coCartesian edges to coCartesian edges; cf. the proof of \cite[Proposition A.5.10]{gepner-haugseng-enriched})
\begin{align*}
	(\nu^*,\id)_*\colon \Alg'_{\mathcal O^\tensor/\opLM}(\mathcal C_{\mathcal V}^\tensor) \rightleftarrows \Alg'_{\mathcal O^\tensor/\opLM}(\mathcal C_{\mathcal V'}^\tensor)\noloc G_{\mathcal O^\tensor}
\end{align*}
for every $\mathcal O^\tensor \to \opLM$ as above. For all large enough regular cardinals $\kappa$, take $\mathcal O^\tensor$ to be the full suboperad of $\mathcal C_{\mathcal V'}^\tensor$ spanned by $\mathcal C$ and $\mathcal V'^\kappa$. Plugging the natural inclusion $\mathcal O^\tensor \injto \mathcal C_{\mathcal V'}^\tensor$ into $G_{\mathcal O^\tensor}$ produces a map $\mathcal O^\tensor \to \mathcal C_{\mathcal V}^\tensor$. This is natural in $\mathcal O^\tensor$, so by taking the union over all $\kappa$ we get a map of $\infty$-operads $\alpha\colon \mathcal C_{\mathcal V'}^\tensor \to \mathcal C_{\mathcal V}^\tensor$. By construction $\alpha$ is the identity on $\mathcal C$. Moreover, plugging $\mathcal O^\tensor = \opTriv \bsum\, \mathcal C_{triv}$ into the above adjunction (and using the naturality of $G_{\mathcal O^\tensor}$) easily shows that the restriction of $\alpha$ to $\mathcal V'$ is right adjoint to $\nu^*$. Hence $(\nu_*, \id) = \alpha$ has the desired form.
\end{proof}

\begin{corollary} \label{rslt:forgetful-functor-of-enriched-presheaves}
Let $\nu^*\colon \mathcal V \to \mathcal V'$ be a colimit-preserving monoidal functor of presentable closed monoidal $\infty$-categories and let $\nu_*\colon \mathcal V' \to \mathcal V$ be a right adjoint of $\nu^*$. Let $\mathcal C$ be a $\mathcal V'$-enriched $\infty$-category. Then $\mathcal C$ is naturally $\mathcal V$-enriched via $\nu_*$ and there is a natural functor
\begin{align*}
	\theta\colon \Fun_{\mathcal V'}(\mathcal C, \mathcal V') \to \Fun_{\mathcal V}(\mathcal C, \mathcal V)
\end{align*}
which acts as the composition with $\nu_*$ on the underlying functor categories. The functor $\theta$ preserves small limits. If $\nu_*$ is conservative or preserves small colimits, then the same holds for $\theta$.
\end{corollary}
\begin{proof}
Let $\mathcal V^\tensor$ and $\mathcal V'^\tensor$ denote the $\infty$-operads exhibiting the monoidal structures on $\mathcal V$ and $\mathcal V'$. By \cref{rslt:automatic-lax-monoidal-structure-of-right-adjoint} $\nu_*$ can naturally be seen as a map of $\infty$-operads $\mathcal V'^\tensor \to \mathcal V^\tensor$ over $\opAssoc$. Taking the pullback along the forgetful functor $\opLM \to \opAssoc$ produces a map of $\infty$-operads
\begin{align*}
	\nu_*^\tensor\colon \mathcal V'^{\tensor\tensor} \to \mathcal V^{\tensor\tensor},
\end{align*}
where $\mathcal V'^{\tensor\tensor} \to \opLM$ and $\mathcal V^{\tensor\tensor} \to \opLM$ are the fibrations of $\infty$-operads exhibiting the self-enrichments of $\mathcal V'$ and $\mathcal V$. Denote $p'\colon \mathcal C^\tensor_{\mathcal V'} \to \opLM$ and $p\colon \mathcal C^\tensor_{\mathcal V} \to \opLM$ the fibrations of $\infty$-operads exhibiting the enrichment of $\mathcal C$ over $\mathcal V'$ and $\mathcal V$ (where the latter enrichment comes from \cref{rslt:forget-enrichment-along-presentable-V}). By \cref{rslt:connecting-map-for-forgetting-enrichment} there is a natural map of $\infty$-operads $(\nu_*, \id)\colon \mathcal C^\tensor_{\mathcal V'} \to \mathcal C^\tensor_{\mathcal V}$ over $\opLM$ which restricts to $\nu_*$ over $\mathfrak a$ and to the identity over $\mathfrak m$. We obtain natural functors
\begin{align*}
	(\nu_*^\tensor)_*\colon \Fun_{\mathcal V'}(\mathcal C, \mathcal V') \to \Alg^{\nu_*}_{\mathcal C^\tensor_{\mathcal V'}/\opLM}(\mathcal V^{\tensor\tensor}) \from \Fun_{\mathcal V}(\mathcal C, \mathcal V) \noloc (\nu_*, \id)^*,
\end{align*}
where $\Alg^{\nu_*}$ denotes the subcategory of $\Alg$ which consists of those functors that restrict to $\nu_*$ over $\mathfrak a$ and where morphisms of functors are required to be the identity over $\mathfrak a$. To construct the desired functor $\theta$, it is enough to show that the functor $(\nu_*, \id)^*$ is an equivalence. Writing $\mathcal C$ as a filtered colimit of small full subcategories we can reduce to the case that $\mathcal C$ is small. In this case \cref{rslt:enriched-presheaves-are-presentable} implies that the enriched functor category $\Fun_{\mathcal V}(\mathcal C, \mathcal V)$ is presentable, with generators given by $\IHom_{\mathcal C, \mathcal V}(X, -) \tensor V$, where $X \in \mathcal C$ and $V \in \mathcal V$. The same arguments as in the proof of these facts apply to $\Alg^{\nu_*}_{\mathcal C^\tensor_{\mathcal V'}/\opLM}(\mathcal V^{\tensor\tensor})$, showing that this $\infty$-category is also presentable, with generators also given by $\IHom_{\mathcal C,\mathcal V}(X, -) \tensor V$ (where $X$ and $V$ are as above).

We now show that $(\nu_*, \id)^*$ is an equivalence. We first show that it is fully faithful, so let $f, g \in \Fun_{\mathcal V}(\mathcal C, \mathcal V)$ be given. By writing $f$ as an iterated small colimit of the generators and using that $(\nu_*, \id)^*$ preserves small colimits (by \cref{rslt:enriched-forgetful-functor-properties}), we can w.l.o.g. assume that $f = \IHom_{\mathcal C, \mathcal V}(X, -) \tensor V$ for some $X \in \mathcal C$ and some $V \in \mathcal V$. But we have
\begin{align*}
	(\nu_*,\id)^* \IHom_{\mathcal C, \mathcal V}(X, -) \tensor V = \IHom_{\mathcal C, \mathcal V}(X, -) \tensor V,
\end{align*}
so that we get
\begin{align*}
	\Hom((\nu_*,\id)^* f, (\nu_*,\id)^* g) = \Hom(V, g(X)) = \Hom(f, g)
\end{align*}
by \cref{rslt:weak-enriched-yoneda}, proving full faithfulness of $(\nu_*,\id)^*$. It is now easy to see that $(\nu_*,\id)^*$ is an equivalence: The image of $(\nu_*,\id)^*$ is closed under colimits and contains the generators $\IHom_{\mathcal C, \mathcal V}(X, -) \tensor V$, thus $(\nu_*,\id)^*$ is essentially surjective. This finishes the construction of the functor $\theta$.

If $\nu_*$ is conservative or preserves small colimits, then it follows immediately from \cref{rslt:enriched-forgetful-functor-properties} that the same is true for $\theta$, as it can be checked on the underlying functor categories (it follows similarly that $\theta$ preserves all small limits because $\nu_*$ does).
\end{proof}

We will now come to the main result of this subsection, which shows an invariance of $\Fun_{\mathcal V}(\mathcal C, \mathcal V)$ under $\mathcal V$ in nice cases.

\begin{proposition} \label{rslt:invariance-of-enriched-presheaves}
Let $\nu^*\colon \mathcal V \to \mathcal V'$ be a colimit-preserving monoidal functor of presentable closed monoidal $\infty$-categories. Assume that the right adjoint $\nu_*\colon \mathcal V' \to \mathcal V$ of $\nu^*$ satisfies the following properties:
\begin{enumerate}[(a)]
	\item $\nu_*$ preserves all small colimits.
	\item $\nu_*$ is conservative.
	\item For all $V \in \mathcal V$ and $V' \in \mathcal V'$ the natural morphism
	\begin{align*}
		\nu_*(V') \tensor V \isoto \nu_*(V' \tensor \nu^*(V))
	\end{align*}
	is an isomorphism.
\end{enumerate}
Let $\mathcal C$ be a $\mathcal V'$-enriched $\infty$-category. Then $\mathcal C$ is naturally $\mathcal V$-enriched via $\nu_*$ and the functor
\begin{align*}
	\Fun_{\mathcal V'}(\mathcal C, \mathcal V') \isoto \Fun_{\mathcal V}(\mathcal C, \mathcal V).
\end{align*}
is an equivalence.
\end{proposition}
\begin{proof}
The functor $\theta\colon \Fun_{\mathcal V'}(\mathcal C, \mathcal V') \to \Fun_{\mathcal V}(\mathcal C, \mathcal V)$ is the one from \cref{rslt:forgetful-functor-of-enriched-presheaves}. To show that it is an equivalence, we can assume that $\mathcal C$ is small (in the general case write $\mathcal C$ as a filtered colimit of its small subcategories). In this case \cref{rslt:enriched-presheaves-are-presentable} implies that both $\Fun_{\mathcal V'}(\mathcal C, \mathcal V')$ and $\Fun_{\mathcal V}(\mathcal C, \mathcal V)$ are presentable, with generators given by $\IHom_{\mathcal C, \mathcal V'}(X, -) \tensor V'$ and $\IHom_{\mathcal C, \mathcal V}(X, -) \tensor V$ respectively, where $X \in \mathcal C$, $V' \in \mathcal V'$ and $V \in \mathcal V$.

We now show that $\theta$ is an equivalence. We first show that it is fully faithful, so let $f, g \in \Fun_{\mathcal V'}(\mathcal C, \mathcal V')$ be given. Since $\nu_*$ preserves small colimits by assumption, so does $(\nu_*^\tensor)_*$, hence we can w.l.o.g. assume that $f = \IHom_{\mathcal C, \mathcal V'}(X, -) \tensor V'$ for some $X \in \mathcal C$ and $V' \in \mathcal V'$. Since $\nu_*$ is conservative, $\mathcal V'$ is generated under small colimits by the image of $\nu^*$ (apply \cref{rslt:compact-generators-implies-compactly-generated} to the family of functors $\Hom(\nu^* V_i, -)\colon \mathcal V' \to \Ani$ for a generating family $(V_i)_i$ of $\mathcal V$). We can thus reduce to the case $V' = \nu^* V$ for some $V \in \mathcal V$. In this case, using the assumed projection formula for $\nu_*$ we get
\begin{align*}
	\theta(f) = \theta \big(\IHom_{\mathcal C, \mathcal V'}(X, -) \tensor \nu^* V \big) = \IHom_{\mathcal C, \mathcal V}(X, -) \tensor V.
\end{align*}
(This identity can be checked on the underlying functors $\mathcal C \to \mathcal V$ because the forgetful functor from $\mathcal V$-enriched functors to the underlying functors is conservative by \cref{rslt:enriched-forgetful-functor-is-conservative}.) In particular we have
\begin{align*}
	\Hom(\theta(f), \theta(g)) = \Hom(V, \nu_* g(X)) = \Hom(\nu^* V, g(X)) = \Hom(f, g),
\end{align*}
proving full faithfulness of $\theta$. Essential surjectivity of $\theta$ follows from the fact that the image is closed under small colimits and contains all generators by the above computation.
\end{proof}

In the main part of this thesis we are mainly interested in the $\infty$-category of enriched endofunctors of a given closed symmetric monoidal $\infty$-category $\mathcal V$ viewed as enriched over itself. These endofunctors come equipped with the composition monoidal structure, which we now show to be functorial in $\mathcal V$. In the following results, $\infcatinf^\ocircle$ denotes the $\infty$-category of monoidal $\infty$-categories.

\begin{lemma} \label{rslt:functoriality-of-composition-monoidal-structure}
Let $S$ be an $\infty$-category and $f\colon S \to \infcatinf$ a functor. Then there is a natural fibration of generalized $\infty$-operads $\mathcal E(S) \to S \cprod \opAssoc$ such that for every $s \in S$ the fiber $\mathcal E(S)_s \to \opAssoc$ is equivalent to the $\infty$-category $\Fun(f(s), f(s))$ equipped with the composition monoidal structure.
\end{lemma}
\begin{proof}
Let $\mathcal T \to S$ be the coCartesian fibration which classifies $f$. View $\mathcal T$ as an object of $(\infcatinf)_{/S}$ and view the latter $\infty$-category as enriched over itself via the cartesian monoidal structure. Consider the associated $\infty$-category $(\infcatinf)_{/S}[\mathcal T]$ as defined in \cite[Definition 4.7.1.1]{lurie-higher-algebra}. Roughly speaking it consists of pairs $(\mathcal T', \eta)$, where $\mathcal T'$ is an $\infty$-category over $S$ and $\eta\colon \mathcal T' \cprod_S \mathcal T \to \mathcal T$ is a functor over $S$. Then $(\infcatinf)_{/S}[\mathcal T]$ contains a final object $\End_S(\mathcal T) \to S$ such that for every simplicial set $K$ over $S$ we have
\begin{align}
	\Hom_S(K, \End_S(\mathcal T)) = \Hom_S(K \cprod_S \mathcal T, \mathcal T). \label{eq:def-of-End-S}
\end{align}
By \cite[Corollary 4.7.1.40]{lurie-higher-algebra} (and \cite[Proposition 4.1.3.19]{lurie-higher-algebra}) we can view $\End_S(\mathcal T)$ as an associative algebra in $(\infcatinf)_{/S}$, given by a functor $\varphi\colon \opAssoc \to (\infcatinf)_{/S}$. Applying the unstraightening functor we obtain a map $\alpha\colon \End_S(\mathcal T)^\ocircle \to S \cprod \opAssoc$ of coCartesian fibrations over $\opAssoc$. We claim that $\alpha$ is a categorical fibration. Namely, since $\opAssoc$ is a 1-category, the unstraightening can be computed by the relative nerve construction (see \cite[Proposition 3.2.5.21]{lurie-higher-topos-theory}), so by \cite[Lemma 3.2.5.11.(3)]{lurie-higher-topos-theory} it is enough to show that for every $x \in \opAssoc$ the map $\varphi(x) \to S$ is a categorical fibration. Now $\phi(x)$ is equivalent to a product of copies of $\End_S(\mathcal T)$ inside $(\infcatinf)_{/S}$, so it is enough to show that $\End_S(\mathcal T) \to S$ is a categorical fibration. To see this, let $\mathfrak P = (M_S, T_S, \emptyset)$ be the categorical pattern on $S$ (see \cite[Definition B.0.19]{lurie-higher-algebra}) where $M_S$ consists of all equivalences in $S$ and $T_S$ consists of all 2-simplices in $S$. Let furthermore $M$ be the class of all equivalences in $\mathcal T$. Then $(\mathcal T, M) \to (S, M_S)$ is easily seen to be $\mathfrak P$-fibered because $\mathcal T \to S$ is a categorical fibration. Moreover, the diagram $(S, M_S) \from (\mathcal T, M) \to (S, M_S)$ (where both maps are the projection $\mathcal T \to S$) satisfies the conditions of \cite[Theorem B.4.2]{lurie-higher-algebra}: Condition (1) follows from the fact that $\mathcal T \to S$ is a coCartesian fibration (see \cite[Example B.3.11]{lurie-higher-algebra}), condition (4) follows from the fact that $\mathcal T \to S$ is a categorical fibration (cf. \cite[Proposition 3.3.1.8]{lurie-higher-algebra}) and the other conditions are obvious. It follows that the functor $F\colon (\catsimpmark)_{/\mathfrak P} \to (\catsimpmark)_{/\mathfrak P}$ given by $F(K) = K \cprod_S \mathcal T$ is a left Quillen functor. It therefore admits a right adjoint $G$, which satisfies $G(\mathcal T) = \End_S(\mathcal T)$. Being a right Quillen functor, $G$ preserves fibrant objects and hence $\End_S(\mathcal T) \to S$ is $\mathfrak P$-fibered; in particular it is a categorical fibration, as desired.

We have now constructed a categorical fibration $\End_S(\mathcal T)^\ocircle \to S \cprod \opAssoc$. Using the fact that the projection $\End_S(\mathcal T)^\ocircle \to \opAssoc$ is coCartesian, it is easy to see that $\End_S(\mathcal T)^\ocircle$ is a generalized $\infty$-operad. We can thus take $\mathcal E(S) = \End_S(\mathcal T)^\ocircle$.
\end{proof}

\begin{proposition} \label{rslt:functoriality-of-enriched-endofunctors}
Let $\mathcal V^\tensor$ be a presentably symmetric monoidal $\infty$-category. Then there is a natural functor
\begin{align*}
	\CAlg(\mathcal V) \to \infcatinf^\ocircle, \qquad A \mapsto \End(\Mod_A(\mathcal V)),
\end{align*}
where for every closed symmetric monoidal $\infty$-category $\mathcal D^\tensor$ we denote $\End(\mathcal D) := \Fun_{\mathcal D}(\mathcal D, \mathcal D)$ equipped with the composition monoidal structure. This functor sends a morphism $f\colon A \to B$ in $\CAlg(\mathcal V)$ to the functor
\begin{align*}
	f^\natural\colon \End(\Mod_A(\mathcal V)) \xto{f^* \comp -} \Fun_{\Mod_A(\mathcal V)}(\Mod_B(\mathcal V), \Mod_A(\mathcal V)) = \End(\Mod_B(\mathcal V)),
\end{align*}
where the equivalence on the right is the one from \cref{rslt:invariance-of-enriched-presheaves}.
\end{proposition}
\begin{proof}
It follows from \cite[Theorem 4.5.3.1]{lurie-higher-algebra} that there is a coCartesian fibration $\mathcal T \to \CAlg(\mathcal V)$ whose fiber $\mathcal T_A$ over every $A \in \CAlg(\mathcal V)$ is the $\infty$-category underlying the $\infty$-operad which exhibits $\Mod_A(\mathcal V)$ as enriched over itself. Applying \cref{rslt:functoriality-of-composition-monoidal-structure} to the associated functor $\CAlg(\mathcal V) \to \infcatinf$ yields a fibration of generalized $\infty$-operads $\mathcal E'^\ocircle \to \CAlg(\mathcal V) \cprod \opAssoc$ whose fiber over every $A \in \CAlg(C)$ is the monoidal $\infty$-category $\mathcal E'^\ocircle_A \to \opAssoc$ with $\mathcal E'_A = \Fun(\mathcal T_A, \mathcal T_A)$. Let $\mathcal E^\ocircle \subset \mathcal E'^\ocircle$ be the subcategory spanned by the enriched functors (and morphisms between enriched functors) inside each $\Fun(\mathcal T_A, \mathcal T_A)$. We claim that the map $p\colon \mathcal E' \to \CAlg(\mathcal A) \cprod \opAssoc$ is a coCartesian fibration. We already know that $p$ is a categorical fibration (and in particular an inner fibration). As a first step we show that $p$ is a locally coCartesian fibration, so we need to verify that for every morphism $e\colon (A, x) \to (B, y)$ in $\CAlg(\mathcal V) \cprod \opAssoc$ and every object $F \in \mathcal E^\ocircle$ lying over $(A, x)$ there is a locally $p$-coCartesian edge $F \to G$ in $\mathcal E^\ocircle$ such that $p(G) = (B, y)$. Since $\mathcal E^\ocircle$ is still a generalized $\infty$-operad and thus admits (coCartesian) inert morphisms, we can reduce to the case $y = \langle 1 \rangle$, $x = \langle n \rangle$ with $x \to y$ being active, w.l.o.g. corresponding to the natural ordering on $\langle n \rangle^\circ$. In this case the fiber $\mathcal E^\ocircle_e$ of $\mathcal E^\ocircle$ over the edge $e$ looks as follows:
\begin{itemize}
	\item The fiber of $\mathcal E^\ocircle_e$ over $(A, x)$ is the $\infty$-category $\End(\Mod_A(\mathcal V))^n$. The fiber of $\mathcal E^\ocircle_e$ over $(B, y)$ is the $\infty$-category $\End(\Mod_B(\mathcal V))$.

	\item Suppose we are given objects $(F_1, \dots, F_n) \in \mathcal E^\ocircle_e$ in the fiber over $(A, x)$ and $G \in \mathcal E^\ocircle_e$ in the fiber over $(B, y)$. The map $f\colon A \to B$ induces a map $f^*\colon \mathcal T_A \to \mathcal T_B$. Then $\Hom_{\mathcal E^\ocircle_e}((F_1, \dots, F_n), G)$ is equal to the subspace of $\Hom(f^* \comp F_1 \comp \dots F_n, G \comp f^*)$ of those natural transformations which are trivial over the associative parts of $\mathcal T_A, \mathcal T_B \to \opLM$.
\end{itemize}
Viewing $\Mod_B(\mathcal V)$ also as enriched over $\Mod_A(\mathcal V)$ together with an enriched morphism $f^*\colon \Mod_A(\mathcal V) \to \Mod_B(\mathcal V)$ we obtain using \cref{rslt:invariance-of-enriched-presheaves} the natural functor
\begin{align*}
	f^\natural\colon \End(\Mod_A(\mathcal V)) \xto{f^* \comp -} \Fun_{\Mod_A(\mathcal V)}(\Mod_B(\mathcal V), \Mod_A(\mathcal V)) = \End(\Mod_B(\mathcal V)).
\end{align*}
Then $\mathcal E^\ocircle_e$ classifies the functor
\begin{align*}
	\End(\Mod_A(\mathcal V))^n \to \End(\Mod_B(\mathcal V)), \qquad (F_1, \dots, F_n) \mapsto f^\natural (F_1 \comp \dots \comp F_n)
\end{align*}
and is in particular a coCartesian fibration. This finishes the proof that $p\colon \mathcal E^\ocircle \to \CAlg(\mathcal V) \cprod \opAssoc$ is a locally coCartesian fibration.

To prove that $p$ is a coCartesian fibration, we invoke \cite[Remark 2.4.2.9]{lurie-higher-topos-theory}, which reduces the claim to showing that the above defined functor $f^\natural$ commutes with the monoidal structures in $\End(\Mod_A(\mathcal V))$ and $\End(\Mod_B(\mathcal V))$. In other words, we need to see that for every map $f\colon A \to B$ in $\CAlg(\mathcal V)$ and all $F_1, F_2 \in \End(\Mod_A(\mathcal V))$ the natural morphism $f^\natural(F_1 \comp F_2) \isoto f^\natural F_1 \comp f^\natural F_2$ in $\End(\Mod_B(\mathcal V))$ (induced by the locally coCartesian fibration $p$) is an isomorphism. This can be checked on the underlying functors and in particular after plugging in any $N \in \Mod_B(\mathcal V)$. One checks easily that $f_* (f^\natural F)(N) = F(f_* N)$, where $f_*\colon \Mod_B(\mathcal V) \to \Mod_A(\mathcal V)$ denotes the forgetful functor. Hence
\begin{align*}
	&f_* (f^\natural F_1 \comp f^\natural F_2)(N) = f_* (f^\natural F_1)((f^\natural F_2)(N)) = F_1(f_* (f^\natural F_2)(N)) = F_1(F_2(f_* N)) =\\&\qquad= f_* (f^\natural(F_1 \comp F_2))(N).
\end{align*}
Since $f_*$ is conservative, we arrive at $f^\natural(F_1 \comp F_2) = f^\natural F_1 \comp f^\natural F_2$, as desired. This finishes the proof that $p\colon \mathcal E \to \CAlg(\mathcal V) \cprod \opAssoc$ is a coCartesian fibration. Therefore $p$ defines a coCartesian family of monoidal $\infty$-categories (see \cite[Definition 4.8.3.1]{lurie-higher-algebra}). By \cite[Example 4.8.3.3]{lurie-higher-algebra} $p$ classifies the claimed functor $\CAlg(\mathcal V) \to \infcatinf^\ocircle$.
\end{proof}

\subsection{Abstract 6-Functor Formalisms} \label{sec:infcat.sixfun}

A central objective of this thesis is to construct a 6-functor formalism. However, in the $\infty$-categorical setting it is not obvious how to encode such a 6-functor formalism and, once one has found a suitable encoding, how to obtain all the desired higher homotopies. The goal of the present subsection is to define 6-functor formalisms in general and provide some tools for constructing them, following \cite{enhanced-six-operations}. Instead of using the notation in the reference we prefer to work with the $\infty$-categories of correspondences, which we find cleaner.

Let us start with the abstract definition of 6-functor formalisms. The goal of a 6-functor formalism on some category $\mathcal C$ of geometric objects (e.g. schemes or diamonds) is to associate to every $X \in \mathcal C$ an $\infty$-category $\D(X)$ and define the six functors $\tensor$, $\IHom$, $f_*$, $f^*$, $f_!$ and $f^!$ in this setting (here the latter four functors are associated to a map $f\colon Y \to X$, where the shriek functors may only be defined for a special class $E$ of maps $f$). Moreover, these six functors should satisfy several compatibilities: the last four functors are functorial in the map $f$, every second functor is right adjoint to the one before it, we have proper base-change (i.e. unrestricted base-change for $f_!$) and the projection formula for $f_!$ holds. It is possible to encode all of these data in a single lax symmetric monoidal functor $\Corr(\mathcal C)_{E,all} \to \infcatinf$ for a certain ``$\infty$-category of correspondences $\Corr(\mathcal C)_{E,all}$''. Our first goal is to introduce the relevant terminology; the following definitions are based on \cite[Definition 6.1.1]{enhanced-six-operations}.

\begin{definition}
A \emph{geometric setup} is a pair $(\mathcal C, E)$, where $\mathcal C$ is an $\infty$-category and $E$ is a collection of homotopy classes of edges in $\mathcal C$ satisfying the following properties:
\begin{enumerate}[(i)]
	\item $E$ contains all isomorphisms and is stable under composition.

	\item Pullbacks of $E$ exist and remain in $E$.
\end{enumerate}
\end{definition}

\begin{definition}
\begin{defenum}
	\item For $n \ge 0$ let $C(\Delta^n)$ be the full subcategory of $\Delta^n \cprod (\Delta^n)^\opp$ spanned by the pairs $([i], [j])$ with $i \le j$. An edge of $C(\Delta^n)$ is called \emph{vertical} (resp. \emph{horizontal}) if its projection to the second (resp. first) factor is degenerate. A square of $C(\Delta^n)$ is called \emph{exact} if it is both a pullback square and a pushout square. One can extend $C$ to a colimit-preserving endofunctor of the category of simplicial sets. Its right adjoint is denoted $B$, so that for every simplicial set $K$ we have $B(K)_n = \Hom(C(\Delta^n), K)$.

	\item \label{def:category-of-correspondences} Let $(\mathcal C, E)$ be a geometric setup. We let
	\begin{align*}
		\Corr(\mathcal C)_{E,all} \subset B(\mathcal C)
	\end{align*}
	be the simplicial subset whose $n$-cells are given by those maps $C(\Delta^n) \to \mathcal C$ which send vertical edges to $E$ and exact squares to pullback squares. By \cite[Lemma 6.1.2]{enhanced-six-operations} $\Corr(\mathcal C)_{E,all}$ is an $\infty$-category, called the \emph{$\infty$-category of correspondences} associated to $(\mathcal C, E)$.
\end{defenum}
\end{definition}

Given a geometric setup $(\mathcal C, E)$ the $\infty$-category $\Corr(\mathcal C)_{E,all}$ can be described as follows: The objects of $\Corr(\mathcal C)_{E,all}$ are the objects of $\mathcal C$. A morphism $Y \to X$ in $\Corr(\mathcal C)_{E,all}$ is given by a diagram
\begin{center}\begin{tikzcd}
	Y & \arrow[l,"h",swap] Y' \arrow[d,"v"]\\
	& X
\end{tikzcd}\end{center}
in $\mathcal C$, where $v$ lies in $E$. Given another morphism $Z \to Y$ in $\Corr(\mathcal C)_{E,all}$, represented by a diagram $Z \from Z' \to Y$, the composed morphism $Z \to X$ is given by the following diagram:
\begin{center}\begin{tikzcd}
	Z & \arrow[l] Z' \arrow[d] & \arrow[l] Z' \cprod_Y Y' \arrow[d]\\
	& Y & \arrow[l] Y' \arrow[d]\\
	&& X
\end{tikzcd}\end{center}
Note that usually $\Corr(\mathcal C)_{E,all}$ is not an ordinary category, even if $\mathcal C$ is. Moreover, we naturally have the subcategories
\begin{align*}
	\mathcal C_E, \mathcal C^\opp \subset \Corr(\mathcal C)_{E,all},
\end{align*}
where $\mathcal C_E \subset \mathcal C$ is the subcategory where we only allow morphisms in $E$. The first of the above inclusions is obtained via mapping $[f\colon Y \to X] \in E$ to $Y \xfrom{\id} Y \xto{f} X$ while the second inclusion is obtained via mapping $[f\colon Y \to X]$ to $X \xfrom{f} Y \xto{\id} Y$.

Additionally, there is a natural operadic structure on $\Corr(\mathcal C)_{E,all}$ which can be described as follows:

\begin{definition}
Let $\mathcal C$ be an $\infty$-category.
\begin{defenum}
	\item We denote $\mathcal C^\amalg \to \catFinAst$ the $\infty$-operad from \cite[Proposition 2.4.3.3]{lurie-higher-algebra} (if $\mathcal C$ admits finite coproducts then $\mathcal C^\amalg$ is the symmetric monoidal $\infty$-category with tensor operation being the coproduct). The fiber of $\mathcal C^\amalg$ over $\langle n \rangle \in \catFinAst$ is canonically identified with $\mathcal C^n$.

	\item Let $E$ be a collection of edges in $\mathcal C$. We define the following classes $E^+$ and $E^-$ of edges in $(\mathcal C^\opp)^{\amalg,\opp}$: Write an edge $f$ of $(\mathcal C^\opp)^{\amalg,\opp}$ lying over $\alpha\colon \langle m \rangle \to \langle n \rangle$ in the form $f\colon (Y_j)_{1\le j \le n} \to (X_i)_{1 \le i \le m}$. Then $E^+$ consists of those $f$ where the induced edge $Y_{\alpha(i)} \to X_i$ belongs to $E$ for all $i \in \alpha^{-1}(\langle n \rangle^\circ)$, while $E^- \subset E^+$ is the subset of those edges where $\alpha$ is the identity.
\end{defenum}
\end{definition}

\begin{definition} \label{def:correspondence-category-operad}
Let $(\mathcal C, E)$ be a geometric setup. We denote
\begin{align*}
	\Corr(\mathcal C)^\tensor_{E,all} := \Corr((\mathcal C^\opp)^{\amalg,\opp})_{E^-,all}.
\end{align*}
Then by the proof of \cite[Proposition 6.1.3]{enhanced-six-operations} there is a canonical functor $\Corr(\mathcal C)^\tensor_{E,all} \to \catFinAst$ making $\Corr(\mathcal C)^\tensor_{E,all}$ an $\infty$-operad with underlying $\infty$-category $\Corr(\mathcal C)_{E,all}$. In the following we will always implicitly endow $\Corr(\mathcal C)_{E,all}$ with this operadic structure.
\end{definition}

\begin{remark}
In \cref{def:correspondence-category-operad}, if we additionally assume that $\mathcal C$ admits all finite products then $\Corr(\mathcal C)_{E,all}^\tensor$ is even a symmetric monoidal $\infty$-category, i.e. the functor to $\catFinAst$ is a coCartesian fibration. This is shown in \cite[Proposition 6.1.3]{enhanced-six-operations}.
\end{remark}

Given a geometric setup $(\mathcal C, E)$ the inclusion $(\mathcal C^\opp)^\amalg \subset \Corr(\mathcal C)_{E,all}^\tensor$ is easily seen to be a map of $\infty$-operads. In particular, assuming the existence of finite products in $\mathcal C$, the ``tensor product'' of two objects $X, Y \in \Corr(\mathcal C)_{E,all}$ is just their product in $\mathcal C$.

We have now introduced the necessary notation to make the main definition of this subsection, namely that of a 6-functor formalism.

\begin{definition}
Let $(\mathcal C, E)$ be a geometric setup. A \emph{pre-6-functor formalism} on $(\mathcal C, E)$ is a map of $\infty$-operads
\begin{align*}
	\D\colon \Corr(\mathcal C)_{E,all} \to \infcatinf,
\end{align*}
where $\infcatinf$ is equipped with the product symmetric monoidal structure. Given a pre-6-functor formalism $\D$ as above we introduce the following additional notation:
\begin{defenum}
	\item By restricting $\D$ to the $\infty$-operad $(\mathcal C^\opp)^\amalg$ we obtain a lax symmetric monoidal functor $\mathcal C^\opp \to \infcatinf$, which by \cite[Theorem 2.4.3.18]{lurie-higher-algebra} corresponds to a functor $\D^*\colon \mathcal C^\opp \to \infcatinf^\tensor$. In particular for every $X \in \mathcal C$ the $\infty$-category $\D(X) = \D^*(X)$ comes equipped with a symmetric monoidal structure, which we usually denote
	\begin{align*}
		- \tensor -\colon \D(X) \cprod \D(X) \to \D(X).
	\end{align*}

	\item For every morphism $f\colon Y \to X$ in $\mathcal C$ we obtain a symmetric monoidal functor
	\begin{align*}
		f^* := \D^*(f)\colon \D(X) \to \D(Y).
	\end{align*}

	\item By restricting $\D$ to the subcategory $\mathcal C_E \subset \Corr(\mathcal C)_{E,all}$ we obtain a functor $\D_!\colon \mathcal C_E \to \infcatinf$. For every $[f\colon Y \to X] \in E$ we denote
	\begin{align*}
		f_! := \D_!(f)\colon \D(Y) \to \D(X).
	\end{align*}
\end{defenum}
\end{definition}

\begin{definition} \label{def:abstract-6-functor-formalism}
Let $(\mathcal C, E)$ be a geometric setup. A \emph{6-functor formalism} on $(\mathcal C, E)$ is a pre-6-functor formalism $\D\colon \Corr(\mathcal C)_{E,all} \to \infcatinf$ such that each symmetric monoidal $\infty$-category $\D(X)$ is closed and all the functors $f^*\colon \D(X) \to \D(Y)$ and $f_!\colon \D(Y) \to \D(X)$ admit right adjoints. Given a 6-functor formalism $\D$ we introduce the following additional notation:
\begin{defenum}
	\item The internal hom in $\D(X)$ (given as the right adjoint of the tensor operation) is denoted
	\begin{align*}
		\IHom(-, -)\colon \D(X)^\opp \cprod \D(X) \to \D(X).
	\end{align*}

	\item For every morphism $f\colon Y \to X$ in $\mathcal C$ the right adjoint of $f^*$ is denoted
	\begin{align*}
		f_*\colon \D(Y) \to \D(X).
	\end{align*}
	We denote by $\D_*\colon \mathcal C \to \infcatinf$ the associated functor $X \mapsto \D(X)$, $f \mapsto f_*$.

	\item For every $[f\colon Y \to X] \in E$ the right adjoint of $f_!$ is denoted
	\begin{align*}
		f^!\colon \D(X) \to \D(Y).
	\end{align*}
	We also denote by $\D^!\colon \mathcal C_E^\opp \to \infcatinf$ the associated functor $X \mapsto \D(X)$, $f \mapsto f^!$.
\end{defenum}
\end{definition}

\begin{proposition} \label{rslt:abstract-6-functor-formalism-properties}
Let $\D\colon \Corr(\mathcal C)_{E,all} \to \infcatinf$ be a 6-functor formalism on a geometric setup $(\mathcal C, E)$. Then the associated six functors $\tensor$, $\IHom$, $f^*$, $f_*$, $f_!$ and $f^!$ satisfy the following compatibilities:
\begin{propenum}
	\item The functors $f^*$, $f_*$, $f_!$ and $f^!$ are natural in $f$, i.e. for composable maps $f, g$ in $\mathcal C$ we have natural isomorphisms $(f \comp g)^* = g^* \comp f^*$ and $(f \comp g)_* = f_* \comp g_*$ and for composable maps $f, g \in E$ we have $(f \comp g)_! = f_! \comp g_!$ and $(f \comp g)^! = g^! \comp f^!$.

	\item For every morphism $f\colon Y \to X$ in $\mathcal C$ the functor $f^*\colon \D(X) \to \D(Y)$ is symmetric monoidal.

	\item (Proper Base Change) For every cartesian diagram
	\begin{center}\begin{tikzcd}
		Y' \arrow[r,"g'"] \arrow[d,"f'"] & Y \arrow[d,"f"]\\
		X' \arrow[r,"g"] & X
	\end{tikzcd}\end{center}
	in $\mathcal C$ with $f \in E$ we have a natural isomorphism
	\begin{align*}
		g^* f_! = f'_! g'^*
	\end{align*}
	of functors $\D(Y) \to \D(X')$.

	\item (Projection Formula) For every $[f\colon Y \to X] \in E$ and all $M \in \D(X)$ and $N \in \D(Y)$ there is a natural isomorphism
	\begin{align*}
		f_!(N \tensor f^* M) = f_! N \tensor M
	\end{align*}
	in $\D(X)$.
\end{propenum}
\end{proposition}
\begin{proof}
Parts (i) and (ii) follow from the construction (using that the formation of adjoint functors is natural, see \cite[Corollary 5.2.2.5]{lurie-higher-topos-theory}).

To prove (iii) let the diagram be given as in the claim. Consider the morphisms $h_1 = Y \xfrom{\id} Y \xto{f} X$ and $h_2 = X \xfrom{g} X' \xto{\id} X'$ in $\Corr(\mathcal C)_{E,all}$. Then $h_2 \comp h_1$ is given by $Y \xfrom{g'} Y' \xto{f'} X'$, hence
\begin{align*}
	f'_! g'^* = \D(h_2 \comp h_1) = \D(h_2) \comp \D(h_1) = g^* f_!,
\end{align*}
as desired.

To prove (iv), let $[f\colon Y \to X] \in E$ be given. One checks immediately that the following diagram in $\Corr(\mathcal C)_{E,all}^\tensor$ commutes:
\begin{center}\begin{tikzcd}
	(X, Y) \arrow[r] \arrow[d] & (X, X) \arrow[d]\\
	Y \arrow[r] & X
\end{tikzcd}\end{center}
Here the top horizontal map is given by $(X, Y) \xto{(\id,f)} (X, X)$, the left vertical map is given by $(X, Y) \xfrom{(f,\id)} Y$, the right vertical map is given by $(X, X) \xfrom{(\id, \id)} X$ and the bottom horizontal map is given by $Y \xto{f} X$. By applying $\D$ to this diagram and using that $\D$ is a map of $\infty$-operads (which we can equivalently view as a lax Cartesian structure via \cite[Proposition 2.4.1.7]{lurie-higher-algebra}) we obtain the following commuting diagram of $\infty$-categories:
\begin{center}\begin{tikzcd}
	\D(X) \cprod \D(Y) \arrow[r,"{(\id, f_!})"] \arrow[d] & \D(X) \cprod \D(X) \arrow[d]\\
	\D(Y) \arrow[r,"f_!"] & \D(X)
\end{tikzcd}\end{center}
Note that by definition of the symmetric monoidal structures, the right vertical map is just $- \tensor -$ and one sees easily that the left vertical maps is $f^* \tensor -$. Therefore the commutativity of the diagram shows the projection formula.
\end{proof}

Having established the general definition of a 6-functor formalism, we will now discuss a way to construct it, following \cite{enhanced-six-operations}. A common situation is the following: We start with an ($\infty$-)category $\mathcal C$ of geometric objects (e.g. schemes, diamonds). By some general constructions one constructs a functor
\begin{align*}
	\D\colon \mathcal C^\opp \to \infcatinf^\tensor, \qquad X \mapsto \D(X),
\end{align*}
associating to every space $X \in \mathcal C$ a symmetric monoidal $\infty$-category $\D(X)$ of ``sheaves on $X$'' and to every map $f\colon Y \to X$ of spaces a symmetric monoidal pullback functor $f^*\colon \D(X) \to \D(Y)$. The hard part is to construct the functor $f_!\colon \D(Y) \to \D(X)$ for a suitable class $E$ of morphisms $f\colon Y \to X$ in $\mathcal C$ and to fit $f_!$ in a 6-functor formalism. One can usually single out the following special cases of $f_!$:
\begin{itemize}
	\item If $f$ belongs to a class of edges $I \subset E$ of ``local isomorphisms'' we want $f_!$ to be a left adjoint of $f^*$.

	\item If $f$ belongs to a class of edges $P \subset E$ of ``proper maps'' we want $f_! = f_*$.
\end{itemize}
Assuming that every morphism in $E$ can be factored into elements of $I$ and $P$ the above two rules determine $f_!$ for all $f \in E$. However, it is not at all clear if the thus constructed functor $f_!$ is well-defined (i.e. independent of the choice of composition of edges in $I$ and $P$) and satisfies the required compatibilities for a 6-functor formalism; even worse, showing the desired compatibilities on a 1-catorical level is a priori not enough to construct all the required higher homotopies constituting a 6-functor formalism. The following result provides a list of conditions one needs to check in order to make the above construction work.

\begin{definition} \label{def:sixfun-suitable-decomposition-of-E}
Let $(\mathcal C, E)$ be a geometric setup such that $\mathcal C$ admits pullbacks. A \emph{suitable decomposition of $E$} is a pair $I, P \subset E$ of subsets satisfying the following properties:
\begin{defenum}
	\item Every $f \in E$ is of the form $f = p \comp i$ for some $i \in I$ and some $p \in P$.
	\item Every morphism $f \in I \isect P$ is $n$-truncated for some $n \ge -2$ (which may depend on $f$).
	\item $I$ and $P$ contain all identity morphisms and are stable under pullback.
	\item Given two morphisms $f\colon Y \to X$ and $g\colon Z \to Y$ in $\mathcal C$ such that $f \in I$ (resp. $f \in P$) then we have $g \in I$ (resp. $g \in P$) if and only if $f \comp g \in I$ (resp. $f \comp g \in P$).
\end{defenum}
\end{definition}

\begin{proposition} \label{rslt:construct-6-functor-formalism-out-of-I-P}
Let $(\mathcal C, E)$ be a geometric setup such that $\mathcal C$ admits pullbacks, let $I, P \subset E$ be a suitable decomposition of $E$ and let $\D\colon \mathcal C^\opp \to \infcatinf^\tensor$ be a functor. For every $f\colon Y \to X$ denote $f^* := \D(f)\colon \D(X) \to \D(Y)$. Assume that the following conditions are satisfied:
\begin{enumerate}[(a)]
	\item For every $[j\colon U \to X] \in I$ the following is true:
	\begin{itemize}
		\item The functor $j^*$ admits a left adjoint $j_!\colon \D(Y) \to \D(X)$.
		\item For every map $g\colon X' \to X$ in $\mathcal C$ with pullbacks $g'\colon U' \to U$ and $j'\colon U' \to X'$, the natural map $j'_! g'^* \isoto g^* j_!$ is an isomorphism of functors $\D(U) \to \D(X')$.
		\item For all $M \in \D(X)$ and $N \in \D(Y)$ the natural morphism $j_!(N \tensor j^* M) \isoto j_! N \tensor M$ is an isomorphism.
	\end{itemize}

	\item For every $[f\colon Y \to X] \in P$ the following is true:
	\begin{itemize}
		\item The functor $f^*$ admits a right adjoint $f_*\colon \D(Y) \to \D(X)$.
		\item For every map $g\colon X' \to X$ in $\mathcal C$ with pullbacks $g'\colon Y' \to Y$ and $f'\colon Y' \to X'$, the natural map $g^* f_* \isoto f'_* g'^*$ is an isomorphism of functors $\D(Y) \to \D(X')$.
		\item For all $M \in \D(X)$ and $N \in \D(Y)$ the natural morphism $f_* N \tensor M \isoto f_*(N \tensor f^* M)$ is an isomorphism.
	\end{itemize}

	\item For every cartesian diagram
	\begin{center}\begin{tikzcd}
		U' \arrow[r,"j'"] \arrow[d,"f'"] & X' \arrow[d,"f"]\\
		U \arrow[r,"j"] & X
	\end{tikzcd}\end{center}
	in $\mathcal C$ such that $j \in I$ and $f \in P$, the natural map $j_! f'_* \isoto f_* j'_!$ is an isomorphism of functors $\D(U') \to \D(X)$.
\end{enumerate}
Then $\D$ can be extended to a pre-6-functor formalism
\begin{align*}
	\D\colon \Corr(\mathcal C)_{E,all} \to \infcatinf
\end{align*}
such that for all $f \in P$ we have $f_! = f_*$ and for all $j \in I$ the functor $j_!$ is left adjoint to $j^*$. Suppose that additionally the following conditions are satisfied:
\begin{enumerate}[(a)]
	\setcounter{enumi}{3}
	\item For every $X \in \mathcal C$, the symmetric monoidal $\infty$-category $\D(X)$ is closed.

	\item For every map $f\colon Y \to X$ in $\mathcal C$ the functor $f^*$ admits a right adjoint $f_*\colon \D(Y) \to \D(X)$.

	\item For every $[f\colon Y \to X] \in P$ the functor $f_*$ admits a right adjoint $f^!\colon \D(X) \to \D(Y)$.
\end{enumerate}
Then the above pre-6-functor formalism $\D$ is a 6-functor formalism.
\end{proposition}
\begin{proof}
The second part of the claim is clear, so we only need to prove the first part. By \cite[Proposition 2.4.1.7]{lurie-higher-algebra} the desired map of $\infty$-operads can equivalently be described by a lax Cartesian structure (see \cite[Definition 2.4.1.1]{lurie-higher-algebra}) $\Corr(\mathcal C)_{E,all}^\tensor \to \infcatinf$. Using notation as in \cite[\S1.3]{enhanced-six-operations}, we deduce from \cite[Example 4.30]{glueing-restricted-nerves} that there is a categorical equivalence of simplicial sets
\begin{align*}
	\delta^*_{2,\{2\}}((\mathcal C^\opp)^{\amalg,\opp})^\cart_{E^-,all} \to \Corr(\mathcal C)^\tensor_{E,all}
\end{align*}
By \cite[Proposition 1.2.7.3]{lurie-higher-topos-theory} it is therefore enough to construct a functor
\begin{align*}
	\delta^*_{2,\{2\}}((\mathcal C^\opp)^{\amalg,\opp})^\cart_{E^-,all} \to \infcatinf
\end{align*}
satisfying the required properties. To simplify notation let us write $\mathcal C_\amalg := (\mathcal C^\opp)^{\amalg,\opp}$ in the following. Using the assumptions on $I$ and $P$, it follows from \cite[Theorem 5.4]{glueing-restricted-nerves} that there is a categorical equivalence
\begin{align*}
	\delta^*_{3,\{3\}}(\mathcal C_\amalg)_{I^-,P^-,all} \to \delta^*_{2,\{2\}}(\mathcal C_\amalg)_{E^-,all}.
\end{align*}
We are thus reduced to constructing a certain functor
\begin{align*}
	F\colon \delta^*_{3,\{3\}}(\mathcal C_\amalg)_{I^-,P^-,all} \to \infcatinf.
\end{align*}
The simplicial set $\delta^*_{3,\{3\}}(\mathcal C_\amalg)_{I^-,P^-,all}$ can roughly be described as follows: Its objects are the objects of $\mathcal C_\amalg$ and its morphisms are cubes made out of cartesian squares, such that all the morphisms in vertical, resp. horizontal direction lie in $I^-$, resp. $P^-$. Moreover, the subscript $\{3\}$ indicates that the morphisms going in the direction of the third dimension of the cube are taken in the opposite direction. The functor $F$ should map such a cube to the morphism which is obtained as the composition of the morphisms $j_!$ along vertical maps $j$ of the cube, $f_*$ along horizontal maps $f$ of the cube and $g^*$ along maps $g$ in the third dimension of the cube.

In order to construct $F$, we start with the composed functor
\begin{align*}
	F_0\colon \delta^*_{3,\{1,2,3\}}(\mathcal C_\amalg)_{I^-,P^-,all} \to (\mathcal C^\opp)^\amalg \to \infcatinf,
\end{align*}
where the first functor is the diagonal functor and the second functor is induced by $\D$ via \cite[Theorem 2.4.3.18]{lurie-higher-algebra}. Here $\delta^*_{3,\{1,2,3\}}(\mathcal C_\amalg)_{I^-,P^-,all}$ is the simplicial set whose objects are the objects of $\mathcal C_\amalg$ and whose morphisms are cubes as above, but with the morphisms in all directions being inverted. The functor $F_0$ maps such a cube to the composition of pullback functors along all three dimensions. Thus in order to get the desired functor $F$ out of $F_0$, we need to invert the directions of the vertical and the horizontal maps in the cubes and replace their image under $F_0$ by their left adjoints resp. right adjoints. This can be done by applying \cite[Lemma 1.4.4]{enhanced-six-operations} twice, as follows.

We first invert the direction of the maps in dimension $1$ (i.e. the vertical maps) by passing to left adjoints. By \cite[Lemma 1.4.4]{enhanced-six-operations} we need to verify that certain diagrams with respect to the dimensions $(1, 2)$ and $(1, 3)$ are adjointable. The former is a special case of the latter, so we only need to check the condition for $(1, 3)$. As in the proof of \cite[Lemma 3.2.5]{enhanced-six-operations} the condition boils down precisely to proper base-change and projection formula for the functors $j_!$ associated to maps $j \in I$; this holds by assumption (a).

We now invert the direction of the maps in dimension $2$ (i.e. the horizontal maps) by passing to right adjoints. We again apply \cite[Lemma 1.4.4]{enhanced-six-operations} which requires us to check certain compatibilities with respect to the dimensions $(2,1)$ and $(2,3)$. The latter condition reduces to proper base-change and the projection formula for the functors $f_*$ associated to maps $f \in P$; this is assumption (b). The former condition reduces to assumption (c).
\end{proof}

Using \cref{rslt:construct-6-functor-formalism-out-of-I-P} one can construct a basic 6-functor formalism in many cases. However, often one wants to have a class of edges $E$ which may not necessarily be decomposed into local isomorphisms $I$ and proper maps $P$. Among others, the following situations may occur:
\begin{enumerate}[(1)]
	\item Given a map $f\colon Y \to X$ in $E$, a decomposition of the form $f = p \comp i$ with $p$ proper and $i$ a local isomorphism is only possible \emph{locally on $Y$}.

	\item Given a map $f\colon Y \to X$ in $E$, it may be the case that $f$ only becomes of a ``nice'' form after pullback along a cover $X' \surjto X$ of $X$.
\end{enumerate}
In the following we provide general results to perform extensions of type (1) and (2) above. We start with extensions of type (1). The general abstract condition one needs to verify is given by the following result. We will afterwards single out two special cases: Firstly, under mild assumptions one can always extend $\D$ from $(\mathcal C, E)$ to the setup $(\mathcal C, E')$ where $E'$ consists of the edges which are of the form $\bigdunion_i Y_i \to X$ with all $Y_i \in X$ belonging to $E$ (see \cref{rslt:extend-6-functor-to-disjoint-unions-on-source}). Secondly, one can extend $\D$ from $(\mathcal C, E)$ to $(\mathcal C, E')$ if every morphism $f\colon Y \to X$ in $E'$ can be covered by a morphism $g\colon Z \to Y$ such that $f \comp g \in E$ and $\D^!$ satisfies descent along $g$ (see \cref{rslt:extend-6-functors-locally-on-source}).

\begin{lemma} \label{rslt:6-functor-left-kan-extension-criterion}
Let $\mathcal C$ be an $\infty$-category equipped with two geometric setups $(\mathcal C, E)$ and $(\mathcal C, E')$ such that $E \subset E'$. Let $\D\colon \Corr(\mathcal C)_{E,all} \to \infcatinf$ be a 6-functor formalism satisfying the following properties:
\begin{enumerate}[(a)]
	\item For every $X \in \mathcal C$, $\D(X)$ is presentable.

	\item For every $[f\colon Y \to X] \in E'$ let $\mathcal I_f \subset \mathcal C_{/Y}$ be the subcategory whose objects are the morphisms $g\colon Z \to Y$ such that $g, f \comp g \in E$ and whose morphisms are the $Y$-morphisms $Z' \to Z$ which lie in $E$. Then the natural functor
	\begin{align*}
		\D^!(Y) \isoto \varprojlim_{Z \in \mathcal I_f^\opp} \D^!(Z)
	\end{align*}
	is an equivalence.
\end{enumerate}
Then $\D$ extends uniquely to a 6-functor formalism
\begin{align*}
	\D\colon \Corr(\mathcal C)_{E',all} \to \infcatinf.
\end{align*}
\end{lemma}
\begin{proof}
By condition (a) and \cite[Remark 4.8.1.9]{lurie-higher-algebra} we can view $\D$ as a map of $\infty$-operads $\Corr(\mathcal C)_{E,all}^\tensor \to (\catPrL)^\tensor$, where $\catPrL$ denotes the $\infty$-category of presentable $\infty$-categories and colimit preserving functors (equipped with the tensor product from \cite[Proposition 4.8.1.15]{lurie-higher-algebra}). By applying the operadic left Kan extensions result \cite[Proposition 3.1.3.3]{lurie-higher-algebra} to $\mathcal A^\tensor := \Corr(\mathcal C)_{E,all}^\tensor$, $\mathcal B^\tensor := \Corr(\mathcal C)_{E',all}^\tensor$ and $\mathcal C^\tensor := (\catPrL)^\tensor$ we see that in order to construct the desired extension of $\D$ we need to check the following: For every $X \in \mathcal C$ the diagram
\begin{align*}
	\D\colon \mathcal K^\triangleright := ((\Corr(\mathcal C)_{E,all}^\tensor)_\act)^\triangleright_{/X} \to (\catPrL)^\tensor
\end{align*}
is an operadic colimit diagram. Here the $\infty$-category $\mathcal K$ has as objects the active morphisms $Y_\bullet \to X$ in $\Corr(\mathcal C)_{E',all}^\tensor$ and as morphisms the active $X$-morphisms in $\Corr(\mathcal C)_{E,all}^\tensor$. Every object $Y_\bullet \to X$ in $\mathcal K$ has the form
\begin{center}\begin{tikzcd}
	(Y_i)_{1 \le i \le n} & \arrow[l] Y' \arrow[d]\\
	& X
\end{tikzcd}\end{center}
where $Y'$ lies in $\mathcal C$ (i.e. in the fiber over $\langle 1 \rangle$), the map $Y' \to Y$ lies in $E'$ and the map $Y_\bullet \from Y'$ lies over the active morphism $\langle n \rangle \to \langle 1 \rangle$ in $\catFinAst$. We let $\mathcal K' \subset \mathcal K$ denote the full subcategory spanned by those objects $Y_\bullet \to X$ as above where $n = 1$ and $Y = Y'$. We claim that the inclusion $\mathcal K' \injto \mathcal K$ is cofinal. Namely, by \cite[Theorem 4.1.3.1]{lurie-higher-topos-theory} we have to check the following: Given any morphism $f\colon Y_\bullet \to X$ as above the $\infty$-category $\mathcal K'_{f/}$ is weakly contractible. But note that $Y' = Y' \to X$ is an initial object of $\mathcal K'_{f/}$, so that this $\infty$-category is indeed weakly contractible. All in all this means that the operadic colimit over $\mathcal K$ is the same as the operadic colimit over $\mathcal K'$. Since the tensor product on $\catPrL$ is compatible with colimits (see \cite[Remark 4.8.1.24]{lurie-higher-algebra}) it follows from \cite[Example 3.1.1.7]{lurie-higher-algebra} that the above claim about an operadic colimit reduces to showing that the diagram $\D\colon \mathcal K'^\triangleright \to \catPrL$ is a colimit diagram.
One checks that there is a natural equivalence $\mathcal K' = \mathcal C_E \cprod_{\mathcal C_{E'}} (\mathcal C_{E'})_{/X}$, where the right-hand side is the $\infty$-category whose objects are morphisms $Y \to X$ in $E'$ and whose morphisms are the $X$-morphisms in $E$. We have thus reduced the claim to showing that the diagram
\begin{align*}
	\D_!\colon (\mathcal C_E \cprod_{\mathcal C_{E'}} (\mathcal C_{E'})_{/X})^\triangleright \to \catPrL
\end{align*}
is a colimit diagram. Consider the full subcategory $(\mathcal C_E)_{/X} \subset \mathcal K'$. This subcategory has the final object $\id_X$, so it is enough to show that the colimit of $\D$ over $\mathcal K'$ is the same as the colimit of $\D_!$ over $(\mathcal C_E)_{/X}$. This follows if we can show that the diagram $\D\colon \mathcal K' \to \catPrL$ is a left Kan extension of the diagram $\D_!\colon (\mathcal C_E)_{/X} \to \catPrL$. To prove that $\D_!$ satisfies this property, fix some $f \in \mathcal K'$, i.e. $f$ is a morphism $f\colon Y \to X$ in $E'$. From the general computation of Kan extensions one sees immediately that with the notation as in the claim, we need to show that the diagram
\begin{align*}
	\D_!\colon \mathcal I_f^\triangleright \to \catPrL
\end{align*}
is a colimit diagram. By passing to right adjoints, we can equivalently show that the diagram $\D^!\colon (\mathcal I_f^\opp)^\triangleleft \to \catPrR$ is a limit diagram. Since the inclusion $\catPrL \injto \infcatinf$ preserves limits (see \cite[Theorem 5.5.3.18]{lurie-higher-topos-theory}) we arrive at assumption (b).
\end{proof}

\begin{proposition} \label{rslt:extend-6-functor-to-disjoint-unions-on-source}
Let $\D\colon \Corr(\mathcal C)_{E,all} \to \infcatinf$ be a 6-functor formalism on a geometric setup $(\mathcal C, E)$ such that the following conditions are satisfied (with an optional fixed regular cardinal $\kappa$):
\begin{enumerate}[(a)]
	\item For every $X \in \mathcal C$, $\D(X)$ is presentable.

	\item The $\infty$-category $\mathcal C$ admits fiber products and ($\kappa$-)small coproducts which are
	\begin{itemize}
	 	\item disjoint, i.e. for $X, Y \in \mathcal C$ the object $X \cprod_{X \dunion Y} Y$ is an initial object of $\mathcal C$, and
	 	\item universal, i.e. if $(X_i)_{i\in I}$ is a ($\kappa$-)small family in $\mathcal C$ and $f\colon \bigdunion_i X_i \to X$ is a map, then the pullback of $f$ along any map $g\colon Y \to X$ is of the form $\bigdunion_i Y_i \to Y$, where $Y_i \to Y$ is the pullback of $X_i \to X$ along $g$.
	 \end{itemize}
	 Moreover, for any $X, Y \in \mathcal C$ the map $X \to X \dunion Y$ lies in $E$.

	 \item The functor $\D^!\colon \mathcal C_E^\opp \to \infcatinf$ preserves ($\kappa$-)small products.
\end{enumerate}
Let $E'$ be the collection of maps in $\mathcal C$ which can be written as $\bigdunion_i Y_i \to X$ for some ($\kappa$-)small collection of objects $Y_i \in \mathcal C$ such that all the maps $Y_i \to X$ lie in $E$. Then $(\mathcal C, E')$ is a geometric setup and $\D$ extends uniquely to a 6-functor formalism
\begin{align*}
	\D\colon \Corr(\mathcal C)_{E',all} \to \infcatinf
\end{align*}
\end{proposition}
\begin{proof}
It is straightforward to see that $(\mathcal C, E')$ is a geometric setup. To see that $\D$ extends we need to check condition (b) of \cref{rslt:6-functor-left-kan-extension-criterion}, so let $f\colon Y = \bigdunion_i Y_i \to X$ be a given morphism in $E'$ and let $\mathcal I_f$ be as in \cref{rslt:6-functor-left-kan-extension-criterion}. Consider the full subcategory $\mathcal J \subset \mathcal I_f$ consisting of those maps $Z \to Y$ which factor over some $Y_i$ and where $Z$ is not initial. We claim that the diagram $\D^!\colon \mathcal I_f \to \infcatinf$ is the right Kan extension of its restriction to $\mathcal J$. To see this, fix any $Z \in \mathcal I_f$. We need to verify that the natural functor
\begin{align*}
	\D^!(Z) \isoto \varprojlim_{Z' \in \mathcal J_{/Z}^\opp} \D^!(Z')
\end{align*}
is an equivalence. Let us denote $Z_i := Y_i \cprod_Y Z$, so that $Z = \bigdunion_i Z_i$. Any $g\colon Z' \to Z$ factors over some $Y_i$ and hence also over some $Z_i$ (for a unique $i$). Thus, if we denote by $\mathcal J_{/Z,i} \subset \mathcal J_{/Z}$ the full subcategory spanned by those maps $g\colon Z' \to Z$ which factor over $Z_i$ then $\mathcal J_{/Z} = \bigdunion_i \mathcal J_{/Z,i}$. Moreover, each of the $\infty$-categories $\mathcal J_{/Z,i}$ has a final object, namely $Z_i$. It follows that
\begin{align*}
	\varprojlim_{Z' \in \mathcal J_{/Z}^\opp} \D^!(Z') = \prod_i \varprojlim_{Z' \in \mathcal J_{/Z,i}^\opp} \D^!(Z') = \prod_i \D^!(Z_i) = \D^!(Z)
\end{align*}
by condition (c). This proves that $\D^!\colon \mathcal I_f \to \infcatinf$ is indeed the right Kan extension of its restriction to $\mathcal J$, hence we are reduced to showing that the natural functor
\begin{align*}
	\D^!(Y) \isoto \varprojlim_{Z \in \mathcal J^\opp} \D^!(Z)
\end{align*}
is an equivalence. But note that $\mathcal J = \bigdunion_i \mathcal J_i$, where $\mathcal J_i \subset \mathcal J$ is the full subcategory spanned by those $Z \to Y$ which factor over $Y_i$. Thus we can argue as before to show that $\D^!$ does indeed satisfy the above limit property.
\end{proof}

\begin{lemma} \label{rslt:6-functor-cech-nerve-is-cofinal}
Let $\mathcal C$ be an $\infty$-category and let $X \in \mathcal C$ be an object satisfying the following properties:
\begin{enumerate}[(a)]
	\item For every $Y \in \mathcal C$ there is a morphism $Y \to X$.
	\item All the products $X^{\cprod n} = X \cprod \dots \cprod X$ exist in $\mathcal C$.
\end{enumerate}
Then the right Kan extension along $\{ [0] \} \injto \Delta^\opp$ produces a simplicial object $X_\bullet\colon \Delta^\opp \to \mathcal C$ with $X_n = X^{\cprod n+1}$ and the functor
\begin{align*}
	X_\bullet\colon \Delta^\opp \to \mathcal C
\end{align*}
is cofinal.
\end{lemma}
\begin{proof}
The existence of the simplicial object $X_\bullet$ follows easily from (b). To show that the functor $\Delta^\opp \to \mathcal C$ is cofinal, we argue as in the proof of \cite[Proposition A.3.3.1]{lurie-spectral-algebraic-geometry}. By \cite[Theorem 4.1.3.1]{lurie-higher-topos-theory} it suffices to show that for every object $Y \in \mathcal C$ the $\infty$-category $\mathcal X := \Delta^\opp \cprod_{\mathcal C} \mathcal C_{Y/}$ is weakly contractible. The projection $\mathcal X \to \Delta^\opp$ is a left fibration and thus classified by a functor $\chi\colon \Delta^\opp \to \Ani$ -- explicitly we have $\chi([n]) = \Hom(Y, X^{\cprod n+1})$. By \cite[Proposition 3.3.4.5]{lurie-higher-topos-theory} it is enough to show that $\varinjlim \chi$ is contractible. But $\chi$ is the Čech nerve of the map $q\colon \Hom(Y, X) \to *$ in $\Ani$, so since $\Ani$ is an $\infty$-topos it suffices to show that $q$ is an effective epimorphism. This is equivalent to saying that $\Hom(Y, X)$ is non-empty, which is true by assumption (a).
\end{proof}

\begin{proposition} \label{rslt:extend-6-functors-locally-on-source}
Let $\mathcal C$ be an $\infty$-category equipped with two geometric setups $(\mathcal C, E)$ and $(\mathcal C, E')$ such that $E \subset E'$, let $\D\colon \Corr(\mathcal C)_{E,all} \to \infcatinf$ be a 6-functor formalism and let $S \subset E$ be a subset of ``special covers'' such that the following conditions are satisfied:
\begin{enumerate}[(a)]
	\item For every $X \in \mathcal C$, $\D(X)$ is presentable.

	\item For every map $f\colon Y \to X$ in $S$ with associated Čech nerve $Y_\bullet \to X$, the natural functor
	\begin{align*}
		\D^!(X) \isoto \varprojlim_{n\in\Delta} \D^!(Y_n)
	\end{align*}
	is an equivalence.

	\item For every map $f\colon Y \to X$ in $E'$ there is a map $g\colon Z \to Y$ in $S$ such that $f \comp g \in E$.

	\item Pullbacks of edges in $S$ exist in $\mathcal C_E$ and remain in $S$. Moreover the inclusion $\mathcal C_E \injto \mathcal C$ preserves these pullbacks.
\end{enumerate}
Then $\D$ extends uniquely to a 6-functor formalism
\begin{align*}
	\D\colon \Corr(\mathcal C)_{E',all} \to \infcatinf.
\end{align*}
\end{proposition}
\begin{proof}
We need to verify condition (b) of \cref{rslt:6-functor-left-kan-extension-criterion}, so let $[f\colon Y \to X] \in E'$ be given and let $\mathcal I_f \subset \mathcal C_{/Y}$ be as in \cref{rslt:6-functor-left-kan-extension-criterion}. By (c) there is some $[g\colon Z \to Y] \in \mathcal I_f$ such that $g \in S$. Let $\mathcal J \subset \mathcal I_f$ be the full subcategory spanned by those maps $Z' \to Y$ which factor as $Z' \to Z \xto{g} Y$, where the first map is in $E$. Let $Z_\bullet \to Y$ be the Čech nerve of $g$. Then by \cref{rslt:6-functor-cech-nerve-is-cofinal} and (b) we have
\begin{align*}
	\varprojlim_{Z' \in \mathcal J^\opp} \D^!(Z') = \varprojlim_{n\in\Delta} \D^!(Z_n) = \D^!(Y).
\end{align*}
It is thus enough to show that the diagram $\D^!\colon \mathcal I_f^\opp \to \infcatinf$ is the right Kan extension of its restriction to $\mathcal J$. To verify this, fix any $W \in \mathcal I_f$; we need to show that the natural functor
\begin{align*}
	\D^!(W) \isoto \varprojlim_{W' \in \mathcal J_{/W}^\opp} \D^!(W')
\end{align*}
is an equivalence. But note that by (d) every $[W' \to W] \in \mathcal J_{/W}$ factors over $W \cprod_Y Z$ and the map $W \cprod_Y Z \to W$ lies in $S$. Thus the above limit for $\D^!(W)$ follows in the same way as before.
\end{proof}

\begin{remarks}
\begin{remarksenum}
	\item The idea behind \cref{rslt:extend-6-functors-locally-on-source} is as follows. Given everything as in that result we fix some edge $f\colon Y \to X$ in $E'$. In order to define the functor $f_!\colon \D(Y) \to \D(X)$, pick some edge $g\colon Z \to Y$ in $S$ such that $f \comp g \in E$ (using assumption (c)) and let $g_\bullet\colon Z_\bullet \to Y$ be the associated Čech nerve. Then each $f \comp g_n$ lies in $E$ and we define
	\begin{align*}
		f_! := \varinjlim_{n\in\Delta} (f \comp g_n)_! g_n^!.
	\end{align*}
	This definition is sensible because by (b) we have $\id = \varinjlim_{n\in\Delta} g_{n!} g_n^!$. In order to show that this definition is independent of $g$ we need assumption (d).

	\item Condition (d) in \cref{rslt:extend-6-functors-locally-on-source} is implied by the following condition: $S$ is stable under composition and pullbacks in $\mathcal C$ and for all composable maps $f, g$ in $\mathcal C$ with $f, f \comp g \in E$ we have $g \in E$.
\end{remarksenum}
\end{remarks}

We now consider extensions of 6-functor formalisms of type (2), i.e. we want to extend a 6-functor formalism $\D$ on $(\mathcal C, E)$ to a new 6-functor formalism on $(\mathcal C', E')$ where $\mathcal C \subset \mathcal C'$ is a full subcategory and every edge $f\colon Y \to X$ in $E'$ is locally on $X$ an edge in $E$. This can be done as follows:

\begin{proposition} \label{rslt:extend-6-functors-locally-on-target}
Let $(\mathcal C, E)$ and $(\mathcal C', E')$ be geometric setups such that $\mathcal C \subset \mathcal C'$ is a full subcategory and $E \subset E'$. Assume that the following conditions are satisfied:
\begin{enumerate}[(a)]
	\item The inclusion $\mathcal C \injto \mathcal C'$ preserves pullbacks of edges in $E$.
	\item If $[f\colon Y \to X] \in E'$ satisfies $X \in \mathcal C$ then $Y \in \mathcal C$ and $f \in E$.
\end{enumerate}
Then every pre-6-functor formalism $\D\colon \Corr(\mathcal C)_{E,all} \to \infcatinf$ extends uniquely to a pre-6-functor formalism
\begin{align*}
	\D'\colon \Corr(\mathcal C')_{E',all} \to \infcatinf
\end{align*}
such that for every $X \in \mathcal C'$ we have
\begin{align*}
	\D'^*(X) = \varprojlim_{Y \in \mathcal C_{/X}^\opp} \D^*(Y).
\end{align*}
Moreover, if $\D$ is a 6-functor formalism such that $\D(X)$ is presentable for all $X \in \mathcal C$ then $\D'$ is a 6-functor formalism.
\end{proposition}
\begin{proof}
The second part of the claim is easy to see because if all $\D(X)$ are presentable then also all $\D'(X)$ are presentable (using the explicit description of $\D'(X)$ in the claim) and one checks easily that all the functors $f^*$, $f_!$ and $\tensor$ preserve colimits and therefore admit right adjoints.

We now prove the first part of the claim, i.e. the existence of $\D'$. In the following we abbreviate $\mathcal X^\tensor := \Corr(\mathcal C)_{E,all}^\tensor$ and $\mathcal X'^\tensor := \Corr(\mathcal C')_{E',all}^\tensor$ and similarly without the superscript $(-)^\tensor$. Suppose we are given a pre-6-functor formalism on $(\mathcal C, E)$, which we can view as a lax Cartesian structure $\D\colon \mathcal X^\tensor \to \infcatinf$ by \cite[Proposition 2.4.1.7]{lurie-higher-algebra}. We now define the functor
\begin{align*}
	\D'\colon \mathcal X'^\tensor \to \infcatinf
\end{align*}
to be the right Kan extension of $\D$ along the obvious functor $\alpha^\tensor\colon \mathcal X^\tensor \to \mathcal X'^\tensor$ ($\alpha$ exists by condition (a)). We need to check that $\D'$ is indeed an extension of $\D$ and that it is a lax Cartesian structure. In order to do that, we will compute $\D'(X_\bullet)$ for every $X_\bullet = (X_i)_{1\le i \le n} \in \mathcal X'^\tensor$. Fixing such an $X_\bullet$, let $\mathcal K := \mathcal X^\tensor_{X_\bullet/}$, so that $\D'(X_\bullet) = \varprojlim_{Y_\bullet \in \mathcal K} \D(Y_\bullet)$. Every object in $\mathcal K$ is of the form
\begin{center}\begin{tikzcd}
	(X_i)_{1 \le i \le n} & \arrow[l] (Y'_j)_{1 \le j \le m} \arrow[d]\\
	& (Y_j)_{1 \le j \le m}
\end{tikzcd}\end{center}
such that $Y_j \in \mathcal C$ for all $j$, the vertical map lies over the identity $\langle m \rangle \to \langle m \rangle$ in $\catFinAst$ and all the maps $Y'_j \to Y_j$ lie in $E'$. Consider the full subcategory $\mathcal K' \subset \mathcal K$ consisting of those objects where the above vertical map is degenerate. We claim that the inclusion $\mathcal K'^\opp \injto \mathcal K^\opp$ is cofinal, so that any limit over $\mathcal K$ is the same as the limit over the restriction to $\mathcal K'$. To prove this cofinality we employ \cite[Theorem 4.1.3.1]{lurie-higher-topos-theory} which reduces the claim to showing that for every $Y_\bullet \in \mathcal K'$ the $\infty$-category $\mathcal K'_{/Y_\bullet}$ is contractible. The object $Y_\bullet$ comes equipped with maps $X_\bullet \from Y'_\bullet \to Y_\bullet$ as above. By condition (b) we have $Y'_j \in \mathcal C$ for all $j$ and each of the maps $Y'_j \to Y_j$ lies in $E$. In particular there is a morphism $Y'_\bullet \to Y_\bullet$ in $\mathcal K$ and one checks that this morphism is a final object of $\mathcal K'_{/Y_\bullet}$, proving that this $\infty$-category is indeed weakly contractible. Note furthermore that $\mathcal K' = (\mathcal C^\opp)^\amalg_{X_\bullet/}$, so that we deduce
\begin{align*}
	\D'(X_\bullet) = \varprojlim_{Y_\bullet \in \mathcal K} \D(Y_\bullet) = \varprojlim_{Y_\bullet \in (\mathcal C^\opp)^\amalg_{X_\bullet/}} \D(Y_\bullet).
\end{align*}
Now let $\mathcal K'' \subset (\mathcal C^\opp)^\amalg_{X_\bullet/}$ be the full subcategory spanned by those morphisms $X_\bullet \from Y_\bullet$ which lie over the identity $\langle n \rangle$ in $\catFinAst$. One checks easily that the inclusion $\mathcal K''^\opp \injto ((\mathcal C^\opp)^\amalg_{X_\bullet/})^\opp$ is cofinal and thus
\begin{align*}
	\D'(X_\bullet) = \varprojlim_{Y_\bullet \in \mathcal K''} \D(Y_\bullet) = \varprojlim_{Y_\bullet \in \mathcal K''} \prod_i \D(Y_i) = \prod_i \varprojlim_{Y \in \mathcal C^\opp_{/X_i}} \D(Y) = \prod_i \D'(X_i),
\end{align*}
as desired.
\end{proof}

\clearpage
\bibliography{bibliography}
\addcontentsline{toc}{part}{References}

\end{document}